\documentclass{amsart}

\usepackage{comment}
\usepackage{amsmath}
\usepackage{amssymb}
\usepackage{etoolbox} 
\usepackage[pdfauthor={Richard J. Mathar},colorlinks=true,citecolor=blue]{hyperref}

\def\Li{\mathop{\mathrm{Li}}\nolimits}
\def\li{\mathop{\mathrm{li}}\nolimits}

\def\Cl{\mathop{\mathrm{Cl}}\nolimits}
\def\Res{\mathop{\mathrm{Res}}\nolimits}
\def\arctanh{\mathop{\mathrm{arctanh}}\nolimits}
\def\csch{\mathop{\mathrm{csch}}\nolimits}
\def\sech{\mathop{\mathrm{sech}}\nolimits}
\def\sh{\mathop{\mathrm{sh}}\nolimits}

\def\Arsh{\mathop{\mathrm{Arsh}}\nolimits}
\def\arcsinh{\mathop{\mathrm{arcinh}}\nolimits}
\def\erfi{\mathop{\mathrm{erfi}}\nolimits}
\def\erf{\mathop{\mathrm{erf}}\nolimits}
\def\Shi{\mathop{\mathrm{Shi}}\nolimits}
\def\Ei{\mathop{\mathrm{Ei}}\nolimits}
\def\Lif{\mathop{\mathrm{Lif}}\nolimits}
\def\sgn{\mathop{\mathrm{sgn}}\nolimits}

\numberwithin{equation}{section}

\newcommand{\be}{\begin{equation}}
\newcommand{\ee}{\end{equation}}

\usepackage{orcidlink}
\begin{document}

\title{Yet Another Table of Integrals}

\author{Richard J. Mathar \orcidlink{0000-0001-6017-6540}}
\urladdr{https://mathar.www3.mpia.de}
\address{Max Planck Institute for Astronomy, K\"onigstuhl 17, 69117 Heidelberg, Germany}


\date{\today}

\begin{abstract}
This collection of sums and integrals has been harvested from the mathematical
and physical literature in unstructured ways.
Its main use is backtracking the original sources whenever an integral
of the reader's application
resembles one of the items in the collection.
[arXiv:1207.5845]

\end{abstract}

\maketitle

\section*{Introduction}
Dealing with the analysis of real numbers in the physical sciences
shows a strange attraction towards integrals.
Closed-form integration
beats numerical integration, and often adaptive series expansion helps
to dissect cumbersome integral kernels to digestable pieces.

The current table started as a incoherent list of bookmarks pointing to
``interesting'' formulas that complement or correct the Gradstein-Rhyshik tables \cite{GR},
\url{http://www.mathtable.com/gr/},
or the Dieckmann tables \cite{Dieckmann} \url{https://www-elsa.physik.uni-bonn.de/~dieckman}.
As such it  does not replicate the original sources in full but
is to be merely regarded
as an aid to find places
at which certain forms and classes of
integrals or sums have been targeted.

The notation is generally not harmonized. Stirling numbers
appear in bracketed and indexed notations, and at least two
different meanings of harmonic numbers $H$ with lower and upper indices
are met.

One hint of use:
The list of references
appears \emph{prior} to each formula.

\subsection{Finite Series}

\cite{Brien}
\be
\sum_{k=1}^n aq^{k-1} = \frac{a(q^n-1)}{q-1}=\frac{ql-a}{q-1}
\ee
where l=$aq^{n-1}$ is the last term.

\cite{SpiveyMM78,SnowAM18,SchumacherJIS19}
\be
\sum_{j=1}^n j^k=\frac{1}{n}\left(\rho(n,i)+\sum_{i=1}^n\sigma_{k+1}(i)\right)
=
\frac{B_{k+1}(n+1)-B_{k+1}}{k+1}
\ee
where $\rho(n,k)\equiv \sum_{d=1}^n d^k (n\mod d)$ is a sum over $n\mod d$
multiplied, then summed, over $d^k$, and $\sigma_k(n)=\sum_{d|n}d^k$.

\cite{NetoJIS19b}
\be
\sum_{j=1}^{\lfloor x \rfloor} i^n
=\frac{x^{n+1}}{n+1}+(-)^n\frac{B_{n+1}}{n+1}
+\frac{1}{n+1}\sum_{j=1}^{n+1}(-)^j\binom{n+1}{j}B_j(\{x\}) x^{n-j+1}.
\ee
where $\{.\}$ denotes the fractional part.

\cite{Linarxiv04,DzhumadilJIS13}
\be
\sum_{k=1}^{n-1} k^r=\sum_{k=0}^r \frac{B_k}{k!}\, \frac{r!}{(r-k+1)!}n^{r-k+1}.
\ee

\cite{SpiveyDM307}
\be
\sum_{k=0}^n k^m=
\sum_{j=0}^m\left\{\begin{array}{c}m\\j\end{array}\right\} \binom{n+1}{j+1}j!.
\ee
By telescoping the last term in this formula becomes \cite{SinghJIS19}
\be
n^m=
\sum_{j=0}^m\left\{\begin{array}{c}m\\j\end{array}\right\} \binom{n}{j}j!.
\ee

\cite{BeraGC36}
Let $c_i$  be a sequence of complex numbers and 
$S(m,k)\equiv \frac{1}{k!}\sum_{j=1}^k (-)^{k-j}\binom{k}{j}j^m$ be the
Stirling Numbers Of The First Kind,
then
\be
\sum_{k=1}^m k^\alpha c_k=\sum_{j=1}^m j! S(\alpha,j)\sum_{k=j}^m\binom{k}{j}c_k.
\ee

\cite{AcuMM61}
Let $S_k(n)\equiv \sum_{l=1}^n l^k$, then
\be
\sum_{i=0}^k \binom{k+1}{i}S_in) = (n+1)^k-1.
\ee

\cite{TakacsMC32}
\be
\sum_{0\le j\le n/m}\binom{n-jm}{k}=P(n+m,k,m)-P(r,k,m)
\ee
for $n\ge 0,$ $k\ge 1$, $m\ge 1$ where
\be
P(x,k,m)\equiv \frac{1}{m}\sum_{j=1}^{k+1}\binom{x}{j}A(m,k+1-j)
\ee
with g.f.\
\be
\frac{mx}{(1+x)^m-1}=\sum_{j=0}^\infty A(m,j)x^j.
\ee

\cite{RoyAMM94}
\be
\sum_{k=0}^n\binom{r+k}{k}
=\binom{r+n+1}{n},\quad n=0,1,2,\ldots
\ee

\cite{SpiveyDM307}
\be
\sum_k \binom{n}{2k}k=n2^{n-3},\quad n\ge 2.
\ee

\cite{ZhaoDM281}
\be
\sum_{k=0}^n[t^k](\frac{1}{1-t} f(\frac{t}{1-t})) = \sum_{j=0}^n\binom{n+1}{j+1}f_j \leadsto_{f=t/(1-t)^2}
\sum_{j=1}^n j \binom{n+1}{j+1} = \sum_{k=0}^n k 2^{k-1}.
\ee

\cite{ZhaoDM281}
\be
\sum_{k=j}^{\lfloor n/2 \rfloor} \binom{n}{2k}\binom{k}{j} = 2^{n-2j-1}\frac{n}{n-j}\binom{n-j}{j}
\leadsto
\sum_{k=j}^{\lfloor n/2 \rfloor} [t^k] (\frac{1}{1-t} f(\frac{t}{1-t})) = \sum_{j=0}^{\lfloor n/2\rfloor} 2^{n-2j-1}\frac{n}{n-j}\binom{n-j}{j}f_j.
\ee

\cite{ZhaoDM281}
\be
\sum_{k=j}^{\lfloor n/2 \rfloor} \binom{n+1}{2k+1}\binom{k}{j} = 2^{n-2j}\binom{n-j}{j}.
\leadsto
\sum_{k=j}^{\lfloor n/2 \rfloor} \binom{n+1}{2k+1} [t^k] (\frac{1}{1-t} f(\frac{t}{1-t})) = \sum_{j=0}^{\lfloor n/2\rfloor} 2^{n-2j}\frac{n}{n-j}\binom{n-j}{j}f_j.
\ee
Similar correspondences to $[t^k]$ coefficients by replacing $\binom{k}{j}$ on the LHS of 
the following binomial identities:
\cite{ZhaoDM281}
\be
\sum_{k=j}^n \binom{n-k}{s}\binom{k}{j} = \binom{n+1}{s+j+1}.
\ee

\cite{ZhaoDM281}
\be
\sum_{k=j}^n \binom{s+k}{s}\binom{k}{j} = \frac{n+1-j}{s+1+j}\binom{n+1}{j}\binom{n+1+s}{s}.
\ee

\cite{ZhaoDM281}
\be
\sum_{k=j}^n (-4)^k \binom{n+k}{2k}\binom{k}{j} = (-)^n 2^{2j} \frac{2n+1}{2j+1}\binom{n+j}{2j}.
\ee

\cite{ZhaoDM281}
\be
\sum_{k=j}^n (-4)^k \frac{n}{n+k}\binom{n+k}{2k}\binom{k}{j} = (-)^n 2^{2j} \frac{n}{n+j}\binom{n+j}{2j}.
\ee

\cite{ZhaoDM281}
\be
\sum_{k=j}^{\lfloor n/2 \rfloor} (-)^k \binom{n-k}{k}2^{n-2k}\binom{k}{j} = (-)^j \binom{n+1}{2j+1}.
\ee

\cite{ZhaoDM281}
\be
\sum_{k=j}^n \binom{\alpha}{k} \binom{\beta}{n-k}\binom{k}{j} = \binom{\alpha}{j}\binom{\alpha+\beta-j}{n-j}.
\ee

\cite{ZhaoDM281}
\be
\sum_{k=j}^n (-)^k \binom{n}{k} \binom{2n-k}{n}\binom{k}{j} = (-)^j \binom{n}{j}^2.
\ee

\cite{ZhaoDM281}
\be
\sum_{k=j}^{\lfloor n/2 \rfloor} \binom{n}{2k} \binom{2k}{k}2^{n-2k}\binom{k}{j} = \binom{2n-2j}{n}\binom{n}{j}.
\ee

\cite{ZhaoDM281}
\be
\sum_{k=0}^{\lfloor n/2 \rfloor} \binom{n}{2k}[t^k] \frac{f(t)}{1-t}
=\sum_{j=0}^{\lfloor n/2\rfloor} 2^{n-2j-1}\frac{n}{n-j}\binom{n-j}{j}[t^j] f(\frac{t}{1+t}).
\ee

\cite{ZhaoDM281}
\be
\sum_{k=0}^{\lfloor n/2 \rfloor} \binom{n+1}{2k+1}[t^k] \frac{f(t)}{1-t}
=\sum_{j=0}^{\lfloor n/2\rfloor} 2^{n-2j}\binom{n-j}{j}[t^j] f(\frac{t}{1+t}).
\ee
and similar relations by replacing the binomial terms of the LHS by $[t^k] \frac{f(t)}{1-t}$.

\cite{SpiveyDM307}
\be
\sum_k \binom{n}{2k+1}k=(n-2)2^{n-3},\quad n\ge 2.
\ee

\cite{SpiveyDM307}
\be
\sum_k \binom{n}{2k}k^{\underline m}=n(n-m-1)^{\underline {m-1}}2^{n-2m-1},
\quad n\ge m+1,
\ee
\cite{SpiveyDM307}
\be
\sum_k \binom{n}{2k+1}k^{\underline m}=n(n-m-1)^{\underline {m}}2^{n-2m-1},
\quad n\ge m+1,
\ee
where $k^{\underline m}$ is the falling factorial $k(k-1)(k-2)\cdots (k-m+1)$.

\cite{ZhaoJIS13}
\be
\sum_{n\ge 0}\binom{n+k-1}{n}x^n=\frac{1}{(1-x)^k}.
\ee

\cite{ZhaoJIS13}
\be
\sum_{n\ge 0}\binom{n+k}{n}x^n=\frac{1}{(1-x)^{k+1}}.
\ee

\cite{ZhaoJIS13}
\be
\sum_{n\ge 0}n\binom{n+k}{n}x^n=\frac{(k+1)x}{(1-x)^{k+2}}.
\ee

\cite{SpiveyDM307}
\be
\sum_{k=0}^n\binom{n}{k}k^m=\sum_{j=0}^m\left\{\begin{array}{c}m\\j\end{array}\right\} \binom{n}{j}j!2^{n-j}.
\ee

\cite{SpiveyDM307}
\be
\sum_{k=0}^n\binom{n}{k}(-1)^kk^m=\left\{\begin{array}{c}m\\n\end{array}\right\} (-1)^n n!.
\ee

\cite{SpiveyDM307}
\be
\sum_{k=0}^n\binom{n}{2k}k^m=
n\sum_{j=1}^{\min(m,n-1)}\left\{\begin{array}{c}m\\j\end{array}\right\}
\binom{n-j-1}{j-1}(j-1)!2^{n-2j-1}.
\ee

\cite{SpiveyDM307}
\be
\sum_{k=0}^n\binom{n}{2k+1}k^m=
\sum_{j=1}^{\min(m,n-1)}\left\{\begin{array}{c}m\\j\end{array}\right\}
\binom{n-j-1}{j-1}j!2^{n-2j-1}.
\ee

\cite{KirschenhoferEJC3}
\be
\sum_{0\le k\le N, k\neq K}\binom{N}{k}
(-1)^k\frac{1}{(k-K)^m}
=
\binom{N}{K}(-1)^{K+1}\frac{1}{m!}
Y_m(\ldots,(i-1)![H_{N-K}^{(i)}+(-1)^iH_K^{(i)}],\ldots)
\ee
where
$0\le K\le N$, $Y_m(\ldots,x_i,\ldots)$ are the Bell polynomials
and $H_r^{(i)}=\sum_{j=1}^r j^{-i}$ Harmonic numbers of the
$i$-th order.

\cite{KirschenhoferEJC3}
\be
\sum_{0\le k\le N, k\neq K}\binom{N}{k}
(-1)^k\frac{1}{k-K}
=
\binom{N}{K}(-1)^{K+1}
(H_{N-K}-H_K)
\ee
where
$0\le K\le N$,
$H_r^{(i)}=\sum_{j=1}^r j^{-i}$ Harmonic numbers of the
$i$-th order.

\cite{KirschenhoferEJC3}
\be
\sum_{0\le k\le N}\binom{N}{k}
(-1)^k\frac{1}{k^m}
=
-
\frac{1}{m!}
Y_m(\ldots,(i-1)!H_N^{(i)},\ldots)
\ee
where
$H_r^{(i)}=\sum_{j=1}^r j^{-i}$ Harmonic numbers of the
$i$-th order.

\cite{KirschenhoferEJC3}
\be
\sum_{0\le k\le N}\binom{N}{k}(-1)^k\frac{1}{(k-\xi)^m}
=
\Gamma(-\xi)\frac{\Gamma(N+1)}{\Gamma(N+1-\xi)}
\frac{1}{(m-1)!} Y_{m-1}(\ldots,(i-1)!\zeta_N(i,-\xi),\ldots),
\ee
where $\zeta_N(i,-\zeta)=\sum_{j=0}^N (j-\xi)^{-j}$.

\cite[\S 4.3]{Riordan}
\be
a_n=\sum_{k=0}^n\binom{p+k}{k}b_{n-k}\Leftrightarrow
b_n=\sum_{k=0}^n(-1)^k\binom{p+k}{k}a_{n-k}.
\ee

\cite[\S 4.3]{Riordan}
\be
a_n=\sum_{k=0}^n\binom{p+k}{k}b_{n-qk}\Leftrightarrow
b_n=\sum_{k=0}^n(-1)^k\binom{p+k}{k}a_{n-qk}.
\ee

\cite[\S 4.3]{Riordan}
\be
a_n=\sum_{k=0}^n\binom{2k}{k}b_{n-k}\Leftrightarrow
b_n=\sum_{k=0}^n\frac{1}{1-2k}\binom{2k}{k}a_{n-k}.
\ee

\cite[\S 4.3]{Riordan}
\be
a_n=\sum_{k=0}^n\frac{1}{k+1}\binom{2k}{k}b_{n-k}\Leftrightarrow
b_n=a_n-\sum_{k=1}^n\frac{1}{k}\binom{2k}{k}a_{n-k}.
\ee

\cite{ChuRDCM}
\be
\sum_{k=0}^n(u+vk)\binom{n}{k}(a+bk)^{k-1}(c+b(n-k))^{n-k-1}
=\frac{(a+c)u+avn}{ac(a+c+bn)}(a+c+bn)^n.
\ee
\cite{ChuRDCM}
\be
\sum_{k=0}^n \binom{n}{k}(a+bk)^k(c-bk)^{n-k}=
\sum_{k=0}^n k!\binom{n}{k}b^k(a+c)^{n-k}.
\ee
\cite{ChuRDCM}
\be
\sum_{k=0}^n\binom{n}{k}(a+bk)^{k-2}(c+b(n-k))^{n-k-1}
=\frac{(a+b)(a+c)^2+b^2cn}{a^2c(a+b)}(a+c+bn)^{n-2}.
\ee
\cite{ChuRDCM}
\be
\sum_{k=0}^n \binom{n}{k}(a+bk)^{k-2}(c+b(n-k))^{n-k}=
\frac{(a+b)(a+c)+b^2n}{a^2(a+b)}(a+c+bn)^{n-1}.
\ee
\cite{ChuRDCM}
\be
\sum_{k=0}^n \binom{n}{k}(a+bk)^{k-2}(c+b(n-k))^{n-k+1}=
\frac{(a+b)c+b^2n}{a^2(a+b)}(a+c+bn)^n.
\ee
\cite{ChuRDCM}
\be
\sum_{k=0}^n(u+vk)\binom{n}{k}(a+bk)^{k-1}(c-bk)^{n-k-1}
=
\frac{(a+c)u+n(av-bu)}{a(a+c)(c-bn)}(a+c)^n
\ee
\cite{ChuRDCM}
\be
\sum_{k=0}^n[(a+c)u+k(av-bu)]\binom{n}{k}(a+c)^{k-1}(-c+bk)^{n-k}
=
(a(u+vn)(a+bn)^n.
\ee
\cite{ChuRDCM}
\be
\sum_{k=0}^n\binom{n}{k}(a+bk)^{k-2}(c-bk)^{n-k}
=
\frac{(a+b)(a+c)-abn}{a^2(a+b)(a+c)}(a+c)^n.
\ee
\ldots plus more convoluted Abel-Gould identities.

\cite{ProdingerIPL46}
\be
x_n=a_n+2^{1-n}\sum_{k=0}^n\binom{n}{k}x_k
\Leftrightarrow
x_n=x_0+\sum_{k=0}^n (-1)^k\binom{n}{k}\frac{k\bar x_1+\hat a_k-\bar x_0}{1-2^{1-k}}
\ee
where $\bar x_0\equiv a_0+x_0$, $\bar x_1=a_1+x_0$ and
\be
\hat a_n =\sum_{k=0}^n\binom{n}{k} (-1)^k a_k.
\ee

\cite{ProdingerIPL46}
\be
x_{n+1}=a_{n+1}+2^{1-n}\sum_{k=0}^n\binom{n}{k}x_k
\Leftrightarrow
x_n=-\sum_{k=0}^n (-1)^k\binom{n}{k}\hat x_{k-2}
\ee
where $\bar x_0=0$,
\be
\hat x_n=Q_n\sum_{i=1}^{n+1}(\hat a_i-\hat a_{i+1}-a_1)/Q_{i-1},\quad
Q_n\equiv \prod_{k=1}^n(1-2^{-k}).
\ee

\cite{KilicJCAM201}
\be
u_n=Au_{n-1}-Bu_{n-2},\, 
v_n=Av_{n-1}-Bv_{n-2}, 
\ee
\be
\leadsto
u_{n+1}=\prod_{j=1}^n\left(
A-2i\sqrt{-B}\cos\frac{\pi j}{n+1}
\right)
=
(i\sqrt{-B})^n\frac{\sin\left(
(n+1)\cos^{-1}\frac{-iA}{2\sqrt{-B}}
\right)}{\sin\left(\cos^{-1}\frac{-iA}{2\sqrt{-B}}\right)}
,
\ee
\be
\leadsto
v_n=\prod_{k=1}^n\left(
A-2i\sqrt{-B}\cos\frac{\pi (k-1/2)}{n}
\right)
=
2(i\sqrt{-B})^n\cos\left(
n\cos^{-1}\frac{-iA}{2\sqrt{-B}}
\right)
,
\ee
where $u_0=0$, $u_1=1$, $v_0=2$, $v_1=A$

\cite{RamirezJIS16}
Let 
\be
t_0=1,\quad t_1=1,\quad t_n=\left\{\begin{array}{ll}
at_{n-1}+t_{n-2},& n \mathrm{even},\\
bt_{n-1}+t_{n-2},& n \mathrm{odd}.
\end{array} \right.
\ee
Then
\be
t_n=a^{\xi(n-1)}\sum_{i=0}^{\lfloor(n-1)/2\rfloor}\binom{n-i-1}{i}(ab)^{\lfloor(n-1)/2\rfloor-i}
\ee
where $\xi(n)=0$ when $n$ even, $\xi(n)=1$ when $n$ odd.

\cite{RamirezJIS16}
Let
\be
F_n^{(a,b)}(q,s) = \chi(n)F_{n-1}^{(a,b)}(q,s)+q^{n-2}sF_{n-2}^{(a,b)}(q,s)
\ee
starting at $F_0^{(a,b)}(q,s)=0, F_1^{(a,b)}(q,s)=1$, then
\be
F_n^{(a,b)}(q,s)=\chi(n)F_{n-1}^{(a,b)}(q,qs)+qsF_{n-2}^{(a,b)}(q,q^2s).
\ee
and
\be
F_n^{(a,b)}(q,s)=a^{\xi(n-1)}\sum_{l=0}^{\lfloor(n-1)/2\rfloor}\left[
\begin{array}{c}n-l-1\\ l\end{array}
\right](ab)^{\lfloor(n-1)/2\rfloor -l}q^{l^2}s^l.
\ee
where the $q$-binomial is 
\be
\left[\begin{array}{c}n\\k\end{array}\right]\equiv \frac{(q;q)_n}{(q;q)_k(q;q)_{n-k}}.
\label{eq.qbin}
\ee

\cite{WilliamsJNT3}
\be
\sum_{x=1}^{p-1}F(x)+\sum_{x=1}^{p-1}\left(\frac{x}{p}\right)
F(x)=\sum_{x=1}^{p-1}F(x^2),
\ee
if $p$ is an odd prime and $F(x)=F(x+p)$, where $(\frac{x}{p})$
is the Legendre Symbol.

\cite{WilliamsJNT3}
\be
\sum_{x=1}^{p-1}\left(\frac{x}{p}\right)
e(k(x+\bar x)) = \left(\frac{k}{p}\right)
i^{(p-1)^2/4}p^{1/2}(e(2k)+e(-2k)).
\ee
where $\bar x$ is the unique solution to $x\bar x\equiv 1$ $\pmod p$,
and $e(t)=e^{2\pi i t/p}$.

\cite{ChenJNT124}
\begin{multline}
\sum_{k=0}^m \binom{m}{k}\frac{(n+k)!}{(n+k+s)!}a_{n+k+s}
=
\sum_{k=0}^n \binom{n}{k}\frac{(-1)^{n-k}(m+k)!}{(m+k+s)!}b_{m+k+s}
\\
+
\sum_{j=0}^{s-1}\sum_{i=0}^{s-1-j}
\binom{s-1-j}{i}
\binom{s-1}{j}
\frac{(-1)^{n+1+i}a_j}{(s-1)!(m+n+1+i)\binom{m+n+i}{n}}
,
\end{multline}
where $b_n\equiv \sum_{k=0}^n \binom{n}{k}a_k$. 

\cite{ChenJNT124}
\be
\sum_{k=0}^m \binom{m}{k}\binom{n+k}{s}a_{n+k+s}
=
\sum_{k=0}^n \binom{n}{k}\binom{m+k}{s}(-1)^{n-k}b_{m+k+s}
,
\ee
where $b_n\equiv \sum_{k=0}^n \binom{n}{k}a_k$. 

\cite{ChenJNT124}
\begin{multline}
\sum_{k=0}^m
\frac{\binom{m}{k}}{\binom{n+k+s}{s}}
x^{m-k}A_{n+k+s}(y)
=
\sum_{k=0}^n
\frac{\binom{n}{k}}{\binom{m+k+s}{s}}
(-1)^{n+m+s}
x^{n-k}A^*_{m+k+s}(z)
\\
+\sum_{j=0}^{s-1}\sum_{i=0}^{s-1-j}
\binom{s-1-j}{i}\binom{s-1}{j}
\frac{(-1)^{n+1+i}x^{m+n+s-j}sA_j(y)}{(m+n+1+i)\binom{m+n+i}{n}}
,
\end{multline}
and
\be
\sum_{k=0}^m
\binom{m}{k}\binom{n+k}{s}
x^{m-k}A_{n+k+s}(y)
=
\sum_{k=0}^n
\binom{n}{k}\binom{m+k}{s}
(-1)^{n+m+s}
x^{n-k}A^*_{m+k+s}(z)
\ee
where $A_n(x)\equiv \sum_{k=0}^n\binom{n}{k}(-1)^k a_k x^{n-k}$
and
$A^*_n(x)\equiv \sum_{k=0}^n\binom{n}{k}(-1)^k a_k^* x^{n-k}$
and
$a^*_n\equiv \sum_{k=0}^n\binom{n}{k}(-1)^k a_k$ and $x+y+z=1$.

\cite{ChenJNT124}
\be
a_{n,m}=\sum_{k=0}^n \binom{n}{k}\alpha^{n-k}\beta^k a_{0,m+k},
\ee
and
\be
\sum_{k=0}^n\binom{n}{k}\alpha^{n-k}\beta^k a_{0,m+k} =
\sum_{k=0}^m\binom{m}{k}(-\alpha)^{m-k}\beta^{-m}a_{n+k,0},
\ee
where $a_{n,m}\equiv \alpha a_{n-1,m}+\beta a_{n-1,m+1}$ for $n\ge 1,m\ge 0$.

\cite{TakacsACMA}
Let $A_m(\alpha,\beta)=\binom{\alpha+\beta m}{m}\frac{\alpha}{\alpha+\beta m}$, then
\begin{equation}
\sum_{i=0}^mA_i(\alpha,\beta) A_{m-i}{\gamma,\beta} = A_m(\alpha+\gamma,\beta).
\end{equation}
and
\begin{equation}
\sum_{i=0}^m iA_i(\alpha,\beta) A_{m-i}{\gamma,\beta} = \frac{m\alpha}{\alpha+\gamma}A_m(\alpha+\gamma,\beta).
\end{equation}
\cite{TakacsACMA}
\begin{equation}
\binom{nm}{n}=\sum_{k=1}^n \frac{b-1}{bk-1}\binom{nk}{k}\binom{bn-bk}{n-k},
\end{equation}
for $b=2,3,\ldots$ and $n=1,2,\ldots$

\cite{GrahamAMM73}
\be
x_{nm}=\sum_{k=0}^m x_{mk}\binom{n+k}{2m}
\leadsto
x_{nm}=\sum_{k=0}^m \frac{2k+1}{m+k+1}\binom{n+k}{m+k}\binom{n-1-k}{m-k}x_{kk},
\ee
for $m<n$, $n\ge 0$, $0\le m\le n$.

\cite{GouldBAMS66}
If 
\be
F_n=\sum_{k=0}^n (-)^k\binom{n}{k}\binom{a+bk}{n}f(k),
\label{eq.Fgould}
\ee
and $f(k)$ is independent of $n$ and $f(0)=1$, then
\be
\sum_{k=0}^\infty \binom{a+bk}{k}z^kf(k)=x^a\sum_{n=0}^\infty (-)^n F(n)u^n,
\ee
where $z=(x-1)/x^b$ and $u=(x-1)/x$. If $F(n)$ is defined as in \eqref{eq.Fgould}, then
\be
\binom{a+bn}{n} f(n) = \sum_{k=0}^n (-)^k\frac{a+bk-k}{a+bn-k}\binom{a+bn-k}{n-k} F(k).
\ee
and
\be
\sum_{k=0}^\infty \frac{w}{w+bk}\binom{a+bk}{k}z^k=x^a\sum_{n=0}^\infty
(-)^n \binom{a-w}{n}\binom{w/b+n}{n}^{-1} u^n.
\ee

\cite{GouldBAMS66}
If
\be
F(n)=\sum_{k=0}^n (-)^k\binom{n}{k}\frac{(a+bk)^n}{n!}f(k),
\label{eq.F2gould}
\ee
where $f(k)$ is independent of $n$ and $f(0)=1$, then
\be
\sum_{k=0}^\infty \frac{(a+bk)^k}{k!}z^kf(k)=x^a\sum_{n=0}^\infty (-)^n F(n)u^n
\ee
where $z=u/x^b$ and $u=\log x$. If $F$ is defined by \eqref{eq.F2gould}, then
\be
\frac{(a+bn)^n}{n!}f(n) = \sum_{k=0}^n (-)^k f(k) \frac{(a+bn)^{n-k}}{(n-k)!}
\, \cdot \frac{a+bk}{a+bn}.
\ee
\be
\sum_{k=0}^\infty \frac{w}{w+bk}\cdot \frac{(a+bk)^k}{k!} z^k
=x^a\sum_{n=0}^\infty (-)^n \frac{(a-w)^n}{n!}\binom{w/b+n}{n}^{-1}u^n.
\ee

\cite{BorosSSA11}
\be
F(x)\equiv\sum_{k\ge 1} f(x/k);\quad G(x)\equiv \sum_{k\ge 1}(-)^kf(x/k)
\leadsto
 F(x)=2^{-n}F(2^nx)+\sum_{k\ge 1}2^{-k}G(2^kx).
\ee

\cite{BorosSSA11}
\be
F(x) = r_1F(m_1x)+r_2G(m_2x)
\leadsto
r_2\sum_{k=1}^n r_1^{k-1}G(m_1^{k-1}m_2x) =
F(x)-r_1^nF(m_1^nx)
.
\ee

\cite{ChenJIS19}
\be
\sum_{k=0}^n (\binom{n}{k})^2H_k=\binom{2n}{n} \left(2H_n-H_{2n}\right)
.
\ee

\cite{SpiveyDM307}
\be
\sum_{k}\binom{n}{2k}\frac{1}{k+1}=\frac{n2^{n+1}+2}{(n+1)(n+2)}
.
\ee

\cite{ChoiCKM14}
\be
\sum_{k=0}^{n-1}(-)^k\binom{n-1}{k}\frac{1}{2^k(n+k+1)}
=\frac{2^n(n-1)!n!}{(2n)!}-\frac{1}{n2^n}.
\ee

\cite{SprugnoliEJCN6}
\be
\sum_{k=0}^n \frac{(-1)^k}{1-2k}\binom{n}{k} = \frac{4^n}{\binom{2n}{n}}.
\ee

\cite{SprugnoliEJCN6}
\be
\sum_{k=0}^n \frac{1}{(1+2k)(1+2n-2k}\binom{2k}{k}\binom{2n-2k}{n-k} = 
\frac{16^n}{(n+1)(2n+1)\binom{2n}{n}}.
\ee

\cite{SpiveyDM307}
\be
\sum_{k}\binom{n}{2k+1}\frac{1}{k+1}=\frac{2^{n+1}-2}{n+1}
.
\ee

\cite{RottbrandITSF10}
\be
\sum_{j=0}^n \binom{m-a+b}{j}
\binom{n+a-b}{n-j}
\binom{a+j}{m+n}
=
\binom{a}{m}
\binom{b}{n}
,
\ee
where $n,m$ are integer and $a,b$ real.

\cite{RottbrandITSF10}
\be
\sum_{m=0}^P (-)^m\binom{P}{m}
\binom{a-m}{M}
=
\binom{a-P}{M-P}
,
\ee
where $P,M$ are integer and $a,b$ real.

\cite{PilehroodEJC18}
\begin{multline}
\sum_{k=0}^{n-1}(-)^{n-1-k}
\binom{2k}{k}^5(205k^2+160k+32)
=
16n\binom{2n}{n}\sum_{k=0}^{n-1}\binom{n+k-1}{k}^4(2k+n)
\\
=
8n^2\binom{2n}{n}^2\,\sum_{k=0}^{n-1}\binom{2n-1}{n+k}
\binom{2n-k-2}{n-k-1}^2.
\end{multline}

\cite{BorosMC71}
\be
\sum_{j=0}^{N-k}\binom{2N+1}{2j}\binom{N-j}{k}=\binom{2N-k}{k}4^{N-k}.
\ee
\cite{BorosMC71}
\be
\sum_{j=0}^k\binom{2N}{2j}\binom{N-j}{N-k}
=\frac{2^{2k-1}N}{k}\binom{k+N-1}{N-k}.
\ee

\cite{BorosMC71}
\begin{multline}
\sum_{k=0}^p d_{p+1-k}z^{2k}
\sum_{j=0}^{p-k}\binom{2p-2k}{2j}(1+z^2)^{p-k-j}
=
\left(\sum_{j=1}^{p+1}d_j\right)z^{2p}
\\
+\sum_{i=1}^p 2^{2i-1}z^{2(p-i)}
\left(
\sum_{j=1}^{p+1-i}\frac{j+i-1}{i}\binom{j+2i-2}{j-1}d_{j+i}
\right)
\end{multline}
with $d_{p+1}\equiv 1$.

\cite{BorosMC71}
\be
\sum_{j=0}^N\binom{2N+1}{2j}(1+z^2)^{N-j}=
\sum_{j=0}^N\binom{N+j}{2j}4^j z^{2(N-j)}.
\ee

\cite{ChuRDCM}
\begin{equation}
\sum_{k=0}^n\frac{u+vk}{(a+bk)(c+b(n-k))}
\binom{a+bk}{k}\binom{c+b(n-k)}{n-k}
=\frac{(a+c)u+avn}{ac(a+c+bn)}\binom{a+c+bn}{n}.
\end{equation}

\cite{ChuRDCM}
\be
\sum_{k=0}^n\binom{a+bk}{k}\binom{c-bk}{n-k}
=\sum_{k=0}^n\binom{a+\epsilon+bk}{k}\binom{c-\epsilon-bk}{n-k}.
\ee

\cite{ChuRDCM}
\be
\sum_{k=0}^n\binom{a+bk}{k}\binom{c-bk}{n-k}
=\sum_{k=0}^n\binom{a+c-k}{n-k}b^k.
\ee

\cite{ChuRDCM}
\be
\sum_{k=0}^n\frac{a(a-1)}{(a+bk)(a+bk-1)}\binom{a+bk}{k}\binom{c+b(n-k)}{n-k}
=\frac{(a+b-1)(a+c)+b(b-1)n}{(a+b-1)(a+c+bn)}\binom{a+c+bn}{n}.
\ee

\cite{ChuRDCM}
\begin{multline}
\sum_{k=0}^n\frac{a(a-1)}{(a+bk)(a+bk-1)}
\binom{a+bk}{k}\frac{1+c+b(n-k)}{1+c}\binom{c+b(n-k)}{n-k}
\\
=\frac{(a+b-1)(1+c)-b(b-1)n}{(a+b-1)(1+c+)}\binom{a+c+bn}{n}.
\end{multline}

\cite{GouldJIS19}
\be
\sum_{k=0}^n \frac{\binom{n}{k}}{\binom{x}{k+b}} = \frac{x+1}{(x-n+1)\binom{x-n}{b}}.
\ee

\cite{GouldJIS19}
\be
\sum_{k=0}^n (-)^k\binom{n}{k}\frac{1}{\binom{b+k}{c}}=\frac{c}{n+c}\frac{1}{\binom{n+b}{b-c}}.
\ee

\cite{FrankelAMM51}
\be
\sum_{r=0}^m (-)^r\binom{m+k+a}{r+k+a}\binom{r+k}{k} = \binom{m+a-1}{m}.
\ee

\cite{WangJIS14}
\be
\sum_{k=0}^{2n}(-)^k\binom{2n}{k}^3 = (-)^n\binom{2n}{n}\binom{3n}{n}.
\ee

\cite{RoyAMM94}
\be
\sum_k \binom{r}{k}\binom{s}{n+k}
=\binom{r+s}{r+n},\quad r=0,1,2,\ldots
\quad,
n=\ldots -2,-1,0,1,2,3,\ldots
\ee

\cite{RottbrandITSF10,RoyAMM94}
\be
\sum_{m=0}^n (-)^m\binom{n}{m}
\binom{a+m}{p}
=
(-)^n \binom{a}{p-n}
,
\ee
where $n,p$ are integer and $a$ is real.

\cite{RoyAMM94}
\be
\sum_{k\ge 0}\binom{n+k}{m+2k}\binom{2k}{k}\frac{(-1)^k}{k+1}
= \binom{n-1}{m-1}.
\ee

\cite{RoyAMM94}
\be
\sum_{\nu=0}^s (-)^\nu \binom{\beta}{\nu}\binom{\beta+s-\nu}{\beta}
\frac{\alpha}{\alpha+s-\nu}
= \frac{(\alpha-\beta)_s}{(\alpha+1)_s}.
\ee

\cite{RoyAMM94}
\be
\sum_{k=-l}^l (-)^k\binom{2l}{l+k}\binom{2m}{m+k}\binom{2n}{n+k}
=
\frac{(l+m+n)!(2l)!(2m)!(2n)!}{(l+m)!(m+n)!(n+l)!l!m!n!},\quad l=\min(l,m,n).
\ee

\cite{FoataLNM138}
\be
\sum_{k\ge 0}t^k\binom{l+k}{l} = \frac{1}{(1-t)^{l+1}}.
\ee

\cite{Wegschaider97}
\be
\sum_k \sum_j\binom{n}{k}\binom{n+k}{k}\binom{k}{j}^3
=
\sum_k \binom{n}{k}^2\binom{n+k}{k}^2.
\ee

\cite{Wegschaider97}
\be
\sum_{k_1} \sum_{k_2\le k_1} \sum_{k_3\le k_2}
(k_1-k_2)(k_1-k_3)(k_2-k_3)\binom{n}{k_1}
\binom{n}{k_2}
\binom{n}{k_3}
=
n^2(n-1)8^{n-2}\frac{(3/2)_{n-2}}{(3)_{n-2}}.
\ee

\cite{Wegschaider97}
\be
\sum_i \sum_j \binom{i+j}{i} \binom{n-i}{j}
\binom{n-j}{n- i- j}
=\sum_{k=0}^n\binom{2k}{k}.
\ee

\cite{Wegschaider97}
\be
\sum_{i=0}^n\sum_{j=0}^m\binom{i+j}{j}^2
\binom{m+m-i-j}{n-j}^2
=\frac{1}{2}\binom{2m+2n+2}{2n+1}.
\ee

\cite{Wegschaider97}
\be
\sum_{s=0}^k \sum_{b\ge 0} (-)^b\binom{s}{b}
\binom{k-s}{2v-b}\binom{k-2v}{s-b}
=\binom{k-v}{k-2v}2^{k-2v},\quad k\ge 2v.
\ee

\cite{Wegschaider97}
\be
\sum_{j,k} (-)^{j+k}\binom{j+k}{k+l}
\binom{r}{j}\binom{n}{k}\binom{s+n-j-k}{m-j}
=(-)^l\binom{n+r}{n+l}\binom{s-r}{m-n-l}.
\ee

\cite{FoataLNM138}
\be
r!\sum_{0\le i\le k} (-)^i(k-i+r)^{n-1}\binom{n+r}{i}\binom{k-i+r}{r}
={}^rA_{n-1+r,k}
\ee
where $^rA_n(s+1)=\sum_{0\le k\le n-r}s^k(n-k)!S(n+1-r,n+1-r-k)$ are Eulerian Numbers.

\cite{ChenJACM196}
\be
\sum_{i=0}^{n} \sum_{j=0}^n\binom{i+j}{i}
\binom{m-i+j}{j}\binom{n-j+i}{i}\binom{m+n-i-j}{m-i}
=
\frac{(m+n+1)!}{m!n!}\sum_k \frac{1}{2k+1}
\binom{m}{k}\binom{n}{k}.
\ee

\cite{ChenJACM196}
\be
\sum_{i=0}^n \sum_{j=0}^n\binom{i+j}{i}^2
\binom{4n-2i-2j}{2n-2i}=(2n+1)\binom{2n}{n}^2
.
\ee

\cite{ChenJACM196}
\be
\sum_{i=0}^{\lfloor m/2\rfloor} \sum_{j=0}^{\lfloor n/2\rfloor}\binom{i+j}{i}^2
\binom{m+n-2i-2j}{n-2i}
=
\frac{\lfloor(m+n+1)/2\rfloor !\lfloor(m+n+2)/2\rfloor !}
{\lfloor m/2\rfloor !\lfloor(m+1)/2\rfloor ! \lfloor(n+1)/2\rfloor !}
.
\ee

\cite{ChenJACM196}
\be
\sum_i \sum_j\binom{n}{j}
\binom{n+j}{j}
\binom{j}{i}^2
\binom{2i}{i}^2
\binom{2i}{j-i}
=
\sum_k \binom{n}{k}^3\binom{n+k}{k}^3.
\ee

\cite{ChenJACM196}
\be
\sum_j \sum_k (-1)^{j+k}\binom{j+k}{k+l}
\binom{r}{j}
\binom{n}{k}
\binom{s+n-j-k}{m-j}
=(-1)^l
\binom{n+r}{n+l}\binom{s-r}{m-n-l}.
\ee

\cite{ChenJACM196,Wegschaider97}
\be
\sum_r \sum_s (-1)^{n+r+s}\binom{n}{r}
\binom{n}{s}
\binom{n+s}{s}
\binom{n+r}{r}
\binom{2n-r-s}{n}
=\sum_k
\binom{n}{k}^4.
\ee

\cite{ShanksAMM58}
\be
\binom{n}{i}^k=\binom{n+ik-i}{ik}+A_2\binom{n+ik-i-1}{ik}
+\cdots A_{ik-i}\binom{n+1}{ik}+\binom{n}{ik},
\ee
\begin{eqnarray*}
S_i^k(n,0)&\equiv& \binom{n}{i}^k;
\nonumber
\\
S_i^k(n,p)&\equiv& S_i^k(1,p-1)+S_i^k(2,p-1)+\cdots S_i^k(n,k-1)
\\
&=&\binom{n+ik-i+p}{ik+p}
+A_2\binom{n+ik-i+p-1}{ik+p}
+\cdots
+A_{ik-i}\binom{n+p+1}{ik+p}
+\binom{n+p}{ik+p}
\end{eqnarray*}
where $A_j=A_{ik-i-j+2}$ for $j=2,3,\ldots,ik-i$,
\be
A_j=\left|
\begin{array}{cccccc}
1 & 0 & 0 & \cdots & 0 & \binom{i}{i}^k \\
\binom{ik+1}{ik} & 1 & 0 & \cdots & 0 & \binom{i+1}{i}^k \\
\binom{ik+2}{ik} & \binom{ik+1}{ik} & 1 & \cdots & 0 & \binom{i+2}{i}^k \\
\vdots  & \vdots & \vdots & \vdots & \vdots & \vdots \\
\binom{ik+j-1}{ik} & \binom{ik+j-2}{ik} & \binom{ik+j-3}{ik} & \cdots & \binom{ik+1}{ik} & \binom{i+j-1}{i}^k \\
\end{array}
\right|
\ee

\cite{TrifFQ38}
\be
\sum_{k=0}^{2n} (-)^k\frac{\binom{2n}{k}}{\binom{4n}{2k}}=\frac{4n+1}{2n+1}.
\ee

\cite{TrifFQ38}
\be
\sum_{k=0}^{2n} (-)^k\frac{\binom{4n}{2k}}{\binom{2n}{k}}= -\frac{1}{2n-1}.
\ee

\cite{SprugnoliEJCN6}
\be
\sum_{k=0}^n (-4)^k\frac{\binom{n}{k}}{\binom{2k}{k}} =\frac{1}{1-2n}.
\ee

\cite{SprugnoliEJCN6}
\be
\sum_{k=1}^n (-4)^k\frac{\binom{n-1}{k-1}}{k^2\binom{2k}{k}} =\frac{H_n-2H_{2n}}{n}.
\ee

\cite{TrifFQ38}
\be
\sum_{k=0}^{m} \frac{\binom{m}{k}}{\binom{n+m}{p+k}}= \frac{n+m+1}{n+1}\binom{n}{p}^{-1},
\ee
with $m$, $n$, $p$ nonnegative integers and $p\le n$.

\cite{TrifFQ38}
\be
\sum_{k=0}^{n} (-)^k \frac{1}{\binom{n+m}{m+k}}= \frac{n+m+1}{m+n+2}\left(\binom{m+n+1}{m}^{-1}
+(-)^n
\right),
\ee
with $m$ and $n$ nonnegative integers.

\cite{SprugnoliEJCN6}
\be
\sum_{k=1}^n \frac{4^k}{k\binom{2k}{k}} = 2\frac{4^n}{\binom{2n}{n}}-2.
\ee

\cite{SprugnoliEJCN6}
\be
\sum_{k=2}^{n-1} \frac{4^k}{k(k-1)\binom{2k}{k}} = 4-\frac{(2n-1)4^n}{n(n-1)\binom{2n}{n}}.
\ee

\cite{SprugnoliEJCN6}
\be
\sum_{k=1}^{n-1} \frac{4^k}{k(2k+1)\binom{2k}{k}} = 2-\frac{4^n}{n\binom{2n}{n}}.
\ee

\cite{SprugnoliEJCN6}
\be
\sum_{k=1}^n (-4)^k\frac{\binom{n}{k}}{k\binom{2k}{k}} = H_n-2H_{2n}.
\ee

\cite{BelbachirJIS14}
Let $S_n^{(m)}\equiv \sum_{k=0}^n k^m/\binom{n}{k}$, then
\be
S_n^{(0)}=\frac{n+1}{2^{n+1}}\sum_{k=1}^{n+1}\frac{2^k}{k};
\ee
\be
\sum_{n\ge 0}S_n^{(0)}x^n = \frac{2}{(x-1)(x-2)}-2\frac{\ln(1-x)}{(x-2)^2};
\ee
\be
\sum_{n\ge 0}S_n^{(1)}x^n = -\frac{x(3x-4)}{(x-1)^2(x-2)^2}+2x\frac{\ln(1-x)}{(x-2)^3};
\ee
\be
\sum_{n\ge 0}S_n^{(2)}x^n = \frac{2x(2x^3-3x^2-3x+5)}{(x-2)^3(x-1)^3}-2(x^2+2x-2)\frac{\ln(1-x)}{(x-2)^4};
\ee
and similar generating functions in terms of Worpitzky numbers\ldots
\be
S_{n+1}^{(2)}=\frac{(n-1)(n+2)^2}{2(n-2)(n+1)^2}S_n^{(2)}+\frac{(n+2)(n^2-2n-2)}{2(n-2)};
\ee
\be
S_n^{(2)}=\frac14 (n+1)(n-2)S_n^{(0)}+\frac12 (n+1)^2;
\ee
for non-negative $n$:
\be
S_n^{(3)}=\frac18 n(n^2-3n-6)S_n^{(0)}+\frac34 n(n+1)^2 ;
\ee
\be
S_{n+j}^{(0)}=(n+j+1)\left( \frac{1}{2^j(n+1)}S_n^{(0)}+\sum_{r=0}^{j-1}\frac{1}{2^r(n+j-r+1)}\right);
\ee
\be
S_n^{(4)}=\frac{1}{16} (n+1)(n^3-7n^2-2n+16)S_n^{(0)}+\frac18 (7n-8)(n+1)^3 ;
\ee
\be
S_{n+1}^{(m)}=\delta_{0m}+\frac{1}{n+1}\sum_{j=0}^{m+1}\binom{m+1}{j}S_n^{(j)}.
\ee

\cite{SpiveyDM307}
\be
\sum_{k=0}^n\left[\begin{array}{c}n\\k\end{array}\right]H_k=
n!\sum_{k=1}^n\frac{c_{k-1}}{k!},
\ee
where $H_k$ are the harmonic numbers
and $[\ldots]$ unsigned Stirling numbers
of the first kind.

\cite{MaierTAAMS358}
\be
\sum_{j=0}^n \frac{(\alpha_1)_j\cdots (\alpha_{r+1})_j}{(\beta_1+1)_j\cdots (\beta_r+1)_j}=
\frac{(\alpha_1+1)_n\cdots (\alpha_{r+1}+1)_n}{(\beta_1+1)_n\cdots (\beta_r+1)_nn!}
\ee
for every integer $n\ge 0$ if the $r$ equations
\begin{eqnarray}
\sum_{1\le i\le r+1}\alpha_i &=& \sum_{1\le i\le r}\beta_i,\label{eq.maiertaams}\\
\sum_{1\le i< j\le r+1}\alpha_i\alpha_j &=& \sum_{1\le i<j \le r}\beta_i\beta_j,\\
\sum_{1\le i< j< k\le r+1}\alpha_i\alpha_j\alpha_k &=& \sum_{1\le i<j<k \le r}\beta_i\beta_j\beta_i,
\end{eqnarray}
etc. are satisfied.

\cite{AndersenACMA}
If $p=b+c+d+e-n+1$ is a non-negative integer then
\begin{multline}
\sum_{k=0}^n \binom{n}{k}[n+2a]_k[b+a]_k[c+a]_k[d+a]_k
[n-2a]_{n-k}[b-a]_{n-k}[c-a]_{n-k}[d-a]_{n-k}[e-a]_{n-k}(n+2a-2k)
\\
=
(-)^n[n+2a]_{2n+1}[b+d-p]_{n-p}[b+e-p]_{n-p}[d+e-p]_{n-p}
\\ \times
\sum_{j=0}^p \binom{p}{j}
[n]_j[c+a]_j[c-a]_j(b+e-n]_{p-j}[d+e-n]_{p-j}[d+e-n]_{p-j}
\end{multline}
where $[x]_n=x(x-1)\cdots (x-n+1)$ is the falling factorial.

\cite{AndersenACMA}
If $p=a_1+a_2+b_1+b_2-n+m+1$ is a non-negative integer then
\begin{multline}
\sum_{k=m}^n\binom{n-m}{k-m}[a_1]_{k-m}[a_2]_{k-m}[b_1]_{n-k}[b_2]_{n-k}(-1)^{k-m}
\\
= 
\sum_{k=m}^n\binom{n-m}{k-m}[a_1]_{k-m}[a_2]_{k-m}[b_1]_{n-k}[n-m+p-a_1-a_2-b_1-1]_{n-k}(-1)^{k-m}
\\
= 
\sum_{k=0}^p\binom{p}{k}[n-m]_k[a_2]_k[a_1+b_1-p]_{n-m-k}[n-m-a_2-b_1-1]_{n-m-k}(-1)^k
\\
= 
[a_1+b_1-p]_{n-m-p}[n-m-a_2-b_1-1]_{n-m-p}
\sum_{k=0}^p\binom{p}{k}[n-m]_k[a_2]_k[a_1+b_1-n+m]_{p-k}[p-a_2-b_1-1]_{p-k}(-1)^k
\end{multline}

\cite{AndersenACMA}
For arbitrary $a, b,c$ and integers $m$ and $n$ with $m\le n$
\begin{multline}
\sum_{k=m}^n\binom{n-m}{k-m} [n-m+2a]_{k-m} [b+a]_{k-m}[c+a]_{k-m}
\\ \times
[b-a]_{n-k}[c-a]_{n-k}[n-m-2a]_{n-k} (n+m+2a-2k)
=
[n-m+2a]_{2n-2m+1}[b+c]_{n-m}.
\end{multline}

\cite{BlanchardRWalk}
\be
1=\sum_{1\alpha_1+2\alpha_2+\cdots N\alpha_N=N} \frac{1}{1^{\alpha_1}2^{\alpha_2}\cdots N^{\alpha_N}\alpha_1!\alpha_2!\cdots}
\ee

\cite{SpiveyDM307}
\be
\sum_{k=0}^ns(n,k)=s(n-1,0)+s(n-1,-1),
\ee
where $s(.,.)$ are the unsigned Stirling numbers of the first kind.

\cite{SpiveyDM307}
\be
\sum_{k=0}^ns(n,k)k=s(n-1,1)+s(n-1,0).
\ee

\cite{YuanJNT128}
Let $S_r(n)\equiv \sum_{k=0}^n \binom{n}{k}^r$. This defines P-finite polynomials
\be
P_0(n)S_r(n+1)+P_1(n)S_r(n)+P_2(n)S_r(n-1)=0, n\ge 0
\ee
where for $r=3$
\be
P_0(n)=(n+1)^2, \, P_1(n)=-(7n^2+7n+2),\, P_2(n)=-8n^2,
\ee
for $r=4$
\be
P_0(n)=(n+1)^3, \, P_1(n)=-2(2n+1)(3n^2+3n+1),\, P_2(n)=-4n(4n+1)(4n-1).
\ee
4-term recurrences for $r=5,6$ are also in the paper.

\cite[(6.16)]{GrahamKnuthPatashnik}
\be
\sum_{k=0}^ns(k,m)\binom{k}{m}=s(n+1,m+1).
\ee
where $s(.,.)$ are the unsigned Stirling Numbers of the First Kind.
\cite{SpiveyDM307,SprugnoliDM142}
\be
n!\sum_{j=0}^n s(j,m)/j!=s(n+1,m+1).
\ee

\cite{SpiveyDM307}
\be
\sum_{k=0}^ns(n,k)k^{\underline m}=m![s(n-1,m)+s(n-1,m-1)],
\ee
where $k^{\underline m}=k(k-1)(k-2)\cdots (k-m+1)$.

\cite{SpiveyDM307}
\be
\sum_{k=0}^ns(n,k)k^m=
\sum_{j=0}^m\left\{\begin{array}{c}m\\j\end{array}\right\}
\left(s(n-1,j)+s(n-1,j-1)\right)j!
.
\ee

\cite{SpiveyDM307}
\be
\sum_{k=0}^n\frac{s(n,k)}{k+1}
=
b_n
,
\ee
where $b_n=\int_0^1 x^{\underline n}dx$ are the Cauchy numbers
of the first type \cite[A006232]{sloane}.

\cite{SpiveyDM307}
\be
n!\sum_{k=0}^n\left[\begin{array}{c}k\\m\end{array}\right]\frac{1}{k!}
=
\left[\begin{array}{c}n+1\\m+1\end{array}\right]
.
\ee
where $[\cdots]$ are the unsigned Stirling numbers of the first kind.

\cite{SpiveyDM307}
\be
\sum_{k=0}^n\left[\begin{array}{c}n\\k\end{array}\right]k
=
\left[\begin{array}{c}n+1\\2\end{array}\right]
.
\ee

\cite{SpiveyDM307}
\be
\sum_{k=0}^n\left[\begin{array}{c}n\\k\end{array}\right]k^{\underline m}
=
\left[\begin{array}{c}n+1\\m+1\end{array}\right]m!
,
\ee
where $k^{\underline m}=k(k-1)(k-2)\cdots (k-m+1)$.

\cite{SpiveyDM307}
\be
\sum_{k=0}^n\left[\begin{array}{c}n\\k\end{array}\right]k^ m
=
\sum_{j=0}^m \left[\begin{array}{c}n+1\\j+1\end{array}\right]
\left\{\begin{array}{c}m\\j\end{array}\right\}j!
.
\ee
\cite{ZhaoDM281}
\be
\sum_{k=0}^n \begin{Bmatrix}k+1\\m+1\end{Bmatrix}
=
\sum_{j=0}^n \begin{Bmatrix}j\\m\end{Bmatrix}
\binom{n+1}{j+1}
.
\ee
where $\{.\}$ denotes the Stirling numbers of the 2nd kind.

\cite{BalakrishnanACMA}
For a vector ${\bf a}=(a_0,a_1,\ldots)$ define $p_0(x,{\bf a})=1$, $p_k(x,{\bf a})=(x-a_0)(x-a_1)\cdots (x-a_{k-1})$,
for $k=1,2,\ldots$. Define the generalized Stirling numbers $s(n,k,{\bf a})$ of the first and $S(n,k,{\bf a})$ of
the second kind as
\begin{equation}
p_n(x,{\bf a})=\sum_{k=0}^n s(n,k,{\bf a})x^k, n=0,1,\ldots
\end{equation}
and
\begin{equation}
x^n=\sum_{k=0}^n S(n,k,{\bf a})p_k(x,{\bf a}).
\end{equation}
Then
\begin{equation}
s(n+1,k,{\bf a})=s(n,k-1,{\bf a})-a_ns(n,k,{\bf a});
S(n+1,k,{\bf a})=S(n,k-1,{\bf a})+a_ks(n,k,{\bf a}).
\end{equation}

\cite{ChuRDCM}
If
\be
\phi(x;n) \equiv \prod_{k=0}^{n-1}(a_k+xb_k); \phi(x;0)=1
\ee
then
\be
f(n)=\sum_{k=0}^n(-)^k\binom{n}{k}\phi(k;n)g(k)
\leftrightarrow
g(n)=\sum_{k=0}^n(-)^k\binom{n}{k}\frac{a_k+kb_k}{\phi(n;k+1)}f(k)
\ee

\cite{BenyiECA2}
Let
\be
\Lif_k(t)\equiv \sum_{n\ge 0}\frac{t^n}{n!(n+1)^k}
\ee
be the k-th polylogarithm factorial function and
\be
\hat c_n^{(k)}=(-1)^n\sum_{m=0}^n
\left[\begin{array}{c}n\\m\end{array}\right]\frac{1}{(m+1)^k}
\ee
be the poly-Cauchy numbers of the second kind \cite[A344639]{sloane}. Then
\be
\Lif_k(-\ln(1+t)) = \sum_{n\ge 0}\hat c_n^{(k)}\frac{t^n}{n!}.
\ee

\cite{SpiveyDM307}
\be
\sum_{k=0}^n\left[\begin{array}{c}n\\k\end{array}\right]\frac{1}{k+1}=c_n,
\ee
where $c_n=\int_0^1 (x)_ndx$ are the Cauchy numbers of the second
type \cite[A002657]{sloane},
using Pochhammer's symbol.

\cite{NetoJIS20}

\be
\sum_{k=0}^n \sum_{j=0}^m (m)_{n-k} \left[\begin{array}{c}m\\ j\end{array}\right] \binom{n}{k}k!
\ee

\cite{Boyadzhievarxiv2012}
\be
\sum_{k=0}^n\left[\begin{array}{c}n\\k\end{array}\right] k^2=n!(H_n+H_n^2-H_n^{(2)}).
\ee

\cite{ChenJNT124}
\begin{multline}
\sum_{k=0}^m
\binom{m}{k}(n+k)!\left\{\begin{array}{c}n+k+s\\ q \end{array}\right\}
=
\sum_{k=0}^n
\binom{n}{k}(m+k)!(-1)^{n-k}
\frac{(m+k+s)^q}{(m+k+s)!}
\\
+\sum_{j=0}^{s-1}\sum_{i=0}^{s-1-j}
\binom{s-1-j}{i}
\frac{(-1)^{n+1+i}
\left\{
\begin{array}{c}j\\q\end{array}
\right\}
}{(s-1-j)!(m+n+1+i)\binom{m+n+i}{n}}
\end{multline}
and
\be
\sum_{k=0}^m
(n+k)!
\left\{\begin{array}{c}n+k-s\\q\end{array}\right\}
=
\sum_{k=0}^n
\binom{n}{k}(m+k)!
(-1)^{n-k}
\frac{(m+k-s)^q}{(m+k-s)!}
.
\ee

\cite{MezoCEJM}
\be
\begin{Bmatrix} n \\ m \end{Bmatrix}_r
=\sum_{k=2}^n \binom{n}{k}\sum_{l=1}^{k-1}(-1)^{l-1}
\binom{l+r-2}{l-1}\begin{Bmatrix} k-l \\ m-1\end{Bmatrix}_{r-1}
\ee
where $\left\{\right\}_r$ are $r$-Stirling numbers of the second kind,
namely \cite{Yaqubiarxiv1412}
\be
\begin{Bmatrix} n\\ m \end{Bmatrix}_r=\begin{cases}
0, & n<r \\
\delta_{mr}& n=r\\
m\begin{Bmatrix}n-1 \\ m\end{Bmatrix}_r+\begin{Bmatrix} n-1 \\ m-1\end{Bmatrix}_r, & n>r.
\end{cases}
\ee

\cite{Boyadzhievarxiv2012}
\be
\sum_{k=0}^n 
\left[\begin{array}{c}n\\k \end{array}\right]
\left\{\begin{array}{c}k+1\\p+1\end{array} \right\}
 = \frac{n!}{p!}\binom{n}{p}.
\ee

\cite{CampbellIJMM16,CampbellIJMM17}
\be
\sum_{k=1}^n a_k\frac{1-\exp(b_kx)}{1-\exp(b_kx/k)}
=(\sum_{k=1}^n a_k)
+\sum_{m=2}^n \sum_{k=1}^{\lfloor n/m\rfloor}
\sum_{1<j\le m, (j,m)=1} a_{mk}\exp(b_{mk}jx/m).
\ee
where $a_k$, $b_k$ are arbitrary sequences.
\cite{CampbellIJMM16}
\be
\sum_{k=1}^\infty a_k \frac{1-q^{kb_k}}{1-q^{b_k}}
=\sum_{k=1}^\infty
[ a_k+\sum_{m=2}^\infty \sum_{1\le j\le m, (j,m)=1} a_{mk}q^{jkb_{mk}}]
.
\ee

\cite{SpiveyDM307}
\be
\sum_{k=0}^n\binom{n}{k}H_k=2^n\left(H_n-\sum_{k=1}^n\frac{1}{k2^k}\right)
.
\ee
where $H_k$ are the harmonic numbers.

\cite{SavioHarm}
\be
\sum_{i=1}^n i^m H_i^{(k)} = H_n^({(k)}\sum_{i=1}^n i^m -\sum_{j=1}^n 1/j^k \sum_{i=1}^{j-1}i^m
\ee

\cite{WangJMAA434}
\be
\sum_{k=1}^n kH_k = \frac{n(n+1)}{2}H_n-\frac{(n-1)n}{4}.
\ee

\cite{WangJMAA434}
\be
\sum_{k=1}^n k^2H_k = \frac{n(n+1)(2n+1)}{6}H_n-\frac{(n-1)n(4n+1)}{36}.
\ee

\cite{WangJMAA434}
\be
\sum_{k=1}^n k^2H_{n+k} = \frac{n(n+1)(2n+1)}{6}(2H_{2n}-H_n)-\frac{n(n+1)(10n-1)}{36}.
\ee

\cite{WangJMAA434}
\be
\sum_{k=1}^n kH_{p+k} = \frac{(n-p)(p+n+1)}{2}H_{p+n}
+\frac{p(p+1)}{2}H_p-\frac{n(n-2p-1)}{4}.
\ee

\cite{WangJMAA434}
\be
\sum_{k=1}^n k^2H_{p+k} = \cdots
\ee

\cite{SavioHarm}
\be
\sum_{i=1}^n \frac{H_i}{i+m} = H_nH_{n+m}-\frac12(H_n^2+H_n^{(2)})-\sum_{i=1}^n\frac1i \sum_{j=1}^{m-1}\frac{1}{i+j}.
\ee

\cite{SavioHarm}
\be
\sum_{i=1+m}^n \frac{H_i}{i-m} = \frac12(H_{n-m}^2+H_{n-m}^{(2)})+\sum_{i=1}^{n-m}\frac1i \sum_{j=1}^{m}\frac{1}{i+j}.
\ee

\cite{SprugnoliEJCN6}
\be
\sum_{k=0}^n (-1)^k(H_k-2H_{2k})\binom{n}{k} = \frac{4^n}{n\binom{2n}{n}}.
\ee

\cite{WangRJ45}
\be
\sum_{k=1}^n H_k^2=(n+1)H_n^2-(2n+1)H_n+2n.
\ee

\cite{WangRJ45}
\be
\sum_{k=1}^n H_k^3=(n+1)H_n^3-\frac{3(2n+1)}{2}H_n^2+3(2n+1)H_n+\frac12 H_n^{(2)}-6n.
\ee

\cite{WangRJ45}
\be
\sum_{k=1}^n H_{p+k}^2=(p+n+1)H_{p+n}^2-(p+1)H_p^2-(2p+2n+1)H_{p+n}+(2p+1)H_p+2n.
\ee

\cite{WangJMAA434}
\be
\sum_{k=1}^n (-)^k\binom{n}{k}\frac{H_{n+k}}{2^k} = \frac{H_n}{2^{n+1}}.
\ee

\cite{WangJMAA434}
\be
\sum_{k=1}^{1+2n} (-)^k\binom{1+2n}{k}\binom{1+2n+k}{k}\frac{H_k}{2^k} = (-)^{n+1}\frac{n!^24^n}{(1+2n)!}.
\ee

\cite{WeiArxiv1709}
\begin{multline}
\sum_{k=0}^n (-)^k\binom{n}{k}\frac{\binom{n+k}{k}}{\binom{x+k}{k}}H_n^{(2)}(x)
=\frac{(-)^n}{2}\frac{\binom{n-x}{n}}{\binom{n+x}{n}}
\{
[H_n^{(2)}(x)-H_n^{(2)}(-x)]
+\frac{4n}{x(x-n)^2}
\\
+[H_n(x)-H_n(-x)-H_n(x/2-n/2)][H_n(x)-H_n(x)-H_n(x/2-n/2)-\frac{2(x+n)}{x(x-n)}]
\}
.
\end{multline}

\cite{WeiArxiv1709}
\be
\sum_{k=0}^n (-)^k\binom{n+k}{k}\binom{p+n}{n-k}
H_{p+k}^{(2)}=\cdots
\ee

\subsection{Numerical Series}

\cite{Silagarxiv1207}\cite[A050970,A068205]{sloane}
\be
\sum_{k=-\infty}^\infty \frac{1}{(4k+1)^n} = (\frac{\pi}{2})^n \delta_n
\ee
where $\delta_n$ is the volume of the n-dimensional convex polytope:
\be
\delta_n= \int_0^1 du_1\int_0^1du_2 K_1(u_1,u_2)\cdots \int_0^1 du_n K_2(u_{n-1},u_n)K_1(u_n,u_1)
\ee
and $K_1(u,v)=\theta(1-u-v)$, $K_n(u,v)=\int_0^1 K_1(u,u_1)K_{n-1}(u1,v) du_1$,
are essentially provided by the Euler polynomials.

\cite{Cvijovicarxiv0911}
Let $\chi_x(z) = \sum_{k=0}^\infty z^{2k+1}/(2k+1)^s$. then
\be
\chi_{2n+1}(z) = (-)^n \frac{\pi^{2n+1}}{\delta 4 (2n)!}
\int_0^\delta E_{2n}(t) \frac{2z(1+z^2)\sin(\pi t)}{1-2z^2\cos(2\pi t)+z^4}dt
\ee
\be
\chi_{2n}(z) = (-)^n \frac{\pi^{2n}}{\delta 4 (2n-1)!}
\int_0^\delta E_{2n-1}(t) \frac{2z(1-z^2)\cos(\pi t)}{1-2z^2\cos(2\pi t)+z^4}dt
\ee
where $|z|<1$, $\delta=1$ or $1/2$, and $E$ are Euler polynomials. Special cases
lead to integral representations of
$\lambda(s)=\sum_{k\ge 0} 1/(2k+1)^s$
and
$\beta(s)=\sum_{k\ge 0} (-)^k/(2k+1)^s$.

\cite{Boyadzhievarxiv2012}
\be
\sum_{n=0}^\infty 
\frac{n!}{z(z+1)\cdots(z+n)} = \frac{1}{z-1}.
\ee

\cite{Boyadzhievarxiv2012}
\be
\sum_{n=0}^\infty 
\frac{n!n}{z(z+1)\cdots(z+n)} = \frac{1}{(z-1)(z-2)}.
\ee

\cite{Boyadzhievarxiv2012}\cite[24.1.4]{AS}
\be
\sum_{n=0}^\infty 
\frac{n!}{(n+1)z(z+1)\cdots(z+n)} = \sum_{n=0}^\infty \frac{1}{(z+n)^2} = \psi'(z)
\ee
where $\psi$ is the digamma function.

\cite{Boyadzhievarxiv2012}
\be
\sum_{n=0}^\infty 
\frac{1}{z(z+1)\cdots(z+n)}\sum_{k=0}^n\left[\begin{array}{c}n\\k\end{array}\right]
a_k = \sum_{k=0}^\infty \frac{a_k}{z^{k+1}}
\ee

\cite{Boyadzhievarxiv2012}
\be
(p-1)!\sum_{n=0}^\infty 
\left[\begin{array}{c}n\\k\end{array}\right]
\frac{1}{n!(n+m+1)(n+m+2)\cdots (n+m+p)}
 = \sum_{j=0}^m \binom{m}{j}\frac{(-)^j}{(j+p)^{k+1}}.
\ee

\cite{Boyadzhievarxiv2012}
\be
\sum_{n=0}^\infty \frac{(-)^nd_n}{z(z+1)\cdots(z+n)}
=\sum_{k=0}^\infty \frac{1}{(k+1)z^{k+1}}
\ee
where $d_n$ are the Cauchy numbers of the second kind \eqref{eq.dnBoyad}.

\cite{MezoArxiv0811}
\be
S(r,m)=S(1,m)+\sum_{k=1}^{r-1}\frac{S(k,m-1)-B(k,m)}{k}.
\ee
where $H_n^{(r)}\equiv \sum_{k=1}^n\frac{1}{k^r}$,
where $S(r,m)\equiv \sum_{n=1}^\infty \frac{H_n^{(r)}}{n^m}$,
where $B(k,m)\equiv \,_{m+1}F_m(1,1,\ldots ,1,k+1;2,2,\ldots,2;1)$ .

\cite{MezoArxiv0811}
\be
S(2,3)=\frac{\pi^4}{72}-\frac{\pi^2}{6}+2\zeta(3)
\approx 2.1120837816098848 .
\ee

\cite{MezoArxiv0811}
\be
S(2,4)=\frac{\pi^4}{72}+3\zeta(5)-\zeta(3)\left(1+\frac{\pi^2}{6}\right).
\ee

\cite{CoffeyJCAM183}
\be
\sum_{n=1}^\infty \frac{(-)^n}{n(n+1)^{j+1}}
=\sum_{m=1}^j(2^{-m}-1)\zeta(m+1)+j+1-2\ln 2.
\ee

\cite{CoffeyJCAM183}
\be
\sum_{n=1}^\infty \frac{(-)^n}{n^{j+1}(n+1)}
=(-)^j [\sum_{m=1}^j(-)^m(2^{-m}-1)\zeta(m+1)+1-2\ln 2].
\ee

\cite{CoffeyJCAM183}
\begin{multline}
\sum_{n=1}^\infty \frac{(-)^n}{n(n+z)^{j+1}}
=\sum_{m=1}^j
\frac{(-)^m}{m!z^{j-m+1}}
[
\frac{1}{2^m}\psi^{(m)}(z/2)+(-1/2)^m m!(z/2)^{-m-1}-\psi^{(m)}(z)-(-)^mm!z^{-1-m}
]
\\
+\frac{1}{z^{j+1}}
[
\psi(z/2)-\psi(z)+1/z
].
\end{multline}

\cite{SaqabiJMAA159}
Suppose $P(x)$ is a polynomial in $x$ of degree $\le r-2$ with partial fractions
\be
\frac{xP(x)}{(\lambda_1+\mu_1 x)\cdots(\lambda_r+\mu_r x)}=\sum_{k=1}^r \frac{A_k}{\lambda_k+\mu_k x},
\ee
then
\be
\sum_{n\ge 1}\frac{P(n)}{(\lambda_1+\mu_1 n)\cdots(\lambda_r+\mu_r n)}
=-\sum_{k=1}^r \mu_k^{r-2}P(-\lambda_k/\mu_k)\psi(1+\lambda_k/\mu_k)\prod_{j=1,j\neq k}^r \frac{1}{\lambda_j\mu_k-\lambda_k\mu_j}
\ee
provided that all the denominators on the rhs are nonzero.
\cite{SaqabiJMAA159}
\be
\sum_{n\ge 1}(-)^n \frac{Q(n)}{(\lambda_1+\mu_1 n)\cdots(\lambda_r+\mu_r n)}
=-\sum_{k=1}^r \mu_k^{r-2}Q(-\lambda_k/\mu_k)\beta(1+\lambda_k/\mu_k)\prod_{j=1,j\neq k}^r \frac{1}{\lambda_j\mu_k-\lambda_k\mu_j}
\ee
where $\beta(z)=\frac12[\psi((1+z)/2)-\psi(z/2)]=\sum_{N\ge 0} (-1)^N/(z+n)$.

\cite{EspinosaMCOM79}
\be
\sum_{-\infty}^\infty \frac{1}{n^2+q^2}=\frac{\pi \coth \pi q}{q}.
\ee

\cite{EspinosaMCOM79}
\begin{multline}
\sum_{-\infty}^\infty \frac{1}{(n^2+q_1^2)((N-n)^2+q_2^2}
=\frac{2\pi}{2q_12q_2}
\big[
\frac{1+n_b(q_1)+n_b(q_2)}{Ni+q_1+q_2}
+\frac{n_b(q_1)-n_b(q_2)}{Ni-q_1+q_2}
\\
-\frac{n_b(q_1)-n_b(q_2)}{Ni+q_1-q_2}
-\frac{1+n_b(q_1)+n_b(q_2)}{Ni-q_1-q_2}
\big]
\end{multline}
where $n_b(z)\equiv 1/(e^{2\pi z}-1)$.

\cite{EspinosaMCOM79}
Define the Matsubara sum associated to the Graph
$G$ (loopless multigraph such that the degree of each 
vertex is at least 2, with $I$ lines) by
\be
S_G\equiv \sum_{n_1=-\infty}^\infty
\sum_{n_2=-\infty}^\infty
\cdots
\sum_{n_I=-\infty}^\infty
\frac{\delta_g(n_1,n_2,\cdots n_I;\{N_v\})}
{(n_1^2+q_1^2)(n_2^2+q_2^2)\cdots (n_I^2+q_I^2)},
\ee
where
\be
\delta_g(n_1,\ldots;\{N_v\})=\prod_{v=1}^V \delta_{T_v,N_v}
\ee
imposes a series of constraints over the vertices $v$,
and $T_v\equiv \sum_i s_i^{v}n_i$ is an algebraic
sum at vertex $v$, with $s_i^{v}$ having values $\pm 1$ or $0$
depending on the orientation of the line $i$ with respect to the vertex $v$.
The $q_i$ are weights associated with the lines $i$.
Then the integral
\be
I(N,q_1,q_2,\ldots)
\int_{-\infty}^\infty dx_1
\int_{-\infty}^\infty dx_2
\cdots \frac{1}{(x_1^2+q_1^2)(x_2^2+q_2^2)\cdots ((N-x_1-x_2\cdots)^2+q_I^2)}
\ee
is related to the sum via
\be
S_G =\hat O_G I_G
\ee
where the operator $\hat O_G=\prod_{i=1}^I [1+n_{b_i}(1-\hat R_i)]$
is composed of the functions $n_b$ of the previous formula and the reflection operator $\hat R_i$
(which switches the sign of the variable $q_i$).

\cite{ErdosSz7}
\be
\sum_{n=2}^\infty \frac{1}{n(n-1)}=\sum_{k\ge 2} [\zeta(k)-1].
\ee

\cite[A152416]{sloane}\cite{MatharArxiv0803}
\be
\sum_{n=2}^\infty \frac{1}{n^s(n-1)}=s-\sum_{l=2}^s \zeta(l).
\ee

\cite{CoffeyJCAM183}
\be
\sum_{n=1}^\infty \frac{1}{n^{j+1}(n+1)}=(-1)^j[\sum_{m=1}^j (-)^m\zeta(m+1)+1].
\ee

\cite{CoffeyJCAM183}
\be
\sum_{n=1}^\infty \frac{1}{n(n+1)^{j+1}}=-\sum_{m=1}^j \zeta(m+1)+j+1.
\ee

\cite{BorweinJMAA316}
\be
\sum_{n\ge 1}\frac{1}{n(n^2+1)} = \gamma +\Re\psi(1+i)
\approx 0.67186598552400983.
\ee

\cite{BorweinJMAA316}
\be
\sum_{n\ge 1}\frac{1}{n^2(n^2+1)} = \frac{\pi^2}{6}-\frac{\pi\coth \pi -1}{2}
\approx 0.56826001937964526.
\ee

\cite{PilehroodEJC18}
\be
\sum_{n\ge 1}\frac{(-)^{n-1}}{n^2-a^2}
=
\frac14\sum_{n\ge 1}\frac{(-1)^{n-1}(10n^2-3n-a^2)}{n(2n-1)(n^2-a^2)\binom{2n}{n}\prod_{j=1}^n(1-a^2/(n+j)^2)}.
\ee

\cite{PilehroodEJC18}
\be
\sum_{k\ge 0}\frac{k+\alpha}{((k+\alpha)^2-a^2)((k+\alpha)^2-b^2)}
=\frac12\sum_{n\ge 1}\frac{(-)^{n-1}(1\pm a\pm b)_{n-1}
(5n^2-6n(1-\alpha)+2(1-\alpha)^2-a^2-b^2)}{n\binom{2n}{n}(\alpha\pm a)_n(\alpha\pm b)_n}
\ee

\cite{PilehroodEJC18}
\be
\sum_{k\ge 0}\frac{k+\alpha}{(k+\alpha)^4-x^2(k+\alpha^2)-y^4} =\ldots
\ee

\cite{CooperBAMS71}
Let
\begin{multline}
a(q)\equiv \sum_m \sum_n q^{m^2+mn+n^2},
b(q)\equiv \sum_m \sum_n q^{m^2+mn+n^2}\omega^{m-n},\\
c(q)\equiv \sum_m \sum_n q^{(m+1/3)^2+(m+1/3)(n+1/3)+(n+1/3)^2},
\end{multline}
where $\omega \equiv \exp(2\pi i/3)$, $q=e^{-2\pi t}$, $\Re t>0$, then
\be
a(q)^3=b(q)^3+c(q)^3,
\ee
\be
b(q)=\prod_{n=1}^\infty \frac{(1-q^n)^3}{1-q^{3n}},
\ee
\be
c(q)=3q^{1/3}\prod_{n=1}^\infty \frac{(1-q^{3n})^3}{1-q^n},
\ee
\be
a(q)=a(q^3)+2c(q^3) = 1+6\sum_{n=1}^\infty[\frac{q^{3n-2}}{1-q^{3n-2}}-\frac{q^{3n-1}}{1-q^{3n-1}}],
\ee
\be
b(q)=a(q^3)-c(q^3),
\ee
\be
a(q)^2 =1+12\sum_{n=1,3\nmid n}^\infty \frac{nq^n}{1-q^n}.
\ee

\cite{MatharArxiv1301}
\be
\sum_{k=1}^\infty \frac{1}{(2k)^{2s}(2k+1)^{2s}}
=
\sum_{t=1}^{2s}\binom{4s-t-1}{2s-1}\left\{
[1-\frac{1-(-)^t}{2^t}]\zeta(t)-1
\right\}
.
\ee

\cite{MatharArxiv1301}
\be
\sum_{k=1}^\infty \frac{1}{(2k)^2(2k+1)^2}
=
-3+\frac{\pi^2}{6}+2\log 2
.
\ee

\cite{MatharArxiv1301}
\be
\sum_{k=1}^\infty \frac{1}{(2k)^4(2k+1)^4}
=
-35+\frac{\pi^4}{90}+3\zeta(3)+\frac{5\pi^2}{3}+20\log 2
.
\ee

\cite{MatharArxiv1301}
\be
\sum_{n=1}^\infty \frac{1}{n^{2s}(n+1)^{2s}}
=
\sum_{t=1}^{2s}\binom{4s-t-1}{2s-1}\left\{
[1+(-)^t]\zeta(t)-1
\right\}
.
\ee

\cite{MatharArxiv1301}
\be
\sum_{n=1}^\infty \frac{1}{n^2(n+1)^2}
=
\frac{\pi^2}{3}-3
.
\ee

\cite{MatharArxiv1301}
\be
\sum_{n=1}^\infty \frac{1}{n^4(n+1)^4}
=
-35+\frac{10\pi^2}{3}+\frac{\pi^4}{45}
.
\ee

\cite{SofoIJPAM50}
\be
\sum_{n=1}^\infty \frac{1}{(4n^2-1)(16n^2-1)(16n^2-9)}
=\frac{22/7-\pi}{60}.
\ee

\cite[(387)]{Jolley}
\be
\sum_{n\ge 1} \frac{(2n-1)!!}{(2n)!!}\frac{1}{2n+1}
=
\frac{1}{6}{}_3F_2(1,3/2,3/2; 2,5/2; 1)
=
\frac{\pi}{2}-1.
\ee

\begin{multline}
\sum_{k\ge 1}\frac{(2k-1)!!}{(2k)!!k}
=
\sum_{k\ge 1}\frac{(2k)!}{(2^kk!)^2k}
=
\sum_{k\ge 1}\frac{\Gamma(2k+1)}{2^{2k}\Gamma(k+1) k!k}
=
\frac{1}{\sqrt{\pi}}\sum_{k\ge 1}\frac{\Gamma(k+1/2)}{ k!k}
=
\frac{1}{\sqrt{\pi}}\sum_{k\ge 0}\frac{\Gamma(k+3/2)}{ (k+1)!(k+1)}
\\
=
\frac{1}{\sqrt{\pi}}\sum_{k\ge 0}\frac{\Gamma(k+3/2)}{ (k+1)^2k!}
=
\frac{\Gamma(3/2)}{\sqrt{\pi}}\sum_{k\ge 0}\frac{(3/2)_k}{ (k+1)^2k!},
\end{multline}
then replace $1/(k+1)^2$ via \eqref{eq.poch1},
insert $b=c=1$, $a=1/2$, $e=3/2$ in \eqref{eq.kornw}, then \eqref{eq.drive} and \cite[(15.1.3)]{AS}
\begin{multline*}
=
\frac{\Gamma(3/2)}{\sqrt{\pi}}\sum_{k\ge 0}\frac{(3/2)_k(1)_k(1)_k}{(2)_k(2)_kk!}
=
\frac{1}{2}{}_3F_2(3/2,1,1;2,2;1)
=2\ln 2.
\end{multline*}

\cite{YangIJMEST23}
\be
\sum_{n\ge 1} [\frac{(2n-1)!!}{(2n)!!}]^2 \frac{1}{2n-6}
=\frac{25}{128}\ln 2+\frac{71}{1536}-\frac{43+50G}{128\pi} = -0.039218...
\ee
where $G$ is Catalan's constant.

\cite{YangIJMEST23}
\be
\sum_{n\ge 1} [\frac{(2n-1)!!}{(2n)!!}]^2 \frac{1}{2n-4}
=\frac{9}{32}\ln 2+\frac{11}{128}-\frac{13+18G}{32\pi} = -0.012431...
\ee

\cite{YangIJMEST23}
\be
\sum_{n\ge 1} [\frac{(2n-1)!!}{(2n)!!}]^2 \frac{1}{2n-2}
=\frac{1}{2}\ln 2+\frac{1}{4}-\frac{1+2G}{2\pi} = 0.1458577...
\ee

\cite{YangIJMEST23}\cite[A370374]{sloane}
\be
\sum_{n\ge 1} [\frac{(2n-1)!!}{(2n)!!}]^2 \frac{1}{2n}
=2\ln 2-\frac{4G}{\pi} = 0.2200507...
\ee

\cite{YangIJMEST23}\cite[A218387]{sloane}
\be
\sum_{n\ge 1} [\frac{(2n-1)!!}{(2n)!!}]^2 \frac{1}{2n+1}
=\frac{4G}{\pi}-1 = 0.166243...
\ee

\cite{YangIJMEST23}
\be
\sum_{n\ge 1} [\frac{(2n-1)!!}{(2n)!!}]^2 \frac{1}{2n+3}
=\frac{2G+1}{2\pi}-\frac13 = 0.1173825...
\ee

\cite{YangIJMEST23}
\be
\sum_{n\ge 1} [\frac{(2n-1)!!}{(2n)!!}]^2 \frac{1}{2n+5}
=\frac{18G+13}{32\pi}-\frac15 = 0.0933164....
\ee

\cite[(1)]{LehmerAMM92}
\be
\sum_{n=0}^\infty \frac{\binom{2n}{n}}{8^n}=\surd 2.
\ee
\cite[(1)]{LehmerAMM92}
\be
\sum_{n=0}^\infty \frac{\binom{2n}{n}}{10^n}=\sqrt{5/3}.
\ee
\cite[(1)]{LehmerAMM92}
\be
\sum_{n=0}^\infty \frac{(-1)^n\binom{2n}{n}}{8^n}=\sqrt{2/3}.
\ee
\cite{LehmerAMM92}\cite[A145439]{sloane}
\be
\sum_{n=0}^\infty \frac{\binom{4n}{2n}}{64^n}=\frac{3\surd 2+\surd 6}{6}
\ee
\cite{LehmerAMM92}
\be
4 \sum_{n=0}^\infty \frac{\binom{8n}{4n}}{8^{4n}}
=\frac{3\surd 2+\surd 6}{3}
+\frac{2\sqrt{2+\surd 5}}{\surd 5}
.
\ee
\cite[(6)]{LehmerAMM92}
\be
\sum_{n=1}^\infty \frac{\binom{2n}{n}}{n4^n}
=\log 4
.
\ee
\cite[(6)]{LehmerAMM92}\cite[A157699]{sloane}
\be
\sum_{n=1}^\infty \frac{(-1)^{n+1}\binom{2n}{n}}{n4^n}
=2\log \frac{1+\surd 2}{2}
,
\ee
The formula above corrects a factor 2 in \cite{LehmerAMM92}, see \cite{MatharArxiv0905}.

\cite[A091648]{sloane}
\be
\sum_{n=1,3,5,7,\ldots}^\infty \frac{\binom{2n}{n}}{n4^n}
=\log (1+\surd 2)
.
\ee
\cite[(7)]{LehmerAMM92}
\be
\sum_{n=1}^\infty \frac{\binom{2n}{n}}{n(n+1)4^n}
=\log 4 -1
.
\ee
\cite[(8)]{LehmerAMM92}
\be
\sum_{n=1}^\infty \frac{n\binom{2n}{n}}{8^n}
=1/\surd 2
.
\ee
\cite{LehmerAMM92}
\be
\sum_{n=1}^\infty \frac{n^2\binom{2n}{n}}{8^n}
=\frac{5\surd 2}{4}
.
\ee
\cite{LehmerAMM92}\cite[A019670]{sloane}
\be
\sum_{n=0}^\infty \frac{\binom{2n}{n}}{(2n+1)16^n}
=\frac{\pi}{3}
.
\ee

\cite[(15)]{LehmerAMM92}
\cite{DysonArXiv1009,TrifFQ38,SprugnoliEJCN6}\cite[A073016]{sloane}
\be
\sum_{m=1}^\infty \frac{1}{\binom{2m}{m}}
=\frac{9+2\pi\surd 3}{27}
.
\label{eq.inv2mm}
\ee

\cite{SofoIJPAM50}
\be
\sum_{n=1}^\infty \frac{1}{\binom{nk}{k}} 
=
\sum_{n=1}^\infty \frac{k!(nk-k+j)!}{(nk+j)!}
=
\sum_{r=1}^k (-)^r \binom{k-1}{r-1}\psi(\frac{r}{k}).
\ee

\cite{SprugnoliEJCN6}\cite[A307086]{sloane}
\be
\sum_{n=0}^\infty \frac{(-1)^n}{\binom{2n}{n}}
=\frac45 - \frac{4\surd 5}{25} \ln \phi;\quad \phi\equiv \frac{1+\surd 5}{2}.
\ee

\cite[A086466]{sloane}
\be
\sum_{m=1}^\infty \frac{(-1)^{m-1}}{m\binom{2m}{m}}
=\frac{2}{\surd 5}\log\frac{1+\surd 5}{2} .
\ee
This corrects a factor 2 in \cite{LehmerAMM92} and two typos
in \cite[4.1.42]{Apelblat2}, see \cite{MatharArxiv0905}.

\cite{SofoIJPAM50}
\be
\sum_{n=1}^\infty \frac{1}{\binom{nk+j}{k}} =
\sum_{r=1}^k (-)^r \binom{k-1}{r-1}\psi(\frac{r+j}{k}).
\ee

\cite{TrifFQ38}
\be
\sum_{k=0}^{\infty} \frac{1}{\binom{mk}{nk}}
= \int_0^1\frac{1+(m-1)t^n(1-t)^{m-n}}{(1-t^n(1-t)^{m-n})^2}dt,
\ee
where $m$ and $n$ are positive integers with $m>n$.

\cite{TrifFQ38,SprugnoliEJCN6}
\be
\sum_{k=0}^{\infty} \frac{1}{\binom{4k}{2k}}
= \frac{16}{15}+\frac{\pi \sqrt{3}}{27}
-\frac{2\sqrt{5}}{25}\ln\frac{1+\sqrt{5}}{2}.
\ee

\cite{LehmerAMM92,SprugnoliEJCN6}\cite[A086465]{sloane}
\be
\sum_{n=1}^\infty \frac{(-1)^{n-1}}{\binom{2n}{n}}
=\frac{1}{5}+\frac{4\surd 5}{25}\log\frac{1+\surd 5}{2}
.
\ee
\cite{LehmerAMM92,DysonArXiv1009}\cite[A145429]{sloane}
\be
\sum_{m=1}^\infty \frac{m}{\binom{2m}{m}}
=\frac{2}{27}(\pi\surd 3+9)
.
\ee
\cite{BorweinCECM,BorweinAM70}
And by linear combination with (\ref{eq.inv2mm}):
\be
\sum_{n\ge 1}\frac{18-9n}{\binom{2n}{n}}=2\frac{\pi}{\surd 3}.
\ee

\cite{LehmerAMM92,DysonArXiv1009}
\be
\sum_{m=1}^\infty \frac{m^2}{\binom{2m}{m}}
=\frac{2}{81}(5\pi\surd 3+54)
.
\ee
\cite{LehmerAMM92,DysonArXiv1009}
\be
\sum_{m=1}^\infty \frac{m^3}{\binom{2m}{m}}
=\frac{2}{243}(37\pi\surd 3+405)
.
\ee
\cite{LehmerAMM92}\cite[A145433]{sloane}
\be
\sum_{m=1}^\infty \frac{(-1)^{m-1}m}{\binom{2m}{m}}
=\frac{2}{125}(2\sigma+15),\quad \sigma\equiv \surd 5\log\frac{1+\surd 5}{2}
.
\ee

\cite{LehmerAMM92}
\be
\sum_{m=1}^\infty \frac{(-1)^{m-1}m^2}{\binom{2m}{m}}
=\frac{2}{125}(5-\sigma)
.
\ee

\cite{SprugnoliEJCN6,LehmerAMM92}
\be
\sum_{n\ge 0}\frac{2^n}{\binom{2n}{n}} = \frac{\pi}{2}+2;\quad
\sum_{m=1}^\infty \frac{2^m}{m\binom{2m}{m}}
=\frac{\pi}{2}
\ee

\cite{SprugnoliEJCN6}
\be
\sum_{n\ge 0}\frac{(-2)^n}{\binom{2n}{n}} = -\frac{2\surd 3}{9}\ln\frac{\surd 3+1}{\surd 2}+\frac{2}{3}.
\ee

\cite{SprugnoliEJCN6}
\be
\sum_{n\ge 1}\frac{4^n}{n^2\binom{2n}{n}} = \frac{\pi^2}{2}.
\ee

\cite{SprugnoliEJCN6}
\be
\sum_{n\ge 2}\frac{4^n}{(n-1)^2\binom{2n}{n}} = \pi^2-4 .
\ee

\cite{SprugnoliEJCN6}
\be
\sum_{n\ge 1}\frac{(-4)^n}{n\binom{2n}{n}} = \sqrt{2}\ln(\sqrt{2}-1).
\ee

\cite{SprugnoliEJCN6}
\be
\sum_{n\ge 1}\frac{(-4)^n}{n^2\binom{2n}{n}} = -2[\ln(\sqrt{2}-1)]^2.
\ee

\cite{SprugnoliEJCN6}
\be
\sum_{n\ge 1}\frac{(-4)^n}{(2n+1)\binom{2n}{n}} = -\frac{\ln(\sqrt{2}-1)}{\surd 2}.
\ee

\cite{SprugnoliEJCN6}
\be
\sum_{n\ge 0}\frac{(-4)^n}{(n+1)\binom{2n}{n}} = -\sqrt{2}\ln(\surd 2-1)-[\ln(\surd 2-1]^2.
\ee

\cite{SprugnoliEJCN6}
\be
\sum_{n\ge 2}\frac{(-4)^n}{(n-1)\binom{2n}{n}} = -3\sqrt{2}\ln(\surd 2-1)-2.
\ee

\cite{SprugnoliEJCN6}
\be
\sum_{n\ge 2}\frac{(-4)^n}{(n-1)^2\binom{2n}{n}} = 4[\surd 2\ln(\surd 2-1)+(\ln(\surd 2-1))^2+1].
\ee

Erratum to \cite[4.1.40]{Apelblat2}:
\be
\sum_{n=0}^\infty\frac{(-1)^{n-1} n^2}{\binom{2n}{n}}=\frac{4}{125}\left[5-\sqrt{5}\ln\left(\frac{1+\surd 5}{2}\right)\right].
\ee
\be
\sum_{m=1}^\infty \frac{(-1)^{m-1}m^3}{\binom{2m}{m}}
=\frac{2}{625}(14\sigma-5)
,
\ee
which corrects a factor 2 in \cite{LehmerAMM92}, see \cite{MatharArxiv0905}.

\cite{LehmerAMM92,DysonArXiv1009}
\be
\sum_{m=1}^\infty \frac{2^m}{m^2\binom{2m}{m}}
=\frac{\pi^2}{8}
,
\ee
a special case of (\ref{eq.2xbin}).
\cite{LehmerAMM92}
\be
\sum_{m=1}^\infty \frac{2^m}{\binom{2m}{m}}
=\frac{\pi}{2}+1
.
\ee
\cite{LehmerAMM92}
\be
\sum_{m=1}^\infty \frac{m2^m}{\binom{2m}{m}}
=\pi+3
.
\ee
\cite{LehmerAMM92,DysonArXiv1009}
\be
\sum_{m=1}^\infty \frac{m^22^m}{\binom{2m}{m}}
=
\frac{1}{5}\sum_{m=1}^\infty \frac{m^32^m}{\binom{2m}{m}}
=\frac{7\pi}{2}+11
.
\ee
\cite{LehmerAMM92,DysonArXiv1009}
\be
\sum_{m=1}^\infty \frac{m^42^m}{\binom{2m}{m}}
=113\pi+355
.
\ee
\cite{LehmerAMM92,DysonArXiv1009}
\be
\sum_{m=1}^\infty \frac{m^{10}2^m}{\binom{2m}{m}}
=229093376\pi+719718067
.
\ee
\cite{LehmerAMM92}
\be
\sum_{m=1}^\infty \frac{3^m}{m^2\binom{2m}{m}}
=\frac{2\pi^2}{9}
,
\ee
a special case of (\ref{eq.2xbin}).

\cite{LehmerAMM92}\cite[A186706]{sloane}
\be
\sum_{m=1}^\infty \frac{3^m}{m\binom{2m}{m}}
=\frac{2\pi}{\surd 3}\equiv \nu
.
\ee
\cite{LehmerAMM92}
\be
\sum_{m=1}^\infty \frac{3^m}{\binom{2m}{m}}
=2\nu+3
.
\ee
\cite{LehmerAMM92}
\be
\sum_{m=1}^\infty \frac{m3^m}{\binom{2m}{m}}
=10\nu+18
.
\ee
\cite{LehmerAMM92}
\be
\sum_{m=1}^\infty \frac{m^23^m}{\binom{2m}{m}}
=2(43\nu+78)
.
\ee
\cite{LehmerAMM92}
\be
\sum_{m=1}^\infty \frac{(-1)^{m-1}2^m}{m\binom{2m}{m}}
=\rho/3,\quad \rho\equiv \sqrt{3}\log(2+\surd 3)
.
\ee

\cite{AmghiJIS07}
\begin{multline}
\sum_{m=1}^\infty \frac{(2t)^{2m+2k}}{m(2m+2k)\binom{2m}{m}}
=
\arcsin^2(t)
\binom{2k}{k}+\sum_{j=1}^k
\binom{2k}{k-j}\frac{(-)^{j+1}}{j^2}
\\
+
\sum_{j=1}^k(-)^j\binom{2k}{k-j}
\left(\frac{2\arcsin t \sin(2j\arcsin t)}{j}
+\frac{\cos(2j\arcsin t)}{j^2}\right)
,
\end{multline}
and a similar logarithmic result for
an alternating sign sum on the left hand side.

\cite{AmghiJIS07}
\begin{multline}
\sum_{m=1}^\infty \frac{(2t)^{2m+2k}}{m^2(2m+2k)
\binom{2m}{m}}
=
-\frac{1}{k}\arcsin^2(t)
+\frac{(2t)^{2k}}k \arcsin^2 t
- \sum_{j=1}^k
\binom{2k}{k-j}\frac{(-)^{j+1}}{kj^2}
\\
-
\sum_{j=1}^k(-)^j\binom{2k}{k-j}
\left(\frac{2\arcsin t \sin(2j\arcsin t)}{kj}
+\frac{\cos(2j\arcsin t)}{kj^2}\right)
,
\end{multline}
and a similar logarithmic result for
an alternating sign sum on the left hand side.

\cite{Dilcherarxiv2512}
\begin{multline}
\sum_{k\ge 0}\frac{\binom{2k}{k}}{4^k}\frac{x^{2k+1}}{2k+n+1}
= \frac{\binom{n}{n/2}}{(2x)^n}\left(\arcsin x -\frac{\sqrt{1-x^2}}{2}\sum_{j=0}^{n/2-1}\frac{(2x)^{2j+1}}{(2j+1)\binom{2j}{j}}\right).
\end{multline}

\cite{LehmerAMM92}
\be
\sum_{m=1}^\infty \frac{(-1)^{m-1}2^m}{\binom{2m}{m}}
=\frac{\rho+3}{9}
.
\ee
\cite{LehmerAMM92}
\be
\sum_{m=1}^\infty \frac{(-1)^{m-1}m2^m}{\binom{2m}{m}}
=\frac{1}{3}
.
\ee
\cite{LehmerAMM92}
\be
\sum_{m=1}^\infty \frac{(-1)^{m-1}m^22^m}{\binom{2m}{m}}
=\frac{1}{27}(3-\rho)
.
\ee
\cite{LehmerAMM92}
\be
\sum_{m=1}^\infty \frac{(-1)^mm^32^m}{\binom{2m}{m}}
=\frac{1}{81}(\rho+15)
.
\ee
\cite{LehmerAMM92}
\be
\sum_{m=1}^\infty \frac{(2-\surd 2)^m}{\binom{2m}{m}}
=\frac{3-2\surd 2}{4}(\pi\surd 2+4)
.
\ee
\cite{LehmerAMM92}
\be
\sum_{m=1}^\infty \frac{(-1)^{m-1}3^{2m}}{4^m\binom{2m}{m}}
=\frac{48}{125}(\log 2+\frac{15}{16})
.
\ee
\cite{LehmerAMM92}\cite[A152422]{sloane}
\be
\sum_{m=1}^\infty \frac{2^m(2-\surd 3)^m}{m^2\binom{2m}{m}}
=
2\left(\arcsin \tau\right)^2,
\quad \tau = \frac{\sqrt{3}-1}{2} = \sqrt{2}\sin\frac{\pi}{12}
=\sin\frac{\pi}{3}
-\sin\frac{\pi}{6},
\ee
where both right hand sides in \cite{LehmerAMM92} are erroneous, see \cite{MatharArxiv0905}.

\cite{BorweinCECM,BorweinAM70,Almkvistarxiv0110}
\be
\sum_{n\ge 0}\frac{50n-6}{\binom{3n}{n}2^n}=\pi.
\ee
\cite{BorweinCECM,BorweinAM70}
\be
\sum_{n\ge 1}\frac{1}{\binom{3n}{n}2^n}=\frac{2}{25}-\frac{6}{125}
\ln 2+\frac{11}{250}\pi.
\ee
\cite{BorweinCECM,BorweinAM70}
\be
\sum_{n\ge 1}\frac{n}{\binom{3n}{n}2^n}=\frac{81}{625}-\frac{18}{3125}
\ln 2+\frac{79}{3125}\pi.
\ee
\cite{BorweinCECM,BorweinAM70}
\be
\sum_{n\ge 1}\frac{n^2}{\binom{3n}{n}2^n}=\frac{561}{3125}+\frac{42}{15625}
\ln 2+\frac{673}{31250}\pi.
\ee
\cite{BorweinCECM,BorweinAM70}
\be
\sum_{n\ge 1}\frac{-150n^2+230n-36}{\binom{3n}{n}2^n}=\pi.
\ee
\cite{BorweinCECM,BorweinAM70}
\be
\sum_{n\ge 1}\frac{575n^2-965n+273}{\binom{3n}{n}2^n}=6\log 2
.
\ee
\cite{BorweinCECM,BorweinAM70}
\be
\sum_{n\ge 1}\frac{(-1)^n}{\binom{3n}{n}4^n}
=-\frac{1}{28}-\frac{3}{32}\ln 2+\frac{13}{112}\frac{\arctan(\surd 7/5)}{\surd 7}
.
\ee
\cite{BorweinCECM,BorweinAM70}
\be
\sum_{n\ge 1}\frac{(-1)^nn}{\binom{3n}{n}4^n}
=-\frac{81}{1568}-\frac{9}{256}\ln 2+\frac{17}{6272}
\frac{\arctan(\surd 7/5)}{\surd 7}
.
\ee
\cite{BorweinCECM,BorweinAM70}
\be
\sum_{n\ge 1}\frac{1}{n\binom{3n}{n}2^n}
=\frac{1}{10}\pi-\frac{1}{5}\ln 2.
\ee
\cite{BorweinCECM,BorweinAM70}
\be
\sum_{n\ge 1}\frac{1}{n^2\binom{3n}{n}2^n}
=\frac{1}{24}\pi^2-\frac{1}{2}\ln^2 2.
\ee

\cite{BorosSci12}
\be
\sum_{j=1}^\infty \frac{(-1)^j}{j!}(j/2)!
= 1-\frac{1}{2}e^{1/4}\sqrt{\pi}(1-\erf(1/2))
.
\ee

Erratum to \cite[4.1.47]{Apelblat2}\cite[A145438]{sloane}:
Using twice 9.14.+13, then \cite[9.121.6]{GR} and \cite[9.121.26]{GR},
then substituting $tt'/4=z^2$, then $y=\sqrt{t'/4}$, then $\arcsin y =v$
\begin{gather*}
\sum_{n=1}^\infty\frac{1}{n^3\binom{2n}{n}}
=
\frac{1}{2}\,_4F_3(1,1,1,1;2,2,3/2;1/4)
=
\frac{1}{2}\int_0^1 dt\,_3F_2(1,1,1;2,3/2;t/4)
\\
=
\frac{1}{4}\int_0^1 dt \int_0^1 dt' \,_2F_1(1,1;3/2;tt'/4)
=
\frac{1}{4}\int_0^1 dt \int_0^1 dt' \frac{1}{\sqrt{1-tt'/4}}\,_2F_1(1/2,1/2;3/2;tt'/4)
\\
=
\frac{1}{4}\int_0^1 dt \int_0^1 dt' \frac{1}{\sqrt{1-tt'/4}}\frac{\arcsin(\sqrt{tt'/4})}{\sqrt{tt'/4}}
=
\int_0^1 \frac{dt'}{t'} \int_0^{\sqrt{t'/4}} dz \frac{\arcsin z}{\sqrt{1-z^2}}
=
\frac{1}{2} \int_0^1 \frac{dt'}{t'} \arcsin^2 (\sqrt{t'/4})
\\
=
\int_0^{\sqrt{1/2}} \frac{dy}{y} \arcsin^2 y
=
\int_0^{\sin \sqrt{1/2}} v^2\cot v dv
.
\end{gather*}

This becomes \cite[(35)]{BorweinCECM}\cite{LehmerAMM92}\cite[Theorem 3.3]{BorweinEM10}
\be
\frac{2\pi}{3} \Im \Li_2(e^{i\pi/4})-\frac{4}{3}\zeta(3) =
\frac{2\pi}{3} \Cl_2(\pi/3)-\frac{4}{3}\zeta(3)
=
\frac{\pi\surd 3}{18}\left\{\psi'(1/3)-\psi'(2/3)\right\}-\frac{4\zeta(3)}{3}
.
\ee
So the r.h.s.\ of \cite{LehmerAMM92} is missing a factor 4, and
\cite[4.1.47]{Apelblat} is in addition misleading to imply that the
digamma (instead of the trigamma) functions are in effect \cite{MatharArxiv0905}.

\cite{SofoJIS13}
\begin{multline}
\sum_{n=1}^\infty\frac{t^n}{n^{k+1}\binom{an+j+1}{j+1}}
=\left\{
\begin{array}{ll}
\frac{(j+1)t(-)^k}{k!}\int_0^1\int_0^1\frac{(1-x)^jx^a(\ln y)^k}{1-tx^ay} dxdy, & k\ge 1\\
at\int_0^1\frac{(1-x)^{j+1}x^{a-1}}{1-tx^a} dx, & k= 1\\
\end{array}
\right.
\\
=T_0\,_{a+k+1}F_{a+k}
\Big(
\begin{array}{cc}
1,1,\ldots ,1;(a+1)/a,\ldots,(2a-1)/a \\
2,2,\ldots ,2;(a+j+2)/a,\ldots,(a+j+a+1)/a \\
\end{array}
\mid t
\Big)
\end{multline}
where $T_0=t(j+1)B(j+1,a+1)$.

\cite{DzhumadilJIS13}
\be
\frac{1}{2^m}\sum_{i=1}^\infty \binom{i+1}{2}^{-m} = (-)^{m-1}\binom{2m-1}{m}+(-)^m2\sum_{i=1}^{\lfloor m/2\rfloor}
\binom{2m-2i-1}{m-1}\zeta(2i).
\ee

\cite{BorweinCECM,LehmerAMM92,SprugnoliEJCN6}\cite[A073010]{sloane}
\be
\sum_{n\ge 1}\frac{1}{n\binom{2n}{n}} = \frac{\pi}{3\surd 3}.
\ee

\cite{SprugnoliEJCN6}
\be
\sum_{n\ge 1}\frac{1}{n^2\binom{2n}{n}} = \frac{1}{3}\zeta(2).
\ee
a special case of (\ref{eq.2xbin}).

\cite{Goyanesarxiv09,SprugnoliEJCN6}
\be
\sum_{n\ge 1}\frac{1}{n^4\binom{2n}{n}} = \frac{17}{36}\zeta(4).
\ee
\cite{Goyanesarxiv09}
\be
\sum_{n\ge 1}\frac{1}{n^5\binom{2n}{n}} = 2\pi\Cl_4(\pi/3)
-\frac{19}{3}\zeta(5)+\frac{2}{3}\zeta(3)\zeta(2).
\ee
\cite{Goyanesarxiv09}
\be
\sum_{n\ge 1}\frac{1}{n^6\binom{2n}{n}} = -\frac{4\pi}{3}\Im L_{4,1}(e^{i\pi/3})
+\frac{3341}{1296}\zeta(6)-\frac{4}{3}\zeta^2(3).
\ee
\cite{Goyanesarxiv09}
\be
\sum_{n\ge 1}\frac{(-1)^n}{n\binom{2n}{n}} = -2\frac{\arctanh(1/\surd 5)}{\surd 5}.
\ee

Correction of a sign error in \cite{LehmerAMM92}, see \cite{MatharArxiv0905,SprugnoliEJCN6}:
\be
\sum_{n\ge 1}\frac{(-1)^n}{n^2\binom{2n}{n}}
= -2\left( \ln \frac{\sqrt{5}-1}{2}\right)^2
= -2\left( \ln \frac{\sqrt{5}+1}{2}\right)^2
.
\ee

\cite{SprugnoliEJCN6}
\be
\sum_{n\ge 0}\frac{(-1)^n}{(2n+1)\binom{2n}{n}} = 4\ln \phi/\surd 5;\quad \phi\equiv (1+\surd 5)/2.
\ee

\cite{SprugnoliEJCN6}
\be
\sum_{n\ge 0}\frac{(-1)^n}{(n+1)\binom{2n}{n}} = 8\ln \phi/\surd 5
-4 (\ln\phi)^2
;\quad \phi\equiv (1+\surd 5)/2.
\ee

\cite{SprugnoliEJCN6}
\be
\sum_{n\ge 2}\frac{(-1)^n}{(n-1)\binom{2n}{n}} = 3\ln \phi/\surd 5
-1/2
;\quad \phi\equiv (1+\surd 5)/2.
\ee

\cite{SprugnoliEJCN6}
\be
\sum_{n\ge 2}\frac{(-1)^n}{(n-1)^2\binom{2n}{n}} = 1-5\ln \phi
+(\ln \phi)^2
;\quad \phi\equiv (1+\surd 5)/2.
\ee

\cite{SprugnoliEJCN6}
\be
\sum_{n\ge 2}\frac{1}{n^2(n^2-1)\binom{2n}{n}} = 11/8-\sqrt{3}\pi/4.
\ee

\cite{SprugnoliEJCN6}
\be
\sum_{n\ge 2}\frac{4^n}{n^2(n^2-1)\binom{2n}{n}} = 4-3\pi^2/8.
\ee

\cite{SprugnoliEJCN6}
\be
\sum_{n\ge 2}\frac{(-1)^n}{n^2(n^2-1)\binom{2n}{n}} = 4(\ln\phi)^2-\frac{\sqrt{5}}{2}\ln\phi -3/8.
\ee

\cite{SprugnoliEJCN6}
\be
\sum_{n\ge 2}\frac{(-4)^n}{n^2(n^2-1)\binom{2n}{n}} = \frac52 [\ln(\surd 2-1)]^2-\sqrt{2} \ln(\surd 2-1)-3.
\ee

\cite{BorweinCECM,BorweinAM70}
\be
\sum_{n\ge 1}\frac{(-1)^n}{n^3\binom{2n}{n}} = -\frac{2}{5}\zeta(3).
\ee
\cite{BorweinCECM,BorweinAM70}
\be
\sum_{n\ge 1}\frac{(-1)^n}{n^4\binom{2n}{n}} = -4K_4(\rho)+4K_4(-\rho)
+\frac{1}{2}\ln^4\rho +7\zeta(4),
\quad
K_k(x)\equiv \sum_{r=0}^{k-1} \frac{(-\ln |x|)^r}{r!}L_{k-r}(x).
\ee

\cite{PilehroodEJC18}
\be
\frac12\sum_{n\ge 1}\frac{(-)^{n-1}(10n-3)}{n^2(2n-1)\binom{2n}{n}}=\zeta(2).
\ee

\cite{PilehroodEJC18}
\be
\sum_{n\ge 1}\frac{21n-8}{n^3\binom{2n}{n}^3}=\zeta(2).
\ee

\cite{BorweinCECM,BorweinAM70}
\be
\sum_{n\ge 1}\frac{1}{n^3\binom{3n}{n}2^n} = -\frac{33}{16}\zeta(3)
+\frac{1}{6}\ln^3 2
-\frac{1}{24}\pi^2 \ln 2
+\pi \Im L_2(i)
.
\ee
\cite{BorweinCECM,BorweinAM70}
\be
\sum_{n\ge 1}\frac{1}{n^4\binom{3n}{n}2^n} = -\frac{143}{16}\zeta(3)\ln 2
+\frac{91}{640}\pi^4
-\frac{3}{8}\ln^4 2
+\frac{3}{8}\pi^2\ln^2 2
-8L_4(1/2)-8\Re L_{3,1}(\frac{1+i}{2})
-8\Re L_4(\frac{1+i}{2})
.
\ee
\cite{BorweinCECM,BorweinAM70}
\begin{gather*}
\sum_{n\ge 1}\frac{(-1)^n}{8^nn\binom{6n}{2n}} =
\left(-\frac{1}{3}+\frac{2}{57}\sqrt{114\surd 57 -342}\right)\ln 2
-\frac{1}{114}\sqrt{114\surd 57-342}
\ln\left(13+\surd 57+\sqrt{-30+26\surd 57}\right)
\\
+\frac{1}{57}\sqrt{114\surd 57 +342}
\arctan\frac{\sqrt{2\surd 57+9}}{7}
.
\end{gather*}
\cite{BorweinCECM,BorweinAM70}
\be
\sum_{n\ge 1}\frac{(-1)^n}{2^nn\binom{4n}{n}} =
\sum_{13+12r+2r^3=0}\frac{\ln(r+2)}{r+3}-\frac{6}{7}\ln 3+\frac{1}{7}\ln 2.
\ee
\cite{BorweinCECM,BorweinAM70}
\be
\sum_{n\ge 1}\frac{(-1)^n}{n\binom{3n}{n}} =
\sum_{8+4r+r^3=0}\frac{\ln(r+2)}{r+3}-\ln 2.
\ee

\cite{YangJIS10}
Let $A_k\equiv \sum_{n=0}^\infty \binom{n+k}{n}/\binom{2n}{n}$, then
\be
A_{k+1}=-\frac{2}{3(k+1)}+\frac{2(2k+3)}{3(k+1)}A_k.
\ee

\cite{YangJIS10}
Let $B_k\equiv \sum_{n=0}^\infty (-1)^n \binom{n+k}{n}/\binom{2n}{n}$, then
\be
B_{k+1}=-\frac{2}{5(k+1)}+\frac{2(2k+3)}{5(k+1)}B_k.
\ee
and similar recurrences with an additional factor $n$ or $n^2$ in the denominator
of the sums.

\cite{GuilleraRJ45}
\be
\sum_{k\ge 0} \frac{(a)_k(b)_k}{(2b)_k (1)_k} (\frac{4t}{(1+t)^2})^k
=
(1+t)^{2a}\sum_{k\ge 0}\frac{(a)_k(a-b+1/2)_k}{(b+1/2)_k(1)_k}t^{2k}.
\ee

\cite{GuilleraRJ45}
\be
\sum_{k\ge 0} \frac{(1/3)_k(2/3)_k}{(1)_k^2} x^{3k}
=
(1+2t)\sum_{k\ge 0}\frac{(1/3)_k(2/3)_k}{(1)_k^2}t^{3k};\quad t=\frac{1-\sqrt[3]{1-x^3}}{1+2\sqrt[3]{1-x^3}}.
\ee

\cite{FurduiJIS14}
\be
\sum_{k=1}^\infty (-)^k\frac{\ln k}{k}=\gamma \ln 2 -\frac12 \ln^2 2.
\ee

\cite{Almkvistarxiv0110}
\be
\frac{1}{3^25^27^2}\sum_{n\ge 0}\frac{-89286+3875948n-34970134n^2+110202472n^3-115193600n^4}{\binom{8n}{4n}(-4)^n}=\pi.
\ee

\cite{BorosSSA11}
\be
\sum_{k=1}^\infty \frac{k^2}{k^4+4x^4}
=
\frac{\pi}{4x}\,\frac{\sin 2\pi x -\sinh 2\pi x}{\cos 2\pi x -\cosh 2\pi x}
.
\ee

\cite{BorosSSA11}
\be
\sum_{k=1}^\infty \frac{k^2}{k^4+4}
=
\frac{\pi}{4}\coth \pi
.
\ee

\cite{BorosSSA11}
\be
\sum_{k=1}^\infty \frac{k^2}{k^4+x^2k^2+x^4}
=
\frac{\pi}{2x\surd 3}\,\frac{\sinh \pi x\surd 3 -\sqrt{3}\sin\pi x}{\cosh \pi x \surd 3-\cos \pi x}
.
\ee

\cite{PilehroodDMTCS10}
\be
\sum_{k=1}^\infty \frac{A_0+B_0k+C_0k^2}{(k^2-a^2)(k^2-b^2)}
=
\sum_{n=1}^\infty \frac{d_n}{ \prod_{m=1}^n (m^2-a^2)(m^2-b^2)}
\ee
for $|a|<1$, $|b|<1$, with $d_n$ defined via a recurrence in the reference.

\be
\sum_{n=0}^\infty nx^n=\frac{x}{(1-x)^2}.
\ee

\be
\sum_{n=0}^\infty n^2 x^n=\frac{x(1+x)}{(1-x)^3}.
\ee

\cite[A000578]{sloane} 
\be
\sum_{n=0}^\infty n^3 x^n=\frac{x(1+4x+x^2)}{(1-x)^4}.
\ee

\cite[A000583]{sloane} 
\be
\sum_{n=0}^\infty n^4 x^n=\frac{x(1+x)(1+10x+x^2)}{(1-x)^5}.
\ee

\cite[A000584]{sloane} 
\be
\sum_{n=0}^\infty n^5 x^n=\frac{x(1+26x+66x^2+26x^3+x^4)}{(1-x)^6}.
\ee

\cite[A001014]{sloane} 
\be
\sum_{n=0}^\infty n^6 x^n=\frac{x(x+1)(1+56x+246x^2+56x^3+x^4)}{(1-x)^7}.
\ee

\cite[A001015]{sloane} 
\be
\sum_{n=0}^\infty n^7 x^n=\frac{x(1+120x+1191x^2+2416x^3+1191x^4+120x^5+x^6)}{(1-x)^8}.
\ee
See \cite[A008292]{sloane} for coefficients in the numerator
polynomial if exponents of $n$ are higher than 7.
The pole at $x=1$ in the generating functions obviously delimits
the radius of convergence.

\cite{HsuFQ38}
Let
\be
S_{a,p}(n)\equiv \sum_{k=0}^na^kk^p ;
\ee
then
\be
S_{a,p}(n)=\frac{a^{n+1}}{a-1}\sum_{r=0}^{p-1}\binom{p}{r}f_r(a)(n+1)^{p-r}
+f_p(a)\frac{a^{n+1}-1}{a-1}
\ee
with
\be
f_0(a)=1,\quad a\sum_{j=0}^r\binom{r}{j}f_j(a)-f_r(a)=0,\quad r=1,2,\ldots
\ee

\cite{Adamchik}
\be
\sum_{k=0}^\infty\frac{(-1)^k}{(2k+1)^2}=G.
\ee

\cite{Adamchik,BradleyCS2001}
\be
2\sum_{k=0}^\infty\frac{1}{(4k+1)^2}=G+\frac{\pi^2}{8}.
\ee

\cite{Adamchik,BradleyCS2001}
\be
-2\sum_{k=0}^\infty\frac{1}{(4k+3)^2}=G-\frac{\pi^2}{8}.
\ee

\cite{BradleyCS2001}
\be
\sum_{n=0}^\infty \frac{1}{2^{n+1}} \sum_{k=0}^\infty\binom{n}{k}\frac{(-1)^k}{(2k+1)^2}
=G.
\ee

\cite{BradleyCS2001}
\be
1-\sum_{n=1}^\infty \frac{n\zeta(2n+1)}{16^n}
=G.
\ee

\cite{BradleyCS2001}
\be
\frac{1}{4}\sum_{n=1}^\infty n 16^{-n} (3^{2n}-1)\zeta(2n+1)
=G-\frac{1}{6}
\ee
(and similar $\zeta$-sums).

\cite{BradleyCS2001}
\be
\sum_{n=0}^\infty \frac{1}{n+1}\binom{n+1/2}{n}^{-2}
= \,_4F_3(1,1,1,1;2,\frac{3}{2},\frac{3}{2};1)
=2\pi G-\frac{7}{2}\zeta(3)
.
\ee

\cite{BradleyCS2001}
\be
\sum_{n=0}^\infty \frac{4^n}{(2n+1)^2\binom{2n}{n}}
=2 G.
\ee

\cite{BradleyCS2001}
\be
\frac{3}{8}\sum_{n=0}^\infty \frac{1}{(2n+1)^2\binom{2n}{n}}
= G-\frac{1}{8}\pi\log(2+\surd 3).
\ee

\cite{BradleyCS2001}
\be
-\frac{1}{32}\pi\sum_{n=0}^\infty \frac{(2n+1)^2}{(n+1)^3 16^n}\binom{2n}{n}^2
= G-\frac{1}{2}\pi\log 2.
\ee

\cite{BradleyCS2001}
\be
\sum_{n=0}^\infty \frac{\surd 2}{(2n+1)^2 8^n}\binom{2n}{n}
= G+\frac{1}{4}\pi\log 2.
\ee

\cite{BradleyCS2001}
\be
\sum_{n=1}^\infty \frac{\sin(n\pi/4)}{n^2 2^{n/2}}
= G-\frac{1}{8}\pi\log 2.
\ee

\cite{BradleyCS2001}
\be
\sum_{n=0}^\infty \frac{2^{n+1}(n!)^2}{(2n+1)! (n+1)^2}
= 2\pi G-\frac{35}{8}\zeta(3)+\frac{1}{4}\pi^2 \log 2.
\ee

\cite{BradleyCS2001}
\be
\frac{1}{4}\pi \, _3F_2(1/2,1/2,n+1/2;1,n+3/2;1) -\frac{1}{2}\sum_{k=0}^{n-1}\frac{(k!)^2}{(\frac{3}{2})_k^2}
= G.
\ee

\cite{PilehroodAA130,OgreidJCAM98}
\be
\sum_{n=0}^\infty \frac{1}{(n+a)(n+b)}=\frac{\psi(a)-\psi(b)}{a-b}.
\ee

\cite{PilehroodAA130}
\be
\sum_{n=0}^\infty \frac{(-)^n}{(n+a)(n+b)}= -\frac{1}{2}\,\frac{\psi(1/2+a/2)-\psi(a/2)-\psi(1/2+b/2)+\psi(b/2)}{a-b}.
\ee

\cite{CoffeyJCAM183}
\be
\sum_{l=1}^\infty \frac{1}{l\prod_{j=0}^N(l+j)}
=\zeta(2)+\sum_{k=1}^N\frac{(-)^k}{k}\binom{N}{k}
[\psi(k+1)+\gamma].
\ee
\cite{CoffeyJCAM183}
\be
\sum_{l=1}^\infty \frac{(-)^l}{l\prod_{j=0}^N(l+j)}
=-\frac{1}{2}\zeta(2)+\sum_{k=1}^N\frac{(-)^k}{k}\binom{N}{k}
[\psi(1+k/2)-\psi(1+k)].
\ee

\cite{SebahGourdonLn2}
\be
\log 2 = \frac{1}{2}-\sum_{k\ge 1}\frac{(-1)^{k}}{k(4k^4+1)}
.
\ee

After inserting $x=1$ in (\ref{1.511.+1}), \cite[4.1.13]{Apelblat2}
\begin{multline}
2\ln 2
=
1
+
\sum_{k=1}^\infty (-1)^{k+1}\frac{1}{k(k+1)}
=
1
+
2 \sum_{k=1,3,5,\ldots }^\infty \frac{1}{k(k+1)(k+2)}
\\
=
1
+
4 \sum_{k=1,5,9,13\ldots }^\infty \frac{6+4k+k^2}{k(k+1)(k+2)(k+3)(k+4)}
.
\end{multline}
\cite[4.1.20]{Apelblat2}
\begin{multline*}
2\ln 2
=
\frac{5}{4}
+
\sum_{k=1}^\infty (-1)^{k+1}\frac{1}{k(k+1)(k+2)}
=
\frac{5}{4}
+
3 \sum_{k=1,3,5,\ldots }^\infty \frac{1}{k(k+1)(k+2)(k+3)}
\\
=
\frac{17}{12}
-
3 \sum_{k=2,4,6,\ldots }^\infty \frac{1}{k(k+1)(k+2)(k+3)}.
\end{multline*}
\be
\frac{4}{3}\ln 2
=
\frac{8}{9}
+
\sum_{k=1}^\infty (-1)^{k+1}\frac{1}{k(k+1)(k+2)(k+3)}
=
\frac{8}{9}
+
4 \sum_{k=1,3,5,\ldots}^\infty \frac{1}{k(k+1)(k+2)(k+3)(k+4)}.
\ee
More irregular denominators follow from hybridization. For
example we can multiply the penultimate  formula by
$4$, the previous formula by $3$, and subtract
\be
4\ln 2
=
\frac{7}{3}
+12
\sum_{k=1,3,5,\ldots}^\infty
\frac{1}{k(k+1)(k+2)(k+4)}.
\ee

\begin{multline*}
\frac{2}{3}\ln 2
=
\frac{131}{288}
+
\sum_{k=1}^\infty (-1)^{k+1}\frac{1}{k(k+1)(k+2)(k+3)(k+4)}
\\
=
\frac{131}{288}
+
5 \sum_{k=1,3,5,\ldots}^\infty \frac{1}{k(k+1)(k+2)(k+3)(k+4)(k+5)}
.
\end{multline*}

Inserting $x=-1/2$ in (\ref{1.511.+1})
\cite{SebahGourdonLn2,BaileyMaCom66}
yields
\be
\ln 2
=
1
-
\sum_{k=1}^\infty \frac{1}{2^k k(k+1)}
.
\ee
\be
\ln 2
=
\frac{1}{2}
+
\sum_{k=1}^\infty \frac{1}{2^{k-1} k(k+1)(k+2)}
.
\ee
\be
\ln 2
=
\frac{5}{6}
-3
\sum_{k=1}^\infty \frac{1}{2^{k-1}k(k+1)(k+2)(k+3)}
.
\ee
More irregular denominators follow from hybridization. For
example we can multiply the penultimate  formula by
$3$ and add to the previous formula,
\be
4\ln 2
=
\frac{7}{3}
+3
\sum_{k=1}^\infty
\frac{1}{2^{k-1}k(k+1)(k+3)}.
\ee

\be
\ln 2
=
\frac{7}{12}
+3
\sum_{k=1}^\infty \frac{1}{2^{k-3}k(k+1)(k+2)(k+3)(k+4)}
.
\ee
\cite{SebahGourdonLn2}
\be
\log 2 = \frac{1327}{1920}+\frac{45}{4}
\sum_{k\ge 4}\frac{(-1)^k}{k(k^2-1)(k^2-4)(k^2-9)}.
\ee

\cite{SebahGourdonLn2}
\be
\log 2 = \sum_{k\ge 1}\left( \frac{1}{3^k}+\frac{1}{4^k}\right)\frac{1}{k}.
\ee
\be
\log 2 = \sum_{k\ge 0}\left( \frac{1}{8k+8}+\frac{1}{4k+2}\right)\frac{1}{4^k}.
\ee
\cite{SebahGourdonLn2}
\be
\log 2 = \frac{2}{3}+
\frac{1}{2}\sum_{k\ge 1}\left( \frac{1}{2k}+\frac{1}{4k+1}+\frac{1}{8k+4}
+\frac{1}{16k+12} \right)\frac{1}{16^k},
\ee
with a factor $1/2$ corrected on the right hand side.

\cite{SebahGourdonLn2}
\be
\log 2 = \frac{2}{3}
\sum_{k\ge 0}\frac{1}{(2k+1)9^k}.
\ee
\be
\log 2 = \frac{3}{4}
\sum_{k\ge 0}\frac{(-1)^k k!^2}{2^k(2k+1)!}.
\ee
\cite{SebahGourdonLn2,LupasRoGer2000}
\be
\log 2 = \frac{3}{4}
+\frac{1}{4}
\sum_{k\ge 1}\frac{(-1)^k (5k+1)}{k(2k+1)16^k}\binom{2k}{k}.
\ee

\cite[A154920]{sloane}\cite{BaileyBBP}
\be
\log 3 =\sum_{k\ge 1}\left[\frac{9}{2k-1}+\frac{1}{2k}\right]\frac{1}{9^k}.
\ee

\begin{table}
\caption{
Formulas of the type
$
s\log p = t\log q+
\sum_{k=1}^\infty \frac{(-1)^{k+1}}{k}r^k.
$
\cite{MatharArxiv0908}
}
\begin{tabular}{rrrrl}
$s$ & $p$ & $t$ & $q$ & $r$ \\
\hline
1 & 3 & 1 & 2 & $ 1 / 2 \approx 0.50000000$\\
1 & 3 & 2 & 2 & $ -1 / 4 \approx -0.25000000$\\
2 & 3 & 3 & 2 & $ 1 / 8 \approx 0.12500000$\\
5 & 3 & 8 & 2 & $ -13 / 256 \approx -0.05078125$\\
12 & 3 & 19 & 2 & $ 7153 / 524288 \approx 0.01364326$\\
\hline
1 & 5 & 2 & 2 & $ 1 / 4 \approx 0.25000000$\\
3 & 5 & 7 & 2 & $ -3 / 128 \approx -0.02343750$\\
\hline
1 & 7 & 2 & 2 & $ 3 / 4 \approx 0.75000000$\\
1 & 7 & 3 & 2 & $ -1 / 8 \approx -0.12500000$\\
5 & 7 & 14 & 2 & $ 423 / 16384 \approx 0.02581787$\\
\hline
1 & 11 & 3 & 2 & $ 3 / 8 \approx 0.37500000$\\
2 & 11 & 7 & 2 & $ -7 / 128 \approx -0.05468750$\\
11 & 11 & 38 & 2 & $ 10433763667 / 274877906944 \approx 0.03795781$\\
\hline
1 & 13 & 3 & 2 & $ 5 / 8 \approx 0.62500000$\\
1 & 13 & 4 & 2 & $ -3 / 16 \approx -0.18750000$\\
3 & 13 & 11 & 2 & $ 149 / 2048 \approx 0.07275391$\\
7 & 13 & 26 & 2 & $ -4360347 / 67108864 \approx -0.06497423$\\
10 & 13 & 37 & 2 & $ 419538377 / 137438953472 \approx 0.00305254$\\
\hline
1 & 17 & 4 & 2 & $ 1 / 16 \approx 0.06250000$\\
\hline
1 & 19 & 4 & 2 & $ 3 / 16 \approx 0.18750000$\\
4 & 19 & 17 & 2 & $ -751 / 131072 \approx -0.00572968$\\
\hline
1 & 23 & 4 & 2 & $ 7 / 16 \approx 0.43750000$\\
1 & 23 & 5 & 2 & $ -9 / 32 \approx -0.28125000$\\
2 & 23 & 9 & 2 & $ 17 / 512 \approx 0.03320312$\\
\hline
1 & 29 & 4 & 2 & $ 13 / 16 \approx 0.81250000$\\
1 & 29 & 5 & 2 & $ -3 / 32 \approx -0.09375000$\\
7 & 29 & 34 & 2 & $ 70007125 / 17179869184 \approx 0.00407495$\\
\hline
1 & 31 & 4 & 2 & $ 15 / 16 \approx 0.93750000$\\
1 & 31 & 5 & 2 & $ -1 / 32 \approx -0.03125000$\\
\hline
1 & 37 & 5 & 2 & $ 5 / 32 \approx 0.15625000$\\
4 & 37 & 21 & 2 & $ -222991 / 2097152 \approx -0.10633039$\\
5 & 37 & 26 & 2 & $ 2235093 / 67108864 \approx 0.03330548$\\
\hline
1 & 41 & 5 & 2 & $ 9 / 32 \approx 0.28125000$\\
2 & 41 & 11 & 2 & $ -367 / 2048 \approx -0.17919922$\\
3 & 41 & 16 & 2 & $ 3385 / 65536 \approx 0.05165100$\\
\hline
1 & 43 & 5 & 2 & $ 11 / 32 \approx 0.34375000$\\
2 & 43 & 11 & 2 & $ -199 / 2048 \approx -0.09716797$\\
5 & 43 & 27 & 2 & $ 12790715 / 134217728 \approx 0.09529825$\\
7 & 43 & 38 & 2 & $ -3059295837 / 274877906944 \approx -0.01112965$\\
\hline
1 & 47 & 5 & 2 & $ 15 / 32 \approx 0.46875000$\\
1 & 47 & 6 & 2 & $ -17 / 64 \approx -0.26562500$\\
2 & 47 & 11 & 2 & $ 161 / 2048 \approx 0.07861328$\\
\hline
1 & 53 & 5 & 2 & $ 21 / 32 \approx 0.65625000$\\
1 & 53 & 6 & 2 & $ -11 / 64 \approx -0.17187500$\\
3 & 53 & 17 & 2 & $ 17805 / 131072 \approx 0.13584137$\\
4 & 53 & 23 & 2 & $ -498127 / 8388608 \approx -0.05938137$\\
\hline
1 & 59 & 5 & 2 & $ 27 / 32 \approx 0.84375000$\\
1 & 59 & 6 & 2 & $ -5 / 64 \approx -0.07812500$\\
\hline
1 & 61 & 5 & 2 & $ 29 / 32 \approx 0.90625000$\\
1 & 61 & 6 & 2 & $ -3 / 64 \approx -0.04687500$\\
\hline
1 & 67 & 6 & 2 & $ 3 / 64 \approx 0.04687500$\\
\hline
\end{tabular}
\end{table}

\cite[A164985]{sloane}\cite{BaileyBBP}
\be
27\log 5 = 4\sum_{k\ge 0}\left[\frac{9}{4k+1}+\frac{3}{4k+2}+\frac{1}{4k+3}\right]\frac{1}{81^k}.
\ee
\be
\log 5 = 2\log 3-\log 2+\sum_{k=1}^\infty \frac{1}{k10^k}.
\ee
\be
4\log 7 = 5\log 2+\log 3+2\log 5+\sum_{k=1}^\infty \frac{(-1)^k}{k2400^k}.
\ee
\cite{BaileyBBP}
\be
3^5\log 7 = \sum_{k\ge 0}\left[\frac{405}{6k+1}+\frac{81}{6k+2}+\frac{72}{6k+3}
+\frac{9}{6k+4}+\frac{5}{6k+5} \right]\frac{1}{3^{6k}}.
\ee

\cite{BaileyBBP}
\begin{multline}
2\times 3^9\log 11 = \sum_{k\ge 0}\big[\frac{85293}{10k+1}+\frac{10935}{10k+2}
+\frac{9477}{10k+3} +\frac{1215}{10k+4}+\frac{648}{10k+5}+\frac{135}{10k+6}
\\
+\frac{117}{10k+7}
+\frac{15}{10k+8}+\frac{13}{10k+9}\big]\frac{1}{3^{10k}}.
\end{multline}

\be
\log 11 = \log 2+\log 5+\sum_{k=1}^\infty \frac{(-1)^{k+1}}{k10^k}.
\ee
\be
\log 11 = 2\log 2+\log 3-\sum_{k=1}^\infty \frac{1}{k12^k}.
\ee
\be
\log 11 = 2(\log 2+\log 5-\log 3)+\sum_{k=1}^\infty \frac{(-)^k}{k100^k}.
\ee
\be
\log 13 = 2\log 2+\log 3+\sum_{k=1}^\infty \frac{(-1)^{k+1}}{k12^k}.
\ee
\be
3\log 17 = 3\log 2+4\log 5-\sum_{k=1}^\infty \frac{1}{k}\left(\frac{87}{5000}\right)^k.
\ee
\be
\log 19 = 2\log 2+\log 5-\sum_{k=1}^\infty \frac{1}{k20^k}.
\ee
\cite{LafontBBP}
\be
\log 11 = \sum_{k\ge 1}\frac{1}{k}(\frac{3}{2^k}+\frac{1}{2^{2k}}+\frac{1}{2^{5k}}
-\frac{1}{2^{10k}})
=\frac{1}{2}
\sum_{k\ge 1}\frac{1}{k}(\frac{13}{3^k}-\frac{4}{3^{2k}}-\frac{1}{3^{5k}}).
\ee
\cite{LafontBBP}
\be
\log 13 = \sum_{k\ge 1}\frac{1}{k}(\frac{3}{2^k}+\frac{1}{2^{2k}}+\frac{1}{2^{3k}}
+\frac{1}{2^{4k}}-\frac{1}{2^{12k}})
=
\sum_{k\ge 1}\frac{1}{k}(\frac{7}{3^k}-\frac{2}{3^{2k}}-\frac{1}{3^{3k}}).
\ee
\cite{LafontBBP}
\be
\log 17 = \sum_{k\ge 1}\frac{1}{k}(\frac{4}{2^k}+\frac{1}{2^{4k}}-\frac{1}{2^{8k}})
.
\ee
\cite{LafontBBP}
\be
\log 19 = \sum_{k\ge 1}\frac{1}{k}(\frac{3}{2^k}+\frac{3}{2^{2k}}+\frac{1}{2^{9k}}-\frac{1}{2^{18k}})
.
\ee
Similar formulas as the four above are obtained by inserting
\be
\log[(2^s-1)^\tau] = s\tau\log 2-\sum_{k\ge 1}\frac{\tau}{k2^{ks}}
= -s\tau \log \frac{1}{2}-\sum_{k\ge 1}\frac{\tau}{k2^{ks}}
= \sum_{k\ge 1}\frac{\tau}{k}\left(\frac{s}{2^k}-\frac{1}{2^{ks}}\right),
\ee
---immediate consequence of putting $x=1/2$ in (\ref{1.511})---into the right
hand sides of \cite{ChamberlJIS6}:
\be
\log 41=\log(2^{20}-1)
-\log(2^{10}-1)
+\log[(2^2-1)^2]
-\log[(2^4-1)^2]
.
\ee
\be
\log 43=\log(2^{14}-1)-\log(2^2-1)-\log(2^7-1).
\ee
\be
\log 73=\log(2^9-1)-\log(2^3-1).
\ee
\be
\log 151=\log(2^{15}-1)-\log(2^3-1)-\log(2^5-1).
\ee
\be
\log 241=\log(2^{24}-1)-\log(2^{12}-1)-\log(2^8-1)+\log(2^4-1).
\ee
\be
\log 257=\log(2^{16}-1)-\log(2^8-1).
\ee
\be
\log 331=\log(2^{30}-1)-\log(2^{15}-1)-\log(2^{10}-1)+\log(2^5-1)-\log(2^2-1).
\ee
\be
\log 337=\log(2^{21}-1)-\log(2^7-1)-\log[(2^6-1)^2]+\log[(2^2-1)^4].
\ee
\be
\log 683=\log(2^{22}-1)-\log(2^{11}-1)-\log(2^2-1).
\ee
\be
\log 2731=\log(2^{26}-1)-\log(2^{13}-1)-\log(2^2-1).
\ee
\be
\log 5419=\log(2^{42}-1)-\log(2^{21}-1)-\log(2^{14}-1)+\log(2^7-1)-\log(2^2-1).
\ee
\be
\log 43691=\log(2^{34}-1)-\log(2^{17}-1)-\log(2^2-1).
\ee
\be
\log 61681=\log(2^{40}-1)-\log(2^{20}-1)-\log(2^8-1)+\log(2^4-1).
\ee
\be
\log 174763=\log(2^{38}-1)-\log(2^{18}-1)-\log(2^2-1).
\ee
\be
\log 262657=\log(2^{27}-1)-\log(2^9-1).
\ee
\be
\log 599479=\log(2^{33}-1)-\log(2^{11}-1)-\log(2^3-1).
\ee

\be
\sum_{k\ge 0}\frac{\log(2k+1)}{(2k+1)^l} = -\frac{\log 2}{2^l}\zeta(l)-(1-2^{-l})\zeta'(l).
\ee
\be
\sum_{k\ge 0}\frac{\log(ak+b)}{(ak+b)^l} = \frac{\log a}{a^l}\zeta(l,b/a)-\frac{1}{a^l}\zeta'(l,b/a).
\ee

\cite{BorosSci12}
\be
\sum_{k=1}^\infty \frac{k k!}{(2k)!}
=
\frac{1}{8}\left(2+3 e^{1/4}\sqrt{\pi} \erf(1/2)\right)
.
\ee
\be
\sum_{k=1}^\infty \frac{(-1)^k \alpha^k x^k}{k k!}
=
-\gamma -\Gamma(0,\alpha x)-\ln \alpha -\ln x
.
\ee

\cite{Adamchik}
\begin{multline*}
\frac{3}{2}\sum_{k=0}^\infty
\frac{1}{16^k}\Big( \frac{1}{(8k+1)^2}
-\frac{1}{(8k+2)^2}
+\frac{1}{2(8k+3)^2}
-\frac{1}{4(8k+5)^2}
+\frac{1}{4(8k+6)^2}
-\frac{1}{8(8k+7)^2}\Big)
\\
-
\frac{1}{4}\sum_{k=0}^\infty
\frac{1}{4096^k}\Big( \frac{1}{(8k+1)^2}
+\frac{1}{2(8k+2)^2}
+\frac{1}{8(8k+3)^2}
-\frac{1}{64(8k+5)^2}
-\frac{1}{128(8k+6)^2}
-\frac{1}{512(8k+7)^2}\Big)=G,
\end{multline*}
which is Catalan's constant \cite[A006752]{sloane}.

\cite{BroadhurstArxiv98}
\be
\sum_{k\ge 0}\frac{1}{16^k}\left(
\frac{4}{8k+1}-\frac{2}{8k+4}-\frac{1}{8k+5}-\frac{1}{8k+6}
\right)=\pi.
\ee

\cite{GoureAME7}
\be
\frac{1}{16807}
\sum_{n=0}^\infty \frac{1}{2^n\binom{7n}{2n}}
\left(
\frac{59296}{7n+1}
-\frac{10326}{7n+2}
-\frac{3200}{7n+3}
-\frac{1352}{7n+4}
-\frac{792}{7n+5}
+\frac{552}{7n+6}
\right)
=\pi
.
\ee

\cite{Adamchik}
\be
\frac{1}{2}\sum_{k=0}^\infty\frac{4^k k!^2}{(2k)!(2k+1)^2}=G.
\ee

\cite{Adamchik}
\be
-\frac{\pi}{32}\sum_{k=0}^\infty\frac{(2k+1)!^2}{16^k k!^4(k+1)^3}
=G-\frac{\pi}{2}\log 2.
\ee

\cite{Adamchik}
\be
\sqrt{2}\sum_{k=0}^\infty\frac{(2k)!}{8^k k!^2(2k+1)^2}=G+\frac{\pi}{4}\log 2.
\ee

\cite{Adamchik,BradleyRJ3}
\be
\frac{3}{8}\sum_{k=0}^\infty\frac{k!^2}{(2k)!(2k+1)^2}=G+\frac{\pi}{8}\log( 2+\sqrt{3}).
\ee

\cite{BradleyRJ3}
\be
\frac{5}{8}\sum_{k=0}^\infty\frac{L(2k+1)}{(2k+1)^2{2k\choose k}}
=G-\frac{\pi}{8}\log( \frac{10+\sqrt{50-22\surd 5}}{10-\sqrt{50-22\surd 5}}).
\ee
where $L$ are the Lucas numbers.

\cite{OgreidJCAM98}
\be
\sum_{n=1}^\infty \frac{1}{n^2}\frac{[\Gamma(n)]^2}{\Gamma(2n)}
=
-\frac{8}{3}\zeta(3)+\frac{4\pi}{3}\Cl_2(\pi/3).
\ee

\cite{OgreidJCAM98}
\be
\sum_{n=1}^\infty \frac{(-)^n}{n^2}\frac{[\Gamma(n)]^2}{\Gamma(2n)}
=
-\frac{4}{5}\zeta(3).
\ee

\cite{OgreidJCAM98}
\be
\sum_{n=1}^\infty \frac{\Gamma(n+k)}{\Gamma(1+n+2k)}
=
\frac{\Gamma(k)}{\Gamma(1+2k)}.
\ee

\cite{OgreidJCAM98}
\be
\sum_{n=1}^\infty \frac{1}{n^2}[\gamma+\psi(1+kn)]
=
(\frac{k^2}{2}+\frac{3}{2k})\zeta(3)+\pi\sum_{j=1}^{k-1}j \mathrm{Cl}_2(\frac{2\pi j}{k}).
\ee

\cite{OgreidJCAM98}
\be
\sum_{n=1}^\infty \frac{(-)^n}{n^2}[\gamma+\psi(1+kn)]
=
(\frac{k^2}{2}-\frac{9}{8k})\zeta(3)+\pi\sum_{j=1}^{k-1}j \mathrm{Cl}_2(\frac{2\pi j}{k}+\frac{\pi}{k}).
\ee
In the two previous formulas, $k=1,2,3,\ldots$ and sums over $j$ are
understood to be zero when $k=1$.

\cite{OgreidJCAM98}
\be
\sum_{n=1}^\infty \frac{1}{n^2}[\gamma+\psi(kn)]
=
(\frac{k^2}{2}+\frac{1}{2k})\zeta(3)+\pi\sum_{j=1}^{k-1}j \mathrm{Cl}_2(\frac{2\pi j}{k}).
\ee

\cite{OgreidJCAM98}
\be
\sum_{n=1}^\infty \frac{(-)^n}{n^2}[\gamma+\psi(kn)]
=
(\frac{k^2}{2}-\frac{3}{8k})\zeta(3)+\pi\sum_{j=1}^{k-1}j \mathrm{Cl}_2(\frac{2\pi j}{k}+\frac{\pi}{k}).
\ee

\cite{OgreidJCAM98}
\be
\sum_{n=1}^\infty \frac{(-)^n}{n}[\gamma+\psi(1+kn)]
=
-\frac{1}{4}(k+1/k)\zeta(2)+\frac{1}{2}\sum_{j=1}^k\left\{\Cl_1(\frac{2\pi j}{k}+\frac{\pi}{k})\right\}^2.
\ee

\cite{OgreidJCAM98}
\be
\sum_{n=1}^\infty \frac{(-)^n}{n}[\gamma+\psi(kn)]
=
\frac{1}{4}(\frac{1}{k}-k)\zeta(2)+\frac{1}{2}\sum_{j=1}^k\left\{\Cl_1(\frac{2\pi j}{k}+\frac{\pi}{k})\right\}^2,
\ee
plus special cases using $\Cl_1(\theta)=-\log|2\sin\frac{\theta}{2}|$.

\cite{CoffeyJCAM183}
\be
\sum_{n=1}^\infty \frac{1}{n^2}\sum_{r=0}^{q-1}[\gamma+\psi(n+r/q)]
=
q\sum_{n=1}^\infty \frac{1}{n^2}[\gamma+\psi(qn)]-q\ln q\zeta(2).
\ee

\cite{CoffeyJCAM183}
\be
\sum_{n=1}^\infty \frac{(-)^n}{n^2}\sum_{r=0}^{q-1}[\gamma+\psi(n+r/q)]
=
q\sum_{n=1}^\infty \frac{(-)^n}{n^2}[\gamma+\psi(qn)]+\frac{q}{2}\ln q\zeta(2).
\ee

\cite{CoffeyJCAM183}
\be
\sum_{n=1}^\infty \frac{1}{n^2}[\gamma+\psi(n+1/2)]
=
\frac{7}{2}\zeta(3)-2\ln 2\zeta(2).
\ee

\cite{MatharVixra2610}
\be
\sum_{n=0}^\infty \frac{1}{(n+1/2)^2}[\psi(n+3/2)+\gamma] 
=
\frac{7}{2}\zeta(3).
\ee

\cite{MatharVixra2610}
\be
\sum_{n=1}^\infty \frac{1}{n^3}\psi(n+1/2) 
=
-\frac{53}{360}\pi^4+16\Li_4(1/2)-\frac12\log^2 2[\pi^2-\log ^2 2]+\zeta(3)(12\log 2-\gamma).
\ee

\cite{CoffeyJCAM183}
\be
\sum_{n=1}^\infty \frac{1}{n^2}[\gamma+\psi(2n+1/2)]
=
14\zeta(3)-4\pi G-2\ln2 \zeta(2).
\ee

\cite{CoffeyJCAM183}
\be
\sum_{n=1}^\infty \frac{(-)^n}{n^2}[\gamma+\psi(2n+1/2)]
=
14\zeta(3)+\zeta)2\ln 2+2\pi G-8\pi\Cl_2(\pi/4).
\ee

\cite{CoffeyJCAM183}
\be
\sum_{n=1}^\infty \frac{1}{n^2}[\gamma+\psi(kn+1/2)]
=
\frac{7}{2}k^2 \zeta(3)+2\pi\sum_{j=1}^{2k-1}j\Cl_2(\pi j/k)
-\pi \sum_{j=1}^{k-1}j\Cl_2(2\pi j/k)-2\zeta(2)\ln 2.
\ee

\cite{CoffeyJCAM183}
\be
\sum_{n=1}^\infty \frac{(-)^n}{n^2}[\gamma+\psi(kn+1/2)]
=
\frac{7}{2}k^2 \zeta(3)+2\pi\sum_{j=1}^{2k-1}j\Cl_2(\frac{\pi j}{k}+\frac{\pi}{2k})
-\pi \sum_{j=1}^{k=1}j\Cl_2(\frac{2\pi j}{k}+\frac{\pi}{k})+\zeta(2)\ln 2.
\ee

\cite{CoffeyJCAM183}
\be
\sum_{n=1}^\infty \frac{1}{n}\psi'(n+1)=\zeta(3).
\ee

\cite{CoffeyJCAM183}
\be
\sum_{n=1}^\infty \frac{(-)^n}{n}\psi'(n+1)=\zeta(3)-\frac{\pi^2}{4}\ln 2.
\ee

\cite{MatharVixra2610}
\be
\sum_{n=0}^\infty \frac{1}{n+1/2}\psi'(n+1/2)
=
\frac{7}{2}\zeta(3)+\pi^2\log 2.
\ee

\cite{MatharVixra2610}
\be
\sum_{n=1}^\infty \frac{1}{n^2}\psi'(n+1/2)
=
-\frac{121}{360}\pi^4+32\Li_4(1/2)+28\zeta(3)\log 2 -\frac43\log^2 2[\pi^2-\log^2 2].
\ee

\cite{MatharVixra2610}
\be
\sum_{n=1}^\infty \frac{1}{n}\psi''(n+1/2)
=
28\log(2)\zeta(3)-\pi^4/4.
\ee

\cite{CoffeyJCAM183}
\be
\sum_{n=1}^\infty (-)^n\psi^{(j)}(n)=(-)^j[\frac{1}{2^{j+1}}-1]j! \zeta(j+1),
\ee
for integer $j$.

\cite{CvijovicJPA41}
\be
\sum_{n\ge 0} [\psi(\lambda+n)-\psi(\lambda)](\lambda)_n\frac{z^n}{n!} = \frac{z}{(1-z)^\lambda}{}_2F_1(1,1;1;z).
\ee

\cite{CvijovicJPA41}
\be
\sum_{n\ge 0} [\psi(\lambda+n)-\psi(\lambda)]z^n =\frac{1}{\lambda} \frac{z}{(1-z)}{}_2F_1(1,\lambda;\lambda+1;z).
\ee

\cite{CvijovicJPA41}
\be
\sum_{n\ge 0} [\psi(\lambda+n)-\psi(\lambda)]\frac{z^n}{n!} = \frac{1}{\lambda} ze^z {}_2F_2(1,1;2,\lambda+1;-z).
\ee

\cite{CvijovicJPA41}
\be
\sum_{n\ge 0} [\psi(\lambda+n)-\psi(\lambda)]\frac{z^n}{(\lambda)_n} = \frac{1}{\lambda^2} ze^z {}_2F_2(\lambda,\lambda;\lambda+1,\lambda+1;-z).
\ee

\cite{Povol0901}
\be
\surd e
=
\frac{16}{31}
\left[
1+\sum_{n=1}^\infty \frac{1+n^3/2+n/2}{2^n n!}
\right]
,
\ee
which is a special case of $\frac{1}{2}P_3(1/2)+\frac{1}{2}P_1(1/2)+P_0(1/2)$ of \cite{Israel0901}
\be
\sum_{n=0}^\infty n^j\frac{t^n}{n!} = P_j(t)\exp(t)
\ee
where
\be
P_0(t)\equiv 1;\quad 
P_j(t)\equiv t\left[P_{j-1}'(t)+P_{j-1}(t)\right], j>0.
\ee
\cite{AdamchikJCAM79}
\be
\sum_{k\ge 1}\left[\begin{array}{c}k\\ p\end{array}\right]
\frac{1}{k!k} = \zeta(p+1).
\ee

\cite{AdamchikJCAM79}
\be
\sum_{k\ge 2}\left[\begin{array}{c}k\\ 2\end{array}\right]
\frac{1}{k!k^q} = \frac{q+1}{2}\zeta(q+2)
-\frac{1}{q}\sum_{k=1}^{q-1}k\zeta(k+1)\zeta(q+1-k)
.
\ee

\cite{AdamchikJCAM79}
\be
\sum_{k\ge 1}\left[\begin{array}{c}k\\ p\end{array}\right]
\frac{z^k}{k!k} = \zeta(p+1)
+\sum_{k=1}^{p}(\frac{(-)^{k-1}}{k!}\Li_{p+1-k}(1-z)\log^k(1-z)
.
\ee

\cite{AdamchikJCAM79}
\be
\sum_{k\ge p}\left[\begin{array}{c}k\\ p\end{array}\right]
\frac{1}{k!k^q} =
\sum_{k\ge q}\left[\begin{array}{c}k\\ q\end{array}\right]
\frac{1}{k!k^p}
=
\frac{(-)^{q-1}}{(q-1)!p!}
\lim_{\beta\to 0}
\lim_{\alpha\to 0}
\frac{d^{q+p-1}}{d\alpha^qd\beta^{q-1}}
\frac{\Gamma(1-\alpha)\Gamma(1+\beta)}{\Gamma(\beta)\Gamma(1-\alpha+\beta)}
,
\ee
which can always be represented in finite terms of zeta functions.

\cite{Boyadzhievarxiv2012}
\be
\sum_{n=k}^\infty \left[\begin{array}{c}n\\k\end{array}\right]
\frac{1}{x(x+1)\cdots(x+n)} = \frac{1}{x^{k+1}}.
\ee

\cite{Boyadzhievarxiv2012}
\be
\sum_{n=k}^\infty \left[\begin{array}{c}n\\k\end{array}\right]
\frac{1}{m(m+1)\cdots(m+n)} = \zeta(k+1).
\ee

\cite{Boyadzhievarxiv2012}
\be
\sum_{p=1}^\infty \frac{x^{p+n}}{p(p+1)\cdots (p+n)}
= P_n(x)+\frac{(-)^{n-1}}{n!}(1-x)^n\ln(1-x)
\ee
where $P_n(x)$ is the polynomial
\be
P_n(x)=\frac{1}{n!}\left( H_n(x-1)^n+\sum_{k=1}^n\binom{n}{k}\frac{(x-1)^{n-k}}{k}\right),
\ee
and $H$ are the harmonic numbers.

\cite{Boyadzhievarxiv2012}
\be
\sum_{p=1}^\infty \frac{(-)^p}{p(p+1)\cdots (p+n)}
= \frac{2^n}{n!}\left(\sum_{k=1}^n \frac{1}{2^kk}-\ln 2\right).
\ee

\cite{Boyadzhievarxiv2012}
\be
\sum_{n=k}^\infty \left[\begin{array}{c}n\\k\end{array}\right]
\frac{\psi'(n)}{n!} = \frac{k+3}{2}\zeta(k+2)-\frac12 \sum_{j=1}^{k-1}\zeta(j+1)\zeta(k-j+1)
\ee
where $\psi$ is the digamma function.

\cite{Boyadzhievarxiv2012}
\be
\sum_{n=0}^\infty \left\{\begin{array}{c}n\\m \end{array}\right\}
x^n = \frac{x^m}{(1-x)(1-2x)\cdots(1-mx)}
\ee

\cite{Boyadzhievarxiv2012}
\be
\frac{1}{p!}\sum_{n=0}^\infty \frac{n!}{z(z+1)\cdots (z+n)}\binom{n}{p}
 = \frac{1}{(z-1)(z-2)\cdots (z-p-1)}
\ee

\cite{MillerFQ}
\be
\sum_{n=1}^\infty \arctan \frac{2\theta^2}{n^2} = (\frac14+m)\pi
-\arctan\frac{\tanh \pi \theta}{\tan \pi\theta}.
\ee
\cite{MillerFQ}
\be
\sum_{n=1}^\infty \arctan \frac{2k^2}{n^2} = (k-\frac14)\pi.
\ee

\cite{MillerFQ}
\be
\sum_{n=1}^\infty \arctan \frac{k^2}{2n^2} = (k-\frac12)\frac{\pi}{2}.
\ee
\cite{MillerFQ}
\be
\sum_{k=1}^\infty \arctan \frac{z^{2n}}{k^{2n}} =
(z\sec \frac{\pi}{4n}-\frac12)\frac{\pi}{2}+\sum_{k=1}^n(-)^k
\arctan\frac{\sin \zeta}{\cos\zeta-\exp\eta}, \quad 0<|z|<\infty; n =1,2,3,\ldots
\ee
where $\zeta\equiv 2\pi z\cos\frac{2k-1}{4n}\pi$, $\eta=2\pi z\sin\frac{2k-1}{4n}\pi$.

\cite{MillerFQ}
\be
\sum_{k=1}^\infty \arctan\frac{x^4}{k^4}
=[\alpha-\arctan\frac{\sin 2\alpha}{\cos 2\alpha-\exp 2\beta}]
-[\beta-\arctan\frac{\sin 2\beta}{\cos 2\beta -\exp 2\alpha}]-\frac{\pi}{4}
\ee
where $\alpha =\pi x\cos \pi/u$ and $\beta = \pi x \sin \pi /8$.

\cite{BorosSSA11}
\be
\sum_{k=1}^\infty \arctan \frac{2}{k^2} = \frac{3\pi}{4}.
\ee

\cite{BorosSSA11}\cite[A091007]{sloane}
\be
\sum_{k=1}^\infty \arctan \frac{1}{k^2}
=
\arctan\frac{\tan(\pi/\surd 2)-\tanh(\pi/\surd 2)}
{\tan(\pi/\surd 2)+\tanh(\pi/\surd 2)}
.
\ee

\cite{BorosSSA11}
\be
\sum_{k=1}^\infty \arctan \frac{a}{a^2k^2+a(a+2b)k+1+ab+b^2}
=
\frac{\pi}{2}-\arctan(a+b).
\ee

\cite{BorosSSA11}
\begin{multline}
\sum_{k=1}^\infty \arctan \frac{a^2k^2+a^2k-1-b^2}
{a^4k^4+2a^3(a+2b)k^3+a^2(2+a^2+6ab+6b^2)k^2+2a(a+2b)(1+ab+b^2)+(1+b^2)(1+[a+b]^2)}
\\
=
\frac{1}{1+(a+b)^2}.
\end{multline}

\cite{BorosSSA11}
\be
\sum_{k=1}^\infty \arctan \frac{2ak+a+b}{a_4k^4+a_3k^3+a_2k^2+a_1k+a_0}
=
\frac{\pi}{2}-\arctan(a+b+c)
.
\ee

\cite{BorosSSA11}
\be
\sum_{k=1}^n \arctan \frac{f(k+1)-f(k-1)}{1+f(k+1)f(k-1)}
=
\arctan f(n+1)-\arctan f(1)+\arctan f(n)-\arctan f(0)
.
\ee

\cite{BorosSSA11}
\be
\sum_{k=1}^\infty \arctan \frac{8k}{k^4-2k^2+5}
=
\pi -\arctan\frac{1}{2}
.
\ee

\cite{BorosSSA11}
\be
\sum_{k=1}^\infty \arctan \frac{4ak}{k^4+a^2+4}
=
\arctan\frac{a}{2}+\arctan a
.
\ee

\cite{BorosSSA11}
\be
\sum_{k=1}^\infty \arctan \frac{2xy}{k^2-x^2+y^2}
=
\arctan\frac{y}{x}-\arctan \frac{\tanh \pi y}{\tan \pi x}
.
\ee

\cite{Adegokearxiv1710}
\be
\sum_{k\ge 1}\arctan \frac{3}{k^2+3k+1}=\frac{\pi}{2}.
\ee

\cite{Adegokearxiv1710}
\be
\sum_{k\ge 1}\arctan \frac{4k}{(2k^2-1)^2}=\frac{\pi}{2}.
\ee

\cite{Adegokearxiv1710}
\be
\sum_{k\ge 1}\arctan \frac{q}{k^2+qk+1}=\sum_{k=1}^q \arctan\frac{1}{k}.
\ee

\cite{Adegokearxiv1710}
\be
\sum_{k\ge 1}\arctan \frac{\alpha m q}{k^2+(2\beta+mq)k+\beta(\beta+mq)+\alpha^2}
\prod_{j=1}^{m-1}\arctan\frac{\alpha}{k+jq+\beta}
=\sum_{k=1}^q \prod_{j=0}^{m-1} \arctan\frac{\alpha}{k+jq+\beta}.
\ee

\cite{Adegokearxiv1710}
\be
\sum_{k\ge 1}\arctan \frac{2q^2}{(k+q-1)^2}
\arctan \frac{q}{k+q-1}=\frac{\pi^2}{8}+\sum_{k=2}^q \arctan\frac{q}{k-1}
\arctan\frac{q}{k+q-1}.
\ee

\cite{SimpsonEJC34}
Let
\be
u(n,k)=u(n-1,k-1)+(a_{n-1}+b_k)u(n-1,k), u(n,0)=\prod_{i=0}^{n-1}(a_i+b_0), u(0,k)=\delta_{0,k}
\ee
where $a_i$ and $b_i$ are sequences of numbers with $b_i\neq b_j$ when $i\neq j$, then
\be
u(n,k) =\sum_{j=0}^k \frac{\prod_{i=0}^{n-1}(b_j+a_i)}{\prod_{i=0,i\neq j}^k(b_j-b_i)}.
\ee

\cite{RamirezJIS16}
Let
\be
f_n(x,s)=xf_{n-1}(x,s)+q^{n-2}sf_{n-2}(x,s), f_0(x,s)=0, f_1(x,s)=1
\ee
be the Carlitz $q$-Fibonacci polnomials. The generating function is
\be
\sum_{n=0}^\infty f_n(q,s)z^n = \sum_{j=0}^\infty \frac{q^{j^2}s^jz^{2j+1}}{(xz;q)_{j+1}}.
\ee
\cite{RamirezJIS16}
Let
\be
D_n(q)=D_{n-1}(q)+q^{n-2}D_{n-2}(q), D_0(q)=0, D_1(q)=1
\ee
then
\be
D_n(q)=\sum_{j=0}^{\lfloor (n-1)/2\rfloor}\left[
\begin{array}{c}n-j-1\\j \end{array}
\right]q^{j^2}
\ee
where the bracket is \eqref{eq.qbin},
and the generating function is
\be
\sum_{n=0}^\infty D_n(q)x^n = \sum_{j=0}^\infty \frac{q^{j^2}x^{2j+1}}{(x;q)_{j+1}}.
\ee

\subsubsection{Numerical double series (including series of Harmonic Numbers)}

\cite{Connonarxiv07II,CoffeyJCAM159,ChenJIS19}\cite[A152649]{sloane}
\be
\sum_{n=1}^\infty \frac{H_n^{(1)}}{n^3}=\frac{1}{2}\zeta^2(2) =\frac{1}{72}\pi^4,
\label{eq.Hnrdef}
\ee
where $H_n^{(r)}\equiv \sum_{k=1}^n\frac{1}{k^r}$.

\cite{SofoJMA2}
\be
\sum_{n=1}^\infty \frac{H_n^{(p)}}{n^p}=\frac{1}{2}(\zeta^2(p)+\zeta(2p)),
\ee
\be
\sum_{n=1}^\infty \frac{H_n^{(p)}}{n^q}=\zeta(p+q)+\sum_{n=1}^\infty \frac{H_n^{(p)}}{(n+1)^q},
\ee
where $H_n^{(r)}\equiv \sum_{k=1}^n\frac{1}{k^r}$.

\cite{SofoJMA2}
\be
\sum_{n=1}^\infty \frac{H_n^{(p)}}{n^2}=
\frac{p^2+3p+4}{4}\zeta(p+2)
+\zeta(2)\zeta(p)-\sum_{j=1}^{p+2}[\binom{2j-2}{p-1}+2j-2]\zeta(2j-1)\zeta(p-2j+3).
\ee

\cite{ValeanJISxx}
\be
\sum_{n\ge 1}(-)^{n-1}\frac{H_n^{(m)}}{n}
=
\frac12\left(
m\zeta(m+1)-2\log(2)(1-2^{1-m})\zeta(m)
-\sum_{k=1}^{m-2}(1-2^{-k})(1-2^{1+k-m})\zeta(k+1)\zeta(m-k)
\right)
\ee

\cite{ValeanJISxx}
\begin{gather}
\sum_{n\ge 1}(-)^{n-1}\frac{H_{2n}^{(m)}}{n}
=
m\zeta(m+1)-2^{-m}(1-2^{-m+1})\log(2)\zeta(m)
\\
-\sum_{k=0}^{m-1}\beta(k+1)\beta(m-k)
-\sum_{k=1}^{m-2}2^{-m-1}(1-2^{-k})(1-2^{1+k-m})\zeta(k+1)\zeta(m-k)
\nonumber
\end{gather}
where $\beta$ is the Dirichlet beta function.

\cite{Koubaarxiv10,SofoMath7}
\be
\sum_{n\ge 1}(-)^{n-1}\frac{H_{kn}}{n}
=
\frac{(k^2+1)\pi^2}{24k}-\frac12 \sum_{j=0}^{k-1}\log^2\left(2\sin\frac{(2j+1)\pi}{2k}\right).
\ee

\cite{ZhaoJIS13}
\be
\sum_{n=1}^\infty H_nt^n=-\frac{\ln(1-t)}{1-t}.
\ee

\cite{ZhaoJIS13}
\be
\sum_{n=1}^\infty H_n^{[r]}t^n=-\frac{\ln(1-t)}{(1-t)^r}.
\ee
where $H_n^{[r]}\equiv \sum_{k=1}^n H_k^{[r-1]}$ for $n,r\ge 1$
with $H_n^{[0]}\equiv 1/n$.

\cite{ZhaoJIS13}
\be
\sum_{n=r+1}^\infty H_{n,r} t^n=(-)^r \frac{[\ln(1-t)]^r}{r!(1-t)}.
\ee
where $H_{n,r}\equiv \sum_{1\le n_1<\cdots <n_r\le n} \frac{1}{n_1n_2\cdots n_r}$ for $n,r,\ge 1$ and $H_{n,0}=1$.

\cite{LiuJDE23}
\be
\sum_{n=1}^\infty \frac{H_n}{(n+1)^2}=\zeta(3).
\ee

\cite{LiuJDE23}
\be
\sum_{n=1}^\infty \frac{H_n}{n^2}=2\zeta(3).
\ee

\cite{LiuJDE23}
\be
\sum_{n=1}^\infty \frac{(1/2)_n}{(n+1)!}H_n=4\ln 2.
\ee

\cite{LiuJDE23}
\be
\sum_{n=1}^\infty \frac{(1/2)_n}{(n+1)!}H_{n+1}=3.
\ee

\cite{LiuJDE23}
\be
\sum_{n=1}^\infty \frac{(1/2)_n}{(n+1)!}O_n=2;\quad O_n\equiv\sum_{k=1}^n\frac{1}{2k-1}.
\ee

\cite{CoffeyJCAM159}
\be
\sum_{n=1}^\infty \frac{H_n^{(1)^2}}{n^4}=\frac{97}{24}\zeta(6)-2\zeta^2(3)
\approx 1.22187994531988 .
\ee

\cite{CoffeyJCAM159}
\be
\sum_{n=1}^\infty \frac{H_n^{(1)^3}}{n^3}=\zeta^2(3)-\frac{1}{3}\zeta(6)
\approx 1.1058264444388.
\ee

\cite{Connonarxiv07II,SofoMath7}
\be
\sum_{n=1}^\infty \frac{H_n^{(1)}}{n^{2p+1}}=\frac{1}{2}\sum_{j=2}^{2p}(-1)^j\zeta(j)\zeta(2p-j+2),
\ee
where $H_n^{(r)}\equiv \sum_{k=1}^n\frac{1}{k^r}$.

\cite{BorweinPAMS123,MezoArxiv0811,SofoAADM,BatirAADM18,SofoMath7}
\be
\sum_{n=1}^\infty \frac{H_n^{(1)}}{n^m}=\frac{m+2}{2}\zeta(m+1)-\frac{1}{2}\sum_{k=1}^{m-2}\zeta(m-k)\zeta(k+1),
\ee
where $H_n^{(r)}\equiv \sum_{k=1}^n\frac{1}{k^r}$.

\cite{SofoMath7}
\begin{multline}
\sum_{n\ge 1} \frac{H_n^{(p)}}{n^q} = (-)^p \sum_{j=1}^{\lfloor q/2\rfloor}
\binom{p+q-2j-1}{p-1}\zeta(p+q-2j)\zeta(2j)
\\
+\frac{1-(-)^p}{2}\zeta(p)\zeta(q)+(-)^p\sum_{j=1}^{\lfloor p/2\rfloor} \binom{p+q-2j-1}{q-1}\zeta(p+q-2j)\zeta(2j)
\\
+\frac{(-)^p}{2}[(-)^p-\binom{p+q-1}{p}-\binom{p+q-1}{q}]\zeta(p+q).
\end{multline}

\cite{SofoMath7}
\begin{multline}
2\sum_{n\ge 1} (-)^{n+1}\frac{H_n^{(p)}}{n^q} = [1-(-)^p] \zeta(p)\eta(q)+\eta(p+q)
\\
+2(-)^p \sum_{i+2k=q} \binom{p+i-1}{p-1}\zeta(p+i)\eta(2k)
-2 \sum_{j+2k=p} \binom{q+j-1}{q-1}(-)^j\eta(q+j)\eta(2k).
\end{multline}
where $\eta(z)\equiv \sum_{n\ge 1} (-)^{n+1}/n^z$ is the Dirichlet eta function.

\cite{SofoAADM}
\be
\sum_{n=1}^\infty \frac{H_n}{n(n+\alpha)}
=\frac{1}{2\alpha}[3\zeta(2)+\psi^2(\alpha)+2\gamma\psi(\alpha)+\gamma^2-\psi'(x)],
\ee
and
\be
\sum_{n=1}^\infty \frac{H_n}{(n+\alpha)^2}
=\gamma\psi'(\alpha)+\psi(\alpha)\psi'(\alpha)-\frac{1}{2}\psi''(\alpha),
\ee
and
\be
\sum_{n=1}^\infty \frac{H_n}{(n+x)^q}
=
\frac{(-)^q}{(q-1)!}
\left[\left(\psi(x)+\gamma)\right)\psi^{(q-1)}(x)
-\frac{1}{2}\psi^{(q)}(x)+
\sum_{m=1}^{q-2}\binom{q-2}{m}\psi^{(m)}\psi^{(q-m-1)}(x)
\right]
.
\ee
The paper also demonstrates a finite expansion of $\sum_{n>=1}H_n/[n^q\binom{an+k}{k}^t]$
in terms of $\zeta$ and $\psi$ functions for $t=1$ and $2$.

\cite{SofoMath7}
\be
\sum_{n\ge 1}(-)^{n+1}\frac{H_n}{n+1} = \frac{\log^2 (2)}{2}.
\ee

\cite{SofoMath7}
\be
\sum_{n\ge 1}(-)^{n+1}\frac{H_n}{n(n+1)} = \frac{\pi^2}{12}-\log^2 (2).
\ee

\cite{SofoMath7}
\be
\sum_{n\ge 1}(-)^{n+1}\frac{H_{pn}}{n(n+1)} = \frac{p}{2}\zeta(2)
+H_{p-1}\log 2 -\sum_{j=0}^{p-1} \log^2[2\sin\frac{(2j+1)\pi}{2p}]
-\frac12 \sum_{j=1}^{p-1}\frac1j (H_{-j/(2p)}-H_{-(p+j)/(2p)}).
\ee

\cite{SofoMath7}
\be
\sum_{n\ge 1}(-)^{n+1}\frac{H_n}{(n+1)^2} = \zeta(3)/8.
\ee

\cite{SofoMath7}
\be
\sum_{n\ge 1}(-)^{n+1}\frac{H_n}{n^2} = \frac58 \zeta(3).
\ee

\cite{SofoMath7}
\be
\sum_{n\ge 1}(-)^{n+1}\frac{H_{2n}}{n^3} = -5\Li_4(1/2)+\frac{13}{192}\pi^4
-\frac{35}{8}\zeta(3)\log 2+\frac{5}{24}\pi^2\log^2 2-\frac{5}{24}\log^4 2.
\ee

\cite{SofoMath7}
\be
\sum_{n\ge 1}(-)^{n+1}\frac{H_{pn}}{n} = \frac{1+p^2}{4p}\zeta(2)
-\frac12 \sum_{j=0}^{p-1}\log^2[2\sin\frac{(2j+1)\pi}{2p}].
\ee

\cite{SofoMath7}
\be
\sum_{n\ge 1}(-)^{n+1}\frac{O_n}{n} = \frac{\pi^2}{16}.
\ee

\cite{SofoMath7}
\be
\sum_{n\ge 1}(-)^{n+1}\frac{O_n}{n+1} = \frac{\pi^2}{16}(8-\pi)-\log 2.
\ee

\cite{SofoMath7}
\be
\sum_{n\ge 1}(-)^{n+1}\frac{O_n}{n(2n+1)} = \frac{\pi^2}{16}+\frac{\pi\log 2}{4}-G.
\ee

\cite{SofoMath7}
\be
\sum_{n\ge 1}(-)^{n+1}\frac{O_n}{2n+1} = \frac{G}{2} -\frac{\pi\log 2}{8}.
\ee

\cite{SofoMath7}
\be
\sum_{n\ge 1}(-)^{n+1}\frac{O_n}{2n-1} = \frac{G}{2} +\frac{\pi\log 2}{8}.
\ee

\cite{SofoMath7}
\be
\sum_{n\ge 1}\frac{H_nO_n}{n(2n-1)} = \frac54 \zeta(3)-2\zeta(2)\log2 +2\zeta(2).
\ee

\cite{SofoMath7}
\be
\sum_{n\ge 1}(-)^{n+1}\frac{nO_n}{4n^2-1} = \frac{G}{4}.
\ee

\cite{SofoMath7}
\be
\sum_{n\ge 1}(-)^{n+1}\frac{O_n}{n^2} = \pi G -\frac74 \zeta(3).
\ee

\cite{SofoMath7}
\be
\sum_{n\ge 1}(-)^{n+1}\frac{O_n}{n^3} = -4\Li_4(1/2)+\frac{151}{2880}\pi^4-\frac72 \zeta(3)\log 2
+\frac{\pi^2}{6}\log^2 2-\frac16\log^4 2.
\ee

\cite{SofoMath7}
\be
\sum_{n\ge 1}(-)^{n+1}\frac{O_n}{(n+1)^2} = \frac74 \zeta(3)-\frac{\pi^2}{12}-\pi (G-1)-2\log 2.
\ee

\cite{SofoMath7}
\be
\sum_{n\ge 1}(-)^{n+1}\frac{O_n^2}{n} = \frac{7}{16}\zeta(3).
\ee

\cite{SofoMath7}
\be
\sum_{n\ge 1}\frac{O_n}{n^{2m}} = \frac{2^{2m+1}-1}{4}\zeta(2m+1)
-\frac12\sum_{j=1}^{m-1}(2^{2m-2j+1}-1)\zeta(2j)\zeta(2m-2j+1).
\ee

\cite{SofoMath7}
\be
\sum_{n\ge 1}(-)^{n+1}\frac{H_n^2}{n} = \frac34 \zeta(3)-\frac{\pi^2\log 2}{12}+\frac13 \log^3(2).
\ee

\cite{LiuJDE23}
\be
\sum_{n\ge 1} \frac{(a)_n(b_n)}{(c)_nn!}H_n(a) = 
\Gamma\left[\begin{array}{c}c,c-a-b\\ c-a,c-b\end{array}\right][\psi(c-a)-\psi(c-a-b)].
\ee

\cite{LiuJDE23}
\be
\sum_{n\ge 1} \frac{(a)_n(b_n)}{(c)_nn!}H_n(c) = 
\Gamma\left[\begin{array}{c}c,c-a-b\\ c-a,c-b\end{array}\right][\psi(c-a)+\psi(c-b)-\psi(c)
-\psi(c-a-b)].
\ee
with special cases
\be
\sum_{n\ge 1} \frac{H_n}{(n+1)(n+2)}=1;\quad 
\sum_{n\ge 1} \frac{H_{n+2}}{(n+1)(n+2)}=5/4.
\ee

\cite{LiuJDE23}
\begin{multline}
\sum_{n\ge 1}\frac{(a)_n(b)_n}{(1+a-b)_nn!}(-)^n[H_n(a)-H_n(1+a-b)]
=
\Gamma\left[\begin{array}{c}1+a-b,1+a/2\\ 1+a,1+a/2-b\end{array}\right]
\\ \times
[\psi(1+a-b)+\frac12 \psi(1+a/2)-\psi(1+a-\frac12 \psi(1+a/2-b)].
\end{multline}

\cite{LiuJDE23}
\be
\sum_{n\ge 1}\frac{(a)_n(b)_n}{(1+a-b)_nn!}(-)^n[H_n(b)-H_n(1+a-b)]
=
\Gamma\left[\begin{array}{c}1+a-b,1+a/2\\ 1+a,1+a/2-b\end{array}\right]
[\psi(1+a/2-b)-\psi(1+a-b)].
\ee

\cite{BatirAADM18}
\begin{multline}
\sum_{k=1}^\infty \frac{H_{k-1}}{k(p+k)^{m+1}}
=\frac12 \{ (\gamma+\psi(p+1))^2+\zeta(2)-\psi'(p+1)\}
\\
+\frac12 \sum_{l=1}^m \frac{(-)^l}{l!p^{m-l+1}} \{ 2\gamma\psi^{(l)}(p+1)-\psi^{(l+1)}(p+1)
+\sum_{j=0}^l\binom{l}{j}\psi^{(j)}(p+1)\psi^{(l-j)}(p+1)\}.
\end{multline}

\cite{SofoJMA2}
\be
\sum_{n\ge 1} \frac{H_n^{(5)}}{(n+1)^2} = 10\zeta(7)-4\zeta(2)\zeta(5)
-2\zeta(4)\zeta(3).
\ee

\cite{Boyadzhievarxiv2012}
\be
\sum_{n=0}^\infty 
\frac{n!H_n}{z(z+1)\cdots(z+n)} = \frac{1}{(z-1)^2}.
\ee

\cite{BaileyExpM3,CoffeyJCAM159}
\be
\sum_{k=1}^{\infty}\frac{H_k^{(1)2}}{k^2}=\frac{17}{4}\zeta(4).
\ee

\cite{BaileyExpM3}
\begin{multline}
\sum_{k=1}^{\infty}\frac{H_k^{(1)2}}{(k+1)^n}=
\frac{1}{3}n(n+1)\zeta(n+2)+\zeta(2)\zeta(n)-\frac{1}{n}\sum_{k=0}^{n-2}\zeta(n-k)\zeta(k+2)
\\
+\frac{1}{3}\sum_{k=2}^{n-2}\zeta(n-k)\sum_{j=1}^{k-1}\zeta(j+1)\zeta(k+1-j)
+\sum_{k=1}^{\infty} \frac{H_k^{(2)}}{(k+1)^n}.
\end{multline}

\cite{BaileyExpM3}
\be
\sum_{k=1}^{\infty}\frac{H_k^{(2)}}{(k+1)^2}=
\frac{1}{2}\zeta^2(2)-\frac{1}{2}\zeta(4)=\frac{1}{120}\pi^4.
\ee

\cite{ZhaoJIS13}
\be
\sum_{n\ge 1}\frac{H_n}{\binom{n+k}{k}}=\frac{k}{(k-1)^2},\quad k>1.
\ee

\cite{BatirAADM18}
\be
\sum_{k\ge 1} \frac{H_k}{(n+k+1)^m\binom{n+k}{k}}
=...
\ee
\be
\sum_{k\ge 1} \frac{H_k}{k(n+k+1)^m\binom{n+k}{k}}
=...
\ee

\cite{RutledgeAMM45,BaileyExpM3}\cite[A214508]{sloane}
\be
A_4\equiv \sum_{k=1}^\infty (-)^{k+1}\frac{1}{(k+1)^2} H_k^{(2)}
=
-4\Li_4(1/2)+\frac{13\pi^4}{288}+\log(2)[-\frac{7}{2}\zeta(3)+\frac{\pi^2}{6}\log 2
-\frac{\log^3 2}{6}].
\label{eq.A4def}
\ee

\cite{ChenJIS19}
\be
\sum_{k=1}^{\infty}\frac{H_k^{(1)^2}}{k^2}=
\frac{17}{360}\pi^4.
\ee

\cite{BaileyExpM3}
\be
\sum_{k=1}^{\infty}\frac{H_k^{(2)}}{(k+1)^4}=
-6\zeta(6)+\frac{8}{3}\zeta(2)\zeta(4)+\zeta^2(3)=\zeta^2(3)-\frac{4}{2835}\pi^6.
\ee

\cite{BaileyExpM3,BatirAADM18}
\be
\sum_{k=1}^{\infty}\frac{H_k^{(1)}}{(k+1)^n}=
\frac{1}{2}n\zeta(n+1)-\frac{1}{2}\sum_{k=1}^{n-2}\zeta(n-k)\zeta(k+1).
\ee

\cite{BaileyExpM3}
\be
\sum_{k=1}^{\infty}\frac{H_k^{(2)}}{(k+1)^{2n-1}}=
-\frac{1}{2}(2n^2+n+1)\zeta(2n+1)+\zeta(2)\zeta(2n-1)+\sum_{k=1}^{n-1}2k\zeta(2k+1)\zeta(2n-2k).
\ee

\cite{BaileyExpM3}
\begin{multline}
\sum_{k=1}^{\infty}\frac{H_k^{(1)2}}{(k+1)^{2n-1}}=
\frac{1}{6}(2n^2-7n-3)\zeta(2n+1)+\zeta(2)\zeta(2n-1)-\frac{1}{2}\sum_{k=1}^{n-2}(2k-1)\zeta(2n-1-2k)\zeta(2k+2)
\\
+\frac{1}{3}\sum_{k=1}^{n-2}\zeta(2k+1)\sum_{j=1}^{n-2-k}\zeta(2j+1)\zeta(2n-1-2k-2j).
\end{multline}

\cite{MezoCEJM}
\be
\sum_{n=1}^\infty \frac{H_n^{(2)}}{n!}=e\left[1+\frac{1}{4}\,
_2F_2\left(\begin{array}{cc}1 & 1 \\ 3 & 3\end{array}\mid -1\right)\right],
\ee
\be
\sum_{n=1}^\infty \frac{H_n^{(4)}}{2^nn!}=\surd e\left[
\frac{H_1^{(3)}}{2^11!}
+\frac{H_2^{(2)}}{2^22!}
+\frac{H_3^{(1)}}{2^33!}
+\frac{3!}{2^4(4!)^2}\,
_2F_2\left(\begin{array}{cc}1 & 1 \\ 5 & 5\end{array}\mid -\frac{1}{2}\right)\right]
,
\ee
where $H_n^{(r)}$ are hyperharmonic numbers (\ref{eq.hyperH}).

\cite{ChenJIS19}
\be
\sum_{n\ge 0} \frac{1}{8^n}\binom{2n}{n}H_n = 2\sqrt{2} \ln\left(\frac{1+\surd 2}{2}\right).
\ee

\cite{ChenJIS19}
\be
\sum_{n\ge 0} \frac{(-)^n}{8^n}\binom{2n}{n}H_n = \frac{2\sqrt{6}}{3} \ln\left(\frac{3+\surd 6}{6}\right).
\ee

\cite{ChenJIS19}
\be
\sum_{n\ge 0} \frac{1}{8^n}\binom{2n}{n}H_{2n} = \frac{\sqrt{2}}{2} \ln\left(\frac{3+2\surd 2}{2}\right).
\ee

\cite{CampbellRMJM49}
\be
\sum_{n\ge 1} (\binom{2n}{n})^2 \frac{H_n}{16^n(2n-1)^2} = \frac{12-16\ln 2}{\pi}.
\ee

\cite{CampbellRMJM49}
\be
\sum_{n\ge 1} (\binom{2n}{n})^2 \frac{H_n}{32^n} = \frac{\Gamma^2(1/4)}{4\surd \pi}
(1-\frac{4\ln 2}{\pi}).
\ee

\cite{CampbellRMJM49}
\be
\sum_{n\ge 1} (\binom{2n}{n})^2 \frac{H_n}{32^n(n+1)} = 
8-\frac{2\Gamma^2(1/4)}{\pi^{3/2}}
-\frac{4\pi^{3/2}+16\sqrt{\pi}\ln 2}{\Gamma^2(1/4)}.
\ee

\cite{WangJIS14}
\be
\sum_{k= 0}^{2n}\frac{(-1)^k}{\binom{2n}{k}}H_k^{(2)} =
\frac{2n+1}{2n+2}H_{n}^{(2)}.
\ee

\cite{SofoJMA2}
\be
\sum_{k=1}^\infty \frac{H_n^{(2)}}{\binom{n+k}{k}}
=\frac{k}{k-1}(\zeta(2)-H_{k-1}^{(2)}).
\ee
\be
\sum_{k=1}^\infty \frac{H_n^{(p)}}{\binom{n+k}{k}}
=\sum_{r=2}^k(-)^r r \binom{k}{r}
[
\sum_{j=1}^{r-1}\frac{(-)^{p+1}H_j}{j^p}
+\sum_{s=2}^p(-)^{p-s}H_{r-1}^{(p-s+1)}\zeta(s)
]
\ee
\be
\sum_{k=1}^\infty \frac{H_n^{(p)}}{n\binom{n+k}{k}}
=\zeta(p+1)+\sum_{r=1}^k (-)^{r+1}\binom{k}{r}
[
\sum_{j=1}^{r-1}\frac{(-)^{p+1}H_j}{j^p}
+\sum_{s=2}^p(-)^{p-s}H_{r-1}^{(p-s+1)}\zeta(s)
]
\ee
where $p\ge 1$ and $k>1$ are positive integers, $H$ as in \eqref{eq.Hnrdef}.

\cite{XuFilo31}
\be
\sum_{n\ge 1} \frac{H_n^2}{\binom{n+k}{k}} = \frac{k}{k-1}\left(\zeta(2)-H_{k-1}^{(2)}+\frac{2}{(k-1)^2}\right),\quad k\ge 2
\ee

\cite{CoppoJNT150}
\be
\sum_{n\ge 1} (-)^{n-1}\frac{O_n}{n} =\pi^2/16.
\ee
\be
\sum_{n\ge 1} (-)^{n-1}\frac{O_n^{(2)}}{n} =\frac74 \zeta(3)-\frac{\pi}{2}G.
\ee
\be
\sum_{n\ge 1} (-)^{n-1}\frac{O_n^{(3)}}{n} =\frac{\pi^4}{64}-G^2.
\ee
\be
\sum_{n\ge 1} \frac{2^n}{\binom{2n}{n}n}O_n =2G.
\ee
\be
\sum_{n\ge 1} \frac{2^n}{\binom{2n}{n}n}(O_n)^3 =12\beta(4).
\ee
\be
\sum_{n\ge 1} \frac{2^{2n-1}}{\binom{2n}{n}n^2}O_n =\frac{7}{2}\zeta(3).
\ee
\be
\sum_{n\ge 1} \frac{2^{n-1}}{\binom{2n}{n}n^2}O_n =\frac{7}{4}\zeta(3)-\frac{\pi}{2}G,
\ee
\be
\sum_{n\ge 1} \frac{2^{2n-1}}{\binom{2n}{n}n^3}O_n =7\zeta(3)\ln 2-\frac{\pi^4}{32}-8G(1),
\ee

\be
\sum_{n\ge 1} \frac{2^{2n-1}}{\binom{2n}{n}n^4} =\frac{\pi^2}{2}(\ln 2)^2
-\frac72 \zeta(3)\ln 2 +\frac{\pi^4}{96}+4G(1),
\ee
were $O_n\equiv \sum_{k=1}^n 1/(2k-1)$, and $G(1)=\sum_{n\ge 1}O_n/(2n)^3$ is Ramanujan's constant \cite[A256576]{sloane}.

\cite{CoppoJNT150}
Let $\tilde H_n\equiv \sum_{k=1}^n (-1)^{k-1}/k$, then
\be
\sum_{n\ge 1}\frac{\tilde H_n}{(2n-1)^2}
=[\frac{\pi^2}{8}+1]\ln 2-\frac{\pi}{2}+2G -\frac74 \zeta(3)+\frac{\pi G}{2}.
\ee
\be
\sum_{n\ge 1}\frac{\tilde H_n}{(2n-1)^3}
=[\frac78 \zeta(3)-1]\ln 2+\frac{\pi}{2}-2G+\frac{\pi^3}{16}+G^2-\frac{\pi^4}{64}.
\ee

\cite{WangJIS14}
\be
\sum_{k\ge 0}\frac{(-1)^k}{\binom{2n}{k}}H_k =
\frac{2n+1}{2n+2}[H_{2n+1}+H_n-H_{n+1}].
\ee

\cite{ChenJIS19}
\be
\sum_{n\ge 1} \binom{2n}{n}H_n x^n = \frac{1}{\sqrt{1-4x}} \ln\left(\frac{1+\sqrt{1-4x}}{2\sqrt{1-4x}}\right).
\ee

\cite{ChenJIS19}
\be
\sum_{n\ge 1} \binom{2n}{n}{H_{2n}-H_n}) x^n = -\frac{1}{\sqrt{1-4x}} \ln\left(\frac{1+\sqrt{1-4x}}{2}\right).
\ee

\cite{ChenJIS19}
\be
\sum_{n\ge 1} \binom{2n}{n}H_{2n} x^n = \frac{1}{\sqrt{1-4x}} \left[ \ln\left(\frac{1+\sqrt{1-4x}}{2}\right)-2\ln \sqrt{1-4x}\right].
\ee

\cite{BradleyCS2001}
\be
2\sum_{n=0}^\infty \frac{(-1)^{n+1}}{n^2}\sum_{k=0}^{n-1}\frac{1}{2k+1}
=2\pi G-\frac{7}{2}\zeta(3)
.
\ee

\cite{BradleyCS2001}
\be
\sum_{n=0}^\infty \frac{2^n}{(2n+1)\binom{2n}{n}}\sum_{k=0}^n\frac{1}{2k+1}
=2 G.
\ee

\cite{BradleyCS2001}
\be
\sum_{n=0}^\infty \frac{(-1)^n}{2n+1}\sum_{k=1}^{n}\frac{1}{k}
= G-\frac{1}{2}\pi\log 2.
\ee

\cite{BradleyCS2001}
\be
-2\sum_{n=0}^\infty \frac{(-1)^n}{2n+1}\sum_{k=0}^{n-1}\frac{1}{2k+1}
= G-\frac{1}{4}\pi\log 2.
\ee

\cite[\S 1.2]{SteinMForms}
\be
\sum_{m,n\in \mathbb{Z}}\frac{1}{(mz+n)^k} = 
2\zeta(k)+2\frac{(2\pi i)^k}{(k-1)!}
\sum_{n\ge 1}\sigma_{k-1}(n)e^{2\pi i z n},
\ee
where $k$ is an even integer $\ge 4$
and where the sum is over all integer pairs such that $mz+n\neq 0$.

\cite{Espinosaarxiv08,EspinosaRJ22}
Let $T$ denote Tornheim double sums
\be
T(a,b,c)\equiv \sum_{r=1}^{\infty} \sum_{s=1}^\infty \frac{1}{r^a s^b(r+s)^c}
\ee
with $a+c>1$, $b+c>1$ and $a+b+c>2$. Then
\be
T(m,k,n)=\sum_{i=1}^m \binom{m+k-i-1}{m-i}T(i,0,N-i)
+
\sum_{i=1}^k \binom{m+k-i-1}{k-i}T(i,0,N-i)
.
\ee

\begin{multline*}
T(m,0,n)=(-1)^m\sum_{j=0}^{\lfloor (n-1)/2\rfloor}\binom{m+n-2j-1}{m-1}
\zeta(2j)\zeta(m+n-2j)
\\
+(-1)^m\sum_{j=0}^{\lfloor m/2\rfloor}\binom{m+n-2j-1}{n-1}
\zeta(2j)\zeta(m+n-2j)-\frac{1}{2}\zeta(m+n),
\end{multline*}
if $N\equiv m+k+n$ is odd.
\be
T(1,0,5)=\frac{3}{4}\zeta(6)-\frac{1}{2}\zeta^2(3).
\ee
\be
T(2,0,4)=\zeta^2(3)-\frac{4}{3}\zeta(6).
\ee
\be
T(3,0,3)=\frac{1}{2}\zeta(3)^2-\frac{1}{2}\zeta(6).
\ee
\be
T(4,0,2)=-\zeta^2(3)+\frac{25}{12}\zeta(6).
\ee
\be
T(0,0,N)=\zeta(N-1)-\zeta(N),\quad N\ge 3;
\quad
T(0,0,4)=\zeta(3)-\frac{\pi^4}{90}.
\ee
\be
T(n,0,n)=\frac{1}{2}\zeta^2(n)-\frac{1}{2}\zeta(2n).
\quad
T(2,0,2)=\frac{\pi^4}{120}.
\ee
\be
T(0,0,8)=\zeta(7)-\zeta(8).
\ee
\be
T(1,0,7)=\frac{5}{4}\zeta(8)-\zeta(3)\zeta(5).
\ee
\be
T(4,0,4)=\frac{1}{12}\zeta(8).
\ee
\be
T(0,0,14)=\zeta(13)-\zeta(14).
\ee
\be
T(1,0,13)=\frac{11}{4}\zeta(14)-\zeta(3)\zeta(11)-\zeta(5)\zeta(9)-\frac{1}{2}\zeta(7)^2.
\ee
\be
T(1,0,2)=\zeta(3).
\ee
\be
T(1,1,1)=2\zeta(3).
\ee
\be
T(1,0,3)=\frac{\pi^4}{360}.
\ee
\be
T(1,1,2)=\frac{\pi^4}{180}.
\ee
\be
T(2,1,1)=\frac{\pi^4}{72}.
\ee
\be
T(2,2,0)=\frac{\pi^4}{36}.
\ee
\be
T(1,0,5)=\frac{3}{4}\zeta(6)-\frac{1}{2}\zeta^2(3).
\ee
\be
T(3,3,0)=\zeta^2(3).
\ee
\be
T(4,1,1)=\frac{7}{6}\zeta(6)-\frac{1}{2}\zeta^2(3).
\ee
\be
T(4,2,0)=\frac{7}{4}\zeta(6).
\ee
\be
T(1,1,4)=\frac{3}{2}\zeta(6)-\zeta^2(3).
\ee

\cite{OgreidJCAM98}
\be
\sum_{n=1}^\infty \sum_{k=1}^\infty \frac{k}{n(1+k)^2(n+k)}
=
\zeta(3).
\ee

\cite{OgreidJCAM98}
\be
\sum_{n=1}^\infty \sum_{k=1}^\infty \frac{1}{k(n+k)(n+2k)}
=
\frac{3}{4}\zeta(3).
\ee

\cite{OgreidJCAM98}
\be
\sum_{n=0}^\infty \sum_{k=1}^\infty \frac{1}{k(n+k)(n+2k)}
=
\frac{5}{4}\zeta(3).
\ee

\cite{OgreidJCAM98}
\be
\sum_{n=1}^\infty \sum_{k=1}^\infty \frac{1}{k!}
\frac{\Gamma(2k)\Gamma(n+k)}{\Gamma(1+n+2k)}
=
\frac{1}{2}\zeta(2).
\ee

\cite{OgreidJCAM98}
\be
\sum_{n=1}^\infty \sum_{k=1}^\infty \frac{1}{k!k}
\frac{\Gamma(2k)\Gamma(n+k)}{\Gamma(1+n+2k)}
=
\frac{1}{2}\zeta(3).
\ee

\cite{OgreidJCAM98}
\be
\sum_{n=1}^\infty \sum_{k=1}^\infty \frac{1}{k!}
\frac{\Gamma(2k)\Gamma(n+k)}{\Gamma(1+n+2k)}[\gamma+\psi(1+2k)]
=
\frac{11}{8}\zeta(3)
\ee
and similar series.

\cite{Boyadzhievarxiv2012}
\be
\sum_{n=k}^\infty \left[\begin{array}{c}n\\k\end{array}\right]
\sum_{m=0}^\infty
\frac{1}{(m+a)(m+a+1)\cdots(m+a+n)} = \zeta(k+1,a).
\ee

\cite{FurduiJIS14}
\be
\sum_{n=1}^\infty \sum_{m=1}^\infty \frac{(-)^{n+m}}{n+m}=\ln 2 -\frac12.
\ee

\cite{Rogersarxiv1303}
\be
\frac14\sum_{(n,m)\neq (0,0)} \frac{(n-im)^4}{(n^2+m^2)^4}
=
\frac14\sum_{(n,m)\neq (0,0)} \frac{1}{(n+im)^4}
=\frac{\Gamma^8(\frac14)}{3840\pi^2}.
\ee

\cite{Rogersarxiv1303}
\be
\frac14\sum_{(n,m)\neq (0,0)} \frac{1}{(m+in)^8}
=\frac{\Gamma^{16}(\frac14)}{2^{10}525 \pi^4}.
\ee

\cite{Rogersarxiv1303}
\be
\sum_{m,n=-\infty}^\infty \frac{(-)^m(2n+1-2mi)}{(2n+1+2im)^3}
=
\frac{\Gamma^8(\frac14)}{256\surd 2 \pi^2}
\ee

\cite{Rogersarxiv1303}
\be
\sum_{(n,m)\neq (0,0)} \frac{(-)^{m+1}(m^2-2n^2)}{(m^2+2n^2)^2}
=\frac{\Gamma^2(\frac18)\Gamma^2(\frac38)}{48\pi}.
\ee

\cite{Rogersarxiv1303}
\be
\sum_{(n,m)\neq (0,0)} \frac{(-)^{m+1}(m^2-4n^2)}{(m^2+4n^2)^2}
=\frac{\Gamma^4(\frac14)}{32\pi}.
\ee

\cite{Rogersarxiv1303}
\be
\sum_{(n,m)\neq (0,0)} \frac{(-)^mn^2}{(m^2+4n^2)^2}
=\frac{\Gamma^4(\frac14)}{256\pi}-\frac{3\pi\log 2}{32}.
\ee

\cite{Rogersarxiv1303}
\be
\sum_{(n,m)\neq (0,0)} \frac{(-)^{m+1}m^2}{(m^2+4n^2)^2}
=\frac{\Gamma^4(\frac14)}{64\pi}+\frac{3\pi\log 2}{8}.
\ee

\cite{Rogersarxiv1303}
\be
\sum_{(n,m)\neq (0,0)} \frac{(-)^{m+n+1}(m^2-3n^2)}{(m^2+3n^2)^2}
=
\frac{\Gamma^6(\frac13)}{2^{14/3}\pi^2}.
\ee

\cite{Rogersarxiv1303}
\be
\sum_{(n,m)\neq (0,0)} \frac{(-)^{m}(2n^2-m^2)}{(m^2+mn+2n^2)^2}
=
\frac{\Gamma^2(\frac17)\Gamma^2(\frac27)\Gamma^2(\frac47)}{56\pi^2}.
\ee

\cite{Rogersarxiv1303}
\be
\sum_{(n,m)\neq (0,0)} \frac{(-)^{m+n}m^2n^2}{(m^2+n^2)^3}
=
\frac{\Gamma^8(\frac14)}{2^93\pi^3}-\frac{\pi\log 2}{8}.
\ee

\cite{Rogersarxiv1303}
\be
\sum_{(n,m)\neq (0,0)} \frac{(-)^{m+n}m^4}{(m^2+n^2)^3}
=
-\frac{\Gamma^8(\frac14)}{2^93\pi^3}-\frac{3\pi\log 2}{8}.
\ee

\cite{Rogersarxiv1303}
\be
\sum_{(n,m)\neq (0,0)} \frac{(-)^{m}m^2n^2}{(m^2+n^2)^3}
=
-\frac{\Gamma^8(\frac14)}{2^{10}3\pi^3}-\frac{\pi\log 2}{16}.
\ee

\subsubsection{Numerical multiseries}

\cite{FurduiJIS14}
\be
\sum_{n_1,n_2,\ldots n_k=1}^\infty \frac{1}{(n_1+n_2+\cdots+n_k)^m}
=\frac{1}{(k-1)!}\sum_{i=1}^k s(k,i)\left(
\zeta(m+1-i)-1-\frac{1}{2^{m+1-i}}-\cdots-\frac{1}{(k-1)^{m+1-i}} \right).
\ee
where $k$ is a fixed positive integer and $m>k$, where $s$ are the stirling numbers of the first kind.

\cite{FurduiJIS14}
\be
\sum_{n_1,n_2,\ldots n_k=1}^\infty \frac{n_i}{(n_1+n_2+\cdots+n_k)^m}
=\frac{1}{k!}\sum_{i=1}^k s(k,i)\left(
\zeta(m-i)-1-\frac{1}{2^{m-i}}-\cdots -\frac{1}{(k-1)^{m-i}}
\right).
\ee
where $k,i$ are fixed positive integers, $1\le i\le k$ and $m-k>1$, where $s$ are the strigling numbers of the first kind.

\cite{FurduiJIS14}
\be
\sum_{n_1,n_2,\ldots n_k=1}^\infty \frac{n_i}{(n_1+n_2+\cdots+n_k)!}
=\frac{e}{k!}.
\ee
where $k\ge 2$ and $i$ are fixed positive integers, $1\le i\le k$

\cite{FurduiJIS14}
\be
\sum_{n=1}^\infty \sum_{m=1}^\infty \frac{(-)^{n+m}}{(\lfloor \sqrt{n+m}\rfloor)^3}
=\frac34 \zeta(3)-\ln 2-\frac{\pi^2}{12}.
\ee
\cite{FurduiJIS14}
\be
\sum_{n=1}^\infty \sum_{m=1}^\infty (-)^{n+m}\frac{H_{n+m}}{n+m} =-\frac{\ln 2}{2}-\frac{\ln ^2 2}{2}+\frac{\pi^2}{12}.
\ee

\subsection{Infinite Products}

\cite{CooperBAMS71}
\be
\prod_{n=1}^\infty \frac{(1-q^{3n})^{10}}{(1-q^n)^3(1-q^{9n})^3} = 1+3\sum_{n=1, 9\nmid n}^\infty 
\frac{nq^n}{1-q^n}.
\ee

\cite{CooperBAMS71}
\be
\int_0^{e^{-2\pi/3}} \prod_{n=1}^\infty \frac{(1-q^{3n})^{10}}{(1-q^n)^6} = \frac{1}{3^{3/2}}.
\ee

\cite{CampbellIJMM16,CampbellIJMM17}
\be
(1-y)\prod_{k=2}^\infty\prod_{1\le j\le k,(j,k)=1}
(1-x^jy^k)^{1/k} = (\frac{1-y}{1-xy})^{1/(1-x)},\quad |y|<1, |xy|<1.
\ee
\cite{CampbellIJMM16}
\be
(1-xy)\prod_{k=2}^\infty\prod_{1\le j\le k,(j,k)=1}
(1-x^jy^k)^{1/j} = (1-xy)^{1/(1-x)},\quad |x|<1, |xy|<1.
\ee
\cite{CampbellIJMM17}
\be
\prod_{a,b<c} (1-x^ay^bz^c)^{-1/c}
=(\frac{(1-xz)(1-yz)}{(1-z)(1-xyz)})^{1/((1-x)(1-y))} 
\ee
for $a,b>0$ and $(a,b,c)=1$, $|z|,|xz|,|yz|,|xyz|<1$.

\cite{CampbellIJMM17}
\be
\prod_{a,c\le b} (1-x^ay^bz^c)^{-1/c}
=(\frac{(1-xyz)^{xy}}{(1-yz)^y)})^{1/((1-x)(1-y))} 
\ee
for $a,b>0$ and $(a,b,c)=1$, $|yz|,|xyz|<1$.

\subsection{Formulae from Differential Calculus}
\cite{KolbigJCAM69}
\be
\left(\frac{d}{dx}\right)^m e^{g(x)}=e^{g(x)}Y_m\left(
g'(x),g''(x),\ldots,g^{(m)}(x)
\right)
\ee
where $Y_m$ is the $m$th exponential complete Bell polynomial,
\be
Y_n(x_1,\ldots,x_n)=\sum_{\pi(n)}\frac{n!}{k_1!\cdots k_n!}
\left(\frac{x_1}{1!}\right)^{k_1}
\cdots
\left(\frac{x_n}{n!}\right)^{k_n}.
\ee

\be
\frac{d}{dx}\left[f(x)^{g(x)}\right]
=
f(x)^{g(x)}\left[
(\log f)g'+g\frac{f'}{f}
\right]
\ee

\cite{DilAADM5}
\be
\sum_{n\ge 0}\frac{g^{(n)}(0)}{n!}f(n)x^n
=\sum_{n\ge 0} \frac{f^{(n)}(0)}{n!}
\sum_{k=0}^n \left\{\begin{array}{c}n\\ k\end{array}\right\} x^k g^{(k)}(x).
\ee
where the braces are Stirling numbers of the second kind.

\cite{RzadkowskiJIS12}
\be
\frac{d^n}{dt^n}(f(e^t)) = \sum_{k=1}^n\left\{\begin{array}{c}n\\k\end{array}\right\} f^{(k)}(e^t)e^{kt}.
\ee

\cite{RzadkowskiJIS12}
\be
\frac{d^n}{dt^n}(f(\log t)) = \frac{1}{t^n}\sum_{k=1}^n (-1)^{n-k}\left[\begin{array}{c}n\\k\end{array}\right] f^{(k)}(\log t).
\ee

\cite{RainaITSF3}
Let 
\be
S_n^{(k)}(r;y)\equiv \frac{\partial ^n}{\partial t^n} \left\{
t^{rk}f(yt)
\right\}_{\mid t=1}
\ee
with $n,k \in N$, $r$ a real or complex number and $f(yt)$ differentiable
at $t=1$,
then
\be
f(y)\sum_{n\ge 0}A_n (-rz)^n
=\sum_{n\ge 0}A_n \frac{x^n}{n!} \sum_{k=0}^n (-)^k\binom{n}{k}
S_n^{(k)}(r;y),
\ee
where $\{A_n\}$ is an arbitrary bounded sequence.

\cite{RainaITSF3}
If $\{\Omega_n\}$ is any bounded squence and the elements $\{a_1,\ldots a_n\}$
are such that $\Re (a_i)>-1$ for $i=1,\ldots,n$ then
\be
\sum_{k=0}^\infty \prod_{j=1}^n \frac{(a_j)_k}{(a_j+1)_k} \Omega_k z^k
=
\sum_{i=1}^n\left(\sum_{k=0}^\infty \frac{(a_i)_k}{(a_k+1)_k}\Omega_k z^k\right)
\frac{1}{a_i}\prod_{\mu=1,\mu\neq i}^n \frac{a_\mu}{a_\mu-a_i}.
\ee

\section{Elementary Functions}

\subsection{Powers of Binomials}

\cite{BerndtBLMS15}
\be
\left(\frac{1+\sqrt{1+4x}}{2}\right)^n
=(1+x)^{n/2}
\sum_{k\ge 0}
\frac{(-)^{k+1}\Gamma(\frac{3k-2}{2})}{\Gamma(\frac{3k-n}{2}-k+1) k!}
\left(\frac{x}{(1+x)^{3/2}}\right)^k
.
\ee

\be
1+\sum_{n\ge 1}\binom{\alpha+n-1}{\alpha-1}x^n = {}_1F_0(\alpha\mid x)
=\frac{1}{(1-x)^\alpha}.
\ee

\cite{LehmerAMM92}
\be
\sum_{n=0}^\infty \binom{2n}{n}x^n=\frac{1}{\sqrt{1-4x}}
.
\ee
\cite{LehmerAMM92}\cite[A000108]{sloane}
\be
\sum_{n=0}^\infty \binom{2n}{n}\frac{x^n}{n+1}
=\frac{1-\sqrt{1-4x}}{2x}
.
\ee

\begin{multline}
\sum_{n\ge 0}\binom{3n}{n}z^n
=
\sum_{n\ge 0}\frac{\Gamma(3n+1)}{\Gamma(2n+1)}\frac{z^n}{n!}
=
\sum_{n\ge 0}\frac{\Gamma(n+1/3)\Gamma(n+2/3)\Gamma(n+1)}
{\Gamma(n+1/2)\Gamma(n+1)}\frac{z^n}{n!}
\frac{(2\pi)^{1/2} 2^{2n+1/2} 3^{3n+1/2}}{2\pi}
\\
=
\frac{\Gamma(1/3)\Gamma(2/3)}{\Gamma(1/2)} \sum_{n\ge 0}\frac{\Gamma(n+1/3)\Gamma(n+2/3)\Gamma(1/2)}
{\Gamma(n+1/2)\Gamma(1/3)\Gamma(2/3)}\frac{z^n}{n!}
\frac{(2\pi)^{1/2} 2^{2n+1/2} 3^{3n+1/2}}{2\pi}
\\
=
\sqrt{6/(2\pi)}
\frac{\Gamma(1/3)\Gamma(2/3)}{\Gamma(1/2)} \sum_{n\ge 0}\frac{\Gamma(n+1/3)\Gamma(n+2/3)\Gamma(1/2)}
{\Gamma(n+1/2)\Gamma(1/3)\Gamma(2/3)}\frac{(2^2 3^3z)^n}{n!}
\\
=
\sqrt{6/(2\pi)}
\frac{\Gamma(1/3)\Gamma(2/3)}{\Gamma(1/2)} \sum_{n\ge 0}\frac{(1/3)_n (2/3)_n}
{(1/2)_n}\frac{(2^2 3^3z)^n}{n!}
\\
=
\sqrt{6/(2\pi)}
\frac{2\pi}{\sqrt{3\pi}} {} _2F_1(1/3, 2/3; 1/2; 2^23^3z)
\end{multline}
This is eventually reduced by \cite[15.1.18]{AS}
\begin{equation*}
_2F_1(a,1-a;1/2;\sin^2z)=\frac{\cos[(2a-1)z]}{\cos z}.
\end{equation*}

\cite{LehmerAMM92}
\be
\sum_{n=1}^\infty \frac{1}{n}\binom{2n}{n}x^n
=2\log\frac{1-\sqrt{1-4x}}{2x}
.
\ee
\cite{LehmerAMM92}
\be
x\sum_{n=1}^\infty \frac{1}{n(n+1)}\binom{2n}{n}x^n
=2x\log\frac{1-\sqrt{1-4x}}{x}
+\frac{\sqrt{1-4x}}{2}
-x(\log 4 -1)-\frac{1}{2}
,
\ee
which corrects a sign error in \cite{LehmerAMM92}, see \cite{MatharArxiv0905}.

\cite{BoyadJIS15}
\be
\sum_{n=0}^\infty \frac{1}{n+1}\binom{2n}{n}H_n x^{n+1}
=\sqrt{1-4x}\log(2\sqrt{1-4x})-(1+\sqrt{1-4x})\log(1+\sqrt{1-4x})+\log 2,
\ee
\be
\sum_{n=0}^\infty \binom{2n}{n}H_n x^n
=\frac{2}{\sqrt{1-4x}}\log(\frac{1+\sqrt{1-4x}}{2\sqrt{1-4x}}),
\ee
\be
\sum_{n=0}^\infty \binom{2n}{n}H_n (-1)^{n-1} x^n
=\frac{2}{\sqrt{1+4x}}\log(\frac{2\sqrt{1+4x}}{1+\sqrt{1+4x}}),
\ee
where $H_n=\sum_{k=1}^n 1/k$.

\cite{BoyadJIS15}
\be
\sum_{n=0}^\infty \binom{2n}{n}P_q(n) x^{n+1}
=\frac{1}{\sqrt{1-4x}}
\sum_{k=0}^q\binom{2k}{k} k!\left(\frac{x}{1-4x}\right)^k\sum_{m=k}^q a_m S(m,k)
\ee
where
$P_q(z)=a_qz^q+a_{q-1}z^{q-1}+\ldots+q_0$ is a polynomial and $S$ are the
Stirling Numbers of the Second Kind.

\cite[A005430]{sloane}
\be
\sum_{n=1}^\infty n\binom{2n}{n}x^n
=2x(1-4x)^{-3/2}
.
\ee

\cite{LehmerAMM92}\cite[A002736]{sloane}
\be
\sum_{n=1}^\infty n^2\binom{2n}{n}x^n
=\frac{2x(2x+1)}{(1-4x)^{5/2}}
.
\ee

\cite{BoyadJIS15}
\be
\sum_{n=0}^\infty n^m\binom{2n}{n}x^n
=\frac{1}{\sqrt{1-4x}} \sum_{k=0}^m S(m,k)\binom{2k}{k}k!\left(\frac{x}{1-4x}\right)^k
.
\ee

\cite{LehmerAMM92}
\be
\sum_{n=0}^\infty \frac{1}{2n+1}\binom{2n}{n}x^n
=\frac{1}{2x}\arcsin(2x)
.
\ee

\cite{BoyadJIS15}
\be
\sum_{n=0}^\infty \frac{1}{n+m+1}\binom{2n}{n}x^n
=\frac{1}{2^{2m+1}x^{m+1}}\sum_{k=0}^m
\binom{m}{k}\frac{(-)^k}{2k+1}[1-(1-4x)^{k+1/2}]
.
\ee

\cite{BoyadJIS15}
\begin{multline}
\sum_{n=0}^\infty \frac{1}{n^2}\binom{2n}{n}x^n
=2\Li_2\left(\frac{1-\sqrt{1-4x}}{2}\right)
\\
-\log^2(1+\sqrt{1-4x})-2\log 2 \log\frac{1-\sqrt{1-4x}}{x}
+3\log^2 2
.
\end{multline}

\cite{LehmerAMM92}
\be
\sum_{m=1}^\infty \frac{(-1)^{m-1}}{\binom{2m}{m}}(2x/y)^{2m}
=\frac{xy^2}{h^3}\left[ \log\frac{x+h}{y}+\frac{xh}{y^2}
\right]
\ee
where $h\equiv\sqrt{x^2+y^2}$.

\cite{LehmerAMM92}
\be
\sum_{m=1}^\infty \frac{(2x)^{2m}}{m\binom{2m}{m}}
=\frac{2x\arcsin x}{\sqrt{1-x^2}}
.
\ee
\cite{LehmerAMM92}
\be
\sum_{m=1}^\infty \frac{(2x)^{2m}}{m^2\binom{2m}{m}}
=2(\arcsin x)^2
.
\label{eq.2xbin}
\ee
\cite{LehmerAMM92}
\be
\sum_{m=1}^\infty \frac{(2x)^{2m}}{\binom{2m}{m}}
=\frac{x^2}{1-x^2}+\frac{x\arcsin x}{(1-x^2)^{3/2}}
.
\ee

\cite{BatirAMC220}
\begin{multline}
\sum_{k\ge 1} \frac{(256x/27)^k}{k^2\binom{4k}{k}}
=
-\frac{3}{2}\log^2|\frac{\varphi-1}{\varphi+1}|
+\frac{3}{4}\left[
\arctan(-\sqrt{\frac{\varphi}{\varphi+2}})
+\arctan(\frac{\varphi+3}{\varphi-1}\sqrt{\frac{\varphi}{\varphi+2}})
\right]^2
\\
+\frac{3}{4}\left[\arctan(\sqrt{\frac{\varphi}{\varphi-2}})
+\arctan(\frac{3-\varphi}{\varphi+1}\sqrt{\frac{\varphi}{\varphi-2}})\right]^2
\end{multline}
where $\varphi(x)=\sqrt{1+\frac{3}{2}(x^2+\sqrt{x^4-x^3})^{-1/3}
+\frac{3}{2x}(x^2+\sqrt{x^4-x^3})^{1/3}}$,
and two similar expressions if the power of $k$ in the denominator is
not 2 but 1 or 0.

\begin{multline}
\frac{1}{(\sigma^2+\sigma_L^2)^{\gamma/2}}
= 
\frac{1}{\sigma_L^\gamma} \exp\left(-\frac{\gamma}{2}\,\frac{\sigma^2}{\sigma_L^2}\right)
\Big(1
+\frac{\gamma}{4}\left(\frac{\sigma}{\sigma_L}\right)^4
-\frac{\gamma}{6}\left(\frac{\sigma}{\sigma_L}\right)^6
+\frac{\gamma(\gamma+4)}{32}\left(\frac{\sigma}{\sigma_L}\right)^8
\\
-\frac{\gamma(5\gamma+12)}{120}\left(\frac{\sigma}{\sigma_L}\right)^{10}
+\frac{\gamma(3\gamma^2+52\gamma+96)}{1152}\left(\frac{\sigma}{\sigma_L}\right)^{12}
-\frac{\gamma(35\gamma^2+308\gamma+480)}{6720}\left(\frac{\sigma}{\sigma_L}\right)^{14}
+\cdots\Big)
\end{multline}

\cite[3.41(a)]{Rivlin}
\be
\frac{1}{a^2-x^2}=\frac{2}{a\sqrt{a^2-1}}{\sum_{j=0}^\infty}' (a-\sqrt{a^2-1})^{2j}T_{2j}(x)
\ee
where the prime at the sum symbols means taking only half of the value at index zero.

\cite{Blasiakarxiv08}
\be
(1+x+x^{-1})^n=\sum_{j=-n}^n \binom{n}{j}_2 x^j,\quad
\binom{n}{m}_2\equiv
\sum_{j\ge 0}\frac{n!}{j!(m+j)!(n-2j-m)!}
.
\ee

\cite{PilehroodDMTCS10}
\be
\sum_{s=0}^\infty
\zeta(4s+3)x^{4s} =
\sum_{k=1}^\infty \frac{k}{k^4-x^4}
=\frac{5}{2}\sum_{k=1}^\infty \frac{(-)^{k-1}}{\binom{2k}{k}}
\,
\frac{k}{k^4-x^4}\prod_{m=1}^{k-1}\frac{ m^4+4x^4}{m^4-x^4}.
\ee

\cite{BalanzarioRJ19}
Let $B_0(x)$ be a period function of period 1. Assume $B_0(x)$
has a continuous derivative in the open interval (0,1).
Let $a_0\equiv\int_0^1 B_0(x)dx$ and define
\be
B_n(x)=\int_0^2 B_{n-1}(y)dy+\int_0^1(y-1)B_{n-1}(y)dy.
\ee
Then
\be
\sum_{k=0}^\infty B_k(x)t^k=\frac{t e^{xt}}{e^t-1}
\left[a_0-\int_0^1 B'_0(1-y)\frac{e^{ty}-1}{t}dy\right]
+\int_0^x e^{t(x-y)}B_0'(y)dy,
\ee
and
\be
B_n(x)=\Re\left[2\sum_{j=1}^\infty\frac{e^{2\pi ijx}}{(2\pi ij)^n}(a_j-a_0-ib_j) \right]
\ee
where $a_j$ and $b_j$ are the Fourier coefficients of $B_0(x)$.
For example $B_0(x)=\cos(\pi x)$ yields
$B_2(x)=(1-2x-\cos(\pi x))/\pi^2$ and
\be
\frac{\pi^3}{2}B_2(x)=-\sum_{j=1}^\infty
\frac{\sin(2\pi jx)}{j(4j^2-1)}.
\ee
Consider also a periodic sequence $\lambda_{j+T}=\lambda_j$,
the Dirichlet series
\be
f(s)=\sum_{j=1}^\infty \frac{a_j-a_0}{j^s}\lambda_j,\quad
g(s)=\sum_{j=1}^\infty \frac{b_j}{j^s}\lambda_j,\quad
\ee
and the Fourier coefficients
\be
\alpha_j\equiv \sum_{k=1}^T \lambda_k\cos(2\pi j k/T),\quad
\beta_j\equiv \sum_{k=1}^T \lambda_k\sin(2\pi j k/T).
\ee
If $\lambda_{T-j}=\lambda_j$ for $1\le j<T$, then
\be
f(2n) = \frac{(-)^n}{2T}(2\pi)^{2n}\sum_{j=1}^T \alpha_j B_{2n}(j/T),
\ee
\be
g(2n+1) = \frac{(-)^{n+1}}{2T}(2\pi)^{2n+1}\sum_{j=1}^T \alpha_j B_{2n+1}(j/T),
\ee
and similar expressions for other even-odd symmetries of the $\lambda$.

\cite{BalanzarioRJ19}
\be
8\log 2-4=\zeta(3)+\sum_{j=1}^\infty \frac{1}{j^2(4j^2-1)}
=\sum_{k=0}^\infty \frac{\zeta(3+2k)}{4^k}.
\ee

\cite{BalanzarioRJ19}
\be
\sum_{j=1}^\infty \frac{\cos(2\pi j/3)}{j^{2k}(4\pi^2j^2+1)}
=\frac{(-)^k}{4}(2\pi)^{2k} \alpha_{2k},
\ee
with g.f.
\be
\sum_{k=0}^\infty \alpha_k t^k=
\frac{e+e^{4/3}-e^{1+t}-e^{t+4/3}+t^3[e^{(2t+5)/3}-e^{(2+2t)/3}
-e^{(2+t)/3}+e^{(5+t)/3}]}{(e-1)e^{2/3}(e^t-1)(t^2-1)}
-\frac{2te^{t/3}}{e^t-2}.
\ee

\be
\sum_{k=1}^\infty \frac{1}{k^j(k^2- s)^l} 
= 
\sum_{i\ge 0}
\frac{(l)_i}{i!}s^i 
\zeta(j+2l+2i)\quad 2l+j> 1.
\label{eq.invk}
\ee

\subsection{The Exponential Function}

\cite{GandhiAMM77}, \cite[A001469]{sloane}
\be
\frac{2t}{e^t+1}=\sum_{N=0}^\infty \frac{G_N}{N!}t^N;
\ee
with $(G+1)^N+G_N=1$ for $N>1$, $G_1=1$, $G_2=-1$, $G_4=1$, $G_6=-3$,
$G_8=17$, $G_{10}=-155$, $G_{12}=2073$, Genocchi numbers.

\cite{Blasiakarxiv08}
\be
\exp(xt+yt^2)=\sum_{n=0}^\infty\frac{t^n}{n!}H_n(x,y),\quad
H_n(x,y)\equiv n!\sum_{k=0}^{[n/2]}\frac{x^{n-2k}y^k}{k!(n-2k)!}.
\ee

\cite{ZhaoDM281}
\be
[t^n](e^t-1)^m = \frac{m!}{n!}\begin{Bmatrix}n \\ m\end{Bmatrix}
\therefore
\sum_{k=0}^n \begin{bmatrix}{c}n\\k \end{bmatrix}
k![t^k]f(e^t-1)=\sum_{j=1}^n n!\binom{n-1}{j-1}f_j .
\ee

\subsection{Trigonometric and Hyperbolic Functions}
\cite{PilehroodDMTCS12}
\be
\sin\frac{\pi a}{6} = \frac{a}{2} \sum_{n=0}^\infty \frac{\binom{2n}{n}}{16^n(2n+1)}\prod_{m=0}^{n-1} (1-\frac{a^2}{(2m+1)^2}).
\ee

\cite{TaylorPEMS32}
\be
\frac{\cos n\alpha -\cos n\theta}{\cos \alpha -\cos \theta}
= 2^{n-1}\prod_{r=1}^{n-1}[\cos \theta -\cos(\alpha+\frac{2r\pi}{n})]
\ee

\cite{NetoJIS17}
\be
\frac{d^n}{dx^n}\cot x = (2i)^n (\cot(x)-i)\sum_{k=1}^n \frac{k!}{2^k}S(n,k)(i\cot(x) -1)^k
\ee
\be
\frac{d^n}{dx^n}\sec x = \sec(x) i^n \sum_{j=0}^n(-)^j j!\sum_{k=j}^n 
\binom{n}{k}S(k,j)2^{k-j}(i\tan(x) +1)^j
\ee
where $i=\sqrt{-1}$ and $S(n,k)$ are the Stirling numbers of the second kind.

\subsection{Fourier Series}\label{sec.Fser}

\cite[B2a]{DaiJOSAA12}
\be
\cos(m\theta)=\sum_{k=0}^{\lfloor m/2\rfloor} (-1)^k\binom{m}{2k}
\cos^{m-2k}\theta \sin^{2k}\theta.
\ee
\cite{TaylorPEMS32}
\be
2\cos n\theta = (2\cos \theta)^n
-\frac{n}{1!}(2\cos\theta)^{n-2}
+\frac{n(n-3)}{2!}(2\cos\theta)^{n-4}+\cdots
\ee
\cite{TaylorPEMS32}
\be
\frac{\cosh n k -\cos n\theta}{\cosh k - \cos\theta}
=\frac{\sh n k}{\frac \sh k}+\frac{2}{\sh k}\sum_{r=1}^{n-1}\sh(n-r) k \cos r\theta.
\ee

\cite{TaylorPEMS32}
\be
\frac{\cos n \alpha -\cos n \theta}{\cos \alpha - \cos \theta}
=\frac{\sin n \alpha}{\sin \alpha}
+\frac{2}{\sin \alpha}\sum_{r=1}^{n-1}\sin(n-r)\alpha \cos r\theta
=\frac{\sin n \theta}{\sin \theta}
+\frac{2}{\sin \theta}\sum_{r=1}^{n-1}\sin(n-r)\theta \cos r\alpha.
\ee
\cite{TaylorPEMS32}
\be
\frac{\cos n \alpha -\cos n \theta}{\cos \alpha - \cos \theta}
=4\sum_{r,s}\cos r\alpha\cos s\theta+2\sum_t (\cos t \alpha+\cos t \theta)+c,
\ee
with sums over all positive integral nonzero values or $r,s,t$ for which  $n-r-s$, $n-t$ are odd and $c$ zero or unity according as $n$ is even or odd.
\cite{TaylorPEMS32}
\be
\frac{1-\cos n\theta}{1-\cos \theta} = n+2\sum_{r=1}^{n-1} (n-r) \cos r\theta.
\ee
\cite{TaylorPEMS32}
\be
\frac{\cos n \alpha -(-)^n\cos n \theta}{\cos\alpha +\cos\theta}
= \frac{\sin n\alpha}{\sin \alpha}+\frac{2}{\sin\alpha}
\sum_{r=1}^{n-1} (-)^r\sin(n-r)\alpha\cos r\theta.
\ee
\cite{TaylorPEMS32}
\be
\frac{1 -(-)^n\cos n \theta}{1 +\cos\theta}
= n+2
\sum_{r=1}^{n-1} (-)^r(n-r)\cos r\theta.
\ee
The hyperbolic analogues are similar
\cite{TaylorPEMS32}, for example:
\be
\frac{\cosh n k -\cosh n u}{\cosh k - \cosh u}
=\frac{\sinh n k}{\sh k}
+\frac{2}{\sinh k}\sum_{r=1}^{n-1}\sinh(n-r)k \cosh ru.
\ee
\cite{TaylorPEMS32}
\be
\frac{\sin n\alpha}{\sin \alpha}=2\cos(n-1)\alpha +2\cos(n-3)\alpha+\cdots,
\ee
where the final term is unit if $n$ is odd.
\cite{TaylorPEMS32}
\be
\frac{\cos (2n-1)\alpha}{\cos \alpha}=2\cos(2n-2)\alpha -2\cos(2n-4)\alpha
+\cdots+ (-1)^{n-1}
\ee
\cite{TaylorPEMS32}
\be
\frac{(-1)^{n-1}+\cos 2n\alpha}{\cos \alpha}=2\cos(2n-1)\alpha -2\cos(2n-3)\alpha
+\cdots+ (-1)^{n-1} 2\cos\alpha.
\ee
\cite{TaylorPEMS32}
\be
\frac{\sin 2n\alpha}{\cos \alpha}=2\sin(2n-1)\alpha -2\sin(2n-3)\alpha
+\cdots.
\ee
\cite{TaylorPEMS32}
\be
\frac{1-\cos 2n\alpha}{\sin \alpha}=2\sin(2n-1)\alpha +2\sin(2n-3)\alpha
+\cdots.
\ee

\cite{TaylorPEMS32}
\be
\frac{d^m}{d\alpha^m} \frac{\sin n\alpha}{\sin \alpha}_{\mid \alpha=0}
= 2n_m \cos(m\pi/2),
\ee
where $n_m=(n-1)^m+(n-3)^m+\cdots$, last term $1^m$ or $2^m$ according as $n$
is even or odd.

\cite{TaylorPEMS32}
\be
\frac{\sin n\alpha}{\sin \alpha}
=n+2\sum_{m\ge 1}(-)^m n_{2m}\frac{\alpha^{2m}}{(2m)!}.
\ee

If $n$ is even \cite{TaylorPEMS32}
\be
\frac{(-1)^{n/2-1}+\cos n\alpha}{\cos \alpha}
=c'+2\sum_{m\ge 1}(-)^m n'_{2m}\frac{\alpha^{2m}}{(2m)!},
\ee
where $n'_m=(n-1)^m-(n-3)^m+\cdots$,
$c'$ being zero or 2 according as $n/2$ is even or odd.
If $n$ is odd \cite{TaylorPEMS32}
\be
\frac{\cos n\alpha}{\cos \alpha}
=1+2\sum_{m\ge 1}(-)^m n'_{2m}\frac{\alpha^{2m}}{(2m)!}.
\ee
\cite{TaylorPEMS32}
\be
\frac{1-\cos n\theta}{1-\cos\theta}=n^2+2\sum_{m\ge 1}(-)^ma_{n,2m}\frac{\theta^{2m}}{(2m)!}
\ee
where $a_{n,m}\equiv \sum_{r=1}^{n-1}(n-r)r^m$.

\cite{TaylorPEMS32}
\be
\frac{\cosh n k-\cos n\theta}{\cosh k-\cos \theta}
= k_{n,0}+2\sum_{m\ge 1} (-)^mk_{n,2m}\frac{\theta^{2m}}{(2m)!}
\ee
where for $m>0$, $k_{n,m}=\sum_{r=1}^{n-1}r^m\frac{\sinh(n-r)k}{\sinh k}
=\sum_{r=1}^{n-1}r^m \frac{p^{n-r}-p^{-n+r}}{p-p^{-1}}$, $p=e^k$,
and $k_{n,0} = (\frac{p^n-1}{p-1})^2/p^{n-1}$.
Definite and indefinite integrals of these rational polynomials follow
by term-by-term intgration \cite{TaylorPEMS32}.

\be
\sin^{2k}\alpha+\cos^{2k}\alpha = 1 -\sum_{l=1}^{k-1}\binom{k}{l}\sin^{2l}\alpha\cos^{2k-2l}\alpha
.
\ee

\cite{CvijovicMathComp64}
\be
S_\nu(\alpha)\equiv\sum_{k=0}^\infty \frac{\sin(2k+1)\alpha}{(2k+1)^\nu};
\ee
\be
C_\nu(\alpha)\equiv\sum_{k=0}^\infty \frac{\cos(2k+1)\alpha}{(2k+1)^\nu};
\ee
\be
S_{2n+1}(\alpha)=\frac{(-1)^n}{4(2n)!}\pi^{2n+1}E_{2n}\left(\frac{\alpha}{\pi}\right) ;
\ee
\be
C_{2n}(\alpha)=\frac{(-1)^n}{4(2n-1)!}\pi^{2n}E_{2n-1}\left(\frac{\alpha}{\pi}\right) ;
\ee
\be
S_\nu(2\pi p/q)=\frac{1}{q^\nu}\sum_{s=1}^{q-1}\zeta(\nu,s/q)
\left[\sin(s2\pi p/q)-\frac{\sin(s4\pi p/q)}{2^\nu}
\right]
;
\ee
\be
C_\nu(2\pi p/q)=\frac{1}{q^\nu}\left\{\zeta(\nu)
\left(1-\frac{1}{2^\nu}\right)
+\sum_{s=1}^{q-1}\zeta(\nu,s/q)
\left[\cos(s2\pi p/q)-\frac{\cos(s4\pi p/q)}{2^\nu}
\right]
\right\}
.
\ee

\cite[1.448.1]{GR}
\be
\sum_{k=1}^{\infty}\frac{p^k\sin (kx)}{k}=\arctan \frac{p\sin x}{1-p\cos x},
\quad 0<x<2\pi, p^2\le 1,
\ee
with special case
\be
\sum_{k=1}^{\infty}\frac{\sin (kx)}{2^k k}=\arctan \frac{\sin x}{2-\cos x},
\quad 0<x<2\pi.
\ee

\cite[1.448.2]{GR}
\be
\sum_{k=1}^{\infty}\frac{p^k\cos (kx)}{k}=\ln \frac{1}{\sqrt{1-2p\cos x+p^2}},
\quad 0<x<2\pi,
\ee
with special case
\be
\sum_{k=1}^{\infty}\frac{\cos (kx)}{2^k k}= -\ln \sqrt{1-\cos x+1/4},
\quad 0<x<2\pi.
\ee

\cite[(1.10)]{HubbertNumAlg42}
\be
-\sqrt{2-2\cos\theta}=-\frac{4}{\pi}+\frac{2}{\pi}\sum_{k=1}^\infty
\frac{\cos k\theta}{(k-1/2)(k+1/2)}
\ee

\cite{AljarrahJCAM143}
\be
\cot\frac{nz}{2} = \frac{1}{n}\sum_{k=0}^{n-1}\frac{\sin z}{\cos(2\pi k/n)-\cos z}.
\ee

\cite{AljarrahJCAM143}
\be
2^{n-1}(\cos nz -1) = \prod_{k=0}^{n-1}[\cos z-\cos(2\pi k/n)].
\ee

\cite{WanAAM48}
\be
K(\sin t)=\sum_{n=0}^\infty
\frac{\Gamma(n+1/2)^2}{\Gamma(n+1)^2}\sin[(4n+1)t].
\ee
\cite{WanAAM48}
\be
E(\sin t)=
\sum_{n=0}^\infty
\frac{\Gamma(n+1/2)^2}{2\Gamma(n+1)^2}\sin[(4n+1)t]
+
\sum_{n=0}^\infty
\frac{(n+1/2)\Gamma(n+1/2)^2}{2(n+1)\Gamma(n+1)^2}\sin[(4n+3)t].
\ee

\cite{SrivastavaZAA19}
\be
\Cl_{2n+1}(\pi/2)=-2^{-2n-1}(1-2^{-2n})\zeta(2n+1), n\in \mathbb{N}.
\ee

\cite{SrivastavaZAA19}
\be
\Cl_{2n+1}(\pi/3)=\frac{1}{2}(1-2^{-2n})(1-3^{-2n})\zeta(2n+1), n\in \mathbb{N}.
\ee

\cite{SrivastavaZAA19}
\be
\Cl_{2n+1}(2\pi/3)=-\frac{1}{2}(1-3^{-2n})\zeta(2n+1), n\in \mathbb{N}.
\ee

\cite{SrivastavaZAA19}
\be
\sum_{k\ge 1}\frac{\cos(k\pi/2)}{k^s}=-2^{-s}(1-2^{1-s})\zeta(s), \Re s > 1.
\ee

\cite{LehmerAMM92}
\be
\sum_{m=1}^\infty \frac{m^{k-2}(2x)^{2m}}{\binom{2m}{m}}
=\frac{x}{2^{k-2}(1-x^2)^{k-1/2}}\left[
\arcsin x V_x(x^2)+x\sqrt{1-x^2}W_k(x^2)
\right]
,\quad k\ge 0
,
\ee
\cite{LehmerAMM92}
\be
\sum_{m=1}^\infty \frac{m^{k-2}4^m(\sin\theta)^{2m}} {\binom{2m}{m}}
=\frac{\sin 2\theta}{(2\cos^2\theta)^k}\left[
2\theta V_k(\sin^2\theta)+\sin 2\theta W_k(\sin^2\theta)
\right]
,
\ee
\cite{LehmerAMM92}
\begin{multline}
\sum_{m=1}^\infty \frac{(-1)^{m-1}m^{k-2}4^m(\sinh z)^{2m}} {\binom{2m}{m}}
\\
=\frac{\sinh 2 z}{(2\cosh^2 z)^k}\left[
2\log\{\sinh z+\cosh z\} V_k(-\sinh^2 z)+\sinh 2z W_k(-\sinh^2z)
\right]
,
\end{multline}
where $V_1(t)=1$, (see \cite[A156919]{sloane}) $W_1(t)=0$,
\begin{gather*}
V_{k+1}(t)=\{(2k-2)t+1\}V_k(t)+2(1-t)\delta V_k(t)
\\
W_{k+1}(t)=\{(2k-4)t+2\}W_k(t)+2(1-t)\delta W_k(t)+V_k(t)
\end{gather*}
and $\delta$ is the operator $x\frac{d}{dx}$.

\subsection{Expansions of Hyperbolic Functions}

\cite{BerndtBLMS15}
\be
\cos x \cosh x = \sum_{k\ge 0}\frac{(-)^k(2x^2)^{2k}}{(4k)!}
.
\ee

\cite{BerndtBLMS15}
\be
\sin x \sinh x = \sum_{k\ge 0}\frac{(-)^k(2x^2)^{2k+1}}{(4k+2)!}
.
\ee

\cite{BerndtBLMS15}
\be
\coth(2 x) \tanh x = 1-\sum_{k\ge 1}
\frac{2^{2k-1}(2^{2k}-1)(2k-1)B_{2k}x^{2k-2}}{(2k)!}
.
\ee

\cite{BorosSSA11}
\be
\sum_{j=0}^\infty \frac{x}{\sinh 2^{-j}x}-2^j = 1-\frac{x}{\tanh x}
.
\ee

\cite{BorosSSA11}
\be
\sum_{j=0}^\infty \frac{2^j-\coth 2^{-j}}{2^j\sinh 2^{-j}} =
\frac{1+4e^2-e^4}{1-2e^2+e^4}
.
\ee

\cite{GlasserMC28}
\be
\sum_{m,n=-\infty}^\infty
\frac{F(|2m+2n+1|)}{\cosh(2m+1)u\,\cosh\, 2nu}
=
2\sum_{n=0}^\infty \frac{(2n+1)F(2n+1)}{\sinh(2n+1)u},
\ee
for any summable function $F(x)$.

\cite{GlasserMC28}
\be
\sum_{m,n=-\infty}^\infty
\frac{F(m+n+1)+F(m-n)}{\sinh(2m+1)u\,\cosh\, (2n+1)u}
=
8\sum_{n=1}^\infty \frac{nF(n)}{\cosh(2nu)},
\ee
and
\be
\sum_{k,m,n=-\infty}^\infty
\frac{F(k+m+n+1)+F(k-m+n)}{\cosh(2ku) \cosh(2m+1)u\,\sinh(2n+1)u}
=
8\sum_{n=1}^\infty \frac{n^2F(n)}{\sinh(2nu)},
\ee
where $F$ is any sine transform (and hence odd).

\subsection{The Logarithmic Function}
\be
\ln(1-x)
= -\sum_{k=1}^\infty \frac{x^k}{k}.
\label{1.511}
\ee
\be
(1+x)\ln(1+x)=x+\frac{1}{1\cdot 2}x^2-\frac{1}{2\cdot 3}x^3
+\frac{1}{3\cdot 4}x^4-\frac{1}{4\cdot 5}x^5+\cdots
=x+\sum_{k=1}^\infty (-1)^{k+1}\frac{x^{k+1}}{k(k+1)}.
\label{1.511.+1}
\ee
By integration of this w.r.t.\ $x$ \cite[2.729]{GR}, \cite{Medinaarxiv1325}
\be
\frac{1}{2}(1+x)^2\ln(1+x)=
\frac{x}{2}+\frac{3x^2}{4}
+
\sum_{k=1}^\infty (-1)^{k+1}\frac{x^{k+2}}{k(k+1)(k+2)}.
\ee
\be
\frac{1}{6}(1+x)^3\ln(1+x)=
\frac{x}{6}+\frac{5x^2}{12}
+\frac{11x^3}{36}
+
\sum_{k=1}^\infty (-1)^{k+1}\frac{x^{k+3}}{k(k+1)(k+2)(k+3)}.
\ee
\be
\frac{1}{24}(1+x)^4\ln(1+x)=
\frac{x}{24}
+\frac{7x^2}{48}
+\frac{13x^3}{72}
+\frac{25x^4}{288}
+
\sum_{k=1}^\infty (-1)^{k+1}\frac{x^{k+4}}{k(k+1)(k+2)(k+3)(k+4)}.
\ee
\be
x+(1-x)^l\ln(1-x)=
\sum_{i=2}^l\tau_{i,l}x^i-(-1)^l l!
\sum_{k=1}^\infty \frac{x^{k+l}}{k(k+1)(k+2)\cdots (k+l)},
\ee
where
\be
\tau_{2,l}=l-1/2,\quad l\ge 2;
\quad
i\tau_{i,l}=(-1)^i\binom{l-1}{i-1}-l\tau_{i-1,l-1},\quad i\ge 3.
\ee

\cite{ZhaoDM281}
\be
[t^k](-\log(1-t))^m = \frac{m!}{n!} \begin{bmatrix}{c}n\\ m\end{bmatrix}
\therefore
\sum_{k=1}^{n+1}(k-1)!\begin{Bmatrix} n \\ k-1\end{Bmatrix}
[t^k]f(-\ln (1-t)) = \sum_{j=0}^{n+1}\frac{(-)^jj!f_j}{n+1}
\binom{n+1}{j}B_{n-j+1}.
\ee

\cite{Boyadzhievarxiv2012}
\be
\frac{-x}{(1-x)\ln(1-x)}=\sum_{n=0}^\infty d_n\frac{x^n}{n!}
\ee
where $d_n$ are the Cauchy numbers of the second kind with
\be
(-)^nd_n=\sum_{k=0}^n \left[\begin{array}{c}n\\k\end{array}\right]\frac{1}{k+1}.
\label{eq.dnBoyad}
\ee

\cite{Boyadzhievarxiv2012}
\be
\frac{x}{\ln(1+x)}=\sum_{n=0}^\infty c_n\frac{x^n}{n!}
\ee
where $c_n$ are the Cauchy numbers of the first kind.

\cite{LafontMer}
\be
\log n = -\sum_{s\ge 1}\frac{\alpha(s,n)}{s}
\ee
with
\be
\alpha(s,n)
\equiv
\left\{
\begin{array}{ll}
n-1,& n\mid s; \\
-1 & n\nmid s.\\
\end{array}
\right.
\ee

Superposition of \cite[1.513.1]{GR} and \cite[1.511]{GR}, see \cite[A165998]{sloane}:
\be
\frac{1}{3x}\ln\frac{1+x}{(1-x)^2}
=
1+\frac{x}{6}+\frac{x^2}{3}+\frac{x^3}{12}+\cdots
+\frac{x^{2j}}{2j+1}+\frac{x^{2j+1}}{3(2j+1)}+\cdots
\ee

\cite{Amdeberhanarxiv3663}
\be
\sum_{j=k}^\infty \frac{|S_j^{(k)}|b^j}{j!}=\frac{\ln^k(1+b)}{k!}.
\ee

\cite{Israel0907}
\be
\sum_{k=0}^\infty \frac{x^k}{\binom{k+L}{L}}
=\,_2F_1(1,1;L+1;x) = L\sum_{j=0}^{L-2}\frac{(x-1)^j}{(L-j-1)x^{j+1}}
-L(x-1)^{L-1}\frac{\ln(1-x)}{x^L}
;\quad L\ge 1.
\ee

\cite{BricePRA6}
\be
\frac{1}{2x}\left\{
1-\ln(1+x)-\frac{1-x}{\sqrt{x}}\arctan \sqrt{x}
\right\}
=\sum_{k=1}^\infty \frac{(-1)^{k+1}x^{k-1}}{(2k-1)2k(2k+1)};
\ee
$0< x\le 1$.

\cite{MezoCEJM}
\be
-\frac{\ln(1-t)}{(1-t)^r}
=
\sum_{n=0}^\infty H_n^{(r)}t^n,
\label{eq.hyperH}
\ee
where the Hyperharmonic numbers are defined by $H_n=\sum_{k=1}^n (1/k)$
and $H_n^{(r)}=\sum_{k=1}^n H_k^{(r-1)} = \binom{n+r-1}{r-1}(H_{n+r-1}-H_{r-1})$.

\subsubsection{Series of logarithms}
\cite[A340440]{sloane}
\be
\sum_{k=2}^\infty \frac{\log k}{k^2-1} = \sum_{k\ge 2} \sum_{i\ge 0}  \frac{\log k}{k^{2(1+i)}} = -\sum_{i \ge 1}\zeta'(2i)
\approx 1.023138726427939\ldots
\ee

\be
\sum_{k=2}^\infty \frac{\log k}{k^2- s} 
= 
-\sum_{i\ge 1}s^{i-1} \zeta'(2i), \quad |s| <4.
\ee

Closely related to \eqref{eq.invk}:
\be
\sum_{k=2}^\infty \frac{\log k}{k^j(k^2- s)^l} 
= 
-\sum_{i\ge 0}
\frac{(l)_i}{i!}s^i 
\zeta'(j+2l+2i).
\ee

\cite[A343920]{sloane}
\be
\sum_{k=4,6,8,10,\ldots} \frac{\log k}{k^2-4}
=
\frac{3}{16}\log 2 
+\frac{1}{4}\sum_{k\ge 2} \frac{\log k}{k^2-1}\approx 0.38574977796197,
\ee
where the remaining sum is given above.

\be
\sum_{k\ge 1}\frac{\log(k+\epsilon)}{k^2}
=
-\zeta'(2) -\sum_{i\ge 1}\frac{(-\epsilon)^i}{i}\zeta(i+2).
\ee

\cite[A343919]{sloane}
\be
\sum_{k\ge 1}\frac{\log(k+1/2)}{k^2}
=
-\zeta'(2) -\sum_{i\ge 1}\frac{(-)^i}{i2^i}\zeta(i+2)
\approx 1.43506248930485670040268217252716365\ldots
\ee

\cite[A343918]{sloane}
\be
\sum_{k=2}^\infty \frac{k\log k}{(k^2-1)^2} 
= \sum_{k\ge 2} \log k \sum_{i\ge 1} \frac{i}{k^{2i+1}}
= -\sum_{i \ge 1}i \zeta'(2i+1) \approx 0.2814827931466416082856987201176296\ldots
\ee

\cite[A340485]{sloane}
\be
\sum_{k=2}^\infty \frac{\log k}{(k^2-1)^2} = -\sum_{i \ge 1}i \zeta'(2i+2)\approx 0.1073253716.
\ee

\cite[A340484]{sloane}
\be
\sum_{k=2}^\infty \frac{\log^2 k}{k^2-1} = \sum_{i \ge 1}\zeta''(2i).
\ee

\subsection{The Inverse Trigonometric and Hyperbolic Function}

\cite{BorosSSA11}
\be
\sum_{k=1}^\infty 2^{-k}\tan\frac{x}{2^k} = \frac{1}{2^n}
\cot\frac{x}{2^n}-\cot x
.
\ee

\cite{BorosSSA11}
\be
\sum_{k=1}^\infty \csc\frac{x}{2^{k-1}} 
=\cot\frac{x}{2^n}-\cot x
.
\ee

\cite{BorosSSA11}
\be
\sum_{k=1}^\infty \arctan\frac{2x^2}{k^2} =
\frac{\pi}{4}-\arctan\frac{\tanh \pi x}{\tan \pi x}
.
\ee

\cite{BorosSSA11}
\be
\sum_{k=1}^\infty (-)^{k-1}\arctan\frac{2x^2}{k^2} =
-\frac{\pi}{4}-\arctan\frac{\sinh \pi x}{\sin \pi x}
.
\ee

\section{Indefinite Integrals of Elementary Functions}

\subsection{Rational Functions}
Aids to partial fraction decompositions:
\cite{WitulaAppMathComp197}
\be
\frac{1}{s^n(s^2+as+b)}
=\frac{-\alpha_{n-1}s+\alpha_n}{b^n(s^2+as+b)}+\sum_{k=0}^{n-1}\frac{\alpha_k}
{b^{k+1}s^{n-k}};
\ee
where $b>0$, $a^2-4b<0$, $\alpha_0\equiv 1$,
\be
\alpha_m\equiv (-1)^m\sqrt{b^m}U_m\left(\frac{a}{2\sqrt b}\right)
=\sum_{k=0}^{\lfloor m/2\rfloor}(-1)^{k+m}\binom{m-k}{k}a^{m-2k}b^k
.
\ee

\cite{WitulaAppMathComp197}
\be
\frac{\alpha s+\beta}{(s^2+as+b)(s^2+kas+kb)}
=\frac{1}{(k-1)b^2}\left( \frac{L(s+a)+\beta b}{s^2+as+b}-\frac{L(s+ka)+\beta b}{s^2+kas+kb}\right),
\ee
$k\neq 1$, $ab\neq 0$, $L\equiv \left|\begin{array}{cc} \alpha & \beta \\ a & c \end{array}\right|$.

\cite{WitulaAppMathComp197}
\be
\frac{1}{(pq+ap+b)(pq+cp+b)}
=\left|\begin{array}{cc}a & c \\ b & b\end{array}\right|^{-1}
\left(\frac{q+a}{pq+ap+b}-\frac{q+c}{pq+cp+b}\right),
\ee
$a\neq c$, $b\neq 0$.

\cite{WitulaAppMathComp197}
\be
\frac{\alpha p+\beta}{(pq+ap+b)(pq+cp+b)}
=\left|\begin{array}{cc}a & c \\ b & b\end{array}\right|^{-1}
\left(\frac{\left|\begin{array}{cc}\alpha & \beta \\ q+c & b\end{array}\right|}
{pq+cp+b}-\frac{
\left|\begin{array}{cc}\alpha & \beta \\ q+a & b\end{array}\right|
}{pq+ap+b}\right),
\ee
$a\neq c$, $b\neq 0$.

\cite{WitulaAppMathComp197}
\be
\frac{as+b}{(s^2+\alpha s+\beta)(s^2+\gamma s+\delta)}
=
\frac{1}{k}\left(
\frac{\left(b-a\frac{\delta-\beta}{\gamma-\alpha}\right)s+b\alpha-a\beta -b\frac{\delta-\beta}{\gamma-\alpha}}
{s^2+\alpha s+\beta}-\frac{
\left(b-a\frac{\delta-\beta}{\gamma-\alpha}\right)s+b\gamma-a\delta -b\frac{\delta-\beta}{\gamma-\alpha}
}{s^2+\gamma s+\delta}
\right),
\ee
$\alpha\neq \gamma$.

\cite{MatharArxiv1301}
\be
\frac{1}{n^{2s}(n+1)^{2s}}
=
\sum_{t=1}^{2s}\binom{4s-t-1}{2s-1}\left[\frac{(-1)^t}{n^t}+\frac{1}{(n+1)^t}\right]
.
\ee

\cite[170.]{Dwight}
\be
\int \frac{dx}{a^4+x^4}=\frac{1}{4a^3\surd 2}\log\frac{x^2+ax\surd 2+a^2}{x^2-ax\surd 2+a^2}
+\frac{1}{2a^3\surd 2}\arctan\frac{ax\surd 2}{a^2-x^2}.
\ee

\cite[171.]{Dwight}
\be
\int \frac{dx}{a^4-x^4}=\frac{1}{4a^3}\log\left|\frac{a+x}{a-x}\right|
+\frac{1}{2a^3}\arctan\frac{x}{a}.
\ee

\cite[170.1]{Dwight}
\be
\int \frac{x dx}{a^4+x^4}=\frac{1}{2a^2}\arctan\frac{x^2}{a^2}
.
\ee

\cite[170.2]{Dwight}
\be
\int \frac{x^2 dx}{a^4+x^4}=-\frac{1}{4a\surd 2}\log\frac{x^2+ax\surd 2+a^2}{x^2-ax\surd 2+a^2}
+\frac{1}{2a\surd 2}\arctan\frac{ax\surd 2}{a^2-x^2}.
\ee

\cite[170.3]{Dwight}
\be
\int \frac{x^3 dx}{a^4+x^4}=\frac{1}{4}\log(a^4+x^4).
\ee

\cite[171.1]{Dwight}
\be
\int \frac{x dx}{a^4-x^4}=\frac{1}{4a^2}\log\left|\frac{a^2+x^2}{a^2-x^2}\right|
.
\ee

\cite[171.2]{Dwight}
\be
\int \frac{x^2 dx}{a^4-x^4}=\frac{1}{4a}\log\left|\frac{a+x}{a-x}\right|
-\frac{1}{2a}\arctan\frac{x}{a}
.
\ee

\cite[171.3]{Dwight}
\be
\int \frac{x^3 dx}{a^4-x^4}=-\frac{1}{4}\log|a^4-x^4|.
\ee

\cite[173]{Dwight}
\be
\int \frac{dx}{x(a+bx^m)}=\frac{1}{am}\log|\frac{x^m}{a+bx^m}|.
\ee

\subsection{Algebraic Functions}
\cite[186.11,188.11]{Dwight}
\be
\int \frac{dx}{(a^2+b^2x)x^{1/2}}=\frac{2}{ab}\arctan\frac{bx^{1/2}}{a}.
\ee
\be
\int \frac{dx}{(a^2-b^2x)x^{1/2}}
= \frac{1}{2ab}\log\left|\frac{a+bx^{1/2}}{a-bx^{1/2}}\right|.
\ee
\be
\int \frac{dx}{(1-a-bx)\sqrt{a+bx}}
= \frac{2}{b}\arctanh\sqrt{a+bx}.
\ee

\cite[185.11,187.11]{Dwight}
\be
\int \frac{x^{1/2}dx}{a^2+b^2x}=\frac{2x^{1/2}}{b^2}-\frac{2a}{b^3}\arctan\frac{bx^{1/2}}{a}.
\ee
\be
\int \frac{x^{1/2}dx}{a^2-b^2x}
=- \frac{2x^{1/2}}{b^2}+\frac{a}{b^3}\log\left|\frac{a+bx^{1/2}}{a-bx^{1/2}}\right|.
\ee

\cite[185.13,187.13]{Dwight}
\be
\int \frac{x^{3/2}dx}{a^2+b^2x}=\frac{2}{3}\,\frac{x^{3/2}}{b^2}-\frac{2a^2x^{1/2}}{b^4}+\frac{2a^3}{b^5}\arctan\frac{bx^{1/2}}{a}.
\ee
\be
\int \frac{x^{3/2}dx}{a^2-b^2x}
=- \frac{2}{3}\,\frac{x^{3/2}}{b^2}-\frac{2a^2x^{1/2}}{b^4}+\frac{a^3}{b^5}\log\left|\frac{a+bx^{1/2}}{a-bx^{1/2}}\right|.
\ee

\cite[186.21,188.21]{Dwight}
\be
\int \frac{dx}{(a^2+b^2x)^2x^{1/2}}
=\frac{x^{1/2}}{a^2(a^2+b^2x)}+\frac{1}{a^3b}\arctan\frac{bx^{1/2}}{a}.
\ee
\be
\int \frac{dx}{(a^2-b^2x)^2x^{1/2}}
= \frac{x^{1/2}}{a^2(a^2-b^2x)}+\frac{1}{2a^3b}\log\left|\frac{a+bx^{1/2}}{a-bx^{1/2}}\right|.
\ee

\cite[185.21,187.21]{Dwight}
\be
\int \frac{x^{1/2}dx}{(a^2+b^2x)^2}
= -\frac{x^{1/2}}{b^2(a^2+b^2x)}+\frac{1}{ab^3}\arctan\frac{bx^{1/2}}{a}.
\ee
\be
\int \frac{x^{1/2}dx}{(a^2-b^2x)^2}
= \frac{x^{1/2}}{b^2(a^2-b^2x)} -\frac{1}{2ab^3}\log\left|\frac{a+bx^{1/2}}{a-bx^{1/2}}\right|.
\ee

\cite[185.23,187.23]{Dwight}
\be
\int \frac{x^{3/2}dx}{(a^2+b^2x)^2}=\frac{2x^{3/2}}{b^2(a^2+b^2x)}+\frac{3a^2x^{1/2}}{b^4(a^2+b^2x)}-\frac{3a}{b^5}\arctan\frac{bx^{1/2}}{a}.
\ee
\be
\int \frac{x^{3/2}dx}{(a^2-b^2x)^2}=\frac{3a^2x^{1/2}-2b^2x^{3/2}}{b^4(a^2-b^2x)}
 -\frac{3a}{2b^5}\log\left|\frac{a+bx^{1/2}}{a-bx^{1/2}}\right|.
\ee

\cite[186.13,188.13]{Dwight}
\be
\int \frac{dx}{(a^2+b^2x)x^{3/2}}=-\frac{2}{a^2x^{1/2}}-\frac{2b}{a^3}\arctan\frac{bx^{1/2}}{a}.
\ee
\be
\int \frac{dx}{(a^2-b^2x)x^{3/2}}=-\frac{2}{a^2x^{1/2}}
 +\frac{b}{a^3}\log\left|\frac{a+bx^{1/2}}{a-bx^{1/2}}\right|.
\ee

\cite[186.23,188.23]{Dwight}
\be
\int \frac{dx}{(a^2+b^2x)^2x^{3/2}}=-\frac{2}{a^2(a^2+b^2x)x^{1/2}}
-\frac{3b^2x^{1/2}}{a^4(a^2+b^2x)}
-\frac{3b}{a^5}\arctan\frac{bx^{1/2}}{a}.
\ee
\be
\int \frac{dx}{(a^2-b^2x)^2x^{3/2}}=-\frac{2}{a^2(a^2-b^2x)x^{1/2}}
+\frac{3b^2x^{1/2}}{a^4(a^2-b^2x)}
 +\frac{3b}{2a^5}\log\left|\frac{a+bx^{1/2}}{a-bx^{1/2}}\right|.
\ee

\cite[189.2,189.4]{Dwight}
\be
\int \frac{dx}{(a^4+x^2)x^{1/2}}=
\frac{1}{2a^3\surd 2}\log\frac{x+a\sqrt{2x}+a^2}{x-a\sqrt{2x}+a^2}
+\frac{1}{a^3\surd 2}\arctan\frac{a\sqrt{2x}}{a^2-x}.
\ee
\be
\int \frac{dx}{(a^4-x^2)x^{1/2}}=
\frac{1}{2a^3}\log\left|\frac{a+x^{1/2}}{a-x^{1/2}}\right|
+\frac{1}{a^3}\arctan\frac{x^{1/2}}{a}.
\ee

\cite[188.23,189.3]{Dwight}
\be
\int \frac{x^{1/2}dx}{a^4+x^2}=
-\frac{1}{2a\surd 2}\log\frac{x+a\sqrt{2x}+a^2}{x-a\sqrt{2x}+a^2}
+\frac{1}{a\surd 2}\arctan\frac{a\sqrt{2x}}{a^2-x}.
\ee
\be
\int \frac{x^{1/2}dx}{a^4-x^2}=
\frac{1}{2a}\log\left|\frac{a+x^{1/2}}{a-x^{1/2}}\right|
-\frac{1}{a}\arctan\frac{x^{1/2}}{a}.
\ee

\be
\int \frac{x^n}{\sqrt{a+bx}}dx
=
\frac{2\sqrt{a+bx}}{b^{n+1}}(-a)^n\sum_{k=0}^n\binom{n}{k}\frac{(-z/a)^k}{2k+1}.
\ee

\cite[2.241.2]{GR}
\be
\int (\alpha+\beta x)^n (a+bx)^{m-1/2} dx
=2\sqrt{z^{2m+1}}\frac{1}{b^{n+1}}\sum_{k=0}^n \binom{n}{k}
(-\Delta)^k
\frac{(\beta z)^{n-k}}{2n+2m-2k+1}
.
\ee
where $z\equiv a+bx$, $\Delta =a\beta-b\alpha$.

\cite[2.252]{GR}
\be
\int \frac{dt}{(t^2+p)^k\sqrt{c(t^2+q)}}
=
\frac{1}{\surd c}\int dv
\frac{(1-v^2)^{k-1}}{[p+v^2(q-p)]^k}
\ee
where
$$
v=\frac{t}{\sqrt{t^2+q}};\quad t=v\sqrt\frac{q}{1-v^2};\quad
\frac{dt}{dv} = \frac{\sqrt{q}}{(1-v^2)^{3/2}}.
$$

\begin{multline}
\int\frac{z^mdx}{t^n\sqrt{z}} = 
-z^m\sqrt{z}\big\{ \frac{1}{(n-1)\Delta}\,\frac{1}{t^{n-1}}
\\
+\sum_{k=2}^{n-1}
\frac{(2n-2m-3)(2n-2m-5)\cdots (2n-2m-2k+1)(-b)^{k-1}}
{2^{k-1}(n-1)(n-2)\cdots (n-k)\Delta^k}\,\frac{1}{t^{n-k}}
\big\}\\
-\frac{(2n-3m-3)(2n-3m-5)\cdots (-2m+3)(-2m+1)(-b)^{n-1}}
{2^{n-1}(n-1)!\Delta^n}\int\frac{z^mdx}{t\sqrt{z}}
\end{multline}
where $z=a+bx$ and $t=\alpha+\beta x$ and $\Delta\equiv a\beta-\alpha b$.

\cite{vanderPoortenLNCS2369}
\be
\int \frac{f_m(z)}{\sqrt{D_m(z)}}dz
=\log(x_m(z)+y_m(z)\sqrt{D_m(z)}),
\ee
if for example
\be
f=4z+2;\quad D=z^4+8(z+1);\quad x=z^4-2z^3+2z^2+4z-4;\quad y=z^2-2z+2.
\ee
or
\be
f=5z+1;\quad D=(z^2+1)^2+4z;\quad x=z^5-z^4+3z^3+z^2+2;\quad y=z^3-z^2+2z.
\ee
or
\be
f=6z+2;\quad D=(z^2+2)^2+8z;\quad x=z^6-2z^5+8z^4-4z^3+8z^2+8z;\quad y=z^4-2z^3+6z^2-4z+4.
\ee
or
\be
f=3z-s;\quad D=(z^2-s^2)^2+t(z-s);\quad x=1+2(z+s)(z^2-s^2)/t;\quad y=2(z+s)/t.
\ee

\subsection{The Exponential Function}

\cite{Brien}
\be
\int\frac{e^{-\beta x}}{x^n}dx = e^{-\beta x}[
\sum_{k=0}^{n-2}(-)^k\frac{(n-k-2)!}{(n-1)!}\frac{\beta ^k}{x^{n-k-1}}
]
+(-\beta)^{n-1}\frac{\Ei(-\beta x)}{(n-1)!}; \beta\neq 0
\ee

\cite{Brien}
\be
\int\frac{e^{-\beta x}}{x^2}dx = -\frac{e^{-\beta x}}{x}
-\beta\Ei(-\beta x);\beta \neq 0
\ee

\be
\int x^{m+2}e^{-ax^2} dx=\frac{m+1}{2a}\left[-x^{m+1}e^{-ax^2}+\int x^me^{-ax^2}dx\right];
\ee
$m\neq -1$.

\be
\int \frac{e^{-x}}{(1+\alpha e^{-x})^s} dx
=
\frac{1}{\alpha(s-1)(1+\alpha e^{-x})^{s-1}},\quad s>1.
\ee

\cite{Brien}
\be
\int x^m e^{\pm ax^n} dx = \frac{1}{na^\gamma}\int s^{\gamma-1}{e^{\pm s}}ds
=\pm \frac{x^{m+1-n}}{na}e^{\pm ax^n} \mp\frac{m+1-n}{na}\int x^{m-n} e^{\pm a^n}dx
\ee
where $s=ax^n$, $\gamma=(m+1)/n$.

\cite{Brien}
\be
\int \frac{e^{\pm ax^n}}{x^m} dx = \frac{1}{na^z}\int s^{z-1}{e^{\pm s}}ds
\ee
where $s=ax^n$, $z=(1-m)/n$.

\cite{Brien}
\be
\int \exp(a^x) dx = \frac{\Ei(a^x)}{\ln a}; \ln a >0.
\ee

\cite{Brien}
\be
\int \frac{a^x}{x^b} dx = \frac{(-1)^b}{(\ln a)^{1-b}}\Gamma(1-b,-x\ln a); \ln a >0.
\ee

\cite{Brien}
\be
\int \frac{2^x}{x^2} dx = (\ln 2)\li(2^x)-\frac{2^x}{x}.
\ee

\cite{Brien}
\be
\int a^{x^n} dx = -\frac{1}{n(-\ln a)^{1/n}}\Gamma(\frac1n,-x^n\ln a); n\neq 0, \ln a <0.
\ee

\cite{Brien}
\be
\int x^m a^{x^n} dx = -\frac{1}{n(\ln a)^{(m+1)/n}}\Gamma(\frac{m+1}{n},-x^n\ln a); n\neq 0, \ln a >0.
\ee

\cite{Brien}
\be
\int \frac{ dx}{a^{x^n}} = -\frac{1}{n(\ln a)^{1/n}}\Gamma(\frac1n,x^n\ln a); n\neq 0, \ln a >0.
\ee

\cite{Brien}
\be
\int a^{x^2}dx = \frac{1}{2}\sqrt{\frac{\pi}{\ln a}}\erfi(\sqrt{\ln a}x);\ln a >0
\ee

\cite{Brien}
\be
\int \frac{dx}{a^{x^2}} = \frac{1}{2}\sqrt{\frac{\pi}{\ln a}}\erf(\sqrt{\ln a}x);\ln a >0
\ee

\cite{Brien}
\be
\int \frac{e^{ax}-1}{e^{ax}+1}dx=\frac{2}{a}\ln(e^{ax}+1)-x=\frac{2}{a}\ln(e^{ax/2}+e^{-ax/2}).
\ee

\cite{Brien}
\be
\int \frac{e^{ax}-m}{e^{ax}+n}dx=\frac{m+n}{na}\ln(e^{nax/(m+n)}+e^{-max/(m+n)}).
\ee

\cite{Brien}
\be
\int \frac{e^{ax}+m}{e^{ax}-n}dx=\frac{m+n}{na}\ln(e^{nax/(m+n)}-ne^{-max/(m+n)}).
\ee

\cite{Brien}
\be
\int \frac{e^{ax}+m}{e^{ax}+n}dx=\frac{1}{na}\ln[(e^{ax}+n)^{n-m}e^{max}].
\ee

\cite{Brien}
\be
\int \frac{e^{ax}-m}{e^{ax}-n}dx=\frac{1}{na}\ln[(e^{ax}-n)^{n-m}e^{max}].
\ee

\cite{Brien}
\be
\int \frac{x+1}{e^{ax}}dx = -e^{-ax}\frac{ax+a+1}{a^2}.
\ee

\cite{Brien}
\be
\int \frac{x}{e^{x}-1}dx = x\log(1-e^x)-\frac{x^2}{2}+\Li_2(e^x).
\ee

\cite{Brien}
\be
\int \frac{xe^{ax}}{(1+ax)^2}dx = \frac{e^{ax}}{a^2(1+ax)}.
\ee

\subsection{Hyperbolic Functions}
\cite{CoffeyJMP49}
\be
\int_0^y \frac{x\cosh x}{\cosh 2x-\cos 2t}dx
=\frac{\csc t}{4}[-4\Cl_2(\pi+t)-\Cl_2(2\omega_1)+Cl_2(2(\omega_1-t))+\Cl_2(2t)-\Cl_2(2\omega_3)
+\Cl_2(2(\omega_3+t))]
\ee
where
$\omega_1\equiv \arctan\frac{e^y\sin t}{1+e^y\cos t}$,
$\omega_3\equiv \arctan\frac{e^y\sin t}{1-e^y\cos t}$.

\be
\int \frac{\sinh(ax)}{\cosh^m x}dx
=
\frac{1}{a}\frac{\cosh(ax)}{\cosh^m x}
+\frac{m}{2a}\int 
\frac{\sinh[(a+1)x]}{\cosh^{m+1} x}dx
-\frac{m}{2a}\int 
\frac{\sinh[(a-1)x]}{\cosh^{m+1} x}dx
.
\ee

\subsection{Trigonometric Functions}

\begin{multline}
\int \sin x \frac{(a+\sin x)^2}{(a+\sin x)^2+(b+\cos x)^2}dx
=
\frac{1}{4}
\left[
\frac{2b^2}{a^2+b^2}+\frac{a^2-b^2}{(a^2+b^2)^2}-3
\right]
\cos(x)
\nonumber
\\
-\frac{ab}{2}
\left[
\frac{1}{a^2+b^2}-\frac{1}{(a^2+b^2)^2}
\right]
\sin (x)
\nonumber
\\
+
\frac{1}{8}\frac{b\cos(2x)-a\sin(2x)}
{a^2+b^2}
+
\frac{b}{4}\left[
1
+\frac{2(a^2-b^2)}{(a^2+b^2)^2}
+\frac{-3a^2+b^2}{(a^2+b^2)^3}
\right]
 \ln\sqrt{(\cos x +b)^2 + (\sin x+a)^2}
\\
+\frac{a}{4}\left[
1
-\frac{4b^2}{(a^2+b^2)^2}
-\frac{a^2-3b^2}{(a^2+b^2)^3}
\right]
\arctan \frac{\sin x +a}{\cos x+b}
+\frac{a}{4}\left[
\frac{4b^2}{(a^2+b^2)^2}+\frac{-3b^2+a^2}{(a^2+b^2)^3}
\right]
x
.
\end{multline}

\cite{ChoiJMAA231}
\be
\int_0^z \pi t \cot \pi t dt = z\log(2\pi)+\log\frac{G(1-z)}{G(1+z)}.
\ee
where $G$ is the reciprocal of the double Gamma function.

\cite{ChoiJMAA231}
\be
\int_0^z \pi t \tan \pi t dt = 
-\frac12 \log\frac{\cos \pi z}{\pi}-z \log(2\pi)
-\log\Gamma(\frac12-z)
+\log\frac{G(1/2+z)}{G(1/2-z)}
=
\log\frac{G(1-z)}{G(1+z)}
-\frac12
\log\frac{G(1-2z)}{G(1+2z)}
.
\ee

\cite{ChoiJMAA231}
\be
\int_0^z(\pi t \tan \pi t)^2 dt =
-\frac{\pi^2 z^3}{3}+\pi z^2\tan \pi z +\log\frac{\cos \pi z}{\pi}
+2z\log(2\pi)+2\log\frac{G(3/2-z)}{G(1/2+z)}.
\ee

\cite{ChoiJMAA231}
\be
\int_0^z(\pi t \cot \pi t)^2 dt =
-\frac{\pi^2 z^3}{3}-\pi z^2\cot \pi z
+2z\log(2\pi)+2\log\frac{G(1-z)}{G(1+z)}.
\ee

\subsection{Rational Functions of Trigonometric Functions}

\cite{TaylorPEMS32}
\be
\int_0 \frac{\sin n\theta}{\sin\theta }d\theta
= 
\frac{2\sin(n-1)\theta}{n-1}
+\frac{2\sin(n-3)\theta}{n-3}+\cdots +c\theta
\ee
where $c$ is zero or unit depending on $n$ being even or odd.

\cite{TaylorPEMS32}
\be
\int_0 \frac{\cos(2n-1) \theta}{\cos\theta }d\theta
= 
\frac{2\sin(2n-2)\theta}{2n-2}
-\frac{2\sin(2n-4)\theta}{2n-4}+\cdots +(-1)^{n-1}\theta.
\ee

\cite{TaylorPEMS32}
\be
\int_{\pi/2} \frac{\sin 2n\theta}{\cos\theta }d\theta
= 
-\frac{2\cos(2n-1)\theta}{2n-1}
+\frac{2\cos(2n-3)\theta}{2n-3}-\cdots .
\ee

\cite{TaylorPEMS32}
\be
\int_{\pi/2} \frac{1-\cos 2n\theta}{\sin\theta }d\theta
= 
-\frac{2\cos(2n-1)\theta}{2n-1}
-\frac{2\cos(2n-3)\theta}{2n-3}-\cdots .
\ee
More of this type by integration of right hand sides of formulas in section \ref{sec.Fser}.

\cite{CoffeyJMP49}
\be
\int_0^b \frac{a}{\sin a}da = \Cl_2(b)-\Cl_2(b+\pi)+i\frac{\pi}{4}(2b-\pi)+b\ln\frac{1-e^{ib}}{1+e^{ib}}+i\frac{\pi^2}{4}.
\ee

\cite{CoffeyJMP49}
\be
\int_0^b \frac{a^2}{\sin^2 a}da = \Cl_2(2b)+ib(\pi-2b)
+b[2\ln(1-e^{2ib})-b\cot b].
\ee

\cite{CoffeyJMP49}
\be
\int_0^b \frac{a}{\tan a}da = \Cl_2(b)+\Cl_2(b+\pi)-i\frac{\pi^2}{4}
+b[\ln(1-e^{2ib})-ib/2]+\frac{i}{2}(\pi b-b^2+\frac{\pi^2}{2}).
\ee

\cite{CoffeyJMP49}
\be
\int_0^b \frac{a^2}{\tan^2 a}da = \Cl_2(2b)+ib(\pi-2b)
-\frac13 b^3+b[2\ln(1-e^{2ib})-b\cot b]
.
\ee

\cite{CoffeyJMP49}
\be
\int_0^b \frac{x}{\sin x+a}dx = \frac{1}{\sqrt{1-a^2}}
[
-2\Cl_2(\pi+\phi)+ib\phi+\Cl_2(\pi-b+\phi)-Cl_2(\pi-b-\phi)
+ib\ln\frac{1-e^{-ib}/u_+}{1-e^{-ib}/u_-}
]
,
\ee
where $\phi=\arctan(a/\sqrt{1-a^2})$, $u_{\pm}=ia\pm \sqrt{1-a^2}$.

\cite{CoffeyJMP49}
\be
\int_0^b \frac{x}{\tan x+a}dx = \frac{1}{8v_+}
\frac{2}{1+a^2}[\Cl_2(\phi_a)-\Cl_2(\phi_a-2b)+2b\ln(1-e^{i(\phi_a-2b)})]\frac{2ib}{1+ia}[b+\frac{b-\phi_a+\pi}{1-ia}]
,
\ee
where $\phi_a=-\arctan(2a/\sqrt{1-a^2})$, $v_+=\sqrt{\frac{1+ia}{1-ia}}$.

\cite{ChoiJMAA231}
\be
\int_0^z (\frac{\pi t}{\cos \pi t})^2dt
=
\pi z^2\tan\pi z+\log\frac{\cos \pi z}{\pi}+2z\log(2\pi)
+2\log\Gamma(\frac12-z)-2\log\frac{G(\frac12+z)}{G(\frac12-z)}.
\ee
where $G$ is the reciprocal of the double Gamma Function.

\cite{ChoiJMAA231}
\be
\int_0^z (\frac{\pi t}{\sin \pi t})^2dt
=
-\pi z^2\cot\pi z+2z\log(2\pi)
+2\log\frac{G(1-z)}{G(1+z)}.
\ee

\cite{ChoiJMAA231}
\be
\int_0^z (\frac{\pi t}{\sin \pi t})^2\cos \pi t dt
=
-\frac{\pi z^2}{\sin \pi z}-2\log\frac{\cos(\pi z/2)}{\pi}
+4\log\frac{G(1-z/2)}{G(1+z/2)}
+4\log\frac{G(1/2+z/2)}{G(3/2-z/2)}
\ee

\cite{ChoiJMAA231}
\be
\int_0^z \frac{t\sin at}{\cos b t -\cos a t}dt
=
\frac{2az}{a^2-b^2}\log(2\pi)
+\frac{2\pi}{(a+b)^2}\log\frac{G(1-(a+b/2\pi)z)}{G(1+(a+b/2\pi)z)}
+\frac{2\pi}{(a-b)^2}\log\frac{G(1-(a-b/2\pi)z)}{G(1+(a-b/2\pi)z)}
.
\ee

\cite{ChoiJMAA231}
\be
\int_0^z \frac{a}{1+a^2+2a\cos t}dt
=
\frac{2a}{1-a^2}\arctan(\frac{1-a}{1+a}\tan\frac{z}{2}).
\ee

\cite{ChoiJMAA231}
\be
\int_0^z \frac{\cos t}{1+a^2+2a\cos t}dt
=
\frac{z}{2a}-\frac{1+a^2}{a(1-a^2)}\arctan(\frac{1-a}{1+a}\tan\frac{z}{2}).
\ee

\cite[(6.6.12)]{Nahin}
\be
\int_0^\phi \frac{\cos^2 x}{\sqrt{1+\cos^2 x}}dx = \sqrt{2}E(1/\surd 2,\phi)-\frac{1}{\surd 2}F(1/\surd 2,\phi),\quad 0\le \phi \le \pi/2
\ee
\subsection{The Logarithm}

\cite{MezoArxiv0811}
\be
\int \frac{\log z}{(1-z)z}dz=\Li_2(1-z)+\frac{1}{2}\log^2 z.
\ee

\cite[3.1.6.]{Apelblat2}
\be
\int \frac{\ln x}{1+ax}dx =\frac{1}{a}[\ln x \ln(1+ax)+\Li_2(-ax)]
.
\ee
Correcting a sign error in \cite[3.1.7]{Apelblat2}:
\be
\int \frac{\ln (a+bx)}{c+hx}dx =\frac{1}{h}[\ln\left(\frac{ah-bc}{h}\right) \ln(c+hx)-\Li_2(\frac{bc+bhx}{bc-ah})]
.
\ee

\cite{Brien}
\be
\int \ln(\ln x)^n dx= xln (\ln x)^2-n \li (x),
\ee
where $\li x=\Ei(\ln x)$.

\cite{Brien}
\be
\int x^n \ln(\ln x)^m dx= \frac{1}{n+1}[x^{n+1}\ln(\ln x)^m-m\li(x^{n+1})].
\ee

\cite{Brien}
\be
\int \frac{\ln(\ln x)^m}{x} dx= \ln x[\ln(\ln x)^m-m].
\ee

\cite{Brien}
\be
\int (a+bx)^m\ln^n(c+k\ln x) dx=\frac{1}{(m+1)b}[
(a+bx)^{m+1}\ln^n(c+k x)-n\int \frac{(a+bx)^{m+1}\ln^{n-1}(c+kx)}{c+kx}dx
]
.
\ee

\cite{Brien}
\be
\int \frac{\ln(c+kx)^n}{a+bx} dx=\frac{1}{b}[
\ln(a+bx)\ln^n(c+kx)- nk\int \frac{\ln(a+bx)\ln^{n-1}(c+kx)}{c+kx}dx
]
.
\ee

\cite{Brien}
\be
\int \frac{(\ln x)^n}{(a+bx)^m} dx=\frac{1}{b(m-1)}[
-\frac{(\ln x)^n}{(a+bx)^{m-1}}+n\int \frac{(\ln x)^{n-1}}{x(a+bx)^{m-1}}dx
]
.
\ee

\cite{Brien}
\be
\int \frac{(\ln x)^n}{a+bx} dx=\frac{1}{b}[
(\ln x)^n\ln(a+bx)-n\int \frac{(\ln x)^{n-1}\ln(a+bx)}{x}dx
]
.
\ee

\cite{Brien}
\be
\int \frac{(\ln x)^n}{x^m} dx=\frac{1}{m-1}[
-\frac{(\ln x)^n}{x^{m-1}}+n\int \frac{(\ln x)^{n-1}}{x^m}dx
]
=-\frac{1}{x^{m-1}} \sum_{k=0}^n \frac{n!}{(n-k)!} \frac{(\ln x)^{n-k}}{(m-1)^{k+1}}
.
\ee

\cite{Brien}
\be
\int x^m\ln^n(a+bx)dx=\frac{1}{m+1}[
x^{m+1}\ln(a+bx)^n-nb\int \frac{x^{m+1}\ln^{n-1}(a+bx)}{a+bx}dx
]
.
\ee

\cite{Brien}
\be
\int \frac{dx}{\ln(a+bx)}=\frac{\li(a+bx)}{b}.
\ee

\cite{Brien}
\be
\int \frac{x}{\ln(a+bx)}dx=\frac{1}{b^2}\{\li[(a+bx)^2]-a \li(a+bx)\}
.
\ee

\cite{Brien}
\be
\int \frac{\ln^m(a+\ln x)}{x^n}dx=-\frac{\ln^m(a+\ln x)}{(n-1)x^{n-1}}
+\frac{m}{n-1}\int \frac{\ln^{m-1}(a+\ln x)}{x^n(a+\ln x)}dx
.
\ee

\cite{Brien}
\be
\int \frac{\ln^m(a+\ln x)}{x}dx=(a+\ln x)\ln^m(a+\ln x)
-m\int \frac{\ln^{m-1}(a+\ln x)}{x}dx
.
\ee

\cite{Brien}
\be
\int \frac{\ln(a+\ln x)}{x}dx=(a+\ln x)\ln(a+\ln x)
-\ln x
.
\ee

\cite{Brien}
\be
\int x^n \ln^m(a+\ln x)dx=\frac{1}{n+1}[
x^{n+1}\ln^m(a+\ln x)-m\int \frac{x^n\ln^{m-1}(a+\ln x)}{a+\ln x}dx
]
.
\ee

\cite{Brien}
\be
\int x^n \ln(a-\ln x)dx=\frac{1}{n+1}[
x^{n+1}\ln^m(a-\ln x)-e^{(n+1)a}\Ei[(n+1)(\ln x -a)]
]
.
\ee

\cite{Brien}
\be
\int \frac{\ln(a-\ln x)}{x}dx=(\ln x-a)\ln(a-\ln x)-\ln x.
.
\ee

\cite{Brien}
\be
\int \frac{dx}{(a\mp\ln x)^n}=\pm\frac{x}{(n-1)(a\mp\ln x)^{n-1}}\mp\frac{1}{n-1}\int \frac{dx}{(a\mp\ln x)^{n-1}}
.
\ee

\cite{Brien}
\be
\int \frac{dx}{a\mp\ln x}=\mp e^a \Ei(\ln x \mp a )
.
\ee

\cite{Brien}
\be
\int \frac{dx}{x^n (\ln x)^m}=-\frac{1}{m-1}[\frac{1}{x^n(\ln x)^{m-1}}
+(n-1)\int\frac{dx}{x^n(\ln x)^{m-1}} 
]
.
\ee

\cite{Brien}
\be
\int \frac{dx}{x^n \ln x}=\li(\frac{1}{x^{n-1}})
.
\ee

\cite{Brien}
\be
\int \frac{dx}{x (\ln x)^m}=-\frac{1}{m-1}\frac{1}{(\ln x)^{m-1}}
.
\ee

\cite{Brien}
\be
\int \frac{x^n}{(a+\ln x)^m}dx=-x^{n+1}\sum_{k=0}^{m-2}\frac{(m-k-2)!}{m-1)!} \frac{(n+1)^k}{(a+\ln x)^{m-k-1}}
+e^{-(n+1)a} \frac{(n+1)^{m-1}}{(m-1)!}\Ei[(n+1)(a+\ln x)]
.
\ee

\cite{Brien}
\be
\int \frac{x^n}{a\pm \ln x}dx=\pm 
e^{-(n+1)a} \Ei[(n+1)(\ln x\pm a)]
.
\ee

\cite{Brien}
\be
\int \frac{\ln x}{(a\pm \ln x)^2}dx=\
\frac{ax}{a\pm \ln x}+(1\mp a)e^{\mp a}\Ei(\ln x \pm a)
.
\ee

\cite{Brien}
\be
\int e^{\pm \mu x}\ln^p x dx
=\frac{1}{\mu}[\pm \ln ^p x e^{\pm \mu x}\mp p\int\frac{\ln^{p-1}x e^{\pm \mu x}}{x}dx]; \mu\neq 0.
.
\ee

\cite{ChoiJMAA231}
\be
\int_0^z \frac{\log t}{1+a^2t^2}dt
=\frac{\log z}{a}\arctan az +\frac{\pi}{2a}\log\frac{\cos(\arctan az)}{\pi}
+\frac{\pi}{a}
\log\frac
{G(1+\frac{1}{\pi}\arctan az)G(3/2-\frac{1}{\pi}\arctan az)}
{G(1-\frac{1}{\pi}\arctan az)G(1/2+\frac{1}{\pi}\arctan az)}
,
\ee
where $G$ is the reciprocal of the double Gamma function.

\cite{ChoiJMAA231}
\be
\int_0^z \frac{\log t}{\sqrt{1-a^2t^2}}dt =
\frac{1}{a}\arcsin az \log\frac{z}{2\pi}
+\frac{\pi}{a}
\log\frac{G(1+\frac{1}{\pi}\arcsin az)}
{G(1-\frac{1}{\pi}\arcsin az)}
.
\ee

\be
\int \frac{\ln x}{\sqrt{1-x^2}}dx = (\ln x -1)\arcsin x + \frac{x^3}{9}{}_3F_2(3/2,3/2,3/2; 5/2,5/2;x^2),
\ee
where the hypergeometric term can be reduced via \eqref{eq.3F2Plus}.

\cite[2.727.5]{GR}
\be
\int \frac{\ln x}{\sqrt{a+bx}}dx=
\left\{\begin{array}{ll}
\frac{2}{b}\left[(\ln x-2)\sqrt{a+bx}+\sqrt{a}\ln\frac{\sqrt{a+bx}+\sqrt{a}}{\sqrt{a+bx}-\sqrt{a}}\right],& a>0 ;\\
\frac{2}{b}\left[(\ln x-2)\sqrt{a+bx}+2\sqrt{-a}\arctan\sqrt{\frac{a+bx}{-a}}\right],& a<0 .
\end{array}
\right.
\ee

\be
\int du \ln u
\arctan\frac{c}{\sqrt{u-c^2}}
= 
(\ln u-3)
c\sqrt{u-c^2}
+u(\ln u-1)
\arctan\frac{c}{\sqrt{u-c^2}}
+2c^2\arctan\frac{\sqrt{u-c^2}}{c}.
\ee

Recursively
\begin{multline}
\int (a+bx)^{k-1/2}\ln x dx
\\
=
\frac{k-1/2}{(k+1/2)^2b}
(a+bx)^{k+1/2}
+\frac{1}{k+1/2}[x(\ln x -1)-a/b](a+bx)^{k-1/2} 
+
\frac{(k-1/2)a}{k+1/2}\int
\ln x (a+bx)^{k-3/2}
.
\end{multline}

\begin{multline}
\int z \ln(c+z+1/z)dz =
\frac{z^2-c^2/2+1}{2}\ln(z^2+cz+1)
-\frac{c\sqrt{c^2/4-1}}{2}
\ln\frac{z+c/2+\sqrt{c^2/4-1}}{z+c/2-\sqrt{c^2/4-1}}
\\
-\frac{(z-c)^2}{4}
-\frac{1}{2}z^2\ln z
.
\end{multline}

\cite{CoffeyJMP49}
\be
\kappa \int_0^u \ln(\sin \kappa x+\sin \alpha)dx =
\Cl_2(\alpha)-\Cl_2(\kappa u+\alpha)
+\Cl_2(\alpha-\kappa u +\pi)-\Cl_2(\alpha+\pi)-\kappa u\ln 2.
\ee
for $\kappa >0$, $|u|\le |\alpha|$.

\cite{CoffeyJMP49}
\be
\int_0^x \ln|\cos A-\cos kt|dt =
-\frac{1}{k}[\Cl_2(kx-A)+\Cl_2(kx+A)+kx\ln 2].
\ee

\cite{CoffeyJMP49}
\be
\int_0^u \ln(\sin^2x -\sin^2\alpha)dx =
\frac12[\Cl_2(2\alpha-2u)-\Cl_2(2\alpha+2u)]-2u\ln 2.
\ee

\cite{CoffeyJMP49}
\be
\int_b^\infty \ln\frac{u+a}{u-a}\frac{du}{1+u^2}=\Cl_2(\pi-\theta)-\frac12[\Cl_2(2\omega)+\Cl_2(2\chi)],
\ee
where
$
\theta=\arccos\frac{1-a^2}{1+a^2}
$
,
$r=(b+a)/(b-a)$, $\tan\omega=r\sin\theta/(1-r\cos\theta)$, $\chi=\pi-\theta-\omega$.

\cite{WatrasieOA14}
\begin{multline}
\int \log^2(\frac{b+x}{b-x})dx
=(b+x)\log^2(b+x)-(b-x)\log^2(b-x)
\\
-2(b+x)\log(b+x)\log(b-x)+4\log(2b)\log(b-x)
-4b\sum_{n=1}^\infty (\frac{b-x}{2b})^n \frac{1}{n^2},\quad |x|<|b|.
\end{multline}
\cite{WatrasieOA14}
\begin{multline}
\int_{-a}^a \log^2(\frac{b+x}{b-x})dx
=2(b+a)\log^2(b+a)-2(b-a)\log^2(b-a)
-4(b+a)\log(b+a)\log(b-a)-4b\log^2b
\\
-8b\log (2)\log (b)
+8b\log(2b)\log(b-a)
+8b[\frac{\pi^2}{12}-\frac12\log^2 2]
-8b\sum_{n=1}^\infty [\frac12(1-a/b)]^n \frac{1}{n^2},\quad a<b.
\end{multline}

\cite{ChoiJMAA231}
\be
\int_0^z \log \sin \pi t dt = z\log\frac{\sin \pi z}{2\pi}+\log\frac{G(1+z)}{G(1-z)}.
\ee
where $G$ is the reciprocal of the double Gamma Function.

\cite{ChoiJMAA231}
\be
\int_0^z \log \sin a t dt = z\log\frac{\sin a z}{2\pi}
+\frac{\pi}{a}\log\frac{G(1+az/\pi)}{G(1-az/\pi)}.
\ee

\cite{ChoiJMAA231}
\be
\int_0^z \log \cos \pi t dt = 
(z-\frac12)\log\frac{\cos \pi z}{2\pi}-\frac12\log 2
-\log\Gamma(\frac12-z)
+\log\frac{G(1/2+z)}{G(1/2-z)}.
\ee

\cite{ChoiJMAA231}
\be
\int_0^z \log \tan \pi t dt = 
z\log \tan \pi z
+\frac12 \log\frac{\cos \pi z}{\pi}
+\log\Gamma(\frac12-z)
+\log\frac{G(1+z)}{G(1-z)}
-\log\frac{G(1/2+z)}{G(1/2-z)}.
\ee

\cite{ChoiJMAA231}
\be
\int_0^z
\log(2+2\cos t)dt = (2z-2\pi)\log\frac{\cos(z/2)}{\pi}
+4\pi\log\frac{G(1/2+z/2\pi)}{G(3/2-z/2\pi)}.
\ee

\section{Definite Integrals of Elementary Functions I.}

\subsection{General formulae}

\cite{GlasserMCom40}
\be
\int_{-\infty}^\infty F(u)dx=\int_{-\infty}^\infty F(x)dx
\ee
where $u=x-\sum_{j=1}^{n-2} a_j/(x-C_j)$ where $a_j$ is any sequence
of positive constants and the $C_j$ are any real constants.

\cite{Amdeberhanarxiv2445}
\be
\int_0^\infty f([ax-b/x]^2) dx = \frac{1}{a}\int_0^\infty f(y^2)dy,\qquad a,b>0.
\ee

\cite{BorosSci12}
Let
\[
\varphi(x)=\sum_{k=0}^\infty A_kx^k
,
\]
then
\be
\int_0^\infty x^{\beta-1}e^{-x}\varphi(x)dx
=
\sum_{k=0}^\infty A_k\Gamma(k+\beta)
,\quad \beta>0
.
\ee
Let
\[
\varphi(x)=\sum_{k=0}^\infty A_kx^{k/p}
,
\]
then
\be
\int_0^\infty x^{\beta-1}e^{-x}\varphi(x)dx
=
\sum_{k=0}^\infty A_k\Gamma(k/p+\beta)
,\quad \beta>0
.
\ee
Let
\[
\varphi(x)=\sum_{k=0}^\infty A_kx^{k}
,
\]
then
\be
p \int_0^\infty x^{\beta-1}e^{-x^p}\varphi(x)dx
=
\sum_{k=0}^\infty A_k\Gamma(\frac{\beta+k}{p})
,\quad \beta>0
.
\ee

\cite{GlasserArxiv1403}
\be
\int_0^\infty x^{s-1} F(x) dx = \Gamma(s)\phi(s).
\ee
where $\phi$ is the analytic continuation of the Taylor coefficient of $F$.

\cite{GlasserArxiv1403}
\be
\int_{-\infty}^\infty \frac{F(iw/u)}{(1+u^2}du = F(w).
\ee

\cite{GlasserArxiv1403}
\be
\int_{-\infty}^\infty \frac{F(ax(x+i))}{\cosh(\pi x)}dx = F(a/4).
\ee

\cite{GlasserArxiv1403}
\be
\int_{-\pi}^\pi \frac{F(e^{ix}\cos x)}{a^2\cos^2(x)+b^2\sin^2 x}dx = \frac{2\pi}{ab}F(b/(a+b)).
\ee

\cite{GlasserArxiv1403}
\be
\int_{-\infty}^\infty \frac{F(ax(x+i))}{\sinh(\pi p(2x+i))}dx = -i\frac{F(a/4)}{2p}.
\ee

\subsection{Power and Algebraic Functions}

\begin{multline}
\int_0^1 dx \frac{x^{k+1}}{(1+\delta^2+\beta^2x^2)^{k+2}}
=
\frac{1}{(1+\delta^2)^{k+2}(k+2)}{}_2F_1(2+k, 1+k/2; 2+k/2; -\frac{\beta^2}{1+\delta^2})
\\
=
\frac{1}{(1+\delta^2)(k+2)(1+\delta^2+\beta^2)^{k+1}}
{}_2F_1(1, -k/2; 2+k/2; -\frac{\beta^2}{1+\delta^2}).
\end{multline}
For even $k$, this is a terminating hypergeometric function, a polynomial
of order $k/2$ in the argument. For odd $k$, this is rewritten as \cite[15.3.26]{AS}
\be
_2F_1(1,-k/2;2+k/2;-\frac{\beta^2}{1+\delta^2})
= \frac{1+\delta^2}{1+\delta^2-\beta^2}
{}_2F_1(1,1/2;2+k/2;-\frac{4\beta^2(1+\delta^2)}{(1+\delta^2-\beta^2)^2})
\ee
and for increasing odd $k$ recursively computed with
the contiguous relations \cite[15.2.27]{AS} anchored at \cite[15.1.5]{AS}
\be
{}_2F_1(1,1/2;3/2;-z^2) = \frac{1}{z}\arctan z;
\quad
{}_2F_1(1,1/2;1/2;-z^2) = \frac{1}{1+z^2}.
\ee

\cite{ReynoldsAIMS6}
\begin{multline}
\int_0^\infty \frac{x^{m-1}}{(x^{2u}+1)(x^{3u}+1)}dx
=\frac{\pi \csc(\pi m/u)}{6u}[
-2\surd 3\sin\frac{2\pi m}{3u}
+3\surd 2\sin\frac{\pi(2m+u)}{4u}
+2\cos\frac{2\pi m}{3u}+1
]
\end{multline}

\cite{LambertMathComp66}
\be
\int_a^u\frac{xdx}{\sqrt{(x-a)(x-b)(x-c)}}
=\frac{2}{\sqrt{a-c}}[(a-b)\Pi(\mu,1,q)+bF(\mu,q)],\, u>a>b>c.
\ee

\cite{LambertMathComp66}
\begin{multline}
\int_u^c\frac{dx}{(r-x)\sqrt{(a-x)(b-x)(c-x)}}
=\frac{2(c-b)}{(r-b)(r-c)\sqrt{a-c}}
\\
\times
\Pi(\beta,\frac{r-b}{r-c},p)+\frac{2}{(r-b)\sqrt{a-c}}F(\beta,p),\, a>b>c>u,r\neq c
.
\end{multline}

\cite{LambertMathComp66}
\begin{multline}
\int_a^u\frac{dx}{(x-r)\sqrt{(x-a)(x-b)(x-c)}}
=\frac{2}{(b-r)(a-r)\sqrt{a-c}}
\\
\times
\left[
(b-a)\Pi(\mu,\frac{b-r}{a-b},q)+(a-r)F(\mu,q)
\right]
,\, u>a>b>c,r\neq a
.
\end{multline}

\cite{LambertMathComp66}
\begin{multline}
\int_u^b\sqrt{\frac{(x-c)(b-x)}{a-x}}dx
=\frac{2}{3}\sqrt{a-c}[2(b-a)F(\delta,q)
\\
+
(2a-b-c)E(\delta,q)]
+\frac{2}{3}(b+c-a-u)\sqrt{\frac{(b-u)(u-c)}{a-u}}
,\, a>b>u\ge c
.
\end{multline}

\cite{LambertMathComp66}
\begin{multline}
\int_a^u\sqrt{\frac{(x-b)(x-c)}{x-a}}dx
=\frac{2}{3}\sqrt{a-c}[2(a-b)F(\mu,q)
\\
+
(b+c-2a)E(\mu,q)]
+\frac{2}{3}(u+2a-2b-c)\sqrt{\frac{(u-a)(u-c)}{u-b}}
,\, u>a>b>c
.
\end{multline}

\cite{LambertMathComp66}
\begin{multline}
\int_a^u\sqrt{\frac{(x-a)(x-c)}{x-b}}dx
=\frac{2}{3}\sqrt{a-c}[(a+c-2b)E(\mu,q)
\\
-
(a-b)F(\mu,q)]
+\frac{2}{3}(u+b-a-c)\sqrt{\frac{(u-a)(u-c)}{u-b}}
,\, u>a>b>c
.
\end{multline}

\cite{CruzPRA17}
\be
\int_0^\infty \frac{u^\alpha du}{(au+1)^\beta(bu+1)^\gamma}
=A_0\log a+B_0\log b
+ \sum_{r=1}^{\beta-1}\frac{A_r}{r}
+ \sum_{s=1}^{\gamma-1}\frac{B_s}{s}
;
\ee
where
\be
A_r=(-1)^{\alpha+\beta+r+1}a^{\gamma-\alpha-1}
\sum_{j=0}^{\beta-\gamma-1}
\binom{\alpha}{j}\binom{\beta+\gamma-r-j-2}{\gamma-1}
\frac{b^{\beta-r-j-1}}{(a-b)^{\beta+\gamma-r-j-1}}
;
\ee
\be
B_s=(-1)^{\beta}b^{\beta-\alpha-1}
\sum_{j=0}^{\gamma-s-1}
\binom{\alpha}{j}\binom{\beta+\gamma-s-j-2}{\beta-1}
\frac{a^{\gamma-s-j-1}}{(a-b)^{\beta+\gamma-s-j-1}}
.
\ee
\be
r=0,1,\ldots,\beta-1;\quad s=0,1,\ldots,\gamma-1.
\ee

\subsection{Powers of $x$, of binomials of the form $\alpha+\beta x^p$, and of polynomials in $x$}

\cite[(C6.2)]{Nahin}
\be
\int_0^\infty \frac{x}{x^n+1}dx = \frac{\pi/n}{\sin\frac{2\pi}{n}}.
\ee

\cite[(C8.11)]{Nahin}
\be
\int_0^\infty \frac{x^m}{x^n+b}dx = \frac{\pi}{nb^{(n-m-1)/n}\sin\frac{(m+1)\pi}{n}}.
\ee

\cite{DriverETNA25}
\be
_3F_2(-n,b/2,(b+1)/2;c/2,(c+1)/2;x)
=\frac{\Gamma(c)}{\Gamma(b)\Gamma(c-b)}
\int_0^1 t^{b-1}(1-t)^{c-b-1}(1-xt^2)^ndt,\quad \Re c>\Re b>0.
\ee

\cite{BorosSIR40}
\begin{multline}
\int_0^\infty\left[\frac{x^2}{x^4+2ax^2+1}\right]^r\cdot \frac{x^2+1}{x^b+1}\frac{dx}{x^2}
\\
=
\int_0^\infty\left[\frac{x^2}{x^4+2ax^2+1}\right]^r\frac{dx}{x^2}
=
\int_0^\infty\left[\frac{x^2}{x^4+2ax^2+1}\right]^r dx
\\
=
\int_0^\infty\left[\frac{x^2}{x^4+2ax^2+1}\right]^r \frac{x^2+1}{x^2} dx
=2^{-1/2-r}(1+a)^{1/2-r} B(r-1/2,1/2),
\end{multline}
with $B$ Euler's beta function,
$a>=1, r>1/2$, any $b$.

From this by specialization \cite{BorosSIR40} 
\be
\int_0^\infty \frac{x^4}{(x^4+x^2+1)^3 }dx=
\frac{\pi}{48\surd 3}.
\ee

\cite{BorosSIR40}
\be
\int_0^\infty \frac{x^3}{(x^4+7x^2+1)^{5/2} }dx=
\frac{2}{243}.
\ee

\cite{BorosSIR40}
\be
\int_0^\infty \frac{\surd x}{(x^4+14x^2+1)^{5/4} }dx=
\frac{\Gamma^2(3/4)}{4\sqrt{2\pi}}.
\ee

\cite{BorosSIR40}
\be
\int_0^\infty\left[\frac{x^2}{bx^4+2ax^2+c}\right]^r dx
=
\frac{B(r-1/2,1/2)}{2^{r+1/2}\sqrt{b}[a+\sqrt{bc}]^{r-1/2}}
\ee
with $b>0$, $c\ge 0$, $a>-\sqrt{bc}$ and $r>1/2$.

\cite{BorosSIR40}
\be
\int_0^\infty\left[\frac{x^2}{x^4-x^2+1}\right]^r\frac{ dx}{x^2}
=
\frac{\surd\pi}{2}\,\frac{\Gamma(r-1/2)}{\Gamma(r)}
\ee

Special cases of \eqref{eq.x4denom} \cite[(2.3.5)]{Nahin}
\be
\int_0^\infty \frac{dx}{x^4+x^2+1}=\frac{\pi}{2\surd 3}.
\ee
\be
\int_0^\infty \frac{dx}{x^4-x^2+1}=\frac{\pi}{2}.
\ee
\be
\int_0^\infty \frac{dx}{x^4+2x^2+1}=\frac{\pi}{4}.
\ee

\cite{BorosMC71}
\be
\int_0^\infty \frac{z^{2n}dz}{(z^4+2d_1z^2+1)^{m+1}}
=\frac{\pi}{2^{3m+3/2}(1+d_1)^{m+1/2}}\sum_{j=0}^{m-n}2^j(1+d_1)^j
\binom{2m-2j-1}{m-j}\binom{m-n+j}{2j}\binom{2j}{j}\binom{m}{j}^{-1}
\ee
for $0\le n\le m$.
and the same expression with upper summation limit $n-m-1$
for $m+1\le n\le 2m+1$.

\cite{BorosMC71}
\begin{multline}
\int_0^\infty \frac{z^{2n}dz}{(bz^4+2az^2+c)^{m+1}}
\\
=
\pi \left\{c(c/b)^{m-h}[8(a+\sqrt{bc})]^{2m+1}\right\}^{-1/2}
\sum_{k=0}^{m-n}
2^k
\binom{2m-2k}{m-k}\binom{m-n+k}{2k}\binom{2k}{k}\binom{m}{k}^{-1}
(\frac{a}{\sqrt{bc}}+1)^k
\end{multline}
for $0\le n\le m$, 
and the same expression with upper summation limit $n-m-1$
for $m+1\le n\le 2m+1$.

\cite{BorosMC71}
\be
\int_0^\infty \frac{z^{2n}dz}{(z^8+a_2z^6+2a_1z^4+a_2z^2+1)^{m+1}}
=\sum_{k=n}^{2m+1}\sum_{j=0}^{k-m-1} t_{k,j}(m,n;a_1,a_2),
\ee
for $m+1\le n\le 2m+1$, $1+a_1+a_2>0$, $a_2+4>-8\sqrt{8(1+a_1+a_2)}$, and 
\be
\ldots 
=\sum_{k=n}^{m}\sum_{j=0}^{m-k} t_{k,j}(m,n;a_1,a_2)
+
\sum_{k=m+1}^{2m+1}\sum_{j=0}^{k-m-1} t_{k,j}(m,n;a_1,a_2)
\ee
for $0\le n\le m$, where $c_1\equiv a_2+4$, $c_2\equiv 1+a_1+a_2$, and
\begin{multline}
t_{k,j}(m,n;a_1,a_2)\equiv \pi 2^{-(3m+2+k+j)/2}
c_2^{(m-k-j)/2}
(c_1+\sqrt{8c_2})^{j-m-1/2}
\\ \times
\binom{4m-n-k+2}{k-n}
\binom{2m-2j}{m-j}
\binom{m-k+j}{2j}
\binom{2j}{j}
\binom{m}{j}^{-1}.
\end{multline}
The denominator polynomials of the previous formulas are symmetric.
There are transformations for even integrals of rational functions $P(x)/Q(x)$
which preserve the values of $\int_0^\infty P(x)/Q(x) dx$, explicit
in \cite{BorosMC71} and transform the set of numerator and denominator
coefficients to $\int P^+(x)/Q^+(x)dx$ with symmetric polynomials.

\cite[2.6.9]{Nahin}
\be
\int_0^\infty \frac{dx}{bx^4+2ax^2+1}=\frac{\pi}{2\surd 2} \,\frac{1}{\sqrt{a+\surd b}}.
\ee

\cite{Amdeberhanarxiv2118,MollNAMS49}
\be
\int_0^\infty \frac{dx}{(x^4+2ax^2+1)^{m+1}}=\frac{\pi}{2}\frac{\sum_{l=0}^md_l(m)a^l}{[2(a+1)]^{m+1/2}}
\ee
where
\be
d_l(m)\equiv 2^{-2m}\sum_{k=l}^m2^k\binom{2m-2k}{m-k}\binom{m+k}{m}\binom{k}{l}.
\ee

\cite{MollNAMS49}
\be
\int_0^\infty \frac{dx}{bx^4+2ax^2+1}
= \frac{\pi}{2\sqrt{2}}\,\frac{1}{\sqrt{a+\surd b}}.
\ee

\cite{LucasAMM}
\be
\frac{(-)^n}{4^{n-1}}\int_0^1
\frac{x^{4n}(1-x)^{4n}}{1+x^2}dx
=\pi-\sum_{k=0}^{n-1}
(-)^k\frac{2^{4-2k}(4k)!(4k+3)!}{(8k+7)!}(820k^3+1533k^2+902k+165).
\ee

\cite{LucasAMM}
\be
\int_0^1
\frac{x^{m}(1-x)^{n}}{1+x^2}dx
=\frac{\surd \pi \Gamma(m+1)\Gamma(n+1)}{2^{m+n+1}}
\,_3F_2\left(\begin{array}{c}1,(m+1)/2,(m+2)/2
\\
(m+n+2)/2,(m+n+3)/2
\end{array}\mid -1\right).
\ee

\cite{BackhouseMG79}
\be
\int_0^1
\frac{x^m(1-x)^n}{1+x^2}dx
= R_{m,n}+A_{m,n}\pi +B_{m,n}\ln \surd 2
\ee
induced by the partial fraction decomposition
\be
\frac{x^m(1-x)^n}{1+x^2}=Q_{m,n}(x)+\frac{A_{m,n}}{1+x^2}
+\frac{B_{m,n}x}{1+x^2}
.
\ee

\cite{Scarpelloarxiv1212}\cite[A093341]{sloane}
\be
\int_0^\infty \frac{dx}{\sqrt{1+x^4}} = \mathbf{K}(\frac{1}{\surd 2}).
\ee

\cite{Scarpelloarxiv1212}
\be
\int_0^\infty \frac{dx}{\sqrt{x+x^4}} = \mathbf{K}(\frac{\surd 6-\surd 2}{4}).
\ee

\cite{Scarpelloarxiv1212}
\be
\int_0^1 \frac{dx}{\sqrt[4]{1-x^2}} = \sqrt{3}\left[2\mathbf{E}(\frac{1}{\surd 2}-\mathbf{K}(\frac{1}{\surd 2})\right].
\ee

\cite{Scarpelloarxiv1212}\cite[A062539]{sloane}
\be
\int_0^1 \frac{dx}{(1-x^2)^{3/4}} = \sqrt{2}\mathbf{K}(\frac{1}{\surd 2}).
\ee

\cite{Scarpelloarxiv1212}
\be
\int_0^1 \frac{x^2}{(1-x^2)^{3/4}}dx = \frac{2^{3/2}}{3}\mathbf{K}(\frac{1}{\surd 2}).
\ee

\cite{Scarpelloarxiv1212}
\be
\int_0^1 \frac{1}{\sqrt{1-x^6}}dx = \frac{1}{\sqrt[4]{3}}\mathbf{K}(\frac{\surd 6-\surd 2}{4}).
\ee

\cite{Scarpelloarxiv1212}
\be
\int_0^\infty \frac{1}{\sqrt{1+x^6}}dx = \frac{2}{\sqrt[4]{27}}\mathbf{K}(\frac{\surd 6-\surd 2}{4}).
\ee

\cite{Scarpelloarxiv1212}
\be
\int_0^1 \frac{1}{\sqrt{1-x^8}}dx = \frac{1}{\sqrt{2}}\mathbf{K}(\surd 2 -1).
\ee

\cite{Scarpelloarxiv1212}
\be
\int_0^1 \frac{x^2}{\sqrt{1-x^8}}dx = (1-\frac{1}{\sqrt{2}})\mathbf{K}(\surd 2 -1).
\ee

\cite{Scarpelloarxiv1212}
\be
\int_0^1 \frac{x^{a-1}}{\sqrt{1-x^n}}dx = \cos(a\pi /n)\int_0^\infty \frac{z^{a-1}}{\sqrt{1+z^n}}dz,\quad 2a<n.
\ee

\cite{Scarpelloarxiv1212}
\be
\int_0^\infty \frac{z^{n-a-1}}{\sqrt{1+z^n}}dz
\cdot 
\int_0^1 \frac{x^{a-1}}{\sqrt{1-x^n}}dx
= \frac{2\pi}{n(2a-n)\sin(\pi a/n)},\quad n/2<a<n.
\ee

\cite{Scarpelloarxiv1212}
\be
\int_0^\infty \frac{1}{\sqrt{1+x^8}}dx 
= \int_0^\infty \frac{x^2}{\sqrt{1+x^8}}dx 
= \sqrt{2-\surd 2}\mathbf{K}(\surd 2-1).
\ee

\cite{Scarpelloarxiv1212}
\be
\int_0^1 \frac{x^4}{\sqrt{1-x^8}}dx = \frac{\pi}{8}\frac{\surd 2}{\mathbf{K}(\surd 2-1)}.
\ee

\cite{Scarpelloarxiv1212}
\be
\int_0^1 \frac{x^6}{\sqrt{1-x^8}}dx = \frac{\pi}{24}\frac{2+\surd 2}{\mathbf{K}(\surd 2-1)}.
\ee

\cite{Scarpelloarxiv1212}
\be
\int_0^\infty \frac{1}{\sqrt[3]{1+x^6}}dx = \frac{\sqrt[3]{4}}{\sqrt[4]{3}}\mathbf{K}(\frac{\surd 6-\surd 2}{4}).
\ee

\cite{Scarpelloarxiv1212}
\be
\int_0^1 \frac{1}{\sqrt[3]{1-x^6}}dx = \frac{\sqrt[3]{4}}{\sqrt[4]{27}}\mathbf{K}(\frac{\surd 6-\surd 2}{4}).
\ee

\cite{Scarpelloarxiv1212}
\be
\int_1^\infty \frac{1}{\sqrt[m]{(x^n+a)(x^n+b)}}dx = \frac{m}{2n-m}
F_1\left(\begin{array}{c}\frac{2n-m}{mn};1/m,1/n\\ \frac{2n-m+mn}{mn}\end{array}
\mid -a,-b\right),\quad 2n-m>0,\quad a,b>0.
\ee

\cite{Scarpelloarxiv1212}
\be
\int_0^1 \frac{1}{\sqrt[m]{(x^n+a)(x^n+b)}}dx = \frac{1}{\sqrt[m]{ab}}
F_1\left(\begin{array}{c}1/n;1/m,1/m\\ 1+1/n\end{array}
\mid -1/a,-1/b\right),\quad a,b>0.
\ee

\subsection{The Exponential function}

\cite{Amdeberhanarxiv3663}
\be
-\int_0^2\frac{x^p-x^q}{1-x}dx=\psi(p+1)-\psi(q+1)
\ee

\cite{KolbigMathComp64_449}
\be
\int_u^\infty \exp(-\frac{x^2}{4\beta}-\gamma x)dx=\sqrt{\pi\beta}e^{\beta\gamma^2}
\left[1-\Phi(\gamma\sqrt{\beta}+\frac{u}{2\sqrt\beta})\right]
,\, \Re\beta>0,u\ge 0
.
\ee

\cite{KolbigMathComp64_449}
\be
\int_{-\infty}^\infty \exp(-p^2x^2\pm qx)dx=\exp(\frac{q^2}{4p^2})
\frac{\sqrt{\pi}}{|p|}
.
\ee

\cite{KolbigMathComp64_449}
\be
\int_0^\infty\frac{x^ne^{-\mu x}}{x+\beta}dx
=(-1)^{n-1}\beta^ne^{\beta\mu}\Ei(-\beta\mu)+\sum_{k=1}^n(k-1)!(-\beta)^{n-k}\mu^{-k},\
\,
|\arg \beta|<\pi,\,\Re\mu>0,\, n\ge 0
.
\ee

\cite{GeddesCM}
\be
\int_0^\infty \frac{x^n\exp(-ux)}{1-\exp{-vx}}dx=\frac{(-)^{n+1}}{v^{n+1}}\psi^{(n)}(u/v).
\ee

\cite{GeddesCM}
\be
\int_0^\infty \frac{t^w\exp(-ut^s)}{1-\exp{-vt^s}}dt=\frac{(-)^m}{|s|v^m}\psi^{(m-1)}(u/v).
\ee

\cite[(3.6)]{MuthumSci22}
\be
\int_a^\infty \frac{e^{-t(n+1/2)}}{1-e^{-t}} dt = -\sum_{k=0}^{n-1} \frac{e^{-a(k+1/2)}}{k+1/2}+\log \coth\frac{a}{4}.
\ee

\cite[(3.7)]{MuthumSci22}
\be
\int_a^\infty \frac{e^{-tz}}{1-e^{-t}} dt = \sum_{k=0}^\infty \frac{e^{-a(z+k)}}{z+k}.
\ee

\cite[(5.4.5)]{Nahin}
\be
\int_0^\infty (e^{-\alpha e^x}+e^{-\alpha e^{-x}}-1)dx = -\gamma-\ln \alpha.
\ee

\cite[(5.4.6)]{Nahin}
\be
\int_0^\infty (e^{-x^a}-e^{-x^b})\frac{dx}{x} = \gamma\frac{a-b}{ab}.
\ee

\cite{CoppoJNT150}
\be
\frac{n\binom{2n}{n}}{2^{2n-1}}
\int_0^\infty e^{-t}(1-e^{-2t})^{n-1} \frac{t^m}{m!} dt
=
P_m(O_n,\ldots O_n^{(m)}),
\ee
where $P$ are the generalized Bell polynomials and $O$ the
sums of odd number inverses as in \eqref{eq.Bell}.

\subsection{Rational functions of powers and exponentials}
\cite{KolbigMathComp64_449}
\be
\int_0^\infty \frac{x^{\nu-1}e^{-\mu x}}{1-\beta e^{-x}}dx=\Gamma(\nu)
\sum_{n=0}^\infty (\mu+n)^{-\nu} \beta^\mu
.
\ee

\cite{KolbigMathComp64_449}
\be
\int_0^\infty \frac{(1+ix)^{2n-1}-(1-ix)^{2n-1}}{i}
\frac{dx}{e^{\pi x}+1} = \frac{1}{2n}[1-2^{2n}B_{2n}]
.
\ee

\cite{KolbigMathComp64_449}
\be
\int_0^\infty \frac{x^q e^{-px}dx}{(1-ae^{-px})^2}
=\frac{\Gamma(q+1)}{ap^{q+1}}\sum_{k=1}^\infty \frac{a^k}{k^q},\,
-1\le a<1,\,q>-1,\,p>0
.
\ee

\cite{KolbigMathComp64_449}
\be
\int_0^\infty \frac{(1+a)e^x+a}{(1+e^x)^2}e^{-ax}x^n dx=
-n!\sum_{k=1}^\infty \frac{(-1)^k}{(a+k)^n}, 
\, a>-1,\, n=1,2,\ldots
\ee

By repeated differentiation of \cite[3.415.1]{GR} w.r.t. $b$
\cite{MatharArxiv2308}
\be
\int_0^\infty \frac{xdx}{(x^2+b^2)^2 (e^{\mu x}-1)}
= 
-\frac{1}{4b^2}
-\frac{\pi}{4b^3\mu}
+\frac{\mu}{8b\pi}\psi'\left(\frac{b\mu}{2\pi}\right)
.
\ee

\be
\int_0^\infty \frac{xdx}{(x^2+b^2)^3 (e^{\mu x}-1)}
= 
-\frac{1}{8b^4}
-\frac{3\pi}{16b^5\mu}
+\frac{\mu}{32b^3\pi}\psi'\left(\frac{b\mu}{2\pi}\right)
-\frac{\mu^2}{64b^2\pi^2}\psi''\left(\frac{b\mu}{2\pi}\right)
.
\ee

\be
\int_0^\infty \frac{xdx}{(x^2+b^2)^4 (e^{\mu x}-1)}
=
-\frac{1}{12b^6}
-\frac{5\pi}{32b^7\mu}
+\frac{\mu}{64b^5\pi}\psi'\left(\frac{b\mu}{2\pi}\right)
-\frac{\mu^2}{128b^4\pi^2}\psi''\left(\frac{b\mu}{2\pi}\right)
+\frac{\mu^3}{768b^3\pi^3}\psi'''\left(\frac{b\mu}{2\pi}\right)
.
\ee

The variant with $+1$ in the denominator is reduced to the variant with $-1$ for arguments $\mu$ and $2\mu$:
\be
\int_0^\infty \frac{xdx}{(x^2+b^2)^m (e^{\mu x}+1)}
=
\int_0^\infty \frac{xdx}{(x^2+b^2)^m (e^{\mu x}-1)}
-2\int_0^\infty \frac{xdx}{(x^2+b^2)^m (e^{2\mu x}-1)}
.
\ee

\cite{BorosSci12}
\be
p\int_0^\infty x^{\beta -1}e^{ax-x^p}dx 
=
\sum_{k=0}^\infty \frac{a^k}{k!}\Gamma(\frac{k+\beta}{p})
.
\ee
\be
p^2\int_0^\infty t^{p -1}e^{at-t^p}dt 
=
\sum_{k=1}^\infty \frac{a^k}{(k-1)!}\Gamma(k/p)
.
\ee
\be
\int_0^\infty x^{\beta -1}e^{-x-x^p}dx 
=
\sum_{j=0}^\infty \frac{(-1)^j}{j!}\Gamma(\beta+jp)
.
\ee

\cite{GlennMaComp20}\cite[A002161]{sloane}
\be
\int_0^\infty (1-e^{-1/x^2})dx =\surd \pi.
\ee

\cite{KolbigMathComp64_449}\cite[A155739]{sloane}
\be
\int_0^\infty \left\{e^{-x^2}-\frac{1}{1+x^{2^n}}\right\}
\frac{dx}{x}=-\frac{1}{2}C \approx -0.288607
.
\ee

\cite{vHaeringenMathComp39}
\be
\int_0^\infty \exp(nx-\beta \sh x)dx=\frac{1}{2}
\left[S_n(\beta)-\pi\mathbf{E}_n(\beta)-\pi N_n(\beta)\right]
.
\ee

\cite{vHaeringenMathComp39}
\be
\int_{-\infty}^\infty\frac{\exp(\nu \Arsh x-iax)}{\sqrt{1+x^2}}=
\left\{
\begin{array}{ll}
2\exp\left(-\frac{i\nu\pi}{2}\right)K_\nu(a) & \mathrm{if}\, a>0,\\
2\exp\left(\frac{i\nu\pi}{2}\right)K_\nu(-a) & \mathrm{if}\, a<0.\\
\end{array}
\right.
\ee
$[|\Re \nu|<1]$

\subsection{Hyperbolic Functions}
\cite{Adamchik,BradleyCS2001}
\be
\frac{1}{2}\int_0^\infty\frac{x}{\cosh x}dx=G.
\ee

\cite{ApelblatJAM27}
\be
\int_0^\infty \frac{du}{(b^2+u^2)\sinh(au)} =
\frac{1}{2b}\left[\psi(\frac{ab}{2\pi}+\frac{3}{4})-\psi(\frac{ab}{2\pi}+\frac{1}{4})\right]
\ee
for $\Re a>0$, $\Re b>\max(-\Re a,-\Re 3a)$.

\cite{ApelblatJAM27}
\be
\int_0^\infty \frac{u du}{(b^2+u^2)\sinh(au)} =
\frac{1}{2}\left[\psi(\frac{ab}{2\pi}+\frac{1}{2})-\psi(\frac{ab}{2\pi})\right]-\frac{\pi^2}{4a^2b}
\ee
for $\Re a>0$, $\Re b>0$.

\cite[p. 203]{Nahin}
\be
\int_0^\infty y^{s-1}[1-\tanh(ay)] ady = \frac{2^{1-s}\Gamma(x)\zeta(s) [1-2^{1-s}]}{a^s}.
\ee

\cite[(5.3.17)]{Nahin}
\be
\int_0^\infty \frac{x^s}{\sinh^2 x}dx = 2^{1-s}\Gamma(s+1)\zeta(s).
\ee

\cite{BorosSci12}
\be
\int_0^\infty x^{\beta-1}e^{-x}\frac{\sinh b\sqrt{x}}{b\sqrt{x}}dx
=
\sum_{j=0}^\infty \frac{\Gamma(j+\beta)}{(2j+1)!} b^{2j}
.
\ee
\be
\int_0^\infty t^{2\beta-2}e^{-t^2/b^2}\sinh t dt
=
\frac{1}{2}\sum_{j=0}^\infty \frac{\Gamma(j+\beta)}{(2j+1)!} b^{2(j+\beta)}
.
\ee
\be
\int_0^\infty x^{\beta-1}e^{-x} \sinh \sqrt{x} dx
=
\sum_{k=1}^\infty \frac{\Gamma(\beta+k)}{\Gamma(2k)}
.
\ee

\cite{BradleyCS2001}
\be
\int_0^{\pi/2} \sinh^{-1}(\sin x)dx=
\int_0^{\pi/2} \sinh^{-1}(\cos x)dx
=G.
\ee

\cite{BradleyCS2001}
\be
\int_0^{\pi/2} \csch^{-1}(\csc x)dx=
\int_0^{\pi/2} \csch^{-1}(\sec x)dx
=G.
\ee

\cite{MuthumSci22}
\be
\int_0^\infty \sin 2a \cos^2 bx \tanh x dx=\frac{1}{4}[\frac{2e^{-2a}}{1-e^{-2a}}+
\frac{e^{-(a+b)}}{1-e^{-2(a+b)}}
-\frac{e^{-2(a-b)}}{1-e^{-2(a-b)}}
], a<b.
\ee

\cite{MuthumSci22}
\be
\int_0^\infty \sin 2a \cos^2 ax \tanh x dx=\frac{1}{4}[\frac{e^{-2a}}{1-e^{-2a}}+
2\frac{e^{-a}}{1-e^{-2a}}
].
\ee

\cite[2.3.2]{Nahin}
\be
\int_0^\infty \frac{dx}{x^4+2x^2\cosh(2\alpha)+1}=\frac{\pi}{4\cosh \alpha}.
\label{eq.x4denom}
\ee

\cite{LiMath10}
Let
\be
H(m,n)\equiv \int_0^1 \frac{\arctanh ^m x}{x^n}dx
\ee
then
\be
H(m,0)=\left\{
\begin{array}{ll}
1,& m=0;\\
\frac{m!}{2^{m-1}}\eta(m), m\ge 1.
\end{array}
\right.
\ee
\be
H(m,1)=
\frac{m!}{2^{m-1}}\lambda(m+1), m\ge 1.
\ee
\be
H(m,n)=
\frac{m}{n-1}H(m-1,n-1)+\frac{n-3}{n-1}H(m,n-2),\quad n\ge n, n\ge 2
\ee
\be
H(m,n)=
\left\{
\begin{array}{ll}
\frac{1}{1-n},& m=0,\\
\frac{n+1}{n-1} H(m,n+2)-\frac{m}{n-1}H(m-1,n+1),& m\ge 1
\end{array}
\right. n<0
\label{eq.LiH}
\ee
where $\lambda(x)=\sum_{n\ge 0} 1/(2n+1)^x$ is Dirichlet's lambda function
and $\eta(x)=\sum_{n\ge 0}(-)^n/(n+1)^x$ is Dirichlet's eta function.

\subsection{Rational Functions of Sines and Cosines}
\cite{DenMComp23}
\be
\int_{-\pi/2}^{\pi/2} \cos^\alpha \theta \cos (2n\theta) d\theta
=
(-)^{n+1}\frac{2}{\surd \pi}\sin(\frac{\alpha\pi}{2})
\Gamma(\frac{1+\alpha}{2})
\Gamma(1-\frac{\alpha}{2})
\frac{1}{2n+\alpha}
\times\left\{
\begin{array}{rl}
1,& n=1\\
\prod_{k=1}^{n-1}\frac{2(n-k)-\alpha}{2(n-k)+\alpha},&n\ge 2.
\end{array}
\right.
\ee

\cite{KolbigMathComp64_449}
\be
\int_0^{\pi/2}\frac{\sin 2nx \cos^{2m+1}x}{\sin x}dx=\frac{\pi}{2},\,
n>m\ge 0.
\ee

\cite[B2b]{DaiJOSAA12}
\be
\int_0^{2\pi}\sin^m\theta \cos^n\theta d\theta = 2\pi
\epsilon_m\epsilon_n \frac{(m-1)!!(n-1)!!}{(m+n)!!}
\ee
where $\epsilon_j=1$ if $j$ is even and $\epsilon_j=0$ otherwise.

\be
\int_{\pi/4}^{\pi/2} \cos^2 t \sin^{2k}t dt=\alpha_k +\beta_k \pi
\ee
where $\alpha_0=-1/4$, $\alpha_1=0$, $\alpha_2=1/48$, $\beta_0=1/8$,
$\beta_1=1/32$,
$4(k+1)\alpha_k +2(1-3k)\alpha_{k-1}+(2k-3)\alpha_{k-2}=0$,
$2(k+1)\beta_k+(1-2k)\beta_{k-1}=0.$
$\sum_{k \ge 0} \beta_k z^k = \frac{1}{4(1+\sqrt{1-x})}.$

\cite{BraggAMM106}
\begin{equation}
\int_0^\pi \cos^{pr}(q\theta)\cos^{qr}(p\theta)d\theta
= \frac{\pi}{2^{r(p+q)}}\sum_{l=0}^r\binom{qr}{ql}\binom{pr}{pl}.
\end{equation}
where $p$ and $q$ are coprime and $r$ is a positive integer.

\cite{BraggAMM106}
\begin{equation}
\int_0^\pi (\cos p\phi)^m(\cos q\phi)^n\cos(mp-nq)\phi d\phi
= \frac{\pi}{2^{m+n}}\sum_{l=0}^H \binom{m}{ql}\binom{n}{pl}.
\end{equation}
where $H=\min\{[m/q],[n/p]\}$.

\cite{BraggAMM106}
\begin{equation}
\frac{2^n}{\pi} \int_0^\pi \cos^n \phi \cos((n-m)\phi)
\frac{\sin((m+1)\phi)}{\sin \phi}d\phi
= \sum_{j=0}^m\binom{n}{j}.
\end{equation}

\cite{BraggAMM106}
\begin{equation}
\int_0^\pi \cos^{pq} \phi 
\frac{\sin q(p+1)\phi}{\sin q\phi}d\phi
= \frac{2\pi}{2^{pq}}\sum_{k=0}^p\binom{pq}{kq}.
\end{equation}

\cite[2.3.8]{Nahin}\cite{TaylorPEMS32}
\be
\int_0^\pi \frac{\cos(nx)-\cos(n\alpha)}{\cos(x)-\cos(\alpha)}dx = \pi \frac{\sin(n\alpha)}{\sin(\alpha)}.
\ee

\cite{TaylorPEMS32}
\be
\int_0^\pi \frac{(\cos n\alpha-\cos n\theta)\cos r \theta}{\cos \alpha-\cos\theta}d\theta=
\pi\frac{\sin(n-r)\alpha}{\sin \alpha}.
\ee
\be
\int_0^\pi \frac{1-\cos n\theta)\cos r \theta}{1-\cos\theta}d\theta=
(n-r)\pi.
\ee
\cite{TaylorPEMS32}
\be
\int_0^\pi \frac{\cos n\theta}{1-2p\cos\theta+p^2}d\theta
=\pi\frac{p^{-n}}{p^2-1},\quad p>1
\ee
\be
\int_0^\pi \frac{\cos n\theta}{1-2q\cos\theta+q^2}d\theta
=\pi\frac{q^n}{1-q^2},\quad 0<q<1.
\ee
\cite{TaylorPEMS32}
\be
\int_0^\pi \frac{\cos n\theta}{1-\cos\beta\cos\theta}d\theta
=\frac{\pi}{\sin\beta}(\frac{1-t}{1+t})^n,\quad t\equiv \tan(\beta/2).
\ee
\cite{TaylorPEMS32}
\be
\int_0^\pi \frac{\cos n\theta}{1+\cos\beta\cos\theta}d\theta
=(-)^n \frac{\pi}{\sin\beta}(\frac{1-t}{1+t})^n,\quad t\equiv \tan(\beta/2).
\ee

\cite[2.3.8]{Nahin}
Let
\be
I_n=\int_0^{\pi/2} \frac{1}{(a\cos^2 x+\sin^2 x)^2}dx,\quad n=1,2,3\ldots
\ee
then 
\be
I_1=\frac{\pi}{2\sqrt{ab}}
\ee
and the others by differentiation 
\be
I_n=\frac{1}{1-n}[\frac{\partial I_{n-1}}{\partial a}+\frac{\partial I_{n-1}}{\partial b}].
\ee

\cite{BaileyPCPS27_3}
\be
\int_0^z \sin^\mu t\sin^\nu(z-t)dt=\frac{\sqrt{\pi}\Gamma(\mu+1)\Gamma(\nu+1)}
{2^{(\mu+\nu+1)/2}\Gamma(\mu/2+\nu/2+1)}\sin^{(\mu+\nu+1)/2}z
P_{(\mu-\nu-1)/2}^{-(\mu+\nu+1)/2}(\cos z),
\Re \mu>-1,\Re \nu>-1.
\ee

\cite{BaileyPCPS27_3}
\be
\frac{2^m \Gamma(m+1/2)}{\sqrt{\pi}\Gamma(m+n+1)\Gamma(m-n)}\int_0^z
\sin^{m+n}t\sin^{m-n-1}(z-t)dt = \sin^mz P_n^{-m}(\cos z).
\ee

\cite{BaileyPCPS27_3}
\be
\int_0^z \left(\frac{\sin t}{\sin(z-t)}\right)^\mu dt = \frac{\pi \sin\mu z}{\sin\mu \pi}, \quad -1<\Re \mu <1.
\ee

\cite{MatharArxiv0706}
\be
\int_0^{\arcsin q} \left(\cos\phi\pm \sqrt{q^2-\sin^2\phi}\right)^{n+2}T_m(\cos \phi)d\phi=\ldots
\ee

\cite{Adamchik}
\be
\int_0^{\pi/2}\sinh^{-1}(\sin x)dx=G.
\ee

\cite{Adamchik}
\be
\int_0^{\pi/2}\sinh^{-1}(\cos x)dx=G.
\ee

\cite{Adamchik}
\be
\int_0^{\pi/2}\csch^{-1}(\csc x)dx=G.
\ee

\cite{Adamchik}
\be
\int_0^{\pi/2}\csch^{-1}(\sec x)dx=G.
\ee

\cite[p40]{MO2Afl}
\be
\int_0^\infty \cos(T_n(t,-x))dt=\frac{\pi\sqrt{x}}{2n\sin\frac{\pi}{2n}}
\left[J_{1/n}(2x^{n/2})-J_{-1/n}(2x^{n/2})\right],
\ee
\be
\int_0^\infty \cos(T_{2m}(t,x))dt=\frac{\pi\sqrt{x}}{4m\sin\frac{\pi}{4m}}
\left[J_{-1/(2m)}(2x^m)-J_{1/(2m)}(2x^m)\right],\quad m=1,2,3\ldots
\ee
\be
\int_0^\infty \cos(T_{2m+1}(t,x))dt=\frac{2\sqrt{x}\cos\frac{\pi}{4m+2}}{2m+1}
K_{1/(2m+1)}(2x^{m+1/2}),\quad m=1,2,3\ldots
\ee
where $x$ real positive, where
\be
T_n(t,x)\equiv t^n\,_2F_1(-\frac{n}{2},\frac{1-n}{2};1-n;-\frac{4x}{t^2}),
\quad
n=2,3,4,\ldots
\ee
for example $T_2=t^2+2x$, $T_3=t^3+3tx$.

\subsection{Trigonometric and Rational Functions}

\cite[3.4.1]{Nahin}
\be
\int_0^\infty \frac{\cos(ax)-\cos(bx)}{x^2}dx=\frac{\pi}{2}(b-a).
\ee

\cite[(4.3.3)]{Nahin}
\be
\int_0^\infty \frac{\sin (x^q)}{x^q}dx=\frac{\pi}{2q\Gamma(2-1/q)\sin(\pi/(2q))}
=\frac{\Gamma(1/q)}{q-1}\cos(\frac{\pi}{2q}),\quad q>1
.
\ee

\cite[(4.3.9)]{Nahin}
\be
\int_0^\infty \frac{\cos (bx)}{x^p}dx=\frac{b^{p-1}\pi}{2\Gamma(p)\cos(\frac{p\pi }{2})}, \quad 0<p<1.
\ee

\cite[(4.3.11)]{Nahin}
\be
\int_0^\infty \cos (bx^k)dx=\frac{\Gamma(1/k)\cos\frac{\pi}{2k}}{kb^{1/k}}, \quad b>0,k>1.
.
\ee

\cite[(4.3.14)]{Nahin}
\be
\int_0^\infty \frac{e^{-rx}\cos(px)-e^{-sx}\cos(qx)}{x}dx
=\frac12 \ln \frac{q^2+s^2}{p^2+r^2}.
.
\ee

\cite[(C4.9)]{Nahin}
\be
\int_0^\infty \frac{\sin (x^2)}{\surd x}dx
=\frac{\pi}{4\Gamma(3/4)\cos(3\pi/8)}
.
\ee

\cite{Adamchik,BradleyCS2001,CoffeyJMP49,ChoiJMAA231}
\be
\frac{1}{2}\int_0^{\pi/2}\frac{x}{\sin x} = G.
\ee

\cite{Adamchik,BradleyCS2001}
\be
\frac{3}{4}\int_0^{\pi/6}\frac{x}{\sin x} = G-\frac{\pi}{8}\log(1+\sqrt{3}).
\ee

\cite{CoffeyJMP49,BradleyCS2001,ChoiJMAA231}
\be
\int_0^{\pi/4}\frac{x^2}{\sin^2 x}dx = G-\frac{\pi}{16}(\pi -4\ln 2).
\ee

\cite{Adamchik,BradleyCS2001}
\be
-\frac{\pi^2}{4}\int_0^{1}(x-\frac{1}{2})\sec(\pi x)dx = G.
\ee

\cite{Amdeberhanarxiv3663}
\be
\frac{2^{s+2}}{\pi}\int_0^{\pi/2}x \cos^s x \sin(sx)dx-\gamma
=\psi(s+1).
\ee

\cite{KolbigMathComp64_449}
\be
\int_0^{\pi/4}x^m\tan x dx = 2\pi G-\frac{7}{2}\zeta(3),\, m=2.
\ee

\cite{CoffeyJMP49}
\be
\int_0^{\pi/2} \frac{a}{\tan a}da=\frac{\pi}{2}\ln 2.
\ee

\cite{CoffeyJMP49}
\be
\int_0^{\pi/4}\frac{a^2}{\tan^2 a}da = G-\frac{\pi}{16}(\pi -4\ln 2)-\frac{\pi^3}{192}.
\ee

\cite{CrandallEM3,CvijovicJMAA351}
\be
\int_0^{1/2}x^n\cot\pi x dx
=
\frac{n!}{2^n}\sum_{k=1,k\mathrm{odd}}^n
\frac{(-)^{(k-1)/2}}{\pi^k}\frac{\eta(k)}{(n-k+1)!}
+\frac{(-)^n+1}{2}\frac{4n!(1-2^{-n-1})}{(2\pi)^{n+1}}\zeta(n+1),
\ee
where $\eta(s)\equiv (1-2^{1-s})\zeta(s)$.

\cite{CvijovicJMAA351}
\be
\int_0^{\pi/4}x^n\cot x dx
=
\frac12(\frac{\pi}{4})^n(\frac{2}{n}-\sum_{k=1}^\infty \frac{\zeta(2k)}{4^{2k-1}(n+2k)}), n\in N
\ee

\cite{CvijovicJMAA351}
\begin{multline}
\int_0^{p\pi/q}x^n\cot x dx
=
(\frac{p\pi}{q})^n
\ln[2\sin(\frac{p\pi}{q})]
+(-)^{\lfloor n/2\rfloor}\{1+(-)^n\}
2^{-n-1}n!\zeta(n+1)+2n!(p\pi/q)^{n+1}
\\
[
\sum_{k=1}^{\lfloor(n+1)/2\rfloor}\sum_{l=1}^q
\frac{(-)^{k-1}}{(n-2k+1)!}\frac{1}{(2p\pi)^{2k}}\sin\frac{2lp\pi}{q}
\zeta(2k,l/q)
+\sum_{k=1}^{\lfloor n/2\rfloor}\sum_{l=1}^q
\frac{(-)^{k-1}}{(n-2k)!}\frac{1}{(2p\pi)^{2k+1}}\cos\frac{2lp\pi}{q}
\zeta(2k+1,l/q)
]
\end{multline}
$0<p/q<1$, where $\zeta$ is the Hurwitz zeta-function.

\cite{CvijovicJMAA351}
\begin{multline}
\int_0^{p\pi/q}x^n\tan x dx
=
-(\frac{p\pi}{q})^n
\ln[2\cos(\frac{p\pi}{q})]
+(-)^{\lfloor n/2\rfloor}\{1+(-)^n\}
2^{-n-1}(1-2^{-n})n!\zeta(n+1)+2n!(p\pi/q)^{n+1}
\\
[
\sum_{k=1}^{\lfloor(n+1)/2\rfloor}\sum_{l=1}^q
\frac{(-)^{k-1}}{(n-2k+1)!}\frac{(-)^{l-1}}{(2p\pi)^{2k}}\sin\frac{2lp\pi}{q}
Z_q(2k,l/q)
\\
+\sum_{k=1}^{\lfloor n/2\rfloor}\sum_{l=1}^q
\frac{(-)^{k-1}}{(n-2k)!}\frac{(-)^{l-1}}{(2p\pi)^{2k+1}}\cos\frac{2lp\pi}{q}
Z_q(2k+1,l/q)
]
\end{multline}
$0<p/q<1/2$, where $Z_q(s,a)$ is $\eta(s,a)$ if $q$ is odd
and $\zeta(s,a)$ if $q$ is even. $\eta(s,a) = [\zeta(s,a/2)-\zeta(s,a/2+1/2)]/2$.

\cite{CvijovicJMAA351}
\begin{multline}
\int_0^{p\pi/q}x^n\csc x dx
=
(\frac{p\pi}{q})^n
\ln[\tan(\frac{p\pi}{2q})]
+(-)^{\lfloor n/2\rfloor}\{1+(-)^n\}
n! (1-2^{-n-1})\zeta(n+1)+2n!(p\pi/q)^{n+1}
\\
[
\sum_{k=1}^{\lfloor(n+1)/2\rfloor}\sum_{l=1}^q
\frac{(-)^{k-1}}{(n-2k+1)!}\frac{1}{(2p\pi)^{2k}}\sin\frac{(2l-1)p\pi}{q}
\zeta(2k,(2l-1)/(2q))
\\
+\sum_{k=1}^{\lfloor n/2\rfloor}\sum_{l=1}^q
\frac{(-)^{k-1}}{(n-2k)!}\frac{1}{(2p\pi)^{2k+1}}\cos\frac{(2l-1)p\pi}{q}
\zeta(2k+1,(2l-1)/(2q))
]
\end{multline}
$0<p/q<1/2$.

\cite[p 6]{OberhettFT}
\be
\int_0^\infty \frac{\cos(xy)}{(a^2+x^2)^{\nu+1/2}}dx = \sqrt{\pi}\left(\frac{y}{2a}\right)^\nu \frac{1}{\Gamma(\nu+1/2)}K_\nu(ay).
\ee

\cite[(C8.5)]{Nahin}
\be
\int_{-\infty}^\infty \frac{\cos(mx)}{ax^2+bx+c}dx = -2\pi \frac{\cos\frac{mb}{2a}\sin\frac{m\sqrt{b^2-4ac}}{2a}}{\sqrt{b^2-4ac}},\quad \mathrm{when} b^2>4ac.
\ee

\cite[p 10]{OberhettFT}
\begin{multline}
\int_0^\infty \frac{x^\nu}{(a^2+x^2)^{\mu+1}}\cos(xy)dx
= \frac{a^{\nu-2\mu-1}}{2}
B\left(\frac{1}{2}+\frac{1}{2}\nu,\mu-\frac{1}{2}\nu+\frac{1}{2}\right)
\, _1F_2(\frac{\nu+1}{2};\frac{\nu+1}{2}-\mu,\frac{1}{2},\frac{a^2y^2}{4})
\\
+\sqrt{\pi}\frac{2^{-2\mu+\nu-2}}{\Gamma(1+\mu-\frac{1}{2}\nu)}
y^{2\mu-\nu+1}\Gamma(\frac{1}{2}\nu-\mu-\frac{1}{2})
\,_1F_2(\mu+1-\frac{\nu}{2};\mu-\frac{\nu}{2}+\frac{3}{2};\frac{a^2y^2}{4})
.
\end{multline}

\cite{KolbigMathComp64_449}
\be
\int_0^\infty \frac{x^{2m}\cos(ax)dx}{(\beta^2+x^2)^{n+1/2}}
=
\frac{(-1)^m\sqrt{\pi}}{2^n\beta^n\Gamma(n+1/2)}
\cdot \frac{d^{2m}}{da^{2m}} \{ a^nK_n(a\beta)\},\,
a>0,\,\Re\beta>0,\, 0\le m\le n
.
\ee

\cite[p 116]{OberhettFT}
\be
\int_0^\infty \frac{\sin(xy)}{x(a^2+x^2)^{\nu+1/2}}dx
=
\frac{\pi y}{2a^{2\nu}}
[K_\nu{\mathbf L}_{\nu-1}(ay)+{\mathbf L}_{\nu}(ay)K_{\nu-1}(ay)]
\ee
where $\mathbf L$ are Struve functions.

\cite{Cvijovicarxiv0911}
\be
\int _0^{\delta \pi} \sin(2k+1)t \frac{2z(1+z^2)\sin t}
{1-2z^2\cos(2t)+z^4}dt = \delta \pi z^{2k+1}
\ee
\be
\int _0^{\delta \pi} \sin(2k+1)t \frac{2z(1-z^2)\cos t}
{1-2z^2\cos(2t)+z^4}dt = \delta \pi z^{2k+1}
\ee
for $\delta=1/2$ and (??) 1, $k$ a nonnegative integer.

\cite{Adamchik}
\be
\frac{1}{2\pi}\int_0^{\pi/2}\frac{x^2}{\sin x} = G-\frac{7}{4\pi}\zeta(3).
\ee

\cite{Adamchik,BradleyCS2001}
\be
\int_0^{\pi/2}\frac{x \csc x}{\cos x + \sin x} = G+\frac{\pi}{4}\log 2.
\ee

\cite{Adamchik,BradleyCS2001}
\be
-2\int_0^{\pi/2}\frac{x \cos x}{\cos x + \sin x} = G-\frac{\pi^2}{8}-\frac{\pi}{4}\log 2.
\ee

\cite{Adamchik,BradleyCS2001}
\be
2\int_0^{\pi/2}\frac{x \sin x}{\cos x + \sin x} = G+\frac{\pi^2}{8}-\frac{\pi}{4}\log 2.
\ee

\cite{ChoiJMAA231}
\be
\int_0^{\pi/4}x\cot x dx = \frac{\pi}{8}\log 2 +\frac{G}{2}.
\ee

\cite{ChoiJMAA231}
\be
\int_0^{\pi/4}\log \cos x dx = -\frac{\pi}{4}\log 2 +\frac{G}{2}.
\ee

\cite{ChoiJMAA231}
\be
\int_0^{\pi/4} x\tan x dx = -\frac{\pi}{8}\log 2+\frac{G}{2}.
\ee

\cite{ChoiJMAA231}
\be
\int_0^{\pi/4}(\frac{x}{\cos x})^2 dx = \frac{\pi^2}{16}+\frac{\pi}{4}\log 2 -G.
\ee

\cite{ChoiJMAA231}
\be
\int_0^{\pi/4} x^2\tan^2 x dx = -\frac{\pi^3}{192}
+\frac{\pi^2}{16}+\frac{\pi}{4}\log 2 -G.
\ee

\cite{ChoiJMAA231}
\be
\int_0^{\pi/4} x^2\cot^2 x dx = -\frac{\pi^3}{192}
-\frac{\pi^2}{16}+\frac{\pi}{4}\log 2 +G.
\ee

\cite{Amdeberhanarxiv2379}
\be
\int_0^{\infty} x^{-p}\cos^{2n+1}(x+b)dx=
\frac{\Gamma(1-p)}{2^{2n}}\sum_{k=0}^n\binom{2n+1}{n-k}
\frac{\sin[\pi p/2-(2k+1)b]}{(2k+1)^{1-p}}
.
\ee

\cite{Amdeberhanarxiv2379}
\be
\int_0^{\infty} x^{-p}\sin^{2n+1}(x+b)dx=
\frac{\Gamma(1-p)}{2^{2n}}\sum_{k=0}^n\binom{2n+1}{n-k}
\frac{\cos[\pi p/2-(2k+1)b]}{(2k+1)^{1-p}}
\ee

\cite{Amdeberhanarxiv2379}
\be
\int_0^{\infty} x^{-p}\cos^{2n+1}x dx=
\frac{\Gamma(1-p)}{2^{2n}}\sin\left(\frac{\pi p}{2}\right)\sum_{k=0}^n\frac{\binom{2n+1}{n-k}}{(2k+1)^{1-p}}
\ee
for $0<p<1$.

\cite[(7.6.4)]{Nahin}
\be
\int_0^\infty \frac{\sin^{2n-1}x}{x}dx = \frac{\pi}{2^{2n-1}}\binom{2n-2}{n-1}.
\ee

\cite[(7.6.5)]{Nahin}
\be
\int_0^\infty \frac{\sin^{2n-1}x \cos x}{x}dx = \frac{\pi}{2^{2n}n}\binom{2n-2}{n-1}.
\ee

\cite{ThompsonMC20}
Let $I_n(b)=\frac{2}{\pi}\int_0^\infty (\frac{sin x}{x})^n\cos(bx)dx$, then
\be
I_n(b)=\frac{1}{2(n-1)}[(n+b)I_{n-1}(b+1)+(n-b)I_{n-1}(b-1)]
\ee
starting at
\be
I_3(b) = \left\{
\begin{array}{ll}
\frac{1}{8}[(b+3)^2-3(b+1)^2], & 0\le b\le 1\\
\frac{1}{8}[(b+3)^2-3(b+1)^2+e(b-1)^2], & 1\le b\le 3\\
0, & b\ge 3.
\end{array}
\right.
\ee
Negative $b$ covered by $I_n(-b)=I_n(b)$.
\cite{MedhurstMC19}
\be
I_n(b)=\frac{n}{2^{n-1}}\sum_{0\le r<(b+n)/2}
\frac{(-1)^r(b+n-2r)^{n-1}}{r!(n-r)!}
,\quad 0\le b<n
\ee
$I_n(b)=0$ for $n\le b<\infty$.

\cite[(7.6.6)]{Nahin}
\be
\int_0^\infty \frac{\sin^{2n}x}{x^2}dx = \frac{\pi}{2^{2n-1}}\binom{2n-2}{n-1}.
\ee

\cite{Amdeberhanarxiv2379}
\be
\int_0^{\infty} x^{-p}\sin^{2n+1}x dx=
\frac{\Gamma(1-p)}{2^{2n}}\cos\left(\frac{\pi p}{2}\right)\sum_{k=0}^n(-1)^k\frac{\binom{2n+1}{n-k}}{(2k+1)^{1-p}}
\ee
for $0<p<1$.

\cite{Amdeberhanarxiv2379}
\be
\int_0^{\infty} \cos^{2n+1} x^p dx=
\frac{1}{2^{2n}}\Gamma\left(\frac{p+1}{p}\right)
\cos\left(\frac{\pi}{2p}\right)\sum_{k=0}^n\frac{\binom{2n+1}{n-k}}{(2k+1)^{1/p}}
\ee
for $p>1$.

\cite{Amdeberhanarxiv2379}
\be
\int_0^{\infty} \sin^{2n+1} x^p dx=
\frac{1}{2^{2n}}\Gamma\left(\frac{p+1}{p}\right)
\sin\left(\frac{\pi}{2p}\right)\sum_{k=0}^n(-1)^k\frac{\binom{2n+1}{n-k}}{(2k+1)^{1/p}}
\ee
for $p>1$.

\cite{Amdeberhanarxiv2379}
\be
\int_0^{\pi/2} x^p\cos^{2n}x dx=
\sum_{j=0}^{\lfloor p/2 \rfloor} a_{n,p,p+1-2j}\pi^{p+1-2j}
+\delta_{odd,p}\cdot a^*_{n,p},
\ee
where for $p\ge 2$ and $0\le j\le \lfloor p/2\rfloor$
\be
a_{n,p,p+1-2j}=\frac{(-)^j\binom{2n}{n}p!}{2^{2n+p+1}(p+1-2j)!}\sum_{1\le k_1\le k_2\le\cdots \le k_j\le n}
\frac{1}{k_1^2k_2^2\cdots k_j^2}.
\ee
and $a_{n,p}^*$ is a similar multinomial sum.
A similar form exists for odd powers of the cosine.

\cite[(5.1.4)]{Nahin}
\be
\int_0^\pi \frac{x\sin x}{a+b\cos^2x}dx = \frac{\pi}{\sqrt{ab}}\tan^{-1}(\sqrt{b}{a})
,\quad a>b.
\ee

\cite{KolbigMathComp64_449}
\be
\int_0^\infty \sin^{2m+1}x\frac{xdx}{a^2+x^2}
=
\frac{\pi}{2^{2m+1}}e^{-(2m+1)a}\sum_{k=0}^m (-1)^{m+k}
\binom{2m+1}{k}e^{2ka}.
\ee

\cite{KolbigMathComp64_449}
\be
\int_0^\infty \cos^{2m}x\frac{xdx}{a^2+x^2}
=
\frac{\pi}{2^{2m+1}a}\binom{2m}{m}+\frac{\pi}{2^{2m}a}\sum_{k=1}^m
\binom{2m}{m+k}e^{-2ka},\, a>0.
\ee

\cite{ZerrAMM6}
\be
\int_0^{\pi/2}\frac{\sin^2\theta d\theta}{\sqrt{1-e^2\sin^2\theta}}
=\frac{1}{e^2}[F(e,\pi/2)-E(e,\pi/2)].
\ee

\cite{ZerrAMM6}
\be
\int_0^{\pi/2}\frac{\cos^2\theta d\theta}{\sqrt{1-e^2\sin^2\theta}}
=\frac{1}{e^2}[E(e,\pi/2)-(1-e^2)F(e,\pi/2)].
\ee

\cite{ZerrAMM6}
\be
\frac{1}{\pi}\int_0^{2\pi}\frac{\cos 2\phi d\phi}{\sqrt{1+e^2-2e\cos\phi}}
=\frac{4}{3\pi e^2}[(2+e^2)F(e,\pi/2)-2(1+e^2)E(e,\pi/2)].
\ee

\cite{ZerrAMM6}
\be
\int_0^{\pi/2}\frac{\sin^4 \theta d\theta}{\sqrt{1-e^2\sin^2\theta}}
=\frac{1}{3 e^4}[(2+e^2)F(e,\pi/2)-2(1+e^2)E(e,\pi/2)].
\ee

\cite{ZerrAMM6}
\be
\int_0^{\pi/2}\frac{\sin^2 \theta \cos^2 \theta d\theta}
{\sqrt{1-e^2\sin^2\theta}}
=\frac{1}{3 e^4}[(2-e^2)E(e,\pi/2)-2(1-e^2)F(e,\pi/2)].
\ee

\cite{ZerrAMM6}
\be
\int_0^{\pi/2}\frac{\cos^4 \theta d\theta}{\sqrt{1-e^2\sin^2\theta}}
=\frac{1}{3 e^4}[(2-5e^2+3e^4)F(e,\pi/2)-2(1-2e^2)E(e,\pi/2)].
\ee

\cite{ZerrAMM6}
\be
\int_0^{\pi/2}\sqrt{1-e^2\sin^2\theta}\sin^2\theta d\theta
=\frac{1}{3 e^2}[(1-e^2)F(e,\pi/2)-(1-2e^2)E(e,\pi/2)].
\ee

\cite{ZerrAMM6}
\be
\int_0^{\pi/2}\sqrt{1-e^2\sin^2\theta}\cos^2\theta d\theta
=\frac{1}{3 e^2}[(1+e^2)E(e,\pi/2)-(1-e^2)F(e,\pi/2)].
\ee

\cite{ZerrAMM6}
\be
\frac{1}{\pi}\int_0^{2\pi}\frac{\cos 3\phi d\phi}{\sqrt{1+e^2-2e\cos\phi}}
=\frac{4}{15\pi e^3}[(8+3e^2+4e^4)F(e,\pi/2)-(8+7e^2+8e^4)E(e,\pi/2)].
\ee

\cite{ZerrAMM6}
\be
\int_0^{\pi/2}\frac{\sin^6 \theta d\theta}{\sqrt{1-e^2\sin^2\theta}}
=\frac{1}{15 e^6}[(8+3e^2+4e^4)F(e,\pi/2)-(8+7e^2+8e^4)E(e,\pi/2)].
\ee

\cite{ZerrAMM6}
\be
\int_0^{\pi/2}\frac{\sin^4 \theta \cos^2\theta d\theta}{\sqrt{1-e^2\sin^2\theta}}
=\frac{1}{15 e^6}[(8-3e^2-2e^4)E(e,\pi/2)-(8-7e^2-e^4)F(e,\pi/2)].
\ee

\cite{ZerrAMM6}
\be
\int_0^{\pi/2}\frac{\sin^2 \theta \cos^4\theta d\theta}{\sqrt{1-e^2\sin^2\theta}}
=\frac{1}{15 e^6}[(8-17e^2+9e^4)F(e,\pi/2)-(8-13e^2+3e^4)E(e,\pi/2)].
\ee

\cite{ZerrAMM6}
\be
\int_0^{\pi/2}\frac{\cos^6\theta d\theta}{\sqrt{1-e^2\sin^2\theta}}
=\frac{1}{15 e^6}[(8-23e^2+23e^4)E(e,\pi/2)-(8-27e^2+34e^4-15e^6)F(e,\pi/2)].
\ee
and similar expressions for combined 8th powers in the numerator.

\cite{ZerrAMM6}
\be
\int_0^{\pi/2}\sqrt{1-e^2\sin^2\theta}\sin^2\theta\cos^2\theta d\theta
=\frac{1}{15 e^4}[2(1-e^2+e^4)E(e,\pi/2)-(2-3e^2+e^4)F(e,\pi/2)].
\ee

\cite{ZerrAMM6}
\be
\int_0^{\pi/2}\sqrt{1-e^2\sin^2\theta}\sin^4\theta d\theta
=\frac{1}{15 e^4}[2(1+e^2-2e^4)F(e,\pi/2)-(2+3e^2-8e^4)E(e,\pi/2)].
\ee

\cite{ZerrAMM6}
\be
\int_0^{\pi/2}\sqrt{1-e^2\sin^2\theta}\cos^4\theta d\theta
=\frac{1}{15 e^4}[2(1-4e^2+3e^4)F(e,\pi/2)-(2-7e^2-3e^4)E(e,\pi/2)].
\ee

\cite{EvansAMC36}
\be
\int_0^{2\pi}\sqrt{a+\sin \sigma}d\sigma = 2\pi \sqrt{a}(1-2\sum_{n\ge 1} \frac{(4n-2)!}{(n!)^2(2n-1)!64^n a^{2n}}).
\ee

\cite{EvansAMC36}
\be
\int_0^{2\pi}\sqrt{a+\sin^2 \sigma}d\sigma = 2\pi \sqrt{a+1/2}(1-2\sum_{n\ge 1} \frac{(4n-2)!}{(n!)^2(2n-1)!256^n (a+1/2)^{2n}}).
\ee

\cite{EvansAMC36}
\be
\int_0^{2\pi}\frac{d\sigma}{\sqrt{a+\sin \sigma}} = \frac{2\pi}{\sqrt{a}}
(1+2\sum_{n\ge 1} \frac{(4n-1)!}{(n!)^2(2n-1)!64^n a^{2n}}).
\ee

\cite{EvansAMC36}
\be
\int_0^{2\pi}\frac{d\sigma}{\sqrt{a+\sin^2 \sigma}} = \frac{2\pi}{\sqrt{a+1/2}}
(1+2\sum_{n\ge 1} \frac{(4n-1)!}{(n!)^2(2n-1)!256^n (a+1/2)^{2n}}).
\ee

\cite{FittMG72}
\be
\int_0^{\pi/2} \frac{\sqrt{\sin x}}{\sqrt{\sin x}+\sqrt{\cos x}}
=\pi/4.
\ee

\cite{FittMG72}
\be
\int_{-\alpha}^{\alpha} \frac{\cos \lambda_m\theta}{(\cos \theta)^{2+\lambda_m}}d\theta
=2(-)^{m+2}\frac{(\cos \alpha)^{-1-\lambda_m}}{1+\lambda_m}.
\ee
where $\lambda_m = \frac{(2m-1)\pi}{2\alpha}-1$.

\cite{ElRabiiSci14}
\[
I(a,b)\equiv \int_0^\infty x^{-a}\left(1-\frac{\sin^b x}{x^b}\right)dx.
\]
then
\be
I(a,b)= \frac{\pi \sec(\pi a/2)}{2^b\Gamma(a+b)}
\sum_{k=0}^{\lfloor (b-1)/2\rfloor}
(-1)^{k+1}\binom{b}{k} (b-2k)^{a+b-1},
\ee
with special cases
\be
I(2,b)= \frac{\pi}{2^b(b+1)!}
\sum_{k=0}^{\lfloor (b-1)/2\rfloor}
(-1)^{k+1}\binom{b}{k} (b-2k)^{b+1},
\ee
\be
I(a,1)= -\frac{\pi \sec(\pi a/2)}{2\Gamma(1+a)}\,\quad
\int_0^\infty x^{-3/2}(1-\frac{\sin x}{x})dx=\frac{2\sqrt{2\pi}}{3},
\ee
\be
I(a,2)= -\frac{\pi 2^{a-1}\sec(\pi a/2)}{\Gamma(2+a)},\quad
\int_0^\infty x^{-3/2}(1-\frac{\sin^2 x}{x^2})dx=\frac{16\sqrt{\pi}}{15},
\ee
\be
I(a,3)= \frac{(3-3^{2+a})\pi \sec(\pi a/2)}{8\Gamma(3+a)}\,\quad
\int_0^\infty x^{-3/2}(1-\frac{\sin^3 x}{x^3})dx=\frac{2}{35}(9\sqrt{3}-1)\sqrt{2\pi}.
\ee

\cite[(7.3.6)]{Nahin}
\be
\int_0^1 \frac{\sin \ln x}{\ln x}dx=\frac{\pi}{4}.
\ee
\subsection{Trigonometric Functions and Exponentials}

\cite[(7.11.28)]{Nahin}
\be
\int_0^\infty e^{-\beta x^2}\frac{\sin x}{x} dx=\frac{\pi}{2}\erf(\frac{1}{2\surd \beta})
\ee

\cite[(8.7.20)]{Nahin}
\be
\int_0^\infty \frac{\cos x}{e^x+e^{-x}}dx =\frac{\pi/2}{e^{\pi/2}+e^{-\pi/2}}. 
\ee

\cite{LeonardMM95}
\be
\int_0^\infty e^{-tx^2}\cos x^2 dx = \sqrt{\pi/8}
\sqrt{\frac{\sqrt{1+t^2}+t}{1+t^2}}.
\ee

\cite{LeonardMM95}
\be
\int_0^\infty e^{-tx^2}\sin x^2 dx = \sqrt{\pi/8}
\sqrt{\frac{\sqrt{1+t^2}-t}{1+t^2}}.
\ee

\cite{BraggAMM106}
\begin{equation}
\frac{2^n}{\pi}\int_0^{\pi}
e^{-x\cos 2\phi}\cos(x\sin 2\phi+n\phi)\cos^n\phi d\phi = L_n(x)
\end{equation}
where $L_n(z)=\sum_{j=0}^n\frac{(-)^j}{j!}\binom{n}{j}z^j$ are Laguerre
polnomials.

\cite{BraggAMM106}
For $n$ even
\begin{equation}
(-1)^{n/2}\frac{2^n}{\pi}\int_0^\pi e^{-x\cos4\phi}
\cos(x\sin 4\phi+n\phi)\sin^n \phi d \phi
= p_n(x,1),
\end{equation}
for $n$ odd
\begin{equation}
(-1)^{(n-1)/2}\frac{2^n}{\pi}\int_0^\pi e^{-x\cos4\phi}
\sin(x\sin 4\phi+n\phi)\sin^n \phi d \phi
= p_n(x,1),
\end{equation}
where
\begin{equation}
p_n(z,\zeta)=\sum_{k=0}^{[n/2]} \frac{(-)^k}{k!}\binom{n}{2k}z^k\zeta^k.
\end{equation}

\cite{MuthumSci22}
\be
\int_0^\infty e^{-xz} (\cos ax-\cos bx)\frac{dx}{x}
=\frac12 \log\frac{1+b^2 /z^2}{1+a^2/z^2}.
\ee

\cite{MuthumSci22}
\be
\int_0^\infty \frac{e^{-xz}}{1-e^{-x}} (\cos ax -\cos bx)dx
=\log \left| \frac{\Gamma(z-ia)}{\Gamma(z-ib)}\right| .
\ee

\cite{MuthumSci22}
\be
\int_0^\infty \frac{x^{n-1} e^{-xz}}{1-e^{-x}} (\cos ax-\cos bx)dx
=-\frac{(n-1)!}{2} [\zeta(n,z-ia)+\zeta(n,z+ia)-\zeta(n,z-ib)-\zeta(n,z+ib)]
\ee

\cite{MuthumSci22}
\be
\int_0^\infty e^{-xz} \sin^{2m}az \frac{dx}{x}
=\frac{(-)^{m+1}}{2^{2m}}\sum_{j=0}^{m-1}(-)^j\binom{2m}{j} \log[z^2+(2m-2j)a^2]-\frac{1}{2^{2m}}\binom{2m}{m}\log z.
\ee

\cite[(2.11)]{MuthumSci22}
\be
\int_0^\infty e^{-xz} x^{n-1} \frac{1}{1-e^{-x}} (\sin^{2m} ax-\sin^{2m} bx)dx
=\frac{(-)^{m+1}}{2^{2m}}(n-1)!\sum_{j=0}^{m-1}(-)^j\binom{2m}{j} [\zeta(n,z-ia(2m-2j))+\zeta(n,z+ia(2m-2j))
-\zeta(n,z-ib(2m-2j))-\zeta(n,z+ib(2m-2j))].
\ee

\section{Definite Integrals of Elementary Functions II}

\subsection{Logarithmic Functions}
Let
\be
I_k\equiv \int_0^1 x^k\frac{\ln x}{\sqrt{1-x}}dx,\quad k=0,1,2,\ldots
\ee
with special values
\be
I_0 = 4(\ln 2 -1);\quad I_1 = \frac{4}{3}\left(-\frac{5}{3}+2\ln 2\right).
\ee
(This corrects a sign error in \cite[1.3.3.11]{Apelblat2}.)
Then recursively
\be
(2k+1)I_k=I_{k-1}+2(k-1)I_{k-2}+4(S_{2k-3}+S_{2k-1}),\quad k>1,
\ee
where $S_k$ is defined
as
\be
S_k\equiv \int_0^{\pi/2} \sin^k \varphi \cos^2\varphi d\varphi,
\ee
such that \cite[2.510]{GR}
\begin{equation*}
S_0=\frac{\pi}{4};\quad S_1=\frac{1}{3};\quad 
(k+2)S_k = (k-1)S_{k-2}; \quad S_{2k+1}=\frac{(2k)!!}{(2k+3)!!}.
\end{equation*}
As a shortcut another representation is
\be
I_k=\alpha_k+\beta_k \ln 2
\ee
where $\alpha_0=-4$, $\alpha_1=-20/9$, $(2k+1)^2\alpha_k-2(4k^2-2k+1)\alpha_{k-1}+4(k-1)^2\alpha_{k-2}=0$,
$\beta_0=4$ and $(2k+1)\beta_k-2k\beta_{k-1}=0$.

\cite[3.4.3]{Nahin}
\be
\int_0^1\frac{x^a-1}{\ln x}dx=\ln(a+1), \quad a\ge 0.
\ee

\cite[3.4.3]{Nahin}
\be
\int_0^1 x^a\ln^2 x dx=\frac{2}{(a+1)^3}.
\ee

\subsection{Logarithms of more complicated arguments}

\cite{SofoMath7}
\be
\int_0^1 \frac{\ln(1+x)}{1+x^2}dx=\frac{\pi}{8}\log 2.
\ee

\cite{SofoMath7}
\be
\int_0^1 \frac{\ln(1-x)}{1+x^2}dx=\frac{\pi}{8}\log 2 -G.
\ee

\cite[C2.5]{Nahin}
\be
\int_0^\infty \frac{\ln(1+x)}{x^{3/2}}dx=2\pi.
\ee

\cite[p. 192]{Nahin}
\be
\int_0^\infty \frac{(\ln x)^2}{1+x^2}dx = \frac{\pi^3}{8}.
\ee

\cite{XuFilo31}
Let $I(\alpha,m,k)\equiv \int_0^1 x^{\alpha-1}\ln^m x \ln^k(1-x) dx$,
\be
I(\alpha,0,0)=\frac{1}{\alpha},\quad I(\alpha,i,0) = (-)^i i!/\alpha^{i+1}
\ee
then
\begin{multline}
I(\alpha,m,k) = \sum_{i=0}^{m-1}\binom{m-1}{i}(m-i-1)!\frac{(-)^{m-i}}{\alpha^{m-i}}I(\alpha,i,k)
\\
+\sum_{j=0}^{k-1}\binom{k}{j}(-)^{m+k-j}(m+k-j-1)! H_\alpha^{(m+k-j)}I(\alpha,0,j)
-\sum_{j=0}^{k-1}\binom{k}{j}(-)^{m+k-j}(m+k-j-1)! \zeta(m+k_j)I(\alpha,0,j)
\\
+\sum_{i=1}^{m-1}\sum_{j=0}^{k-1}\binom{m-1}{i}\binom{k}{j}(-)^{m+k-i-j}(m+k-i-j-1)! H_\alpha^{(m+k-i-j)}I(\alpha,i,j)
\\
-\sum_{i=1}^{m-1}\sum_{j=0}^{k-1}\binom{m-1}{i}\binom{k}{j}(-)^{m+k-i-j}(m+k-i-j-1)! \zeta(m+k-i-j)I(\alpha,i,j).
\end{multline}

\cite{RutledgeAMM45}
In terms of the constant (\ref{eq.A4def})
we have
\be
\int_0^1 \frac{\log^2 u}{u}\log(1+u)du = \frac{7\pi^4}{360},
\ee
\be
\int_0^1 \frac{\log^2 u}{u}\log(1-u)du = -\frac{\pi^4}{45},
\ee
\be
\int_0^1 \frac{\log^2 u}{u}\log\frac{1+u}{1-u}du = \frac{\pi^4}{24},
\ee
\be
\int_0^1 \frac{1}{u}\log^{s-1}\frac{1+u}{1-u}du = 2\Gamma(s)[1+\frac{1}{3^s}+\frac{1}{5^s}+\cdots],
\ee
\be
\int_0^1 \frac{\log u}{u}\log(1+u)\log(1-u)du = -\frac12 A_4+\frac{11}{2880}\pi^4.
\ee
\be
\int_0^1 \frac{\log u}{u}\log^2(1+u)du = A_4-\frac{\pi^4}{288},
\ee
\be
\int_0^1 \frac{\log u}{u}\log^2(1-u)du = -\frac{\pi^4}{180},
\ee
\be
\int_0^1 \frac{\log u}{u}\log^2\frac{1+u}{1-u}du = 2A_4-\frac{\pi^4}{60},
\ee
\be
\int_0^1 \frac{1}{u}\log^3(1+u)du = \frac{3}{2}A_4-\frac{\pi^4}{960},
\ee
\be
\int_0^1 \frac{1}{u}\log^3(1-u)du = -\frac{\pi^4}{15},
\ee
\cite{RutledgeAMM45}
\be
\int_0^1 \frac{\log u}{u}\log(1+u)du = -\frac{3}{4}\zeta(3),
\ee
\be
\int_0^1 \frac{\log u}{u}\log(1-u)du = \zeta(3),
\ee
\be
\int_0^1 \frac{\log u}{u}\log\frac{1+u}{1-u}du = -\frac{7}{4}\zeta(3),
\ee
\be
\int_0^1 \frac{1}{u}\log^2(1+u)du = \frac{1}{4}\zeta(3),
\ee
\be
\int_0^1 \frac{1}{u}\log^2(1-u)du = 2\zeta(3).
\ee

\cite[(5.2.3)]{Nahin}
\be
\int_0^1 \frac{1}{x}\ln\frac{1+x}{1-x}=\frac{\pi^2}{4}.
\ee

\cite[(5.3.11)]{Nahin}
\be
\int_0^1 \frac{1}{x}\ln^2\frac{1-x}{1+x}=\frac72 \zeta(3).
\ee

\cite{KolbigMathComp64_449}
\be
2\int_0^{\pi/2}\ln|1-\sin x|dx
= -\pi\ln 2-4G.
\ee

\cite{BradleyCS2001,YangIJMEST23}
\be
-2\int_0^{\pi/4}\log(2\sin x)dx=G.
\ee
A factor 2 is missing in \cite{Adamchik}.

\cite{ChaudhuriAMM74}
\be
(-)^n \int_0^{\pi/2}[\log \sin^p\theta \cos^q\theta]^nd\theta
=\frac{\pi}{2}
\left|
\begin{array}{cccccc}
t_1 & -1 & 0 & 0 & \ldots & 0\\
t_2 & t_1 & -2 & 0 & \ldots & 0 \\
t_3 & t_2 & t_1 & -3 & \ldots & 0 \\
\ldots & \ldots & \ldots \\
t_{n-1} & t_{n-2} & t_{n-3} & t_{n-4} & \ldots & -(n-1) \\
t_2 & t_{n-1} & t_{n-2} & t_{n-3} & \ldots & t_1
\end{array}
\right|
\ee
where 
\be
t_r\equiv \sum_{m=1}^\infty \frac{p^r+q^r}{(2m+1)^r}
-\sum_{m=0}^\infty \frac{(p+q)^r}{(2m+2)^r}.
\ee

\cite{YangIJMEST23}
\be
\int_0^{3\pi/2}\log(2\sin \frac{x}{2})dx=G.
\ee

\cite{BradleyCS2001}
\be
\frac{1}{4}\int_0^{\pi/2} \log\frac{1+\cos x}{1-\cos x}dx=G.
\ee

\cite{BradleyCS2001}
\be
\frac{1}{4}\int_0^{\pi/2} \log\frac{1+\sin x}{1-\sin x}dx=G.
\ee

\cite{KolbigMathComp64_449}
\be
\int_0^{\pi/2} (\ln\tan x)^{2n}dx
=
\left(\pi/2\right)^{2n+1}|E_{2n}|.
\ee

\cite{KolbigMathComp64_449}
\be
\int_0^\infty \frac{\ln x dx}{(x+a)^2} =\frac{\ln a}{a},\, 0<a.
\ee

\cite{Adamchik,BradleyCS2001}
\be
2\int_0^{\pi/4}\log(2\cos x)dx=G.
\ee

\cite{Amdeberhanarxiv3663,MollSci14}\cite[A115252]{sloane}
\be
\int_{\pi/4}^{\pi/2}\ln\ln\tan x dx
=
\int_0^1\ln\ln\left(\frac{1}{x}\right)\frac{dx}{1+x^2}
=
\frac{\pi}{2}\ln \left(\frac{\Gamma(3/4)\sqrt{2\pi}}{\Gamma(1/4)}\right)
\ee

\cite{Vilceanu2009}
If $\chi_{-\Delta}$ is an odd primitve character ($\mod \Delta$)
then
\be
\int_0^1\frac{\sum_{n=1}^{\Delta-1} \chi_{-\Delta}(n)x^{n-1}}
{1-x^\Delta}\ln \ln \frac{1}{x}dx=
\left\{
\begin{array}{ll}
\frac{\pi}{\surd 3}\ln\frac{\sqrt{3} \Gamma^2(2/3)}{(2\pi)^{2/3}}, & \mathrm{if} \Delta =3;\\
\pi\ln\frac{\Gamma(3/4)}{\pi^{1/4}}, & \mathrm{if} \Delta =4;\\
\frac{\pi}{\sqrt \Delta}
[\ln 2\pi-\sum_{r=1}^{\Delta-1}\chi_{-\Delta}\ln\Gamma\frac{r}{\Delta}],
& \mathrm{if} \Delta >4.
\end{array}
\right.
\ee

\cite{MillerJCAM100}
\be
\frac{1}{4}\int_0^1 \frac{x^4-6x^2+1}{(1+x^2)^3}\log\log(1/x) dx
=
-(\gamma + \log 4)[\zeta(-1,1/4)-\zeta(-1,3/4)]
+\zeta'(-1,1/4)-\zeta'(-1,3/4).
\ee

\cite{Adamchik}\cite{ChoiJMAA231}
\be
-\int_0^1\frac{\log(x)}{x^2+1}dx=G.
\ee

\cite[p. 130]{Nahin}
\be
\int_0^\infty \frac{\ln x}{1+x^3}dx = -\frac{2\pi^2}{27}.
\ee

\cite[(5.5.12)]{Nahin}
\be
\int_0^\infty \frac{\ln(1+tx)}{1+x^2}dx
=\frac{\pi}{4}\ln(1+t^2)-\ln(t) \tan^{-1}t+G(t),\quad t\ge 0,
\ee
where a generalized Catalan's function is defined as
\be
G(x) = \left\{
\begin{array}{ll}
t-\frac{t^3}{3^2}+\frac{t^5}{5^2}-\frac{t^7}{7^2}+\ldots,  & 0\le t\le 1;\\
\frac{\pi}{2}\ln t+\frac{1}{t}-\frac{1}{3^2t^3}+\frac{1}{5^2t^5}-\frac{1}{7^2t^7}+\ldots,  & t\ge 1.
\end{array}
\right.
\ee

\cite{Vilceanu2009}
\be
\int_0^1 \frac{1}{1-x^2}(-\log x)^ldx = l!(1-\frac{1}{2^{1+l}})\zeta(l+1).
\ee

\cite{Vilceanu2009}
\be
\int_0^1 \frac{x}{1-x^2}(-\log x)^ldx = \frac{1}{2^{l+1}}l!\zeta(l+1).
\ee

\cite{Connonarxiv07II}
\be
\int_0^1 \frac{\log t}{1-(1-t)(1-u)}dt
= \frac{1}{1-u}\Li_2\left(-\frac{u-1}{u}\right)
\ee

\cite{Amdeberhanarxiv3663}
\be
\int_0^b \frac{\ln t}{(1+t)^{n+1}}dt
= \frac{1}{n}[1-(1+b)^{-n}]\ln b -
\frac{1}{n}\ln(1+b)-\frac{1}{n(1+b)^{n-1}}\sum_{j=1}^{n-1}
\frac{1}{j!}\binom{n-1}{j}|S_{j+1}^{(2)}|b^j
\ee

\cite{SofoJIS15}
\be
\int_0^1 \frac{\ln (1-x)}{x}((1-x)^p-1)dx
= \sum_{n=1}^p \frac{(-)^{n+1}\binom{p}{n}H_n^{(1)}}{n}.
\ee

\cite{BarbieriNC11A}
\be
\int_0^1 dy \log(1+y)(\frac{1}{y+1/x}-\frac{1}{y+x})
=
\Li_2(x)+\Li_2(-x)+\log x \log(1-x)-\frac12 \zeta(2).
\ee

\cite{BarbieriNC11A}
\be
\int_0^1 \frac{\log^2(x)}{x} =\infty
\ee

\cite{BarbieriNC11A}
\be
\int_0^1 \frac{\log^2(x)}{1+x} =\frac32 \zeta(3)
\ee

\cite{ReynoldsAIMS7}
\be
\int_0^\infty \frac{y^m\log^k(ay)}{1+\alpha y^n} dy
=\frac{1}{n^2} e^{\pi i (m+1)/n} (i/n)^{k-1} (2\pi)^{k+1} \alpha^{-(m+1)/n}
\Phi( e^{2\pi i(m+1)/n},-k,\frac{\pi -in \log a +i\log\alpha}{2\pi}).
\ee

\cite{BarbieriNC11A}
Table of $\int _0^1 f(x) = \alpha_1\zeta(3)+\alpha_2\zeta(2)\log 2+\alpha_3 \log^3 2$:
Free entries at the $\alpha$-coefficients indicate 0.

\begin{tabular}{lrrr}
$f(x)$ & $\alpha_1$ & $\alpha_2$ & $\alpha_3$ \\
\hline
$\log^2(x)/(1+x)$ & 3/2 \\
$\log^2(x)/(x-1)$ & -2 \\
$\log(x)\log(1+x)/x$ & -3/4 \\
$\log(x)\log(1+x)/(x+1)$ & -1/8 \\
$\log(x)\log(1+x)/(x-1)$ & -1 & 3/2 \\
$\log(x)\log(1-x)/x$ & 1 \\
$\log(x)\log(1-x)/(x+1)$ & 13/8 & -3/2 \\
$\log(x)\log(1-x)/(x-1)$ & -1 \\
$\log^2(1+x)/x$ & 1/4 \\
$\log^2(1+x)/(x+1)$ & & & 1/3 \\
$[\log^2(1+x)-\log^2 2]/(x-1)$ & -1/4 & 1 & -2/3 \\
$\log(1+x)\log(1-x)/x$ & -5/4\\
$\log(1+x)\log(1-x)/(x+1)$ & 1/8 & -1/2 & 1/3\\
$[\log(1+x)-\log 2]\log(1-x)/(x-1)$ & -7/8 & 1/2 & -1/6\\
$\log^2(1-x)/x$ & 2\\
$\log^2(1-x)/(x+1)$ & 7/4 & -1 & 1/3\\
$\Li_2(x)/x$ & 1\\
$\Li_2(x)/(x+1)$ & -5/8 & 1\\
$[\Li_2(x)-\zeta(2)]/(x-1)$ & 2\\
$\Li_2(-x)/x$ & -3/4\\
$\Li_2(-x)/(x+1)$ & 1/4 & -1/2\\
$[\Li_2(-x)+\frac12 \zeta(2)]/(x-1)$ & -5/8\\
\end{tabular}

\cite{BarbieriNC11A}
Table of $\int _0^1 f(x) = \alpha_1\zeta^2(2)+\alpha_2 a_4
+\alpha_3\zeta(3)\log 2+\alpha_4 \zeta(2)\log^2 2+\alpha_5\log^4 2$
Free entries at the $\alpha$-coefficients indicate 0,
and $a_n\equiv \Li_n(1/2)=\sum_{k=1}^\infty 1/(2^kk^n)$.

\begin{tabular}{lrrrrr}
$f(x)$ & $\alpha_1$ & $\alpha_2$ & $\alpha_3$ &$\alpha_4$&$\alpha_5$\\
\hline
$\log^3(x)/(1+x)$ & -21/10 \\
$\log^3(x)/(x-1)$ & 12/5\\
$\log^2(x)\log(1+x)/x$ & 7/10\\
$\log^2(x)\log(1+x)/(x+1)$ & -3/2 & 4 & 7/2 & -1 & 1/6\\
$\log^2(x)\log(1+x)/(x-1)$ & 19/20 & & -7/2\\
$\log^2(x)\log(1-x)/x$ & -4/5\\
$\log^2(x)\log(1-x)/(x+1)$ & 2/5 & -4 & & 1 & -1/6\\
$\log x\Li_2(-x)/x$ & 7/20\\
$\log x\Li_2(-x)/(x+1)$ & 13/8 & -4 & -7/2 & 1 & -1/6\\
\end{tabular}

many more $f(x)$ which are triple products of logarithms
divided through $x$, $x+1$ or $x-1$, or are prodcuts of logarithms
multiplied by $\Li_2$, or are $\Li_3$. Some cases where $f(x)$ is a quadruple product of logarithms are in \cite{BlumleinCPC180}.

K\"olbig \cite{KolbigJCAM18} provides extensive tables of 
\be
\int_0^1 \frac{x^n}{(1-x)^l}\ln^m x dx =
\ee
and
\be
\int_0^1 \frac{x^{n-1/2}}{(1-x)^l}\ln^m x dx =
\ee
\be
\int_0^1 \frac{x^n}{(1-x)^{l+1/2}}\ln^m x dx =
\ee
for small integer $m$, $n$, $l$.

\cite{AdamchikJCAM79}
\be
\frac{z^p}{(p-1)!}
\int_0^1\frac{t^{z-1}}{(1-t)^z}\log^{p-1}\frac{1}{t}dt
=
\,_{p+1}F_p\left(\begin{array}{c}z,z,\ldots,z\\
z+1,\ldots,z+1
\end{array}\mid 1\right)
.
\ee

\cite{AdamchikJCAM79}
\be
\frac{1}{\Gamma(1-z)\Gamma(p)}
\int_0^1\frac{t^{z-1}}{(1-t)^z}\log^{p-1}(t)dt
=
\left[\begin{array}{c}z\\p \end{array}\right].
\ee

\cite{Amdeberhanarxiv3663}
\be
\int_0^x \frac{\ln t}{(1+t^2)^{n+1}}dt
=\frac{\binom{2n}{n}}{2^{2n}}
\left[g_0(x)+p_n(x)\ln x-\sum_{k=0}^{n-1}\frac{\tan^{-1}x+p_k(x)}{2k+1}\right]
\ee
where
\be
g_0(x)\equiv \ln x \tan^{-1}x-\int_0^x\frac{\tan^{-1}t}{t}dt
.
\ee

\cite{ReynoldsAIMS6}
\begin{multline}
\int_0^\infty
\frac{x^{m-1}\log^k(ax)}{(x^{2u}+1)(x^{3u}+1)}dx
=
\frac{2^{k-1/2}\pi^{k+1}(i/u)^k e^{i\pi(2m-u)/(4u)}}{u}
\\ \times
[\Phi(e^{2im\pi/u},-k,\frac{\pi-2iu\log a}{4\pi})
-ie^{i\pi m/u}\Phi( e^{2im\pi/u},-k,3/4-\frac{iu\log a}{2\pi}
)]
\\
+\frac{e^{-i\pi/6}(2\pi)^{k+1}(i/u)^k e^{i\pi m/(3u)}}{3u}
[e^{4\pi m/(3u)}\Phi(e^{2im\pi/u},-k,5/6-\frac{iu\log a}{2\pi})
-e^{i\pi 3}\Phi( e^{2im\pi/u},-k,\frac{\pi-3iu\log a}{6\pi}
)]
\\
-\frac13 2^k\pi^{k+1}(i/u)^{k+1} e^{i\pi m/u}
\Phi(e^{2im\pi/u},-k,\frac{\pi-iu\log a}{2\pi}),
\end{multline}
where $\Phi$ is the Lerch function.

\cite{MollSci14}
\begin{multline}
\int_0^\infty \frac{\ln^{n-1}x dx}{(x-1)(x+a)}
=
\frac{(-)^n(n-1)!}{1+a}
\big\{
[1+(-)^n]\zeta(n)
\\
-\sum_{j=0}^{\lfloor n/2\rfloor}
\binom{n}{2j}(2^{2j}-2)(-)^jB_{2j}\pi^{2j}\log^{n-2j} a
\big\}
\end{multline}
for $n\ge 2$, $a>0$.

\cite{SofoJIPAM10}
\be
-2t\int_0^1 \frac{(1-x)^{j+1}\log(1-x)}{(1-tx(1-x))^3}dx= \sum_{n\ge 1}\frac{t^n}{C_n(j)}\sum_{r=1}^n\frac{1}{r+j+n}
\ee
where $C_n(j)=\binom{2n+j}{n}/(n+1)$ are Catalan related numbers.
For $t=2$ the integral becomes a sum of $G$, $\zeta(2)$, $\pi\ln 2$, $\pi$
and 1 with rational coefficients, and similar results are given for $t=1/2$.

\cite{BaileyJPA39}
\be
\int_0^\infty \frac{\log(1+x)}{1+x+x^2}dx
= -\int_0^1\frac{(1+t)\log t}{1+t^3}dt
= -\frac{3}{2}L_{-3}(2),
\label{eq.C3}
\ee
where $L_{-3}(s)$ is a Dirichlet series \cite[A086724]{sloane}.

\cite{Amdeberhanarxiv3663}
\be
\int_0^1 \frac{\ln t}{(1+t^2)^{n+1}}dt
=-2^{-2n}\binom{2n}{n}\left(G+\sum_{k=0}^{n-1}\frac{\frac{\pi}{4}+p_k(1)}{2k+1}\right)
\ee
where
\be
p_k(1)=\sum_{j=1}{k}\frac{2^j}{2j\binom{2j}{j}}.
\ee

\cite{Adamchik}
\be
\int_0^1 \left( \frac{2}{x^2-4x+8}-\frac{3}{x^2+2x+2}\right)\log x dx=C.
\ee

\cite{Connonarxiv07II,SofoJIS15}
\be
(-1)^{p+1}n\int_0^1 (1-t)^{n-1} \log^pt dt=p!\sum_{k=1}^n\binom{n}{k} \frac{(-1)^k}{k^p}
.
\ee

\cite{Connonarxiv07II}\cite[A152648]{sloane}
\be
\int_0^1 \frac{\log^2 t}{1-t}dt = \sum_{n=1}^\infty \frac{H_n^{(1)}}{n^2}
=2\zeta(3)
=2\Li_3(t)-2\Li_2(t)\log t-\log(1-t)\log^2t+c
\ee
where $H_n^{(r)}\equiv \sum_{k=1}^n\frac{1}{k^r}$.

\cite{Connonarxiv07II}
\be
\int_0^1 \frac{\log^3 u}{1-u}du
=-6\zeta(4).
\ee

\cite[(7.3.1)]{Nahin}
\be
\int_0^{\pi/2}x \ln(\sin x)dx = \frac{7}{16}\zeta(3)-\frac{\pi^2}{8}\ln 2.
\ee

\cite[(7.3.2)]{Nahin}
\be
\int_0^{\pi/2} \ln^2(a \sin x)dx = \frac{\pi^3}{24}+\frac{\pi}{2}\ln^2\frac{2}{a}.
\ee

\cite[(7.3.3)]{Nahin}
\be
\int_0^{\pi/2} \ln(a \sin x)\ln(a\cos x)dx = \frac{\pi}{2}\ln^2\frac{2}{a}-\frac{\pi^3}{48}.
\ee

\cite{MezoRJ46}
Let $S_n=(-)^n\int_0^1 \log^n(\sin(\pi x))dx$ then
$S_1=\log 2$,
\be
S_{n+1}=S_n\log 2+\sum_{j=1}^n(1-1/2^j)\zeta(j+1)\frac{n!}{(n-j)!}S_{n-j},\quad n\ge 1.
\ee

\cite{Vildanovarxiv1007}
\be
\int_{0}^{\pi/2}
\frac{\ln(2\cos x)}{x^2+\ln^2(2\cos x)}dx =\pi/4.
\ee

\cite{Vildanovarxiv1007}
\be
\int_{0}^{\pi/2}
\ln[x^2+\ln^2(2\cos x)]dx =0.
\ee

\cite{Vildanovarxiv1007}
\be
\int_{0}^{\pi/2}
\ln[x^2+\ln^2(2e^{-a}\cos x)]dx =x\ln\frac{a}{e^b-1},
\ee
and
\be
\int_{0}^{\pi/2}
\ln[x^2+\ln^2(2e^{-a}\cos x)]\cos 2x dx =
\frac{\pi}{2}\left(1-\frac{1}{a}-e^b+\frac{1}{e^b-1}\right)
\ee
where $b=\min(a,\ln 2)$.

\cite{Vildanovarxiv1007}
\be
\int_{0}^{\pi/2}
\ln[x^2+\ln^2(\cos x)]dx =\frac{\pi}{2}\ln\ln 2.
\ee

\cite{Vildanovarxiv1007}
\be
\int_{0}^{\pi/2}
\ln[x^2+\ln^2(\cos x)]\cos 2x dx =-\frac{\pi}{\ln 2}.
\ee

\cite{Vildanovarxiv1007}
\be
\int_{0}^{\pi/2}
\frac{\ln \cos x}{x^2+\ln^2(\cos x)} dx
=\frac{\pi}{2}\left(1-\frac{1}{\ln 2}\right).
\ee

\cite{Vildanovarxiv1007}
\be
\int_{0}^{\pi/2}
\frac{x\sin 2x}{x^2+\ln^2(\cos x)} dx
=\frac{\pi}{4\ln^2 2}.
\ee

\cite{Vildanovarxiv1007}
\be
\int_{0}^{\pi/2}
\frac{x\sin 2x}{x^2+\ln^2(2\cos x)} dx
=\frac{13\pi}{48}.
\ee

\cite{Vildanovarxiv1007}
\be
\int_{-\pi /2}^{\pi/2}
\frac{(1+e^{-2ix})^\beta}{\ln(1+e^{-2ix})-a} dx
=-\frac{\pi}{a}+\pi\frac{e^{(\beta+1)a}}{e^a-1}H(\ln 2 -a),
\ee
and
\be
\int_{0}^{\pi/2}
\frac{x\sin x}{x^2+\ln^2(2e^{-a}\cos x)} dx
=\frac{\pi}{4a^2}+\frac{\pi e^a}{4}
\left(1-\frac{1}{(e^a-1)^2}\right)
H(\ln 2 -a),
\ee
where $H$ is the unit step function.

\cite[(5.2.4)]{Nahin}
\be
\int_0^{\pi/2} \cot x \ln(\sec x) dx =\frac{\pi^2}{24}.
\ee

\cite{Amdeberhanarxiv3663,Vildanovarxiv1007}
\be
\frac{4}{\pi}\int_0^{\pi/2}\frac{x^2dx}{x^2+\ln^2(2\cos x)}
=
\frac{1}{2}(1+\ln(2\pi)-\gamma).
\ee

\cite{Amdeberhanarxiv3663}
\be
\frac{4}{\pi}\int_0^{\pi/2}\frac{x^2dx}{x^2+\ln^2(2e^{-a}\cos x)}
=
\frac{\gamma}{a}+\frac{a+\ln(1-e^{-a})-\gamma-\ln a}{1-e^{-a}}
+\frac{a}{1-e^{-a}}\int_0^1e^{-at}\ln \Gamma(t)dt
\ee
\be
=
\frac{\gamma}{a}+\frac{a+\ln(1-e^{-a})+\Gamma(0,a)}{1-e^{-a}}
+\frac{1}{1-e^{-a}}\int_0^1e^{-at}\psi(t+1) dt
\ee
where $0<a<\ln 2$.

\cite{Amdeberhanarxiv3663}
\be
\frac{4}{\pi}\int_0^{\pi/2}\frac{x^2dx}{x^2+\ln^2(2e^{-a}\cos x)}
=
\frac{\gamma}{a}
+\int_0^\infty e^{-at}\psi(t+1)dt
\ee
where $a>\ln 2$.

\cite{Amdeberhanarxiv3663}
\be
\int_0^{\pi/2}\frac{x^2\ln(2\cos x)dx}{(x^2+\ln^2(2\cos x))^2}
=
\frac{7\pi}{192}+\frac{\pi\ln 2\pi}{96}-\frac{\zeta'(2)}{16\pi}.
\ee

\cite{Vildanovarxiv1007}
\be
\frac{1}{2i}\int_{-\pi/2}^{\pi/2}
\frac{x(1+\exp(-2ix))^\beta}{\ln(1+\exp(-2ix))}dx
=
\frac{\pi}{8}
[1+\ln 2\pi -\gamma(2\beta+1)-2\ln\Gamma(\beta+1)]
\ee
with $\Re \beta>-1$.

\cite{Vildanovarxiv1007}
\be
\int_{-\pi/2}^{\pi/2}
\frac{(1+\exp(-2ix))^\beta}{\ln(1+\exp(-2ix))}dx
=
\frac{\pi}{2}
(1+2\beta)
.
\ee

\cite{ChoiJMAA231}
\be
\int_0^{\pi/2}
\log(1+a^2+2a\cos t)dt
= -\pi\log\frac{\cos(\arctan a)}{\pi}
-2\pi\log\frac{G(1+\frac{1}{\pi}\arctan a)G(\frac32-\frac{1}{\pi}\arctan a)}
{G(1-\frac{1}{\pi}\arctan a)G(\frac12+\frac{1}{\pi}\arctan a)}
,
\ee
where $G$ is the reciprocal of the double Gamma function.

\subsection{Logarithmic functions of compound arguments and powers}

\cite[2.4.2]{Nahin}
\be
\int_0^{\pi/2} \log\frac{\sin x}{x}dx = \frac{\pi}{2}[1-\log \pi].
\ee

\cite{LiMath10}
\be
\int_0^{\pi/4}y^m\ln\sin y dy = \Theta(m)+\Lambda(m).
\ee
\be
\int_0^\infty y^m\ln\tanh y dy = -\frac{1}{m+1}H(m+1,1)
\ee
where $\Theta$, $\Lambda$ and $H$ are defined in \eqref{eq.LiTheta},
\eqref{eq.LiLambda}, \eqref{eq.LiH}.

\cite{Adamchik}
\be
-\int_0^{1}\frac{\log(\frac{1}{\sqrt{2}}(1-x))}{x^2+1} = G.
\ee

\cite{Adamchik}
\be
-\int_0^{1}\frac{\log[\frac{1}{2}(1-x^2)]}{x^2+1} = G.
\ee

\cite[(C5.9)]{Nahin}
\be
\int_0^\infty \ln \frac{e^x+1}{e^x-1}dx = \pi^2/4.
\ee

\cite{Vilceanu2009}
\be
\int_0^1 \frac{(-\log(1\pm x))^k}{1\pm x}x^r(-\log x)^ldx
=\mp k!l!\sum_{n_1>n_2>n_3>\ldots >n_{k+1}\ge 1}
\frac{(\mp 1)^{n_1}}{(n_1+r)^{l+1}n_2n_3\cdots n_{k+1}}.
\ee

\cite{Vilceanu2009}
For integers $k$, $r>0$, $l\ge 1$
\be
\int_0^1 (-\log(1\pm x))^kx^r(-\log x)^ldx
= k!l!\sum_{n_1>n_2>n_3>\ldots >n_{k}\ge 1}
\frac{(\mp 1)^{n_1}}{(n_1+r+1)^{l+1}n_2n_3\cdots n_k}.
\ee

\cite{Vilceanu2009}
For integers $k$, $l$, $r\ge 0$, $m$, $l\ge 1$
\be
I_{k,l,r,m}^{-}\equiv \int_0^1 \frac{(-\log(1- x))^k}{(1-x)^m}x^r(-\log x)^ldx
= \sum_{i=0}^r(-)^i\binom{r}{i}
I_{l,k,i-m,0}^{-}.
\ee

\cite[(C5.15)]{Nahin}
\be
\int_0^2 \frac{\ln(1-x)\ln(1+x)}{x}dx
=-\frac{5}{8}\zeta(3).
\ee

\cite{Vilceanu2009}
\be
-\int_0^1 \log(1- x)(-\log x)^ldx
= l!\sum_{n_1\ge 1}\frac{1}{(n_1+1)^{l+1}n_1}
=l!\left(l+1-\sum_{i=1}^{l+1}\zeta(i)\right),\quad l\ge 1.
\ee

\cite{Vilceanu2009}
\be
-\int_0^1 \log(1+ x)(-\log x)^ldx
= -l!\sum_{n_1\ge 1}\frac{(-)^{n_1}}{(n_1+1)^{l+1}n_1}, \quad l\ge 0.
\ee

\cite{ReynoldsAIMS6}
\be
\int_0^\infty \frac{\log x \log(x/a)}{(x^2+1)(x^3+1)}dx
=\frac{\pi^2}{432}(27\pi -37\log\frac{1}{a}).
\ee

\cite{ReynoldsAIMS6}
\be
\int_0^\infty \frac{x\log x \log(x/a)}{(x^2+1)(x^3+1)}dx
=\frac{\pi^2}{3888}[(160\surd 3-243)\pi -45\log\frac{1}{a}].
\ee

\cite{ReynoldsAIMS6}
\begin{multline}
\int_0^\infty \frac{\log x }{\surd x (x^{2u}+1)(x^{3u}+1)}dx
=-\frac{\pi^2\csc^2\frac{\pi}{2u}}
{72u^2}
[
20\cos\frac{\pi}{6u}
+27\cos\frac{\pi}{4u}
+12\cos\frac{\pi}{2u}
+9\cos\frac{3\pi}{4u}
\\
+4(
-9\sin^3\frac{\pi}{4u}
+5\surd 3\sin\frac{\pi}{6u}
-\surd 3\sin\frac{5\pi}{6u}
+\cos\frac{5\pi}{6u}
)
]
.
\end{multline}

\cite{BarbieriNC11A}
\be
\int \frac{\log(a+by)}{c+ey}dy=\frac{1}{e}
[
\log\frac{ae-bc}{e}\log(c+ey)-\Li_2(-b\frac{c+ey}{ae-bc})
]
\ee

\cite{BarbieriNC11A,KolbigJSC1}
Let
\be
S_{n,p}(y)\equiv \frac{(-1)^{n+p-1}}{(n-1)!p!}
\int_0^2 \frac{\log^{n-1}t \log^p(1-yt)}{t}dt
\ee
then
\be
S_{1,1}(y)=\Li_2(y)=-\int_0^1 \frac{\log(1-yt)}{t}dt
=-\int_0^y \frac{\log(1-t)}{t}dt.
\ee
\be
S_{2,1}(y)=\Li_3(y)=-\int_0^1 \frac{\log t \log(1-yt)}{t}dt
=\int_0^y \frac{\Li_2(t)}{t}dt.
\ee
\be
S_{1,2}(y)=\frac12 \int_0^1 \frac{\log^2(1-yt)}{t}dt
=\frac12 \int_0^y \frac{\log^2(1-t)}{t}dt.
\ee
\be
\int_0^y \frac{dt}{t} S_{n,p}(t) = S_{n+1,p}(y);
S_{0,p}\equiv \frac{(-1)^p}{p!}\log^p(1-y).
\ee
\be
S_{n,p}(1-y)=\sum_{s=0}^{n-1}\frac{\log^s(1-y)}{s!}
[
S_{n-s,p}(1)-\sum_{r=0}^{p-1}\frac{(-1)^r}{r!}\log^r y S_{p-r,n-s}(y)
]
+\frac{(-1)^p}{n!p!}\log^n(1-y)\log^p y;
\ee
\begin{multline}
S_{n,p}(-\frac{1}{y})=
(-)^n\sum_{s=0}^{p-1}(-)^s\sum_{r=0}^s\frac{(-)^r}{r!}
\log^r y\binom{n+s-r-1}{s-r}S_{n+s-r,p-s}(-y)
\\
+
(-)^p[\sum_{r=0}^{n-1}\frac{(-1)^r}{r!}\log^r y C_{n-r,p}
+\frac{(-)^{n+p}}{(n+p)!}\log y],
\end{multline}
where 
\be
C_{n,p}\equiv (-)^p(1-(-)^n)S_{n,p}(-1)-(-)^n\sum_{r=1}^{p-1}(-)^{p-r}S_{n+r,p-r}(-1).
\ee
\cite{BlumleinCPC180}
\be
S_{n,p}(x^2) = 2^n\frac{(-)^{n+p-1}}{(n-1)!p!} \sum_{l=0}^p \binom{p}{l}\int_0^x \frac{dz}{z} \ln^{n-1}(\frac{z}{x})
\ln^{p-l}(1-z)\ln^l(1+z).
\ee
\cite{BlumleinCPC180}
\be
\frac12 S_{1,p}(x^2) = S_{1,p}(x)+S_{1,p}(-x)+ \frac{(-)^p}{p!} \sum_{l=0}^{p-1} \binom{p}{l}\int_0^x \frac{dz}{z}
\ln^{p-l}(1-zx)\ln^l(1+zx).
\ee

\cite{GeddesCM}
\be
\int_0^\infty t^w \ln^m (t) \exp(-t^s) dt = sig(s) s^{-m-1}[\frac{d^m}{dz^m}\Gamma(z)]_{\mid{z=(w+1)/s}}.
\ee

\cite{GeddesCM}
\be
\int_0^\infty t^w \ln^m (t) \frac{\exp(-t^s)}{1-\exp(-t^s)} dt = 
\sgn(s) s^{-m-1}[\frac{d^m}{dz^m}(\Gamma(z)\zeta(z))]_{\mid{z=(w+1)/s}}.
\ee

\cite{KolbigJSC1}
\begin{multline}
\int_0^1 t^{-1/2}\ln^n t\ln^p(1-t) dt
=
2^{n+p+1}\int_0^{\pi/2}\ln^n \sin t \ln^p \cos t dt
\\
=
\pi n! p!\sum_{\nu=0}^n C_{n-\nu}(3;2,0,1)
\sum_{\rho=0}^p \binom{\nu+\rho}{\rho}C_{p-\rho}(4;2,0,1)
C_{\nu+\rho}(1;0,0,1),\quad ,n,p=0,1,2,\ldots
\end{multline}
where the $C$ are given as
\be
C_0(a;\alpha_1,\alpha_2,\alpha_3)=1,
\ee
\be
C_k(a;\alpha_1,\alpha_2,\alpha_3)=\frac{1}{k}\sum_{\kappa=}^k
c_k(a,\alpha_1,\alpha_2,\alpha_3)C_{k-\kappa}(a,\alpha_1,\alpha_2,\alpha_3),
\quad k\ge 1
\ee
\cite{KolbigJSC1}
\begin{multline}
\int_0^\infty e^{-\mu t}t^n\ln^m t
=
2^{m+1}\int_0^\infty e^{\mu t^2}t^{2n+}\ln^m t dt
\\
=
(-)^n\mu^{-n-1}m!\sum_{\rho=0}^{\min(m,n)}
(-)^\rho S_{n+1}^{(\rho+1)}C_{m-\rho}(\mu;1,0,0),\quad \Re \mu>0, n,m=0,1,2\ldots
\end{multline}
where $S_k^{(j)}$ are the Stirling Number of the First Kind.

\cite{KolbigJSC1}
\begin{multline}
\int_0^\infty e^{-\mu t}t^{n-1/2} \ln^m t
=2^{m+1}\int_0^\infty e^{-\mu t^2}t^{2n}\ln^m t
\\
=\sqrt{\pi}(-)^n \mu^{-n-1/2}m!\sum_{\rho=0}^{\min(m,n)}
\tilde S(\rho,n)C_{m-\rho}(4\mu;2,0,1),\quad \Re \mu>0, n,m=0,1,2\ldots
\end{multline}
where $\tilde S(p,q)=\sum_{j=p}^q (-)^j S_q^{(j)}\binom{j}{p}2^{p-j}$.

\cite{KolbigJSC1}
\begin{multline}
\int_0^\infty \frac{x^n \ln^m x}{(1+\beta x)^l}dx
=\frac{(-)^{l+m} m!}{\beta^{n+1}(l-1)!}
\sum_{m_1=0}^{\min(m,n)} S_{n+1}^{(m_1+1)}
\sum_{m_2=0}^{\min(m-m_1,l-n-2)} (-)^{m_2} S_{l-n-1}^{(m_2+1)}
\\
\times
\sum_{m_3=0}^{m-m_1-m_2} \frac{|2^{m_3}-2| |B_{m_3}|}{m_3!(m-m_1-m_2-m_3)!}
\pi^{m_3} \ln^{m-m_1-m_2-m_3}\beta,
\end{multline}
for $|\arg\beta|<\pi, n\ge 0, l\ge n+2, m\ge 0$.

\cite{KolbigJSC1}
\begin{multline}
\int_0^\infty \frac{x^n \ln^m x}{(1-x)^l}dx
=\frac{(-)^{m+n+l} m!}{(l-1)!}
\sum_{m_1=0}^{n} S_{n+1}^{(m_1+1)}
\sum_{m_2=0}^{l-n-2} (-)^{m_2} S_{l-n-1}^{(m_2+1)}
\\
\times
S_{l-n-1}^{(m_2+1)} \frac{(2\pi)^{m-m_1-m_2} |B_{m-m_1-m_2}|}{(m-m_1-m_2)!}
\pi^{m_3} \ln^{m-m_1-m_2-m_3}\beta,
\end{multline}
for $n\ge 0, l\ge n+2, m\ge l-1$.
The integral to be understood as the PV if $m=l-1$.

\cite{KolbigJSC1}
\begin{multline}
\int_0^\infty \frac{x^{n-1/2} \ln^m x}{(1+\beta x)^l}dx
=\frac{(-)^{l+m+1} m!}{\beta^{n+1/2}(l-1)!}
\sum_{m_1=0}^{\min(m,n)} (-)^{m_1} \tilde S(m_1,n)
\sum_{m_2=0}^{\min(m-m_1,l-n-1)} \tilde S(m_1,l-n-1)
\\
\times
\sum_{m_3=0}^{m-m_1-m_2} \frac{|E_{m_3}|\pi^{m_3+1}}{m_3!(m-m_1-m_2-m_3)!}
\ln^{m-m_1-m_2-m_3}\beta,
\end{multline}
for $|\arg\beta|<\pi, n\ge 0, l\ge n+1, m\ge 0$, where $E$ are Euler numbers.

\cite{KolbigJSC1}
\begin{multline}
\int_0^\infty \frac{x^{n-1/2} \ln^m x}{(1-x)^l}dx
=\frac{(-)^{m+n+l+1} m!}{(l-1)!}
\sum_{m_1=0}^{n} (-)^{m_1} \tilde S(m_1,n)
\sum_{m_2=0}^{l-n-1}\tilde S(m_2,l-n-1) (2\pi)^{m-m_1-m_2+1}
\\
\times
\frac{(2^{m-m_1-m_2+1}-1) |B_{m-m_1-m_2+1}|}{(m-m_1-m_2+1)!}
\end{multline}
for $n\ge 0, l\ge n+1, m\ge l-1$.
The integral to be understood as the PV if $m=l-1$.

\cite{KolbigJSC1}
\begin{multline}
\int_0^\infty \frac{x^{n} \ln^m x}{(1+\beta x)^{l-1/2}}dx
=\frac{(-)^l2^{l-1} m!}{\beta^{n+1}(2l-3)!!}
\sum_{m_1=0}^{\min(m,n)} (-)^{m_1} S_{n+1}^{(m_1+1)}
\\
\times
\sum_{m_2=0}^{\min(m-m_1,l-n-2)}(-)^{m_2} \tilde S(m_2,l-n-2)
C_{m-m_1-m_2}(\beta/4;1,-2,-1),
\end{multline}
for $|\arg\beta|<\pi, n\ge 0, l\ge n+2, m\ge 0$. The previous integral
also evaluates
\be
\int_0^\infty \frac{x^{n-1/2}\ln^m x}{(1+\beta x)^{l-1/2}}dx
=(-)^m \beta^{1/2-l}\int_0^\infty \frac{x^{l-n-2}\ln^m x}{(1+x/\beta)^{l-1/2}}dx.
\ee

\cite{KolbigJCAM18}
\be
\int_0^1 \frac{x^{\nu-1}}{(1-x)^\lambda}\ln^m x dx = 
(\partial/\partial \nu)^m B(1-\lambda,\nu)
=\Gamma(1-\lambda)(\partial/\partial \nu)^m \frac{\Gamma(\nu)}{\Gamma(1+\nu-\lambda)}.
\ee

\cite{KolbigJCAM18}
\be
\int_0^1 \frac{x^{\nu-1}}{(1-x)^l}\ln^m x dx = 
\frac{(-)^l}{(l-1)!}(d/d\nu)^m \psi(\nu)\prod_{j=1}^{l-1}(\nu-j).
\ee

\cite{KolbigJCAM18}
\be
\int_0^\infty \frac{X^{\nu-1}}{(X^\mu+z^\mu)^\lambda}\ln^m X dx = \frac{1}{\mu}
(\partial/\partial \nu)^m z^{\nu-\lambda \mu}B(\frac{\nu}{\mu},\lambda-\frac{\nu}{\mu}).
\ee

\cite{BorosSIR40}
\be
\int_0^\infty \left[\frac{1}{x^4-x^2+1}\right]^r
\ln\frac{x^2}{x^4-x^2+1}
\,
\frac{ dx}{x^2}
=
\frac{\surd\pi}{2}\,\frac{\Gamma(r)\Gamma'(r-1/2)-\Gamma(r-1/2)\Gamma'(r)}{\Gamma^2(r)}
\ee

\cite{BorosSIR40}
\be
\int_0^\infty \frac{1}{x^4-x^2+1}
\ln^2\frac{x^2}{x^4-x^2+1}
\,
\frac{ dx}{x^2}
=
\frac{\pi}{2}(\frac{\pi^2}{3}+4\ln^2 2).
\ee

\cite{BorosSIR40}
\be
\int_0^\infty \left[\frac{x}{x^2+1}\right]^{2r}
\ln\frac{x}{x^2+1}
\,
\frac{ dx}{x^2}
=
\frac{\sqrt{\pi}}{4}G(r)[\psi(r-1/2)-\psi(r)-2\ln 2],
\ee
where $G(r)\equiv 2^{1-2r}B(r-1/2,1/2)/\surd \pi$.

\cite{ChoiEAJM29}
\be
\int_0^1 x^{-1/2} (1-x)^{2a+1}(4-x)^a (\ln\frac{(1-x)^2}{4-x})^ndx
=
\frac{2\surd \pi}{3}D_a^n[\frac{2^{2a}\Gamma(a+1)}{\Gamma(a+3/2)}],\quad \Re a>-1, n\in N.
\ee

\cite{ChoiEAJM29}
\be
\int_0^1 x^{-1/2} (1-x)^{2a+1}(4-x)^a \ln\frac{(1-x)^2}{4-x}dx
=
\frac{\surd \pi}{3}\frac{2^{2a}\Gamma(a+1)}{\Gamma(a+3/2)}[2\ln 2+\psi(a+1)-\psi(a+3/2)],\quad \Re a>-1.
\ee

\cite{ChoiEAJM29}
\begin{multline}
\int_0^1 x^{-1/2} (1-x)^{2a+1}(4-x)^a (\ln\frac{(1-x)^2}{4-x})^2dx
=
\frac{\surd \pi}{3}\frac{2^{2a}\Gamma(a+1)}{\Gamma(a+3/2)}
\\ \times
[\{2\ln 2+\psi(a+1)-\psi(a+3/2)\}^2 +\psi'(a+1)-\psi'(a+3/2)],\quad \Re a>-1.
\end{multline}

\cite{BorosSIR40}
\be
\int_0^\infty \left[\frac{x}{x^2+1}\right]^{2r}
\ln\frac{x}{x^2+1}
\,
\frac{x^2+1}{x^2(x^b+1)}dx
=
\frac{\sqrt{\pi}}{4}G'(r)
\ee
where $G\equiv 2^{1-2r}B(r-1/2,1/2)/\surd \pi$.

\cite{BorosSIR40}
\be
\int_0^\infty \left[\frac{x}{x^2+1}\right]^{2r}
\ln^2 \frac{x}{x^2+1}
\,
\frac{dx}{x^2}
=
\frac{\sqrt{\pi}}{8}G''(r)
\ee
where $G''\equiv G(r)[\psi'(r-1/2)-\psi'(r)+(\psi(r-1/2)-\psi(r)-2\ln 2)^2]$

\cite{Amdeberhanarxiv3663}
\be
\int_0^\infty e^{-ax}\ln x\, dx = -\frac{\gamma+\ln a}{a}.
\ee

\cite{Amdeberhanarxiv3663}
\be
\int_0^\infty \frac{\ln x}{e^x+e^{-x}-1}dx
=
\int_0^1\ln\ln\left(\frac{1}{x}\right)\frac{dx}{1-x+x^2}
=
\frac{2\pi}{\sqrt{3}}\left(\frac{5}{6}\ln 2\pi-\ln\Gamma\left(\frac{1}{6}\right)\right)
.
\ee

\cite{KolbigJCAM69}
\be
\frac{a^\nu}{\Gamma(\nu)}
\int_0^\infty
e^{-ax}x^{\nu-1}\ln^m x dx
=
\phi^m(a,\nu)+\sum_{j=1}^m\binom{m}{j}\eta_j\phi^{m-j}(a,\nu),
\ee
where
\be
\phi(a,\nu)\equiv\psi(\nu)-\ln a,
\ee
\be
\eta_j(\nu)\equiv (-1)^j\sum_{\pi_0(j)}(j;0,k2,\ldots,k_j)^*\zeta^{k_2}(2,\nu)
\cdots \zeta^{k_j}(j,\nu),
\ee
\be
\zeta(k,\nu)\equiv \sum_{l=0}^\infty\frac{1}{(l+\nu)^k},\quad k\ge 2.
\ee

\cite{Amdeberhanarxiv2379}
\begin{multline*}
\int_0^\infty \log x\cos^{2n+1} x^2 dx =
-\frac{\sqrt{\pi}}{2^{2n+3}}(\pi+2\gamma+4\log 2)\sum_{k=0}^n\binom{2n+1}{n-k}\frac{1}{\sqrt{4k+2}}
\\
-\frac{\sqrt{\pi}}{2^{2n+2}}\sum_{k=0}^n\binom{2n+1}{n-k}\frac{\log(2k+1)}{\sqrt{4k+2}}
.
\end{multline*}

\cite[3.4.9]{Nahin}
\be
\int_0^\pi \ln(a+b \cos x) dx= \pi \ln\frac{a+\sqrt{a^2-b^2}}{2}, \quad a>b.
\ee

\cite[3.4.10]{Nahin}
\be
\int_0^\pi \frac{\ln(1+b \cos x)}{\cos x} dx= \pi \sin^{-1} b.
\ee

\cite{Amdeberhanarxiv3663}
\be
\int_0^{\pi/2}x\ln(2\cos x)dx=-\frac{7}{16}\zeta(3).
\ee

\cite{Amdeberhanarxiv3663}
\be
\int_0^{\pi/2}x^2\ln(2\cos x)dx=-\frac{\pi}{4}\zeta(3).
\ee

\cite{Amdeberhanarxiv3663}
\be
\int_0^{\pi/2}x^2\ln^2(2\cos x)dx=\frac{11\pi}{16}\zeta(4)
=\frac{11\pi^5}{1440}.
\ee

\cite{CoffeyJCAM159}
\be
\int_0^{\pi/2}x^4\ln^2(2\cos x)dx=\frac{5\pi^7}{8064}+\frac{3\pi}{4}\zeta^2(3)
.
\ee

\cite{CoffeyJCAM159}
\be
\int_0^{\pi/2}x^2\ln^4(2\cos x)dx=\frac{33\pi^7}{4480}+\frac{3\pi}{2}\zeta^2(3)
.
\ee

\cite{DolderJCAM11,CoffeyJCAM183}
\be
\Cl_2(\theta)=-\int_0^\theta \ln(2\sin\frac{t}{2})dt
= -\sin\theta\int_0^2 \frac{\ln \theta}{\rho^2-2\rho\cos\theta+1}d\rho.
\ee

\cite{DolderJCAM11}
\be
\Cl_2(\pi/3)
= \frac{\surd 3}{6}[\psi'(\frac{1}{3})-\frac{2}{3}\pi^2].
\ee

\cite{DolderJCAM11}
\be
\Cl_2(p \pi/q)
= -\frac{1}{4q^2} \sum_{k=1}^{q-1}
[\psi'(1-\frac{k}{2q})-\psi'(\frac{1}{2}-\frac{k}{2q})]\sin k\frac{p}{q}\pi.
\ee
for $p$ odd.

\cite{DolderJCAM11}
\be
\Cl_2(p \pi/q)
= -\frac{1}{4q^2} \sum_{k=1}^{q-1}
[\psi'(1-\frac{k}{2q})+\psi'(\frac{1}{2}-\frac{k}{2q})]\sin k\frac{p}{q}\pi.
\ee
for $p$ even and $q\ge 3$ odd.

\cite{CoppoJNT150}
Let $L_k(\alpha)=\pi^2\int_0^\alpha x \ln^{k-1}(2\sin(\pi x/2))dx$, then
\be
L_1(\alpha)=(\pi\alpha)^2/2.
\ee
\be
L_2(\alpha)=\zeta(3)
-\sum_{n=1}^\infty \frac{\cos(\pi n\alpha)}{n^3}
-\alpha\pi \sum_{n=1}^\infty \frac{\sin(\pi n\alpha)}{n^2}
.
\ee
\be
L_3(\alpha)=\frac{\pi^4}{16}
(\alpha^4-8\alpha^3/3+2\alpha^2)
+2\sum_{n=1}^\infty \frac{H_n}{(n+1)^3}\cos(\pi(n+1)\alpha)
+2\pi\alpha \sum_{n=1}^\infty \frac{H_n}{(n+1)^2}\sin(\pi(n+1)\alpha)
-\frac{1}{2}\zeta(4)
.
\ee

\cite{DolderJCAM11}
\be
\mathrm{Ti}_2(\tan \theta)=
\theta \ln(\tan \theta)+\frac{1}{2}\Cl_2(2\theta)+\frac{1}{2}\Cl_2(\pi-2\theta),\quad 0\le \theta<\pi/2.
\ee
where
\be
\mathrm{Ti}_2(x)=\int_0^x \frac{\arctan t}{t}dt.
\ee

\cite{CoffeyJCAM159}
\begin{multline}
\int_0^{\pi/2}x^3\ln(\cos x)\sin[(p-1)x]\cos^{p-1}xdx=
-\frac{\pi}{15}2^{-(p+6)}[60\gamma^4 -60\gamma^2\pi^2+\pi^4
+60\gamma(\pi^2-2\gamma^2)\ln 2
\\
+60\psi'''(p)+60\{-(\pi^2+6\gamma(-\gamma+\ln 2))\psi^2(p)
+(4\gamma-2\ln 2)\psi^3(p)+\psi^4(p)+6\gamma \ln 2 \psi'(p)
\\
-3[\psi'(p)]^2-2(\gamma+\ln 2)\psi''(p)+(8\gamma-4\ln 2)\zeta(3)
+\psi(p)(4\gamma^3-2\gamma\pi^2+(\pi^2-6\gamma^2)\ln 2
\\
+6\ln 2\psi'(p)-2\psi''(p)+8\zeta(3))\}]
.
\end{multline}

\cite{Adamchik,BradleyCS2001}
\be
\int_0^{\pi/4}\log(\cot x)dx=
-\int_0^{\pi/4} \log(\tan x)dx
=G.
\ee

\cite{EspinosaRJ6,MedinaRJ20}
\be
\int_{\pi/4}^{\pi/2}
\ln\ln\tan x \, dx=\frac{\pi}{2}
\ln\frac{\Gamma(3/4)\sqrt{2\pi}}{\Gamma(1/4)}
.
\ee

\cite{MedinaRJ20,AdamchikICSAC}
\be
\int_0^1 x^j\log \log\frac{1}{x}dx
=
-\frac{\gamma+\log(j+1)}{j+1}.
\ee

\cite{AdamchikICSAC}
\be
\int_0^1\frac{x^{p-1}}{1+x^n}\log\log\frac{1}{x}dx
=
\frac{1}{2n}[\log(2n)+\gamma]
\left(
\psi(\frac{p}{2n})-\psi(\frac{n+p}{2n})
\right)
+
\frac{1}{2n}
\left(
\zeta'(1,\frac{p}{2n})-\zeta'(1,\frac{n+p}{2n})
\right)
\ee
for $\Re p>0$ and $\Re n>0$,
and a similar expression if $(1+x^n)^2$ or $(1+x^n)^3$ are in the denominator.

\cite{ReynoldsAIMS7}
\be
\int_0^\infty \frac{\log \log y}{1+y^2}dy
=\frac{\pi}{4}(\pi i +\log\frac{4\pi^2\Gamma^4(3/4)}{\Gamma^4(1/4)}).
\ee

\cite{AdamchikICSAC}
\be
\int_0^1\frac{x^{nr-1}}{1+x^n}\log\log\frac{1}{x}dx
=
\frac{1}{2n}[\log(2n)+\gamma]
\left(
\psi(r/2)-\psi(\frac{r+1}{2})
\right)
+
\frac{1}{2n}
\left(
\zeta'(1,r/2)-\zeta'(1,(r+1)/2)
\right)
.
\ee

\cite{AdamchikICSAC}
\be
\int_0^1 x^{p-1}\frac{1-x}{1-x^n}\log\log\frac{1}{x}dx
=
\frac{1}{n}[\log(n)+\gamma]
\left(
\psi(p/n)-\psi(\frac{p+1}{n})
\right)
+
\frac{1}{n}
\left(
\zeta'(1,p/n)-\zeta'(1,(p+1)/n)
\right)
.
\ee

\cite{AdamchikICSAC}
\be
\int_0^1 \frac{1}{1-x+x^2}\log\log\frac{1}{x}dx
=
\frac{\pi}{3\sqrt{3}}
[5\log(2\pi)-6\log\Gamma(1/6)]
.
\ee

\cite{MedinaRJ20}
Let
\be
R_{m,j}(a)\equiv \int_0^1 \frac{x^j\log \log 1/x}{(x+a)^{m+1}}dx,
\ee
and (Eulerian numbers)
\be
A_{m,j}=\sum_{k=0}^j(-)^k\binom{m+1}{k}(j-k)^m
\ee
and
\be
E_m\equiv \int_0^1 \frac{T_{m-1}(x)\log \log 1/x}{(x+1)^{m+1}}dx
\ee
defined via polynomials
\be
T_m(x)\equiv \sum_{j=0}^m(-)^jA_{m+1,j+1}x^j,
\ee
then
\be
E_m=(1-2^m)\zeta'(1-m)-(-)^m[\gamma(2^m-1)+2^m\log 2]\frac{B_m}{m}
\ee
where $B_m$ are the Bermoulli numbers. The $R_{m,j}$ are then recursively
\be
R_{0,0}(1)=-(\log^2 2)/2;\quad R_{m,0}(1)=\frac{E_m}{b_0(m)}-\sum_{k=1}^{m-1}\frac{b_k(m)}{b_0(m)}R_{m-k,0}(1),
\ee
\be
R_{0,0}(a)=-\gamma\log(1+1/a)-\Li'_1(-1/a);\quad \Li'_c(x) = -\sum_{n>=1}\frac{\log n}{n^c}x^n,
\ee
and for $m>0$
\be
R_{m,0}(a)=-\frac{\gamma}{a^m(1+a)m}-\frac{\gamma}{a^{m+1}m!}
\sum_{j=2}^m\frac{S_1(m,j)T_{j-2}(1/a)}{(1+1/a)^j}
-\frac{1}{a^mm!}\sum_{j=1}^mS_1(m,j)\Li'_{1-j}(-1/a)
,
\ee
where
\be
b_k(m)\equiv (-)^k\sum_{j=0}^{m-1}\binom{j}{k}A_{m,j+1},
\ee
and the unsigned Stirling numbers of the first kind are $S_1(m,j)$ as in
$(t)_m=\sum_{j=1}^m S_j(m,j)t^j$. For larger parameters $m$ then
\be
R_{m,0}(a)=\sum_{j=0}^r \alpha_{j,r}(a)R_{m-r+j,j}(a),\quad \alpha_{j,r}(a)\equiv (-)^j\binom{r}{j}a^{-r}.
\ee

\cite{MedinaRJ20}
Let
\be
D_{m,j}(r,\theta)=\equiv \int_0^1\frac{x^j\log\log 1/x}{(x^2-2rx\cos\theta+r^2)^{m+1}}dx,
\ee
then
\be
D_{0,0}(1,\theta)=\frac{\pi}{2\sin\theta}\left[(1-\theta/\pi)\log 2\pi+\log\frac{\Gamma(1-\theta/2\pi)}{\Gamma(\theta/2\pi)}\right]
,
\ee
\be
D_{0,0}(r,\theta)=-\frac{\gamma}{r\sin\theta}\tan^{-1}\frac{\sin\theta}{r-\cos\theta}
+\frac{1}{2ri\sin\theta}(\Li'_{1}(e^{i\theta}/r)-\Li'_{1}(e^{-i\theta}/r))
,
\ee
\be
D_{0,1}(r,\theta)=-\frac{\gamma}{2}\log\frac{r^2-2r\cos\theta+1}{r^2}-\gamma\cot\theta
\tan^{-1}\frac{\sin\theta}{r-\cos\theta}
+\frac{1}{2ri\sin\theta}[\Phi'(e^{i\theta}/r,1,1)-\Phi'(e^{-i\theta}/r,1,1)]
,
\ee
\be
D_{m,j}(r\theta)=-\frac{1}{2rm\sin\theta}\frac{\partial}{\partial \theta}D_{m-1,j-1}(r,\theta),\quad
m,j>0.
\ee

\cite{MedinaRJ20}
\be
\int_0^1\frac{\log(1-x)}{x}\log\log 1/x dx=\int_0^\infty\log t\log(1-e^{-t})dt=\frac{\gamma\pi^2}{6}-\zeta'(2).
\ee

\cite{MedinaRJ20}
\be
\int_0^1\frac{\log(1+x)}{x}\log\log 1/x dx=\frac{\pi^2}{12}(\log 2-\gamma)+\zeta'(2)/2,
\ee
and other examples involving the kernel $\log\log 1/x$.

\cite{EspinosaRJ6}
\be
\int_0^1 q^n\ln(\sin\pi q)dq
=
-\frac{\ln 2}{n+1}+n!
\sum_{k=1}^{\lfloor n/2\rfloor}
\frac{(-1)^k(\zeta(2k+1)}{(2\pi)^{2k}(n+1-2k)!}
.
\ee

\cite{Fichtenholz,MatharArxiv0912} by differentiation of  \cite[3.761.4]{GR} w.r.t.\ the parameter:
\be
\int_0^\infty \frac{\sin x}{x^s}\ln x dx =
\frac{\pi}{2}\frac{1}{[\Gamma(s)\sin\frac{\pi s}{2}]^2}
\left[\Gamma'(s)\sin\frac{s\pi}{2}+\frac{\pi}{2}\Gamma(s)\cos\frac{\pi s}{2}\right]
.
\ee
\begin{multline}
\int_0^\infty \frac{\sin x}{x^s}\ln^2 x dx =
\pi\frac{1}{[\Gamma(s)\sin\frac{\pi s}{2}]^3}
\left\{\Gamma'(s)\sin\frac{s\pi}{2}+\frac{\pi}{2}\Gamma(s)\cos\frac{\pi s}{2}\right\}
\\
-\frac{\pi}{2}
\frac{1}{[\Gamma(s)\sin\frac{\pi s}{2}]^3}
\left\{\Gamma''(s)\sin\frac{s\pi}{2}+\pi\Gamma'(s)\cos\frac{\pi s}{2}-\frac{\pi^2}{4}\Gamma(s)\sin\frac{s\pi}{2}\right\}
.
\end{multline}

\cite{CoffeyJCAM183}
\be
\int_0^1 \ln(1-t)\frac{\ln^2 t}{1-t}dt=-\frac{\pi^4}{180}.
\ee

\cite{CoffeyJCAM183}
\be
\int_0^1 \ln(1-t)\frac{\ln^3 t}{1-t}dt=-\pi^2\zeta(3)+12\zeta(5).
\ee

\cite{CoffeyJCAM183}
\be
\int_0^1 \ln(1-t)\frac{\ln^4 t}{1-t}dt=12\zeta^2(3)-\frac{2\pi^6}{105}.
\ee

\cite{KolbigMathComp39}
Define
\be
s_{n,p}\equiv \frac{(-)^{n+p-1}}{(n-1)!p!}\int_0^1
t^{-1}\log^{n-1}t \log^p(1-t)dt
\ee
then
\be
s_{n,p}=s_{p,n}=\sum_{k=1}^p
\frac{(-)^{k+1}}{k!}
\sum_{m_i}\frac{H_p(m_1,\ldots m_k)}{m_1\cdots m_k}\zeta(m_1)\cdots \zeta(m_k),
\ee
where
\be
H_p(m_1,\ldots,m_k)=\sum_{p_i}\binom{m_1}{p_1}
\cdots\binom{m_k}{p_k},
\ee
the sum over $m_i$ over all sets of integers which 
satisfy $m_i\ge 2$, $\sum_{i=1}^k=n+p$, and the
sum over $p_i$ over all sets of integers which
satisfy $1\le p_i\le m_i-1$, $\sum_{i=1}^kp_i=p$.
Examples with $s_{n,p}=\sum_{k=1}^p (-)^{k+1}\alpha_k(n,p)/k!$
are $\alpha_1(n,p)=(n+p-1)!\zeta(n+p)/(n!p!)$
or $\alpha_2(n,2)=\sum_{\nu=2}^n \zeta(\nu)
\zeta(n-\nu+2)$. The reference provides
an explicit table for $n,p\le 4$.

\cite{KolbigMathComp40}
\begin{gather}
r_{np}\equiv \int_0^{\pi/2}\log^n\cos x\log^p\sin x dx;
\\
r_{10}=-\frac{\pi}{2}\log 2,
\\
r_{11}=\frac{\pi}{2}\left(-\frac{\pi^2}{24}+\log^2 2\right),
\\
r_{20}=\frac{\pi}{2}\left(\frac{\pi^2}{12}+\log^2 2\right),
\\
r_{21}=\frac{\pi}{2}\left(-\log^3 2+\frac{1}{4}\zeta(3)\right),
\\
r_{22}=\frac{\pi}{2}\left(\frac{\pi^4}{160}+\log^4 2-\zeta(3)\log 2\right).
\end{gather}
Due to the symmetry of the integrand upon $x\to \pi/2-x$, $r_{np}=r_{pn}$,
\begin{multline}
2r_{n,n+1}= \int_0^{\pi/2}\log^n\cos x\log^n\sin x[\log \sin x+\log \cos x] dx
=
\int_0^{\pi/2}\log^n\cos x\log^n\sin x \log[ \sin x \cos x] dx
\\
=
\int_0^{\pi/2}\log^n\cos x\log^n\sin x \log[\frac{1}{2} \sin 2x] dx
=
-\ln 2\, r_{n,n} + \int _0^{\pi/2} \log^n\cos x\log ^n\sin x \log (\sin 2x) dx
.
\end{multline}

\cite{Connonarxiv07II}
\be
\int_0^1 \frac{\log u}{1-u}\Li_2\left(\frac{u-1}{u}\right)du
=\frac{17}{4}\zeta(4).
\ee

\cite{Connonarxiv07II}
\be
-\int_0^1 \frac{\Li_{q-1}(1-t)\log t}{1-t}dt = \sum_{n=1}^\infty \frac{H_n^{(1)}}{n^q},
\ee
where $H_n^{(r)}\equiv \sum_{k=1}^n\frac{1}{k^r}$.

\cite{Connonarxiv07II}
\be
\int\frac{\Li_2(1-t)\log t}{1-t}dt=\frac{1}{2}\left[\Li_2(1-t)\right]^2+c.
\ee

\cite{Connonarxiv07II}
\be
\int_0^1\frac{\Li_{2p}(1-t)\log t}{1-t}dt=\frac{1}{2}
\sum_{j=2}^{2p}(-1)^j\zeta(j)\zeta(2p-j+2),
\ee

\cite{BradleyCS2001}
\be
-\frac{1}{4}\int_0^1 \frac{\log x}{(x+1)\sqrt{x}}dx
=
\frac{1}{4}\int_1^\infty \frac{\log x}{(x+1)\sqrt{x}}dx
=G.
\ee

\cite{BradleyCS2001}
\be
\frac{1}{\sqrt{2}}\int_0^{\pi/2} \log(\frac{1+\frac{1}{\surd 2}\sin x}{1-\frac{1}{\surd 2}\sin x})\frac{dx}{1+\cos^2x}
=G.
\ee

\cite{BradleyCS2001}
\be
\frac{1}{2}\int_0^{\pi/4} \log(\frac{1+\sin x}{1-\sin x})\frac{dx}{\cos x\sqrt{\cos 2x}}
=G.
\ee

\cite{BradleyCS2001}
\be
\int_0^{\frac{\sqrt{2}+1}{\sqrt{2}-1}} \frac{(x+1)\log x}{4x\sqrt{6x-x^2-1}}dx
=G.
\ee

\cite{BradleyCS2001}
\be
-\int_0^1 \frac{\log x}{1+x^2}dx = 
\int_1^\infty \frac{\log x}{1+x^2}dx
=G.
\ee

\cite{BradleyCS2001}
\be
-\int_0^1 \log(\frac{1-x}{\sqrt{2}})\frac{dx}{1+x^2}
=G.
\ee

\cite{BradleyCS2001}
\be
-\int_0^1 \log(\frac{1-x^2}{2})\frac{dx}{1+x^2} 
=G.
\ee

\cite{BradleyCS2001}
\be
\int_1^\infty \log(\frac{x+1}{\surd 2})\frac{dx}{1+x^2}
=G.
\ee

\cite{BradleyCS2001}
\be
\int_0^\infty \frac{\log(1+x)}{1+x^2}dx 
=G+\frac{1}{4}\pi\log 2.
\ee

\cite{BradleyCS2001}
\be
-\int_0^1 \frac{\log(1+x^2)}{1+x^2}dx 
=G-\frac{1}{2}\pi\log 2.
\ee

\cite{BradleyCS2001}
\be
-\int_1^{\sqrt{2}} \frac{2\log x}{x\sqrt{x^2-1}}dx 
=G-\frac{1}{2}\pi\log 2.
\ee

\cite{BradleyCS2001}
\be
\int_0^{\pi/2} \log(\cos x+\sin x)dx
=G-\frac{1}{4}\pi\log 2.
\ee

\cite{BradleyCS2001}
\be
-\frac{3}{2}\int_0^{2-\surd 3} \frac{\log x}{1+x^2}dx
=G.
\ee

\cite{BradleyCS2001}
\be
\frac{3}{2}\int_{2+\surd 3}^\infty \frac{\log x}{1+x^2}dx
=G.
\ee

\cite[4.322.8]{GR8}
\be
\int_0^{\pi/2} \ln(\sin x) x^{\mu-1} dx =
-\frac{(\pi/2)^\mu}{\mu} \left [
\frac{1}{\mu}
-\ln\frac{\pi}{2}
+
\mu \sum_{k\ge 1} 
\frac{\zeta(2k)}{k4^k(\mu+2k)}
\right].
\ee

\be
I_{2k}=\int_0^{\pi/4}x^{2k}\frac{\ln x}{\sqrt{1-x^2}}dx
=\int_0^{\pi/4}\sin^{2k} t \ln(\sin t)dt
=\frac{1}{2^{k+1}}\ln 2+S_{2k-2}+(2k-1)I_{2k-2},
\ee
where $I_0=-\frac{\pi}{4}\ln 2-\frac{G}{2}$ with $G\approx 0.9159$
the Catalan constant \cite[4.241.6]{GR}\cite[A006752]{sloane},
where
\be
S_{2k}\equiv \int_0^{\pi/4}\sin^{2k}t\cos^2 t dt
= -\frac{1}{2^{k+2}(k+1)}+\frac{k-1/2}{k+1}S_{2k-2}
\ee
and where $S_0=\frac{1}{4}+\frac{\pi}{8}$.

\be
n\int_0^1 x^n\frac{\log^2 x}{\sqrt{1-x^2}}dx
=(n-2)K(n)+J(n)-M(n)+L(n)
\ee
where $J(n)=\int_0^1 x^n/\sqrt{1-x^2}dx$, $K(n)=\int_0^1 x^n\log(x)/\sqrt{1-x^2}$,
$L(n) = \int_0^1 x^n\log(x)/(1-x^2)^{3/2}$,
$M(n) = \int_0^1 x^n\log^2 (x)/(1-x^2)^{3/2}$, such that
$$
J(n) = \int_0^{\pi/2} \sin^n(t) dt = 2^{n-1}B(\frac{n+1}{2},\frac{n+1}{2});
$$
$K(n)$ known via \cite[4.241.1,4.242]{GR},
$$
L(0)=-\frac{\pi}{2};\quad L(1)=-\ln 2; L(n+2)=-K(n) +L(n);
$$
$$
M(0) 
= \pi \log 2;\quad M(1) = \frac{\pi^2}{12}-\log^2 2;
M(n+2)=-I(n)+M(n).
$$

\cite[(C5.4)]{Nahin}
\be
\int_0^\infty \frac{dx}{(x+a)[\ln^2 x+\pi^2]} =\frac{1}{1-a}+\frac{1}{\ln a},\quad a>0.
\ee

\subsection{Inverse Trigonometric Functions}

\cite{Adamchik,BradleyCS2001}
\be
\frac{2}{\pi}\int_0^{1}\frac{(\tan^{-1}x)^2}{x} = G-\frac{7}{4\pi}\zeta(3).
\ee

\cite{BradleyCS2001}
\be
\frac{3}{2}\int_0^{2-\surd 3}\frac{\tan^{-1}x}{x}dx
=G-\frac{1}{8}\pi\log(2+\surd 3).
\ee

\cite{BradleyCS2001}
\be
-\int_0^1 (\tan^{-1}x)^2dx
=G-\frac{1}{16}\pi^2-\frac{1}{4}\pi\log 2.
\ee

\cite{Adamchik,BradleyCS2001}
\be
\int_0^1\frac{\tan^{-1}x}{x}dx=G.
\ee

\cite{SofoMath7}
\be
\int_0^1 \frac{1}{x}\arctan^3(x)dx = \frac98 \zeta(2)G-\frac32 \beta(4).
\ee

\cite{SofoMath7}
\be
\int_0^1 \frac{1}{x}\arctan^4(x)dx = \frac{1}{16}\pi^3 G-\frac{93}{32}\zeta(5)+\frac32 \pi \beta(4).
\ee

\cite[3.3.1]{Nahin}
\be
\int_0^\infty \frac{\tan^{-1}(ax)-\tan^{-1}(bx)}{x}dx=\frac{\pi}{2}\ln\frac{a}{b}.
\ee

\cite{BradleyCS2001}
\be
2\int_0^1(\frac{1}{4}\pi-\tan^{-1}x)\frac{dx}{1-x^2}=G.
\ee

\cite{BradleyCS2001}
\be
-\int_0^1 \frac{\sin^{-1}x}{\sqrt{1+x^2}}dx
=G-\frac{1}{2}\pi\log(1+\surd 2)
.
\ee

\cite{GlasserMC22}
Let
\be
I_n(a) = \int_0^\infty \frac{(\arctan(z))^{2n}}{z^2+a^2}
\ee
then
\be
I_1(a) = \frac{\pi}{4a}[\frac{\pi^2}{6}-\ln^2(\frac{1+a}{2})
-2\Li_2(\frac{1-a}{2})],\quad a>0.
\ee
where some known special cases of the dilogarithm lead to
\be
I_1(\surd 5) = \frac{\pi}{4\surd 5}[\frac{3\pi^2}{10}-2\ln^2(\frac{1+\surd 5}{2})]
\ee
etc.
\be
\int_0^\infty \frac{(\arctan(z))^4}{z^2+a^2}= \frac{\pi}{4a}[\frac{\pi^4}{40}
+\pi^2\Li_2(\frac{a-1}{a+1})
-6\Li_4(\frac{a-1}{a+1})].
\ee

\cite{GlasserMC22}
\be
\int_0^\infty \frac{(\arctan z)^2}{z^4-6z^2+25}dz
=\frac{\pi}{40}[\frac{7\pi^2}{48}-\frac14 \ln^2 2+\frac{\pi\ln2}{4}-G].
\ee
\be
\int_0^\infty \frac{z^2(\arctan z)^2}{z^4-6z^2+25}dz
=\frac{\pi}{8}[\frac{7\pi^2}{48}-\frac14 \ln^2 2+\frac{\pi\ln2}{4}+G].
\ee

\cite{GlasserMC22}
\be
\int_0^{\pi/2}x^{p+1}\csc^2x dx 
 =(p+1)(\pi/2)^p[p^{-1}-2\sum_{k=1}^\infty (p+2k)^{-1}2^{-2k}\zeta(2k)]
.
\ee

\cite{LiMath10}
Let 
\be
T(m,n)\equiv \int_0^1 \frac{\arctan^m x}{x^n}dx
\ee
then
\be
T(m,0) = \left\{
\begin{array}{ll}
1,&m=0;\\
\frac{\pi}{4}-\frac{\ln 2}{2},& m=1 \\
(\frac{\pi}{4})^m-\frac{m\ln 2}{2}(\frac{\pi}{4})^{m-1}-m(m-1)\Theta(m-2),& m\ge 2 \\
\end{array}
\right.
\ee
\be
T(m,1) = -m\Lambda(m-1), m\ge 1.
\ee
\be
T(m,n) = \frac{m}{n-1}T(m-1,n-1)-\frac{n-3}{n-1}T(m,n-2)-\frac{2}{n-1}(\frac{\pi}{4})^m, m\ge n, n\ge 2
\ee
\be
T(m,n) = \left\{
\begin{array}{ll}
\frac{1}{1-n},& m=0,\\
\frac{m}{1-n}T(m-1,n+1)-\frac{n+1}{n-1}T(m,n+2)-\frac{2}{n-1}(\frac{\pi}{4})^m,& m\ge 1,\\
\end{array}
n<0
\right.
\ee
The auxiliary functions are
\begin{multline}
\Theta(m) = \frac{m!}{2^{m+1}}\eta(m+2)\sin(\frac{m+2}{2}\pi)
+\sum_{j=0}^{\lfloor m/2\rfloor} (-)^j \frac{<m>_{2j}}{2^{2j+1}}
(\frac{\pi}{4})^{m-2j}\beta(2j+2)
\\
-\frac{\ln 2}{m+1}(\frac{\pi}{4})^{m+1}
-\sum_{j=1}^{\lfloor (m+1)/2\rfloor} (-1)^j 
\frac{<m>_{2j-1}}{2^{4j+1}}(\frac{\pi}{4})^{m+1-2j}\eta (2j+1);
\label{eq.LiTheta}
\end{multline}
\begin{multline}
\Lambda(m) = \int_0^{\pi/4} y^m\ln \tan y dy = \frac{m!}{2^m}\lambda(m+2)\sin(\frac{m}{2}\pi)
-\sum_{j=0}^{\lfloor m/2\rfloor} (-)^j \frac{<m>_{2j}}{2^{2j}}
(\frac{\pi}{4})^{m-2j}\beta(2j+2)
\label{eq.LiLambda}
\end{multline}
where $<x>_n=x(x-1)\cdots (x-n+1)$ is the falling factorial,
$\eta(x)=\sum_{n\ge 0}(-)^n/(n+1)^x$ the Dirichlet eta-function
and
$\beta(x)=\sum_{n\ge 0}(-)^n/(2n+1)^x$ the Dirichlet beta-function.

\cite{LiMath10}
Let 
\be
S(m,n)\equiv \int_0^1 \frac{\arcsin^m x}{x^n}dx
\ee
then
\be
S(m,0) = 
\sum_{k=0}^{\lfloor m/2\rfloor -1}(-)^k <m>_{2k} (\frac{\pi}{2})^{m-2k}
+(-)^{\lfloor m/2\rfloor} m!\times \left\{
\begin{array}{ll}
1, m \mathrm{even}\\
\pi/2-1, m \mathrm{odd}
\end{array}
\right.
\ee
\be
S(m,1) = 
-2^m\ln 2(\frac{\pi}{4})^m -m 2^m \Lambda(m-1)-m2^{m+1}\Theta(m-1), m\ge 1;
\ee
\be
S(m,2) = 
-(\frac{\pi}{2})^m -m(m-1)2^{m-1}\Lambda(m-2),m \ge 2;
\ee
\be
S(m,n) = 
\frac{(n-3)^2}{(n-1)(n-2)}
S(m,n-2)+\frac{m(m-1)}{(n-1)(n-2)}
S(m-2,n-2)-\frac{(\pi/2)^m}{(n-1)(n-2)}, m\ge n, n\ge 3;
\ee
\be
S(m,n) = \left\{
\begin{array}{ll}
1/(1-n),& m=0;\\
\frac{\pi}{2-2n}[1-\frac{(n/2)_{-\lfloor n/2\rfloor}}
{(\frac{n-1}{2})_{-\lfloor n/2\rfloor}}]
-\frac{(1)_{-\lfloor n/2\rfloor}(2\pi)}{(3/2)_{-\lfloor n/2\rfloor}(2-2n)} 
\chi(n\equiv_2 0), & m=1;\\
\frac{n(n+1)}{(n-1)^2}S(m,n+2)
-\frac{m(m-1)}{(n-1)^2}S(m-2,n)+\frac{(\pi/2)^m}{(n-1)^2},& m\ge 2
\end{array}
\right.
n<0
\ee
with $\Lambda$ and $\Theta$ defined above, and the Symbol $\chi$
equals 1 if its argument is true, equals 0 if its argument
is false, and $\equiv_2$ indicating evaluation modulo 2.

\cite{GlasserSIAMR26}
\be
\int_0^b\tan^{-1}\frac{a}{\sqrt{1+x^2}}\cdot \frac{dx}{\sqrt{1+x^2}}
=
\frac{1}{2}[g(a,b)+g(b,a)]
\ee
where
\begin{multline}
g(a,b)=\tan^{-1}\frac{b}{a}\ln[4(a^2+1)-\frac{4a\sqrt{a^2+1}}{2a^2+1}]
-2\eta\ln(\sqrt{a^2+1}-a)
\\
-\Cl_2(2\tan^{-1}\frac{b}{a})
+\frac{1}{2}\Cl_2(4\tan^{-1}\frac{b}{a}-2\eta)
+\frac{1}{2}\Cl_2(2\eta)
\end{multline}
and
\be
\eta=\tan^{-1}\frac{ab}{(\sqrt{a^2+1}-a+1)a^2+(\sqrt{a^2+1}-a-1)b^2}.
\ee

\cite{ChoiJMAA231}
\be
\int_0^z \frac{\arcsin at}{t}dt
= \arcsin az \log(2\pi)+\pi \log\frac{G(1-\frac{1}{\pi}\arcsin az)}
{G(1+\frac{1}{\pi}\arcsin az)}
\ee
where $G$ is the reciprocal of the double Gamma function.

\cite[(6.2.6)]{Nahin}
\be
\int_0^1 \frac{\tan^{-1}\sqrt{2+x^2}}{(1+x^2)\sqrt{2+x^2}}dx = \frac{5\pi^2}{96}.
\ee

\cite[(6.3.13)]{Nahin}
\be
\int_0^{\pi/2} \arccos\frac{\cos x}{1+2\cos x}dx=\frac{5\pi^2}{24}.
\ee

\cite{ReynoldsSciA37}
\be
\int_0^\infty \cos(mx)\tan^{-1}[\sinh(\alpha)\sech(bx)]dx
=\frac{\pi}{2m}\sech\frac{\pi m}{2b}\sin\frac{\alpha m}{b}.
\ee

\cite{ReynoldsSciA37}
\be
\int_0^\infty \frac{x^2-1}{(x^2+1)^2} \tan^{-1}[\sech(\pi x)]dx
=-\frac12 \pi [\log(4) -1].
\ee

\cite{ReynoldsSciA37}
\be
\int_0^\infty \frac{1}{(x^2+\pi^2)^2} \tan^{-1}[\sech(x)]dx
=\frac{\log(2\pi/e)}{4\pi^2}.
\ee

\subsection{Multiple Integrals}

\cite{ZhaoJIS13}
\be
\int_0^1\int_0^1 [xy(1-x)(1-y)]^k[(1+xy(1-x)(1-y)]^p dxdy
=\sum_{n=0}^p\frac{\binom{p}{n}}{(2n+2k+1)^2\binom{2n+2k}{n+k}^2}.
\ee

\cite{ZhaoJIS13}
\be
\int_0^1\int_0^1 \frac{[xy(1-x)(1-y)]^k}{1-xy(1-x)(1-y)} dxdy
=\sum_{n\ge 0}\frac{1}{(2n+2k+1)^2\binom{2n+2k}{n+k}^2}.
\ee

\cite{ZhaoJIS13}
\be
\int_0^1\int_0^1 \frac{[xy(1-x)(1-y)]^j}{[1-xy(1-x)(1-y)]^k} dxdy
=\sum_{n\ge 0}\frac{\binom{n+k+1}{n}}{(2n+2j+1)^2\binom{2n+2j}{n+j}^2}.
\ee

\cite{WangJIS13}
\be
jk\int_0^1\int_0^1 \frac{(1-x)^{j-1}(1-y)^{k-1}}{1-tx^ay^b} dxdy
=\sum_{n\ge 0}\frac{t^n}{\binom{an+j}{j}\binom{bn+k}{k}}.
\ee

\cite{BaileyJPA39}
Let
\be
C_n=\frac{4}{n!}\int_0^\infty\cdots \int_0^\infty \frac{1}{[\sum_{j=1}^n(u_j+1/u_j)]^2}
\frac{du_1}{u_1}
\cdots
\frac{du_n}{u_n}
\ee
then
\be
C_2=1;\quad C_3=L_{-3}(s);\quad C_4=7\zeta(3)/12
\ee
see (\ref{eq.C3}).

\cite{BaileyJPA39}
Let
\be
D_n=\frac{4}{n!}\int_0^\infty\cdots \int_0^\infty \frac{\prod_{i<j}
(\frac{u_i-u_j}{u_i+u_j})^2
}{[\sum_{j=1}^n(u_j+1/u_j)]^2}
\frac{du_1}{u_1}
\cdots
\frac{du_n}{u_n}
=
\frac{1}{n!}\int d^nx \frac{\prod_{i<j} \tanh^2[(x_i-x_j)/2]}{(\cosh x_1+\cdots +\cosh x_n)^2}
\ee
then
\be
D_1=2;\quad D_2=1/3;\quad D_3=8+4\pi^2/3-27L_{-3}(2);\quad
D_4=4\pi^2/9-1/6-7\zeta(3)/2,
\ee
see (\ref{eq.C3}).

\cite{Koholdarxiv1201}
\be
\int_{-\infty}^\infty ds_1 e^{-is_1^2}
\int_{-\infty}^{s_1} ds_2 e^{-is_2^2}
\int_{-\infty}^{s_2} ds_3 e^{-is_3^2}
\cdots 
\int_{-\infty}^{s_{2n-1}} ds_{2n} e^{-is_{2n}^2}
=\frac{1}{n!}(\pi/2)^n.
\ee

\cite{Koholdarxiv1201}
\be
\int_{-\infty}^\infty ds_1 
\int_{-\infty}^{s_1} ds_2 
\int_{-\infty}^{s_2} ds_3
\cdots 
\int_{-\infty}^{s_{2n-1}} ds_{2n}
\cos(s_1^2-s_2^2)
\cdots
\cos(s_{2n-1}^2-s_{2n}^2)
=\frac{2}{n!}(\pi/4)^n.
\ee

\cite{BraggAMM106}
\begin{equation}
\frac{1}{(2\pi)^2} \int_0^{2\pi} \int_0^{2\pi}
[(1+ze^{i(\phi+\theta)})(1+ze^{i(\phi-\theta)})(1+ze^{-i\phi})]^n d\phi d\theta
=\sum_{j=0}^{[n/2]}\binom{n}{j}^2\binom{n}{2j}z^{4j}.
\end{equation}

\cite{BraggAMM106}
\begin{equation}
\frac{1}{(2\pi)^2} \int_0^{2\pi} \int_0^{2\pi}
\cos(n\phi/2)\cos^n(\phi/2)[\cos\phi+\cos\theta]^n d\phi d\theta
=\frac{4\pi^2}{2^{2n}}\sum_{j=0}^{[n/2]}\binom{n}{j}^2\binom{n}{2j}.
\end{equation}

\cite{BraggAMM106}
\begin{equation}
\int_0^{2\pi} \int_0^{2\pi}\int_0^{2\pi}
[(\cos\theta_1+\cos\theta_2)(\cos\theta_1+\cos\theta_3)]^nd\theta_1
d\theta_2 d\theta_3
=\frac{(2\pi)^3}{2^{2n}}\sum_{k=0}^{n}\binom{n}{k}^4.
\end{equation}

\cite{BraggAMM106}
\begin{equation}
\int_0^{\pi} \int_0^{\pi}
\cos[(m-n)\theta+(m+n-2p)\phi]\cos^m(\theta+\phi)\cos^n(\theta-\phi)
\cos^p(2\phi)d\theta d\phi
=\frac{\pi^2}{2^{m+n+p}}\sum_{j=0}^{\min(m,n,p)}\binom{m}{j}\binom{n}{j}\binom{p}{j}.
\end{equation}

\cite{BraggAMM106}
\begin{equation}
\int_0^{2\pi} \int_0^{2\pi}
[\cos\theta(\cos\theta+\cos\phi)]^md\theta d\phi
=\frac{4\pi^2}{2^{2m}}\sum_{j=0}^{m}\binom{m}{j}^3.
\end{equation}

\cite{BradleyCS2001}
\be
\frac{2}{\pi}\int_0^{\pi/2} \int_0^{\pi/2} \tan^{-1}(\sin x\sin y)\frac{dx dy}{\sin x}
=G.
\ee

\cite{BradleyCS2001}
\be
\frac{1}{2}\int_0^1 \int_0^{\pi/2} \frac{d\theta dx}{\sqrt{1-x^2\sin^2\theta}}
=G.
\ee

\cite{BradleyCS2001}
\be
\int_0^1 \int_0^{\pi/2} \sqrt{1-x^2\sin^2\theta}d\theta dx
=G+\frac{1}{2}.
\ee

\cite{BradleyCS2001}
\be
8\int_0^1 \int_0^1 \frac{\tan^{-1}(xy)dx dy}{1+x^2y^2}
=2\pi G-\frac{7}{2}\zeta(3).
\ee

\cite{BradleyCS2001}
\be
4\int_0^1 \int_0^1 \frac{\tan^{-1}x}{1+x^2y^2}dxdy
=2\pi G-\frac{7}{2}\zeta(3).
\ee

\cite{BradleyCS2001}
\be
-\int_0^1 \int_0^1 \frac{\log(1-x^2y^2)}{xy\sqrt{(1-x^2)(1-y^2)}}dxdy
=2\pi G-\frac{7}{2}\zeta(3).
\ee

\cite[1945637]{stackexch}
\be
\int_0^\infty dx \int_0^\infty dy \frac{\log|(x+y)(1-xy)|}{(1+x^2)(1+y^2)}
= \frac{\pi^2}{2}\log 2
.
\ee

\cite{BerndtBLMS15}
\be
\int_0^{\pi/2}
\int_0^{\pi/2}
\frac{\tan(\phi/2)d\theta d\phi}{\sqrt{1-x\cos^2\theta\cos^2\phi}}=
\frac{\pi}{4}\int_0^{\pi/2}
\frac{d\phi}{\sqrt{1-(1-x)\sin^2\phi}}
+
\frac{1}{4}\log x\int_0^{\pi/2}\frac{d\phi}{\sqrt{1-x\sin^2\phi}}.
\ee

\cite{Adamchik,BradleyCS2001}
\be
\frac{1}{4}\int_0^1\int_0^1 \frac{1}{\sqrt{(1-x)(1-y)}}\,\frac{dx dy}{x+y}=G.
\ee

\cite{GustafsonBAMS22}
\be
\int_0^1\cdots \int_0^1 \prod_{1\le i<j\le n}|t_i-t_j|^{2z}
\prod_{j=1}^n t_j^{x-1}(1-t_j)^{y-1}dt_j
=
\prod_{j=1}^n\frac{\Gamma(x+(j-1)z)\Gamma(y+(j-1)z)\Gamma(jz+1)}
{\Gamma(x+y+(n+j-2)z)\Gamma(z+1)},
\ee
where $n$ is a positive integer, $x$, $y$, $z$ are in $\mathbb{C}$,
and $\Re x$, $\Re y>0$, $\Re z> -\max\{1/n,\Re x/(n-1), \Re y/(n-1)\}$.

\cite{Connonarxiv07II}
\be
\int_0^1\int_0^1 \frac{\Li_{q-2}[(1-t)(1-u)]\log t\log u\,dudt}{(1-t)(1-u)}
= \sum_{n=1}^\infty\frac{\left[H_n^{(1)}\right]^2}{n^q}
\ee
where $H_n^{(r)}\equiv \sum_{k=1}^n\frac{1}{k^r}$.

\cite{Connonarxiv07II}
\be
\int_0^1\int_0^1 \frac{\log t\log u\,dudt}{1-(1-t)(1-u)}
= \sum_{n=1}^\infty\frac{\left[H_n^{(1)}\right]^2}{n^2}
\ee
where $H_n^{(r)}\equiv \sum_{k=1}^n\frac{1}{k^r}$.

\cite{Connonarxiv07II}
\be
-\int_0^1\int_0^1 \frac{\log[1-(1-t)(1-u)]\log t\log u\,dudt}{(1-t)(1-u)}
= \sum_{n=1}^\infty\frac{\left[H_n^{(1)}\right]^2}{n^3}
\ee
where $H_n^{(r)}\equiv \sum_{k=1}^n\frac{1}{k^r}$.

\cite{Connonarxiv07II}
\be
n^2\int_0^1(1-t)^{n-2}\log t\,dt \int_0^1(1-u)^{n-1}\log u\,du
= \left[H_n^{(1)}\right]^2,
\ee
where $H_n^{(r)}\equiv \sum_{k=1}^n\frac{1}{k^r}$.

\cite{Amdeberhanarxiv2379}
\be
\int_0^\infty \int_0^\infty \frac{\cos^{2n+1}(x+y)}{x^py^q}dxds=
-\Gamma(1-p)\Gamma(1-q)\cos\frac{\pi(p+q)}{2}\sum_{k=0}^n\binom{2n+1}{n-k}\frac{(2k+1)^{p+1-2}}{2^{2n}}
\ee

\cite{Amdeberhanarxiv2379}
\be
\int_0^\infty \int_0^\infty \frac{\cos(x+y)}{x^py^q}dxdy=
-\Gamma(1-p)\Gamma(1-q)\cos\frac{\pi(p+q)}{2}.
\ee

\cite{Amdeberhanarxiv2379}
\be
\int_0^\infty \int_0^\infty \frac{\log x \log y}{\sqrt{xy}}\cos(x+y)dxdy=
(\gamma+2\log 2)\pi^2.
\ee

\cite{Amdeberhanarxiv2379}
\be
\int_{R_+^n}(\cos ||x||^2)\cdot \prod_{j=1}^n\log x_j dV=
\frac{(-)^{\Delta_n}\pi^{n/2}}{2^{2n}}
\times
\begin{cases}
\Re \psi_n& n\, \mathrm{even}\\
\Im\psi_n& n\, \mathrm{odd}
\end{cases}
\ee
where $\Delta_n\equiv n(n+1)/2$ and $\psi_n\equiv (\gamma+2\log 2+\pi i/2)^n e^{\pi i n/4}$.

\cite{GuilleraRJ16}
\be
\int_0^1
\int_0^1
\frac{x^{u-1}y^{v-1}}{1-xyz}(-\ln xy)^sdxdy = \Gamma(s+1)
\frac{\Phi(z,s+1,v)-\Phi(z,s+1,u)}{u-v}
,
\ee
\cite{GuilleraRJ16}
\be
\int_0^1
\int_0^1
\frac{(xy)^{u-1}}{1-xyz}(-\ln xy)^sdxdy = \Gamma(s+1)
\Phi(z,s+2,u)
,
\ee
where
\be
\Phi(z,s,u)=\sum_{k=0}^\infty \frac{z^k}{(u+k)^s}
=
\frac{1}{\Gamma(s)}\int_0^\infty \frac{e^{-(u-1)t}}{e^t-z}t^{s-1}dt
\ee
is the Lerch transcendent.

\cite[(5.3.1)]{Nahin}
\be
\int_0^1\int_0^1 \frac{(xy)^2}{1-xy}dxdy=\sum_{n\ge 1} \frac{1}{(n+a)^2}.
\ee

\cite[(5.3.2)]{Nahin}
\be
\int_0^1\int_0^1 \frac{(xy)^2 [\ln(xy)]^{s-2}}{1-xy}dxdy=(-)^s (s-1)!\sum_{n\ge 1} \frac{1}{(n+a)^s}.
\ee

\cite{GuilleraRJ16}
\be
\int_0^1
\int_0^1
\frac{1}{1+xyi}dxdy = G-\frac{\pi^2i}{48}
,
\ee
\cite{GuilleraRJ16}
\be
\int_0^1
\int_0^1
\frac{-x\ln xy}{1+x^2y^2}dxdy = G-\frac{\pi^2}{48}
,
\ee
\cite{GuilleraRJ16}
\be
\int_0^1
\int_0^1
\frac{-x\ln xy}{1-x^2y^2}dxdy = \frac{\pi^2}{12}
,
\ee
\cite{GuilleraRJ16}
\be
\int_0^1
\int_0^1
\frac{(-\ln xy)^n}{1-xyz}dxdy = \frac{(n+1)!\Li_{n+2}(z)}{z}
,
\ee
\cite{GuilleraRJ16}
\be
\int_0^1
\int_0^1
\frac{-1}{(1-xyz)\ln xy}dxdy = -\frac{\ln(1-z)}{z}
,
\ee
\cite{GuilleraRJ16}
\be
\int_0^1
\int_0^1
\frac{-1}{(2-xyz)\ln xy}dxdy = \ln 2
,
\ee

\cite{GuilleraRJ16}
\be
\int_0^1
\int_0^1
\frac{1}{2-xy}dxdy = \frac{\pi^2}{12}-\frac{\ln^2 2}{2}
.
\ee
\cite{GuilleraRJ16}
\be
\int_0^1
\int_0^1
\frac{-\ln xy}{2-xy}dxdy = \frac{7\zeta(3)}{4}-\frac{\pi^2\ln 2}{6}+\frac{\ln^3 2}{3}
.
\ee
\cite{GuilleraRJ16}
\be
\int_0^1
\int_0^1
\frac{-1}{(\varphi-xy)\ln xy}dxdy = \ln \varphi
,
\ee
\cite{GuilleraRJ16}
\be
\int_0^1
\int_0^1
\frac{1}{\varphi-xy}dxdy = \frac{\pi^2}{10}-\ln^2\varphi
,
\ee
\cite{GuilleraRJ16}
\be
\int_0^1
\int_0^1
\frac{-1}{(\varphi^2-xy)\ln xy}dxdy = \ln\varphi
,
\ee
etc
where
\be
\varphi\equiv (1+\surd 5)/2.
\ee
\cite{GuilleraRJ16}
\be
\int_0^1
\int_0^1
\frac{1-2xy}{(8+xy)(9-xy)}dxdy = \frac{1}{2}\ln^2\frac{9}{8}
.
\ee
\cite{GuilleraRJ16}
\be
\int_0^1
\int_0^1
\frac{52-7xy}{(2+xy)(9-xy)}dxdy = \frac{\pi^3}{3}+3\ln^2 2+2\ln^2 3 -6\ln 2\ln 3
.
\ee
etc
\cite{GuilleraRJ16}
\be
\int_0^1
\int_0^1
\frac{(-\ln xy)^s}{1-xy}dxdy = \Gamma(s+2)\zeta(s+2), \Re s>-1
.
\ee
\cite{GuilleraRJ16}
\be
\int_0^1
\int_0^1
\frac{(-\ln xy)^s}{1+xy}dxdy = \Gamma(s+2)\zeta^*(s+2), \Re s>-2
,
\quad
\zeta^*(s)\equiv (1-2^{1-s})\zeta(s)
.
\ee
\cite{GuilleraRJ16}
\be
\int_0^1
\int_0^1
\frac{(-\ln xy)^s}{1+x^2y^2}dxdy = \Gamma(s+2)\beta(s+2), \Re s>-2
.
\ee
\cite{GuilleraRJ16}
\be
\int_0^1
\int_0^1
\frac{1}{1-xy}dxdy = \zeta(2)
.
\ee
\cite{GuilleraRJ16}
\be
\int_0^1
\int_0^1
\frac{-\ln xy}{1-xy}dxdy = 2\zeta(3)
.
\ee
\cite{GuilleraRJ16}
\be
\int_0^1
\int_0^1
\frac{-1}{(1+x^2y^2)\ln xy}dxdy = \pi/4
.
\ee
\cite{GuilleraRJ16,BradleyCS2001}
\be
\int_0^1
\int_0^1
\frac{1}{1+x^2y^2}dxdy = G
.
\ee

\cite{GuilleraRJ16}
\be
\int_0^1
\int_0^1
\frac{\ln xy}{1+x^2y^2}dxdy = \frac{\pi^3}{16}
.
\ee
\cite{GuilleraRJ16}
\be
\int_0^1
\int_0^1
\frac{(-\ln xy)^s}{1+x^2y^2z^2}dxdy = \Gamma(s+1)\frac{\chi_{s+2}(z)}{z}
,
\quad \Re s>-2, \Re z\neq 0.
\ee
\cite{GuilleraRJ16}
\be
\int_0^1
\int_0^1
\frac{1}{1-x^2y^2\tan^2(\pi/8)}dxdy
= \frac{\pi^2}{16\tan\frac{\pi}{8}}-
\frac{\ln^2 \tan\frac{\pi}{8}}{4\tan\frac{\pi}{8}}
.
\ee
\cite{GuilleraRJ16}
\be
\int_0^1
\int_0^1
\frac{x^{u-1}y^{v-1}}{-\ln xy}dxdy
=
\frac{1}{u-v}\ln\frac{u}{v}
.
\ee
\cite{GuilleraRJ16}
\be
\int_0^1
\int_0^1
\frac{(xy)^{u-1}}{-\ln xy}dxdy
=
\frac{1}{u}
.
\ee
\cite{GuilleraRJ16}
\be
\int_0^1
\int_0^1
\frac{-x^{u-1}y^{v-1}}{(1+xy)\ln xy}dxdy
=
\frac{1}{u-v}
\ln\frac{\Gamma(u/2)\Gamma(\frac{u+1}{2})}
{\Gamma(v/2)\Gamma(\frac{v+1}{2})}
,
\quad u>0, v>0.
\ee
\cite{GuilleraRJ16}
\be
\int_0^1
\int_0^1
\frac{-(xy)^{u-1}}{(1+xy)\ln xy}dxdy
=
\frac{1}{2}
\left[\psi(\frac{u+1}{2}-\psi(\frac{u}{2})\right]
.
\ee
\cite{GuilleraRJ16}\cite[A053510]{sloane}
\be
\int_0^1
\int_0^1
\frac{1+x}{-(1+xy)\ln xy}dxdy
=
\ln\pi
.
\ee
\cite{GuilleraRJ16}\cite[A094640]{sloane}
\be
\int_0^1
\int_0^1
\frac{1-x}{-(1+xy)\ln xy}dxdy
=
\ln\frac{4}{\pi}
.
\ee
\cite{GuilleraRJ16}
\be
\int_0^1
\int_0^1
\frac{-x}{(1+x^2y^2)\ln xy}dxdy
=
\ln\frac{\sqrt{2\pi}}{\Gamma^2(3/4)}
.
\ee
\cite{GuilleraRJ16}
\be
\int_0^1
\int_0^1
\frac{x^{u-1}y^{v-1}}{1-xy}dxdy
=
\frac{\psi(u)-\psi(v)}{u-v}
;
\quad
\int_0^1
\int_0^1
\frac{(xy)^{u-1}}{1-xy}dxdy
=
\psi'(u).
\ee
\cite{GuilleraRJ16}\cite[A073010]{sloane}
\be
\int_0^1
\int_0^1
\frac{y}{(1-x^3y^3)}dxdy
=
\frac{\pi}{3\surd 3}
.
\ee
\cite{GuilleraRJ16}
\be
\int_0^1
\int_0^1
x^{u-1}y^{v-1}(-\ln xy)^sdxdy
=
\Gamma(s+1)\frac{v^{-s-1}-u^{-s-1}}{u-v}
\ee
and others of similar shape.

\cite{TaylorPEMS32}
\be
\int_0^\pi \int_0^\pi \frac{\cos n\alpha-\cos n\theta}{\cos \alpha-\cos \theta}d\alpha d\theta =0\,\mathrm{or}\,\pi^2
\ee
according as $n$ is even or odd.
\cite{TaylorPEMS32}
\be
\int_0^\pi \int_0^\pi \frac{(\cos n\alpha-\cos n\theta)\cos r\alpha \cos s\theta}{\cos \alpha-\cos \theta}d\alpha d\theta =0\,\mathrm{or}\,\pi^2
\ee
according as $n-r-s$ is even or odd.

\cite{OkuiNBS79B} representations as complete elliptic Intgrals:
\be
\int_0^\infty \int_0^\infty e^{-px^2} e^{-qy^2}I_0(2axy) y^l dxdy
=\ldots
\ee
for even $l=0,\ldots,6$.
\cite{OkuiNBS79B}
\be
\int_0^\infty \int_0^\infty e^{-px^2} e^{-qy^2}I_0(2axy)x^2 y^l dxdy
=\ldots
\ee
for even $l=2,\ldots,8$.
\cite{OkuiNBS79B}
\be
\int_0^\infty \int_0^\infty e^{-px^2} e^{-qy^2}I_1(2axy)x^k y^l dxdy
=\ldots
\ee
for odd $k \ge 1$, odd $l \ge 1$ but even $k+l<12$. Similar integrals
also where $I_2$ and $I_3$ represent the Besel function.
\cite{OkuiNBS79B}
\be
\int_0^\infty \int_0^\infty e^{-px^2} e^{-qy^2}I_k^2(axy) x^s y^l dxdy
=\ldots
\ee
for even small integer $k$ and symmetry-adapted exponents $s$ and $l$.
\cite{OkuiNBS79B}
\be
\int_0^\infty \int_0^\infty e^{-px^2} e^{-qy^2}I_0(ax^2)I_0(by^2)I_0(2cxy) xy dxdy
=\ldots
\ee
\cite{OkuiNBS79B}
\be
\int_0^\infty \int_0^\infty e^{-px^2} e^{-qy^2}I_0(ax^2)I_0(by^2)I_0(2cxy) xy^3 dxdy
=\ldots
\ee
and more double integrals of triple products of $I_k$ times $xy$.
\cite{OkuiNBS79B}
\be
\int_0^\infty \int_0^\infty e^{-px^2} e^{-qy^2}I_k(2axy)\cosh(2bxy) x^{-l} y^l dxdy
=\ldots
\ee
for $0\le k\le 4$ and $1\le l\le 4$. Additional formulas where $\sinh(2bxy)$
takes the place of $\cosh(2bxy)$.

\cite[(6.4.7)]{Nahin}
\be
\int_0^\infty dx \int_x^\infty dt \int_x^\infty du \cos(t^2-u^2) = \frac12 \sqrt{\pi/2}.
\ee

For $r\equiv \sqrt{x^2+y^2+z^2}$ via multinomial expansion
\begin{multline}
R(2n)\equiv \int_0^1 dz\int_0^1 dy \int_0^1 dx r^{2n} 
= \sum_{i=0}^n\sum_{j=0}^{n-i}\binom{n}{i,j,n-i-j} \frac{1}{(2i-1)(2j-1)(2n-2i-2j-1)}\\
= \sum_{l=0}^n\binom{n}{l}\frac{1}{2n-2l-1}\sum_{i=0}^l \binom{l}{i}\frac{1}{(2i-1)(2l-2i-1)} \\
= -\sum_{l=0}^n\binom{n}{l}\frac{1}{(2n-2l-1)(2l-1)}{}_3F_2\left(\begin{array}{c} -1/2,-l,1/2-l\\1/2,3/2-l\end{array}\mid -1\right)\\
=1,1,\frac{19}{15},\frac{583}{315},\frac{1573}{525},\frac{2599}{495},\ldots (n=0,1,2,\ldots)
\end{multline}

\cite[A130590]{sloane}
\begin{multline}
R(1)=\int_0^1 dz\int_0^1dy \int_0^1 dx r 
=
2\int_0^1 dz \int_0^1 dx \int_0^x dy \sqrt{x^2+y^2+z^2}\\
=
2\int_0^1 dz \int_0^{\pi/4} d\phi \int_0^{1/\cos\phi} \rho d\rho \sqrt{\rho^2+z^2}
=
\frac23 \int_0^1 dz \int_0^{\pi/4} d\phi [(z^2+\frac{1}{\cos^2\phi})^{3/2}-z^3]\\
=
-\frac{\pi}{24}
+\frac23 \int_0^1 dz \int_{z^2+1}^{z^2+2} \frac{\cos^3\phi}{2\sin\phi}dt t^{3/2}
=
-\frac{\pi}{24}
+\frac13 \int_0^1 dz \int_{z^2+1}^{z^2+2} \frac{1}{\sqrt{1-\frac{1}{t-z^2}}(t-z^2)^{3/2}}dt t^{3/2}\\
=
-\frac{\pi}{24}
+\frac13 \int_0^1 dz \int_{z^2+1}^{z^2+2} \frac{1}{\sqrt{t-z^2-1}(t-z^2)}dt t^{3/2}\\
=
-\frac{\pi}{24}
+\frac13 \int_0^1 dz \int_{\sqrt{z^2+1}}^{\sqrt{z^2+2}} \frac{1}{\sqrt{u^2-z^2-1}(u^2-z^2)}2u du u^3
\\
=
-\frac{\pi}{24}
+\frac23 \int_0^1 dz \int_{\sqrt{z^2+1}}^{\sqrt{z^2+2}} \frac{1}{\sqrt{u^2-z^2-1}}(u^2+z^2+z^4\frac{1}{u^2-z^2})du
\end{multline}
where
\begin{multline}
\int \frac{u^2+z^2}{\sqrt{u^2-z^2-1}} du = \frac12 u\sqrt{u^2-z^2-1}+\frac12(3z^2+1)\ln(u+\sqrt{u^2-z^2-1}),
\end{multline}
\begin{multline}
\int \frac{1}{\sqrt{u^2-z^2-1}(u^2-z^2)} du 
=\frac{1}{2z}\arctan\frac{2zu\sqrt{u^2-z^2-1}}{u^2-z^2(u^2-z^2-1)}.
\end{multline}
\begin{multline}
\int_{\sqrt{z^2+1}}^{\sqrt{z^2+2}} \frac{1}{\sqrt{u^2-z^2-1}}(u^2+z^2+z^4\frac{1}{u^2-z^2})du
\\
=
\frac12 \sqrt{2+z^2}
+\frac12 (3z^2+1)\ln\frac{1+\sqrt{2+z^2}}{\sqrt{1+z^2}}
+\frac12 z^3\arctan(z\sqrt{2+z^2}).
\end{multline}
Note that expressing $r$ as the divergence of a vector field and using the Gauss' law ends
up with the same type of double integrals.
\be
\int dz \sqrt{2+z^2}=\frac{z}{2}\sqrt{2+z^2}+\Arsh\frac{z}{\surd 2}.
\ee
\begin{equation}
\int z^3 \arctan(z\sqrt{2+z^2}) dz = 
-\frac14 z\sqrt{2+z^2}+\Arsh\frac{z}{\surd 2}
+\frac14 z^4\arctan(z\sqrt{2+z^2})-\frac12 \arctan\frac{z}{\sqrt{2+z^2}}.
\end{equation}
\be
\int (3z^2+1)\ln(1+z^2) dz = z(1+z^2)\ln(1+z^2)-\frac23 z^3.
\ee
\be
-\int (3z^2+1)\ln\sqrt{1+z^2} dz = -z(1+z^2)\ln\sqrt{1+z^2}+\frac13 z^3.
\ee
\begin{multline}
\int \ln(1+\sqrt{2+z^2}) 
= z[\ln(1+\sqrt{2+z^2})-1]+\Arsh\frac{z}{\surd 2}-\arctan\frac{z}{\sqrt{2+z^2}}+\arctan z
.
\end{multline}
\begin{multline}
3\int z^2\ln(1+\sqrt{2+z^2}) 
=
z^3\ln(1+\sqrt{2+z^2})-\frac13 z^3
-2\Arsh\frac{z}{\surd 2}+\frac{z}{2}\sqrt{2+z^2}
+\arctan\frac{z}{\sqrt{2+z^2}}
+z-\arctan z.
\end{multline}
The sum of the previous two integrals is
\be
\int(3z^2+1) \ln(1+\sqrt{2+z^2}) dz
=
z(1+z^2)\ln(1+\sqrt{2+z^2})
-\Arsh\frac{z}{\surd 2}+\frac{z}{2}\sqrt{2+z^2}-\frac13 z^3.
\ee
\begin{multline}
\int(3z^2+1) \ln\frac{1+\sqrt{2+z^2}}{\sqrt{1+z^2}} dz
=
z(1+z^2)\ln\frac{1+\sqrt{2+z^2}}{\sqrt{1+z^2}}
-\Arsh\frac{z}{\surd 2}+\frac{z}{2}\sqrt{2+z^2}.
\end{multline}
\begin{multline}
\frac23 \int dz \int_{\sqrt{z^2+1}}^{\sqrt{z^2+2}} \frac{1}{\sqrt{u^2-z^2-1}}(u^2+z^2+z^4\frac{1}{u^2-z^2})du
\\
=
\frac13[
\frac{3z}{4}\sqrt{2+z^2}+\Arsh\frac{z}{\surd 2}
+z(1+z^2)\ln\frac{1+\sqrt{2+z^2}}{\sqrt{1+z^2}}
\\
+\frac14 z^4\arctan(z\sqrt{2+z^2})-\frac12 \arctan\frac{z}{\sqrt{2+z^2}}
]
\end{multline}
\begin{multline}
\frac23 \int_0^1 dz \int_{\sqrt{z^2+1}}^{\sqrt{z^2+2}} \frac{1}{\sqrt{u^2-z^2-1}}(u^2+z^2+z^4\frac{1}{u^2-z^2})du
=
\frac{\surd 3}{4}+\frac13 \Arsh\frac{1}{\surd 2}+\frac23 \ln\frac{1+\surd 3}{\surd 2}.
=
\frac{\surd 3}{4}+ \ln\frac{1+\surd 3}{\surd 2}.
\end{multline}
\begin{equation}
R(1)=
\frac{\surd 3}{4}+ \ln\frac{1+\surd 3}{\surd 2}-\frac{\pi}{24}
\approx 0.960591956455052959\ldots
.
\end{equation}

\cite{WanAAM48}
\begin{multline}
\int_{[0,1]^3} \frac{x^{h_2-1}y^{h_3-1}z^{h_4-1}
(1-x)^{h_0-h_2-h_3}(1-y)^{h_0-h_3-h_4}(1-z)^{h_0-h_4-h_5}
}{(1-x(1-y(1-z)))^{h_1}}dxdydz\\
=
\frac{\Gamma(h_0+1)\prod_{j=2}^4 \Gamma(h_j)
\prod_{j=1}^4\Gamma(h_0+1-h_j-h_{j+1}}{\prod_{j=1}^5
\Gamma(h_0+1-h_j)}\\ \times
{}_7F_6\left(
\begin{array}{c}
h_0,1+\frac{h_0}{2},h_1,h_2,h_3,h_4,h_5 \\
\frac{h_0}{2},1+h_0-h_1,1+h_0-h_2,1+h_0-h_3,1+h_0-h_4,
h+h_0-h_5
\end{array}\mid 1
\right).
\end{multline}

\cite{WanAAM48}
\be
\frac18\int_{[0,1]^3} \frac{\surd y dx dy dz}{\sqrt{x(1-x)(1-y)z(1-z)(1-x(1-y(1-z)))}}=\int_0^1E'(x)K'(x)dx.
\ee

\cite{WanAAM48}
\be
\frac18\int_{[0,1]^3} \frac{dx dy dz}{\sqrt{x(1-x)y(1-y)z(1-z)(1-x(1-y(1-z)))}}=\int_0^1(x)K'(x)^2dx.
\ee

\cite[(2.5.1)]{Nahin}
Let
\be
R(m,n)=\frac{1}{(2\pi)^2}\int_{-\pi}^{\pi}\int_{\pi}^{\pi}\frac{-1\cos(mx+ny)}{2-\cos x -\cos y}dxdy.
\ee
then 
\be
R(m+1,n)+R(m-1,n)+R(m,n+1)+R(m,n-1)-4R(m,n) = \left\{
\begin{array}{rr}
2, & \mathrm{ if} m=n=0\\
0, & \mathrm{ otherwise}
\end{array}
\right.
\ee
with $R(1,2)=4/\pi -1/2$. Others by recursion.

\cite[(6.5.1)]{Nahin}
\be
\frac{1}{\pi^3}
\int_0^\pi 
\int_0^\pi 
\int_0^\pi 
\frac{du dv dw}{1-\cos u\cos v\cos w}=\frac{\Gamma^4(1/4)}{4\pi^3}
\ee

\cite[p. 246]{Nahin}\cite[A091671]{sloane}
\be
\frac{1}{\pi^3}
\int_0^\pi 
\int_0^\pi 
\int_0^\pi 
\frac{du dv dw}{3-\cos u\cos w-\cos w \cos u -\cos u \cos v}=\frac{3\Gamma^6(1/3)}{2^{14/3}\pi^4}.
\ee

\cite[p. 246]{Nahin}\cite[A091672]{sloane}
\be
\frac{1}{\pi^3}
\int_0^\pi 
\int_0^\pi 
\int_0^\pi 
\frac{du dv dw}{3-\cos u-\cos v-\cos w}=\frac{\Gamma(1/24)\Gamma(5/24)\Gamma(7/24)\Gamma(11/24)}{16\sqrt{6} pi^3}.
\ee

\cite{BorweinNuy2011}
\be
W_n(2k)=\sum_{a_1+a_2+\cdots+a_n=k}\left(
\begin{array}{c}k \\ a_1,\ldots,a_n
\end{array}\right)^2,
\ee
where the sum is over all compositions (unordered partitions) with $n$ terms, and
\be
W_n(s)\equiv \int_{[0,1]^n} \left|\sum_{k=1}^n e^{2\pi ix_k}\right|^s d^nx
.
\ee
\cite{BorweinNuy2011}
\be
W_3(k)=\Re \,_3F_2\left(\begin{array}{c}1/2,-k/2,-k/2\\ 1,1\end{array}\mid 4\right).
\ee
\cite{BorweinNuy2011}
\be
W_n(s)=n^s\sum_{m\ge 0}(-1)^m\binom{s/2}{m}\sum_{k=0}^m\frac{(-)^k}{n^{2k}}
\binom{m}{k}
\sum_{a_1+a_2+\cdots+a_n=k}
\left(
\begin{array}{c}k \\ a_1,\ldots,a_n
\end{array}\right)^2.
\ee

\section{Indefinite Integrals of Special Functions}
\subsection{Elliptic Integrals and Functions}
\subsection{The Exponential Integral}
\subsection{The Sine Integral and Cosine Integral}
\cite{WatrasieOA14}
\be
\int Si(ax+b)dx = \frac{ax+b}{a}Si(ax+b)+\frac{\cos(ax+b)}{a}.
\ee
\cite{WatrasieOA14}
\be
\int Ci(ax+b)dx = \frac{ax+b}{a}Ci(ax+b)-\frac{\sin(ax+b)}{a}.
\ee
\cite{WatrasieOA14}
\be
\int Si(ax)\sin(bx) dx = \frac{1}{2b}[Si(a+b)x +Si(a-b)x -2Si(ax)\cos(bx)].
\ee
\cite{WatrasieOA14}
\be
\int Si(ax)\cos(bx) dx = \frac{1}{2b}[2Si(ax)\sin(bx) +Ci(a+b)x -Ci(b-a)x].
\ee
\cite{WatrasieOA14}
\be
\int Ci(ax)\sin(bx) dx = \frac{1}{2b}[Ci(a+b)x +Ci(b-a)x -2Ci(ax)\cos(bx)].
\ee
\cite{WatrasieOA14}
\be
\int Ci(ax)\cos(bx) dx = \frac{1}{2b}[2Ci(ax)\sin(bx) -Si(a+b)x +Si(a-b)x].
\ee
\cite{WatrasieOA14}
\be
\int x^n Si(ax)dx = \frac{x^n}{n+1} I_0 -\frac{n}{a(n+1)}C_{n-1}
\ee
where $I_0=\int Si(ax)dx$, $C_n=\int x^n \cos(ax)dx$.

\cite{WatrasieOA14}
\be
\int x^n Ci(ax)dx = \frac{x^n}{n+1} \int Ci(ax)dx +\frac{n}{a(n+1)}S_{n-1}
\ee
where $I_0=\int Ci(ax)dx$, $S_n=\int x^n \sin(ax)dx$.

\cite{WatrasieOA14}
\be
I_n=\int x^nSi(ax+b)dx =\frac{1}{n+1}[x^nI_0-\frac{nb}{a}I_{n-1}-\frac{n}{a}C_{n-1}],
\ee
\be
I_n=\int x^nCi(ax+b)dx =\frac{1}{n+1}[x^nI_0-\frac{nb}{a}I_{n-1}+\frac{n}{a}S_{n-1}].
\ee
Plus integrals over products of $Si$ and $Ci$\ldots
\subsection{The Error Funcation and Fresnel Integrals}
\subsection{Cylinder Functions}

\be
\int x^2 Z_{\nu+1}(x)dx = -2x^2Z_\nu(x)+4\int xZ_\nu(x)dx +\int x^2 Z_{\nu-1}(x)dx,
\ee
by partial integration of $\int x Z_\nu dx$ with \cite[8.471.2]{GR}, where $Z$
is a Bessel Function.

\be
\int x^{\mu+1} Z_{\nu-1}(x)dx = x^{\mu+1}Z_\nu(x)+(\nu-\mu-1)\int x^\mu Z_\nu(x)dx,
\ee
$\mu \neq -1$, by partial integration of $\int x^\mu Z_\nu dx$ with \cite[8.472.1]{GR}, where $Z$
is a Bessel Function.
Equivalent formula for spherical Bessel functions $j_n(z)\equiv \sqrt{\pi/(2z)}J_{n+1/2}(z)$:
\be
(n-m)\int x^m j_n(x)dx = \int x^{m+1}j_{n-1}(x)dx - x^{m+1}j_n(x).
\ee

\be
\int x^{\mu+1} Z_{\nu+1}(x)dx = -x^{\mu+1}Z_\nu(x)+(\mu+1+\nu)\int x^\mu Z_\nu(x)dx,
\ee
$\mu \neq -1$, by partial integration of $\int x^\mu Z_\nu dx$ with \cite[8.472.2]{GR}, where $Z$
is a Bessel Function.

\cite{PiquetteACM24}
\be
\int\frac{ \sin(x)Z_\nu(x)}{x^{3/2}}dx =
\frac{2[(2\nu+1)\sin(x)-2x\cos(x)]}{x^{1/2}(2\nu-1)(2\nu+1)}Z_\nu(x)
-\frac{4x^{1/2}\sin(x)}{(2\nu-1)(2\nu+1)}Z_{\nu+1}(x)
\ee
where $Z$ is a Bessel function $J$ or $Y$.

\cite{PiquetteACM24}
\be
\int\frac{ \cos(x)Z_\nu(x)}{x^{3/2}}dx =
\frac{2[(2\nu+1)\cos(x)+2x\sin(x)]}{x^{1/2}(2\nu-1)(2\nu+1)}Z_\nu(x)
-\frac{4x^{1/2}\cos(x)}{(2\nu-1)(2\nu+1)}Z_{\nu+1}(x)
\ee
where $Z$ is a Bessel function $J$ or $Y$.

\cite{PiquetteACM24}
\begin{multline}
\int\frac{ Z_\mu(x)\bar Z_{\nu}(x)}{x^2}dx =
-\frac{1+\mu+\nu+2\mu\nu+\mu\nu^2+\mu^2\nu-\mu^2-\mu^3-\nu^2-\nu^3+2x^2}
{x(-1+\mu-\nu)(-1+\mu+\nu)(1+\mu-\nu)(1+\mu+\nu)}Z_{\mu}\bar Z_{\nu}
\\
+\frac{1}{(1-\mu-\nu)(1-\mu+\nu)}Z_\mu\bar Z_{\nu+1}
+\frac{1}{(1-\mu-\nu)(1+\mu-\nu)}Z_{\mu+1}\bar Z_{\nu}
\\
-\frac{2x}{(1-\mu-\nu)(1-\mu+\nu)(1+\mu-\nu)(1+\mu+\nu)}Z_{\mu+1}\bar Z_{\nu+1}
.
\end{multline}

\cite{PiquetteACM24}
\begin{multline}
\int\frac{ Z_\mu(x)\bar Z_{\nu}(x)}{x^3}dx =
-\frac{(2+\mu+\nu)(4\mu+4\nu+\mu\nu^2+\mu^2\nu-\mu^3-\nu^3+4x^2}
{x^2(\mu+\nu)[\nu^2-(2-\mu)^2][\nu^2-(2+\mu)^2]}Z_{\mu}\bar Z_{\nu}
\\
-\frac{4\mu\nu^2+2\mu^2\nu^2-4\mu^2-4\nu^3-\mu^4+4\nu^2-\nu^4+8x^2}
{x(\mu^2-\nu^2)[\nu^2-(2-\mu)^2][\nu^2-(2+\mu)^2]}Z_\mu\bar Z_{\nu+1}
\\
+\frac{4\mu^2\nu+2\mu^2\nu^2+4\mu^2-\mu^4-4\nu^2-4\nu^3-\nu^4+8x^2}
{x(\mu^2-\nu^2)[\nu^2-(2-\mu)^2][\nu^2-(2+\mu)^2]}Z_{\mu+1}\bar Z_{\nu}
\\
-\frac{4}{[-\nu^2+(2-\mu)^2][-\nu^2+(2+\mu)^2]}Z_{\mu+1}\bar Z_{\nu+1}
.
\end{multline}

\cite{PiquetteACM24}
\begin{multline}
\int\frac{ Z^2_\nu(x)}{x^2}dx =
\frac{1+2\nu+2x^2}{(4\nu^2-1)x}Z^2_\nu(x)
-\frac{2}{-1+2\nu}Z_\nu(x)Z_{\nu+1}(x)
-\frac{2x}{1-4\nu^2}Z^2_{\nu+1}(x)
.
\end{multline}

\cite{PiquetteACM24}
\begin{multline}
\int\frac{ Z^2_\nu(x)}{x^4}dx =
\frac{-9-6\nu+x^2(6+16\nu+8\nu^2)+36\nu^2+24\nu^3+16x^4}{3x^3(1-4\nu^2)(9-4\nu^2)}Z^2_\nu(x)
\\
-\frac{2(-3+4\nu+4\nu^2+8x^2)}{3x^2(1-2\nu)(9-4\nu^2)}Z_\nu(x)Z_{\nu+1}(x)
-\frac{2(1-4\nu^2-8x^2)}{3x(1-4\nu^2)(9-4\nu^2)}Z^2_{\nu+1}(x)
.
\end{multline}

\cite{PiquetteACM24}
\begin{multline}
\int x^2 Z^3_{1/3}(x)dx =
(-\frac{4}{9}x-\frac{16}{81x})Z^3_{1/3}(x)
-\frac{4x}{3}Z_{1/3}(x)Z^2_{4/3}(x)
+(\frac{8}{9}+x^2)Z^2_{1/3}(x)Z_{4/3}(x)
+\frac{2}{3}x^2 Z^3_{4/3}(x)
.
\end{multline}

\cite{PiquetteACM24}
\begin{multline}
\int \frac{Z^4_1(x)}{x}dx =
\frac{x^2}{4}Z^4_2(x)+(\frac{3}{4}+\frac{x^2}{4})Z^4_1(x)
-\frac{3x}{2}Z_1(x)Z^3_2(x)
+6(\frac{1}{2}+\frac{x^2}{12})Z^2_1(x)Z^2_2(x)
+4(-\frac{3x}{8}-\frac{1}{2x})Z^3_1(x)Z_2(x)
.
\end{multline}

\cite{PiquetteACM24}
\begin{multline}
\int \frac{Z^4_3(x)}{x^3}dx =
(\frac{1}{24}+\frac{1}{2x^2}+\frac{2}{x^4}+\frac{x^2}{378})Z^4_3(x)
+(\frac{5}{216}+\frac{2}{27x^2}+\frac{x^2}{378})Z^4_4(x)
\\
+4(-\frac{x}{108}-\frac{5}{54x}-\frac{1}{3x^3})Z_3(x)Z^3_4(x)
+6(\frac{7}{216}+\frac{1}{3x^2}+\frac{4}{3x^4}+\frac{x^2}{1134})Z^2_3(x)Z^2_4(x)
\\
+4(-\frac{x}{108}-\frac{1}{8x}-\frac{1}{x^3}-\frac{4}{x^5})Z^3_3(x)Z_4(x)
.
\end{multline}
where $Z$ and $\bar Z$ are Bessel functions $J$ or $Y$.

\cite{PiquetteACM24}
\begin{multline}
\int x^l Z_\mu(x)\bar Z_\nu(x) dx=
A_{00}(x)Z_\mu(x)\bar Z_\nu(x)
+A_{01}(x)Z_\mu(x)\bar Z_{\nu+1}(x)
+A_{10}(x)Z_{\mu+1}(x)\bar Z_{\nu}(x)
\\
+A_{11}(x)Z_{\mu+1}(x)\bar Z_{\nu+1}(x)
,
\end{multline}
where $Z$ and $\bar Z$ are Bessel functions $J$ or $Y$,
where
\begin{multline*}
A_{00}=\frac{x}{2(\mu+\nu)}D^3 A_{11}
+\frac{3+\mu+\nu}{2(\mu+\nu)}D^2 A_{11}
+\frac{-7-3\mu-3\nu-2\mu\nu+\mu^2+\nu^2-4x^2}{2x(\mu+\nu)}D A_{11}
\\
+\frac{(-2-\mu-\nu)(-4-2\mu\nu+\mu^2+\nu^2-2x^2)}{2x^2(\mu+\nu)} A_{11}
+\frac{x^{l+1}}{\mu+\nu}
,
\end{multline*}
\begin{multline*}
A_{01}=\frac{-x^2}{2(\mu^2-\nu^2)}D^3 A_{11}
+\frac{3x}{2(\mu^2-\nu^2)}D^2 A_{11}
-\frac{7-3\mu^2-\nu^2+4x^2}{2(\mu^2-\nu^2)}D A_{11}
\\
+\frac{4+\mu\nu^2-3\mu^2-\mu^3-\nu^2+2x^2}{x(\mu^2-\nu^2)} A_{11}
+\frac{x^{l+2}}{\mu^2-\nu^2}
,
\end{multline*}
\begin{multline*}
A_{10}=\frac{x^2}{2(\mu^2-\nu^2)}D^3 A_{11}
-\frac{3x}{2(\mu^2-\nu^2)}D^2 A_{11}
+\frac{7-\mu^2-3\nu^2+4x^2}{2(\mu^2-\nu^2)}D A_{11}
\\
-\frac{4+\mu^2\nu-\mu^2-3\nu^2-\nu^3+2x^2}{x(\mu^2-\nu^2)} A_{11}
-\frac{x^{l+2}}{\mu^2-\nu^2}
,
\end{multline*}
\[
A_{11}(x)=x^{l+3}\sum_{n=0}^{n'} d_n x^{2n}
,
\]
\[
d_0=\frac{2(l+1)}{(l+3)^4-8(l+3)^3+2(12-\mu^2-\nu^2)(l+3)^2
-8(l+3)(4-\mu^2-\nu^2)+[(2-\mu)^2-\nu^2][(2+\mu)^2-\nu^2]
}
,
\]
\begin{multline*}
\{
(3+2n+l)^4-8(3+2n+l)^3+2(12-\mu^2-\nu^2)(3+2n+l)^2
-8(3+2n+l)(4-\mu^2
\\
-\nu^2)+[(2-\mu)^2-\nu^2][(2+\mu)^2-\nu^2]
\}d_n=-4(1+2n+l)(2n+l)d_{n-1}, n'\ge n>0,
\end{multline*}
and $d_n=0$ if $n>n'$.
\[
n'=\left\{
\begin{array}{ll}
0, & l=-1\\
\frac{|l|}{2}-1, & l<-1, even\\
\frac{|l|-3}{2}, & l<-1, odd\\
\infty, & l\ge 0
\end{array}
\right.
\]

\cite{PiquetteACM24}\cite{Auluckarxiv10}
\begin{multline}
\int x^l Z_\nu^2(x) dx =
A_{00}(x)Z_\nu^2(x)+2A_{01}(x)Z_\nu(x) Z_{\nu+1}(x)
+A_{11}(x)Z_{\nu+1}^2(x)
\end{multline}
where
\begin{equation}
A_{00}(x)=\frac{1}{2}D^2A_{11}(x)-\frac{3+2\nu}{2x}DA_{11}(x)+\frac{(2+2\nu+x^2)}{x^2}A_{11}(x),
\end{equation}
\begin{equation}
A_{01}(x)=\frac{1}{2}DA_{11}(x)-\frac{1+\nu}{x}A_{11}(x),
\end{equation}
\begin{equation}
A_{11}(x)=x y(x),
\end{equation}
\begin{equation}
y=\sum_{n=0}^{(l-1)/2}b_nx^{2n+1},
\end{equation}
\begin{equation}
b_{(l-1)/2}=1/(2l),
\quad 
b_n=\frac{2(n+1)[\nu^2-(n+1)^2]}{2n+1}b_{n+1},
\end{equation}
if $0\le n \le(l-3)/2$, and if $l\ge 3$ a positive odd integer.

\cite[136,1]{Watson}
\be
\int z^{-\mu-\nu-1}{\mathcal C}_{\mu+1}(z){\mathcal D}_{\nu+1}(z)dz
= -\frac{z^{-\mu-\nu}}{2(\mu+\nu+1)}\left\{
{\mathcal C}_{\mu}(z){\mathcal D}_\nu(z)
+{\mathcal C}_{\mu+1}(z){\mathcal D}_{\nu+1}(z)
\right\}
\ee
where ${\mathcal C}$ and $\mathcal D$ are arbitrary Bessel functions.

\cite[136,2]{Watson}
\be
\int z^{\mu+\nu+1}{\mathcal C}_{\mu}(z){\mathcal D}_{\nu}(z)dz
= \frac{z^{\mu+\nu+2}}{2(\mu+\nu+1)}\left\{
{\mathcal C}_{\mu}(z){\mathcal D}_\nu(z)
+{\mathcal C}_{\mu+1}(z){\mathcal D}_{\nu+1}(z)
\right\}
\ee
where ${\mathcal C}$ and $\mathcal D$ are arbitrary Bessel functions.

\cite{Watson}
\begin{multline}
(\rho+\mu+\nu)\int z^{\rho-1}{\mathcal C}_{\mu}(z){\mathcal D}_{\nu}(z)dz
+(\rho-\mu-\nu-2)
\int z^{\rho-1}{\mathcal C}_{\mu+1}(z){\mathcal D}_{\nu+1}(z)dz
=
\nonumber
\\
z^{\rho}\left\{
{\mathcal C}_{\mu}(z){\mathcal D}_\nu(z)
+{\mathcal C}_{\mu+1}(z){\mathcal D}_{\nu+1}(z)
\right\}
\end{multline}
where ${\mathcal C}$ and $\mathcal D$ are arbitrary Bessel functions.

\cite[136,5]{Watson}
\begin{multline*}
(\mu+\nu)\int {\mathcal C}_{\mu}(z){\mathcal D}_{\nu}(z)\frac{dz}{z}
-(\mu+\nu+2n)
\int {\mathcal C}_{\mu+n}(z){\mathcal D}_{\nu+n}(z)\frac{dz}{z}
=
\\
{\mathcal C}_{\mu}(z){\mathcal D}_\nu(z)
+2\sum_{m=1}^{n-1}{\mathcal C}_{\mu+m}(z){\mathcal D}_{\nu+m}(z)
+{\mathcal C}_{\mu+n}(z){\mathcal D}_{\nu+n}(z)
\end{multline*}
where ${\mathcal C}$ and $\mathcal D$ are arbitrary Bessel functions.

\cite[137,1]{Watson}
\be
\int {\mathcal C}_{n}(z){\mathcal D}_n(z)\frac{dz}{z}
=
-\frac{1}{2n}\left[
{\mathcal C}_0(z){\mathcal D}_0(z)
+2\sum_{m=1}^{n-1}{\mathcal C}_{m}(z){\mathcal D}_{m}(z)
+{\mathcal C}_{n}(z){\mathcal D}_{n}(z)
\right]
\ee
where ${\mathcal C}$ and $\mathcal D$ are arbitrary Bessel functions.

\cite[138]{Watson}
\begin{multline}
(\mu+2)\int z^{\mu+2}{\mathcal C}_{\nu}^2(z) dz
=
(\mu+1)\left\{\nu^2-\frac{1}{4}(\mu+1)^2\right\}
\int z^\mu {\mathcal C}_\nu^2(z) dz
\\
+\frac{1}{2}\left[
z^{\mu+1}\left\{
z{\mathcal C}_\nu'(z)-\frac{1}{2}(\mu+1){\mathcal C}_\nu(z)
\right\}^2
+z^{\mu+1}
\left\{z^2-\nu^2+\frac{1}{4}(\mu+1)^2\right\}{\mathcal C}_\nu^2(z)
\right]
\end{multline}
where ${\mathcal C}$ and $\mathcal D$ are arbitrary Bessel functions.

With \cite[10.121]{AS}, then partial integration
for a product of three spherical Bessel functions
\begin{multline}
\left(n+m+l+2\right) \int j_n(x)j_m(x)j_l(kx)dx
=
-\left(2n-1\right) j_{n-1}(x) j_m(x) j_l(kx)
\\
+\left(n-3-m-l\right) \int j_{n-2}(x) j_m(x) j_l(kx)dx
+\left(2n-1\right) \int j_{n-1}(x) j_{m-1}(x) j_l(kx) dx
\\
+\left(2n-1\right)k \int j_{n-1}(x) j_m(x) j_{l-1}(kx) dx
.
\end{multline}

\be
\int x J_1(x)dx = \frac{\pi}{2}x[ J_1(x) {\mathbf H}_0(x)-J_0(x){\mathbf H}_1(x)]
\ee
where $\mathbf H$ are Struve functions \cite[\S 12]{AS}.

\cite{ConollyGMJ2}
\be
\frac{1}{2\pi i}\int k^{-2\zeta}\frac{K_{\mu+\zeta}(a)}{a^{\mu+\zeta}}\frac{K_{\nu-\zeta}(b)}{b^{\nu-\zeta}}d\zeta
=\frac12 k^{\mu-\nu}(\frac{k+1/k}{a^2k+b^2/k})^{(\mu+\nu)/2} K_{\mu+\nu}(\sqrt{(k+1/k)(a^2k+b^2/k)}).
\ee

\cite{OkuiMC41}
Let $I_e(k,z)=\int_0^z e^{-x}I_0(kx)dx$, then
\be
\int_0^z e^{-px}I_e(k,x)dx= \frac{1}{p}\left\{
\frac{1}{p+1}I_e(\frac{k}{p+1},[p+1]z)-e^{-pz}I_e(k,z)
\right\}.
\ee

\cite{OkuiMC41}
\begin{multline}
\int_0^z e^{-px}I_e(k,x)xdx= \frac{1}{p}[
\left\{
\frac{1}{p(p+1)}+\frac{1}{(p+1)^2-k^2}\right\}
I_e(\frac{k}{p+1},[p+1]z)-(z+\frac{1}{p})e^{-pz}I_e(k,z)
\\
-\frac{z}{(p+1)^2-k^2}e^{-(p+1)z}\left\{
kI_1(k,z)+(p+1)I_0(kz)
\right\}
].
\end{multline}

\cite{OkuiMC41}
\be
\int_0^z I_e(k,x)dx=
(z-\frac{1}{1-k^2})
I_e(i,z)+\frac{ze^{-z}}{1-k^2}\left\{
I_0(kz)+kI_1(kz)
\right\}
.
\ee

\cite{OkuiMC41}
\be
\int_0^z I_e(k,x)x dx=
\frac12
[
\left\{
z^2-\frac{2+k^2}{(1-k^2)^2}\right\}I_e(k,z)+\frac{ze^{-z}}{1-k^2}
\left\{
(z+\frac{2+k^2}{1-k^2})I_0(kz)
+k(z+\frac{3}{1-k^2})I_1(kz)
\right\}
]
.
\ee

\cite{OkuiMC41}
\be
\int_0^z I_e(k,x)(1-x)^{\nu -1} dx=
\frac{1}{\nu(\nu+1)}
\Phi_2[1/2,1/2;\nu+2;-(1-k),-(1+k)],\quad \nu>0.
.
\ee
plus similar integrals of products of $I_e$ with powers, exponentials, rational
functions, Bessel functions etc.

\subsection{Orthogonal Polynomials}
\cite{PiquetteACM24}
\be
\int\ln(1\pm x)P_\nu(x)dx=
\left[-\frac{x}{\nu}\ln(1\pm x)\pm\frac{1}{\nu(\nu+1)}-\frac{x}{\nu^2}  \right]P_\nu
+\left[ \frac{1}{\nu}\ln(1\pm x)+\frac{1}{\nu^2(\nu+1)}\right]P_{\nu+1}
.
\ee
where $P$ are Legendre functions.

\cite{PiquetteACM24}
\begin{multline}
\int P_\mu(x)\bar P_\nu(x)dx=
-\frac{x}{1+\mu+\nu}P_\mu(x)\bar P_\nu(x)
-\frac{1+\nu}{(\mu-\nu)(1+\mu+\nu)}P_\mu(x)\bar P_{\nu+1}(x)
\\
+\frac{1+\mu}{(\mu-\nu)(1+\mu+\nu)}P_{\mu+1}(x)\bar P_\nu(x)
.
\end{multline}

\cite{PiquetteACM24}
\be
\int x [P_\nu(x)]^2dx=
-\frac{1+\nu}{2\nu}\left[ \frac{x^2+\nu}{1+\nu}[P_\nu(x)]^2
-2xP_\nu(x)P_{\nu+1}(x)+[P_{\nu+1}(x)]^2
\right]
.
\ee

\cite{PiquetteACM24}
\be
\int [P_{1/2}(x)]^2dx=
\frac{9}{2}x[P_{3/2}(x)]^2
+\frac{x(-7+16x^2)}{2}[P_{1/2}(x)]^2
-3(-1+2x)(1+2x)P_{1/2}(x)P_{3/2}(x)
.
\ee

\cite{PiquetteACM24}
\begin{multline}
\int [P_{3/2}(x)]^2dx=
\frac{x(93-480x^2+512x^4)}{18}x[P_{3/2}(x)]^2
+\frac{25x(-3+8x^2)}{18}[P_{5/2}(x)]^2
\\
-\frac{5}{9}(3-42x^2+64x^4)P_{3/2}(x)P_{5/2}(x)
.
\end{multline}

\cite{PiquetteACM24}
\begin{multline}
\int x^5 P_{1/3}(x) P_{2/3}(x) dx=
(-\frac{335}{2352}-\frac{512x^2}{392}+\frac{235x^4}{336}+\frac{x^6}{3})P_{1/3}(x)P_{2/3}(x)
\\
+(\frac{685x}{784}-\frac{573x^3}{1176}-\frac{5x^5}{48})P_{1/3}(x)P_{5/3}(x)
+(\frac{235x}{196}-\frac{295x^3}{588}-\frac{2x^5}{21})P_{2/3}(x)P_{4/3}(x)
\\
+(-\frac{5}{6}+\frac{65x^2}{196}+\frac{25x^4}{588})P_{4/3}(x)P_{5/3}(x)
.
\end{multline}

\cite{PiquetteACM24}
\begin{multline}
\int x [P_{1/3}(x)]^3 dx=
(\frac{125x^4}{12}-\frac{14}{3}x^2-\frac{5}{12})[P_{1/3}(x)]^3
+(-4+20x^2)P_{1/3}(x)[P_{4/3}(x)]^2
\\
+(9x-25x^3)[P_{1/3}(x)]^2 P_{4/3}(x)
-\frac{16}{3}x [P_{4/3}(x)]^3
.
\end{multline}

\cite{PiquetteACM24}
\begin{multline}
\int x [P_{1/2}(x)]^4 dx=
(-\frac{5}{16}-\frac{19}{4}x^2)[P_{1/2}(x)]^4
+\frac{81}{4}xP_{1/2}(x)[P_{3/2}(x)]^3
\\
+6( -\frac{9}{16}-\frac{9}{2}x^2)[P_{1/2}(x)]^2[P_{3/2}(x)]^2
+4( \frac{33x}{16}+3x^3)[P_{1/2}(x)]^3P_{3/2}(x)
-\frac{81}{16}[P_{3/2}(x)]^4
.
\end{multline}

\cite{PiquetteACM24}
\begin{multline}
\int x^l P_\mu(x) P_\nu(x)dx = A_{00}(x)P_\mu(x)P_\nu(x)
+ A_{01}(x)P_\mu(x)P_{\nu+1}(x)
\\
+ A_{10}(x)P_{\mu+1}(x)P_{\nu}(x)
+ A_{11}(x)P_{\mu+1}(x)P_{\nu+1}(x),
\end{multline}

\cite{PiquetteACM24}
\be
\int e^{-x^2}H_\nu(x)dx=-e^{-x^2}H_{\nu-1}(x).
\ee

\cite{PiquetteACM24}
\be
\int H_\nu(x) x^{-(\nu+3)}dx=\left[\frac{2x^{-\nu}}{(\nu+1)(\nu+2)}
-\frac{x^{-(\nu+2)}}{\nu+2}\right]H_{\nu}(x)
- \frac{2\nu}{(\nu+1)(\nu+2)}x^{-(\nu+1)}H_{\nu-1}(x).
\ee

\cite{PiquetteACM24}
\be
\int xH_\nu(x) dx=\frac{1+2x^2}{2(\nu+2)}H_{\nu}(x)
- \frac{\nu x}{\nu+2}H_{\nu-1}(x).
\ee

\cite{PiquetteACM24}
\be
\int e^{-x^2} H_\mu(x) \bar H_{\nu} dx=\frac{e^{-x^2}}{2(\mu-\nu)}
[-H_{\mu}(x)\bar H_{\nu+1}(x)+H_{\mu+1}(x)\bar H_{\nu}(x)
]
.
\ee

\cite{PiquetteACM24}
\begin{multline}
\int xe^{-x^2} H_\mu(x) \bar H_{\nu} dx=\frac{e^{-x^2}}{2}
[
-\frac{1+\mu+\nu}{(1-\mu+\nu)(1+\mu-\nu)}H_{\mu}(x)\bar H_{\nu+1}(x)
+\frac{x}{1+\mu-\nu}H_{\mu+1}(x)\bar H_{\nu}(x)
\\
+\frac{x}{1-\mu+\nu}H_{\mu}(x)\bar H_{\nu+1}(x)
-\frac{1}{(1-\mu+\nu)(1+\mu-\nu)}H_{\mu+1}(x)\bar H_{\nu+1}(x)
]
.
\end{multline}

\cite{PiquetteACM24}
\begin{multline}
\int x^2 e^{-x^2} H_\mu(x) \bar H_{\nu} dx=
-e^{-x^2}\frac{H_\mu(x)\bar H_{\nu}(x) x(\mu+\nu)}{(2-\mu+\nu)(2+\mu-\nu)}
\\
+H_{\mu+1}(x)\bar H_\nu(x)
\frac{2+\mu+3\nu+2\mu x^2-2\nu x^2-\mu^2x^2-\nu^2x^2+2\mu\nu x^2}{2(\mu-\nu)
(2-\mu+\nu)(2+\mu-\nu) }
\\
+H_{\mu}(x)\bar H_{\nu+1}(x)
\frac{2+3\mu+\nu-2\mu x^2+2\nu x^2-\mu^2x^2-\nu^2x^2+2\mu\nu x^2}{2(\mu-\nu)
(2-\mu+\nu)(2+\mu-\nu) }
\\
-H_{\mu+1}(x)\bar H_{\nu+1}(x)
\frac{x}{
(2-\mu+\nu)(2+\mu-\nu) }
.
\end{multline}

\cite{PiquetteACM24}
\begin{multline}
\int e^{-3x^2} x^2 H^3_{2/3}(x) dx=e^{-3x^2}
[
-\frac{1}{12}x(5+6x^2)H^3_{2/3}(x)
+\frac{1}{8}(1+6x^2)H^2_{2/3}(x)H_{5/3}(x)
\\
-\frac{3}{8}xH_{2/3}(x)H^2_{5/3}(x)
+\frac{1}{16}H^3_{5/3}(x)
]
.
\end{multline}
where $H$ and $\bar H$ are Hermite functions.

\cite{PiquetteACM24}
\be
\int x e^{-(\nu+1)x}L_\nu(x) dx=\frac{e^{-(\nu+1)x}}{\nu+1}
[-(1+x)L_\nu(x)+L_{\nu-1}(x)].
\ee
\cite{PiquetteACM24}
\be
\int x (1+x)^{-(\nu+3)}L_\nu(x) dx=\frac{(1+x)^{-(\nu+1)}}{\nu+2}
\left[ \left(\frac{x-\nu}{\nu+1}-\frac{x}{1+x}\right)L_\nu(x)+\frac{\nu}{\nu+1}L_{\nu-1}(x)\right]
.
\ee

\cite{PiquetteACM24}
\be
\int e^{-x} L_\mu(x) \bar L_{\nu}(x) dx=
e^{-x}[L_\mu(x)\bar L_\nu(x)
+\frac{1+\nu}{\mu-\nu}L_\mu(x)\bar L_{\nu+1}(x)
-\frac{1+\mu}{\mu-\nu}L_{\mu+1}(x)\bar L_{\nu}(x)
]
.
\ee

\cite{PiquetteACM24}
\begin{multline}
\int x e^{-x} L_\mu(x) \bar L_{\nu}(x) dx=
e^{-x}[
-\frac{1+\mu+\nu-x+2\mu\nu+\mu^2x+\nu^2 x-2\mu\nu x}{(1-\mu+\nu)(1+\mu-\nu)} L_\mu(x)\bar L_\nu(x)
\\
+\frac{(1+\nu)(1+2\mu-\mu x+\nu x)}{(\mu-\nu)(1-\mu+\nu)} L_\mu(x)\bar L_{\nu+1}(x)
-\frac{(1+\mu)(1+2\nu+\mu x-\nu x)}{(\mu-\nu)(1+\mu-\nu)} L_{\mu+1}(x)\bar L_{\nu}(x)
\\
-\frac{2(1+\mu)(1+\nu)}{(1-\mu+\nu)(1+\mu-\nu)} L_{\mu+1}(x)\bar L_{\nu+1}(x)
]
.
\end{multline}

\cite{PiquetteACM24}
\begin{multline}
\int e^{-3x} x L^3_{2/3}(x) dx=
e^{-3x}\Big\{
L^3_{2/3}(x)[
\frac{125}{24}-\frac{625x}{24}+\frac{853x^2}{16}-\frac{675x^3}{16}+\frac{225x^4}{16}-\frac{27x^5}{16}]
\\
+L^3_{5/3}(x)[
-\frac{125}{24}+\frac{125x}{12}-\frac{125x^2}{16}]
+3[\frac{125}{24}-\frac{125x}{8}+\frac{275x^2}{16}-\frac{75x^3}{16}]L_{2/3}(x)L^2_{5/3}(x)
\\
+3[-\frac{125}{24}+\frac{125x}{6}-\frac{515x^2}{16}+\frac{135x^3}{8}-\frac{45x^4}{16}]
L^2_{2/3}(x)L_{5/3}(x)
\Big\}
.
\end{multline}
where $L$ and $\bar L$ are Laguerre functions.

\section{Definite Integrals of Special Functions}

\subsection{Elliptic Integrals and Functions}

\cite{WanAAM48}
\be
\int_0^1 x^nK'(x)dx=\frac{\pi}{4} \frac{\Gamma(\frac{n+1}{2})^2}{\Gamma(\frac{n+2}{2})^2}.
\ee
\cite{WanAAM48}
\be
\int_0^1 x^nE'(x)dx=\frac{\pi}{2(n+1)}
\frac{\Gamma(\frac{n+3}{2})^2}{\Gamma(\frac{n+2}{2})\Gamma(\frac{n+4}{2})}.
\ee

\cite{WanAAM48}
\be
\int_0^1 x^{'n}x^{m}K(x)dx = \frac{\pi}{4}
\frac{\Gamma(\frac{m+1}{2})\Gamma(\frac{n+2}{2})}{\Gamma(\frac{m+n+3}{2})}
{}_3F_2(\frac12,\frac12,\frac{m+1}{2};1,\frac{m+n+3}{2}\mid 1).
\ee
where $x'=\sqrt{1-x^2}$, $K'(x)=K(x')$, $E'(x)=E(x')$.

\cite{WanAAM48}
\be
\int_0^1 x^{'n}x^{m}E(x)dx = \frac{\pi}{4}
\frac{\Gamma(\frac{m+1}{2})\Gamma(\frac{n+2}{2})}{\Gamma(\frac{m+n+3}{2})}
{}_3F_2(-\frac12,\frac12,\frac{m+1}{2};1,\frac{m+n+3}{2}\mid 1).
\ee

\cite{WanAAM48}
\be
\int_0^1 x^{2n+1}K(x)^aE(x)^bK'(x)^cE'(x)^ddx = 
\int_0^1 x(1-x^2)^n K'(x)^aE'(x)^bK(x)^cE(x)^ddx
\label{eq.WanEI}
\ee

\cite{WanAAM48}
\be
\int_0^1 x'K(x)^2dx =
\int_0^1\frac{\surd x}{x+1}K'(x)^2dx
=\frac{\pi^3}{16}[2{}_4F_3(\frac12,\frac12,\frac12,\frac12;1,1,1\mid 1)-1].
\ee

\cite{WanAAM48}
\be
\int_0^1 x^nE'(x)^2dx = \frac{2^{4n}(n+1)^3(n+3)^2}{16(n+2)^3(n+4)}
\frac{\Gamma^8(\frac{n+1}{2})}{\Gamma^4(n+1)}
{}_7F_6\left(
\begin{array}{c}
-\frac12,\frac12,\frac12,\frac32,\frac{n+3}{2},\frac{n+3}{2},\frac{n+7}{4}\\
1,\frac{n+3}{4},\frac{n+2}{2},\frac{n+4}{2},\frac{n+4}{2},\frac{n+6}{2}
\end{array}\mid 1
\right).
\ee

\cite{WanAAM48}
\be
\int_0^1 x^nK(x)K'(x)dx =\frac{\pi^2}{8} \frac{\Gamma(\frac{n+1}{2})^2}{\Gamma(\frac{n+2}{2})^2}
{}_4F_3(\frac12,\frac12,\frac{n+1}{2},\frac{n+1}{2};
1\frac{n+2}{2},\frac{n+2}{2}\mid 1).
\ee

\cite{WanAAM48}
\be
\int_0^1 xK(x)K'(x)dx
=\int_0^1 2x^3K(x)K'(x)dx
=\int_0^1\frac{1-x}{1+x}K(x)K'(x)dx=\frac{\pi^3}{16}.
\ee

\cite{WanAAM48}
\be
\int_0^1 x^nE(x)K'(x)dx =\frac{\pi^2}{8} \frac{\Gamma(\frac{n+1}{2})^2}{\Gamma(\frac{n+2}{2})^2}
{}_4F_3(-\frac12,\frac12,\frac{n+1}{2},\frac{n+1}{2};
1\frac{n+2}{2},\frac{n+2}{2}\mid 1).
\ee

\cite{WanAAM48}
\be
\int_0^1 x^nK(x)E'(x)dx =\frac{\pi^2}{8} \frac{(n+1)\Gamma(\frac{n+1}{2})^2}{(n+2)\Gamma(\frac{n+2}{2})^2}
{}_4F_3(\frac12,\frac12,\frac{n+1}{2},\frac{n+3}{2};
1\frac{n+2}{2},\frac{n+4}{2}\mid 1).
\ee

\cite{WanAAM48}
\be
\int_0^1 x^nE(x)E'(x)dx =\frac{\pi^2}{8} \frac{(n+1)\Gamma(\frac{n+1}{2})^2}{(n+2)\Gamma(\frac{n+2}{2})^2}
{}_4F_3(-\frac12,\frac12,\frac{n+1}{2},\frac{n+3}{2};
1\frac{n+2}{2},\frac{n+4}{2}\mid 1).
\ee

\cite{WanAAM48}
\be
\int_0^1 x^nE'(x)K'(x)dx = \frac{2^{4n}(n+1)^2}{16(n+2)}
\frac{\Gamma^8(\frac{n+1}{2})}{\Gamma^4(n+1)}
{}_7F_6\left(
\begin{array}{c}
-\frac12,\frac12,\frac12,\frac12,\frac{n+1}{2},\frac{n+1}{2},\frac{n+5}{4}\\
1,\frac{n+1}{4},\frac{n+2}{2},\frac{n+2}{2},\frac{n+2}{2},\frac{n+4}{2}
\end{array}\mid 1
\right).
\ee

Recurrences for moments:
\cite{WanAAM48}
\be
(n+1)^3K_{n+2}-2n(n^2+1)K_n+(n-1)^3K_{n-2}=2,
\ee
where $K_n\equiv \int_0^1 x^nK(x)^2dx$.
\cite{WanAAM48}
\be
(n+1)(n+3)(n+5)E_{n+2}-2(n^3+3n^2+n+1)E_n+(n-1)^3E_{n-2}=8
\ee
where $E_n\equiv \int_0^1 x^nE(x)^2dx$.

\cite{WanAAM48}
\be
\int_0^1 xK'(x)^2dx = \frac74\zeta(3).
\ee
\cite{WanAAM48}
\be
\int_0^1 x^3K(x)^2dx = \frac18(2+7\zeta(3)).
\ee
and further odd moments via \eqref{eq.WanEI}.
\cite{WanAAM48}
\be
\int_0^1 xK'(x)K(ix)dx = \frac12\pi G,
\ee
where $G$ is Catalan's constant.

\cite{WanAAM48}
\be
\int_0^1 x^nK'(x)^2dx = \frac{2^{4n}(n+1)}{16}
\frac{\Gamma^8(\frac{n+1}{2})}{\Gamma^4(n+1)}
{}_7F_6\left(
\begin{array}{c}
\frac12,\frac12,\frac12,\frac12,\frac{n+1}{2},\frac{n+1}{2},\frac{n+5}{4}\\
1,\frac{n+1}{4},\frac{n+2}{2},\frac{n+2}{2},\frac{n+2}{2},\frac{n+4}{2}
\end{array}\mid 1
\right).
\ee

\cite{WanAAM48}
\be
\int_0^1 \frac{K(x)}{1+x}dx = \frac{\pi^2}{8}.
\ee

\cite{WanAAM48}
\be
\int_0^1 \frac{1}{1-t^2x^2}K'(x)dx = \frac{\pi}{2}K(t).
\ee
\cite{WanAAM48}
\be
\int_0^1 \frac{1}{1-t^2x^2}E'(x)dx = \frac{\pi}{2t^2}[K(t)-E(t)].
\ee
\cite{WanAAM48}
\be
\int_0^1 \frac{x}{1-t^2x^2}K(x)K'(x)dx = \frac{\pi}{4}K(t)^2.
\ee

\cite{WanAAM48}
\be
\int_0^1 \frac{E(x)}{(1+x)^3}dx = \frac{\pi^2}{32}+\frac14.
\ee
\cite{WanAAM48}
\be
\int_0^1 \frac{E'(x)}{(1+x)^3}dx = \frac{G}{8}+\frac{6}{16}.
\ee

\cite{WanAAM48}
\be
\int_0^1 \frac{K(x)-E(x)}{xx'}dx =\frac{\pi}{2}.
\ee

\cite{Adamchik}
\be
\frac{1}{2}\int_0^{1}K(x^2)dx=G.
\ee

\cite{Adamchik}
\be
\int_0^{1}E(x^2)dx=G+\frac{1}{2}.
\ee

\cite{Rogersarxiv1303}
\be
\int_0^1 \frac{K(k)}{\sqrt{1-k^2}}dk=\frac{\Gamma^4(\frac14)}{16\pi}.
\ee

\cite{Rogersarxiv1303}
\be
\int_0^1 \frac{K'(x)}{\sqrt{1+x}}dx
=
\frac{1}{2}\left[
\frac{\Gamma^2(\frac18)\Gamma^2(\frac38)}{16\pi}-{}_4F_3(\frac34,1,1,\frac54;\frac32,\frac32,\frac32;1)
\right].
\ee

\cite{Rogersarxiv1303}
\be
\int_0^1 \frac{K'(x)}{\sqrt{1-x}}dx
=
\frac{1}{2}\left[
\frac{\Gamma^2(\frac18)\Gamma^2(\frac38)}{16\pi}+{}_4F_3(\frac34,1,1,\frac54;\frac32,\frac32,\frac32;1)
\right].
\ee

\cite{Rogersarxiv1303,RogersRJ37}
\be
\int_0^1 K'(k)^3dk
=
3\int_0^2 K(k)^2K'(k)dk
=
\frac{\Gamma^8(\frac14)}{128\pi^2}.
\ee

\cite{Rogersarxiv1303,RogersRJ37}
\be
\int_0^1 \frac{K'(k)^3}{\sqrt k(1-k^2)^{3/4}}dk=\frac{3\Gamma^8(\frac14)}{32\surd 2 \pi^2}.
\ee

\cite{Rogersarxiv1303}
\be
\int_0^1 E(k)K'(k)^2dk = \frac{\pi^3}{12}+\frac{\Gamma^8(\frac14)}{384\pi^2}.
\ee

\cite{Rogersarxiv1303}
\be
\int_0^1 \frac{K(k)^2}{\sqrt{1-k^2}}dk
=
\frac{\pi^3}{4}{}_4F_3(\frac12,\frac12,\frac12,\frac12;1,1,1;1).
\ee

\cite{Rogersarxiv1303}
\be
\int_0^1 K(k)^2dk
=
\frac{\pi^4}{32}{}_7F_6(\frac54,\frac12,\frac12,\frac12,\frac12,\frac12,\frac12;\frac14,1,1,1,1,1;1).
\ee

\cite{WanAAM48}
\be
(n+k+1)\int_0^1 [x^kE(x)^n-nx^kE(x)^{n-1}K(x)]dx = 1.
\ee

\cite{WanAAM48}
\be
\int_0^1 [5x^2E(x)^2-2E(x)K(x)]dx = 1.
\ee
\cite{WanAAM48}
\be
\int_0^1 [2E'(x)K(x)-(1-x^2)K(x)K'(x)]dx = \frac{\pi}{2}.
\ee
and similar partial integrations of products.

\cite{WanAAM48}
\be
\int_0^1 \frac{\arctan x}{x}K'(x)dx = \pi G,
\ee
where $G$ is Catalan's constant.

\cite{WanAAM48}
\be
\int_0^1 -\frac{\log(1-x^2)}{x}K'(x)dx = \frac78 \pi \zeta(3).
\ee
\cite{WanAAM48}
\begin{multline}
\int_0^{\pi/2}
K(\sin t)^2dt
=2\int_0^{\pi/2}
K(\sin t)E(\sin t)dt
=\int_0^1\frac{K(x)^2}{\sqrt{1-x^2}}dx\\
=2\int_0^1\frac{E(x)^2}{\sqrt{1-x^2}}dx
=2\int_0^1 K(x)K'(x)dx
=
\frac{\pi^3}{4}{}_4F_3
(\frac12,\frac12,\frac12,\frac12;
1,1,1\mid 1).
\end{multline}

\cite{WanAAM48}
\begin{multline}
\int_0^{\pi /2}K(\sin t)dt
=\int_0^1 \frac{K(x)}{\sqrt(1-x^2)}dx =\int_0^1 \frac{K'(x)}{\sqrt{1-x^2}}dx
=\int_0^1 \frac{K(x)}{\surd x}dx\\
=\frac12 \int_0^1 \frac{K'(x)}{\surd x}dx = \frac{1}{16\pi}\Gamma(1/4)^4.
\end{multline}

\cite{MacrobertGJM2}
\begin{multline}
\int_0^\infty \lambda^{l-1}K_m(\lambda)K_n(z/\lambda)d\lambda
=
\sum_{n,-n} \Gamma(\frac{l+m+n}{2})\Gamma(\frac{l-m+n}{2})\Gamma(n)2^{l+2n-3}z^{-n}
\\ \times
F(;1-n,1-\frac{l}{2}-\frac{m}{2}-\frac{n}{2},
1-\frac{l}{2}+\frac{m}{2}-\frac{n}{2};z^2/16)
\\
+\sum_{m,-m} \Gamma(\frac{-l-m-n}{2})\Gamma(\frac{-l-m+n}{2})\Gamma(-m)2^{-l-2m-3}z^{l+m}
F(;1+m,1+\frac{l}{2}+\frac{m}{2}+\frac{n}{2},
1+\frac{l}{2}+\frac{m}{2}-\frac{n}{2};z^2/16).
\end{multline}

\subsection{The Exponential Integral and Related Functions}
\cite{BorosSci12}
\be
\int_0^\infty t^{2\beta}e^{-3t^2}\erf(t)dt
=
\frac{\Gamma(1+\beta)}{\sqrt{\pi}3^{\beta+1}}\,_2F_1(\frac{1}{2},1+\beta;\frac{3}{2};-\frac{1}{3})
.
\ee

\subsection{The Gamma Function and Related Functions}
\cite{Amdeberhanarxiv3663}\cite[A075700]{sloane}
\be
\int_0^1\ln \Gamma(t)dt=\frac{1}{2}\ln 2\pi.
\ee

\cite{Amdeberhanarxiv3663}
\be
\int_0^1 t\ln \Gamma(t)dt=\frac{\zeta'(2)}{2\pi^2}+\frac{1}{6}\ln 2\pi-\frac{\gamma}{12}.
\ee

\cite{Amdeberhanarxiv3663}
\be
\int_0^\infty 2^{-t}\ln \Gamma(t)dt=
2\int_0^1 2^{-t}\ln \Gamma(t)dt-\frac{\gamma+\ln\ln 2}{\ln 2}.
\ee

\cite{Amdeberhanarxiv3663}
\be
\int_0^\infty 2^{-t}t\ln \Gamma(t)dt=
2\int_0^1 2^{-t}(t+1)\ln \Gamma(t)dt-\frac{(\gamma+\ln\ln 2)(1+2\ln 2)-1}{\ln^2 2}.
\ee

\cite{EspinosaRJ6}
\be
\int_0^1 \ln \Gamma(q) \cos((2n+1)\pi q)dq
=
\frac{2}{\pi^2}
\left(
\frac{\gamma+2\ln\sqrt{2\pi}}{(2n+1)^2}
+2\sum_{k=1}^\infty \frac{\ln k}{4k^2-(2n+1)^2}
\right)
.
\ee

\cite{EspinosaRJ6}
\be
\int_0^1 B_{2m}(q)\ln\Gamma(q) dq =
(-)^{m+1}\frac{(2m)!\zeta(2m+1)}{2(2\pi)^{2m}}
.
\ee

\cite{EspinosaRJ6}
\be
\int_0^1 B_{2m-1}(q)\ln\Gamma(q) dq =
\frac{B_{2m}}{2m}\left[
\frac{\zeta'(2m)}{\zeta(2m)}-A
\right]
,
\ee
where $A=2\ln\sqrt{2\pi}+\gamma$.

\cite{EspinosaRJ6}
\begin{multline}
\int_0^1 q^n \ln\Gamma(q) dq =
\frac{1}{n+1}\sum_{k=1}^{\lfloor (n+1)/2 \rfloor}
(-)^k\binom{n+1}{2k-1}\frac{(2k!)}{k(2\pi)^{2k}}
[A\zeta(2k)-\zeta'(2k)]
\\
-
\frac{1}{n+1}\sum_{k=1}^{\lfloor n/2 \rfloor}
(-)^k\binom{n+1}{2k}\frac{(2k!)}{2(2\pi)^{2k}}
\zeta(2k+1)
+\frac{\ln\sqrt{2\pi}}{n+1}
,
\end{multline}
where $A=2\ln\sqrt{2\pi}+\gamma$.

\cite{EspinosaRJ6}
\be
\int_0^1 (q-1/2)\ln\Gamma(q) dq =
\frac{1}{12}\left(
\frac{6\zeta'(2)}{\pi^2}-2\ln\sqrt{2\pi}-\gamma
\right)
.
\ee

\cite{EspinosaRJ6}
\be
\int_0^{1/2}\ln\Gamma(q+1)dq
=\frac{\gamma}{8}+\frac{3\ln\sqrt{2\pi}}{4}
-\frac{13\ln 2}{24}-\frac{3\zeta'(2)}{4\pi^2}-\frac{1}{2}
.
\ee

\cite{EspinosaRJ6}
\begin{multline}
\int_0^1 q^n \psi^{(m)}(q)dq
=(-1)^m\frac{n!}{(n-m)!}
\big[
\frac{\gamma}{n-m+1}+(n-m)!\sum_{k=0}^{m-2}
\frac{\Gamma(m-k)\zeta(m-k)}{(n-k)!}
\\
+\sum_{k=0}^{n-m-1}(-)^k\binom{n-m}{k}[H_k \zeta(-k)+\zeta'(-k)]
\big]
,
\end{multline}
where $H_k$ are harmonic sums.

\cite{EspinosaRJ6}
\be
\int_0^1 q^n \psi(q)dq
=\zeta'(0)
+\sum_{k=1}^{n-1}(-)^k\binom{n}{k}[H_k \zeta(-k)+\zeta'(-k)]
,
\ee
where $H_k$ are harmonic sums.

\cite{MuthumSci22}
\be
\int_0^\infty \sin(ax)[\psi(z+ix)-\psi(z-ix)]dx=i\pi \frac{e^{-az}}{1-e^{-a}}.
\ee

\cite{MuthumSci22}
\be
\int_0^\infty \sin 2a \cos^2 bx [\psi(z+ix)-\psi(z-ix)]dx=\frac{i\pi}{4}[\frac{2e^{-2az}}{1-e^{-2a}}+
\frac{e^{-2(a+b)z}}{1-e^{-2(a+b)}}
-\frac{e^{-2(a-b)z}}{1-e^{-2(a-b)}}
], a<b.
\ee

\cite{MuthumSci22}
\be
\int_0^\infty \sin 2a \cos^2 ax [\psi(z+ix)-\psi(z-ix)]dx=\frac{i\pi}{4}[\frac{e^{-4az}}{1-e^{-4a}}+
2\frac{e^{-2az}}{1-e^{-2a}}
].
\ee

\cite{ApelblatJAM34}
\be
\int_0^\infty \frac{t^{\alpha+z}\psi(\alpha+z+1)}{\Gamma(\alpha+z+1)}dz
=
\frac{t^\alpha}{\Gamma(\alpha+1)}+\nu(t,\alpha)\ln t,\quad
\Re \alpha>-1.
\ee
\be
\int_0^\infty \frac{\psi(\alpha+z+1)}{\Gamma(\alpha+z+1)}dz
=
\frac{1}{\Gamma(\alpha+1)},\quad
\Re \alpha>-1.
\ee
\be
\int_1^\infty \frac{\psi(z)}{\Gamma(z)}dz = 
1.
\ee
\be
\int_0^\infty \frac{t^z \psi(z+1)}{\Gamma(z+1)}dz = 
1+\nu(t)\ln t, t\neq 1.
\ee

\cite{ApelblatJAM34}
\be
\int_0^\infty \frac{t^{\alpha+z}}{\Gamma(\alpha+z+1)}
\{\psi(\alpha+z+1)^2-\psi^{(1)}(\alpha+z+1)\}
dz = 
\frac{2t^\alpha\ln t}{\Gamma(\alpha+1)}
+(\ln t)^2 \nu(t,\alpha)+L^{-1}\left\{\frac{\ln s}{s^{\alpha+1}}\right\}
,\quad \Re \alpha>-1
\ee
where $L^{-1}$ is the inverse Laplace transform.
\be
\int_0^\infty \frac{t^z}{\Gamma(z+1)}
\{\psi(z+1)^2-\psi^{(1)}(z+1)\}
dz = 
-\gamma +\ln t (1+\nu(t)\ln t)
.
\ee
\be
\int_1^\infty \frac{1}{\Gamma(z)}
\{\psi^{(1)}(z)-\psi(z)^2\}
dz = 
\gamma
.
\ee
\be
\int_0^\infty \frac{t^{\alpha+z}}{\Gamma(\alpha+z+1)}
\{\psi(\alpha+z+1)^2-\psi^{(1)}(\alpha+z+1)\}
dz = 
\frac{t^\alpha}{\Gamma(\alpha+1)}
[\psi(\alpha+1)+\ln t]+(\ln t)^2\nu(t,\alpha)
.
\ee
\be
\int_0^\infty \frac{t^{\alpha+z}z^\beta }{\Gamma(\alpha+z+1)}
\,\frac{\psi(\alpha+z+1)}{\Gamma(\beta+1)}
dz = 
\mu(t,\beta-1,\alpha)+\ln t \mu(t,\beta,\alpha),
\quad
\Re \alpha>-1,\Re \beta>-1
.
\ee
\begin{multline}
\int_0^\infty \frac{t^{\alpha+z}z^\beta }{\Gamma(\alpha+z+1)}
\,\frac{\psi(\alpha+z+1)^2-\psi^{(1)}(\alpha+z+1)}{\Gamma(\beta+1)}
dz
\\
= 
\mu(t,\beta-2,\alpha)+2\ln t \mu(t,\beta-1,\alpha)
+(\ln t)^2\mu(t,\beta,\alpha),
\quad
\Re \alpha>-1,\Re \beta>1
.
\end{multline}

\cite{ApelblatJAM34}
Let $L\{f(t)\}\equiv \int_0^\infty e^{-st}f(t)dt = F(s)$ be the Laplace transform
and
\be
\nu(z)\equiv \int_0^\infty \frac{z^t dt}{\Gamma(t+1)},
\ee
\be
\nu(z,\alpha)\equiv \int_0^\infty \frac{z^{\alpha+t} dt}{\Gamma(\alpha+t+1)},
\ee
\be
\mu(z,\beta)\equiv \int_0^\infty \frac{z^t t^\beta dt}
{\Gamma(\beta+1)\Gamma(t+1)},
\ee
\be
\mu(z,\beta,\alpha)\equiv \int_0^\infty \frac{z^{\alpha+t} t^\beta dt}
{\Gamma(\beta+1)\Gamma(\alpha+t+1)},
\ee
then
\be
L\{\nu (t)\}=\frac{1}{s\ln s}.
\ee
\cite{ApelblatJAM34}
\be
L\{\nu (t,\alpha)\}=\frac{1}{s^{1+\alpha}\ln s}.
\ee
\cite{ApelblatJAM34}
\be
L\{\mu (t,\beta)\}=\frac{1}{s(\ln s)^{\beta+1}}.
\ee
\cite{ApelblatJAM34}
\be
L\{\mu (t,\beta,\alpha)\}=\frac{1}{s^{\alpha+1}(\ln s)^{\beta+1}}.
\ee
Above $\Re \alpha>-1$, $\Re\beta>-1$, $\Re s>1$.
\cite{ApelblatJAM34}
\be
\frac{F(\ln s)}{s\ln s}=L\{\int_0^\infty \nu(t,x)f(x)dx\}
.
\ee
\cite{ApelblatJAM34}
\be
\frac{F(\ln s)}{s^{\alpha+1}\ln s}=L\{\int_0^\infty \nu(t,\alpha+x)f(x)dx\}
.
\ee
\cite{ApelblatJAM34}
\be
\frac{F(\ln s)}{s(\ln s)^{\beta+1}}=L\{\int_0^\infty \mu(t,\beta,x)f(x)dx\}
.
\ee
\cite{ApelblatJAM34}
\be
\frac{F(\ln s)}{s^{\alpha+1}(\ln s)^{\beta+1}}=L\{\int_0^\infty \mu(t,\beta,\alpha+x)f(x)dx\}
.
\ee
\cite{ApelblatJAM34}
\be
\int_0^\infty x^{\beta-1}\nu(t,\alpha+x)dx
=
\Gamma(\beta)\mu(t,\beta,\alpha),
\quad \Re\alpha>-1,\Re \beta>0
.
\ee
\cite{ApelblatJAM34}
\be
\int_0^\infty x^{\beta-1}\nu(t,x)dx
=
\Gamma(\beta)\mu(t,\beta)
.
\ee
\cite{ApelblatJAM34}
\be
\int_0^\infty \nu(t,\alpha+x)dx
=
\mu(t,1,\alpha)
.
\ee
\cite{ApelblatJAM34}
\be
\int_0^\infty \nu(t,x)dx
=
\mu(t,1)
.
\ee
\cite{ApelblatJAM34}
\be
\int_0^\infty t^{\lambda-1}\mu(t,\beta,\alpha+x)dx
=
\Gamma(\lambda)\mu(t,\beta+\lambda,\alpha),
\quad \Re\alpha>-1,\Re \beta>-1,\Re\lambda >0
.
\ee
\cite{ApelblatJAM34}
\be
\frac{1}{3\pi}\int_0^\infty (\frac{x}{t})^{1/2}
K_{1/2}(\phi)\mu(x,\beta,\alpha)dx
=
3^{\beta+1}\mu(t,\beta,\alpha/3), \phi\equiv 2(x^3/27t)^{1/2},
\ee
and others,

\subsection{Cylinder Functions}

\cite{NollJOSA66}
\be
2\pi(-1)^{(n-m)/2}\int_0^\infty dk J_{n+1}(2\pi k)J_m(2\pi kr)
=
R_n^m(r)
=
\sum_{s=0}^{(n-m)/2} (-1)^s\binom{n-s}{s}\binom{n-2s}{(n-m)/2-s}r^{n-2s},
\ee
for $n\ge 0$, $0\le m\le n$, $n-m$ even.

\cite{Reynoldsarxiv2502}
\be
\int_0^{\pi/2} J_0(\lambda \sin\theta)^2 d\theta
=\surd\pi \sum_{j\ge 0}
\frac{(-\lambda^2/4)^j \Gamma((3+2j)/2)}{(1+2j)\Gamma^3(1+j)}
{}_1F_2( 1/2+j; 1,1+j;-\lambda^2/4).
\ee

\cite{Reynoldsarxiv2502}
\be
\int_0^{\pi/2} J_\nu(\lambda \cos\theta)^2 d\theta
= \surd \pi\sum_{j\ge 0}
\frac{(-)^j 2^{-1-2j-2\nu}\lambda^{2j+2\mu} \Gamma((1+2j+2\nu)/2)}
{\Gamma(1+j)\Gamma(1+\nu)\Gamma^2(1+j+\nu)}
{}_1F_2( 1/2+j+\nu; 1+\nu,1+j+\nu;-\lambda^2/4).
\ee

\cite{SchulzMC23}
\be
\int_0^{\pi/2} J_1(x\sin\theta)I_1(x\cos\theta) d\theta
=\int_0^x\frac{I_1(\sqrt{x^2-\xi^2}J_1(\xi)}{\sqrt{x^2-\xi^2}}d\xi
=\frac{I_1(x)-J_1(x)}{x}.
\ee

\cite{SchulzMC23}
\be
x\int_0^{\pi/2} I_1(x\sin\theta)I_1(x\cos\theta) d\theta
=\sqrt{2}I_1(\sqrt{2}x)-2I_1(x).
\ee

\cite{SchulzMC23}
\be
\int_0^{\pi/2} J_0(x\cos\theta)I_1(x\sin\theta) x\cos\theta  d\theta
=1-J_0(x).
\ee

\cite{SchulzMC23}
\be
\int_0^{\pi/2} I_0(x\cos\theta)J_1(x\sin\theta) x\cos\theta  d\theta
=I_0(x)-1.
\ee

\cite{SchulzMC23}
\be
\int_0^{\pi/2} I_0(x\cos\theta)I_1(x\sin\theta) x\cos\theta  d\theta
=I_0(\sqrt{2}x)-I_0(x).
\ee

\cite{SchulzMC23}
\be
\int_0^{\pi/2} J_0(x\cos\theta)J_1(x\sin\theta) x\cos\theta  d\theta
=J_0(x)-J_0(\sqrt{2}x).
\ee

\cite{BaileyPLMS40}
\be
\int_0^\infty J_\mu(ct\sin\phi)J_\nu(ct\sin\Phi) K_\rho(ct\cos\phi\cos\Phi)dt
=
...
\ee

\cite{BaileyPLMS40}
\begin{gather}
\int_0^\infty J_\mu(ct\sin\phi\sin\Phi)J_\nu(ct\cos\phi\sin\Phi) J_\rho(ct)dt
=
\frac{\Gamma(\frac{1+\mu+\nu+\rho}{2})\sin^\mu\phi\cos^\mu\Phi \cos^\nu\phi\sin^\nu\Phi}
{c^{\mu+\nu+1}\Gamma(\mu+1)\Gamma(\nu+1)\Gamma(\frac{1-\mu-\nu+\rho}{2})}
\\
\times
\,_2F_1(\frac{1+\mu+\nu-\rho}{2},\frac{1+\mu+\nu+\rho}{2};\mu+1;\sin^2\phi)
\,_2F_1(\frac{1+\mu+\nu-\rho}{2},\frac{1+\mu+\nu+\rho}{2};\nu+1;\sin^2\Phi)
\end{gather}
where $\phi$ and $\Phi$ are positive angles whose sum is acute.

\cite{Reynoldsarxiv2502}
\begin{multline}
\int_0^{\pi/2}
J_\nu(\lambda \sin\theta)
J_\mu(\lambda \sin\theta)
\cos^\beta \theta
\sin^m \theta d\theta
\\
=\sum_{j\ge 0}
\frac{(-)^j 2^{-1-2j-\nu-\mu}\lambda^{2j+\nu+\mu}\Gamma(\frac{1+\beta}{2})\Gamma(\frac{1+2j+m+\nu+\mu}{2})}
{\Gamma(1+j)\Gamma(1+j+\nu)\Gamma(1+\mu)\Gamma(\frac{2+2j+m+\nu+\beta+\mu}{2})}
\\ \times
{}_1F_2(j+(1+m+\nu+\mu)/2; 1+j+(m+\nu+\beta+\mu)/2,1+\mu;-\lambda^2/4)
\end{multline}
where $\Re\beta>0$, $\Re m>0$.

\subsubsection{Cylinder Functions combined with $x$ and $x^2$}

\cite{WingAJM72}
\be
\int_0^1 x J_\nu(\lambda_nx)J_\nu(\lambda_mx)dx
=
\frac{1}{2}\left\{(1-\nu^2/\lambda_n^2)J_\nu^2(\lambda_n)+J_{\nu}'^2(\lambda_n)\right\}\delta_{nm}
\ee
for $\{\lambda_n\}$ ($n=1,2,\cdots$) a sequence of succesive positive
roots of the equation
$xJ_\nu'(x)+HJ_\nu(x)=0$, where $H$ is a real number and $\nu\ge -1$.

\cite{MatharBA17}
\be
\int_0^1 x^{2+l+2n}j_l(2\pi\sigma x)dx
=
\frac{1}{2\pi\sigma}\sum_{k=0}^n \frac{(-n)_k}{(\pi\sigma)^k}j_{l+k+1}(2\pi\sigma);
\quad l,n=0,1,2,3,\ldots
\ee

\cite{PeeblesApJ185}
\be
\int_0^\infty k^2 j_l(kr) j_l(kr')dk = \frac{\pi \delta(r-r')}{2r^2}.
\ee

\cite{GrantJPA26,MehremJPA24}
\begin{multline}
\frac{4p_1p_2p_3}{\pi}
\int_0^\infty
x^2 j_1(p_1x)
j_2(p_2x)
j_3(p_3x) dx
= \Delta(p_1,p_2,p_3)(-1)^{(l_1+l_2+l_3)/2}
\\
 \times
\frac{1}{2}
\sum_{k_1=0}^{l_1}
\sum_{k_2=0}^{l_2}
\sum_{k_3=0}^{l_3}
\frac{(-1)^{k_1+k_2+k_3}}{(k_1+k_2+k_3)!}
\prod_{i=1}^3 \frac{(l_i+k_i)!}{k_i!(l_i-k_i)!}(2p_i)^{-k_i}
\\
\times
[
(-1)^{l_1+k_1}(p_2+p_3-p_1)^{k_1+k_2+k_3}
+(-1)^{l_2+k_2}(p_3+p_1-p_2)^{k_1+k_2+k_3}
\\
+(-1)^{l_3+k_3}(p_1+p_2-p_3)^{k_1+k_2+k_3}
-(p_1+p_2+p_3)^{k_1+k_2+k_3}
]
\end{multline}
supposed that $p_1$, $p_2$ and $p_3$ can be the sides of a plane triangle, that
is where $\Delta(.)=1$ if they form a non-denerate triangle, $\Delta(.)=1/2$
if they form a degenerate triangle, and $\Delta(.)=0$ otherwise.

\begin{equation}
\int _0^{\infty}
x
j_0(ax)
j_0(bx)
j_1(cx) dx
=
\frac{\pi}{8abc^2} (c^2-(a-b)^2)
\end{equation}
for $c>0$, $a>0$, $|c-a|<b<c+a$.

\cite{JacksonJMA3}
\be
\int_0^\infty
J_{n_1}(k_1\rho)
J_{n_2}(k_2\rho)
J_{n_3}(k_3\rho)
\rho d\rho
=\frac{\Delta}{6\pi A}
[
\cos(n_1\alpha_2-n_2\alpha_1)
+\cos(n_2\alpha_3-n_3\alpha_2)
+\cos(n_3\alpha_1-n_1\alpha_3)
]
\ee
if $n_1+n_2+n_3=0$, $\Delta$ as above, the area of the triangle of $k_1$, $k_2$
and $k_3$ given by $2A=k_1k_3\sin\theta_{13}$, and $\alpha$ three external angles
in that triangle.

\cite[2.8.3]{Apelblat2}
\begin{equation}
\int _0^\infty K_0(x)^3 dx 
=\frac{\pi}{2}[
K^2(e^{-i\phi_a/2})
+
K^2(e^{i\phi_a/2})
]
\approx 6.948822781079629789
\end{equation}
where $\phi_a=\pi/3$ and
$
e^{\pm i\phi_a/2} 
= \frac{\sqrt{3}\pm i}{2}
$
with Complete Elliptic Integrals
\begin{multline}
K(e^{-i\pi/6})
=K(\frac{\sqrt{3}-i}{2})
=\frac{\pi}{2}{}_2F_1(1/2,1/2;1;e^{-i\pi/3})
\\
\approx \frac{\pi}{2}(0.98274143349164123997896661 - 0.2633247734726891555758756i).
\end{multline}

\begin{equation}
\int_0^\infty x K_0(x)^3 dx 
\approx 0.5859768 .
\end{equation}

\begin{equation}
\int_0^\infty x K_0^2(x)K_1(x)dx \approx 2.31627 .
\end{equation}

\begin{equation}
\int_0^\infty x^2 K_0(x)^2 K_1(x) dx 
\approx 0.39065120.
\end{equation}

\cite{AlbrightJPA19}
\be
\int \frac{Ai^{n-1}(x)Bi^{n-1}(x)}{(Ai^n(x)+Bi^n(x))^2}dx=\frac{\pi}{n} \frac{Bi^n(x)}{Ai^n(x)+Bi^n(x)}.
\ee

\cite{AlbrightJPA19}
\be
\int \frac{dx}{Ai^2(x)+2Ai(x)Bi(x)+Bi^2(x)}=\frac{\pi Bi(x)}{Ai(x)+Bi(x)}.
\ee

\cite{AlbrightJPA19}
\be
\int \frac{dx}{Ai^2(x)}= \pi \frac{Bi(x)}{Ai(x)}.
\ee

\cite{AlbrightJPA19}
\be
\int \frac{dx}{Bi^2(x)}= -\pi \frac{Ai(x)}{Bi(x)}.
\ee

\cite{AlbrightJPA19}
\be
\int \frac{dx}{Ai(x) Bi(x)}= \pi \ln \frac{Bi(x)}{Ai(x)}.
\ee

\cite{AlbrightJPA19}
\be
\int \frac{Bi^n(x) dx}{Ai^{n+2}(x)}= \frac{\pi}{n+1} (\frac{Bi(x)}{Ai(x)})^{n+1}.
\ee

Triple products of mixed Modified Bessel Functions of index 0 or 1
are covered by Bailey's formula \cite[(3.2)]{BaileyJLMS11}
\begin{multline}
\int_0^\infty t^{\lambda-1} I_\mu(t) K_\nu(t) K_\rho(t)dt
\\
= 
\frac{\sqrt{\pi} 
\Gamma(\frac{\lambda+\mu+\nu-\rho}{2})
\Gamma(\frac{\lambda+\mu+\nu+\rho}{2})
\Gamma(\frac{\lambda+\mu-\nu-\rho}{2})
\Gamma(\frac{\lambda+\mu-\nu+\rho}{2})
}
{\Gamma(1+\mu)2^{\mu+2}
\Gamma(\frac{\lambda+\mu}{2}) \Gamma(\frac{1+\lambda+\mu}{2})
} 
{}_4F_3(\begin{array}{c}
\frac{\lambda+\mu+\nu-\rho}{2},
\frac{\lambda+\mu+\nu+\rho}{2},
\frac{\lambda+\mu-\nu-\rho}{2},
\frac{\lambda+\mu-\nu+\rho}{2}\\
\frac{\lambda+\mu}{2}, \frac{1+\lambda+\mu}{2},
\mu+1
\end{array}\mid \frac{1}{4}).
\end{multline}

Special cases are:
\begin{equation}
\int_0^\infty I_0(x)K_0(x)^2dx = \frac{\pi^2}{4}
{}_3F_2(\begin{array}{c}1/2,1/2,1/2\\
1,1
\end{array}\mid \frac14)\approx 2.554057858916278267;
\end{equation}
where (inserting $\gamma=1/2$, $\alpha=\beta=-1/4$ in (\ref{eq.rama3f2}))
$$
{}_3F_2(\begin{array}{c}1/2,1/2,1/2\\
1,1
\end{array}\mid \frac14)={}_2F_1^2(\begin{array}{c}
1/4,1/4\\ 1
\end{array}\mid \frac14
)
= \frac{2}{\surd 3}[P_{-1/4}(5/3)]^2\approx 1.01740879759^2.
$$
\begin{equation}
\int_0^\infty I_1(x)K_0(x)^2dx = \frac{1}{4}
{}_3F_2(\begin{array}{c}1,1,1\\
3/2,2
\end{array}\mid \frac14)\approx 0.2741556778080377394;
\end{equation}
where (inserting $\gamma=1$, $\alpha=\beta=-1/2$ in (\ref{eq.rama3f2}))
$$
{}_3F_2(\begin{array}{c}1,1,1\\
3/2,2
\end{array}\mid \frac14)=\frac{\pi^2}{9}
$$
\begin{equation}
\int_0^\infty I_1(x)K_0(x)K_1(x)dx = 
\frac{\pi^2}{16}{}_3F_2(\begin{array}{c}1/2,1/2,3/2\\
1,2
\end{array}\mid \frac14)\approx 0.6497774937995530258;
\end{equation}
\begin{equation}
\int_0^\infty x I_0(x)K_0(x)^2 dx = \frac{1}{2}\frac{2\pi}{3^{3/2}}\approx 0.6045997880780726168;
\end{equation}
\begin{equation}
\int_0^\infty x I_0(x)K_0(x)K_1(x) dx =
\frac{\pi^2}{8}{}_3F_2(\begin{array}{c}1/2,1/2,3/2\\
1,1
\end{array}\mid \frac14)\approx 1.374951395237058295;
\end{equation}
\begin{equation}
\int_0^\infty x I_1(x)K_0(x)^2 dx =
\frac{\pi^2}{64}{}_3F_2(\begin{array}{c}3/2,3/2,3/2\\
2,2
\end{array}\mid \frac14)\approx 0.1958449315578383236;
\end{equation}
\begin{equation}
\int_0^\infty x I_1(x)K_0(x)K_1(x) dx =
\frac{1}{4}{} \frac{2\pi}{3^{3/2}}\approx 0.3022998940390363084;
\end{equation}
\begin{equation}
\int_0^\infty x I_1(x)K_1(x)^2 dx =
\frac{3\pi^2}{64}{}_3F_2(\begin{array}{c}1/2,3/2,5/2\\
2,2
\end{array}\mid \frac14)\approx 0.529328969879666946;
\end{equation}
\begin{equation}
\int_0^\infty x^2 I_0(x)K_0(x)^2 dx =
\frac{\pi^2}{32}{}_3F_2(\begin{array}{c}3/2,3/2,3/2\\
1,2
\end{array}\mid \frac14)\approx 0.4984666397189814584;
\end{equation}
\begin{equation}
\int_0^\infty x^2 I_0(x)K_0(x)K_1(x) dx =
\frac{1}{2}{}_2F_1(\begin{array}{c}1,2\\
3/2
\end{array}\mid \frac14)\approx 0.7363998587187150779;
\end{equation}
\begin{equation}
\int_0^\infty x^2 I_0(x)K_1(x)^2 dx =
\frac{3\pi^2}{32}{}_3F_2(\begin{array}{c}1/2,3/2,5/2\\
1,2
\end{array}\mid \frac14)\approx  1.223640541156486728;
\end{equation}
\begin{equation}
\int_0^\infty x^2 I_1(x)K_0(x)^2 dx =
\frac{1}{6}{}_2F_1(\begin{array}{c}2,2\\
5/2
\end{array}\mid \frac14)\approx  0.263600141281284922 ;
\end{equation}
\begin{equation}
\int_0^\infty x^2 I_1(x)K_0(x)K_1(x) dx =
\frac{3\pi^2}{128}{}_3F_2(\begin{array}{c}3/2,3/2,5/2\\
2,2
\end{array}\mid \frac14)\approx  0.3471557856384098910 ;
\end{equation}
\begin{equation}
\int_0^\infty x^2 I_1(x)K_1(x)^2 dx =
\frac{1}{3}{}_2F_1(\begin{array}{c}1,3\\
5/2
\end{array}\mid \frac14)\approx  0.4727997174374301558 .
\end{equation}

\subsubsection{Cylinder Functions and Rational Functions}

\cite{FikiorisMathComp67}
\be
\int_0^\infty \frac{J_\nu(x)}{x^2+a^2}dx = \frac{i}{a}[S_{0,\nu}(ia-e^{-i\nu\pi/2}K_\nu(a)]
=
\frac{1}{a}[is_{0,\nu}(ia)+\frac{\pi}{2}\sec\frac{\nu\pi}{2}I_\nu(a)].
\ee

\cite{SoltMathComp47}
\begin{multline}
\int_0^\infty x^{1-2n}J_\nu(ax)J_\nu(bx)\frac{dx}{x^2+c^2}
=
(-1)^nc^{-2n}\Big\{
I_\nu(bc)K_\nu(ac)
\\
-\frac{1}{2}\left(\frac{b}{a}\right)^\nu
\frac{\pi}{\sin \pi\nu}\sum_{p=0}^{n-1}
\frac{(ac/2)^{2p}}{p!\Gamma(1-\nu+p)}\sum_{k=0}^{n-1-p}
\frac{(bc/2)^{2k}}{k!\Gamma(1+\nu+k)}
\Big\}
\end{multline}
for $0<b<a$, $\Re c>0$, $\Re\nu>n-1$, $n=1,2,\ldots$.
For $0<a<b$, the arguments $a$ and $b$ should be interchanged.

\subsubsection{Cylinder Functions and Powers}

\cite{KolbigMathComp64_449,KrupnikovMathComp41}
\be
\int_0^1 x^\mu J_\nu(ax)dx =
2^\mu\frac{\Gamma(\frac{1}{2}+\frac{1}{2}\mu+\frac{1}{2}\nu)}
{a^{\mu+1}\Gamma(\frac{\nu}{2}+\frac{1}{2}-\frac{1}{2}\mu)}
+a^{-\mu}
\left\{
(\mu+\nu-1)J_\nu(a) S_{\mu-1,\nu-1}(a)
-J_{\nu-1}(a)S_{\mu,\nu}(a)\right\}
\ee
\[
[a>0, \Re(\mu+\nu)>-1].
\]

\begin{multline}
\int_0^u z^\nu K_\nu(x)dx=
\frac{2^\nu-1}{1+2\nu}
[
(1+2\nu)\Gamma(\nu)u\,_1F_2(1/2;3/2,1-\nu;u^2/4)
\\
+2(u/2)^{1+2\nu}\Gamma(-\nu)
\,_1F_2(\nu+1/2;1+\nu,3/2+\nu;u^2/4)
]
.
\end{multline}

\cite{NollJOSA66}
\be
\int_0^\infty x^{-P}[1-J_0(bx)]dx = \frac{\pi b^{P-1}}{2^P\Gamma^2(\frac{P+1}{2}) \sin[\pi(P-1)/2]}.
\ee

\cite[p 22]{ErdelyiT2}
\begin{multline}
\int_0^1 x^\mu J_\nu(xy)\sqrt{xy} dx =
y^{-\mu-1}\big[
(\nu+\mu-\frac{1}{2})yJ_\nu(y)S_{\mu-1/2,\nu-1}(y)
-yJ_{\nu-1}(y)S_{\mu+1/2,\nu}(y)
\\
+2^{\mu+1/2}\frac{\Gamma(\frac{\mu}{2}+\frac{\nu}{2}+\frac{3}{4})}
  {\Gamma(\frac{\nu}{2}-\frac{\mu}{2}+\frac{1}{4})}
\big]
\end{multline}
for $\Re(\mu+\nu)>-3/2$.

\cite{GlasserMC33}
Let 
\be
I_\nu^{\alpha,\beta}(x)\equiv \int_0^\infty t^{2\alpha-1}(1+t^2)^{1-\alpha-\beta}J_\nu(x\sqrt{1+t^2})dt
\ee
then
\be
I_0^{1,1/2}(x)=\frac{1}{x}-J_0(x)+\frac{\pi}{2}[J_0(x){\mathbf H}_1(x)
-J_1(x){\mathbf H}_0(x)];
\ee
\be
I_1^{1,1/2}(x)=J_0(x)/x ;
\ee
\be
I_1^{1/2,1}(x)=x^{-1}\sin x;
\ee
\be
I_0^{1/2,1/2}(x)=x^{-1}\cos x;
\ee
\be
I_{1/2}^{\alpha,5/2-2\alpha}(1)=-\sqrt{\pi/8}
[J_{3/2-2\alpha}(1/2)Y_{1/2}(1/2)
+Y_{3/2-2\alpha}(1/2)Y_{1/2}(1/2)
],\quad 1/2<\alpha<1;
\ee
\be
I_0^{1/2,3/2}(1)=-si(1);
\ee
\be
I_\nu^{\nu-1/2,3/2-\nu/2}(x)=\frac{\Gamma(\nu-1/2)}{2^{1+\nu}\sqrt x}x^{-\nu}\sin x,\quad 1/2<\nu<5/2;
\ee
\be
I_\nu^{\nu+1/2,1/2-\nu/2}(x)=\frac{\Gamma(\nu+1/2)}{2^\nu\sqrt x}x^{-\nu-1}\cos x,\quad 1/2<\nu<5/2;
\ee
\be
I_\nu^{1/2-\nu,1+3\nu/2}(x)=-\frac{\surd \pi}{2^{\nu+1}}\Gamma(1/2-\nu)
x^{\nu}J_\nu(x/2)Y_\nu(x/2),\quad |\nu|<1/2;
\ee
\be
I_\nu^{1,1/4}(x)=x^{-1/2}[J_{\nu-1}(x)S_{1/2,\nu}(x)
+(1/2-\nu)J_\nu(x)S_{-1/2,\nu-1}(x) ];
\ee
\be
\frac{x}{2\nu}[I_{\nu-1}^{\alpha,\beta}(x)+I_{1+\nu}^{\alpha,\beta}]
=I_{\nu}^{\alpha,\beta+1/2}(x);
\ee
\be
x^{-\nu}[x^\nu I_{\nu}^{\alpha,\beta}(x)]'
=I_{\nu-1}^{\alpha,\beta-1/2}(x);
\ee
\be
x^{\nu}[x^{-\nu} I_{\nu}^{\alpha,\beta}(x)]'
=-I_{\nu+1}^{\alpha,\beta-1/2}(x).
\ee

\cite{FikiorisMathComp67}
\begin{multline}
\int_0^\infty \frac{x^{\rho-1}J_\nu(ax)}{(x^2+k^2)^{\mu+1}}dx
=
\frac{a^\nu k^{\rho+\nu-2\mu-2}\Gamma(\rho/2+\nu/2)\Gamma(\mu+1-\rho/2-\nu/2)}
{2^{\nu+1}\Gamma(\mu+1)\Gamma(\nu+1)}
\\
\times
\,_1F_2\left(\frac{\rho+\nu}{2};\frac{\rho+\nu}{2}-\mu,\nu+1;\frac{a^2k^2}{4}\right)
+
\frac{a^{2\mu+2-\rho} \Gamma(\nu/2+\rho/2-\mu-1)}
{2^{2\mu+3-\rho}\Gamma(\mu+2+\nu/2-\rho/2)}
\\
\times
\,_1F_2\left(\mu+1;\mu+2+\frac{\nu-\rho}{2},\mu+2-\frac{\nu+\rho}{2};\frac{a^2k^2}{4}\right)
,
\, a>0, -\Re \nu<\Re\rho < 2\Re \mu+\frac{7}{2}, \Re k >0
.
\end{multline}

\cite{WinkerJOSAA8}
\begin{multline*}
\int_0^\infty \frac{J_{n+1}(k)J_{n'+1}(k)}{k(k^2+k_0^2)^{1+\gamma2}}
dk
=
\\
k_0^{n+n'-\gamma}
\frac{\Gamma(\frac{\gamma-n-n'}{2})\Gamma(\frac{n+n'}{2}+1)}
{2^{n+n'+3}\Gamma(n+2)\Gamma(n'+2)\Gamma(1+\gamma/2)}
\\ \times
\,_3F_4\left(
\begin{array}{cc}
\frac{n+n'}{2}+1, \frac{n+n'}{2}+2, \frac{n+n'+3}{2} \\
n+2,n'+2,n+n'+3,1+\frac{n+n'-\gamma}{2}
\end{array}\mid
k_0^2
\right)
\\
+
\frac{\Gamma(\frac{n+n'-\gamma}{2})\Gamma(3+\gamma)}
{2^{3+\gamma}\Gamma(3+\frac{\gamma+n+n'}{2})
\Gamma(2+\frac{\gamma+n-n'}{2})\Gamma(2+\frac{\gamma+n'-n}{2})}
\,_3F_4\left(
\begin{array}{cc}
2+\frac{\gamma}{2},1+\frac{\gamma}{2},\frac{3+\gamma}{2} \\
2+\frac{\gamma+n-n'}{2},2+\frac{\gamma+n'-n}{2},
3+\frac{\gamma+n+n'}{2},1+\frac{\gamma-n-n'}{2}
\end{array} \mid
k_0^2
\right)
\end{multline*}

\cite{Reynoldsarxiv2502}
\begin{multline}
\int_0^b x^m(1+cx^p)^\beta J_\nu(x\lambda)J_\mu(x\rho) \log^k(\frac{1}{ax})dx
=
\frac{\lambda^\nu \rho^\mu}{2^{\nu+\mu}a^{1+m+\nu+\mu}}
\\ \times
\sum_{j\ge 0}
\sum_{l\ge 0}
\sum_{h\ge 0}
\frac{(-)^{j+l}c^h\lambda^{2j}\rho^{2l}\binom{\beta}{h}}
{2^{2j+2l}a^{2j+2l+hp}(1+2j+2l+m+hp+\nu+\mu)^{k+1}}
\frac{\Gamma(1+k,(1+2j+2l+m+hp+\nu+\mu)\log\frac{1}{ab})}{j!l!\Gamma(1+j+\nu)\Gamma(1+l+\mu)}
\end{multline}
where $\Re \beta >0$, $\Re b>0$.

\cite[p50]{MO2Afl}
\be
\int_0^\infty J_{\nu+n}(at)J_{\nu-n-1}(bt)dt
=
\left\{
\begin{array}{ll}
\frac{b^{\nu-n-1}\Gamma(\nu)}{a^{\nu-n}n!\Gamma(\nu-n)} \,_2F_1(\nu,-n;\nu-n;\frac{b^2}{a^2}), & 0<b<a, \\
(-1)^n/(2a), & 0<b=a, \\
0, & 0<a<b.
\end{array}
\right.
\ee
where $n=0,1,2,\ldots$, $\Re \nu>0$.

\cite[p50]{MO2Afl}
\be
\int_0^\infty J_{\nu}(at)J_{\nu+1}(bt)dt
=
\left\{
\begin{array}{ll}
a^\nu b^{-\nu-1}, & 0<b<a \\
1/(2a), & 0<b=a \\
0, & 0<a<b
\end{array}
\right.
\ee
where $a$, $b$ real positive, $\Re \nu>-1$.

\cite[p50]{MO2Afl}
\be
\int_0^\infty J_{\mu}(at)J_{\nu}(at)dt
=
\frac{2}{\pi}\frac{\sin(\frac{\nu-\mu}{2}\pi)}{\nu^2-\mu^2},
\ee
$\Re(\nu+\mu)>0$, $a>0$.

\cite{NollJOSA66}
\be
\int_0^\infty x^{-P}\left( 1-\frac{4J_1^2(x)}{x^2}\right)dx =
\frac{\pi \Gamma(P+2)}{2^P
\Gamma^2(\frac{P+3}{2})
\Gamma(\frac{P+5}{2})
\Gamma(\frac{P+1}{2})
\sin[\pi(P-1)/2]}.
\ee

\cite{FabrikantZAMM83}
\be
\int_0^\infty t^{\rho-\mu-\nu-3}J_\mu(at)J_\nu(bt)J_\rho(ct)dt
=
\frac{2^{\rho-\mu-\nu-3}a^\mu b^\nu \Gamma(\rho-1)}{c^{\rho-2}\Gamma(\mu+1)\Gamma(\nu+1)}\left(1-\frac{\rho-1}{\mu-1}\,\frac{a^2}{c^2}
-\frac{\rho-1}{\nu-1}\,\frac{b^2}{c^2}
\right)
.
\ee

\cite{TylerJOSA7}
Let
\be
G_{lmn}\equiv \int_0^\infty dx x^\alpha J_l(x)J_m(x)J_n(\beta x)
\ee
with $\Re \alpha<3/2$ and $\Re(\alpha+l+m+n+1)>0$, then for $\beta/2<1$
\begin{multline}
G_{lmn}
=
\frac{\Gamma(\frac{\alpha+n)}{2})(\beta/2)^{-\alpha}}
{
2\Gamma(\frac{-l+m+1}{2})\Gamma(\frac{l-m+1}{2})\Gamma(\frac{-\alpha+n+2}{2})
}
\,_4F_3\left(
\begin{array}{c}
\frac{l+m+1}{2},
\frac{l-m+1}{2},
\frac{-l-m+1}{2},
\frac{-l+m+1}{2}\\
\frac{-\alpha+n+2}{2},
\frac{1}{2},
\frac{-\alpha-n+2}{2}
\end{array}
\mid\beta^2/4\right)
\\
-\frac{\Gamma(\frac{l+m+2}{2})\Gamma(\frac{\alpha+n-1}{2}) (\beta/2)^{-\alpha+1}}
{
\Gamma(\frac{-l+m}{2})\Gamma(\frac{l+m}{2})\Gamma(\frac{l-m}{2})\Gamma(\frac{-\alpha+n+3}{2})
}
\,_4F_3\left(
\begin{array}{c}
\frac{l+m+2}{2},
\frac{l-m+2}{2},
\frac{-l-m+2}{2},
\frac{-l+m+2}{2}\\
\frac{-\alpha+n+3}{2},
\frac{3}{2},
\frac{-\alpha-n+3}{2}
\end{array}
\mid\beta^2/4\right)
\\
+\frac{2^{\alpha+n}\Gamma(-\alpha-n)\Gamma(\frac{\alpha+l+m+n+1}{2}) (\beta/2)^n}
{
\Gamma(\frac{-\alpha-l+m-n+1}{2})\Gamma(\frac{-\alpha+l+m-n+1}{2})
\Gamma(\frac{-\alpha+l-m-n+1}{2})\Gamma(n+1)
}
\\ \times
\,_4F_3\left(
\begin{array}{c}
\frac{\alpha+l+m+n+1}{2},
\frac{\alpha+l-m+n+1}{2},
\frac{\alpha-l-m+n+1}{2},
\frac{\alpha-l+m+n+1}{2}\\
n+1,
\frac{n+\alpha+2}{2},
\frac{\alpha+n+1}{2}
\end{array}
\mid\beta^2/4\right)
,
\label{eq.Glmn}
\end{multline}
and for $\beta/2>1$
\begin{multline}
G_{lmn}(\alpha,\beta)=
\frac{2^{-l-m-1}\Gamma(\frac{\alpha+l+m+n+1}{2}) (\beta/2)^{-\alpha-l-m-1}}
{
\Gamma(\frac{-\alpha-l-m+n+1}{2}) \Gamma(l+1)\Gamma(m+1)
}
\,_4F_3\left(
\begin{array}{c}
\frac{\alpha+l+m+n+1}{2},
\frac{l+m+1}{2},
\frac{l+m+2}{2},
\frac{\alpha+l+m-n+1}{2}\\
m+1,
l+1,
m+l+1
\end{array}
\mid 4/\beta^2\right)
.
\end{multline}
In \eqref{eq.Glmn}, any of the three terms on the right hand side for which
any of the $\Gamma$-functions in the denominator has a pole is to be
interpreted as its limiting value of zero (unless $\alpha$ is integer which
may create a counter-balance in the numerator).

The case $a=c=1$ in \cite[6.578.1]{GR} leads to the same expression
\begin{multline}
\int_0^\infty x^{\rho-1}J_\lambda(x)J_\nu(x)J_\mu(bx)dx
=
\frac{2^{\rho-1}b^\mu\Gamma(\frac{\lambda+\nu+\mu+\rho}{2})}{\Gamma(\nu+1)\Gamma(\mu+1)
\Gamma(1-\frac{\lambda+\mu-\nu+\rho}{2})}
\\
\times F_4\left(\frac{\lambda+\mu-\nu+\rho}{2},\frac{\lambda+\mu+\nu+\rho}{2};\lambda+1,\mu+1;1,b^2\right)
\end{multline}
if the Apell series \cite[9.180.4]{GR} is inserted and one of the double sums
summed by the formula \cite[15.1.20]{AS} for the Gaussian Hypergeometric Function
of unit argument.

\cite{TezerJCAM28}
\be
\int_0^\infty t\left(\frac{t}{\sqrt{u^2+t^2}}-1\right)J_0(\gamma t)
dt
=\frac{u^2}{2}[I_1(u\gamma /2)K_1(u\gamma/2)-I_0(u\gamma/2)K_0(u\gamma/2)].
\ee

\subsubsection{Cylinder Functions and Exponentials}
\subsubsection{Cylinder Functions, Exponentials and Powers}

With 
\[
l_1(a,b,c)\equiv \frac{1}{2}\left[\sqrt{(a+b)^2+c^2}-\sqrt{(a-b)^2+c^2}\right]
\]
\[
l_2(a,b,c)\equiv \frac{1}{2}\left[\sqrt{(a+b)^2+c^2}+\sqrt{(a-b)^2+c^2}\right]
\]
at $a>0$, $b>0$, $c>0$ \cite{FabrikantZAMM83}
\be
\int_0^\infty e^{-cx}J_1(ax)J_{1/2}(bx)\frac{dx}{x^{3/2}}
=
\frac{\sqrt{2}}{\sqrt{\pi b}a}
\left[\frac{l_1}{2}\sqrt{a^2-l_1^2}+\frac{a^2}{2}\sin^{-1}(\frac{l_1}{a})
+c(\sqrt{b^2-l_1^2}-b)\right]
.
\ee
\be
\int_0^\infty e^{-cx}J_1(ax)J_{1/2}(bx)\frac{dx}{\sqrt{x}}
=
\frac{\sqrt{2}}{\sqrt{\pi b}a}
(b-\sqrt{b^2-l_1^2})
.
\ee
\be
\int_0^\infty e^{-cx}J_1(ax)J_{1/2}(bx)\sqrt{x}dx
=
\frac{\sqrt{2}}{\sqrt{\pi b}a}
\frac{l_1\sqrt{a^2-l_1^2}}{l_2^2-l_1^2}
.
\ee
\be
\int_0^\infty e^{-cx}J_1(ax)J_{3/2}(bx)\sqrt{x} dx
=
\frac{2l_1^2\sqrt{b^2-l_1^2} }
{\sqrt{2\pi} b^{3/2}a(l_2^2-l_1^2)}
.
\ee
\be
\int_0^\infty e^{-cx}J_1(ax)J_{3/2}(bx)\frac{dx}{\sqrt{x}}
=
\frac{1}{\sqrt{2\pi}b^{3/2}a}
\left(-l_1\sqrt{a^2-l_1^2}+a^2\sin^{-1}(l_1/a)\right)
.
\ee
\be
\int_0^\infty e^{-cx}J_1(ax)J_{5/2}(bx)\frac{dx}{\sqrt{x}}
=
\frac{2^{-1/2}c}{\sqrt{\pi}b^{5/2}a}
\left(l_1\sqrt{a^2-l_1^2}+\frac{2a^2l_1}{\sqrt{a^2-l_1^2}}-3a^2\sin^{-1}(l_1/a)\right)
.
\ee
and similar combinations of even and odd-indexed $J(ax)$ and $J(bx)$.

\cite{AlhaidariAML20}
\be
\int_0^\infty x^\nu e^{-x/2}J_\nu(\mu x)L_n^{2\nu}(x)dx
=
2^\nu \Gamma(\nu+\frac{1}{2})\frac{1}{\sqrt{\pi\mu}}
(\sin\theta)^{\nu+\frac{1}{2}} C_n^{\nu+\frac{1}{2}}(\cos\theta),
\ee
$\mu\ge 0$, $\nu>-\frac{1}{2}$,
$\cos\theta\equiv\frac{\mu^2-1/4}{\mu^2+1/4}$, $C_n^\lambda(x)$ ultraspherical polynomial.

\cite{OberhettBT}
\be
\int_0^\infty \frac{x^{\nu+1}}{a^2+x^2} J_{\nu}(xy)dx
=\frac{a^\nu}{\sqrt y}K_\nu(ay),\quad 1< \Re \nu < 3/2.
\ee
\cite{OberhettBT}
\be
\int_0^\infty \frac{x^{\nu+1}}{(a^2+x^2)^\mu} J_{\nu}(xy)dx
=\frac{2^{1-\mu}a^{\nu-\mu+1}y^{\mu-1}}{\Gamma(\mu)}K_{\nu-\mu+1}(ay),\quad
\Re \nu>-1, \Re(2\mu-\nu)>1/2.
\ee

\cite{OberhettBT}
\begin{multline}
\int_0^\infty \frac{x^{1-\nu}}{(a^2+x^2)^\mu} J_{\nu}(xy)dx
=
a^{-\mu-\nu+1}y^{\mu-1}
\big[
2^{-\mu}\frac{\Gamma(1-mu)}{1-\nu}I_{\nu+\mu-1}(ay)
\\
-2\frac{1-\mu}{\Gamma(\nu)}
e^{-i\frac{\pi}{2}(\nu-\mu+1)}s_{-\mu+\nu,-\mu-\nu+1}(iay)
\big]
,\quad
\Re(\nu+2\mu)>1/2.
\end{multline}

\cite{BorosSci12}
\be
\int_0^\infty x^{\beta-1}e^{-x}I_n(x)dx
=
\sum_{j=0}^\infty \frac{\Gamma(\beta+2j+n)}{2^{2j+n}j!(j+n)!}
.
\ee
\cite{BorosSci12}
\be
\int_0^\infty x^{2n+1}e^{-x^2/4a}I_0(x)dx = 2^{2n+1} a^{n+1} n! e^a L_n(-a)
=
2^{2n+1}a^{n+1}n!\, _1F_1(n+1;1;a)
.
\ee
\be
\int_0^\infty x^3e^{-x^2/4a}I_0(x)dx = 8a^2(1+a)e^a .
\ee

\cite{AdamchikJCAM64}
\be
\int_0^\infty x e^{-x^2/z}J_2(x)Y_2(x)dx = -\frac{2}{\pi}+\frac{4}{\pi z}
-
\frac{zK_2(z/2)}{2\pi \exp(z/2)}.
\ee

\cite{AdamchikJCAM64}
\be
\int_0^\infty x e^{-x^2/z}I_3(x)K_3(x)dx =
-\frac{32+16z+3z^2}{2z^2}+\frac{\exp(z/2)zK_3(z/2)}{4}.
\ee
\cite{AdamchikJCAM64}
\be
\int_0^\infty x^3 e^{-x^2/z}J_2(x)Y_2(x)dx =
-\frac{4}{\pi}+\frac{z^2(2+z)K_0(z/2)}{4\pi\exp(z/2)}+
\frac{z(8+4z+z^2)K_1(z/2)}{4\pi\exp(z/2)}.
\ee

\cite{AdamchikJCAM64}
\begin{multline}
\int_0^\infty x^5 e^{-x^2/z}I_3(x)K_3(x)dx =
-32+\frac{1}{8}e^{z/2}z^2(32-16z+5z^2-z^3)K_0(z/2)
\\
+\frac{1}{8}
e^{z/2}z(128-64z+24z^2-6z^3+z^4)K_1(z/2).
\end{multline}

\cite{EasonPTRSL247}
Let
\be
I(\mu,\nu,\lambda)\equiv \int_0^\infty J_\mu(at)
J_\nu(bt)e^{-ct} t^\lambda dt,
\ee
then
\be
I(n,n;0) = \frac{(-)^nk}{\pi\sqrt{ab}}\int_0^{\pi/2}
\frac{\cos(2n\psi)d\psi}{\sqrt{1-k^2\sin^2\psi}},
\ee
\be
I(n,n;1) = \frac{(-)^nck^3}{4\pi(ab)^{3/2}}\int_0^{\pi/2}
\frac{\cos(2n\psi)d\psi}{(1-k^2\sin^2\psi)^{3/2}},
\ee
in particular
\be
I(0,0;0)=\frac{k}{2\sqrt{ab}}F_0(k),
\ee
\be
I(1,1;0)=\frac{1}{k\sqrt{ab}}[(1-k^2/2)F_0(k)-E_0(k)],
\ee
\be
I(0,0;1)=\frac{ck^3E_0(k)}{8k'^2(ab)^{3/2}},
\ee
and more results on $I(n+1,n;\pm 1)$ and $I(n+1,n;0)$
in terms of Elliptic Integrals.
Associated recurrences:
\be
a[I(\mu+1,\nu;\lambda)+I(\mu-1,\nu;\lambda)]
=2\mu I(\mu,\nu;\lambda-1) ;
\ee
\be
b[I(\mu,\nu+1;\lambda)+I(\mu,\nu-1;\lambda)]
=2\nu I(\mu,\nu;\lambda-1) ;
\ee
\be
aI(\mu+1,\nu;\lambda)-bI(\mu,\nu-1;\lambda)
=C_{\mu,\nu}+(\mu-\nu+\lambda) I(\mu,\nu;\lambda-1)
-c I(\mu,\nu;\lambda),
\ee
with
\be
C_{\mu,\nu}\equiv\begin{cases}
\frac{a^\mu b^\nu}{2^{\mu+\nu}\Gamma(\mu+1)\Gamma(\nu+1)},&\mathrm{if}\, \lambda+\mu+\nu=0;\\
0& \mathrm{if}\, \lambda+\mu+\nu>0.
\end{cases}
\ee

\cite{RagabGMJ2}
\begin{multline}
\int_0^\infty e^{-\lambda}\lambda^{k-2} K_m(\lambda)K_n(z/\lambda)d\lambda
\\
=
\sum_{n,-n} \frac{\Gamma(\frac12)\Gamma(k+m+n)\Gamma(k-m+n)}{\Gamma(k+n+\frac12) 2^{k+1}}\Gamma(n)z^{-n}\\
\times F(\begin{array}{c} \frac34-\frac{k}{2}-\frac{n}{2}, \frac14-\frac{k}{2}-\frac{n}{2}\\
1-n,1-\frac{k}{2}-\frac{m}{2}-\frac{n}{2}, 1-\frac{k}{2}+\frac{m}{2}-\frac{n}{2},
\frac12-\frac{k}{2}-\frac{m}{2}-\frac{n}{2},\frac12-\frac{k}{2}+\frac{m}{2}-\frac{n}{2} \\
\end{array};\frac{z^2}{4})
\\
+\sum_{m,-m} \Gamma(-\frac{k}{2}-\frac{m}{2}+\frac{n}{2})
\Gamma(-\frac{k}{2}-\frac{m}{2}-\frac{n}{2})
\Gamma(-m)2^{-m-3}(z/2)^{m+k}
\\
\times F(\begin{array}{c} \frac34+\frac{m}{2}, 
\frac14+\frac{m}{2}\\
1+\frac{k}{2}+\frac{m}{2}-\frac{n}{2},
1+\frac{k}{2}+\frac{m}{2}+\frac{n}{2},
\frac12,\frac12+m,1+m
\end{array};\frac{z^2}{4})\\
-\sum_{m,-m} \Gamma(-\frac{k}{2}-\frac{m}{2}-\frac{n}{2}-\frac12)
\Gamma(-\frac{k}{2}-\frac{m}{2}+\frac{n}{2}-\frac12)
\Gamma(-m)2^{-m-3}(z/2)^{m+k+1}
\\
\times F(\begin{array}{c} \frac54+\frac{m}{2}, 
\frac34+\frac{m}{2}\\
\frac32+\frac{k}{2}+\frac{m}{2}-\frac{n}{2},
\frac32+\frac{k}{2}+\frac{m}{2}+\frac{n}{2},
\frac32,1+m,\frac32+m
\end{array};\frac{z^2}{4})
,
\end{multline}
where $\Re z>0$.

\cite{RagabGMJ2,VyasMN40}
\begin{multline}
\int_0^\infty \lambda^{k-1}K_\nu(\lambda)K_{\mu}(x\lambda^{-n})d\lambda
=2^{k-n-3}\pi^{-n}n^{k-1}\\
\times \sum_{i,-i}\frac{1}{i}
E\left\{
\Delta(n;\frac{k}{2}+\frac{\nu}{2}),
\Delta(n;\frac{k}{2}-\frac{\nu}{2}),
\frac{\mu}{2},
-\frac{\mu}{2},
1::\frac14 e^{i\pi}(2n)^{-2n}x^2
\right\},
\end{multline}
if $\Re (k\pm \nu+\frac{n}{2})>0$, $n$ a positive integer,
$x$ real and positive,
where $\Delta(n;\alpha)$ represents the set of
parameters $\alpha/n,(\alpha+1)/n,\ldots (\alpha+n-1)/n$, and where
$E$ is MacRobert's E-function.

\cite{RagabGMJ2,VyasMN40}
\begin{multline}
\int_0^\infty \lambda^{k-1}K_\nu(\lambda)K_{\mu}(x\lambda^n)d\lambda
=-2^{k-n-2}\pi^{2-n}n^{k-1}
 \sum_{\mu,-\mu}\mathrm{cosec} \mu\pi\left\{\frac14(2n)^{2n}x^2\right\}^{\mu/2}
\\
\times
E\left\{
\Delta(n;\frac{k+n\mu+\nu}{2})
\Delta(n;\frac{k+n\mu-\nu}{2}):
1+\mu
:\frac14 e^{\pm i\pi}(2n)^{-2n}x^{-2}
\right\},
\end{multline}
if $\Re (k\pm \nu+n\mu)>0$, $n$ a positive integer,
$x$ real and positive.

\cite{RagabRenc14}
\begin{multline}
\int_0^\infty \lambda^{k-1}K_\nu(\lambda)K_\mu(x\lambda)
d\lambda =
2^{k-3}\frac{1}{\Gamma(k)}
x^{-\nu-k}
\Gamma(\frac{k+\mu+\nu}{2})
\Gamma(\frac{k+\mu-\nu}{2})
\Gamma(\frac{k+\nu-\mu}{2})
\Gamma(\frac{k-\mu-\nu}{2})
\\ \times
{}_2F_1((k+\mu+\nu)/2,(k+\nu-\mu)/2;k;1-1/x^2).
\end{multline}
where $\Re(k\pm \mu\pm \nu)>0$ and $x$ real and positive.

\cite{RagabRenc14}
\begin{multline}
\prod_{r=1}^{m-1}\int_0^\infty \lambda_r^{k_r-1}K_{\mu_r}(\lambda_r)
d\lambda_r
K_\mu(\frac{x}{\lambda_1\lambda_2\cdots \lambda_{m-1}})
\\
=2^{\sum_{r=1}^{m-1}k_r+1-2m}
\frac{1}{2\pi}\sum_{i,-i}\frac{1}{i}E[
\frac{\nu_1+k_1}{2},
\frac{k_1-\nu_1}{2},
\cdots
\frac{\nu_{m-1}+k_{m-1}}{2},
\frac{k_{m-1}-\nu_{m-1}}{2},
\frac{\mu}{2},-\frac{\mu}{2},1::\frac{e^{i\pi}}{2^{2m}}x^2
]
,
\end{multline}
$m=2,3,4\ldots$, where $E$ is MacRobert's E-function,
where the notation ${i,-i}$ in the sum limits means
to evaluate at $i$ and add the evaluation at $-i$.

\cite{RagabRenc14}
\begin{multline}
\prod_{r=1}^{m-1}\int_0^\infty \lambda_r^{k_r-1}K_{\mu_r}(\lambda_r)
d\lambda_r
K_\mu(x\lambda_1\lambda_2\cdots \lambda_{m-1})
\\
=2^{\sum_{r=1}^k k_r+1-2m}
\pi\sum_{\mu,-\mu}\frac{1}{\sin\pi\mu}
(2^{2m-4}x^2)^{-\mu/2}
E[
\frac{
\frac{k_1+\nu_1-\mu}{2},
\frac{k_1-\nu_1-\mu}{2},
\cdots
\frac{k_{m-1}+\nu_{m-1}-\mu}{2},
\frac{k_{m-1}-\nu_{m-1}-\mu}{2}
:\frac{e^{\pm ix}}{2^{2m-4}x^2}
}{1-\mu}
]
,
\end{multline}
$m=2,3,4\ldots$.

\cite{RagabRenc14}
\be
\prod_{r=1}^{m-1}\int_0^\infty \lambda_r^{2r/m-1}K_{\mu}(\lambda_r)
d\lambda_r
K_\mu(\frac{x}{\lambda_1\lambda_2\cdots \lambda_{m-1}})
=\pi^{\mu-1}K_{m\mu}(mx^{1/m}).
\ee

\subsubsection{Cylinder and Trigonometric Functions and Powers}

\cite{FabrikantZAMM83}
\be
\int_0^\infty \sin(cx) x^{\nu-\mu-4}J_\mu(ax)J_\nu(bx)dx
=
\frac{\Gamma(\nu)a^\mu b^{-\nu} c}{2^{\mu-\nu+3}\Gamma(\mu+1)}\left(
\frac{b^2}{\nu-1}-\frac{a^2}{\mu+1}-\frac{2c^2}{3}
\right)
.
\ee
\cite{FabrikantZAMM83}
\be
\int_0^\infty \cos(cx) x^{\nu-\mu-3}J_\mu(ax)J_\nu(bx)dx
=
\frac{\Gamma(\nu)a^\mu b^{-\nu} }{2^{\mu-\nu+3}\Gamma(\mu+1)}\left(
\frac{b^2}{\nu-1}-\frac{a^2}{\mu+1}-2c^2
\right)
.
\ee

\cite{OkuiNBS79B} representations as complete elliptic Intgrals:
\be
\int_0^\infty J_k(ax^2)\left\{\begin{array}{c}\cos(bx^2)\\ \sin(bx^2)\end{array}\right\} dx=\ldots 
\ee
for $k=0,\ldots,3$.
\be
\int_0^\infty J_k(ax^2)\left\{\begin{array}{c}\cos(bx^2)\\ \sin(bx^2)\end{array}\right\} x^{-2}dx=\ldots 
\ee
for $k=1,\ldots,3$.
\be
\int_0^\infty J_k(ax^2)\left\{\begin{array}{c}\cos(bx^2)\\ \sin(bx^2)\end{array}\right\} x^{-4}dx=\ldots 
\ee
for $k=2,3$.
\be
\int_0^\infty J_k^2(ax^2)\cos(2bx^2)dx=\ldots 
\ee
for $k=0,\ldots 2$.
\be
\int_0^\infty J_k^2(ax)\sin(2bx)x^{-2}dx=\ldots 
\ee
for $k=1,2$.
\be
\int_0^\infty J_3(ax)J_1(ax)\sin(2bx)dx=\ldots 
\ee
\be
\int_0^\infty J_3(ax)J_2(ax)\sin(2bx)x^{-1}dx=\ldots 
\ee
\be
\int_0^\infty J_k(ax)J_k(bx)\cos(cx)dx=\ldots 
\ee
for $k=0,\ldots 3$.
\be
\int_0^\infty Y_k(ax)J_k(bx)\cos(cx)dx=\ldots 
\ee
for $k=0,\ldots 3$.
\be
\int_0^\infty K_k(ax)\cos(bx^2)x^l dx=\ldots 
\ee
for $k=0,\ldots 2$ and even $2k\le l\le 6$.
\be
\int_0^\infty K_k(ax)\sin(bx^2)x^l dx=\ldots 
\ee
for $k=0,\ldots 2$ and even $2k-2\le l\le 6$.
Additionally integrals of products of $K_k(ax)I_l(bx) sin(cx)$
and similar products where Hyperbolic Functions and exponentials
appear au lieu of the sines and cosines.

\cite{RottbrandITSF9}
\begin{multline}
\int_0^{\infty} J_{2n}(a\sqrt{t^2+2bt})e^{-pt}dt
= (-1)^n\bigg\{
\frac{e^{b(p-\sqrt{p^2+a^2})}}{\sqrt{p^2+a^2}}
+\frac{2e^{bp}}{a}
\\
\times
\sum_{k=1}^n(-)^k(2k-1)I_{k-1/2}[b(\sqrt{p^2+a^2}-a)/2)]
K_{k-1/2}[b(\sqrt{p^2+a^2}+a)/2]
\bigg\}
.
\end{multline}

\cite{RottbrandITSF9}
\begin{multline}
\int_0^{\infty} J_{2n}(a\sqrt{t^2+2bt})e^{-pt}dt
= (-1)^n\bigg\{
\frac{e^{b(p-\sqrt{p^2+a^2})}}{\sqrt{p^2+a^2}}
+2n\sum_{\lambda=1}^n\frac{(n-1+\lambda)!}{\lambda!(n-\lambda)!}
\frac{1}{b^{2\lambda-1}a^{2\lambda}}
\\
\times
\sum_{\mu=0}^{\lambda-1}\frac{(2\lambda-2-\mu)!}{\mu !(\lambda-1-\mu)!}
(2b\sqrt{p^2+a^2})^\mu
\left(e^{b(p-\sqrt{p^2+a^2})}-\sum_{\nu=0}^{2\lambda-2-\mu}
\frac{[b(p-\sqrt{p^2+a^2})]^\nu}{\nu!}
\right)
\bigg\}
.
\end{multline}

\cite{BorosSci12}
\be
\int_0^\infty x^{\beta-1}e^{-x} I_0(2\sqrt{ax})dx
=
\sum_{k=0}^\infty \frac{a^k}{k!^2}\Gamma(\beta+k)
.
\ee

\section{Definite Integrals of Special Functions II.}
\subsection{Associated Legendre Functions and Powers}

\cite{RashidJPA19}
\begin{multline*}
\int_{-1}^1 P_{l+\alpha}^{-\alpha}(x)P_{k+\beta}^{-\beta}(x)
(1-x^2)^{-p-1}dx
=
2^{-(\alpha+\beta)}
\\
\times
\frac{
\Gamma[\frac{1}{2}(k+1)]\Gamma(-\frac{1}{2}k-\beta)
\Gamma[\frac{1}{2}(\alpha+\beta)-p]\Gamma[\frac{1}{2}(\alpha-\beta-p)]
\Gamma[\frac{1}{2}(l-k)]\Gamma[\frac{1}{2}(l+k+1)+\beta]
}{
\Gamma(\beta+1)\Gamma(-\beta)\Gamma(\alpha+1)
\Gamma(-\frac{1}{2}k)\Gamma[\frac{1}{2}(k+1)+\beta]
\Gamma[\frac{1}{2}(l-k+\alpha-\beta)-p]\Gamma[\frac{1}{2}(l+k+\alpha+\beta+1)-p]
}
\\
\times
\,_4F_3\left(
\begin{array}{cccc}
\frac{1}{2}(\alpha-\beta)+p+1, & 
\frac{1}{2}(\alpha-\beta)-p, & 
-\frac{1}{2}(l-1), & 
-\frac{1}{2}l \\ 
\frac{1}{2}(k-l)+1, & 
-\frac{1}{2}(l+k-1)-\beta, & 
\alpha+1 & 
\end{array};1
\right)
\end{multline*}
for $\Re[\frac{1}{2}(\alpha+\beta)-p]>0$, for $k,l$ both even or both odd.

\cite{RashidJPA19}
\begin{multline}
\int_{-1}^1 P_l^m(x)P_k^n(x)
(1-x^2)^{-p-1}dx
=
(-1)^{(l-k+m-n)/2}
2^{-(m+n)}
\\
\times
\frac{
(l+m)!\Gamma[\frac{1}{2}(m+n)-p]
\Gamma[\frac{1}{2}(m-n)-p]\Gamma[\frac{1}{2}(l+k+1-m+n)]
}{
(l-m)! m![\frac{1}{2}(k-l+m-n)]!
\Gamma[\frac{1}{2}(l-k)-p]
\Gamma[\frac{1}{2}(l+k+1)-p]
}
\\
\times
\,_4F_3\left(
\begin{array}{cccc}
\frac{1}{2}(m-n)+p+1, & 
\frac{1}{2}(m-n)-p, & 
-\frac{1}{2}(l-m-1), & 
-\frac{1}{2}(l-m) \\ 
\frac{1}{2}(k-l+m-n)+1, & 
-\frac{1}{2}(l+k-1-m+n), & 
m+1 & 
\end{array};1
\right)
\end{multline}
for $\Re[\frac{1}{2}(m+n)-p]>0$, for $k-n\ge l-m$.

\cite{RashidJPA19}
\be
\int_{-1}^1 P_l^m(x)P_k^m(x)
(1-x^2)^{-1}dx
=
\frac{1}{m}\frac{(l+m)!}{(l-m)!},\quad k\ge l\ge m>0.
\ee

\cite{RashidJPA19}
\begin{multline*}
\int_{-1}^1 P_l^m(x)P_k^m(x)
(1-x^2)^{-2}dx
=
\frac{(l+m)!}{2(l-m)!(m-1)m(m+1)}
\\ \times
[l(l+1)+(m-1)(m+1)+\frac{1}{2}(k-l)(m+1)(k+l+1)],
\quad k\ge l \ge m >1.
\end{multline*}

\cite{RashidJPA19}
\be
\int_{-1}^1 P_l^m(x)P_k^m(x)
=
\frac{2}{2l+1}\,\frac{(l+m)!}{(l-m)!}\delta_{kl}
\ee

\cite{RashidJPA19}
\begin{multline*}
\int_{-1}^1 P_l^m(x)P_k^m(x)
(1-x^2)dx
=
\\
(-1)^{(k-l)/2}
\frac{(l+m)![l(l+1)+(m-1)(m+1)+\frac{1}{2}(k-l)(m+1)(k+l+1)]}
{(l-m)![\frac{1}{2}(k-l)+1]!\Gamma[\frac{1}{2}(l-k)+2]
(k+l-1)(k+l+1)(k+l+3)
}
\end{multline*}
for $k\ge l \ge m \ge 0$, zero for $k-l\ge 4$.

\subsection{Associated Legendre functions, powers, and trigonometric functions}

\cite{BaileyPCPS27_3,BaileyPCPS27_2}
\be
\int_0^z P_n[\cos(z-t)]P_n^{-m}(\cos t)\frac{dt}{\sin t}=\frac{P_n^{-m}(\cos z}{m},\quad \Re m>0.
\ee
\be
\int_0^z \sin^m t P_{n-m-1}[\cos(z-t)]P_n^{-m}(\cos t)dt=\frac{\Gamma(m+1/2)}{2^{1/2}\Gamma(m+1)} \sin^{m+1/2}z P_{n-1/2}^{-m-1/2}(\cos z),\quad \Re m>-1/2.
\ee
\be
\int_0^z \sin^{m+1} t P_{n-m-2}[\cos(z-t)]P_n^{-m}(\cos t)dt=\frac{\Gamma(m+3/2)}{2^{1/2}\Gamma(m+2)} \sin^{m+3/2}z P_{n-1/2}^{-m-1/2}(\cos z),\quad \Re m>-1.
\ee
\be
\int_0^z \sin^{-1-k}t \sin^k(z-t) P_n^{-m}(\cos t)P_n^{-k}[\cos (z-t)]
= \frac{2^k\Gamma(m-k)\Gamma(k+1/2)}{\sqrt{\pi} \Gamma(k+m+1)} \sin^kz P_n^{-m}(\cos z),\quad \Re m>\Re k > -1/2.
\ee
\be
\int_0^z \sin^mt \sin^{-m}(z-t) P_n^{-m}(\cos t)P_{n-1}^m[\cos (z-t)]
= \frac{\sin[(m+n)z]}{(m+n)\cos m\pi},\quad -1/2<\Re m<1/2.
\ee
\begin{multline}
\int_0^z P_{n+1}^{-m}(\cos t)\sin[n(z-t)]dt
= 2^{3/2}\sin^{1/2}z\sum_{r=0}^\infty \frac{(-1)^r(m+n+2)_{2r}(3/2)_r\Gamma(m+r+3/2)}
{r!\Gamma(m+r+1)}
\\ \times
\frac{m+2r+3/2}{(m+2r+1)(m+2r+2)} P_{n-1/2}^{-(m+3/2+2r)}(\cos z).
\end{multline}

\cite{BaileyPCPS27_2}
\be
\int_0^z \sin^m t P_{n-m-1}[\cos(z-t)] P_n^{-m}(\cos t) dt
=\frac{\Gamma(m+1/2)}{\sqrt 2 \Gamma(m+1)}\sin^{m+1/2}z P_{n-1/2}^{-m-1/2}(\cos z),\quad \Re m >-1/2.
\ee
\cite{BaileyPCPS27_2}
\be
\int_0^z \sin^{m+1} t P_{n-m-2}[\cos(z-t)] P_n^{-m}(\cos t) dt
=\frac{\Gamma(m+1/2)}{\sqrt 2 \Gamma(m+2)}\sin^{m+3/2}z P_{n-1/2}^{-m-1/2}(\cos z),\quad \Re m >-1.
\ee
\cite{BaileyPCPS27_2}
\be
\int_0^z P_n[\cos(z-t)] P_{n-1}(\cos t) dt
=\frac{\sin nz}{n}.
\ee
\cite{BaileyPCPS27_2}
\be
\int_0^z P_{n-2}[\cos(z-t)] \sin nt dt
=n \sin z P_{n-1}^{-1}(\cos z).
\ee
\cite{BaileyPCPS27_2}
\be
\int_0^z \sin tP_{n-2}[\cos(z-t)]P_n(\cos t)dt
=\frac{\sin z sin nz}{2n}.
\ee
\cite{BaileyPCPS27_2}
\be
\int_0^z \sin tP_{n-3}[\cos(z-t)]\sin nt dt
=\frac23 n\sin^2 z P_{n-1}^{-1}(\cos z).
\ee
\cite{BaileyPCPS27_2}
\be
\int_0^z P_{n-2}[\cos(z-t)]\cos nt dt
=\sin^z P_{n-1}(\cos z).
\ee

\cite{NevesJPA39}
\be
\int_0^\pi d\theta \sin^{|m|+1}\theta \exp(\pm i R\cos\theta)P_n^{|m|}(\cos\theta)
=
2(\pm i)^{n+|m|}\frac{ (n+|m|)!}{(n-|m|)!}\,\frac{j_n(R)}{R^{|m|}}
.
\ee

\cite[(23)]{NevesJPA39}
\be
\int_0^\pi d\theta \sin\theta \exp( i R\cos\alpha\cos\theta)P_n^m(\cos\theta)
J_m(R\sin\alpha\sin\theta)
=
2i^{n-m}P_n^m(\cos\alpha)j_n(R)
.
\ee

\cite[\S 11.4,13]{Apelblat}\cite{MatharIJQC90}
\begin{multline}
\int_0^\infty e^{-c^2x^2}H_{2m}(ax)H_{2k}(bx)dx =
\\
 (-)^{m+k}
2^{2m+2k-1}\frac{\Gamma(m+1/2)\Gamma(k+1/2)}{\pi c^{2m+2k+1}}
(c^2-a^2)^m(c^2-b^2)^k\,_2F_1(-m,-k;\frac{1}{2};\frac{a^2b^2}{(c^2-a^2)(c^2-b^2}),
\end{multline}
$c^2-a^2-b^2>0$.

\cite[\S 11.4,14]{Apelblat}
\begin{multline}
\int_0^\infty e^{-c^2x^2}H_{2m+1}(ax)H_{2k+1}(bx)dx =
 (-)^{m+k}
2^{2m+2k+1}\frac{\Gamma(m+3/2)\Gamma(k+3/2)}{\Gamma(3/2)}
\\
\times
\frac{ab(c^2-a^2)^m(c^2-b^2)^k}{c^{2m+2k+3}}
\,_2F_1(-m,-k;\frac{3}{2};\frac{a^2b^2}{(c^2-a^2)(c^2-b^2}).
\end{multline}

\cite[\S 11.4,11]{Apelblat}
\be
\int_{-\infty}^\infty e^{-2x^2}\left[H_n(x)\right]^2dx
=2^{n-1/2}\Gamma(n+1/2).
\ee

\cite[\S 11.4,20]{Apelblat}
\be
\int_{-\infty}^\infty e^{-(a^2+b^2)x^2} H_{2m}(ax)H_{2k}(bx)dx
=
(-)^{m+k}2^{2(m+k)}\Gamma(m+k+\frac{1}{2})\frac{a^{2k}b^{2m}}{(a^2+b^2)^{m+k+1/2}}
.
\ee

\cite[\S 11.4,21]{Apelblat}
\be
\int_{-\infty}^\infty e^{-(a^2+b^2)x^2} H_{2m+1}(ax)H_{2k+1}(bx)dx
=
(-)^{m+k}2^{2(m+k+1)}\Gamma(m+k+\frac{3}{2})\frac{a^{2k+1}b^{2m+1}}{(a^2+b^2)^{m+k+3/2}}
.
\ee

\cite{KolbigJCAM71}
\be
\int_0^\infty x^{\nu+1}
e^{-\alpha x^2} L_m^{\nu-\sigma}(\alpha x^2)L_n^{\sigma}(\alpha x^2) J_{\nu}(xy)dx
=
\frac{(-1)^{m+n}}{2\alpha}
\left(\frac{y}{2\alpha}\right)^\nu \exp\left(-\frac{y^2}{4\alpha}\right)
L_m^{\sigma-m+n}(\frac{y^2}{4\alpha})L_n^{\nu-\sigma+m-n}(\frac{y^2}{4\alpha}),
\ee
for
$y>0$, $\Re\alpha>0$, $Re\nu>-1$.

\be
\int_0^\infty x^{(\alpha+\beta)/2}
e^{-x} L_m^\alpha(x)L_n^\beta(x) J_{\alpha+\beta}(2\sqrt{ax})dx
=
(-1)^{m+n}a^{(\alpha+\beta)/2}L_m^{\beta-m+n}(a)L_n^{\alpha+m-n}(a);
\quad \Re(\alpha+\beta)>-1.
\ee

\cite{AfanasievJCP69}
\be
\int_0^\infty \frac{dx}{(1+x^2)^{3/4}} Q_{1/2}(\frac{1+r^2+x^2}{2r\sqrt{1+x^2}})
=\frac{\pi}{\surd 2}(1-\sqrt{1-r^2}),
\ee
where $Q$ are Legendre Functions of the 2nd Kind.

\subsection{Hypergeometric Functions}

\cite{YangIJMEST23}
\be
\int_0^1 {}_2F_1(1/2,1;3/2;-x^2) dx = G
\ee
where $G$ is Catalan's constant.

\cite[p 238]{ErdelyiT1}
\be
\frac{1}{\Gamma(2\lambda+2\nu)}
\int_0^\infty dt
e^{-pt}
t^{2\lambda+2\nu-1}
\, _1F_2(\nu;\lambda+\nu,\lambda+\nu+\frac{1}{2};-\frac{1}{4}a^2t^2)
=\frac{1}{p^{2\lambda}}\,\frac{1}{(p^2+a^2)^\nu}
;
\quad \Re(\lambda+\nu)>0
.
\ee

previous formula at $\lambda=0$ with \cite[9.1.69]{AS}
gives \cite[6.623.1]{GR}
\begin{multline}
\frac{1}{\Gamma(2\nu)}
\int_0^\infty dt
e^{-pt}
t^{2\nu-1}
\, _0F_1(;\nu+\frac{1}{2};-\frac{1}{4}a^2t^2)
=
\frac{\Gamma(\nu+1/2)}{\Gamma(2\nu)}
\left(\frac{a}{2}\right)^{1/2-\nu}
\int_0^\infty
e^{-pt}
t^{\nu-1/2}
J_{\nu-1/2}(at)
\\
=
\frac{1}{(p^2+a^2)^\nu}
;
\quad \Re(\nu)>0
.
\end{multline}

\cite{MavromatisJCAM59}
\be
\int_0^\infty J_\mu^2(\omega\rho)\,_3F_2
\left(\begin{array}{c}3/2,(\sigma+\nu)/2,(\sigma-\nu)/2
\\ 1-\mu,1+\mu\end{array}\mid -4\omega^2\right)\omega d\omega
= -\frac{\mu\rho^{\sigma-2}K_\nu(\rho)}
{2^{\sigma-1}\Gamma[(\sigma+\nu)/2]\Gamma[(\sigma-\nu)/2]}
\ee
for $\Re\sigma >1+|\Re\nu|$, $\mu=l+1/2$ with $l\in \mathbb{N}$,
$\Re(\sigma+\nu)>0$, $\Re\rho>0$.

\cite{MavromatisJCAM59}
\be
\int_0^\infty J_\mu^2(\omega\rho) J_{\mu+1}(\omega\rho) \,_3F_2
\left(\begin{array}{c}3/2,(\sigma+\nu)/2,(\sigma-\nu)/2
\\ 1-\mu,2+\mu\end{array}\mid -4\omega^2\right)\omega^2 d\omega
= -\frac{\mu(\mu+1)\rho^{\sigma-3}K_\nu(\rho)}
{2^{\sigma-1}\Gamma[(\sigma+\nu)/2]\Gamma[(\sigma-\nu)/2]}
\ee
for $\Re\sigma >2+|\Re\nu|$, $\mu=l+1/2$ with $l\in \mathbb{N}$,
$\Re(\sigma+\nu)>0$, $\Re\rho>0$.

\subsection{Confluent Hypergeometric Functions}
\subsection{Parabolic Cylinder Functions}
Corrigendum to \cite[7.741.5]{GR}

\cite{MatharVixra2207}
\begin{equation}
\int_0^\infty e^{-x^2/4}\cos(bx) \left[D_{2\nu-1/2}(x)+D_{2\nu-1/2}(-x)\right]dx
= 
\sqrt{2\pi} 
e^{-b^2/2}
b^{2\nu-1/2}
\sin[\pi(\nu+1/4)]
.
\end{equation}
This is the limit $a\to 1/4^+$ of the more general
\begin{multline}
\int_0^\infty e^{-ax^2}\cos(bx) \left[D_{2\nu-1/2}(x)+D_{2\nu-1/2}(-x)\right]dx
\\
=
\frac{2^{\nu-1/4}\pi}{\Gamma(3/4-\nu)}
\frac{(a-1/4)^{\nu-1/4}}{(a+1/4)^{\nu+1/4}}
e^{-b^2/(4a+1)}
{}_1F_1\left(-\nu+\frac14;\frac12;-\frac{2b^2}{(4a+1)(4a-1)}\right),\quad a>1/4,\Re\nu>1/4,b>0.
\end{multline}

\subsection{Laguerre Polynomials}
\cite{LordMC14}
\be
C_{r,s,t}\equiv = \int_0^\infty e^{-x} L_r(x)L_s(x)L_t(x) dx
=
(-2)^{r+s+t}\sum_n 2^{-2n}\frac{n!}{(n-r)!(n-s)!(n-t)!(r+s+t-2n)!}
\ee
with recurrence
\be
C_{r,s,t} =
C_{r,s-1,t-1} 
+C_{r-1,s,t-1} 
+C_{r-1,s-1,t} 
+2C_{r-1,s-1,t-1},
\ee
where $L_n(x)=\sum_{\nu=0}^n (-1)^\nu\binom{n}{\nu}\nu^r/r!$.

\section{Special Functions}

\subsection{Elliptic Integrals and Functions}
\cite{LavoieMC49}
\be
E(2^{-1/2})=\frac{\Gamma^2(1/4)}{8\surd \pi}+\frac{\pi^{3/2}}{\Gamma^2(1/4)}
=\frac12 K(2^{-1/2})+\frac{\pi}{2K(2^{-1/2})}.
\ee

\cite{ZuckerMPCPS131}
\be
K(1/\surd 2) = \frac{\Gamma^2(1/4)}{4\surd \pi}.
\ee

\cite{ZuckerMPCPS131}
\be
K(\sqrt{2}-1) = \frac{2^{3/4}}{16\surd \pi} \sqrt{\sqrt{2}+1}\Gamma(1/8)\Gamma(3/8).
\ee

\cite{ZuckerMPCPS131}
\be
K(\frac14\sqrt{2}(\sqrt{3}-1)) = \frac{2^{2/3}3^{1/4}}{8\pi}\Gamma^3(1/3) .
\ee

\cite{ZuckerMPCPS131}
\be
K((\sqrt{3}-\sqrt{2})(2-\sqrt{3})) = \frac{2^{1/12}3^{1/4}}{48}(\sqrt{2}-1)(\sqrt{3}+1)\frac{\Gamma^2(1/24}){\Gamma(1/12)}.
\ee

\cite{ZuckerMPCPS131}
\be
K(\frac{1}{8}\sqrt{2}(3-\sqrt{7})) = \frac{2^{5/2}7^{3/4}}{14} \frac{\Gamma^2(1/7)\Gamma^2(2/7)\Gamma(3/4)}{\Gamma^2{3/8}\Gamma(2/7)\Gamma(4/7)}.
\ee

\subsection{The exponential integral and related functions}

\cite{BerndtBLMS15}
If
\be
\int_0^x \frac{\sin u}{u}du=\frac{\pi}{2}-r\cos(x-\theta);\,
\int_0^x \frac{1-\cos u}{u}du=\gamma +\log x-r\sin(x-\theta)
\ee
then
\be
r\cos\theta \sim \sum_{k\ge 0}\frac{(-)^k(2k)!}{x^{2k+1}};\,
r\sin\theta \sim \sum_{k\ge 1}\frac{(-)^{k+1}(2k-1)!}{x^{2k}};\,
r \sim \sum_{k\ge 1}\frac{(-)^{k+1}(2k-1)!}{kx^{2k}}
\ee
for $x\to\infty$.

\cite{WatrasieOA14}
Derivatives:
\be
Si^{(n+1)}(ax) = \sum_{k=0}^n a^{n+k}\frac{n!}{(n-k)!}\frac{\sin[ax+(n+k)\pi/2]}{x^{k+1}}.
\ee
\be
Ci^{(n+1)}(ax) = \sum_{k=0}^n a^{n-k}\frac{n!}{(n-k)!}\frac{\cos[ax+(n+k)\pi/2]}{x^{k+1}}.
\ee

\subsection{The error function and Fresnel integrals}

\cite{BorosSci12}
\be
2\sqrt{\pi}e^{b^2}\erf(b)
= \sum_{j=0}^\infty \frac{j!b^{2j+1}}{(2j+1)!2^{2j+2}}
.
\ee

\subsection{The gamma function}

\be
\frac{1}{a+r}=\frac{1}{a}\,\frac{(a)_r}{(a+1)_r}.
\label{eq.poch1}
\ee

\cite[2.2.3.1]{SlaterHyp}\cite{RoyAMM94,ChoEAMJ15}
\be
(a)_{m-r}=\frac{(-1)^r (a)_m}{(1-a-m)_r}.
\ee

\cite[(2.4.5.2.)]{SlaterHyp}\cite{DriverETNA25,RoyAMM94}
\be
(a)_{2r}=(a/2)_r (\frac{a+1}{2})_r 2^{2r}
;
\quad
(a)_{qr}=(a/q)_r (\frac{a+1}{q})_r\cdots (\frac{a+q-1}{q})_r  q^{qr}
\ee

\cite{SlaterHyp,ChoEAMJ15}
\be
\Gamma(n+1-k)=\frac{(-1)^k\Gamma(n+1)}{(-n)_k}.
\label{eq.pochmk}
\ee

\cite[(I.4)]{SlaterHyp}
\be
(a+kn)_n=\frac{(a)_{(k+1)n}}{(a)_{kn}}
.
\ee

\cite[(I.6)]{SlaterHyp}
\be
(a-kn)_n=\frac{(-1)^n (1-a)_{kn}}{(1-a)_{(k-1)n}}
.
\ee

\cite[(I.11)]{SlaterHyp}
\be
(a+kn)_{N-n}=\frac{(a)_N (a+N)_{(k-1)n}}{(a)_{kn}}
.
\ee

\cite[(I.13)]{SlaterHyp}
\be
(a-kn)_{N-n}=\frac{(-1)^n (a)_N (1-a)_{kn}}{(1-a-N)_{(k+1)n}}
.
\ee

\cite[(\S 1.2)]{ErdelyiI}
\be
\frac{\Gamma(n+\frac{1}{2}+z)\Gamma(n+\frac{1}{2}-z)}
{\Gamma^2(n+\frac{1}{2})}
=\frac{1}{\cos(\pi z)}\prod_{l=1}^n\left[1-\frac{4z^2}{(2l-1)^2}\right]
,
\ee
for $n=1,2,3,\ldots$.

\cite{Vidunasarxiv04}\cite[A073006]{sloane}
\be
\Gamma(2/3)=\frac{2\pi}{\sqrt 3}\,\frac{1}{\Gamma(1/3)}
.
\ee
\cite{Vidunasarxiv04}\cite[A068465]{sloane}
\be
\Gamma(3/4)=\pi\sqrt{2}\,\frac{1}{\Gamma(1/4)}
.
\ee
\cite{Vidunasarxiv04}\cite[A175379]{sloane}
\be
\Gamma(1/6)=\frac{\sqrt 3}{\sqrt \pi 2^{1/3}}\Gamma^2(1/3)
.
\ee
\cite{Vidunasarxiv04}\cite[A340721]{sloane}
\be
\Gamma(3/5)=\frac{\pi\sqrt{2}\sqrt{\phi^*}}{\sqrt 5}
\frac{1}{\Gamma(2/5)}
,
\ee
where $\phi^*=5-\surd 5$.

\cite{Vidunasarxiv04}\cite[A203145]{sloane}
\be
\Gamma(5/6)=\frac{\pi^{3/2}2^{4/3}}{\sqrt 3}
\frac{1}{\Gamma^2(1/3)}
.
\ee

\cite{Vidunasarxiv04}\cite[A340722]{sloane}
\be
\Gamma(4/5)=\frac{\pi\sqrt{2}\sqrt{\phi}}{\sqrt 5}
\frac{1}{\Gamma(1/5)}
,
\ee
where $\phi=5+\surd 5$.

\cite{Vidunasarxiv04}\cite[A203143]{sloane}
\be
\Gamma(3/8)=\sqrt{\pi} \sqrt{\surd 2 -1}
\frac{\Gamma(1/8)}{\Gamma(1/4)}
.
\ee

\cite{Vidunasarxiv04}\cite[A203144]{sloane}
\be
\Gamma(5/8)=\sqrt{\pi} 2^{3/4}
\frac{\Gamma(1/4)}{\Gamma(1/8)}
.
\ee

\cite{Vidunasarxiv04}\cite[A203146]{sloane}
\be
\Gamma(7/8)=\pi 2^{3/4}
\sqrt{\surd 2 +1}
\frac{1}{\Gamma(1/8)}
.
\ee

\cite{Vidunasarxiv04}
\be
\Gamma(1/10)=\frac{\surd \phi}
{\sqrt \pi 2^{7/10}}
\Gamma(1/5)\Gamma(2/5)
,
\ee
where $\phi=5+\surd 5$.

\cite{Vidunasarxiv04}\cite[A340723]{sloane}
\be
\Gamma(3/10)=\frac{\sqrt\pi \phi^*}
{2^{3/5}\sqrt 5}
\frac{\Gamma(1/5)}{\Gamma(2/5)}
,
\ee
where $\phi^*=5-\surd 5$.

\cite{Vidunasarxiv04}\cite[A340724]{sloane}
\be
\Gamma(7/10)=\sqrt\pi 2^{3/5}
\frac{\Gamma(2/5)}{\Gamma(1/5)}
.
\ee

\cite{Vidunasarxiv04}\cite[A340725]{sloane}
\be
\Gamma(9/10)=\frac{\pi^{3/2}2^{7/10}\surd \phi}{\surd 5}
\frac{1}{\Gamma(2/5)\Gamma(1/5)}
.
\ee

\cite{Vidunasarxiv04}\cite[A203140]{sloane}
\be
\Gamma(1/12)=\frac{3^{3/8}\sqrt{\surd 3+1}}{\sqrt \pi 2^{1/4}}
\Gamma(1/3)\Gamma(1/4)
.
\ee

\cite{Vidunasarxiv04}
\be
\Gamma(5/12)=\frac{\sqrt \pi 2^{1/4}\sqrt{\surd 3-1}}{3^{1/8}}
\frac{\Gamma(1/4)}{\Gamma(1/3)}
.
\ee
\cite{ZuckerMPCPS131}
\be
\frac{\Gamma(1/20)\Gamma(9/20)}{\Gamma(3/20)\Gamma(7/20)}=\frac{5^{1/4}}{2}(\sqrt{5}+1).
\ee
\cite{ZuckerMPCPS131}
\be
\frac{\Gamma(1/24)\Gamma(11/24)}{\Gamma(5/24)\Gamma(7/24)}=\sqrt{3}(\surd 3 +2)^{1/2}.
\ee

\cite{RagabGMJ6}
\be
\Gamma(\frac{-\nu-l}{n})\ldots
\Gamma(\frac{-\nu-l+n-1}{n})=
(-)^{l-1}\frac{(2\pi)^{(n-1)/2}n^{1/2+\nu+l}\pi}{\Gamma(\nu+l+1)\sin\nu\pi}
\ee

\cite{RagabGMJ6}
\be
\Gamma(\frac{-l}{n})\ldots
\Gamma(-\frac{1}{n})
\Gamma(\frac{1}{n})
\ldots
\Gamma(\frac{-l+n-1}{n})=
(-)^{l}\frac{(2\pi)^{(n-1)/2}n^{l-1/2}}{\Gamma(l+1)}
\ee

\cite{RottbrandITSF10}
\be
(N+\frac{1}{2},K)
(N+\frac{1}{2},L)
=\sum_{j=\max(K,L)}^{\min(N,K+L)}
\binom{2j-K-L}{j-K}
(j+\frac{1}{2},K+L-j)
(N+\frac{1}{2},j)
,
\ee
where $(n+\frac{1}{2},j)=\frac{(n+j)!}{j!(n-j)!}$ is Hankel's symbol.

\cite{RottbrandITSF10}
\be
(-)^N
\frac{
(N+\frac{1}{2},K)
(N+\frac{1}{2},L)
}{\binom{K+L}{K}}
=\sum_{j=\max(K,L)}^{\min(N,K+L)}
(-)^j
\binom{2j-K-L}{j-K}
\frac{(j+1/2,K+L-j)(N+1/2,j)}{\binom{K+L}{j}}
,
\ee
where $(n+\frac{1}{2},j)=\frac{(n+j)!}{j!(n-j)!}$ is Hankel's symbol.

\cite{JordanSIAMR26}
\be
\log\frac{\Gamma(z+1/2)}{\sqrt{z}\Gamma(z)}
=
-\sum_{r=1}^k\frac{(1-2^{-2r})B_{2r}}{r(2r-1)z^{2r-1}}+O(z^{-2k+1/2}).
\ee

\cite{JordanSIAMR26}
\be
\log\frac{\Gamma(z+3/4)}{\sqrt{z}\Gamma(z+1/4)}
=
-\sum_{r=1}^k\frac{E_{2r}}{4r(4z)^{2r}}+O(z^{-2k+1/2}).
\ee

\cite{JordanSIAMR26}
\be
\frac{1}{z}\left[ \frac{\Gamma(z+3/4)}{\Gamma(z+1/4)}\right]^2
=1+\frac{2u}{1+}\frac{9u}{1+}\frac{25u}{1+}\frac{49u}{1+}\cdots
\ee
where $u=1/(64z^2)$.

\cite{TricomiPJM1}
\begin{equation}
\frac{\Gamma(z+\alpha)}{\Gamma(z+\beta)} 
\sim \sum C_n(\alpha-\beta,\beta) z^{\alpha-\beta-n},\quad z\to\infty.
\end{equation}
where
\begin{equation}
C_n(\alpha',\beta)=\frac{1}{n}\sum_{m=0}^{n-1}[\binom{\alpha'-m}{n-m+1}-(-)^{n+m}\alpha'\beta^{n-m}] C_m(\alpha',\beta),\quad n=1,2,3\ldots
\end{equation}
\begin{equation}
C_0=1,\quad
C_1=\frac12 \alpha'(\alpha'+2\beta-1) = \frac12 (\alpha-\beta)(\alpha+\beta-1),
\end{equation}
\begin{equation}
C_1=\frac{1}{12} [(\alpha'-2)(3\alpha'-1)+12\beta(\alpha'+beta-1)].
\end{equation}

\cite{ChoiJMAA231}
Let $1/G(z)$ be the reciprocal of the double Gamma function,
\be
\Gamma(1)=G(1)=1;\quad G(z+1)=\Gamma(z)G(z).
\label{eq.A}
\ee
Let $A\approx 1.28242713$ be Glaisher's (Kinkelin's) constant \cite[A074962]{sloane}. Then
there is a Stirling formula
\be
\log G(x+a+1)=\frac{x+a}{2}\log(2\pi)-\log A+\frac{1}{12}-\frac{3x^2}{4}-ax
+(\frac{x^2}{2}-\frac{1}{12}+\frac{a^2}{2}+ax)\log x+O(1/x),\quad x\to \infty.
\ee
Let $G\approx 0.91596$ be the Catalan constant. Then
\be
G(3/4)=2^{-1/8}\pi^{-1/4}e^{G/2\pi}\Gamma(1/4)G(1/4).
\ee
The duplication formula is 
\be
G(a)G^2(a+1/2)G(a+1)=e^{1/4}A^{-3}2^{-2a^2+3a-11/12}\pi^{a-1/2}G(2a).
\ee
\cite{ChoiJMAA231}
\be
G(1/4)=e^{3/32-G/4\pi}A^{-9/8}\Gamma^{-3/4}(\frac14).
\ee
\cite{ChoiJMAA231}
\be
G(3/4)=2^{-1/8}\pi^{-1/4}e^{3/32+G/4\pi}A^{-9/8}\Gamma^{1/4}(\frac14).
\ee
\cite{ChoiJMAA231}
\be
\frac{d}{dz}\log G(z)=\frac{1+\log(2\pi)}{2}-z+(z-1)\psi(z).
\ee

\subsection{The Digamma function $\psi$}

\cite{RamMurtyJNT125}
\be
-\psi(a/q)-\gamma = \log(q)-\sum_{b=1}^{q-1}
\zeta_q^{-ba}\log(1-\zeta_q^b),
\ee
where $\zeta_q$ is the primitive $q$th root of unity, $e^{2\pi i/q}$.

\cite{SofoIJPAM50}
\be
\psi(a/b) = -\gamma -\ln(2b)-\frac{\pi}{2}\cot\frac{\pi a}{b}
+2\sum_{\mu=1}^{[b/2]-1}\cos\frac{2\pi a\mu}{b}\ln \sin \frac{\pi\mu}{b}.
\ee
for $0<a<b$.

\cite[(2.11)]{MuthumSci22}
\be
\psi(1/2+z)=\psi(1/2-z)+\pi \tan \pi z.
\ee

\cite[(2.11)]{MuthumSci22}
\be
\psi(n+1/2+z)-\psi(n+1/2-z) = \pi \tan \pi z-2\sum_{k=0}^{n-1} \frac{z}{(k+1/2)^2-z^2}.
\ee
The derivative of this is
\be
\psi'(n+1/2+z)+\psi'(n+1/2-z) = (\frac{\pi}{\cos(\pi z)})^2-2\sum_{k=0}^{n-1} \frac{(k+1/2)^2+z^2}{[(k+1/2)^2-z^2]^2}.
\ee

A consequence of the duplication formula for the $\Gamma$-function is \cite[8.365.6]{GR}\cite{BlumleinCPC180}
\be
\psi(\frac{z}{2}) = 2\psi(z-1)-\psi(\frac{z-1}{2}) -2 \ln 2.
\ee
\be
\psi'(\frac{z}{2}) = 4\psi'(z-1)-\psi'(\frac{z-1}{2}) .
\ee

\be
\psi'(x+1)=\psi'(x)-\frac{1}{x^2}.
\ee

\cite{KolbigJCAM18}
\be
\psi^{(k)}(1+n)=(-)^{k+1}k![\zeta(k+1)-\sum_{j=1}^n\frac{1}{j^{k+1}}].
\ee

\cite{KolbigJCAM18}
\be
\psi^{(k)}(n+1/2)=(-)^{k+1}k![(2^{k+1}-1)\zeta(k+1)-\sum_{j=1}^n(\frac{2}{2j-1})^{k+1}].
\ee

\cite[B.83]{BlumleinCPC180}\cite{SofoIJPAM50}
\be
\psi^{(n)}(1+z)=\psi^{(n)}(z) +(-)^n\frac{n!}{z^{n+1}}.
\ee
\cite[B.110]{BlumleinCPC180}
\be
\psi^{(n)}(z)\sim 
(-)^{n-1}\left[\frac{(n-1)!}{z^n}+\frac{n!}{2z^{n+1}}
+\sum_{k=1}^\infty B_{2k}\frac{(2k+n-1)!}{(2k-1)!z^{2k+n}}
\right]
.
\ee

\cite{KrupnikovJCAM78}
\be
\psi'(1/4)=\pi^2+8G ; \quad \psi'(3/4) = \pi^2-8G.
\ee
\be
\psi''(1/4)=-2(\pi^3+28\zeta(3)) ; \quad \psi''(3/4) = 2(\pi^3-28\zeta(3)).
\ee
\be
\psi'''(1/4)=8(\pi^4+96\beta^*(4)) ; \quad \psi'''(3/4) = 8(\pi^4-96\beta^*(4)).
\ee
where $\beta^*(n)\equiv \sum_{k\ge 0}(-1)^k/(2k+1)^n$.
\be
\psi(1/6)=-\gamma-2\ln 2-\frac32 \ln 3 -\frac12 \sqrt{3}\pi ;
\ee
\be
\psi(1/3)=-\gamma-\frac32 \ln 3-\frac16 \sqrt{3}\pi ;
\ee
\be
\psi'(1/6)=5\psi'(1/3)-\frac43 \pi^2;
\ee
\be
\psi''(1/6)=-182\zeta(3)-4\sqrt{3}\pi^3 ;
\ee
\be
\psi''(1/3)=-26\zeta(3)-\frac49 \sqrt{3}\pi^3 ;
\ee

\cite{SofoIJPAM50}
\be
\psi^{(m)}(1+z)=(-)^mm!(H_z^{(m+1)}-\zeta(1+z)).
\ee
where $H_n^{(m)}=\sum_{r=1}^n 1/r^m$ are harmonic numbers.

\cite{SavioHarm}
\be
\psi_m(n-1)=\psi_m(0)+(-)^{m-1}(m-1)!H_n^{(m)}.
\ee

\cite{SofoIJPAM50}
\be
\sum_{r=1}^k (-)^r\binom{k-1}{r-1}\psi^{(m)}(\frac{r+j}{k})
=m!\sum_{r=1}^k (-)^{r+m}\binom{k-1}{r-1}H_{(r+j)/k-1}^{(m+1)}.
\ee

\cite{SofoIJPAM50}
\be
\sum_{r=1}^k (-)^r\binom{k-1}{r-1}\psi^{(m)}(\frac{r}{k})
=m!\sum_{r=1}^{k-1} (-)^{r+m}\binom{k-1}{r-1}H_{r/k-1}^{(m+1)}.
\ee

\cite{Adamchik}
\be
\sum_{k=0}^\infty\frac{(-1)^k}{2k+1}[\psi(k+1)+\gamma]=G-\frac{\pi}{2}\log 2.
\ee

\cite{Adamchik}
\be
\sum_{k=0}^\infty\frac{(-1)^k}{2k+1}[\psi(k+3/2)+\gamma]=G-\frac{\pi}{4}\log 2.
\ee

\cite{Adamchik}
\be
\frac{(-1)^n}{q^n(n-1)!}\sum_{k=1}^{q-1}e^{2\pi ikp/q}\psi^{n-1}(k/q)
=
\Li_n(e^{2\pi ip/q})-\frac{\zeta(n)}{q^n}.
\ee

\cite{Adamchik}
\be
\frac{1}{4}\sum_{k=0}^\infty\frac{2^k k!^2}{(2k+1)!}[\psi(k+3/2)+\gamma]=G-\frac{\pi}{4}\log 2.
\ee

\cite{GeddesCM}
\be
\psi^{(n)}(z)=(-v)^{n+1}\int_0^\infty \frac{x^n\exp(-zvx)}{1-\exp(-vx)}dx.
\ee

\cite{Adamchik}
\be
q\sum_{k=0}^{q-1}\Li_2(e^{2\pi ik/q}t)
=
\Li_2(t^q).
\ee

\cite{YangIJMEST23}
\be
G=\frac{1}{4}\Phi(-1,2,1/2) = \frac{1}{16}[\zeta(2,1/4)-\zeta(2,3/4)].
\ee
where $\Phi(z,s,\alpha)$ is Lerch's function and $G$ the Catalan constant.

\subsection{Bessel Functions and Functions Associated with Them}
\subsection{Integral representations of the functions $J_\nu$ and $N_\nu$}
\cite[La. 4.13]{Stein}
\be
J_{\mu+\nu+1}(t)=\frac{t^{\nu+1}}{2^\nu\Gamma(\nu+1)}
\int_0^1 J_{\mu}(ts)s^{\mu+1}(1-s^2)^{\nu}ds.
\ee

\cite[7.2.7]{Erdelyi}\cite{LandauJCAM132}\cite[9.1.14]{AS}
\be
\Gamma(\nu+1)\Gamma(\mu+1)J_\nu(z)J_\mu(z)
=
\left(z/2\right)^{\nu+\mu}
\,_2F_3\left(
\begin{array}{c}\frac{1+\nu+\mu}{2},1+\frac{\nu+\mu}{2}\\ 1+\nu,1+\mu,1+\nu+\mu
\end{array}\mid-z^2\right).
\ee
Special case of this:
\be
\Gamma^2(\nu+1)J_\nu^2(z)
=
\left(z/2\right)^{2\nu}
\,_1F_2\left(\begin{array}{c}\nu+1/2\\ 1+\nu,1+2\nu\end{array}\mid-z^2\right).
\ee
Neumann's expansion of the previous equation \cite[9.1.86]{AS}
\be
J_\nu^2(z)
= \left(z/2\right)^{2\nu}\sum_{k\ge 0}a_{2k} J_{2k}(z)
\ee
where
\begin{multline}
a_2(\nu) = 2\frac{\nu^2-1}{\Gamma^2(\nu+2)};
a_4(\nu) = 2\frac{(\nu+2)(\nu^3-4\nu^2+5\nu+22)}{\Gamma^2(\nu+3)};
\\
a_6(\nu) = 2\frac{(\nu+2)(\nu+3)(\nu^4-11\nu^3+53\nu^2-13\nu-654)}{\Gamma^2(\nu+4)};
\\
a_8(\nu) = 2\frac{(\nu+3)(\nu+4)(\nu^6-19\nu^5+163\nu^4-385\nu^3-3836\nu^2+18764\nu+56592)}{\Gamma^2(\nu+5)};
\end{multline}
This shows quickly increasing $a_{2k}$ as $k$ grows because the expansion
tries to approximate the positive $J_\nu^2$ by a series of oscillating functions.

\cite{AljarrahJCAM143}
\be
\sum_{k=0}^\infty I_{nk}(\xi)
=\frac12 I_0(\xi)+\frac{1}{2n}\sum_{k=0}^{n-1} e^{\xi\cos(2\pi k/n)}.
\ee

\cite{AljarrahJCAM143}
\be
\sum_{k=0}^\infty (\pm i)^{nk} J_{nk}(z)
=\frac12 J_0(z)+\frac{1}{2n}\sum_{k=0}^{n-1} e^{\pm iz\cos(2\pi k/n)}\, ,n=1,2,\ldots.
\ee

\cite{AljarrahJCAM143}
\be
\sum_{k=0}^\infty (-)^{nk} I_{nk}(z)
=\frac12 I_0(z)+\frac{1}{2n}\sum_{k=0}^{n-1} e^{-z\cos(2\pi k/n)}\, ,n=1,2,\ldots.
\ee

\cite{AljarrahJCAM143}
\be
\sum_{k=0}^\infty J_{nk}(z)\cos(nk\pi/2)
=\frac12 J_0(z)+\frac{1}{2n}\sum_{k=0}^{n-1} \cos(z\cos(2\pi k/n))\, ,n=1,2,\ldots.
\ee

\cite{AljarrahJCAM143}
\be
\sum_{k=0}^\infty J_{nk}(z)\sin(nk\pi/2)
=\frac{1}{2n}\sum_{k=0}^{n-1} \sin(z\cos(2\pi k/n))\, ,n=1,2,\ldots.
\ee

\cite{AljarrahJCAM143}
\be
\sum_{k=0}^\infty (-)^k J_{2mk}(z)
=\frac12 J_0(z)+\frac{1}{2m}\cos z+
\frac{1}{m}\sum_{k=1}^{(m-1)/2} \cos(z\cos(2\pi k/m))\, ,m=1,3,\ldots.
\ee

\cite{AljarrahJCAM143}
\be
\sum_{k=0}^\infty J_{4mk}(z)
=\frac12 J_0(z)+\frac{1}{4m}(1+\cos z)+
\frac{1}{2m}\sum_{k=0}^{m-1} \cos(z\cos(\pi k/2m))\, ,m=1,2,\ldots.
\ee
There is no summation on the right hand side for the previous
two equations if $m=1$.

\cite{AljarrahJCAM143}
\be
\sum_{k=0}^\infty (-)^k J^2_{mk}(z)
=\frac12 J_0^2(z)+\frac{1}{2m}J_0(2z)
+\frac{1}{m}\sum_{k=1}^{(m-1)/2} J_0(2z\cos(\pi k/m))\, ,m=1,3,\ldots.
\ee

\cite{AljarrahJCAM143}
\be
\sum_{k=0}^\infty J^2_{2mk}(z)
=\frac12 J_0^2(z)+\frac{1}{4m}[1+J_0(2z)]
+\frac{1}{2m}\sum_{k=1}^{m-1} J_0(2z\cos(\pi k/2m))\, ,m=1,2,\ldots.
\ee

\cite{AljarrahJCAM143}
\be
\sum_{k,l=0}^\infty (-)^l I_{nk+2l+1}(z)
=\frac12 \sum_{l=0}^\infty (-)^l I_{2l+1}(z)+\frac{1}{4n}
\left\{
\sum_{k=0,k\neq n/4,3n/4}^{n-1}
\frac{e^{z\cos(2\pi k/n)}-1}{\cos(2\pi k/n)}+Nz
\right\}
\ee
whre $N=2$ if $4\div n$, otherwise $N=0$.

\cite{AljarrahJCAM143}
For $a>0$, $y>0$, and $m=1,3,\ldots$
\be
\sum_{k=0}^\infty (-)^k J_{mk}(y/2)Y_{mk}(y/2)
=\frac12 J_0(y/2)Y_0(y/2)+\frac{1}{2m}Y_0(y)
+\frac{1}{m}\sum_{k=1}^{(m-1)/2} Y_0(y\cos(\pi k/m)].
\ee
and
\be
\sum_{k=0}^\infty (-)^k I_{mk}(ay/2)K_{mk}(ay/2)
=\frac12 I_0(ay/2)K_0(ay/2)+\frac{1}{2m}K_0(ay)
+\frac{1}{m}\sum_{k=1}^{(m-1)/2} K_0(ay\cos(\pi k/m)].
\ee

\cite{BaileyPCPS25}
\begin{multline}
\frac{z}{2}J_\mu(z\cos\phi\cos\Phi)J_\nu(z\sin\phi\sin\Phi)
=(\cos\phi\cos\Phi)^\mu(\sin\phi\sin\Phi)^\nu
\\
\times
\sum_{n=0}^\infty (-1)^n(\mu+\nu+2n+1)J_{\mu+\nu+2n+1}(z)
\frac{\Gamma(\mu+\nu+n+1)\Gamma(\nu+n+1)}{n!\Gamma(\mu+n+1)\Gamma^2(\nu+1)}
\\
\times
F(-n,\mu+\nu+n+1;\nu+1;\sin^2\phi)
F(-n,\mu+\nu+n+1;\nu+1;\sin^2\Phi)
,
\end{multline}
where $\nu$ and $\mu$ are not negative integers.

\cite{DixonQJM4}
\be
J_\mu(Xz)J_\nu(xz)
=
\frac{X^\mu x^\nu}{\pi}
\int_{-\pi/2}^{\pi/2}
e^{i(\mu-\nu)\theta}
(\lambda_1/\lambda_2)^{\mu+\nu}J_{\mu+\nu}(z\lambda_1\lambda_2)d\theta
\ee
where
$\lambda_1 = +\surd(e^{i\theta}+e^{-i\theta})$,
$\lambda_2 = (X^2e^{i\theta}+x^2e^{-i\theta})^{1/2}$.

\cite{ConollyGMJ2}
\be
K_m(a)K_n(b)=\int_{-\infty}^{\infty} e^{-u(m-n)}
\left(\frac{ae^u+be^{-u}}{ae^{-u}+be^u}\right)^{\frac{m+n}{2}}
K_{m+n}(\sqrt{(ae^u+be^{-u})(ae^{-u}+be^u)})
du
\ee

\cite[p 99]{Erdelyi}
\begin{gather}
\left(z/2\right)^{\gamma-\mu-\nu}
J_\mu(\alpha z)J_\nu(\beta z)
=
\frac{\alpha^\mu \beta^\nu}{\Gamma(\mu+1)\Gamma(\nu+1)}\sum_{m=0}^\infty
\frac{(\gamma+2m)\Gamma(\gamma+m)}{m!}
\\
\times
F_4(-m,\gamma+m;\mu+1,\nu+1;\alpha^2,\beta^2)
J_{\gamma+2m}(z)
\nonumber
\\
=
\frac{\alpha^\mu\beta^\nu}{\Gamma(\nu+1)}
\sum_{m=0}^\infty (\gamma+2m)J_{\gamma+2m}(z)
\sum_{n=0}^\infty \frac{(-1)^n\Gamma(\gamma+m+n)\alpha^{2n}}{n!(m-n)!\Gamma^2(n+\mu+1)}
\,_2F_1(-n,-n-\mu;\nu+1;\frac{\beta^2}{\alpha^2}).
\end{gather}

Application of \cite[(3.5)]{Fields} to \cite[(9.1.14)]{AS}
\begin{multline}
J_\rho(cz)J_\nu(z)J_\mu(z) = \frac{c^\rho (z/2)^{\nu+\mu+\rho}}{
\Gamma(\rho+1)\Gamma(\nu+1)\Gamma(\mu+1)}
\sum_{n\ge 0}\frac{(-c^2 z^2/4)^n}{(\rho+1)_nn!}
\\
\times
\,_4F_3\left(-n,-n-\rho,1+\frac{\nu+\mu}{2},
\frac{1+\nu+\mu}{2};\nu+1,\mu+1,\nu+\mu+1;\frac{4}{c^2
}
\right)
.
\end{multline}

\cite{Blasiakarxiv08}
\be
x^{-\alpha}I_\alpha(2\sqrt{x})=\sum_{r=0}^\infty \frac{x^r}{r!\Gamma(r+\alpha+1)}
.
\ee

\cite{Blasiakarxiv08}\cite[A002426]{sloane}
\be
\exp(t)I_0(2t)=\sum_{n=0}^\infty \frac{t^n}{n!}c_n;
\quad
c_n\equiv \sum_{k=0}^{[n/2]}\frac{n!}{(k!)^2(n-2k)!}
=i^n\sqrt{3^n}P_n\left(-\frac{i}{\sqrt{3}}\right).
\ee

\cite{Blasiakarxiv08}
\be
\exp(yt)\exp(xt^2)^{-\alpha/2}I_\alpha(2t\sqrt{x})
=\sum_{n=0}^\infty \frac{t^n}{n!}\Pi_n^\alpha(x,y);
\quad
\Pi_n^\alpha(x,y)\equiv
n!\sum_{k=0}^{[n/2]}\frac{x^ky^{n-2k}}{(n-2k)!k!\Gamma(k+\alpha+1)}
.
\ee
\be
\frac{\exp(t)}{t}I_1(t)
=\sum_{n=0}^\infty \frac{t^n}{n!}\Pi_n^1(1,1).
\ee

\cite{NevesJPA39}
\be
j_s(R)=\frac{R^s}{2^{s+1}s!}\int_0^\pi d\theta \sin\theta \cos(R\cos\theta)
(\sin\theta)^{2s}
.
\ee

\cite[A122848]{sloane}
\begin{eqnarray}
\frac{d^s}{dt^s} t^\nu Z_\nu(t) &=& \sum_{l=\lfloor(s+1)/2\rfloor}^s
\alpha_{l,s} t^{\nu-(s-l)}Z_{\nu-l}(t).
\\
\frac{d^s}{dt^s} t^\nu K_\nu(t) &=& (-)^s\sum_{l=\lfloor(s+1)/2\rfloor}^s
\alpha_{l,s} t^{\nu-(s-l)}K_{\nu-l}(t)
\end{eqnarray}
with Bessel polynomial coefficients \cite{HanEJC29}
\[
\alpha_{l,s} = \frac{s!}{(s-l)!(2l-s)!2^{s-l}}
.
\]

\cite{RottbrandITSF10}
\be
\sum_{k=1}^n\frac{2k-1}{\sqrt\pi}K_{k-1/2}(x)K_{k-1/2}(y)
=
\sum_{k=1}^n \frac{n(n-1+k)!}{k!(n-k)!}
\left(\frac{x+y}{2xy}\right)^{k-1/2}
K_{k-1/2}(x+y)
.
\ee

\cite{RottbrandITSF10}
\be
\frac{1}{\sqrt\pi}K_{n+1/2}(x)K_{n+1/2}(y)
=
\sum_{\mu=0}^n \frac{(n+\mu)!}{\mu!(n-\mu)!}
\left(\frac{x+y}{2xy}\right)^{\mu+1/2}
K_{k+1/2}(x+y)
.
\ee

\cite{GlasserSIAMR22}
\be
\sum_{n=1}^\infty (-1)^n\frac{J_{2m}(n\pi)}{a^2-n^2} = \frac{\pi J_{2m}(a\pi)}{2a\sin a\pi}.
\ee

\cite{GlasserSIAMR22}
\be
\sum_{n=1}^\infty (-1)^n\frac{nJ_{2m-1}(n\pi)}{a^2-n^2} = \frac{\pi J_{2m-1}(a\pi)}{2\sin a\pi}.
\ee

\cite{SamkoITSF11}
\be
z^m[K_{\nu+m}(z)-K_{\nu-m}(z)]
=\sum_{j=0}^{m-1}(-)^{m-j-1}b_m(j)z^jK_{\nu+j}(z),
\ee
\be
z^m[I_{\nu-m}(z)-I_{\nu+m}(z)]
=\sum_{j=0}^{m-1}b_m(j)z^jI_{\nu+j}(z),
\ee
\be
z^m[J_{\nu-m}(z)-(-1)^m I_{\nu-m}(z)]
=\sum_{j=0}^{m-1}b_m(j)(-z)^jJ_{\nu+j}(z),
\ee
with $m=1,2,3,.\dots$ and
\[
b_m(j)=2^{m-j}\binom{m}{j}\nu(\nu-1)\cdots(\nu+j-m+1)
= \Gamma(\nu+1)\frac{2^{m-j}\binom{m}{j}}{\Gamma(\nu+j-m+1)}
.
\]

\cite[\S 7.10.1]{Erdelyi}
\be
f(z)=\frac{1}{z^\nu}\sum_{n=0}^\infty a_nJ_{\nu+n}(z);\quad
a_n=(\nu+n)2^{\nu+n}\sum_{s=0}^{\lfloor n/2\rfloor}
\frac{\Gamma(\nu+n-s)}{s!2^{2s}}b_{n-2s},
\ee
where $f(z)=\sum_{n=0}^\infty b_nz^n$.
Special case where only the term $s=n/2$ contributes:
\be
1 = \sum_{n\ge 0} \frac{(2n-1)!!}{(2n)!}J_n(1).
\ee

\cite[\S 7.10.1]{Erdelyi}
\be
f(z)=\frac{1}{z^\nu}\sum_{n=0}^\infty a_n z^n J_{\nu+n}(z);\quad
a_n=\sum_{s=0}^n
\frac{\Gamma(\nu+s+1)}{(n-s)!}2^{2s-n+\nu}b_s,
\ee
\be
\Gamma(\nu+n+1)b_n=\sum_{s=0}^\infty (-1)^s 2^{-\nu-n-s}\frac{a_{n-s}}{s!}.
\ee
where $f(z)=\sum_{l=0}^\infty b_lz^{2l}$.

\cite{BaileyPCPS26}
\be
[z^{2n}]\, \frac{1}{z^{2\nu}}\sum_{r=0}^\infty A_r J^2_{\nu + r}(z)
=
\frac{(-)^n\Gamma(\nu+n+1/2)}
{n!\surd\pi \Gamma(2\nu+n+1)\Gamma(\nu+n+1)}\sum_{r=0}^n\frac{A_r(-n)_r}{(2\nu+n+1)_r}
.
\ee

\cite{BaileyPCPS26}
\be
\frac{1}{z^{2\nu}}\sum_{r=0}^\infty
\frac{(2\nu)_r(\nu+1)_r(a_3)_r (a_4)_r}{r!(\nu)_r(1+2\nu-a_3)_r
(1+2\nu-a_4)_r
}J^2_{\nu+r}(z)
=
\frac{1}{4^\nu\Gamma^2(1+\nu)}\,_2F_3
\left(\begin{array}{c}\nu+1/2,1+2\nu-a_3-a_4\\ \nu+1,1+2\nu-a_3,
1+2\nu-a_4\end{array}\mid -z^2\right)
.
\ee

\cite{BaileyPCPS26}
\be
z^\nu J_\nu(2z)=\frac{2\surd \pi}{\Gamma(\nu+1/2)}
\sum_{r=0}^\infty(-)^r(\nu+r)\frac{\Gamma(2\nu+r)}{r!}J^2_{\nu+r}(z).
\ee

\cite{BaileyPCPS26}
\be
z^{-1/2} J_{2\nu+1/2}(2z)=-\frac{\Gamma(\nu+1)}{\Gamma(\nu+1/2)}
\sum_{r=0}^\infty(\nu+r)\frac{\Gamma(2\nu+r)\Gamma(r-1/2)}
{r!\Gamma(2\nu+r+3/2)}J^2_{\nu+r}(z).
\ee

\cite{BaileyPCPS26}
\be
[z^{2n}]\, \frac{1}{z^{2\nu}}\sum_{r=0}^\infty A_r J^2_{\nu +2r}(z)
=
\frac{(-)^n\Gamma(\nu+n+1/2)}
{n!\surd\pi \Gamma(2\nu+n+1)\Gamma(\nu+n+1)}\sum_{r=0}^n\frac{A_r(-n)_{2r}}{(2\nu+n+1)_{2r}}
.
\ee

\cite{BaileyPCPS26}
\be
\frac{1}{z^{2\nu}}\sum_{r=0}^\infty
\frac{(\nu)_r(1+\nu/2)_r(a_3)_r}{r!(\nu/2)_r(1+\nu-a_3)_r
}J^2_{\nu+2r}(z)
=
\frac{1}{4^\nu\Gamma^2(1+\nu)}\,_1F_2
\left(\begin{array}{c}1/2+\nu-a_3\\ \nu+1,1+2\nu-2a_3,
\end{array}\mid -z^2\right)
.
\ee

\cite{BaileyPCPS26}
\be
\frac{J_\lambda(2z\sin\frac{1}{2}\phi)}{(2z\sin\frac{1}{2}\phi)^{\lambda-2\nu}}
=
\sum_{n=0}^\infty \frac{2^{1-\lambda+2\nu}\surd \pi \Gamma(\nu+1)}{n!\Gamma(\nu+1/2)\Gamma(\lambda+1)}
\sin^{2\nu}\frac{1}{2}\phi (\nu+n)\Gamma(2\nu+n)J^2_{\nu+n}(z)
\,_3F_2\left(\begin{array}{c}\nu+1,2\nu+n,-n\\ \nu+1/2,\lambda+\end{array}\mid \sin^{2}\frac{1}{2}\phi\right)
\ee

\be
J_0(z\cos a)+J_0(z\sin a)\equiv 2\sum_{n=0}^\infty b_{2n}(a) J_{2n}(z);
\ee
where

\begin{tabular}{ll}
$2n$ & $b_{2n}$ \\
\hline
0 & 1 \\
2 & 1 \\
4 & $1+6\cos^4a-6\cos^2a = \frac{1}{4}+\frac{3}{4}\cos(4a)$ \\
6 & $1+6\cos^4a-6\cos^2a = \frac{1}{4}+\frac{3}{4}\cos(4a)$ \\
8 & $1+70\cos^8a-20\cos^2a-140\cos^6a+90\cos^4a
  = \frac{9}{64}+\frac{35}{64}\cos(8a)+\frac{5}{16}\cos(4a)$ \\
10 & $1+70\cos^8a-20\cos^2a-140\cos^6a+90\cos^4a
  = \frac{9}{64}+\frac{35}{64}\cos(8a)+\frac{5}{16}\cos(4a)$ \\
12 & $1-42\cos^2a+420\cos^4a-1680\cos^6a+3150\cos^8a-2772\cos^{10}a+924\cos^{12}a$
  \\
  & $= \frac{25}{256}+\frac{104}{512}\cos(4a)+\frac{63}{256}\cos(8a)+\frac{231}{512}\cos(12a)$ \\
\end{tabular}
Apparently, the $b_{2n}$ count numbers of Delannoy paths,
A109983 in \cite{sloane}.

\cite[p31]{MO2Afl}
\be
e^{ik\cos\varphi}=2^\nu\Gamma(\nu)\sum_{m=0}^\infty
(\nu+m)i^m J_{\nu+m}(k\rho) (k\rho)^{-\nu} C_m^{(\nu)}(\cos\varphi),
\ee
$\nu\neq 0,-1,-2,-3,\ldots$.

\cite{FikiorisMathComp67}
\be
S_{\mu,\nu}(z)\sim z^{\mu-1}
\sum_{m=0}^\infty (-1)^m\left(\frac{1-\mu+\nu}{2}\right)_m
\left(\frac{1-\mu-\nu}{2}\right)_m
\left(z/2\right)^{-2m},
\,
|z|\to \infty, |\arg z|<\pi.
\ee
The series terminates and is equal to $S_{\mu,\nu}(z)$ when $\mu\pm\nu$
is a positive odd integer.

\cite{RangarajanQJM15}
\be
\left(\frac{\sin\beta}{\sin\alpha}\right)^{m+n}
P_{m+n}^m(\cos\alpha)
=
\sum_{r=0}^n\binom{2m+n}{r}\left(\frac{\sin(\beta-\alpha)}{\sin\alpha}\right)^{r}
P_{n+m-r}^m(\cos\beta)
.
\ee

\subsection{Orthogonal Polynomials}
\cite{WatsonJLMS13}
\be
L_m^{(\alpha)}
L_n^{(\alpha)}
=\sum_{M=|m-n|}^{m+n}
\frac{L_M^{(\alpha)}(z)(-)^{m+n-M}2^{m+n-M}M!}
{(M-m)!(M-n)!(m+n-M)!}
\,_3F_2\left(\begin{array}{c}
\alpha+M+1,\frac{M-m-n}{2},\frac{(M-m-n+1}{2}\\
M-m+1,M-n+1
\end{array}\mid 1\right).
\ee

\cite{RainaITSF3}
\be
\frac{1}{(1+rx)^\lambda}
=\sum_{n\ge 0}\frac{(\lambda)_n}{n!}x^n
\sum_{k=0}^n (-)^k\binom{n}{k}L_n^{(rk-n)}(y).
\ee

\cite{CohenMC31}
\be
\sum_{n\ge 0} t^n L_n^{a+bn}(x)=\frac{(1+v)^{a+1}}{1-bv}\exp(-xv),
\ee
where $v=t(1+v)^{b+1}$, $v(0)=0$.

\cite{CohenMC31}
\be
\sum_{n\ge 0} t^n L_n^{v+bn}(x(1+an))=
\frac{(1-z)^{1-v}}
{1-z(b+2-ax)+z^2(b+1)}
e^{xz/(z-1)},
\ee
where $t=z(1-z)^b\exp[axz/(1-z)]$ and $|t|<1$.

\cite{CohenMC31}
\be
\sum_{n\ge 0} \frac{t^n}{v+bn+n} L_n^{v+bn}(x(1+an))=
\frac{\exp[xz/(z-1)]}{v(1-z)^v}
\,_1F_1\left(\begin{array}{c}
1\\ \frac{v+1+b}{1+b}
\end{array}\mid \frac{xz(1+b-av)}{(1-z)(1+b)}\right)
,
\ee
where $t=z(1-z)^b\exp[axz/(1-z)]$ and $|t|<1$.

\cite{CohenMC31}
\be
\sum_{n\ge 0} \frac{t^n}{1+an} L_n^{v+avn-n}(x(1+an))=
\frac{1}{(1-z)^v}
e^{xz/(z-1)}
,
\ee
where $t=z(1-z)^{av-1}\exp[axz/(1-z)]$ and $|t|<1$.

\cite{CohenMC31}
\be
\sum_{n\ge 0} \frac{t^n}{1+an} L_n^{v+bn}(x(1+an))=
\frac{\exp[xz/(z-1)]}{(1-z)^v}
\,_2F_1\left(\begin{array}{c}
1, \frac{1+b-av}{a}\\ \frac{a+1}{a}
\end{array}\mid z\right)
,
\ee
where $t=z(1-z)^b\exp[axz/(1-z)]$ and $|t|<1$.

\cite{CohenMC31}
\begin{multline}
\sum_{n\ge 0} \frac{t^n}{1+n} L_n^{v+bn}(x(1+an))=
\frac{\exp[x(1-a-z)/(1-z)]}{(1-z)^v z(v-b)}
\big\{
\,_1F_1\left(\begin{array}{c}
v-b\\ v-b+1
\end{array}\mid \frac{x(a-1)}{1-z}\right)
\\
-(1-z)^{v-b}
\,_1F_1\left(\begin{array}{c}
v-b\\ v-b+1
\end{array}\mid x(a-1)\right)
\big\}
,
\end{multline}
where $t=z(1-z)^b\exp[axz/(1-z)]$ and $|t|<1$.

\cite{CohenMC31}
\be
\sum_{n\ge 0}\frac{(1+bn)^{n/2}}{n!}t^n H_n[x(1+an)/(1+bn)^{1/2}]
=
\frac{e^{-z^2-2xz}}{1+2bz^2+2axz},
\ee
where $t=(-z)e^{bz^2+2axz}$ and $|2azx\exp[bz^2+2axz+1]|<1$.

\cite{CohenMC31}
\be
\sum_{n\ge 0}\frac{(1+bn)^{n/2-1}}{n!}t^n H_n[x(1+an)/(1+bn)^{1/2}]
=
e^{-z^2-2xz}
\,_1F_1\left(\begin{array}{c}
1\\ 1+1/b
\end{array}\mid 2xz(b-a)/b\right)
\,
\ee
where $t=(-z)e^{bz^2+2axz}$ and $|2azx\exp[bz^2+2axz+1]|<1$.

\cite{CohenMC31}
\be
\sum_{n\ge 0}\frac{(1+an)^{n/2-1}}{n!}t^n H_n[x(1+an)]
=
e^{-z^2-2xz}
\,
\ee
where $t=(-z)e^{az^2+2axz}$ and $|2azx\exp[az^2+2axz+1]|<1$.

\cite{DixonJPA6}
\be
\frac{d^s}{dx^s} P_n(x) = \sum_{p=0}^{\lfloor (n-s)/2\rfloor} \binom{s+p-1}{p} (2n-2s-4p+1)\prod_{k=0}^{s-2}
(2n-1-2p-2k)P_{n-s-2p}(x).
\ee

\cite{BaileyPCPS27_2}
\be
\sin^p tP_n^{-m}(\cos t)
=2^p \sum_{r=0}^\infty \frac{(-)^r(m+n+1)_{2r}\Gamma(p+m+r)}
{r!\Gamma(m+r+1)}
(p+m+2r)P_{n-p}^{-p-m-2r}(\cos t).
\ee
\cite{BaileyPCPS27_2}
\be
\sin^m t
=\sum_{r=0}^\infty \frac{(-)^r 2^{m+2r}((m+n+1)/2)_r\Gamma(m+r)}
{r!}
(m+2r)P_n^{-m-2r}(\cos t).
\ee
\cite{BaileyPCPS27_2}
\be
2m \frac{P_n^{-m}(\cos\theta)}{\sin \theta}
=
P_{n+1}^{-m_1}(\cos\theta)+(m+n+1)(m+n+2)P_{n+1}^{-m-1}(\cos \theta).
\ee
\cite{BaileyPCPS27_2}
\be
(\frac12 \sin\theta)^{\mu-\nu} P_\kappa^{-\nu}(\cos \theta)
=
\Gamma(\mu+1)\sum_{r=0}^\infty \frac{(\nu-\mu)_r(\kappa+\nu+1)_{2r}}
{r!\Gamma(\nu+r+1)}(\frac12 \sin\theta)^r P_{\kappa+\nu-\mu+r}^{-\mu -r}(\cos \theta).
\ee

\cite{DixonJPA6}
\begin{multline}
|{\mathbf r}_1-{\mathbf r}_2|^l
Y_l^m(\theta_R,\phi_R)
=\sqrt{\frac{4\pi(l+m)!(2l+1)!}{(2l)!(l-m)!}}
\sum_{f=0}^l\sum_{w=0}^{l-m}
\sqrt{\frac{(2f)!w!(2l-2f)!(l-m-w)!}{(2f+1)!(2l-2f+1)!(2f-w)!(l-2f+m+w)!}}
\\
\times
(-1)^f\binom{l-m}{w}r_1^f Y_f^{f-w}(\theta_1,\phi_1)r_2^{l-f}Y_{l-f}^{m+w-f}(\theta_2,\phi_2)
\end{multline}
$f$ and $w$ are integers and an particular term is considered zero if at least one of the relations $|f-w|\le f$
and $|m+w-f|\le l-f$ is violated. $\theta_R$ and $\phi_R$ describe the orientation fo the difference vector
$\mathbf r_1-\mathbf r_2$.

\cite{FoxJLMS13}
\be
\int_0^1 f(t)t^rdt =c_r, \,(r=0,1,\ldots n)
\leadsto 
f(t)=\sum_{i=0}^n (2i+1)\sum_{r=0}^i
\left\{(-)^r\frac{(i+r)!}{(i-r)!}\,\frac{1}{(r!)^2}c_r \right\}P_i(1-2t).
\ee

\section{Special Functions II.}
\subsection{Hypergeometric Functions}
\cite{BaileyPLMS28}
\be
e^x {}_0F_1(1/2;-x^2/4)=\sum_{n=0}^\infty 2^{n/2}\cos\frac{n\pi}{4} \frac{x^n}{n!}.
\ee
\cite{BaileyPLMS28}
\be
e^x {}_0F_1(3/2;-x^2/4)=\sum_{n=1}^\infty 2^{n/2}\sin\frac{n\pi}{4} \frac{x^{n-1}}{n!}.
\ee

\cite{RoachSSAC96}
\be
_1F_1(1/2;9/2;z)=
-\frac{525+280z+140z^2}{128z^3}e^z+
\frac{525+630z+420z^2+280z^3}{256z^{7/2}}
\sqrt \pi \erfi(\sqrt z).
\ee

\cite{SlaterPCPS50}
\be
_1F_1(\frac{1}{2}+\frac{1}{2}a-b,1+a-b;x)
=
e^{x/2}
\sum_{r=0}^\infty
\frac{(a)_r(b)_r (-x/4)^r}{r!(a/2)_r(1+a-b)_r} \,_0F_1(;\frac{1}{2}a+r+1;(x/4)^2) .
\ee
\cite{BaileyPLMS28}
\be
e^{-x/2}{}_1F_1(a;2a;x) = {}_0F_1(a+1/2;x^2/16).
\ee

\cite{AncaraniJMP49}
\begin{multline}
\frac{d}{da^n} {}_1F_1(a;b;z)=
\frac{z^n}{(b)_n} \sum_{m_1=0}^\infty \cdots \sum_{m_{n+1}=0}^\infty
\frac{(1)_{m_1}(1)_{m_2}\cdots (1)_{m_{n+1}}}{(n+1)_{m_1+m_2+\cdots m_{n+1}}
(b+n)_{m_1+m_2+\cdots m_{n+1}}}
\\
\times
\frac{(a)_{m_1}(a+1)_{m_2+m_2}\cdots (a+n)_{m_1+\cdots m_{n+1}}}
{(a+1)_{m_1}\cdots(a+n)_{m_1+\cdots+m_n}}
\frac{z^{m_1+m_2+\cdots +m_{n+1}}
}{m_1!m_2!\cdots m_{n+1}!}
.
\end{multline}
\cite{BaileyPLMS28}
\be
e^x {}_0F_2(1/3,2/3;-x^3/27)=1+2\sum_{n=1}^\infty 3^{n/2-1}\cos\frac{n\pi}{6} \frac{x^n}{n!}.
\ee
\cite{BaileyPLMS28}
\be
e^x {}_0F_2(2/3,4/3;-x^3/27)=2\sum_{n=0}^\infty 3^{n/2-1/2}\sin\frac{(n+2)\pi}{6} \frac{x^n}{(n+1)!}.
\ee
\cite{BaileyPLMS28}
\be
e^x {}_0F_2(4/3,5/3;-x^3/27)=4\sum_{n=0}^\infty 3^{n/2}\sin\frac{(n+1)\pi}{6} \frac{x^n}{(n+2)!}.
\ee
\cite{BaileyPLMS28}
\be
e^x {}_0F_{\lambda-1}(\frac{1}{\lambda},\frac{2}{\lambda},\ldots \frac{\lambda-1}{\lambda};-(x/\lambda)^\lambda)=
1+\frac{1}{\lambda}\sum_{n=1}^\infty[\sum_{r=1}^\lambda
\cos^n\frac{(2r-1)\pi}{2\lambda} \cos\frac{n(2r-1)\pi}{2\lambda}]\frac{2^nx^n}{n!}.
\ee

\cite{RoachSSAC96}
\begin{multline*}
_1F_2(3/2;5/2,5;z)
=
-\frac{432-24z+96z^2}{5z^3}I_0(2\sqrt z)
+\frac{432+192z+48z^2}{5z^{7/2}}I_1(2\sqrt z)
\\
-\frac{48}{5z}\pi\left(I_0(2\sqrt z)L_1(2\sqrt z)-I_1(2\sqrt z)
L_0(2\sqrt z)\right)
.
\end{multline*}

\cite{RainaITSF3}
\be
\frac{1}{(1-cy)^\alpha (1+rz)^\lambda}
=
\sum_{n\ge 0}\frac{(\lambda)_n}{n!}x^n
\sum_{k=0}^n (-)^k\binom{n}{k}\binom{rk}{n}{}_2F_1(\alpha,1+rk; 1+rk-n; cy)
\ee
where $|cy|<1,|rx|<1$.

\cite{KimarXiv1205}
\be
_2F_1(m+1/2,1/2-m;3/2+m;1) = \frac{\pi}{2^{2m+1}}\frac{(3/2)_m}{m!}.
\ee

\cite{KrupnikovJCAM78}
\be
_2F_1(ai,-ai;1/2; z) = \cosh(2a\arcsin\surd z).
\ee
\cite{KrupnikovJCAM78}
\be
_2F_1(\frac12+ai,\frac12-ai;3/2; z) = \frac{\sinh(2a\arcsin\surd z}{2a\surd z}.
\ee
\cite{KrupnikovJCAM78}
\be
_2F_1(ai,-ai;1/2; -z) = \cos(2a\arcsinh\surd z).
\ee
\cite{KrupnikovJCAM78}
\be
_2F_1(\frac12+ai,\frac12-ai;3/2; -z) = \frac{\sin(2a\arcsinh\surd z}{2a\surd z}.
\ee
\cite{KrupnikovJCAM78}
\be
_2F_1(\frac12,\frac12;3/2; -1) = \ln(1+\surd 2).
\ee
\cite{KrupnikovJCAM78}
\be
_2F_1(\frac12,1;3/2; 3/4) = \frac{2}{\surd 3}\ln(1+\surd 3).
\ee

\cite{ZuckerMPCPS131}
\be
_2F_1(1/6,1/3 ; 1/2; 25/27) = 3^{3/2}/4.
\ee
\cite{ZuckerMPCPS131}
\be
_2F_1(1/4,1/2 ; 3/4; 80/81) = 9/5 .
\ee
\cite{ZuckerMPCPS131}
\be
_2F_1(1/10,3/10; 9/20;1) = \frac{5^{1/4}}{2}(\sqrt 5-1).
\ee
\cite{ZuckerMPCPS131}
\be
_2F_1(1/6,1/4; 11/24;1) = \sqrt{3}(\sqrt 3+2)^{1/2}.
\ee
\cite{ZuckerMPCPS131}
\be
_2F_1(1/7,2/7; 1/2;1) = \frac{1}{2\sin(\pi/14)}.
\ee

\cite{BeckenJCAM126}
\be
_2F_1(a,b;c;z) = (a)_m (b)_m z^m \,_2F_1(a+m,b+m;1+m;z).
\ee

\cite{BeckenJCAM126}
\be
_2F_1(a,b;c;z) = (1-z)^{c-a-b} \,_2F_1(c-a,c-b;c;z).
\ee

\cite{BeckenJCAM126}
\be
_2F_1(a,b;c;z) = (a)_m(b)_m z^m(1-z)^{1-a-b-m}
 \,_2F_1(1-b,1-a;1+m;z).
\ee

\cite{BeckenJCAM126}
\be
_2F_1(a,b;c;z) = (1-z)^{-a}
 \,_2F_1(a,c-b;c;z/(z-1)).
\ee
\cite{BaileyPLMS28}
\be
(1-x)^{-2a} {}_2F_1(a,a+1/2;\rho+1/2;\frac{x^2}{(1-x)^2})
={}_2F_1(2a,\rho; 2\rho;2x).
\ee
\cite{BaileyPLMS28}
\be
(1-x)^{-2a} {}_2F_1(a,a+1/2;1/2;-\frac{x^2}{(1-x)^2})
=\sum_{n=0}^\infty \cos\frac{n\pi}{4}\frac{(2a)_n 2^{n/2}}{n!}x^n.
\ee
\cite{BaileyPLMS28}
\be
(1-x)^{-2a} {}_2F_1(a,a+1/2;3/2;-\frac{x^2}{(1-x)^2})
=\sum_{n=0}^\infty \sin\frac{(n+1)\pi}{4}\frac{(2a)_n 2^{n/2+1/2}}{(n+1)!}x^n.
\ee
\cite{BaileyPLMS28}
\be
(1-x)^{-2a} {}_2F_1(a,a+1/2;2a+1;\frac{4x^3}{27(1-x)^3})
={}_2F_1(3a,3a+1;6a+1;4x/3).
\ee
\cite{BaileyPLMS28}
\be
(1-x)^{-a} {}_2F_1(a,\beta;a+\beta+1/2;-\frac{x^2}{4(1-x)})
={}_2F_1(2a,a+\beta;2a+2\beta;x).
\ee
\cite{BaileyPLMS28}
\be
(1-x)^{a+\beta-1/2} {}_2F_1(a,\beta;a+\beta+1/2;4x(1-x))
={}_2F_1(a-\beta+1/2,\beta-a+1/2;a+\beta+1/2;x).
\ee

\cite{NagelJMP42}
\be
\frac{(z/2)^c}{\Gamma(1+c)} \,_2F_1(a,b;c+1;-z^2/(4ab))
=\sum_{\nu\ge 0} \frac{1}{\nu!}\, _3F_0(-\nu,a,b;1/(ab)) (z/2)^\nu J_{c+\nu}(z).
\ee

\cite{BeckenJCAM126}
\begin{multline}
_2F_1(a,b;c;z) = (1-z)^{-a}
\left(\frac{z}{z-1}\right)^{1-c}
\bigg\{
\frac{\Gamma(|m|)}{\Gamma(a+\bar m)\Gamma(c-a-\underline m)}
\left(\frac{1}{1-z}\right)^{\underline m}
\\
\sum_{n=0}^{|m|-1}
\frac{(1-a-\bar m)_n(1-c+a+\underline m)_n}
{(1-|m|)_n\Gamma(1+n)}\left(\frac{1}{1-z}\right)^n
-\frac{(-)^m}{\Gamma(a+\underline m)\Gamma(c-a-\bar m)}
\left(\frac{1}{1-z}\right)^{\bar m}
\\
\sum_{n=0}^\infty
\frac{(1-a-\underline m)_n(1-c+a+\bar m)_n}
{\Gamma(1+|m|+n)\Gamma(1+n)}\left(\frac{1}{1-z}\right)^n
\\
[-\ln(1-z)-\pi\cot \pi(c-a)-\pi \cot \pi a+\psi(1-a-\underline m+n)
+\psi(1-c+a+\bar m+n)-
\psi(1+|m|+n)-\psi(1+n)]
\bigg\}
,
\end{multline}
where $\bar m \equiv \max(0,m)$, $\underline m\equiv \min (0,m)$.

\cite[1.7.7]{SlaterHyp}
\be
_2F_1(a,-m;c,1)=\frac{(c-a)_m}{(c)_m}
.
\ee

\cite{GottschalkJPA21}
\be
_2F_1(1,a;a+1;z)=-az^{-a}\ln(1-z)-az^{-a}\sum_{i=0}^{a-2}\frac{z^{i+1}}{i+1}; \quad a\in \mathbb{Z}^+, |z|<1\, \mathrm{or }\,z=-1.
\ee

\cite{GottschalkJPA21}
\be
_2F_1(1,a;a+1;-z^2)=2az^{-2a}(-)^{a-1/2} \tan^{-1}z + 2a(-z^2)^{1/2-a}{\sum}_{i=a-1/2}^{'-1}\frac{(-z^2)^i}{2i+1},\quad 2a\in \mathbb{Z}, a\not\in \mathbb{Z}, |z|\le 1
\ee

\cite{GottschalkJPA21}
\be
_2F_1(1,a;a+1;z^2)=az^{-2a}\ln\frac{1+z}{1-z}+ 2az^{-2a}{\sum}_{i=a-1/2}^{'-1}\frac{z^{2i+1}}{2i+1},\quad 2a\in \mathbb{Z}, a\not\in \mathbb{Z}, |z|< 1
\ee
and the primes at the sum symbol means ${\sum'}_{r=p}^{q-1}$ is $\sum_{r=p}^{q-1}$ for $p<q$, is $-\sum_{r=q}^{p-1}$ for $p>q$
and zero for $p=q$.

\cite{GottschalkJPA21}
\be
_2F_1(a,b;a+1;z)= -\frac{\Gamma(a+1)\Gamma(1-b)}{\Gamma(a+1-b)z^a}
[(1-z)^{1-b}-1]
-\frac{\Gamma(a+1)(1-z)^{1-b}}{\Gamma(a+1-b)z^a}
\sum_{i=1}^{a-1}\frac{\Gamma(i+1-b)z^i}{\Gamma(i+1)},a\in \mathbb{Z}^+
\ee

\cite{GottschalkJPA21}
\begin{multline}
_2F_1(a,b;a+1;z)= -\frac{\Gamma(a+1)(1-z)^{1-b}}{\Gamma(a+1-b)z^{a+1/2}}
{\sum}_{i=1}^{'a-1/2}
\frac{\Gamma(i+1/2-b)z^i}{\Gamma(i+1/2)}
+\frac{\Gamma(a+1)\Gamma(3/2-b)\Gamma(b-1/2)}{\Gamma(1+a+b)\Gamma(b)z^{a-1/2}\surd \pi}
\\
\left(
\sum_{i=1}^{b-1}\frac{\Gamma(i)}{(1-z)^i\Gamma(i+1/2)}+\frac{1}{\sqrt{|z|}\omega(z)}
\right)
\end{multline}
where $2a\in \mathbb{Z}$, $a\not \in \mathbb{Z}$, $b\in \mathbb{Z}^+$ and 
\begin{equation}
\omega(z)=\left\{
\begin{array}{rl}
2\tan^{-1}\sqrt{-z},& z<0\\
\ln\frac{1+\surd z}{1-\surd z},& z>0
\end{array}
\right.
\end{equation}
and similar reductions for $2a$ and/or $2b$ in $\mathbb{Z}$\ldots

A special case of \eqref{eq.hypRev} \cite{ChaundyQJM14}
\be
_2F_1(a,-n;c;p) = \frac{(a)_n}{(c)_n}(-p)^n
\,_2F_1(1-c-n,-n;1-a-n;-1/p).
\ee

\cite{VidunasRMJ32}
\be
_2F_1(a+n,b;a-b;-1) = P(n)\frac{\Gamma(a-b)\Gamma(\frac{a+1}{2})}{\Gamma(a)\Gamma(\frac{a+1}{2}-b)}
+ Q(n)\frac{\Gamma(a-b)\Gamma(\frac{a}{2})}{\Gamma(a)\Gamma(\frac{a}{2}-b)}
\ee
where
\begin{eqnarray*}
P(n)=\frac{1}{2^{n+1}}\,_3F_2(-n/2,-(n+1)/2,a/2-b; 1/2,a/2; 1);\quad
\\
Q(n)=\frac{n+1}{2^{n+1}}\,_3F_2(-(n-1)/2,-n/2,(a+1)/2-b; 3/2,(a+1)/2; 1).
\end{eqnarray*}

\cite{MoritaIIS2}
Let $N_n(z)\equiv{}_2F_1(\frac12+n,\frac12+n;1+2n;z)$, then
\begin{equation}
N_{n-1}(z)-(1-z/2)N_n(z)+\frac{(n+1/2)^2}{16n(n+1)}{}z^2N_{n+1}(z)=0.
\end{equation}

\cite{MoritaIIS2}
Let $M_m(z)=F(a,b;c+m;z)$ then
\begin{multline}
-(c+m)(c+m-1)(1-z)M_{m-1}(z)
+(c+m)[c+m-1-(2c+2m-a-b-1)z]M_m(z)
\\
+(c+m-a)(c+m-b)zM_{m+1}(z)=0.
\end{multline}

\cite{ApagoduInt8}
\be
_2F_1(-2n,b;-2n+2r-b;-1) =
\frac{(1/2)_n(b+1-r)_n}{(b/2+1-r)_n(b/2+1/2-r)_n}
\sum_{i=0}^{r-1}\frac{2^{2i}i!\binom{r+i-1}{2i}}{(b-r+1)_i}\binom{n}{i}
.
\ee
\cite{BaileyPLMS28}
\be
_2F_1(a,b;1+a-b;-1) = \frac{\Gamma(1+a-b)\Gamma(1/2)}{2^a\Gamma((1+a)/2)\Gamma(a/2-b+1)}.
\ee

\cite{MaierTAAMS358}
\be
_2F_1(a, b;1+a-b;-1) =
\Gamma\left[
\begin{array}{c}
1+2a-2b, 1+a-b\\
1+a-2b,1+2a-b
\end{array}
\right]
\Gamma\left[
\begin{array}{c}
1+\frac12 a, \frac12+\frac12 a-b,1+a-\frac12 b,\frac12+a-\frac12 b\\
1+a,\frac12+a-b,1+\frac12 a-\frac12 b,\frac12 +\frac12 a-\frac12 b
\end{array}
\right]
.
\ee

\cite{MaierTAAMS358}
\be
_2F_1(a-4,\frac23 a-1;1+\frac13 a;-1) =
\Gamma\left[
\begin{array}{c}
\frac12, \frac13 a, 1+\frac13 a, \frac23 a-\frac32, \frac23 a-2\\
\frac12 a-\frac32, \frac16 a, \frac12 +\frac16 a, 2-\frac16 a,\frac16 a,\frac43 a-3
\end{array}
\right]
,
\ee
where the parametric excess of the $_2F_1(-1)$ has real part greater than $-1$.

\cite{MaierTAAMS358}
\be
_2F_1(a-\frac 32,\frac14 a+\frac14;\frac34 a-\frac14 ;-1) =
\Gamma\left[
\begin{array}{c}
\frac12, \frac32, \frac34+\frac14 a, \frac12 a, \frac34 a-\frac14\\
\frac78 , \frac98, a, \frac18 +\frac14 a, \frac38 +\frac14 a
\end{array}
\right]
,
\ee
where the parametric excess of the $_2F_1(-1)$ has real part greater than $-1$.

\cite{MaierTAAMS358}
\be
_2F_1(a-\frac 12,\frac14 a;\frac34 a;-1) =
\Gamma\left[
\begin{array}{c}
\frac12, \frac12, \frac12+\frac14 a, \frac12+\frac12 a, \frac34 a\\
\frac38 , \frac58, a, \frac38 +\frac14 a, \frac58 +\frac14 a
\end{array}
\right]
,
\ee
where the parametric excess of the $_2F_1(-1)$ has real part greater than $-1$.

\cite{MaierTAAMS358}
\be
_2F_1(a+2,\frac23 a;\frac13 a;-1) =
\Gamma\left[
\begin{array}{c}
\frac12, \frac13 a, 1+\frac13 a, \frac32+\frac23 a, 2+\frac23 a\\
\frac32+\frac12 a , \frac12 +\frac16 a, 1+\frac16 a, -\frac16 a, 2 +\frac43 a
\end{array}
\right]
,
\ee
where the parametric excess of the $_2F_1(-1)$ has real part greater than $-1$.

\cite{VidunasRMJ32}
\begin{multline}
_2F_1(-a,1/2;2a+3/2+n;1/4) = 
\frac{2^{n+3/2}}{3^{n+1}}
\frac{\Gamma(a+5/4+n/2)\Gamma(a+3/4+n/2)\Gamma(a+1/2)}
{\Gamma(a+7/6+n/3)\Gamma(a+5/6+n/3)\Gamma(a+1/2+n/3)}
K(n)
\\
-(-3)^{n-2}
2^{3/2}
\frac{\Gamma(a+5/4+n/2)\Gamma(a+3/4+n/2)\Gamma(a+1)}
{\Gamma(a+3/2)\Gamma(a+1/2+n/2)\Gamma(a+1+n/2)}
L(n)
\end{multline}
where $K(n)$ and $L(n)$ are defined in the reference.

\cite{BaileyPLMS28}\cite{ChoEAMJ15}
\be
_2F_1(a,b;(a+b+1)/2; 1/2) = \frac{\Gamma(1/2)\Gamma( (1+a+b)/2)}{\Gamma((1+a)/2)\Gamma((1+b)/2)}.
\ee

\cite{ChoEAMJ15}
\be
_2F_1(a,b;(a+b)/2; 1/2) = \Gamma(1/2)\Gamma( (a+b)/2)
[
\frac{1}{\Gamma(a/2)\Gamma((1+b)/2)}
+\frac{1}{\Gamma(b/2)\Gamma((1+a)/2)}
]
.
\ee

\cite{ChoEAMJ15}
\be
_2F_1(a,b;1+(a+b)/2; 1/2) = \frac{a+b}{a-b}\Gamma(1/2)\Gamma( (a+b)/2)
[
\frac{1}{\Gamma(a/2)\Gamma((1+b)/2)}
-\frac{1}{\Gamma(b/2)\Gamma((1+a)/2)}
]
.
\ee

\cite{LavoieMC49}
\be
_2F_1(a,-a;e;1/2)
=\frac{\Gamma(e)}{2^{a+1}\Gamma(e-a)}\big\{\frac{\Gamma(\frac{e-a}{2})}{\Gamma(\frac{e+a}{2})}+\frac{\Gamma(\frac{e-a+1}{2})}{\Gamma(\frac{e+a+1}{2})}\big\}.
\ee

\cite{LavoieMC49,ChoiCKM14}
\be
_2F_1(a,2-a;e;1/2)
=\frac{\Gamma(e)}{2^{a-1}(a-1)\Gamma(e-a)}\big\{\frac{\Gamma(\frac{e-a}{2})}{\Gamma(\frac{e+a-2}{2})}+\frac{\Gamma(\frac{e-a+1}{2})}{\Gamma(\frac{e+a+1}{2})}\big\}.
\ee

\cite{LavoieMC49}
\be
_2F_1(a,1-a;e;1/2)
=\frac{\Gamma(e/2)\Gamma(\frac{1+e}{2})}{\Gamma(\frac{e+a}{2})\Gamma(\frac{1+e-a}{2})}.
\ee

\cite{ChoiCKM14}
\begin{multline}
_2F_1(1+a,2-a;b;1/2)
=
\frac{2(a-b+1)(a+b-2)}{a(a-1)}\frac{\Gamma(b/2)\Gamma((b+1)/2)}
{\Gamma((a+b)/2)\Gamma((b-a+1)/2)}
\\
+
\frac{4(b-2)}{a(a-1)}
\frac{\Gamma(b/2)\Gamma((b+1)/2)}{\Gamma((b-a)/2)\Gamma((b+a-1)/2)}
\end{multline}
and similarly contiguously extended relations. For $_2F_1(a,1+n-a;b;1/2)$
and $|n|\le 5$ see the table in \cite{LavoieJCAM72}, and the generic formula for other offsets $n$.

\cite{BerndtBLMS15}
\be
\pi{}_2F_1(1/2,1/2;1;1-x)=\log\frac{16}{x}{}
_2F_1(1/2,1/2;1;x)-4\sum_{k\ge 1}\frac{(1/2)_k^2}{(k!)^2}
\sum_{j=1}^k\frac{x^k}{(2j-1)(2j)}.
\ee

\cite{BerndtBLMS15}
\be
_2F_1(\frac{1}{2},\frac{1}{2};1;\frac{1+x}{2})
=
\frac{\sqrt{\pi}}{\Gamma^2(3/4)}{}_2F_1(1/4,1/4;1/2;x^2)
+
\frac{\Gamma^2(3/4)}{\pi^{3/2}}{}_2F_1(3/4,3/4;3/2;x^2)
.
\ee

\cite{ZuckerMPCPS131}
\be
_2F_1(a,a;2a;z) = (\frac12+\frac12\sqrt{1-z})^{-2a}
{}_2F_1[a,1/2; a+1/2; (\frac{1-\sqrt{1-z}}{1+\sqrt{1-z}})^2].
\ee

\cite{BerndtBLMS15}
\begin{multline}
_2F_1(\frac{1}{2},\frac{1}{2};1;\frac{1}{2}+\frac{x}{1+x^2})
=
\frac{\sqrt{\pi}}{\Gamma^2(3/4)}\sqrt{1+x^2} {}_2F_1(1/4,1/2;3/4;x^4)
\\
+
\frac{\Gamma^2(3/4)}{\pi^{3/2}}x(1+x^2)^{3/2} {}_2F_1(1/2,3/4;5/4;x^4)
.
\end{multline}

\cite{BerndtBLMS15}
\be
_2F_1(n,-n;1/2;x^2) = \cos(2n\sin^{-1}x).
\ee

\cite{BerndtBLMS15}
\be
2nx {} _2F_1(\frac{1}{2}+n,\frac{1}{2}-n;3/2;x^2) = \sin(2n\sin^{-1}x).
\ee

\cite{BerndtBLMS15}
\be
 _2F_1(\frac{1}{2}+n,\frac{1}{2}-n;1/2;x^2) = (1-x^2)^{-1/2} \cos(2n\sin^{-1}x).
\ee

\cite[15.1.6]{AS}
\be
 _2F_1(\frac{1}{2},\frac{1}{2};\frac{3}{2};x^2) = \frac{1}{x}\arcsin x.
\ee
\cite{ZuckerMPCPS131}
\be
\frac{2 \Gamma(\frac12)\Gamma(a+b+1/2)}{\Gamma(a+1/2)\Gamma(b+1/2)} {}_2F_1(a,b;1/2;z)
={}_2F_1(2a,2b;a+b+1/2;\frac{1-\sqrt{z}}{2})
+{}_2F_1(2a,2b;a+b+1/2;\frac{1+\sqrt{z}}{2}).
\ee
\cite{ZuckerMPCPS131}
\be
\frac{\Gamma^2(1/4)}{2\sqrt \pi} {}_2F_1(1/4,1/4;1/2;z)
=K(k)+K'(k), \quad z=(1-2k^2)^2.
\ee
\cite{ZuckerMPCPS131}
\be
{}_2F_1(1/8,3/8;1/2;z)
=\frac{\sqrt{\surd 2+1}\sqrt{1+k^2}}{2^{3/4}K(k_2)}[\sqrt{2}K(k)+K'(k)],\quad k_2=\sqrt{2}-1,
\ee
where $z=(1-6k^2+k^4)^2/(1+k^2)^4$, $K'(k)=K(k'), k'=\sqrt{1-k^2}$.
\cite{ZuckerMPCPS131}
\be
{}_2F_1(1/6,1/3;1/2;z)
=\frac{3^{-3/4}(1+\alpha+\alpha^2)}{K(k_3)(1+2\alpha)^{3/2}}
[(1+2\alpha)K(k)+\sqrt{3}K'(l)]
\ee
where $z=(2+2\alpha-\alpha^2)^2(1+4\alpha+\alpha^2)^2(1-2\alpha-2\alpha^2)^2/4/(1+\alpha+\alpha^2)^6$,
$k^2=\alpha^2(2+\alpha)/(1+2\alpha)$, $l^2=\alpha(2+\alpha)^3/(1+2\alpha)^3$, $k_3=\frac14\sqrt{2}(\surd 3-1)$.
\cite{ZuckerMPCPS131}
\be
{}_2F_1(3/4,3/4 ; 3/2; 3/4) = \frac{2^{2/3}3^{1/4}}{12\pi}(\surd 3-1)\frac{\Gamma^2(1/4)\Gamma^2(1/3)}{\Gamma(1/2)\Gamma(2/3)}
\ee
\cite{ZuckerMPCPS131}
\be
{}_2F_1(1/4,1/2 ; 3/4; 3/4) = (4/3)^{3/4}.
\ee
\cite{ZuckerMPCPS131}
\be
{}_2F_1(3/4,1/2 ; 5/4; 3/4) = \frac{\Gamma^4(1/4)}{12\pi\Gamma^2(1/2)}.
\ee

Via \cite[15.3.1]{AS}, the substitution 
$t=z^3$  gives
\be
 _2F_1(\frac{2}{3},1;\frac{5}{3};x) 
=
\frac{2}{3}\int_0^1 t^{-1/3}(1-xt )^{-1} dt
=
2\int_0^1 \frac{z}{1-xz^3} dz,
\ee
and this is evaluated via \cite[2.126]{GR}
\be
\int \frac{z}{a+bz^3} dz 
= -\frac{1}{3b\alpha}\left[\frac12 \ln\frac{(z+\alpha)^2}{z^2-\alpha z+\alpha^2}
-\surd 3 \arctan
\frac{2z-\alpha}{\surd 3 \alpha}\right],
\ee
where $\alpha\equiv \sqrt[3](a/b)$.

\be
 _2F_1(-\frac{1}{2},-m-\frac{1}{2};\frac{1}{2};z) = 
\sum_{s=0}^m \beta_{s,m}(1-z)^{s+1/2}+\frac{(2m+1)!!}{(2m)!!}\sqrt z \arcsin \surd z
\ee
where
\be
\beta_{s,m}=\left\{ \begin{array}{ll} \frac{(2m+1)!!}{(2m)!!}, & s=0 ;\\ 
-\frac{1}{2s}\frac{(s+3/2)(s+5/2)(s+7/2)\cdots (m+1/2)}{(s+1)(s+2)(s+3)\cdots m}, & 0<s\le m.
\end{array} \right.
\ee
\cite{RoachSSAC96}
\be
 _2F_1(-\frac{1}{2},\frac{1}{2};\frac{3}{2};x^2) = \frac{1}{2}[\frac{\arcsin x}{x}+\sqrt{1-x^2}].
\label{eq.roach1}
\ee

Differentiation of \eqref{eq.roach1} w.r.t. $x$ using \cite[15.2.2]{AS}
\be
{}_2F_1(\frac{1}{2},\frac{3}{2};\frac{5}{2};x^2) = \frac{3}{2x^2}[
\frac{\arcsin x}{x}
-\sqrt{1-x^2}
].
\ee

\cite{DetrichMC33}
\be
_2F_1(l+1,m+1/2;n+3/2;z) = A_{l,m,n}(z)\log\frac{1+z^{1/2}}{1-z^{1/2}}+B_{l,m,n}(z),\quad l\ge 0, n-m\ge 0,
\ee
where
\be
A_{l,m,n}(z)=\frac{(-)^{n-mn-l}\Gamma(n+1/2)\Gamma(n+3/2)z^{-n-1/2}}{
\Gamma(m+1/2)\Gamma(n-l+1/2)\Gamma(l+1)\Gamma(n-m+1)
}{}_2F_1(l-n+1/2,m-n;1/2-n;z),
\ee
\be
B_{l,m,n}(z)=(-1)^{n-m-l}B_{n-m,n-l,n}(z)=
\sum_{p=0}^{-m-1}f_{l,m,n,p}z^p
+\sum_{p=-m}^{-l-1}g_{l,m,n,p}z^p
-\sum_{p=-n}^{-m-1}h_{l,m,n,p}z^p,
\ee
where $f$, $g$ and $h$ are finite sums of $\Gamma$-products.

\cite{MoritaIIS2}
\be
_2F_1(-1/2,-1/2;1;k^2)=2[2E(k)-k'^2K(k)]/\pi .
\ee
\cite{MoritaIIS2}
\be
_2F_1(1/2,1/2;1;k^2)=2K(k)/\pi.
\ee

\cite{MoritaIIS2}
\begin{multline}
-(c-a-n)(c-b-n){}_2F_1(a+n-1,b+n-1;c;z)
\\
+(c-a-b-2n)[c-(a+b+2n-1)z+\frac{2(a+n)(b+n)}{c-a-b-2n-1}(1-z)]
{}_2F_1(a+n,b+n;c;z)
\\
-(1+\frac{2}{c-a-b-2n-1})(a+n)(b+n)(1-z)^2{}_2F_1(a+n+1,b+n+1;c;z)=0.
\end{multline}

\cite{MoritaIIS2}
\begin{multline}
_2F_1(a,b;c;z)
-(1-z){}_2F_1(a+b,b+1;c;z) -\frac{c-a-b-1}{c}z{}_2F_1(a+1,b+1;c+1;z) =0 .
\end{multline}

\cite{MoritaIIS2}
\begin{multline}
_2F_1(a+n-1,b+n-1;c+2n;z)
\\
-\left\{1-z+\frac{1}{c+2n-1}[\frac{(a+n-1)(b+n-1)}{c+2n-2}-\frac{(c-a+n)(c-b+n)}{c+2n}]z
\right\}
{}_2F_1(a+n,b+n;c+2n;z)
\\
+\frac{(a+n)(b+n)(c-a+n)(c-b+n)}
{(c+2n)^2[(c+2n)^2-1]}
z^2 {}_2F_1(a+n+1,b+n+1;c+2n+2;z)=0.
\end{multline}

\cite{RoachSSAC96}
\be
_2F_1(-3/2,-1/2;1/2,z)=\frac{2+z}{2}\sqrt{1-z}+\frac{3\sqrt z}{2}\arcsin\sqrt{z}.
\ee

\cite{RoachSSAC96}
\be
\frac{z^{\mu+1}}{2^\mu\sqrt \pi \Gamma(\frac{3}{2}+\mu)}
\,_2F_1(1,3/2;3/2+\mu,z^2/4)=L_\mu(z) .
\ee

\cite{RoachSSAC96}
\be
_2F_1(-3/2,1/2;3/2,z)
=\frac{5-2z}{8}\sqrt{1-z}+\frac{3}{8\sqrt z}\arcsin\sqrt{z}.
\ee

\cite{RoachSSAC96}
\be
_2F_1(-3/2,-1/2;3/2,z)
=\frac{13+2z}{16}\sqrt{1-z}+\frac{3+12z}{16\sqrt z}\arcsin\sqrt{z}.
\ee

\cite{RoachSSAC96}
\be
_1F_2(-3/2;-1/2,1/2;z)=
(1+2z)\cosh (2\sqrt z)+\sqrt z\sinh(2\sqrt z)
-4z^{3/2}\Shi(2\sqrt z)
.
\ee

\cite{RoachSSAC96}
\be
_1F_2(-3/2;-1/2,2;z)=
-\frac{4+24z-28z^2}{15z\pi}K(\sqrt z)
+\frac{4+56z+4z^2}{15z\pi}E(\sqrt z)
.
\ee

\cite{Vidunasarxiv05}
\be
_2F_1\left(1/4,-1/12;2/3;\frac{x(4+x)^3}{4(2x-1)^3}\right)=
(1-2x)^{-1/4}.
.
\ee

\cite{Vidunasarxiv05}
\be
_2F_1\left(5/4,-1/12;5/3;\frac{x(4+x)^3}{4(2x-1)^3}\right)=
\frac{1+x}{(1+\frac{1}{4}x)^2}(1-2x)^{-1/4}
.
\ee

\cite{Vidunasarxiv05}
\be
_2F_1\left(1/4,7/12;4/3;\frac{x(4+x)^3}{4(2x-1)^3}\right)=
\frac{1}{1+\frac{1}{4}x}(1-2x)^{3/4}
.
\ee

\cite{Vidunasarxiv05}
\be
_2F_1\left(1/4,-5/12;1/3;\frac{x(4+x)^3}{4(2x-1)^3}\right)=
(1+\frac{5}{2}x)(1-2x)^{-5/4}
.
\ee

\cite{Vidunasarxiv05}
\be
_2F_1\left(1/2,-1/6;2/3;\frac{x(2+x)^3}{(2x+1)^3}\right)=
(1+2x)^{-1/2}
.
\ee

\cite{Vidunasarxiv05}
\be
_2F_1\left(1/2,5/6;2/3;\frac{x(2+x)^3}{(2x+1)^3}\right)=
\frac{1}{(1-x)^2}
(1+2x)^{3/2}
.
\ee

\cite{Vidunasarxiv05}
\be
_2F_1\left(1/6,5/6;4/3;\frac{x(2+x)^3}{(2x+1)^3}\right)=
\frac{1}{1+\frac{1}{2}x}
(1+2x)^{1/2}(1+x)^{1/3}
.
\ee

\cite{Vidunasarxiv05}
\be
_2F_1\left(1/6,-1/6;1/3;\frac{x(2+x)^3}{(2x+1)^3}\right)=
(1+2x)^{-1/2}(1+x)^{1/3}
.
\ee

\cite{Vidunasarxiv05}
\be
_2F_1\left(7/24,-1/24;3/4;\frac{108x(x-1)^4}{(x^2+14x+1)^3}\right)=
(1+14x+x^2)^{-1/8}
.
\ee

\cite{Vidunasarxiv05}
\be
_2F_1\left(7/24,23/24;7/4;\frac{108x(x-1)^4}{(x^2+14x+1)^3}\right)=
\frac{1+2x-\frac{1}{11}x^2}{(1-x)^2}
(1+14x+x^2)^{7/8}
.
\ee

\cite{Vidunasarxiv05}
\be
_2F_1\left(5/24,13/24;5/4;\frac{108x(x-1)^4}{(x^2+14x+1)^3}\right)=
\frac{1}{1-x}
(1+14x+x^2)^{5/8}
.
\ee

\cite{Vidunasarxiv05}
\be
_2F_1\left(5/24,-11/24;1/4;\frac{108x(x-1)^4}{(x^2+14x+1)^3}\right)=
\frac{1-22x-11x^2}{(1+14x+x^2)^{11/8}}
.
\ee

\cite{Vidunasarxiv05}
\be
_2F_1\left(19/60,-1/60;4/5;\varphi_1(x)\right)=
(1-228x+494x^2+228x^3+x^4)^{-1/20}
.
\ee
\be
_2F_1\left(19/60,59/60;4/5;\varphi_1(x)\right)=
\frac{(1+66x-11x^2)(1-228x+494x^2+228x^3+x^4)^{19/20}}
{(1+x^2)(1+522x-10006x^2-522x^3+x^4)}
.
\ee
\be
_2F_1\left(11/60,31/60;6/5;\varphi_1(x)\right)=
\frac{(1-228x+494x^2+228x^3+x^4)^{11/20}}
{1+11x-x^2}
.
\ee
\be
_2F_1\left(11/60,-29/60;1/5;\varphi_1(x)\right)=
\frac{1+435x-6670x^2-3335x^4-87x^5}
{(1-228x+494x^2+228x^3+x^4)^{29/20}}
.
\ee
\be
_2F_1\left(13/60,-7/60;3/5;\varphi_1(x)\right)=
\frac{1-7x}
{(1-228x+494x^2+228x^3+x^4)^{7/20}}
.
\ee
\be
_2F_1\left(13/60,53/60;3/5;\varphi_1(x)\right)=
\frac{(1+119x+187x^2+17x^3)(1-228x+494x^2+228x^3+x^4)^{13/20}}
{(1+x^2)(1+522x-10006x^2-522x^3+x^4)}
.
\ee
\be
_2F_1\left(17/60,37/60;7/5;\varphi_1(x)\right)=
\frac{(1+\frac{1}{7}x)(1-228x+494x^2+228x^3+x^4)^{17/20}}
{(1+11x-x^2)^2}
.
\ee
\be
_2F_1\left(17/60,-23/60;2/5;\varphi_1(x)\right)=
\frac{(1+107x-391x^2+1173x^3+46x^4)}
{(1-228x+494x^2+228x^3+x^4)^{23/20}}
.
\ee
Where
\be
\varphi_1(x)=\frac{1728x(x^2-11x-1)^5}{(x^4+228x^3+494x^2-228x+1)^3}.
\ee

\cite{Vidunasarxiv05}
\be
_2F_1\left(7/20,-1/20;4/5;\varphi_2(x)\right)=
\frac{(1+x)^{7/20}}
{(1-x)^{1/20}(1-4x-x^2)^{1/4}}
.
\ee
\be
_2F_1\left(7/20,19/20;4/5;\varphi_2(x)\right)=
\frac{(1+3x)(1+x)^{7/20}(1-x)^{19/20}(1-4x-x^2)^{7/4}}
{(1+x^2)(1+22x-6x^2-22x^3+x^4)}
.
\ee
Where
\be
\varphi_2(x)=\frac{64x(x^2-x-1)^5}{(x^2-1)(x^2+4x-1)^5}.
\ee

\cite{VidunasJCAM178}
\be
_2F_1(a,(a-b+1)/2; a-b+1 ;x )= (1-x/2)^{-a} {}_2F_1(a/2,(a+1)/2 ; (a-b)/2+1; x^2/(2-x)^2) .
\ee

\cite{VidunasJCAM178}
\be
_2F_1(a,(1-a)/3; (4a+5)/6 ;x )= (1-4x)^{-a} {}_2F_1(a/3,(a+1)/3 ; (4a+5)/6; 27x/(4x-1)^3) .
\ee

\cite{VidunasJCAM178}
\be
_2F_1(4a/3,(4a+1)/3; (4a+5)/6 ;x )= (1+8x)^{-a} {}_2F_1(a/3,(a+1)/3 ; (4a+5)/6; 64x(1-x)^3/(1+8x)^3) .
\ee

\cite{VidunasJCAM178}
\be
_2F_1(c,(c+1)/3; (2c+2)/3 ;x )= (1+\omega^2 x)^{-c} {}_2F_1(c/3,(c+1)/3 ; (2c+2)/3; 3(2\omega+1)x(x-2)/(x+\omega)^3) .
\ee
where $\omega^2+\omega+1=0$, $=1/4\pm 3/2m$, $c=1/2\pm 3/2m$.

\cite{VidunasJCAM178}
\be
_2F_1(2/21,5/21; 2/3 ;x )= (1-\frac{33+39\omega}{49}x)^{-1/12} {}_2F_1(1/84,13/84 ; 2/3; \varphi_1(x)) ;
\ee
\be
_2F_1(2/21,3/7; 6/7 ;x )= (1-x)^{-1/84}(1-\frac{241+464\omega}{9}x-\frac{8(62-87\omega)}{27}x^2)^{-1/28} 
{}_2F_1(1/84,29/84 ; 6/7; \frac{1}{\varphi_1(x)}) ;
\ee
where
\be
\varphi_1=\frac{x(x-1)[27x^2-(723+1392\omega)x-496+696\omega]^3}{64[(6\omega+3)x-8-3\omega]^7}.
\ee

\cite{IsmailJAT83,IncePLMS18}
Let
\be
R_{II}(z)=\left\{\begin{array}{ll}
-\frac{(1+a-b)(1-z)^{b-1}}{a}{}_2F_1(a,b;1+a;z), & \Re z<1/2,\\
-\frac{1+a-b}{b-1}z^{-a}{}_2F_1(1-a,1-b;2-b;1-z), & \Re z>1/2.\\
\end{array}
\right.
\ee
This has the continued fraction representation
\be
R_{II}(z)=\frac{1}{z-c_1+K_{n=2}^\infty \frac{z(z-1)\lambda_n}{z-c_n}}
\ee
where
\be
c_n\equiv \frac{n+a-1}{2n+a-b-1};\quad
\lambda_n\equiv \frac{(n-1)(n+a-1-b)}{(2n+a-1-b)(2n+a-3-b)}.
\ee

\cite{ChaundyQJM14}
\be
[\delta(\delta+c-1)-px(\delta+a)(\delta-n)]F_n=0,
\ee
where $F_n\equiv \frac{(c)_n}{n!}\,_2F_1(a,-n;c;px)$.

\cite{ChaundyQJM14}
\be
nF_n-[2n+c-2-p(n+a-1)]F_{n-1}
+(1-p)(n+c-2)F_{n-2}=0
\ee
where $F_n\equiv \frac{(c)_n}{n!}\,_2F_1(a,-n;c;p)$.

\cite{VidunasJMAA355}
\begin{multline}
\frac{(\frac{a+1}{2})_l}{(1/2)_l}{}_2F_1
(a/2, (a+1)/2+l; 1/2-k ;z)
=
\frac{(1+\surd z)^{-a}}{2}
F_2(\begin{array}{c}
a; k,-l\\
-2k,-2l
\end{array}
\mid\frac{2\surd z}{1+\surd z},\frac{2}{1+\surd z})
\\
+\frac{(1-\surd z)^{-a}}{2}
F_2(\begin{array}{c}
a; k,-l\\
-2k,-2l
\end{array}
\mid\frac{2\surd z}{\surd z-1},\frac{2}{1-\surd z})
\end{multline}
where the $F_2$ series on the r.h.s. are finite sums.

\cite{VidunasJMAA355}
\be
_2F_1(a/2,(a+1)/2 ; 1/2;z) = \frac{(1-\surd z)^{-a}+(1+\surd z)^{-a}}{2}.
\ee

\cite{VidunasJMAA355}
\begin{multline}
\frac{(\frac{a+1}{2})_k(a/2)_{k+l+1}}{(1/2)_k (1/2)_{k+1} (1/2)_l}
(-)^k z^{k+1/2} {}_2F_1( (a+1)/2+k,a/2+k+l+1 ; 3/2+k ; z)
\\
=
\frac{(1+\surd z)^{-a}}{2}
F_2(\begin{array}{c}
a;-k,-l \\
-2k,-2l
\end{array}\mid \frac{2\surd z}{1+\surd z},\frac{2}{1+\surd z})
\\
-\frac{(1-\surd z)^{-a}}{2}
F_2(\begin{array}{c}
a;-k,-l \\
-2k,-2l
\end{array}\mid \frac{2\surd z}{\surd z-1},\frac{2}{1-\surd z})
.
\end{multline}

\cite{VidunasJMAA355}
\be
_2F_1((a+1)/2,(a+2)/2 ; 3/2;z) = \frac{(1-\surd z)^{-a}-(1+\surd z)^{-a}}{2a\surd z}, a\neq 0.
\ee
\cite{VidunasJMAA355}
\be
_2F_1(1/2,1; 3/2;z) = \frac{\log(1+\surd z)-\log(1-\surd z)}{2\surd z}.
\ee

\cite{VidunasJMAA355}
\begin{multline}
_2F_1( (a-l)/2,(a+l+1)/2; a+k+1 ; z)
=
(\frac{1+\sqrt{1-z}}{2})^{-a-k}
(1-z)^{k/2}
\\
F_3(\begin{array}{c}
k+1,l+1;-k-l\\
a+k=1
\end{array}\mid \frac{\sqrt{1-z}-1}{2\sqrt{1-z}},\frac{1-\sqrt{1-z}}{2})
.
\end{multline}

\cite{VidunasJMAA355}
\be
_2F_1(a/2,(1+a)/2; a+1;z) = (\frac{1+\sqrt{1-z}}{2})^{-a}.
\ee

\cite{WilletQJM2}
\be
1={}_2F_1(\alpha,\beta;\gamma;x){}_2F_1(-\alpha,-\beta;-\gamma;x)
+\frac{\alpha\beta(\alpha-\gamma)(\beta-\gamma)}{\gamma^2(1-\gamma^2)}
x^2 
{}_2F_1(1+\alpha,1+\beta;2+\gamma;x)
{}_2F_1(1-\alpha,1-\beta;2-\gamma;x).
\ee

\cite{WilletQJM2}
\be
_2F_1(\alpha,\beta;\gamma;px)
{}_2F_1(\alpha',\beta';\gamma';qx)
=\sum_{n}\frac{(\alpha)_n(\beta)_n}{(\gamma)_nn!}{}_4F_3(\alpha',\beta',1-\gamma-n,-n;1-\alpha-n,1-\beta-n,\gamma';q/p)
(px)^n.
\ee
\cite{DriverNumAlg26}
\be
_2F_1(a,b;c;z)
{}_2F_1(a,b;a+b-c+1;z)
={}_4F_3(a,b,(a+b)/2,(a+b+1)/2 ; a+b, a+b-c+1; 4z(1-z)).
\ee
\cite{BaileyPLMS28}
\be
_0F_3(\rho,\rho/2,(1+\rho)/2;-x^2/4)
={}_0F_1(\rho;x){}_0F_1(\rho;-x).
\ee

\cite{MillerJPA38,KimMMN10}
\be
_2F_2(a,c+1;b,c;x) = e^x {}_2F_2(b-a-1,\alpha+1;b,\alpha;-x),\quad
\alpha\equiv c(1+a-b)/(a-c). \label{eq.Millalpha}
\ee

\cite{ParisJCAM173,ChuJCAM216}
\be
_2F_2(a,d;b,c;x)
=e^x\sum_{n\ge 0}\frac{(c-d)_n}{(c)_nn!} (-x)^n\,
_2F_2(b-a,d;b,c+n;-x)
.
\ee

\cite[(9a)]{ChuJCAM216}
\be
_2F_2(a,c+m;b,c;x)
=e^x 
\sum_{n=0}^m \binom{m}{n} \frac{(a)_n}{(b)_n(c)_n}x^n
{}_1F_1(b-a;b+n;-x)
.
\ee

\cite[(12)]{ChuJCAM216}
\be
_2F_2(a,d;b,c;x)
=\frac{\Gamma(b)\Gamma(c)}{\Gamma(a)\Gamma(b+c-a)} e^x
\sum_{k=0}^\infty \frac{(b-a)_k(c-a)_k}{(1)_k(b+c-a)_k}
{}_1F_1(b+c-a-d+k;b+c-a+k;-x)
.
\ee

\cite[(13)]{ChuJCAM216}
\be
_2F_2(a,c+m;b,c;x)
=(-)^m \frac{(1-b)_m}{(c)_m} e^x
\sum_{k=0}^m \frac{(b-c-m)_k(-m)_k}{(1)_k(b-m)_k}
{}_1F_1(b-a-m+k;b-m+k;-x)
.
\ee

\cite{RoachSSAC96}
\be
_2F_2(-3/2,-1/2;-5/2,1;z)
=
\frac{5-4z}{5}e^{z/2}I_0(z/2)+\frac{4z}{5}e^{z/2}I_1(z/2).
\ee

\cite{BaileyPLMS28}
\be
_2F_3((\rho+\sigma)/2, (\rho+\sigma-1)/2; \rho,\sigma,\rho+\sigma-1; 4x)
={}_0F_1(\rho;x){}_0F_1(\sigma;x).
\ee

\cite{BaileyPLMS28}
\be
_2F_3(a,\rho-a;\rho,\rho/2,(1+\rho)/2; x^2/4)
={}_1F_1(a;\rho;x){}_1F_1(a;\sigma;-x).
\ee

\cite{BaileyPLMS28}
\be
_2F_3((a+\beta)/2,(a+\beta+1)/2;a+1/2,\beta+1/2,a+\beta; x^2/4)
={}_1F_1(a;2a;x){}_1F_1(\beta;2\beta;-x).
\ee

\cite{Qureshiarxiv18}
\be
_3F_2(1,1,(p+q)/q ; 2, (p+2q)/2 ;1)
=\frac{q+p}{p}
\left[
\frac{q}{p}-\ln 2q-\frac{\pi}{2}\cot\frac{\pi p}{q}
+2\sum_{j=1}^{[q/2]}\left\{ \cos\frac{2\pi pj}{q}\ln \sin \frac{\pi j}{q}\right\}
\right]
.
\ee

\cite{GottschalkJPA21}
\be
_3F_2(1,a,b; a+1,b+1;1) = \frac{ab}{a-b}[\psi(a)-\psi(b)].
\ee

\cite{GottschalkJPA21}
\begin{multline}
_3F_2(a,a,b; a+1,a+1;z) = \frac{\Gamma(a+1)a\Gamma(2-b)[(1-z)^{2-b}-1]}{\Gamma(a+1-b)z^2}
\sum_{i=1}^{a-1}\frac{1}{i(i+1-b)}
\\
+\frac{a\Gamma(a+1)(1-z)^{2-b}}{\Gamma(a+1-b)z^a}\sum_{i=1}^{a-1}\frac{1}{i(i+1-b)}
\sum_{j=1}^{i-1}\frac{\Gamma(j+2-b)z^j}{\Gamma(j+1)}
-\frac{a\Gamma(a+1)\Gamma(1-b)}{\Gamma(a+1-b)z^a}\theta(b,z)
\end{multline}
where
\begin{equation}
\theta(b,z)=\left\{
\begin{array}{rl}
\sum_{j=1}^{b-2}\frac{(1-z)^{-j}-1}{j}-\ln(1-z),& a\in \mathbb{Z}^+, b\in \mathbb{Z}^+, b\ge a+1\\
\sum_{j=1}^{'b-3/2}\frac{(1-z)^{1/2-j}-1}{j-1/2}+2\ln|1-\sqrt{1-z}|-2\ln|z|+2\ln 2,& a\in \mathbb{Z}^+, 2b\in \mathbb{Z}^+, b\not \in \mathbb{Z}.
\end{array}
\right.
\label{eq.3F2Plus}
\end{equation}

\cite{Amdeberhanarxiv3663}
\be
_3F_2(1,1,2-t;2,3;1)=\frac{2(1-\gamma-\psi(t+1))}{1-t}.
\ee

\cite{RoyAMM94}
\be
_3F_2(-n,-a,-b;c,2-n-a-b-c;1)
=
\frac{(c+b-1)_n(c+a)_n}{(c+a+b-1)_n(c)_n}\left[1-\frac{a}{(c+b-1)(a+c+n-1)}\right].
\ee
for $n=1,2,3,\ldots$

\cite{Koornwinderarxiv98}\cite[(2.3.1.3)]{SlaterHyp}
\be
_3F_2(-m,b,c;e,-m+b+c-e+1;1)
=
\frac{(e-b)_m(e-c)_m}{(e)_m(e-b-c)_m};
\quad
m=0,1,2,\ldots,\quad e,e-b-c \neq 0,-1,-2,\ldots
\ee

\cite{KimarXiv1205}
\be
_3F_2(m+1/2,1/2-m,d+1 ; m+5/2,d;1) = (1+\frac{1-2m}{2d})\frac{3\pi}{2^{2m+3}m!}(5/2)_m.
\ee

\cite{BaileyPLMS28}
\be
_3F_2(a,b,c;d,e;1) = \frac{\Gamma(d)\Gamma(1-e+a)\Gamma(1-e+b)\Gamma(1-e+c)}{\Gamma(1-e)\Gamma(d-a)\Gamma(d-b)\Gamma(d-c)},
\ee
where $a$, $b$ or $c$ a negative integer and $d+e=a+b+c+1$.

\cite{BaileyPLMS28}
\be
_3F_2(a,a-d+1,a-2+1;d,e;1) = \frac{\Gamma(d)\Gamma(e)\Gamma(d+e-3a/2-1)\Gamma(1/2)}
{2^n\Gamma((1+a)/2)\Gamma(d-a/2)\Gamma(e-a/2)\Gamma(d+e-a+1)}.
\ee

\cite{BaileyPLMS28}
\be
_3F_2(a,b,c;(a+b+1)/2,2c;1) =\frac{\Gamma(1/2)\Gamma(1/2+c)\Gamma((1+a+b)/2
\Gamma((1-a-b)/2+c)}
{\Gamma((1+a)/2)\Gamma((1+b)/2)\Gamma((1-a)/2+c)\Gamma((1-b)/2+c)}
\ee
provided that $a+b=1$ and $d+2=2c+1$.

\cite{ExtonJCAM83}
\be
_3F_2(a,1+a/2,-n; a/2,b;1) = \frac{(2+a-b)_n (b-a-1)_n}{(b)_n (1+a-b)_n}
.
\ee

\cite{ExtonJCAM83}
\be
_3F_2(a,b,-n; 1+a-b,1+2b-n;1) = \frac{(a-2b)_n(1+a/2-b)_n(-b)_n}{(1+a-b)_n(a/2-b)_n(-2b)_n}.
\ee

\cite{ExtonJCAM83}
\be
(1-x)^e {}_3F_2(a+e,1+a/2+e/2,e/2; 1+a+e/2,a/2+e/2;x)
={}_2F_1(a,-e/2,1+a+e/2;x).
\ee
\cite{BaileyPLMS28}
\be
_3F_2(2a,2\beta,a+\beta ; 2a+2\beta, a+\beta+1/2;x)
=[{}_2F_1(a,\beta;a+\beta+1/2;x)]^2.
\ee
\cite{BaileyPLMS28}
\be
_3F_2(2a,2\beta,a+\beta ; 2a+2\beta-1, a+\beta+1/2;x)
={}_2F_1(a,\beta;a+\beta-1/2;x)
{}_2F_1(a,\beta;a+\beta+1/2;x).
\ee
\cite{BaileyPLMS28}
\be
_3F_2(2a-1,2\beta,a+\beta-1 ; 2a+2\beta-2, a+\beta-1/2;x)
={}_2F_1(a,\beta;a+\beta-1/2;x)
{}_2F_1(a,\beta-1;a+\beta-1/2;x).
\ee

\cite{MillerJPA38}
\be
_3F_2(-n,a,c+1; b,c;1) = \frac{(b-a-1)_n(\alpha+1)_n}{(b)_n(\alpha)_n}
\ee

\cite{MillerJPA38,KimarXiv1205}
\be
_3F_2(f,a,c+1; b,c;1) = \frac{(c-a)(\alpha-f)}{c}
\frac{\Gamma(b)\Gamma(b-a-f-1)}{\Gamma(b-a)\Gamma(b-f)}
\ee
where $\alpha$ is defined in \eqref{eq.Millalpha}.

\cite{KimMMN10}
\be
_3F_2(-n,b-a-1,\alpha +1;b,\alpha ;1)=\frac{(a)_n(c+1)_n}{(b)_n(c)_n}
\ee
where $\alpha$ is defined in \eqref{eq.Millalpha}.

\cite{BaileyPLMS28}
\be
_3F_2(a,\rho-\beta,-n; \rho,a-\gamma;1)
=(-1)^n\frac{(1+\gamma-n)_n}{(a-\gamma)_n}
{}_3F_2(a,\beta,-n;\rho,1+\gamma-n;1).
\ee

\cite{KimMMN10}
\be
_3F_2(a,b,d+1;c,d;1) = \frac{\Gamma(c)\Gamma(c-a-b-1)}{\Gamma(c-a)\Gamma(c-b)}[c-a-b-1+\frac{ab}{d}].
\ee

\cite{LavoieMC49}
\begin{multline}
_3F_2(a,b,c;\frac{a+b+1}{2},2c-1;1)
=
\frac{2^{a+b}\Gamma(\frac{a+b+1}{2})\Gamma(c-1/2)\Gamma(c-\frac{a+b+1}{2})}
{\Gamma(1/2)\Gamma(a+1)\Gamma(b+1)}\\ \times
\left(
\frac{\Gamma(1+a/2)\Gamma(1+b/2)}{\Gamma(c-a/2-1/2)\Gamma(c-b/2-1/2)}
+\frac{ab\Gamma(\frac{a+1}{2})\Gamma(\frac{b+1}{2})}{4\Gamma(c-a/2)\Gamma(c-b/2)}
\right),\quad \Re(2c-a-b)>1
\end{multline}

\cite{LavoieMC49}
\begin{multline}
_3F_2(a,b,c;\frac{a+b+1}{2},2c+1;1)
=
\frac{2^{a+b}\Gamma(\frac{a+b+1}{2})\Gamma(c+1/2)\Gamma(c-\frac{a+b-1}{2})}
{\Gamma(1/2)\Gamma(a+1)\Gamma(b+1)}\\ \times
\left(
\frac{\Gamma(1+a/2)\Gamma(1+b/2)}{\Gamma(c-a/2+1/2)\Gamma(c-b/2+1/2)}
-\frac{ab\Gamma(\frac{a+1}{2})\Gamma(\frac{b+1}{2})}{4\Gamma(c-a/2+1)\Gamma(c-b/2+1)}
\right),\quad \Re(2c-a-b)>-3.
\end{multline}

\cite{LavoieMC49}
\begin{multline}
_3F_2(a,b,c;\frac{a+b+1}{2},2c;1)
=
\frac{\Gamma(1/2)\Gamma(\frac{a+b+1}{2})\Gamma(c+1/2)\Gamma(c-\frac{a+b-1}{2})}
{\Gamma(\frac{a+1}{2})\Gamma(\frac{b+1}{2})\Gamma(c-a/2+1/2)\Gamma(c-b/2+1/2)}
,\quad \Re(2c-a-b)>-1.
\end{multline}

\cite{LavoieMC49}
\begin{multline}
_3F_2(a,b,c;e,f;1)
=
\frac{\Gamma(e)\Gamma(f)}
{2^{2a+1}\Gamma(e-a)\Gamma(f-a)}
\left(
\frac{\Gamma(\frac{e-a}{2})\Gamma(\frac{f-a}{2})}{\Gamma(\frac{e-b}{2})\Gamma(\frac{f-b}{2})}+\frac{\Gamma(\frac{e-a+1}{2})\Gamma(\frac{f-a+1}{2})}{\Gamma(\frac{e-b+1}{2})\Gamma(\frac{f-b+1}{2})}
\right)\\
,\, a+b=0, e+f=1+2c,\Re c >-1.
\end{multline}

\cite{LavoieMC49}
\begin{multline}
_3F_2(a,b,c;e,f;1)
=
\frac{\Gamma(e)\Gamma(f)}
{2^{2a-1}(a-1)(c-1)\Gamma(e-a)\Gamma(f-a)}
\left(
\frac{\Gamma(\frac{e-a}{2})\Gamma(\frac{f-a}{2})}{\Gamma(\frac{e-b}{2})\Gamma(\frac{f-b}{2})}-\frac{\Gamma(\frac{e-a+1}{2})\Gamma(\frac{f-a+1}{2})}{\Gamma(\frac{e-b+1}{2})\Gamma(\frac{f-b+1}{2})}
\right)\\
,\, a+b=2, e+f=1+2c,\Re c >1.
\end{multline}

\cite{LavoieMC49}
\begin{multline}
_3F_2(a,b,c;e,f;1)
=
\frac{\pi\Gamma(e)\Gamma(f)}
{2^{2c-1}\Gamma(\frac{e+a}{2})\Gamma(\frac{f+a}{2})\Gamma(\frac{e+b}{2})\Gamma(\frac{f+b}{2})},\, a+b=1, e+f=1+2c,\Re c >0.
\end{multline}

\cite{LavoieMC49}
\begin{multline}
_3F_2(1,1+a,a+b;1+b,1+a+b;1)
=\frac{b(b+1)}{a}\big\{ \psi(b+a)-\psi(b-a)-\psi(\frac{1+b+a}{2})+\psi(\frac{1+b-a}{2})\big\}
\end{multline}
for $Re (b-a)>0$.

\cite{LavoieJCAM72}
\begin{multline}
_3F_2(a,b,c;e,f;1)
=
\frac{\Gamma(e)\Gamma(f)\Gamma(c-\frac12(j+|j|))
\Gamma(e-c-\frac12(i+|i|))
\Gamma(a-\frac12(i+j+|i+j|))
}{2^{2a-i-j}\Gamma(e-a)\Gamma(f-a)\Gamma(e-c)\Gamma(a)\Gamma(c)}
\\
\times
\big\{
A_{i,j}\frac{\Gamma(e/2-a/2+\frac14(1-(-)^i)\Gamma(f/2-a/2))}
{\Gamma(e/2+a/2-i/2+[-j/2]\Gamma(f/2+a/2-i/2+((-)^j/4)(-)^i-1)+[-j/2]}
\\
+B_{i,j}\frac{\Gamma(e/2-a/2+\frac14(1+(-)^i))}{\Gamma(e/2+a/2-1/2-i/2+[-j/2+1/2])}
\frac{\Gamma(f/2-a/2+1/2}{\Gamma(f/2+a/2-1/2-i/2+((-)^j/4)(1-(-)^j+[-j/2+1/2])}
\big\}
\end{multline}
where $A_{i,j}$ and $B_{i,j}$ are tabulated functions of $a,b,c$ for $0\le |i|,|j|\le 3$,
$[.]$ is the greatest integer less than or equal to the argument
and $a+b=1+i+j$, $e+f=2c+1+i$.

\cite{DriverETNA25}
\be
_3F_2(a,b/2,(b+1)/2;c/2,(c+1)/2;1)
=\frac{\Gamma(c)\Gamma(c-a-b)}{\Gamma(c-b)\Gamma(c-a)}\,_2F_1(a,b;c-a;-1),
\quad \Re c> \Re b>0, \Re(c-a-b)>0.
\label{eq.drive}
\ee

\cite{DriverETNA25}
\be
_3F_2(a,b/2,(b+1)/2;c/2,(c+1)/2;1/2)
=2^a\sum_{k=0}^\infty \binom{-a}{k}\frac{(c-b)_k}{(c)_k}\,_2F_1(-k,b;c+k;-1),
\quad \Re c> \Re b>0.
\ee

\cite{Baileyarxiv08,BaileyJPA41}
\be
8\pi^3 {} _3F_2(-1/2,1/2,1/2;1,1;1)
=\Gamma^4(1/4)+16\Gamma^4(3/4).
\ee

\cite{ChaundyQJM14}
\be
_0F_1(a;px)\,_0F_1(c';qx)
=\sum_{n\ge 0}\frac{(px)^n}{n!(c)_n}
\,_2F_1(1-c-n,-n;c';q/p).
\ee

\cite{ChaundyQJM14}
\be
_1F_1(a;c;px)\,_1F_1(a';c';qx)
=\sum_{n\ge 0}\frac{(a)_n(px)^n}{n!(c)_n}
\,_3F_2(a',1-c-n,-n;c',1-a-n;-q/p).
\ee

\cite{ChaundyQJM14}
\be
_2F_0(a,b;px)\,_2F_0(a',b';qx)
=\sum_{n\ge 0}\frac{(a)_n(b)_n(px)^n}{n!}
\,_3F_2(a',b',-n;1-a-n,1-b-n;-q/p).
\ee

\cite{ChaundyQJM14}
\be
_2F_1(a,b;c;px)\,_2F_1(a',b';c';qx)
=\sum_{n\ge 0}\frac{(a)_n(b)_n(px)^n}{n!(c)_n}
\,_4F_3(a',b',1-c-n,-n;c',1-a-n,1-b-n;q/p).
\ee

\cite{RuehrSIAMR24}
\be
(a+b+1)_3F_2[-c,-a,1;b+1,\frac{1-a-c}{2};\frac{1}{2}]
=(b+1)_3F_2[-\frac{a}{2},\frac{1-a}{2},1;\frac{1-a-b}{2},\frac{1-a-c}{2};1]
\ee
where $a$, $b$ and $c$ are positive integers of the same parity.

\cite{KimarXiv1205}
\begin{multline}
_3F_2(a,b,d+1; (a+b+3)/2,d;1/2) = \frac{\Gamma(1/2)\Gamma(a/2+b/2+3/2)\Gamma(a/2-b/2-1/2)
}{\Gamma(a/2-b/2+3/2)}
\\ \times
[\frac{(a+b-1)/2-ab/d}{\Gamma((a+1)/2)\Gamma((b+1/)/2)}+\frac{(a+b+1)/d-2}{\Gamma(a/2)\Gamma(b/2)}].
\end{multline}

\cite{KimarXiv1205}
\begin{multline}
_3F_2(a,1-a,d+1; c+1,d;1/2) = 2^{-c}\Gamma(1/2)\Gamma(c+1)
\\ \times
[\frac{2/d}{\Gamma(a/2+c/2)\Gamma(c/2-a/2+1/2)}+\frac{1-c/d}{\Gamma(a/2+c/2+1/2)\Gamma(c/2-a/2+1)}].
\end{multline}

\cite{WeiArxiv1709}
\begin{multline}
b {} _3F_2(a,1-a,1+b; c+1,1+2b-c;1) = 
\\
\Gamma(1/2+c/2)\Gamma(1+c/2)\Gamma(b+1/2-c/2)\Gamma(b+1-c/2)
\\
\times
[
\frac{1}
{\Gamma(a/2+c/2)\Gamma(1/2-a/2+c/2)\Gamma(b+1/2+a/2-c/2)\Gamma(b+1-a/2-c/2)}
\\
+\frac{1}
{\Gamma(1/2+a/2+c/2)\Gamma(1-a/2+c/2)\Gamma(b+a/2-c/2)\Gamma(b+1/2-a/2-c/2)}
]
\end{multline}

\cite{KimarXiv1205}
\begin{multline}
_3F_2(1+2m,1+2n,d+1; m+n+5/2,d;1/2) = \frac{2\pi(3/2)_{m+n+1}}{(2m-2n+1)(2m-2n-1)}
\\ \times
[\frac{(m+n+1/2)-(2m+1)(2n+1)/2}{m!n!}+\frac{(2m+2n+3)/d-2}{\pi(1/2)_m(1/2)_n}].
\end{multline}

\cite{KimarXiv1205}
\begin{multline}
_3F_2(1/2+m,1/2-m,d+1; m+2n+5/2,d;1/2) = \frac{\pi}{2^{m+2n+5/2}}(3/2)_{m+2n+1}
\\ \times
[\frac{2/d}{(m+n)!n!}+\frac{1-(3/2+m+2n)/d}{\pi(1/2)_{m+n+1}(1/2)_{n+1}}].
\end{multline}

\cite{DriverETNA25}
\be
_3F_2(-n,b/2,(b+1)/2;c/2,(c+1)/2;1)
=\frac{(c-b)_n}{(c)_n}\,_2F_1(-n,b;c+n;-1),
\quad \Re c> \Re b>0.
\ee

\cite{Koornwinderarxiv98}
\begin{multline}
_3F_2(a,b,c;e,a+b+c-e+1;1)
+\frac{\Gamma(e-1)\Gamma(a-e+1)\Gamma(b-e+1)\Gamma(c-e+1)\Gamma(a+b+c-e+1}
{\Gamma(1-e)\Gamma(a)\Gamma(b)\Gamma(c)\Gamma(a+b+c-2e+2)}
\nonumber
\\
 \quad \times
_3F_2(a-e+1,b-e+1,c-e+1;2-e,a+b+c-2e+2;1)
\nonumber
\\
 =
\frac{\Gamma(a-e+1)\Gamma(b-e+1)\Gamma(c-e+1)\Gamma(a+b+c-e+1)}
{\Gamma(1-e)\Gamma(b+c-e+1)\Gamma(a+c-e+1)\Gamma(a+b-e+1)}
.
\end{multline}

\cite{Koornwinderarxiv98}
\be
_3F_2(a,b,c;e,a+b+c-e+1;1)
=
\frac{\Gamma(e)\Gamma(a+b+c-e+1)}
{\Gamma(a)\Gamma(b+1)\Gamma(c+1)}
\,_3F_2(e-a,b+c-e+1,1;b+1,c+1;1);
\quad \Re a>0.
\label{eq.kornw}
\ee

\cite{RathieJCAM167}\cite{Milgramarxiv1011}\cite{RaoJPA25}
\begin{gather}
_3F_2(a,b,c;f,e;1)
\nonumber
\\
=
\frac{\Gamma(f)\Gamma(e)\Gamma(f+e-a-b-c)}
{\Gamma(c)\Gamma(e-b+f-c)\Gamma(f+e-a-c)}
\,_3F_2(e+f-a-b-c,f-c,e-c;e-b+f-c,e+f-a-c;1)
\\
=
\frac{\Gamma(f)\Gamma(e)\Gamma(f+e-a-b-c)}
{\Gamma(b)\Gamma(e-b+f-c)\Gamma(e-b+f-a)}
\,_3F_2(e+f-a-b-c,f-b,e-b;e-b+f-c,e-b+f-a;1)
\\
=
\frac{\Gamma(e)\Gamma(f+e-a-b-c)}
{\Gamma(e-a)\Gamma(e-b+f-c)}
\,_3F_2(f-c,f-b,a;e-b+f-c,f;1)
\\
=
\frac{\Gamma(f)\Gamma(f+e-a-b-c)}
{\Gamma(f-a)\Gamma(e-b+f-c)}
\,_3F_2(e-c,e-b,a;e-b+f-c,e;1)
\\
=
\frac{\Gamma(e)\Gamma(f)\Gamma(f+e-a-b-c)}
{\Gamma(a)\Gamma(e+f-a-c)\Gamma(e-b+f-a)}
\,_3F_2(e+f-a-b-c,f-a,e-a;e+f-a-c,e-b+f-a;1)
\\
=
\frac{\Gamma(e)\Gamma(f+e-a-b-c)}
{\Gamma(e-b)\Gamma(e+f-a-c)}
\,_3F_2(f-c,f-a,b;e+f-a-c,f;1)
\\
=
\frac{\Gamma(f)\Gamma(f+e-a-b-c)}
{\Gamma(f-b)\Gamma(e+f-a-c)}
\,_3F_2(e-c,e-a,b;e+f-a-c,e;1)
\\
=
\frac{\Gamma(e)\Gamma(f+e-a-b-c)}
{\Gamma(e-c)\Gamma(e-b+f-a)}
\,_3F_2(f-b,f-a,c;e+f-a-b,f;1)
\\
=
\frac{\Gamma(f)\Gamma(f+e-a-b-c)}
{\Gamma(f-c)\Gamma(e-b+f-a)}
\,_3F_2(e-b,e-a,c;e+f-a-b,e;1)
.
\end{gather}

\cite{BaileyPGMA2,KarpMdpiM8}
\be
_3F_2(a,b,c;d,e;1)
= \frac{\Gamma(d)\Gamma(e)\Gamma(s)}{\Gamma(c)\Gamma(s+a)\Gamma(s+b)}
{}_3F_2(d-c,e-c,s ; s+a,s+b;1)
\ee
where $s\equiv d+e-a-b-c$.

\cite{BaileyPGMA2}
\begin{multline}
a(d+e-a-b-c-1){}_3F_2(a+,b,c;d,e;1)
-(d-a)(e-a){}_3F_2(a-,b,c;d,e;1)
\\
= 
[a(2d+2e-2a-b-c-1)-de+bc]{}_3F_2(a,b,c;d,e;1).
\end{multline}

\cite[(6)]{RajeswariJPA22}
\be
_3F_2(-n,\alpha,\beta;\gamma,\delta;1)
=
\frac{\Gamma(\gamma)\Gamma(\gamma+n-\alpha)}
{\Gamma(\gamma+n)\Gamma(\gamma-\alpha)}
\,_3F_2(-n,\alpha,\delta-\beta;1+\alpha-\gamma-n,\delta;1).
\ee

\cite{ChuMC81}
\begin{multline}
_3F_2(a,b,c; (1+a+b+m)/2,d;1)
=
\frac{(b)_m}{(b-a)_m}\sum_{i=0}^m \binom{m}{i}
{}_3F_2(a+i,b+m-i,c; (1+a+b+m)/2,d;1)
\\
\times
\frac{a-b-m+2i}{a-b-m}\left[\begin{array}{c} a,a-b-m\\ 1+a-b,1-b-m\end{array}\right]_i
\end{multline}
plus similar reductions of $_3F_2(a,b,c;d,2c+n;1)$ to $_3F_2(a,b,c;d,2c;1)$
and 
$_3F_2(a,b,c;(1+a+b+m)/2,2c+n;1)$ to $_3F_2(a,b,c;(1+a+b)/2,2c;1)$.

\cite{ChuMC81}
\begin{multline}
_3F_2(a,b,c; (1+a+b-m)/2,d;1)
=
\left[\begin{array}{c} b,(1-a+b-m)/2\\ b-a,(1+a+b-m)/2\end{array}\right]_m
\sum_{i=0}^m 
(-)^i
\binom{m}{i}
\\
\times
{}_3F_2(a+i,b+m-i,c; (1+a+b+m)/2,d;1)
\frac{a-b-m+2i}{a-b-m}
\left[\begin{array}{c} a,a-b-m\\ 1+a-b,1-b-m\end{array}\right]_i
\end{multline}

\cite{BaileyPGMA2}
\begin{multline}
b(a-b+1)(d-a)(e-a){}_3F_2(a-,b+,c;d,e;x)
-a(b-a+1)(d-b)(e-b){}_3F_2(a+,b-,c;d,e;x)
\\
=
(a-b)[c(b-a+1)(a-b+1)x+de(a+b-1)-2(d+e)ab+ab(a+b+1)]{}_3F_2(a,b,c;d,e;x)
\\
+(a-b+1)(a-b)(a-b-1)(1-x)\frac{abcx}{de}{}_3F_2(a+,b+,c+;d+,e+;x)
.
\end{multline}

\cite{EbisuJMAA463,KarpMdpiM8}
\begin{equation}
_3F_2(a,b,c;d,e;1)
=
u({\bf r})
{}_3F_2(a+k_1,b+k_2,c+k_1;d+k_4,e+k_5;1)
+u({\bf r})
{}_3F_2(a+m_1,b+m_2,c+m_3;d+m_4,e+m_5;1)
\end{equation}

\cite[Entry 13]{BerndtBLMS15}\cite{ClausenCrelle3}
If $\alpha+\beta+\gamma=0$ then
\be
_3F_2(-2\alpha,-2\beta,\gamma;\gamma+1/2,2\gamma;x) ={}_2F_1^2(-\alpha,-\beta;\gamma+1/2;x).
\label{eq.rama3f2}
\ee

\cite{BerndtBLMS15}
\begin{multline}
_3F_2(\alpha,\beta,\gamma;\delta,\epsilon;1)
=
\frac{\Gamma(\delta)\Gamma(\delta-\alpha-\beta)}{
\Gamma(\delta-\alpha)\Gamma(\delta-\beta)}
{}_3F_2(\alpha,\beta,\epsilon-\gamma;\alpha+\beta-\delta+1,\epsilon;1)
\\
+
\frac{
\Gamma(\delta)\Gamma(\epsilon)\Gamma(\alpha+\beta-\gamma)
\Gamma(\delta+\epsilon-\alpha-\beta-\gamma)
}{
\Gamma(\alpha)\Gamma(\beta)\Gamma(\epsilon-\gamma)
\Gamma(\delta+\epsilon-\alpha-\beta)
}
{}_3F_2(\delta-\alpha,\delta-\beta,
\delta+\epsilon-\alpha-\beta-\gamma;\delta-\alpha-\beta+1,
\delta+\epsilon-\alpha-\beta;1)
.
\end{multline}

\cite{KrupnikovJCAM78}
\be
_3F_2(a,a+1/2,b; c,2a+b-c+1;z)
=
[\frac{2}{z}(1-\sqrt{1-z})]^{2a}
{}_3F_2(2a,c-b,2a-c+1; c,2a+b-c+1 ; 1-\frac{2}{z}(1-\sqrt{1-z})).
\ee

\cite{KimMMN10}
\be
_3F_2(h,a,c+1;b,c;x)
=(1-x)^{-h}{}_3F_2(h,b-a-1,\alpha+1 ; b,\alpha ; -\frac{x}{1-x})
\ee
where $\alpha$ is given in \eqref{eq.Millalpha}.

\cite{BaileyPLMS28}
\be
(1-x)^{-3a}{}_3F_2(a,a+1/3,a+2/3 ; \rho_1,\rho_2 ; \frac{-27x}{4(1-x)^3})
=
{}_3F_2(3a,1/2+\rho_1-\rho_2,1/2+\rho_2-\rho_1;\rho_1,\rho_2;x/4),
\ee
and
\be
(1-x)^{-3a}{}_3F_2(a,a+1/3,a+2/3 ; \rho_1,\rho_2 ; \frac{27x^2}{4(1-x)^3})
=
{}_3F_2(3a,\rho_1-1/2,\rho_2-1/2;2\rho_1-1,2\rho_2-1;4x).
\ee
provided $\rho_1+\rho_2=3a+3/2$.

\cite{BaileyPLMS28}
\be
_3F_2(2a-1,a+1/2,a-b-1/2;a+b+1/2,a-1/2;x)
= (1-x)^{1-2a}{}_2F_1(a,b;a+b+1/2;\frac{-4x}{(1-x)^2}).
\ee
\cite{BaileyPLMS28}
\be
(1-x)^{1-2a}{} _3F_2(a,a-1/2,\rho_1+\rho_2-2a;\rho_1,\rho_2;\frac{-4x}{(1-x)^2})
= {}_3F_2(2a-1,2a-\rho_1,2a-\rho_2;\rho_1,\rho_2;x)
\ee
\cite{BaileyPLMS28}
\be
(1-x)^{-3a}{} _3F_2(a,a+1/3,a+2/3;1/3,2/3;\frac{-x^3}{(1-x)^3})
= 1+\frac{2}{3}\sum_{n=1}^\infty \frac{(3a)_n}{n!}\cos\frac{n\pi}{6} 3^{n/2}x^n,
\ee
plus 2 similar formulas for denominators $(2/3,4/3)$ or $(4/3,5/3)$ on the LHS.

\cite{BaileyPLMS28}
\begin{multline}
_3F_2(2a-1,a+1/2,a-b-1/2; a+b+1/2,a-1/2;x)
\\
=(1-x)^{2b-2a}{}_3F_2(2b-1,b+1/2,b-a-1/2;a+b+1/2,b-1/2;x).
\end{multline}

From \eqref{eq.3F2Plus}
\be
{}_3F_2(\frac12,\frac12,\frac12;\frac32,\frac32;\frac12)
=
\frac{1}{\surd 2}[\frac{\pi}{4}\ln 2+G]
\ee
and
\be
{}_3F_2(\frac12,\frac32,\frac32;\frac52,\frac52;\frac12)
=
\frac{9}{4 \surd 2}[4G-2+\pi(\ln2 -1)]
\ee
where $G$ is Catalan's constant.

\cite{KrupnikovJCAM78}
\be
{}_3F_2(\frac14,\frac14,\frac34;\frac12,\frac54;z)
=
\frac12 z^{-1/4}(\arcsin z^{1/4}+arcsinh z^{1/4}).
\ee

\cite{KrupnikovJCAM78}
\be
{}_3F_2(\frac12,\frac34,\frac54;\frac32,\frac32;z)
=
\frac{2}{\surd z} \ln\frac{1+\sqrt{1+\surd z}}{1+\sqrt{1-\surd z}}.
\ee

\cite{KrupnikovJCAM78}
\be
{}_3F_2(\frac34,\frac34,\frac54;\frac32,\frac74;z)
=
3z^{-3/4}(\arcsin z^{1/4}-\arcsinh z^{1/4}).
\ee

\cite{KrupnikovJCAM78}
\be
{}_3F_2(1,1,1;\frac32,2;-z)
=
\frac{\arcsinh^2\surd z}{z}
\ee

\cite{KrupnikovJCAM78}
\be
_3F_2(1,1,1; (2k+1)/4, 2;1)
= \frac{2k-3}{4}\psi'(\frac{2k-3}{4}),\quad k\ge 2.
\ee

\cite{KrupnikovJCAM78}
\be
_3F_2(\frac16 ,\frac16 ,\frac56;\frac76,\frac76 ;1)
= \frac{\Gamma^2(1/6)}{108\Gamma(1/3)}(\pi \sqrt{3}+6\ln 2)
\ee

\cite{KrupnikovJCAM78}
\be
_3F_2(\frac14 ,\frac14 ,\frac34;\frac12,\frac54 ;1)
= \frac14[\pi+2\ln(1+\surd 2)].
\ee

\cite{KrupnikovJCAM78}
\be
_3F_2(\frac34 ,\frac34 ,\frac54;\frac32,\frac74 ;1)
= \frac32[\pi-2\ln(1+\surd 2)].
\ee

\cite{KrupnikovJCAM78}
\be
_3F_2(\frac14 ,\frac14 ,1;\frac54,\frac54 ;1)
= \frac{1}{16}[\pi+8G].
\ee

\cite{KrupnikovJCAM78}
\be
_3F_2(\frac34 ,\frac34 ,1;\frac74,\frac74 ;1)
= \frac{9}{16}(\pi^2-8G).
\ee

\cite{KrupnikovJCAM78}
\be
_3F_2(\frac13 ,\frac13 ,\frac23;\frac43,\frac43 ;1)
= \frac{1}{18}\Gamma^3(1/3).
\ee

\cite{KrupnikovJCAM78}
\be
_3F_2(\frac13 ,\frac56 ,\frac56;\frac{11}{6},\frac{11}{6} ;1)
= \frac{25}{27}\frac{\pi^{3/2}}{\Gamma(1/3)\Gamma(1/6)}[3\sqrt{3}\ln 3+4\sqrt{3}-3\pi]
\ee

\cite{KrupnikovJCAM78}
\be
_3F_2(\frac12 ,\frac34 ,\frac54;\frac32,\frac32 ;1)
= 2\ln(1+\surd 2).
\ee

\cite{KrupnikovJCAM78}
\be
_3F_2(1/2,1,1; 3/2, 3/2;1)
= 2G.
\ee

\cite{KrupnikovJCAM78}
\be
_3F_2(1,1,1; 5/4, 2;1)
= \frac14(\pi^2+8G).
\ee

\cite{KrupnikovJCAM78}
\be
_3F_2(1,1,1; 7/4, 2;1)
= \frac34(\pi^2-8G).
\ee

\cite{KrupnikovJCAM78}
\be
_3F_2(1,5/4,5/4; 9/4, 9/4;1)
= \frac{25}{16}(\pi^2+8G-16).
\ee

\cite{KrupnikovJCAM78}
\be
_3F_2(1,3/2,3/2; 5/2, 5/2;1)
= \frac{9}{8}(\pi^2-8).
\ee

\cite{KrupnikovJCAM78}
\be
_3F_2(1,3/2,2; 5/2, 5/2;1)
= \frac{9}{2}(2G-1).
\ee

\cite{KrupnikovJCAM78}
\be
_3F_2(1,7/4,7/4; 11/4, 11/4;1)
= \frac{49}{144}(9\pi^2-72G-16).
\ee

\cite{KrupnikovJCAM78}
\be
_3F_2(1,2,2;\frac{1}{4}(2k+5),3;1)
=\frac12 (2k+1)[\frac14 (2k-3)\psi'(\frac{2k-3}{4})-1], k\ge 2
\ee

\cite{KrupnikovJCAM78}
\be
_3F_2(1,2,2; 9/4, 3;1)
= \frac{5}{8}(\pi^2+8G-4).
\ee

\cite{KrupnikovJCAM78}
\be
_3F_2(1,2,2; 11/4, 3;1)
= \frac{7}{8}(3\pi^2-24G-4).
\ee

\cite{KrupnikovJCAM78}
\be
_3F_2(1,1,1; 3/2, 2;-1)
= \ln^2(1+\surd 2).
\ee

\cite{KrupnikovJCAM78}
\be
_3F_2(1,1,3/2; 2, 2;-1)
= 4\ln\frac{1+\surd 2}{2}.
\ee

\cite{KrupnikovJCAM78}
\be
_3F_2(1,3/2,3/2; 5/2, 5/2;-1)
= 9(1-G).
\ee

\cite{KrupnikovJCAM78}
\be
_3F_2(1/2,1,1; 3/2, 3/2;1/4)
= \frac13[8G-\pi\ln(2+\surd 3)].
\ee
\ldots plus special values of $_4F_3$ and $_5F_4$.

\cite{AdamchikZAA21}
\be
2{}_3F_2(\frac12,\frac12,\frac12;1,\frac32;1)=\frac{8G}{\pi}.
\ee

\cite{AdamchikZAA21}
\be
4{}_3F_2(\frac12,\frac12,\frac14;1,\frac54;1)=\frac{\Gamma(1/4)^4}{4\pi^2}.
\ee

\cite{BaileyPGMA2}
\be
(2\sigma-n)(\sigma-n)(\beta+\delta-n-1)(\alpha+\beta-n-1)L_{n+2}
+Q_nL_{n+1}+(n+1)(\sigma-n-1)(\gamma+\delta-n)(\alpha+\gamma-n)L_n=0
\ee
where $L_n={}_3F_2(-n,-\beta,n-2\sigma-1 ; -\beta-\delta,-\alpha-\beta ;1)$,
$Q_n=(2\sigma-2n-1)[(\beta\gamma-\alpha\delta)(\sigma+1)+\theta_3(n+1)(2\sigma-n)]$,
$2\sigma=\alpha+\beta+\gamma+\delta$, $2\theta_3=\alpha-\beta-\gamma+\delta$.
$n$ does not need to be a positive integer.

\cite{BaileyPGMA2}
\begin{multline}
d(d-1)(e-d-1)(a-e)(b-e)(c-e)
{}_3F_2(a,b,c;d-1,e+1;1)
\\
-e(e-1)(d-e-1)(a-d)(b-d)(c-d)
{}_3F_2(a,b,c;d+1,e-1;1)
\\
=
de(d-e)[(d+e-1)(bc+ca+ab)-(2de-d-e+1)(a+b+c)+(d+e)(d-1)(e-1)+de-2abc]
{}_3F_2(a,b,c;d,e;1)
.
\end{multline}

\cite{BaileyPGMA2}
\begin{multline}
a(d+e-a-b-c-1)
{}_3F_2(a+1,b,c;d,e;1)
-(d-a)(e-a)
{}_3F_2(a-1,b,c;d,e;1)
\\
=
[a(2d+2e-2a-b-c-1)-de+bc]
{}_3F_2(a,b,c;d,e;1)
.
\end{multline}

\cite{BaileyPGMA2}
\begin{multline}
e(e-1)(d+e-a-b-c-1)
{}_3F_2(a,b,c;d,e-1;1)
-(a-e)(b-e)(c-e)
{}_3F_2(a,b,c;d,e+1;1)
\\
=
e[(e-1)(d+e-1)-(a+b+c)(2e-1)+bc+ca+ab+e^2]
{}_3F_2(a,b,c;d,e;1)
.
\end{multline}

\cite{BaileyPGMA2}
\begin{multline}
d(d-1)(e-a)
{}_3F_2(a-1,b,c;d-1,e;1)
-a(b-d)(c-d)
{}_3F_2(a+1,b,c;d+1,e;1)
\\
=
d[(d-1)(e-a)-bc+ca+ab-ad]
{}_3F_2(a,b,c;d,e;1)
.
\end{multline}

\cite{ChaundyQJM14}
\be
[\delta(\delta+c-1)(\delta+c'-1)-(\delta+a)
(\delta+a')(\delta-n)]F_n=0,
\ee
where $F_n\equiv \frac{(c)_n}{n!}
\,_3F_2(a,a',-n;c,c';x)$.

\cite{ChaundyQJM14}
\be
n(n+c'-1)F_n -
[2(n-1)^2+(2c+2c'-a-a'-1)(n-1)+cc'-aa']F_{n-1}
+(n+c-2)(n+c+c'-a-a'-2)F_{n-2}=0
\ee
where $F_n\equiv \frac{(c)_n}{n!}
\,_3F_2(a,a',-n;c,c';1)$.

\cite{RoachSSAC96}
\begin{multline*}
_2F_3(-1/2,1;1/4,1/2,3/4;z)
=
1+z^{1/4}\sqrt 2 \sqrt \pi e^{2\sqrt z} \erf (\sqrt 2 z^{1/4})
\\
-z^{1/4}\sqrt 2\sqrt \pi e^{-2\sqrt z}\erfi(\sqrt 2 z^{1/4})
-2\sqrt z \pi \erf (\sqrt 2 z^{1/4})\erfi(\sqrt 2 z^{1/4}).
\end{multline*}

\cite{RoachSSAC96}
\begin{multline*}
_3F_2(-1/2,1,2;3,4;z)
=
-\frac{480+3472z-2100z^2}{525z^3}
+\frac{480+3712z-1024z^2+192z^3}{525z^3}\sqrt{1-z}
\\
-\frac{32}{5z^2}\log\left(\frac{1}{2}+\frac{\sqrt{1-z}}{2}\right)
.
\end{multline*}

\cite{BerndtBLMS15}
\be
_2F_3(-\beta,\beta+\gamma;\gamma,\gamma/2,\frac{1+\gamma}{2};x^2/4)
={}
_1F_1(-\beta;\gamma;-x)
{}_1F_1(-\beta;\gamma;x)
.
\ee

\cite{BerndtBLMS15}
\be
_2F_3(1,n;n+1;(n+1)/2,{2+n}{2};x^2/4)
={}
_1F_1(1;n+1;-x)
{}_1F_1(1;n+1;x)
.
\ee

\be
_{3}F_{2}\left(\begin{array}{c}1,1,1\\ 2,2\end{array}\mid \frac12\right)
=
\frac16 \pi^2-\log^2(2).
\ee

\cite{ChuRDCM}
\be
_3F_2(-n,-1+u,1+v-uv;1+u,-n(1+v)-uv;1)=\frac{u+un}{u+n}\frac{(1+c+vn)_n}{(1+uv+vn)_n}.
\ee

\cite{MaierTAAMS358}
\be
_3F_2(-n,a,b;c,1+a+b-c-n;1)
=
\Gamma\left[
\begin{array}{c}c-a+n,c-b+n,c,c-a-b\\c-a,c-b,c+n,c-a-b+n\end{array}
\right]
\ee
where no lower parameter is a nonpositive integer. 

\cite{MaierTAAMS358}
\be
_3F_2(a,b,c;a-n,1+b;1)
=
\Gamma\left[
\begin{array}{c}1+b-a+n,1-a,1-c,1+b\\1+b-a,1-a+n,1+b-c\end{array}
\right]
\ee
where the parametric excess has positive real part.

\cite{MaierTAAMS358}
If $ab+bc+ca=(d-1)(e-1)$ and $d+e-a-b-c=2$ then
\be
_3F_2(a,b,c;d,e;1)
=
\Gamma\left[
\begin{array}{c}d,e\\a+1,b+1,c+1\end{array}
\right]
.
\ee

\cite{MaierTAAMS358}
If $(a-1)(b-1)=c[(a-1)+(b-1)-(e-1)]$ and $\Re (e-a-b+2)>0$ then
\be
_3F_2(a,b,c;c+2,e;1)
=
\Gamma\left[
\begin{array}{c}e,e-a-b+2,c+2\\ e-a+1,e-b+1,c+1\end{array}
\right]
.
\ee

\cite{MaierTAAMS358}
If $(a-1)(b-1)=(d-2)(e-2)$ and $\Re (d+e-a-b-2)>0$ then
\be
_3F_2(a,b,2;d,e;1)
=
\frac{(d-1)(e-1)}{d+e-a-b-2}
.
\ee

\cite{Milgramarxiv1105}
\begin{multline}
_3F_2(a,m,b;c,m-n;1)
=
\Gamma(a+n-m+1)
\\
\times
\left(\sum_{L=0}^{m-1}
\frac{(-1)^L\Gamma(-b+c-a-n+L)\Gamma(1-b+L)}{\Gamma(m-L)\Gamma(-b-n+1+L)
\Gamma(-b+c-m+1+L)\Gamma(L+1)}
\right)
\\
\times
\frac{\Gamma(c)\Gamma(m-n)(-1)^n}{\Gamma(c-a)\Gamma(a)}
,\quad
0<n<m.
\end{multline}

\begin{multline}
_3F_2(a,b,-n;c,m-n;1)
=
\Gamma(a+n-m+1)
\\
\times
\left(\sum_{L=0}^{m-1}
\frac{(-1)^L\Gamma(c-a+L)\Gamma(1+b-m+n+L)}{\Gamma(m-L)\Gamma(b-m+1+L)
\Gamma(c+n-m+1+L)\Gamma(L+1)}
\right)
\\
\times
\frac{\Gamma(c)\Gamma(m-n)(-1)^n}{\Gamma(c-a)\Gamma(a)}
,\quad
0<n<m.
\end{multline}

\begin{multline}
_3F_2(a,-n,b;c,a-c+b-n+m;1)
=
\Gamma(c-b+n-m+1)
\Gamma(-a+c-b-m+1)
\\
\times
\left(\sum_{L=0}^{m-1}
\frac{(-1)^L\Gamma(b+L)\Gamma(1-a+c-m+n+L)}{\Gamma(m-L)\Gamma(-a+c-m+1+L)
\Gamma(c+n-m+1+L)\Gamma(L+1)}
\right)
\\
\times
\frac{\Gamma(c)\Gamma(m)}{\Gamma(n+1-a+c-b-m)\Gamma(b)\Gamma(c-b)}
,\quad
0<n,m.
\end{multline}

\begin{multline}
_3F_2(a,1,b;n+1,c;1)
=
\frac{\Gamma(a-n)\Gamma(b-n)\Gamma(n+1)\Gamma(c)\Gamma(-b+c-a+n)\Gamma(1-b)\Gamma(1+n)\Gamma(c)}
{\Gamma(b)\Gamma(c-b)\Gamma(a)\Gamma(c-a)\Gamma(a)}
\\
\times
\left(\sum_{L=0}^{n-1}
\frac{(-1)^L\Gamma(a-1-L)}{\Gamma(c-1-L)\Gamma(n-L)\Gamma(2+L-b)}
\right)
\\
,\quad
0<n,m.
\end{multline}
and 60 others.

\cite{LievensJCAM169,GasperCM254}
\begin{multline}
_3F_2(\alpha,\beta,-n ; \gamma,\delta;1)
=
\sum_{l=0}^n \binom{n}{l}
\frac{(\mu)_l (\lambda+\delta-\alpha-\beta)_l(\gamma-\mu)_{n-l}}
{(\gamma)_n(\delta)_l}
\\
\times
{}_3F_2(\lambda-\alpha,\lambda-\beta,-l ;
\mu,\lambda+\delta-\alpha-\beta;1)
\,{}_3F_2(\alpha+\beta-\lambda,\lambda-\mu,l-n;
\gamma-\mu,l+\delta;1)
\end{multline}
where $\Re \lambda>0$,
$\Re \nu>0$,
$\Re(\gamma+\mu-\lambda-\nu)>0$.

\cite{LievensJCAM169,GasperCM254}
\begin{multline}
_3F_2(\alpha,\beta,-n ; \gamma,\delta;1)
=
\sum_{l=0}^n \binom{n}{l}
\frac{(\lambda)_l (\beta)_l(\mu)_l}
{(\delta)_l(\mu)_{2l}}
\\
\times
{}_3F_2(\alpha,l+\mu,-l ; \gamma,\lambda ;1)
\,{}_3F_2(l+\beta,l+\lambda,l-n ; 1+2l+\mu,l+\delta;1) ;
\end{multline}

\cite{BerndtBLMS15,BaileyPLMS28}
\be
_4F_1(-\alpha,-\beta,-\frac{\alpha+\beta}{2}-\frac{\alpha+\beta-1}{2};-\alpha-\beta;4x^2)
=
{}_2F_0(-\alpha,-\beta;x)
{}_2F_0(-\alpha,-\beta;-x)
.
\ee
if $\alpha$ or $\beta$ a nonnegative integer.

\cite{MilgramANE18}
\be
_4F_1(3/2,1,1,\epsilon;2;t)=\frac{2}{\surd \pi \Gamma(\epsilon)}G^{4,1}_{2,4}\left(-1/t\mid
\begin{array}{l} 1;2\\ 3/2,1,1,\epsilon 
\end{array}
\right).
\ee

\cite[(4.3.5.1)]{SlaterHyp}\cite{WhipplePLMS25}\cite{KarlssonPA87}
\be
_4F_3(x,y,z,-n;u,v,w;1)
=
\frac{(v-z)_n(w-z)_n}{(v)_n(w)_n}
\,_4F_3(u-x,u-y,z,-n;1-v+z-n,1-w+z-n,w;1) ,
\ee
if $u+v+w=1+x+y+z-n$.

\cite{BerndtBLMS15}
\be
_{4}F_{3}\left(\begin{array}{c}1,1,1,1\\ 2,2,2\end{array}\mid \frac12\right)
=
\frac74 \zeta(3)-\frac16\pi^2 \log 2 +\frac13\log^3(s).
\ee

\cite{BaileyPLMS28}
\be
_4F_3(1-n-c,1-n-a-b,-n/2,(1-n)/2; 1-n-(a+c)/2, 3/2-n-(a+c)/2,1-n-b ;1)
=\frac{(c)_n}{(a+c+n-1)_n}{}_3F_2(a,b,-n;c,1-n-b;1) .
\ee

\cite{BaileyPLMS28}
\be
_4F_3(a,a+1/2-n/2,(1-n)/2 ; \rho+1/2,\rho/2,(p+1)/2;1)
=\frac{(p-2a)_n}{(p)_n}{}_3F_2(2a,\rho,-n; 2\rho,1+2a-n-p;2).
\ee

\cite{GottschalkJPA21}
\begin{multline}
\frac{\Gamma(a+1/2)n2^{2b-n}}{\Gamma(2b+1-n)\surd \pi}
{}_4F_3(1,a+1/2,1-n/2,1/2-n/2; 1/2+b-n/2,1+b-n/2,3/2 ;1)
\\
=\sum_{p=1,3,5,\ldots \mathrm{odd}}^n\binom{n}{p}\frac{\Gamma(1/2+p/2)\Gamma(a+p/2)}{\Gamma(b-n/2+p/2)\Gamma(b-n/2+p/2+1/2)}
.
\end{multline}

\cite{KarlssonQJM22}
\begin{multline}
_4F_3(a,b,c+d+1,d;a+d+1,b+d+1,c;1)
=\frac{\Gamma(a+d+1)\Gamma(b+d+1)\Gamma(c+1)}{\Gamma(a+1)\Gamma(b+1)\Gamma(c+d+1)\Gamma(d+1)}
\\
+
\frac{d(d+1)(a-c)(b-c)}{(a+d+1)(b+d+1)c(c+1)}
{}_4F_3(a+1,b+1,c+d+1,d+2;a+d+2,b+d+2,c+2;1)
.
\end{multline}

\cite{LakinQJM21}
\begin{multline}
_4F_3(a,b,c,d;e,c+d-e,a+b+1;1)
\\
=\frac{(c-e)(d-e)}{e(e-c-d)}
{}
_4F_3(a+1,b+1,c,d;e+1,c+d-e+1,a+b+1;1)
,
\end{multline}
where one of the numerator parameters is a negative integer.

\cite{BaileyPLMS28}
\be
_4F_3(2b,b+1,b-a-1/2,-n; b,b+a+3/2,1-2a+2b-n;1)
=\frac{(2a)_n (a+1)_n (a-b-1/2)_n}{(2a-2b)_n(a+b+3/2)_n(a)_n}.
\ee

\cite{ExtonJCAM83}
\be
_4F_3(a,1+a/2,b,-n ; a/2,1+a-b,1+2b-n;1) = \frac{(1-2b)_n(-b)_n}{(1+a-b)_n(-2b)_n}.
\ee

\cite{ExtonJCAM83}
\be
(1-y)^{-a} {}_4F_3(a+d-1,-a,2-d-2a,d+2a-1;d,d+a,1-d-2a;y)
={}_2F_1(a+d-1,-2a;d;y).
\ee

\cite{AbiodunJMAA70}
\begin{multline}
_4F_3(\gamma,2\beta-\gamma,1+\alpha/2,1/2+\alpha/2; \beta+1/2, 1+\alpha, 1+\beta;1)
=
\frac{\beta\Gamma(2\beta)\Gamma(2\beta-\alpha-\gamma)}{(\beta-\gamma)\Gamma(2\eta-\gamma)\Gamma(2\beta-\alpha)}
\end{multline}
where $\Re(2\beta-\alpha-\gamma)>0$.

\cite{LakinQJM21}
\begin{multline}
_4F_3(a+1,b,c,d;a-d,b+d+1,c+d+1;1)
\\
=\frac{(b+d-a)(c+d-a)d(d+1)}{(a-d)(1+a-d)(1+b+d)(1+c+d)}
{}
_4F_3(a+1,b+1,c+1,d+2;2+a-d,2+b+d,2+c+d;1)
,
\end{multline}
where one of the numerator parameters is a negative integer.

\cite{LakinQJM21}
\begin{multline}
_4F_3(a+1,b,c,d;a+c,b+c+1,d-c+1;1)
\\
=\frac{(b+c)(d-c-a)}{(a+c)(d-c-b)}
{}
_4F_3(a,b+1,c,d;a+c+1,b+c,d-c+1;1)
,
\end{multline}
where one of the numerator parameters is a negative integer.

\cite{LakinQJM21}
\begin{multline}
_4F_3(a+1,b,c,a+b+c;a+b+1,b+c+1,c+a+1;1)
\\
=-\frac{(b+c)}{a}
{}
_4F_3(a,b,c,a+b+c;a+b+1,b+c,c+a+1;1)
,
\end{multline}
where one of the numerator parameters is a negative integer.

\cite{WilletQJM2}
\begin{multline}
_4F_3(-\alpha,-\beta,1-\gamma-n,-n; 1-\alpha-n,1-\beta-n,-\gamma;1)
= -\frac{n(n-1)(\alpha-\gamma)(\beta-\gamma)}{(\alpha+n-1)(\beta+n-1)\gamma(1-\gamma)}\\ \times
_4F_3(1-\alpha,1-\beta,1-\gamma-n,2-n; 2-\alpha-n,2-\beta-n,2-\gamma;1).
\end{multline}

\cite{DriverETNA25}
\begin{multline}
_4F_3(a,b/3,(b+1)/3,(b+2)/3;c/3,(c+1)/3,(c+2)/3;1)
\\
=\frac{\Gamma(c)\Gamma(c-b-a)}{\Gamma(c-a)\Gamma(c-b)}
\sum_{k=0}^\infty \frac{(a)_k(-1)^k(b)_k}{k!(c-a)_k}\,_2F_1(-k,b+k;c-a+k;-1)
.
\end{multline}

\cite{AdamchikZAA21}
\be
\frac{1}{r}{}_4F_3(\frac12,\frac12,r;1,r+1;1)=\sum_{k=0}^\infty \frac{\binom{2k}{k}^2}{(k+r)16^k}.
\ee

\cite{AdamchikZAA21}
\begin{multline}
(n-\frac12){}_4F_3(1,1,n+1/2,n+1/2;2,n+1,n+1;1)
\\
=\frac{4n^2}{2n-1}(H_{n-1}+\log 4)-\frac{16^n}{\binom{2n}{n}^2}{}_3F_2(1/2,1/2,n-1/2;1,n+1/2;1).
\end{multline}

\cite{AdamchikZAA21}
\begin{multline}
\frac{2n^2}{(2n-1)^2}{}_4F_3(1,1,1-n,1-n;2,3/2-n,3/2-n;1)
\\
=\psi(n+1/2)-\psi(1/2)=\sum_{k=0}^{n-1}\frac{2}{2k+1}.
\end{multline}

\cite{KarlssonPA87}
\be
_4F_3(\frac{\alpha}{2},\frac{\alpha+1}{2},\beta+n,-n;1+\alpha,\frac{\beta}{2},
\frac{\beta+1}{2};1)
=\frac{(\beta-\alpha)_n}{(\beta)_n}
.
\ee

\cite{KarlssonPA87}
\be
_4F_3(\alpha,-\alpha,-\frac{m}{2},\frac{1-m}{2};\frac{1}{2},\beta,1-m-\beta;1)
=\frac{(\alpha+\beta)_m+(\beta-\alpha)_m}{2(\beta)_m}
.
\ee

\cite{KarlssonPA87}
\be
_4F_3(a,b,\frac{1}{2}-a-b-n,-n; a+b-\frac{1}{2},1-a-n,1-b-n;1)
=\frac{(2a)_n(2b)_n(a+b)_n}{ (2a+2b-1)_n(a)_n(b)_n}
.
\ee

\cite{KarlssonPA87}
\be
_4F_3(a-\frac{1}{2},b-\frac{1}{2},\frac{1}{2}-a-b-n,-n; a+b-\frac{1}{2},\frac{1}{2}-a-n,\frac{1}{2}-b-n;1)
=\frac{(2a)_n(2b)_n(a+b)_n}{ (2a+2b-1)_n(a+\frac{1}{2})_n(b+\frac{1}{2})_n}
.
\ee

\cite{KarlssonPA87}
\be
_4F_3(a+1,b,\frac{1}{2}-a-b-n,-n; a+b+\frac{1}{2},1-a-n,1-b-n;1)
=\frac{(2a+1)_n(2b)_n(a+b)_n}{ (2a+2b)_n(a)_n(b)_n}
.
\ee

\cite{KarlssonPA87}
\be
_4F_3(a+\frac{1}{2},b-\frac{1}{2},\frac{1}{2}-a-b-n,-n; a+b+\frac{1}{2},\frac{1}{2}-a-n,\frac{1}{2}-b-n;1)
=\frac{(2a+1)_n(2b)_n(a+b)_n}{ (2a+2b)_n(a+\frac{1}{2})_n(b+\frac{1}{2})_n}
.
\ee

\cite{KarlssonPA87}
\be
_4F_3(a,b,\frac{1}{2}-a-b-n,-n; a+b+\frac{1}{2},1-a-n,1-b-n;1)
=\frac{(2a)_n(2b)_n(a+b)_n}{ (2a+2b)_n(a)_n(b)_n}
.
\ee

\cite{KarlssonPA87}
\be
_4F_3(a,b,a+b-\frac{1}{2}-n,-n; a+b+\frac{1}{2},\frac{1}{2}+a-n,\frac{1}{2}+b-n;1)
=\frac{(1/2)_n(a-b+\frac{1}{2})_n(b-a+\frac{1}{2})_n}{ (a+b+\frac{1}{2})_n(\frac{1}{2}-a)_n(\frac{1}{2}-b)_n}
.
\ee

\cite{KarlssonPA87}
\be
_4F_3(a,-a,\frac{1}{2}-n,-n; \frac{1}{2},\frac{1}{2}+b-n,\frac{1}{2}-b-n;1)
=\frac{(\frac{1}{2}+a+b)_n(\frac{1}{2}-a-b)_n
+(\frac{1}{2}+a-b)_n(\frac{1}{2}-a+b)_n}
{2(\frac{1}{2}+b)_n(\frac{1}{2}-b)_n}
.
\ee

\cite{KarlssonPA87}
\be
_4F_3(a,-a,-\frac{1}{2}-n,-n; \frac{1}{2},b-n,-b-n;1)
=\frac{(a+b)_{n+1}(1-a-b)_n
+(b-a)_{n+1}(1+a-b)_n}
{2b(1+b)_n(1-b)_n}
.
\ee

\cite{KarlssonPA87}
\be
_4F_3(\frac{1}{2}+a,\frac{1}{2}-a,-\frac{1}{2}-n,-n; \frac{1}{2},\frac{1}{2}+b-n,\frac{1}{2}-b-n;1)
=\frac{(a+b)_{n+1}(1-a-b)_n
+(b-a)_{n+1}(1+a-b)_n}
{2b(\frac{1}{2}+b)_n(\frac{1}{2}-b)_n}
.
\ee

\cite{KarlssonPA87}
\begin{multline}
_4F_3(a,-a,-\frac{1}{2}-n,-n; \frac{1}{2},\frac{1}{2}+b-n,\frac{1}{2}-b-n;1)
\\
=\frac{(a+b)(\frac{1}{2}+a+b)_{n}(\frac{1}{2}-a-b)_n
+(b-a)(\frac{1}{2}+b-a)_{n}(\frac{1}{2}+a-b)_n}
{2b(\frac{1}{2}+b)_n(\frac{1}{2}-b)_n}
.
\end{multline}

\cite{KarlssonPA87}
\be
_4F_3(a,-a,\frac{1}{2}-n,-n; \frac{1}{2},1+b-n,1-b-n;1)
=\frac{(a+b)_{n}(-a-b)_n
+(a-b)_{n}(b-a)_n}
{2(b)_n(-b)_n}
.
\ee

\cite{BerndtBLMS15}
\be
_4F_3(\alpha,\beta,\frac{\alpha+\beta}{2},\frac{\gamma+\delta}{2};
\gamma,\delta,\alpha+\beta;x)
=
{}
_2F_1(\alpha,\beta;\gamma;\frac{1-\sqrt{1-x}}{2})
_2F_1(\alpha,\beta;\delta;\frac{1-\sqrt{1-x}}{2})
\ee
if $\alpha+\beta+1=\gamma+\delta$.

\cite{BuhringJAMSA8}
\begin{multline}
\frac{1}{\Gamma(b_1)\Gamma(b_2)\Gamma(b_3)}
\,_4F_3\left(
\begin{array}{c}
a_1,a_2,a_3,a_4\\
b_1,b_2,b_3
\end{array}\mid 1
\right)
=\frac{\Gamma(s)}{\Gamma(a_1+s)\Gamma(a_2+s)
\Gamma(a_3)\Gamma(a_4)}
\\
\times
\sum_{k=0}^\infty
\frac{(b_1+b_3-a_3-a_4)_k
(b_2+b_3-a_3-a_4)_k(s)_k}{
(a_1+s)_k(a_2+s)_k k!}
\,
_3F_2
\left(
\begin{array}{c}
b_3-a_3,b_3-a_4,-k\\
b_1+b_3-a_3-a_4,b_2+b_3-a_3-a_4
\end{array}
\mid 1 \right)
,
\end{multline}
where $s=b_1+b_2+b_3-a_1-a_2-a_3-a_4>0$.

\cite{ExtonJCAM79}
\begin{multline}
_4F_3(a/2,a/2+1/2,b/2,b/2+1/2;1/2+a/2-b/2,1+a/2-b/2,1/2;1)
\\
=\frac{\Gamma(1+a-b)\Gamma(1-2b)}{\Gamma(1-b)\Gamma(1+a-2b)}
+\frac{\Gamma(1+a-b)\Gamma(1+a/2)}{\Gamma(1+a)\Gamma(1+a/2-b)}.
\end{multline}
where $\Re b < 1/2$.

\cite{ExtonJCAM79}
\begin{multline}
\frac{2ab}{1+a-b}\,_4F_3(a/2+1/2,a/2+1,b/2+1/2,b/2+1;1+a/2-b/2,
3/2+a/2-b/2,3/2;1)
\\
=\frac{\Gamma(1+a-b)\Gamma(1-2b)}{\Gamma(1-b)\Gamma(1+a-2b)}
+\frac{\Gamma(1+a-b)\Gamma(1+a/2)}{\Gamma(1+a)\Gamma(1+a/2-b)}.
\end{multline}
where $\Re b < 1/2$.

\cite{LavoieMC49}
\begin{multline}
_4F_3(1+a,1+b,1,1; \frac{a+b+3}{2},2,2;1)
=\frac{1+a+b}{2ab}[\psi(1)+\psi(\frac{1-a-b}{2})-\psi(1-a)-\psi(1-b)\\
+\frac{\cos \pi(a-b)/2}{2\cos\pi(a+b)/2}(\psi(1-a/2)+\psi(1-b/2)-\psi(\frac{1-a}{2})-\psi(\frac{1-b}{2}))
]
\end{multline}
where $\Re(a+b)<3$.

\cite{LavoieMC49}
\begin{multline}
_4F_3(1+a,1+c,1,1; \frac{a+3}{2},2c,2;1)
=\frac{(1+a)(2c-1)}{4ac}[\psi(\frac{a+1}{2})-\psi(1/2)+\psi(c-1/2)-\psi(c-a/2-1/2)]\\
+\frac{\Gamma(1/2)\Gamma(\frac{a+3}{2})\Gamma(c-a/2-1/2)\Gamma(c+1/2)}{2\Gamma(1+a/2)\Gamma(c-a/2)\Gamma(c+1)}
\end{multline}
where $\Re(2c-a)>1$.

\cite{LavoieMC49}
\begin{multline}
_4F_3(1+a,1+c,1,1; \frac{a+3}{2},2+2c,2;1)
=\frac{(1+a)(2c+1)}{4ac}[\psi(\frac{a+1}{2})-\psi(1/2)+\psi(c+1/2)-\psi(c-a/2+1/2)]\\
-\frac{\Gamma(1/2)\Gamma(\frac{a+3}{2})\Gamma(c-a/2+1/2)\Gamma(c+3/2)}{2c\Gamma(1+a/2)\Gamma(c-a/2+1)\Gamma(c+1)}
\end{multline}
where $\Re(2c-a)>-3$.

\cite{LavoieMC49}
\begin{multline}
_4F_3(1+a,1+c,1,1; \frac{a+3}{2},1+2c,2;1)
=\frac{1+a}{2a}[\psi(\frac{a+1}{2})-\psi(1/2)+\psi(c+1/2)-\psi(c-a/2+1/2)]
\end{multline}
where $\Re(2c-a)>-1$.

\cite{LavoieMC49}
\begin{multline}
_4F_3(1+c,3,1,1; 1+e,1+f,2;1)
=\frac{ef}{2(c-1)}\big\{1-\frac{(e-1)(f-1)-c+1}{2c}\\
[\psi(\frac{e+1}{2})-\psi(e/2)+\psi(\frac{f+1}{2})-\psi(f/2)]
\big\}
\end{multline}
where $f=1+2c-e$, $\Re c>1$.

\cite{ChuRDCM}
\be
_4F_3(-n,1-\lambda,1-\mu,-a-c ;-\lambda,-\mu,1-c+bn-n;1)
=\frac{a(a-1)}{(a+bn)(a+bn-1)} \frac{(-a-bn)_n}{(c-bn)_n}.
\ee
where $\lambda$ and $\mu$ are the zeros of the quadratic
polynomial $(a+b-1)(a+c-bx)(a+c-bx-1)+(b-1)b(c-bx)x$.

\cite{ChuRDCM}
\begin{multline}
_4F_3(-n,1-\frac{(a+b-1)(a+c)}{ab},1-c/b,-a-c
; -\frac{(a+b-1)(a+c)}{ab},-c/b,1-c+bn-n;1)
\\
=\frac{a(a-1)(c-bn)}{c(a+bn)(a+bn-1)} \frac{(-a-bn)_n}{(c-bn)_n}.
\end{multline}

\cite{ChuRDCM}
\begin{multline}
_4F_3(-n,1-\frac{(a+b-1)(1+c)}{ab},1-c/b,-a-c
; -\frac{(a+b-1)(1+c)}{ab},-c/b,1-c+bn-n;1)
\\
=\frac{a(a-1)}{c(c+1)} \frac{(2-a-bn)_{n-2}}{(2+c-bn)_{n-2}}.
\end{multline}

\cite{ChuRDCM}
\begin{multline}
_4F_3(-n,1+\frac{(a+b-1)(1+c)}{b(b-c-2)},-1+\frac{1-a}{b},-a-c
; -\frac{(a+b-1)(1+c)}{b(b-c-2)},1+\frac{1-a}{b},1-a+bn-n;1)
\\
=\frac{1-a}{c+1} \frac{(-1-c-bn)_{n+1}}{(-1+a-bn)_{n+1}}.
\end{multline}

\cite{ChuRDCM}
\begin{multline}
_5F_4(-n,1-\alpha,1-\beta,-1+\frac{1-a}{b},-a-c
; -\alpha,-\beta, 1+\frac{1-a}{b},1-c+bn-n;1)
\\
=\frac{a-1}{a-bn-1} \frac{(-c-bn)_n}{(a-bn)_n}.
\end{multline}
where $\alpha$ and $\beta$ are the zeros of the quadratic
polynomial $(a+c)(b-1)+(a-bx)(a+c-bx)x$.

\cite{FieldsJAT1}
\be
_4F_3(-n,\beta+1,1,z+2\beta; n+2\beta+1,\beta,1-z ;1)
=\frac{z(n+2\beta)}{2\beta(z-n)},\quad n=0,1,2,\ldots ; \beta(z-n)\neq 0.
\ee

\cite{GottschalkJPA21}
\begin{multline}
_4F_3(1,a,1/2-n/2,-n/2; 1/4+a/2-n/2,3/4+a/2-n/2,3/2;1)
\\
=
\frac{\Gamma(n+a+1/2)}{(n+1)\Gamma(n-a+1/2)}\big\{
\frac{(-)^{n+a}\pi 4^{a-1}}{\Gamma(2a)}
+\sum_{j=1}^a\frac{(-)^j}{\Gamma(2a+1-j)\Gamma(j)}
\\ \times
\left[
\psi(a/2-j/2+n/2+5/4)
-\psi(j/2+n/2-a/2+3/4)
\right]
\big\}
\end{multline}
where $a$ and $n$ are positive integers.

\cite{LievensJCAM169}
\begin{multline}
_4F_3(-n,a,b,c; d,e,f;1) = \frac{(a)_n(e+f-a-b)_n(e+f-a-c)_n}{(e)_n(f)_n(e+f-a-b-c)_n}{}
\\
\times
{}_4F_3(-n,e-a,f-a,e+f-a-b-c; e+f-a-b,e+f-a-c,1-n-a ;1)
\end{multline}
provided $d+e+f=1-n+a+b+c$.

\cite{KarpMdpiM8}
\begin{equation}
_4F_3(a,b,c,f+1;d,e,f;1) = \frac{\Gamma(e)\Gamma(s)}{\Gamma(e-c)\Gamma(s+c)}
{}
_4F_3(d-a-1,d-b-1,c,\eta+1;d,d+e-a-b-1,\eta;1)
\end{equation}
where $s=d+e-a-b-c-1$ is the excess and $\eta=(d-a-1)(d-b-1)f/[ab+(d-a-b-1)f]$.

\cite{KarpMdpiM8}
\begin{equation}
_4F_3(a,b,c,f+1;d,e,f;1) = \frac{[sf-c(d-a-b)]\Gamma(e)\Gamma(s)}{f\Gamma(e+d-a-b)\Gamma(e-c)}
{}
_4F_3(d-a,d-b,c,\xi+1;d,d+e-a-b,\xi;1)
\end{equation}
where $s$ is the excess and $\xi=f+(d-a-b)(f-c)/(e-c-1)$.

\cite{KarpMdpiM8}
\begin{equation}
_4F_3(a,b,c,f+1;d,e,f;1) = \frac{(fs+bc)\Gamma(s)\Gamma(d)\Gamma(e)}{f\Gamma(a)\Gamma(s+b+1)\Gamma(s+c+1)}
{}
_4F_3(s,d-a,e-a,\zeta+1; d+e-a-c,d+e-a-b,\zeta;1)
\end{equation}
where $\zeta=s+bc/f$ and $s$ the excess.

\cite{KarpMdpiM8}
\begin{equation}
_4F_3(a,b,c,f+1;d,e,f;1) = \frac{(abc+fds)\Gamma(s)\Gamma(e)}{fd\Gamma(e-a)\Gamma(s+a+1)}
{}
_4F_3(a,d-b,d-c,\nu+1;d+1,s+a+1,\nu;1)
\end{equation}
where $\nu=(abc+fds)/(bc+fs)$ and $s$ the excess.

\cite{LievensJCAM169,GasperCM254}
\begin{multline}
_4F_3(\alpha,\beta,n+\nu,-n ; \gamma,\delta,\epsilon;1)
=
\frac{(\gamma-\mu)_n(1+\nu-\gamma)_n}{(\gamma)_n (1+\nu+\mu-\gamma)_n}
\\
\times \sum_{l=0}^n
\frac{(\nu+\mu-\gamma)_l}{l!}
\,
\frac{(\mu)_l(\lambda+\delta-\alpha-\beta)_l(\lambda+\epsilon-\alpha-\beta)_l(n+\nu)_l
(-n)_l(1+\nu+\mu-\gamma)_{2l}}
{(1+\nu-\gamma)_l(\epsilon)_l(\delta)_l (1+\mu-\gamma-n)_l(n+1+\nu+\mu-\gamma)_l
(\nu+\mu-\gamma)_{2l}}
\\
{}_4F_3(\lambda-\alpha,\lambda-\beta,l+\nu+\mu-\gamma,-l ; \mu,
\lambda+\delta-\alpha-\beta,\lambda+\epsilon-\alpha-\beta;1)
\,{}_4F_3(\alpha+\beta-\lambda,\lambda-\mu,n+l+\nu,l-n ;
\gamma-\mu,l+\delta,l+\epsilon;1)
\end{multline}
when $\alpha+\beta+\nu+1=\gamma+\delta+\epsilon$.

\cite{LievensJCAM169}
\begin{multline}
_4F_3(\alpha,\beta,n+\nu,-n ; \gamma,\delta,\epsilon;1)
=
\sum_{l=0}^n (-)^l
\frac{(-n)_l}{l!} \frac{(\lambda)_l (1+\alpha+\beta+\mu-\gamma-\lambda)_l
(n+\nu)_l(\mu)_l}{(\delta)_l(\epsilon)_l(\mu)_{2l}}
\\
\times
{}_4F_3(\alpha,\beta,l+\mu,-l;
\gamma,\lambda,1+\alpha+\beta+\mu-\gamma-\lambda;1)
\,{}_4F_3(l+n+\nu,l+\lambda,1+l+\alpha+\beta+\mu-\gamma-\lambda,l-n;
1+2l+\mu,l+\delta,l+\epsilon;1)
\end{multline}
when $\alpha+\beta+\nu+1=\gamma+\delta+\epsilon$, $n$ a nonnegative integer..

\cite{GottschalkJPA21} 
\begin{multline}
_4F_3(1,a,1/2-n/2,-n/2; 3/4+a/2-n/2,5/4+a/2-n/2,3/2;1)
\\
=
\frac{\Gamma(n+a-1/2)}{(n+1)\Gamma(n-a-1/2)}\big\{
-\frac{(-)^{n+a}\pi 4^{a-1}}{\Gamma(2a)}
-\frac{\pi(-1)^{n+a}}{\Gamma^2(a+1/2)}
+\sum_{j=1}^a\frac{(-)^j}{\Gamma(2a+1-j)\Gamma(j)}
\\ \times
\left[
\psi(a/2-j/2+n/2+3/4)
-\psi(j/2-a/2+n/2+1/4)
\right]
\\
-2a\sum_{j=1}^{2a+2}\frac{(-)^j(n+a+3/2-j)}{\Gamma(2a+3-j)\Gamma(j)}
\\ \times
\left[
(-)^{n+1}\psi(a/2+5/4-j/2)
-\psi(a/2+n/2-j/2+7/4)
\right]
\big\}
\end{multline}
where $a$ and $n$ are positive integers. 

\cite{KimarXiv1205}
\begin{multline}
_4F_3(a,b,c,d+1;(a+b+1)/2, 2c+1,d;1)
=\frac{2^{a+b-2}\Gamma(c+1/2)\Gamma(a/2+b/2+1/2)\Gamma(c-a/2-b/2+1/2)}{\Gamma(1/2)\Gamma(a)\Gamma(b)}
\\ \times
[\frac{\Gamma(a/2)\Gamma(b/2)}{\Gamma(c-a/2+1/2)\Gamma(c-b/2+1/2)}+\frac{2c-d}{d}\,
\frac{\Gamma(a/2+1/2)\Gamma(b/2+1/2)}{\Gamma(c-a/2+1)\Gamma(c-b/2+1)}]
\end{multline}
if $\Re d>0$, $\Re(2c-a-b)>-1$.

\cite{MillerJPA36}
\begin{multline}
_4F_3(a,b,c,d;e,f,g;1)
=\Gamma(e)\Gamma(f)\Gamma(g)\Gamma(1-d)
\\
\times
\big\{
\frac{\Gamma(b-a)\Gamma(c-a)}{\Gamma(b)\Gamma(c)\Gamma(e-a)\Gamma(f-a)\Gamma(g-a)\Gamma(1+a-d)}{}_4F_3(a,1+a-e,1+a-f,1+a-g; 1+a-b,1+a-c,1+a-d;1)
\\
+\frac{\Gamma(a-b)\Gamma(c-b)}{\Gamma(a)\Gamma(c)\Gamma(e-b)\Gamma(f-b)\Gamma(g-b)\Gamma(1+b-d)}{}_4F_3(b,1+b-e,1+b-f,1+b-g; 1+b-a,1+b-c,1+b-d;1)
\\
+\frac{\Gamma(a-c)\Gamma(b-c)}{\Gamma(a)\Gamma(b)\Gamma(e-c)\Gamma(f-c)\Gamma(g-c)\Gamma(1+c-d)}{}_4F_3(c,1+c-e,1+c-f,1+c-g; 1+c-a,1+c-b,1+c-d;1)
\big\}
.
\end{multline}

\cite{MillerJPA36}
\begin{multline}
_4F_3(a,b,c,d;e,f,g;1)
=\Gamma(e)\Gamma(f)\frac{\Gamma(1-c)\Gamma(1-d)}{\Gamma(1-g)}
\\
\times
\big\{
\frac{\Gamma(b-a)\Gamma(1+a-g)}{\Gamma(b)\Gamma(e-a)\Gamma(f-a)\Gamma(1+a-c)\Gamma(1+a-d)}{}_4F_3(a,1+a-e,1+a-f,1+a-g; 1+a-b,1+a-c,1+a-d;1)
\\
+\frac{\Gamma(a-b)\Gamma(1+b-g)}{\Gamma(a)\Gamma(e-b)\Gamma(f-b)\Gamma(1+b-c)\Gamma(1+b-d)}{}_4F_3(b,1+b-e,1+b-f,1+b-g; 1+b-a,1+b-c,1+b-d;1)
\\
-\frac{\Gamma(g-1)\Gamma(1+a-g)\Gamma(1+b-g)}{\Gamma(a)\Gamma(b)\Gamma(g-c)\Gamma(1+e-g)\Gamma(1+f-g)}{}_4F_3(1+a-g,1+b-g,1+c-g,1+d-g; 2-g,1+e-g,1+f-g;1)
\big\}
.
\end{multline}

\cite{Saigoarxiv98}
\begin{multline}
\Gamma\left[\begin{array}{c}
a_1,a_1+m,a_3,a_4\\ b_1,b_2,b_3
\end{array}\right]
{}_4F_3(a_1,a_1+m,a_3,a_4; b_1,b_2,b_3; z) = (-z)^{-a_1} \sum_{\nu=0}^{m-1}
\frac{z^{-\nu}}{\nu!} 
\\
\times
\Gamma\left[\begin{array}{c}
m-\nu,a_3-a_1-\nu,a_4-a_1-\nu,a_1+\nu \\ b_1-a_1-\nu,b_2-a_2-\nu,b_3-a_1-\nu
\end{array}\right]
z^{-\nu}
+(-z)^{-a_1}\sum_{\nu=m}^\infty \frac{(-1)^m(-z)^{-\nu}}{\nu! (\nu-m)!} 
\\
\times
\Gamma\left[\begin{array}{c}
a_1+\nu,a_3-a_1-\nu,a_4-a_1-\nu \\ b_1-a_1-\nu,b_2-a_1-\nu,b_3-a_1-\nu
\end{array}\right]
[A+\log(-z)]
\end{multline}
where $m$ is a non-negative integer $a_1\neq 0,-1,-2,\ldots$; $|\arg(-z)|<\pi$,
$A=\psi(\nu+1)
+\psi(\nu+1-m)
+\psi(a_3-a_1-\nu)
+\psi(a_4-a_1-\nu)
-\psi(a_1+\nu)
-\psi(b_1-a_1-\nu)
-\psi(b_2-a_1-\nu)
-\psi(b_3-a_1-\nu)$, $\psi\equiv d \log\Gamma(z)/dz$.
Similar formulas where not 2 but even 3 or 4 of the upper parameters differ by integers.

\cite{KarpMdpiM8}
\be
(\gamma+1){}_4F_3(a-1,b,c,\xi+1;d,e,\xi;1)
=
{}_3F_2(a,b,c;d,e;1)
+\gamma
{}_3F_2(a-1,b,c;d,e;1)
\ee
where $\xi=(\gamma+1)(a-1)$.

\cite{KarpMdpiM8}
\be
(\gamma+1){}_4F_3(a,b,c,\nu+1;d+1,e,\nu;1)
=
{}_3F_2(a,b,c;d,e;1)
+\gamma
{}_3F_2(a+1,b,c;d+1,e;1)
\ee
where $\nu=(\gamma+1)ad/(\gamma d+a)$.

\cite{KarpMdpiM8}
\be
(\gamma+1){}_4F_3(a-1,b,c,\lambda+1;d,e,\lambda;1)
=
{}_3F_2(a,b,c;d,e;1)
+\gamma
{}_3F_2(a,b+1,c+1;d+1,e+1;1)
\ee
where $\lambda=(a-1)bc/(bc+\gamma de)$.

\cite{BuhringJAMSA8}
\be
\frac{1}{\Gamma(b_1)\Gamma(b_2)}\,_4F_3
\left(\begin{array}{c}a_1,a_2,a_3,-m\\b_1,b_2,1-s-m\end{array}\mid 1\right)
=
\frac{(a_1+s)_m(a_2+s)_m(a_3)_m}{(s)_m\Gamma(b_1+m)\Gamma(b_2+m)}
\,_4F_3
\left(\begin{array}{c}b_1-a_3,b_2-a_3,s,-m\\a_1+s,a_2+s,1-a_3-m\end{array}\mid 1\right)
\ee
where $s=b_1+b_2-a_1-a_2-a_3$ not a negative integer or zero, $m=0,1,2\ldots$.

\cite{MillerJCAM231,MillerRMJ43}
\be
_{r+1}F_{r+1}(a,(f_r+1);b,(f_r);y)
=e^y
\,
_{r+1}F_{r+1}(b-a-r,(\xi_r+1);b,(\xi_r);-y)
,
\ee
where $\xi_r$ are nonvanishing zeros of an associated parametric polynomial
of $Q$ degree $r$,
\be
Q_r(t)=\sum_{j=0}^r s_{r-j}\sum_{l=0}^j \left\{\begin{array}{c}j\\l\end{array}\right\}
(a)_l (t)_l(b-a-r-t)_{r-l},
\label{eq.MillerQ}
\ee
and the $s_{r-j}$ are determined by
\be
(f_1+x)\cdots(f_r+x)=\sum_{j=0}^r s_{r-j}x^j.
\ee

\cite{FieldsJAT1}
\be
(-)^m(a_p)_1{}_{p+1}F_p(-m,1+a_p;a_p;1)
=\sum_{j=0}^{p-m}\frac{(j+m)!}{j!}B_j^{(-m)}S_{p-m-j}(a_p),
\ee
where $S_r(a_p)$ are the symmetric polynomials defined by
\be
\prod_{r=1}^p(x+a_r)=\sum_{r=0}^pS_r(a_p)x^{p-r}
\ee
and $B_j^{(-m)}$ are the generalized Bernoulli numbers defined  by
\be
(\frac{e^t-1}{t})^m = \sum_{j\ge 0} \frac{t^j}{j!}B_j^{(-m)}, B_0^{(-m)}=1.
\ee

\cite{LewanowiczMC45}
\be
_{r+2}F_{r+1}(-m,m+\alpha,1+\omega_r; r-h+\alpha,\omega_r;1)
=\frac{(-)^m m!\Gamma(r-h+\alpha)}{\Gamma(m+\alpha)\omega_r} s_h
\ee
for $h=0,1$ where $s_0=1$, $s_1=\frac12(m)_2+\frac12(m-r+1)(m+r+2\alpha-2)+\sum_{i=1}^r \omega_i$.

\cite{SlaterHyp,NorlundAM94}
\be
\left[
\sum_{\nu=1}^B z^{\nu-1}(a_\nu z-b_\nu)\frac{d^\nu}{dz^\nu}+a_0+z^B(1-z)
\frac{d^{B+1}}{dz^{B+1}}
\right] \,_{B+1}F_B((a),(b),z)=0
.
\ee

\cite{GottschalkJPA21}
\be
_pF_{p-1}(1,a,a,\ldots; a+1,a+1,\ldots ;-1) = \frac{(-a/2)^{p-1}}{(p-2)!}[\psi^{(p-2)}(a/2)-\psi^{(p-2)}(a/2+1/2)].
\ee

\cite{GottschalkJPA21,RainvilleBAMS51}
\be
(a_i-b_j+1){}_pF_q(...a_i..;...,b_j...;z)
=
a_i\,{}_pF_q(...a_i+1..;...,b_j...;z)
-(b_j-1){}_pF_q(...a_i..;...,b_j-1...;z).
\ee

\cite[A2.1]{GottschalkJPA21}
\be
\therefore F=\frac{(a_i)_k}{(1+a_i-b_j)_k}
\sum_{n=0}^k \binom{k}{n}\frac{(-)^n(1-b_j)_n}{(1-a_i-k)_n}F(a_i+k-n,b_j-n)
\ee

\cite{GottschalkJPA21,RainvilleBAMS51}
\be
(a_i-a_j){}_pF_q(...a_i..a_j..;.....;z)
=
a_i\,{}_pF_q(...a_i+1..a_j..;......;z)
-a_j{}_pF_q(...a_i..a_j+1...;z).
\ee

\cite[A2.2]{GottschalkJPA21}
\be
\therefore F=-\sum_{n=0}^k \frac{(a_j)_n(1-a_i+a_j)_{n-k-1}(a_i-a_j+k-2n)}{(1-a_i)_{n-k}(1-a_i+a_j)_n} 
F(a_i+k-n,a_j+n),\quad k<|a_i-a_j|
\ee

\cite[A2.3]{GottschalkJPA21}
\be
\therefore F=\frac{(a_i)_k}{(a_i-a_j)_k}F(a_i+k)
-a_j\sum_{n=1}^k \frac{(a_i)_{n-1}}{(a_i-a_j)_n} F(a_i+n-1,a_j+1).
\ee

\cite{GottschalkJPA21}
\be
b_j{}_pF_q(...a_i....;..b_j...;z)
=
a_i\,{}_pF_q(...a_i+1....;..b_j+1...;z)
+(b_j-a_i){}_pF_q(...a_i....;..b_j+1...;z)
.
\ee

\cite[A2.4]{GottschalkJPA21}
\be
\therefore F=\frac{(a_i)_k}{(b_j)_k}\sum_{n=0}^k \binom{k}{n}
\frac{(-1)^n(b_j-a_i)_n}{(1-a_i-k)_n}
F(a_i+k-n,b_j+k).
\ee

\cite{GottschalkJPA21}
\be
(a_i-1){}_pF_q(...a_i..a_j;...;z)
=
(a_i-a_j-1){}_pF_q(...a_i-1..a_j;...;z)
+a_j{}_pF_q(...a_i-1..a_j+1;...;z)
.
\ee

\cite[A2.5]{GottschalkJPA21}
\be
\therefore F=\frac{(1-a_i+a_j)_k}{(1-a_i)_k}\sum_{n=0}^k \binom{k}{n}
\frac{(-1)^n(a_j)_n}{(1-a_i+a_j)_n}
F(a_i-k,a_j+n).
\ee

\cite[A2.6]{GottschalkJPA21}
\be
\therefore F=\frac{(1-a_i+a_j)_k}{(1-a_i)_k}F(a_i-k)
-a_j\sum_{n=1}^k \frac{(1-a_i+a_j)_{n-1}}{(1-a_i)_n}
F(a_i-n,a_j+1).
\ee

\cite{GottschalkJPA21,KrupnikovJCAM78}
\be
_pF_q(a,b... ; a+1,b+1,...;z)
=
\frac{a}{a-b}{}_{p-1}F_{q-1}(b,...;b+1,...;z)
-\frac{b}{a-b}{}_{p-1}F_{q-1}(a,...;a+1,....;z)
.
\ee

\cite{RainvilleBAMS51}
\be
a_1{}_pF_q(a_1,... ; b_1,b_2,...;z)
=
a_1{}_pF_q(a_1+1,... ; b_1,b_2,...;z)
- z
\sum_{j=1}^q U_j\, {}_pF_q(... ; ...b_j+1...;z)
;\quad p<q
\ee
where
\be
U_j = \frac{\prod_{s=1}^p (a_s-b_j)}{b_j\prod_{s=1,s\neq j}^q (b_s-b_j)}.
\ee

\cite{RainvilleBAMS51}
\be
(a_1+z){}_pF_q(a_1,... ; b_1,b_2,...;z)
=
a_1{}_pF_q(a_1+1,... ; b_1,b_2,...;z)
- z
\sum_{j=1}^q U_j\, {}_pF_q(... ; ...b_j+1...;z)
;\quad p=q
\ee

\cite{RainvilleBAMS51}
\be
[(1-z)a_1+(A-B)z]{}_pF_q(a_1,... ; b_1,b_2,...;z)
=
(1-z)a_1{}_pF_q(a_1+1,... ; b_1,b_2,...;z)
- z
\sum_{j=1}^q U_j\, {}_pF_q(... ; ...b_j+1...;z)
;\quad p=q+1
\ee
where $A\equiv \sum_{s=1}^p a_s$, $B\equiv \sum_{s=1}^q b_s$.

\cite{RainvilleBAMS51}
\be
{}_pF_q(a_1,... ; b_1,b_2,...;z)
=
{}_pF_q(a_1,...,a_k-1,... ; b_1,b_2,...;z)
+ z
\sum_{j=1}^q W_{j,k}\, {}_pF_q(... ; ...b_j+1...;z)
;\quad p\le q
\ee
where
\be
W_{j,k} = \frac{\prod_{s=1,s\neq k}^p (a_s-b_j)}{b_j\prod_{s=1,s\neq j}^q (b_s-b_j)}.
\ee

\cite{RainvilleBAMS51}
\be
(1-z) {}_pF_q(a_1,... ; b_1,b_2,...;z)
=
{}_pF_q(a_1,...,a_k-1,... ; b_1,b_2,...;z)
+ z
\sum_{j=1}^q W_{j,k}\, {}_pF_q(... ; ...b_j+1...;z)
;\quad p= q+1
\ee

\cite{WimpMathComp22}
Let $R\ge 0$, $T\ge 0$ be integers, $\tau=\max(T+1,R+2)$, $\gamma,c_i,d_j$ for $i=1\ldots R$,
$j=1\ldots T+1$ be complex constants such that none of $\gamma$, $\gamma+1-d_j$, $j=1,\cdots,T$
are negative integers or zero. Define hypergeometric Polynomials $P_n(z)$
\be
P_n(z)=\frac{1}{\Gamma(d_T)}{}_{R+2}F_T(-n,n+\gamma,c_r ; d_T; z).
\ee
Then
\be
\sum_{\nu=0}^\tau [C_\nu+zD_\nu]P_{n-\nu}(z)=0,\quad n=\tau, \tau+1,\tau+2,\ldots
\ee
where
\begin{multline}
C_\nu=\frac{(-)^\nu(n+1-\nu)_\nu (1-\gamma-2n)_{2\nu} (n-\nu-1+d_{T+1})}
{\nu!(n+\gamma-\nu)_\nu (\tau+1-\gamma-2n)_\nu(n+d_{T+1}-1)}
\\ \times
{}_{T+3}F_{T+2}(-\nu,2n+\gamma-\tau-\nu,n-\nu+d_{T+1};2n+\gamma+1-2\nu,n-\nu-1+d_{T+1};1) 
\end{multline}
and
\begin{multline}
D_\nu=\frac{(-)^{\nu+1}(n+1-\nu)_\nu (1-\gamma-2n)_{2\nu} (n-\nu+c_R)}
{\Gamma(\nu)(n+\gamma-\nu)_\nu (\tau+1-\gamma-2n)_{\nu-1}(n+d_{T+1}-1)}
\\ \times
{}_{R+2}F_{R+1}(1-\nu,2n+\gamma+1-\tau-\nu,n+1-\nu+c_R;2n+\gamma+1-2\nu,n-\nu+c_R;1) 
\end{multline}
and $D_0=D_\tau=0$.

\cite{FieldsJAT1,LewanowiczMC45}
The generalized Jacobi functions
\be
E_n(z,\lambda)\equiv {}_{r+3}F_s(-n,n+\lambda,\alpha_r,1; \beta_s;z)
\ee
with $r+3\le 2$, or $r+2=s$ and $|\arg(1-z)|< \pi$ and
\be
H_n(z,\lambda)\equiv \frac{(\beta_s-1)/z}{(n+1)(n+\lambda-1)(\alpha_r-1)}
{}_{s+1}F_{r+2}(2-\beta_s,1;2+n,2-n-\lambda,2-\alpha_r;\frac{(-1)^{s+r}}{z})
\ee
with $r+1\ge s$, or $r+2=s$ and $|\arg(1-1/z)<\pi$
satisfy (under conditions given in the reference)
\be
\Phi_n(z,\lambda)+\sum_{m=1}^t [A_m(n,\lambda)+zB_m(n,\lambda] 
\Phi_{n-m}(z,\lambda)
=\frac{(\beta_s-1)(n+\lambda)_n)}{(n+\beta_s-1)(n+\lambda-t)_n}
\ee
where $t=\max(r+2,s)$, $B_t(n,\lambda)=0$ and where
\begin{multline}
A_m(n,\lambda)=\frac{(n+1-m)_m(2n+\lambda-2m)_{2m}(n-m-1+\beta_s}
{m!(n+\lambda-m)_m(2n+\lambda-t-m)_m(n-1+\beta_s)}
\\ \times
{}_{s+2}F_{s+1}(-m,2n+\lambda-t-m,n-m+\beta_s;2n+\lambda+1-2m,n-m-1+\beta_s;1),
\end{multline}
\begin{multline}
B_m(n,\lambda)=\frac{(n+1-m)_m(2n+\lambda-2m)_{2m}(n-m+\alpha_r}
{(m-1)!(n+\lambda-m)_m(2n+\lambda-t-m+1)_{m-1}(n-1+\beta_s)}
\\ \times
{}_{r+2}F_{r+1}(1-m,2n+\lambda-t-m+1,n-m+1+\alpha_r;2n+\lambda+1-2m,n-m+\alpha_r;1).
\end{multline}

\cite{Mishevarxiv2020}
\be
_{A+1}F_B[(a),-m; (b);z]
=
\frac{((a))_m(-z)^m}{((b))_m}
\,_{B+1}F_A[(1-b-m),-m;(1-a-m);\frac{(-1)^{A+B}}{z}]
.
\label{eq.hypRev}
\ee
This corrects equation \cite[(2.2.3.2)]{SlaterHyp}.

\cite{SlaterHyp,TakayamaJSC20}
\be
[\delta(\delta+b_1-1)(\delta+b_2-1)\cdots (\delta+b_B-1)
-z(\delta+a_1)(\delta+a_2)\cdots(\delta+a_A)] \,_AF_B((a),(b),z)=0
,
\ee
where $\delta =z\frac{d}{dz}$.

\cite{ChaundyQJM14}
\be
\frac{(\delta +h)_m}{(\delta +k)_m}
\,_rF_s\left(
\begin{array}{c}a,\ldots\\ c,\ldots\end{array}
\mid x\right)
=
\frac{(h)_m}{(k)_m}
\,_{r+2}F_{s+2}\left(
\begin{array}{c}a,\ldots,h+m,k\\ c,\ldots,h,k+m\end{array}
\mid x\right)
.
\ee

\cite{ChaundyQJM14,NorlundAM94}
\be
\sum_{n\ge 0}\frac{(a)_n}{n!} F\left(
\begin{array}{c}-n,A,\ldots\\ C\ldots\end{array}\mid p
\right)x^n
=(1-x)^{-a}F\left(
\begin{array}{c}a,A,\ldots \\ C,\ldots
\end{array}\mid
-\frac{px}{1-x}
\right)
.
\ee

\cite{ChaundyQJM14}
\be
\sum_{n\ge 0}\frac{(a)_n}{n!} F\left(
\begin{array}{c}-n,a+n,A,\ldots\\ C\ldots\end{array}\mid p
\right)x^n
=(1-x)^{-a}F\left(
\begin{array}{c}a/2,(1+a)/2,A,\ldots \\ C,\ldots
\end{array}\mid
-\frac{4px}{(1-x)^2}
\right)
.
\ee

\cite{ChaundyQJM14}
\be
\sum_{n\ge 0}\frac{(a)_n}{n!} F\left(
\begin{array}{c}-n,A,\ldots\\ 1-a-n,C,\ldots\end{array}\mid p
\right)x^n
=(1-x)^{-a}F\left(
\begin{array}{c}A,\ldots \\ C,\ldots
\end{array}\mid
px
\right)
.
\ee

\cite{SlaterHyp}
\be
_{A+1}F_{B+1}\left( c,(a);d,(b);z\right)=
\frac{\Gamma(d)}{\Gamma(c)\Gamma(d-c)}
\int_0^1t^{c-1}(1-t)^{d-c-1}\,_AF_B\left((a),(b),tz\right)dt
.
\ee

\cite{CoffeyJCAM233,DriverETNA25}
\begin{multline}
_{q+1}F_q\left(
\begin{array}{c}
a,\frac{b}{q},\frac{b+1}{q},
\cdots
,\frac{b+q-1}{q}
\\
\frac{c}{q},\frac{c+1}{q},
\cdots
,\frac{c+q-1}{q}
\end{array}\mid x
\right)
=
\frac{\Gamma(c)}{\Gamma(b)\Gamma(c-b)}
\int_0^1 t^{b-1}(1-t)^{c-b-1}(1-xt^q)^{-a}dt
\\
=
\frac{\Gamma(c)}{\Gamma(b)\Gamma(c-b)}
\sum_{k\ge 0}
\binom{c-b-1}{k}
\frac{(-)^k}{b+k}
\,_2F_1(a,(b+k)/q ; 1+(b+k)/q; x)
\end{multline}
\cite{BaileyPLMS28}
\begin{multline}
_\lambda F_{\lambda-1}(-n/\lambda,(1-n)/\lambda,\ldots, (\lambda-1-n)/\lambda;
1/\lambda,
2/\lambda,
\ldots, (\lambda-1)/\lambda;1)
\\
=\frac{2^n}{\lambda}
\sum_{r=1}^\lambda
\cos^n\frac{(2r-1)\pi}{2\lambda}\cos\frac{n(2r-1)\pi}{2\lambda},
\end{multline}
where $n>1$ is an integer and $\lambda$ an odd prime.

\cite{SofoIJPAM50}
\begin{multline}
k\int_0^1 \frac{x^j(1-x)^{k-1}}{1-x^k}dx=
\frac{1}{\binom{k+j}{k}}
\\ \times
{}_{1+k}F_k( (1+j)/k, (2+j)/k,\ldots (k+j)/k,1 ;
(1+j+k)/k, (2+j+k)/k,\ldots (2k+j)/k ;1)
.
\end{multline}

\cite{DriverETNA25}
\be
_{p+k}F_{q+k}\left(
\begin{array}{c}
a_1,\ldots,a_p,\frac{\alpha}{k},\frac{\alpha+1}{k},\ldots,\frac{\alpha+k-1}{k}\\
b_1,\ldots,b_q,\frac{\alpha+\beta}{k},\frac{\alpha+\beta+1}{k},\ldots,\frac{\alpha+\beta+k-1}{k}\\
\end{array}
| ct^k
\right)
=\frac{t^{1-\alpha-\beta}}{B(\alpha,\beta)}
\int_0^t x^{\alpha-1}(t-x)^{\beta-1}
\,_pF_q\left(
\begin{array}{c}
a_1,\ldots,a_p\\
b_1,\ldots,b_q\\
\end{array}
| cx^k
\right)dx.
\ee

\cite{MillerZAMP}
\be
_{r+2}F_{r+1}\left(\begin{array}{cc}a,b,& (f_r+1)\\ c,& (f_r)\end{array}\mid x\right)
=(1-x)^{-a}\,_{r+2}F_{r+1}\left(\begin{array}{cc}a,\lambda,&(\xi_r+1)\\c,&(\xi_r)\end{array}\mid\frac{x}{x-1}\right)
\ee
where $\lambda=c-b-r$, $|x|<1$, $\Re x<1/2$ and the $(\xi_r)$
are the nonvanishing zeros of the parametric polynomial $Q_r(t)$ given
in (\ref{eq.MillerQ}).

\cite{MillerZAMP}
\be
_{r+2}F_{r+1}\left(\begin{array}{cc}a,b,& (f_r+1)\\ c,& (f_r)\end{array}\mid x\right)
=(1-x)^{c-a-b-r}\,_{r+2}F_{r+1}\left(\begin{array}{cc}\lambda,\lambda',&(\xi_r+1)\\c,&(\xi_r)\end{array}\mid x\right)
\ee

\cite{MaierTAAMS358}
\be
_{r+1}F_r\left(\begin{array}{ccc}\alpha_1,&\ldots,&\alpha_{r+1}\\
\beta_1+1,&\ldots,&\beta_r+1\end{array}\mid x\right)
=
(1-x){}_{r+1}F_r\left(\begin{array}{ccc}\alpha_1+1,&\ldots,&\alpha_{r+1}+1\\
\beta_1+1,&\ldots,&\beta_r+1\end{array}\mid x\right),
\ee
for all $r\ge 1$ for which no $\beta_i+1$ is a nonpositive integer and the equations (\ref{eq.maiertaams}) are satisfied.

\cite{GottschalkJPA21,KrupnikovJCAM78}
\be
_pF_{p-1}(1,a,a,\ldots; a+1,a+1,\ldots ;1) = a^{p-1}\zeta(p-1,a) = \frac{(-a)^{p-1}}{(p-2)!}\psi^{(p-2)}(a).
\ee

\cite{GottschalkJPA21}
\be
_{p}F_{p-1}\left(\begin{array}{c}1,a,a,\ldots,a\\ a+1,\ldots,a+1\end{array}\mid z\right)
=
a^{p-1}\Phi(z,p-1,a).
\ee

\cite{AdamchikJCAM79}
\be
_{p+1}F_p\left(\begin{array}{c}a,a,\ldots,a\\ a+1,\ldots,a+1\end{array}\mid 1\right)
=
\frac{a^p}{(a-1)!}\sum_{k=0}^{a-1} (-)^{a-k-1}\zeta(p-k)
\left[\begin{array}{c}a\\k+1\end{array}\right]
\ee
for $a$ a positive integer,
where the $[]$ is Stirling numbers of the first kind.

\cite{GottschalkJPA21}
\be
_{p}F_{p-1}\left(\begin{array}{c}1,a,a,\ldots,a\\ a+1,\ldots,a+1\end{array}\mid -1\right)
=
\frac{(-a/2)^{p-1}}{(p-2)!}[\psi^{(p-2)}(a/2)-\psi^{(p-2)}(a/2+1/2)]
.
\ee

\cite{NorlundAM94}
\be
_{p+1}F_p\left(a_1,\ldots,a_{p+1}; b_1\ldots b_p\mid z\zeta\right)
=
(1-z)^{-a_1}\sum_{k\ge 0}\frac{ (a_1)_k}{k!} \,
_{p+1}F_p\left(-k,a_2,\ldots,a_{p+1}; b_1,\ldots, b_p\mid \zeta\right)
\left(\frac{z}{z-1}\right)^k
.
\ee

\cite{AdamchikJCAM79}
\be
_{p+1}F_p\left(\begin{array}{c}a,a,\ldots,a\\ a+1,\ldots,a+1\end{array}\mid 1\right)
=
\frac{(-)^{p-1}\pi a^p}{\sin(a\pi)(p-1)!}w(a,p-1)
\ee
where
\be
w(n,m)=\frac{1}{(n-1)!}\sum_{i=m+1}^n \left[
\begin{array}{c}n\\ i\end{array}
\right] (i-m)_m (-1)^{n-i}n^{i-m-1}
,
\ee
recursively
\be
w(n,0)=1,\quad
w(n,m)=\sum_{k=0}^{m-1} (1-m)_k H_{n-1}^{(k+1)} w(n,m-1-k)
\ee
and the Harmonic numbers defined in (\ref{eq.Hnrdef}).

\cite{MaierTAAMS358}
\be
_{r+1}F_r\left(
\begin{array}{ccc}
\alpha_1&,\ldots,&\alpha_{r+1}\\
\beta_1+1&,\ldots,&\beta_r+1\\
\end{array}
\mid 1
\right)=\Gamma\left[
\begin{array}{ccc}
\beta_1+1&,\ldots,&\beta_r+1\\
\alpha_1+1&,\ldots,&\alpha_{r+1}+1\\
\end{array}
\right].
\ee
for all $r\ge 1$ for which no $\beta_i+1$ is a nonpositive integer
and the equations (\ref{eq.maiertaams}) are satisfied.

\cite{ChuRDCM}
\begin{multline}
_{1+m}F_m(1+\frac{u}{v},\frac{a}{m},\frac{a+1}{m},\ldots \frac{a+m-1}{m};
u/v, \frac{a+1}{m-1},\ldots,\frac{a+m-1}{m-1}; \frac{m^m(t-1)}{(m-1)^{m-1}t^m})
\\
=t^a \frac{u+(av-mu+u)(t-1)}{u(m+t(1-m))}.
\end{multline}

\cite{ChuRDCM}
\begin{multline}
_{1+2m}F_{2m}(-n,1+\frac{u}{v},\frac{a}{m},\frac{a+1}{m},\ldots \frac{a+m-1}{m}
, \frac{c+n}{m-1},\frac{c+n+1}{m-1}\cdots \frac{c+n+m-2}{m-1};
\\
u/v, \frac{c+1}{m},\ldots,\frac{c+m}{m}
,
\frac{a+1}{m-1},\frac{a+2}{m-1}\ldots,\frac{a+m-1}{m-1}
; 1)
\\
=\frac{c}{c-a} \frac{c+mn-a(1+vn/u)}{c+mn}\frac{(c-a)_n}{(c)_n}.
\end{multline}

\cite{ChuRDCM}
\begin{multline}
_{2m}F_{2m-1}(-n,\frac{a-1}{m},\frac{a}{m},\ldots \frac{a+m-2}{m}
,\frac{c+n}{m-1},\frac{c+n+1}{m-1},\cdots \frac{c+n+m-2}{m-1}
\\
; \frac{c+1}{m},\frac{c+2}{m}\cdots \frac{c+m}{m},
\frac{a+1}{m-1},\frac{a+2}{m-1}\ldots,\frac{a+m-1}{m-1}
; 1)
\\
=\frac{(c-a)_n}{(c)_n} \frac{c\{(a+m-1)(c-a+mn)(c-a+mn+1)-(m-1)(c+mn)mn\}}
{(c+mn)(a+m-1)(c-a+1)(c-a)}
.
\end{multline}

\cite{ChuRDCM}
\begin{multline}
_{2m}F_{2m-1}(-n,\frac{a-1}{m},\frac{a}{m},\ldots \frac{a+m-2}{m}
,\frac{c+n}{m-1},\frac{c+n+1}{m-1},\cdots \frac{c+n+m-2}{m-1}
\\
; \frac{c}{m},\frac{c+1}{m}\cdots \frac{c+m-1}{m},
\frac{a+1}{m-1},\frac{a+2}{m-1}\ldots,\frac{a+m-1}{m-1}
; 1)
\\
=[1+\frac{amn}{(c-a)(a+m-1)}]\frac{(c-a)_n}{(c)_n}
.
\end{multline}

\cite{ChuRDCM}
\begin{multline}
_{2m}F_{2m-1}(-n,\frac{a-1}{m},\frac{a}{m},\ldots \frac{a+m-2}{m}
,\frac{c+n}{m-1},\frac{c+n+1}{m-1},\cdots \frac{c+n+m-2}{m-1}
\\
; \frac{c-1}{m},\frac{c}{m}\cdots \frac{c+m-2}{m},
\frac{a+1}{m-1},\frac{a+2}{m-1}\ldots,\frac{a+m-1}{m-1}
; 1)
\\
=[1+\frac{amn}{(c-1)(a+m-1)}]\frac{(c-a)_n}{(c)_n}
.
\end{multline}

\cite{ExtonJCAM79}
\begin{multline}
_{q+1}F_q\left(\begin{array}{c}a_1,\ldots,a_{q+1} \\ b_1,\ldots,b_q\end{array}\mid 1\right)
+\,_{q+1}F_q\left(\begin{array}{c}a_1,\ldots,a_{q+1} \\ b_1,\ldots,b_q\end{array}\mid -1\right)
\\
=
2\,_{2q+2}F_{2q+1}\left(\begin{array}{c}a_1/2,a_1/2+1/2\ldots,a_{q+1}/2+1/2 
\\ b_1/2,b_1/2+1/2,\ldots,b_q/2+1/2,1/2\end{array}\mid 1\right)
.
\end{multline}

\cite{ExtonJCAM79}
\begin{multline}
_{q+1}F_q\left(\begin{array}{c}a_1,\ldots,a_{q+1} \\ b_1,\ldots,b_q\end{array}\mid 1\right)
-\,_{q+1}F_q\left(\begin{array}{c}a_1,\ldots,a_{q+1} \\ b_1,\ldots,b_q\end{array}\mid -1\right)
\\
=
2\frac{a_1a_2\cdots a_{q+1}}
{b_1b_2\cdots b_q}\,_{q+1}F_{q}
\left(\begin{array}{c}a_1/2+1/2,a_1/2+1\ldots,a_{q+1}/2+1/2 ,a_{q+1}/2+1
\\ b_1/2+1/2,b_1/2+1,\ldots,b_q/2+1,3/2\end{array}\mid 1\right)
.
\end{multline}

\cite{ExtonJCAM79}
\begin{multline}
_{q+1}F_q\left(\begin{array}{c}a,a+2,a+4\ldots,a+2q 
\\ a+1,a+3,a+5\ldots a+2q-1\end{array}\mid 1\right)
+\,_{q+1}F_q\left(\begin{array}{c}a,a+2,a+4,\ldots,a+2q 
\\ a+1,a+3,a+5,a+2q-1\end{array}\mid -1\right)
\\
=
\,_2F_1(a/2,a+q+1/2;1/2;1)
.
\end{multline}

\cite{ExtonJCAM79}
\begin{multline}
_{q+1}F_q\left(\begin{array}{c}a,a+2,a+4\ldots,a+2q 
\\ a+1,a+3,a+5\ldots a+2q-1\end{array}\mid 1\right)
-\,_{q+1}F_q\left(\begin{array}{c}a,a+2,a+4,\ldots,a+2q 
\\ a+1,a+3,a+5,a+2q-1\end{array}\mid -1\right)
\\
=
\frac{2a(a+2)(a+4)\cdots (a+2q)}{(a+1)(a+3)\cdots (a+2q-1)}\,
\,_2F_1(a/2+1/2,a/2+q+1;3/2;1)
.
\end{multline}

\cite{ExtonJCAM79}
\begin{multline}
_{q+1}F_q\left(\begin{array}{c}a,a+2,a+4\ldots,a+2q 
\\ a+1,a+3,a+5\ldots a+2q-1\end{array}\mid 1\right)
\\
=
\frac{\sqrt{\pi}\Gamma(-a-q)}{\Gamma(1/2-a)\Gamma(-a/2-q)}\,
+
\frac{2(a/2)_q}{(a/2+1/2)_{q-1}\Gamma(1-a/2)\Gamma(1/2-a/2-q)},
\end{multline}
if $\Re(a+q)<0$, and a similar expression for argument $-1$.

\cite{KarlssonJMP12}
\begin{equation}
_pF_q\left(\begin{array}{c} b_1+m_1,a_2\ldots a_p
\\
b_1,\ldots b_q
\end{array}\mid z\right)
=
\sum_{j=0}^{m_1}
z^j \binom{m_1}{j} \frac{(a_2)_j\cdots (a_p)_j}{(b_1)_j\cdots (b_q)_j}
{}_{p-1}F_{q-1}
\left(\begin{array}{c}a_2+j,\ldots ,a_p+j\\
b_2+j,\ldots ,b_q+j\\
\end{array}\mid z\right)
.
\end{equation}

\cite{KarlssonJMP12}
\begin{equation}
_{n+2}F_{n+1}\left(\begin{array}{c} b_1+m_1,\ldots,b_n+m_n,b,a
\\
b_1,\ldots b_n,b+1
\end{array}\mid 1\right)
=
\frac{\Gamma(b+1)\Gamma(1-a)}{\Gamma(b+1-a)}\prod_{k=1}^n \frac{(b_k+m_k)_{-b}}{(b_k)_{-b}}
.
\end{equation}

\cite{GottschalkJPA21,KarlssonJMP12}
\begin{multline}
_pF_q\left(\begin{array}{c} b_1+m_1,\ldots,b_n+m_n,a_{n+1},\ldots a_p
\\
b_1,\ldots b_q
\end{array}\mid z\right)
=
\sum_{j_1=0}^{m_1}
\cdots
\sum_{j_n=0}^{m_n}
A(j_1,\ldots j_n)z^{J_n}
\\ \times
{}_{p-n}F_{q-n}
\left(\begin{array}{c}a_{n+1}+J_n,\ldots ,a_p+J_n\\
b_{n+1}+J_n,\ldots ,b_q+J_n\\
\end{array}\mid z\right)
,
\end{multline}
where $m$ are positive integers, where $J_n=j_1+\cdots +j_n$ and
\begin{equation*}
A(j_1,\ldots j_n)=\binom{m_1}{j_1}
\cdots
\binom{m_n}{j_n}
\frac{(b_2+m_2)_{J_1}(b_3+m_3)_{J_2}\cdots (b_n+m_n)_{J_{n-1}}(a_{n+1})_{J_n}\cdots (a_p)_{J_n}}
{(b_1)_{J_1}(b_2)_{J_2}\cdots (b_n)_{J_{n}}(b_{n+1})_{J_n}\cdots (b_q)_{J_n}}
.
\end{equation*}
Also: a similar reduction in the degree if the denominators have positive offset relative to the numerators.

\cite{GottschalkJPA21} 
\begin{multline}
_pF_q(a_1,\ldots,a_n,a_{n+1},\ldots a_p ; a_1+m_1+1,\ldots a_n+m_n+1,b_{n+1},\ldots b_q;z)
\frac{1}{(a_r)_{m_r+1}}
\\
=
\sum_{i=1}^n \sum_{j=0}^{m_i} \frac{(-m_i)_j}{j!(a_i+j)m_i!}
\prod_{r=1,r\neq i}^n \frac{1}{(a_r-a_i-j)_{m_r+1}}
{}_{p-n+1}F_{q-n+1}(a_i+j,a_{n+1},\ldots,a_p ; a_i+j+1,b_{n+1},\ldots b_q;z).
\end{multline}
where $\{a_1,\ldots a_n\}$ are all different and, if $a_r-a_i=N\in N$, then $N>m_i$ is satisfied.

\cite{KarlssonJMP12}
\be
_{n+2}F_{n+1}(b_1+m_1,\cdots,b_n+m_n,b,a ;
b_1,\cdots,b_n,b+1;1)
=\frac{\Gamma(b+1)\Gamma(1-a)}{\Gamma(b+1-a)}\prod_{k=1}^n
\frac{(b_k+m_k)_{-b}}{(b_k)_{-b}}.
\ee
for $\Re (-a)>m_1+m_2+\cdots m_n-1$.

\cite{RainaITSF3}
\begin{multline}
_pF_q(a_1,a_2,\ldots a_p ; a_1+\alpha+1,b_2,\ldots b_q;z)
\frac{\Gamma(1+\alpha)\Gamma(a_1)}{\Gamma(1+\alpha+a_1)}
\\
=\sum_{k=0}^\infty \binom{\alpha}{k}\frac{(-)^k}{a_1+k}\,{}_pF_q
(a_1+k,a_2,\ldots, a_p ; a_1+k+1,b_2,\ldots b_q;z).
\end{multline}

\cite{RainaITSF3}
Let $p\le q$, $|z|<\infty$, and for $p=q+1$ $|\arg(1-z)|<\pi$,
the parameters $a_i$, $i=1\ldots n$ be all different and $\Re(a_i+\alpha_i)>-1$,
then for arbitrary $\alpha_i$
\begin{multline}
_pF_q(a_1,a_2,\ldots a_n,a_{n+1},\ldots a_p ; 
a_1+\alpha+1,\ldots,a_n+\alpha_n+1,b_{n+1},\ldots b_q;z)
\prod_{\mu=1}^n \frac{\Gamma(1+\alpha_\mu)\Gamma(a_\mu)}{\Gamma(1+\alpha_\mu+a_1)}
\\
=
\sum_{i=1}^n 
\sum_{j_1=0}^\infty \ldots
\sum_{j_n=0}^\infty 
\frac{1}{a_i+j_i}
\prod_{\mu=1,\mu\neq i}
\binom{\alpha_\mu}{j_\mu}
\frac{(-)^{j_\mu}}{a_\mu+j_\mu-a_i-j_i}\\
\times {}_LF_M
(a_i+j_i,a_{n+1},\ldots, a_p ; a_i+j_i+1,b_{n+1},\ldots b_q;z).
\end{multline}
where $L=p-n+1$, $M=q-n+1$ provided that no $b_j$, $j=n+1,\ldots q$, 
is zero or a negative integer.

\cite{FieldsJAT1}
\be
_{p+r}F_{q+s}(\alpha_p,\sigma_r; \beta_q,\rho_s;\omega z)
=
\sum_{n=0}^\infty \frac{(\alpha_p)_n(-\omega)^n}{(\beta_q)_n(n+\lambda)_nn!}
{}_pF_{q+1}(n+\alpha_p ; n+\beta_q,2n+\lambda+1;\omega)
{}_{r+2}F_s(-n,n+\lambda,\sigma_r ; \rho_s;z).
\ee

\cite{FieldsJAT1}
\be
_{p+r}F_{q+s}(\alpha_p,\sigma_r; \beta_q,\rho_s;\omega z)
=
\sum_{n=0}^\infty \frac{(\alpha_p)_n(-\omega)^n}{(\beta_q)_n(n+\lambda)_nn!}
{}_pF_q(n+\alpha_p ; n+\beta_q;\omega)
{}_{r+1}F_s(-n,\sigma_r ; \rho_s;z).
\ee

\cite{KrupnikovJCAM78}
\begin{multline}
_{p+1}F_p(c_1+n,c_2+n,\cdots c_p+n,1; d_1+n,d_2+n,\ldots d_q+n; z
\\
=\frac{\prod_{h=1}^q(d_h)_n}{\prod_{h=1}^p(c_h)_n}z^{-n}[_{p+1}F_q(c_1,c_2,\ldots c_p,1 ; d_1,d_2,\ldots d_q;z)
-\sum_{k=0}^{n-1} \frac{\prod_{h=1}^p (c_h)_k}{\prod_{h=1}^q (d_h)_k}z^k
] ; n\ge 1
\end{multline}

\cite{BuhringJAMSA8}
\begin{multline}
\frac{1}{\prod_{j=1}^p\Gamma(b_j)}
\,
_{p+2}F_{p+1}
\left(
\begin{array}{c}
a_1,a_2,\ldots,a_{p+1},-m\\
b_1,b_2,\ldots,b_p,1-s-m\\
\end{array}\mid 1
\right)
\\
=
\frac{(a_1+s)_m(a_2+s)_m\prod_{j=3}^{p+1}(a_j)_m}
{(s)_m\prod_{j=1}^p \Gamma(b_j+m)}
\sum_{k=0}^m\frac{(s)_k(-m)_k}{(a_1+s)_k(a_2+s)_k k!}
\\
\times (S)_k\, _{p+1}F_p\left(
\begin{array}{c}
1-b_1-m,1-b_2-m,\ldots,1-b_p-m,-k\\
1-a_3-m,1-a_4-m,\ldots 1-a_{p+1}-m,1-S-k
\end{array}\mid 1
\right)
\end{multline}
where
$s=\sum_{j=1}^p b_j-\sum_{j=1}^{p+1}a_j$
and
$S=a_1+a+2+s+m-1$ and $s$ not a negative integer or zero.

\cite{GottschalkJPA21} 
\begin{multline}
_{p+1}F_p(a_1,\ldots a_p,b ; a_1+m_1+1,\ldots a_p+m_p+1;z)
\prod_{r=1}^p \frac{\Gamma(a_r)}{\Gamma(a_r+m_r+1)}
\\
=
\sum_{i=1}^n \sum_{j=0}^{m_i} \frac{(-m_i)_j}{j!m_i!(a_i+j)}
\prod_{r=1,r\neq i}^p \frac{\Gamma(a_r-a_i-j)}{\Gamma(a_r+m_r-a_i-j+1)}
{}_2F_1(a_i+j,b ; a_i+j+1;z).
\end{multline}
where $\{a_1,\ldots a_n\}$ are all different and, if $a_r-a_i=N\in N$, then $N>m_i$ is satisfied.

\cite{GottschalkJPA21} 
\begin{multline}
_pF_q(a_1,\ldots a_p,b ; a_1+m_1+1,\ldots a_p+m_p+1;1)
=
\Gamma(1-b)\prod_{r=1}^p \frac{\Gamma(a_r+m_r+1)}{\Gamma(a_r)}
\\
\times
\sum_{i=1}^p \sum_{j=0}^{m_i} \frac{(-)^j \Gamma(a_i+j)}{\Gamma(1-j+m_i)j!\Gamma(a_i+j+1-b)}
\prod_{r=1,r\neq i}^p \frac{\Gamma(a_r-a_i-j)}{\Gamma(a_r+m_r-a_i-j+1)}
\quad m_i\in \mathbb{Z}^+, b\not \in \mathbb{Z}.
\end{multline}

\cite{GottschalkJPA21}
\be
_{p+1}F_p(a,a,\ldots,b ; a+1,a+1,\ldots a+1; 1)
=
-\frac{\Gamma(1-b)(-a)^p}{\Gamma(p)}\left[\frac{d^{(p-1)}}{dr^{(p-1)}}\frac{\Gamma(r)}{\Gamma(r+1-b)}\right]_{r=a}, b\not \in \mathbb{Z}^+
\ee

\cite{BerndtBLMS15}
\be
_0F_2(m+n+1,n+1;x){}_0F_2(m+1,1-n;-x)
=
1+\sum_{k\ge 1}\frac{\alpha_k(\frac{2m+n+k+2}{2})_k(2x)^k}{(m+n+1)_k(m+1)_kk!}
\ee
where
\be
\alpha_k=\begin{cases}
\frac{n}{(n^2-1^2)(n^2-3^2)\cdots(n^2-k^2)},& \mathrm{if} k \mathrm{odd};\\
\frac{1}{(n^2-2^2)(n^2-4^2)\cdots(n^2-k^2)},& \mathrm{if} k \mathrm{even}
\end{cases}
\ee

\cite{CohenMC31}
\be
\sum_{k\ge 0} t^k\frac{(\alpha+1+sk)_k}{k!}
\,_rF_r
\left(\begin{array}{c} \Delta(-k,r) \\ \Delta(-\alpha-sk-k,r)\end{array}
\mid x(\beta+sk)
\right)
=
(1-z)^{\alpha+1}
\frac{1}{1+sz+rsy} e^{-\beta y},
\ee
where $t=(-z)(1-z)^{-s-1}\exp(sy)$, $x=(-y)[1-1/z]^r$,
$(a)_k=\Gamma(a+k)/\Gamma(a)$ and
$\Delta(-k,r)\equiv -k/r, (-k+1)/r,\ldots (-k+r-1)/r$, $|t|<1$.

\cite{CohenMC31}
\begin{multline}
\sum_{k\ge 0} t^k
\frac{(\alpha+1+sk)_k}{(\alpha+l+1+sk)k!((\alpha+l+1+s''+sk)/s'')_{l'}}
\,_rF_r\left(
\begin{array}{c}\Delta(-k,r) \\ \Delta(-\alpha-sk-k,r)
\end{array}
\mid x(\beta+sk)
\right)
\\
=
\frac{(1-z)^{\alpha+l+1}e^{-\beta y}}{l'!}
\sum_{q=0}^{l'}
\sum_{k=0}^{l+s''q}
\frac{(-z)^k(1-z)^{s''q-k}(-l')_q (-l-s''q)_k}{k!q!(\alpha+l+1+s''q+sk)}
\\ \times
\,_1F_1\left(
\begin{array}{c}
1 \\ (\alpha+l+1+s'+s''r+sk)/s'
\end{array}\mid y(\beta-\alpha-l-1-s''q)
\right),
\end{multline}
where $t=(-z)(1-z)^{-s-1}e^{sy}$, $x=-y[1-1/z]^r$ and $|t|<1$.

\cite{CohenMC31}
\be
\sum_{p\ge 0}t^p\frac{(\beta+s'p)^p}{p!}
\,_rF_0\left(
\begin{array}{c} \Delta(-p,r) \\ -
\end{array}
\mid\frac{xr^r(\alpha+s'p)}{(\beta+s'p)^r}
\right)
=\frac{e^{-\beta y -\alpha z}}{1+s'y+rs'z},
\ee
where $t=-ye^{s'y+s'z}$, $x=-z/y^r$ and $|s'y\exp(s'y+s'z+1)|<1$.

\cite{CohenMC31}
\be
\sum_{p\ge 0}t^p\frac{(\beta+s'p)^p}{p!(\alpha+s'p)}
\,_rF_0\left(
\begin{array}{c} \Delta(-p,r) \\ -
\end{array}
\mid\frac{xr^r(\alpha+s'p)}{(\beta+s'p)^r}
\right)
=\frac{e^{-\beta y -\alpha z}}{\alpha}
\,_1F_1\left(\begin{array}{c}1 \\ \alpha/s'+1\end{array}\mid y(\beta-\alpha)\right)
\ee
where $t=-ye^{s'y+s'z}$, $x=-z/y^r$ and $|s'y\exp(s'y+s'z+1)|<1$.

\cite{MalaniCKMS15,LinMCM36}
\begin{multline}
\sum_n \frac{(d)_n (d+1/2)_n }{(3/2)_n}\frac{x^{2n}}{n!}{}_{A+2}F_H
((a),-1/2-n,-n;(h);y)
=\frac{1}{2x(1-2d)}
\\ \times
\{
(1+x)^{1-2d}{}_{A+2}F_H((a),d,d-1/2;(h);\frac{x^2y}{(1+x)^2})
-(1-x)^{1-2d}{}_{A+2}F_H((a),d,d-1/2;(h);\frac{x^2y}{(1-x)^2})
\}
.
\end{multline}

\cite{Scarpelloarxiv1212}
\be
F_1\left(
\begin{array}{c} a;b,b\\
1+a-b\end{array}
\mid e^{2\pi i/3},e^{-2\pi i/3}
\right)
=\frac{\Gamma(1+a-b)\Gamma(1+a/3)}{\Gamma(1+a)\Gamma(1+a/3-b)}
.
\ee

\cite{Scarpelloarxiv1212}
\be
F_D\left(
\begin{array}{c} a;b,b,\ldots b\\
c\end{array}
\mid \omega_{1,n},\ldots \omega_{n-1,n}
\right)
=\frac{\Gamma(a-b+1)\Gamma(1+a/n)}{\Gamma(1+a)\Gamma(1+a/n-b)}
,
\ee
where $\omega_{k,n}\equiv e^{2k\pi i/n}$.

\cite{Scarpelloarxiv1212}
\be
F_D\left(
\begin{array}{c} 2mb-a;b,b,\ldots b\\
2mb \end{array}
\mid x_1,\ldots x_{2m}
\right)
=\frac{1}{2m}\, \frac{\Gamma(a/(2m))\Gamma(2mb)\Gamma((2mb-a)/(2m))}{\Gamma(a)\Gamma(b)\Gamma(2mb-a)}
,
\ee
where $x_k=1+e^{(2k-1)\pi i/(2m)}$, for $k=1,\ldots 2m$, and $2m$ is an even integer, $a>0$, $b>0$, $nb>a$.

\cite{Scarpelloarxiv1212}
\begin{multline}
F_D\left(
\begin{array}{c} (2m-1)b-a;b,b,\ldots b\\
(2m-1)b \end{array}
\mid y_1,\ldots y_{2m-1}
\right)
\\
=\frac{1}{2m-1}\, \frac{\Gamma(a/(2m-1))\Gamma((2m-1)b)\Gamma(((2m-1)b-a)/(2m-1))}{\Gamma(a)\Gamma(b)\Gamma((2m-1)b-a)}
,
\end{multline}
where $y_k=1+e^{(2k-1)\pi i/(2m-1)}$, and $2m-1$ is an odd integer, $a>0$, $b>0$, $nb>a$.

\cite{VidunasJMAA355}
\be
x(1-x)\frac{\partial^2 F_2}{\partial x^2}
-xy\frac{\partial^2 F_2}{\partial x\partial y}
+(c_1-(a+b_1+1)x)\frac{\partial F_2}{\partial x}
-b_1y\frac{\partial F_2}{\partial y}
-ab_1 F_2=0.
\ee
\be
y(1-y)\frac{\partial^2 F_2}{\partial y^2}
-xy\frac{\partial^2 F_2}{\partial x\partial y}
+(c_2-(a+b_2+1)y)\frac{\partial F_2}{\partial y}
-b_2x\frac{\partial F_2}{\partial x}
-ab_2 F_2=0.
\ee

\cite{VidunasJMAA355}
\be
x(1-x)\frac{\partial^2 F_3}{\partial x^2}
+y\frac{\partial^2 F_3}{\partial x\partial y}
+(c-(a_1+b_1+1)x)\frac{\partial F_3}{\partial x}
-a_1b_1 F_3=0.
\ee
\be
y(1-y)\frac{\partial^2 F_3}{\partial y^2}
+x\frac{\partial^2 F_3}{\partial x\partial y}
+(c_2-(a_2+b_2+1)y)\frac{\partial F_3}{\partial x}
-a_2b_2 F_3=0.
\ee

\cite{VidunasJMAA355}
\be
x(1-x)\frac{\partial^2 F_4}{\partial x^2}
-y^2\frac{\partial^2 F_4}{\partial y^2}
-2xy\frac{\partial^2 F_4}{\partial x \partial y}
+(c_1-(a+b+1)x)\frac{\partial F_4}{\partial x}
-(a+b+1)y\frac{\partial F_4}{\partial y}
-ab F_4=0.
\ee
\be
y(1-y)\frac{\partial^2 F_4}{\partial y^2}
-x^2\frac{\partial^2 F_4}{\partial x^2}
-2xy\frac{\partial^2 F_4}{\partial x \partial y}
+(c_2-(a+b+1)y)\frac{\partial F_4}{\partial y}
-(a+b+1)x\frac{\partial F_4}{\partial x}
-ab F_4=0.
\ee

\cite[p 20]{Appell}
\begin{multline}
\alpha F_4(\alpha+1,\beta;\gamma,\gamma';x,y)
-\beta F_4(\alpha,\beta+1;\gamma,\gamma';x,y)
\\
=
(\alpha-\beta)
F_4(\alpha,\beta;\gamma,\gamma';x,y)
.
\end{multline}

\cite[p 21]{Appell}
\begin{multline}
\frac{\beta}{\gamma}xF_4(\alpha+1,\beta+1;\gamma+1,\gamma';x,y)
+\frac{\beta}{\gamma'}yF_4(\alpha+1,\beta+1;\gamma,\gamma'+1;x,y)
\\
=
F_4(\alpha+1,\beta;\gamma,\gamma';x,y)
-F_4(\alpha,\beta;\gamma,\gamma';x,y)
.
\end{multline}

\cite[p 26]{Appell}
\begin{multline}
F_4(\alpha,\beta;\gamma,\gamma';x,y)
=
\sum\frac{(\alpha)_m(\beta)_m}{(\gamma)_m m!}
\frac{\Gamma(\gamma')\Gamma(\beta-\alpha)}{\Gamma(\gamma'-\alpha-m)\Gamma(\beta+m)}
(-y)^{-\alpha-m}
\\
\times \,_2F_1(\alpha+m,\alpha+m+1-\gamma';\alpha+1-\beta;\frac{1}{y})x^m
\\
+
\sum\frac{(\alpha)_m(\beta)_m}{(\gamma)_m m!}
\frac{\Gamma(\gamma')\Gamma(\alpha-\beta)}{\Gamma(\gamma'-\beta-m)\Gamma(\alpha+m)}
(-y)^{-\beta-m}\,
\\
\times
\,_2F_1(\beta+m,\beta+m+1-\gamma';\beta+1-\alpha;\frac{1}{y})x^m
.
\end{multline}

\cite{Appell}
\begin{multline}
F_4(\alpha,\beta;\gamma,\gamma';x,y)
=
\frac{\Gamma(\gamma')\Gamma(\beta-\alpha)}{\Gamma(\gamma'-\alpha)\Gamma(\beta)}
(-y)^{-\alpha}
 F_4(\alpha,\alpha+1-\gamma';\gamma,\alpha+1-\beta;\frac{x}{y},\frac{1}{y})
\\
+
\frac{\Gamma(\gamma')\Gamma(\alpha-\beta)}{\Gamma(\gamma'-\beta)\Gamma(\alpha)}
(-y)^{-\beta}
 F_4(\beta+1-\gamma',\beta;\gamma,\beta+1-\alpha;\frac{x}{y},\frac{1}{y})
.
\end{multline}

\cite{BaileyPLMS40,GasperPJM56,VidunasJMAA355}
\be
F_4\left[\alpha,\beta;\gamma,1+\alpha+\beta-\gamma;x(1-y),y(1-x)\right]
=
F(\alpha,\beta;\gamma;x)F(\alpha,\beta;1+\alpha+\beta-\gamma;y)
.
\ee
\be
F_4\left[\alpha,\beta;\gamma,\beta;-\frac{x}{(1-x)(1-y)},-\frac{y}{(1-x)(1-y)}\right]
=
(1-x)^\alpha (1-y)^\alpha F_1[\alpha;\gamma-\beta,1+\alpha-\gamma;\gamma;x,xy]
.
\ee
\be
F_4\left[\alpha,\beta;\alpha,\beta;-\frac{x}{(1-x)(1-y)},-\frac{y}{(1-x)(1-y)}\right]
=
(1-xy)^{-1} (1-x)^\beta (1-y)^\alpha.
\ee
\be
F_4\left[\alpha,\beta;\beta,\beta;-\frac{x}{(1-x)(1-y)},-\frac{y}{(1-x)(1-y)}\right]
=
(1-x)^\alpha (1-y)^\alpha F[\alpha,1+\alpha-\beta;\beta;xy]
.
\ee
\be
F_4\left[\alpha,\beta;1+\alpha-\beta,\beta;-\frac{x}{(1-x)(1-y)},-\frac{y}{(1-x)(1-y)}\right]
=
(1-y)^\alpha F[\alpha,,\beta;1+\alpha-\beta;-\frac{x(1-y)}{1-x}]
.
\ee

\cite{KarlssonPA87}
\be
F^{p:2;2}_{q:1;1}
\left[
\begin{array}{ccc}
\alpha_1,\ldots \alpha_p: &
a_1,b_1 &
a_2,b_2 \\
\gamma_1,\ldots \gamma_q: &
c_1 &
c_2 \\
\end{array}\mid X,X
\right]
=
_{p+3}F_{q+2}
(\alpha_1,\ldots,\alpha_p,\beta_1,\beta_2,\beta_3;
\gamma_1,\ldots,\gamma_q,\delta_1,\delta_2;X
)
\ee

\begin{tabular}{ccc|ccc|ccc|cc}
$a_1$ & $b_1$ & $c_1$ & $a_2$ & $b_2$ & $c_2$ &
$\beta_1$ & $\beta_2$ & $\beta_3$ & $\delta_1$ & $\delta_2$ \\
\hline
$a$ & $b$ & $a+b-\frac{1}{2}$ & $a$ & $b$ & $a+b+\frac{1}{2}$ &
$2a$ & $2b$ & $a+b$ & $2a+2b-1$ & $a+b+\frac{1}{2}$ \\
$a-\frac{1}{2}$ & $b-\frac{1}{2}$ & $a+b-\frac{1}{2}$ & $a+\frac{1}{2}$ & $b+\frac{1}{2}$ & $a+b+\frac{1}{2}$ &
$2a$ & $2b$ & $a+b$ & $2a+2b-1$ & $a+b+\frac{1}{2}$ \\
\end{tabular}
and more of that format.

\cite{SrivastavaJMAA107}
\be
F^{1:2;2}_{1:1;1}
\left[
\begin{array}{ccc}
\alpha: & -M,\beta: & -N,\beta': \\
\beta+\beta': & \alpha-\beta'-M+1 & \alpha-\beta-N+1 \\
\end{array}\mid 1,1
\right]
=
\frac{(\beta+\beta'-\alpha)_{M+N}(\beta')_M(\beta)_N}
{(\beta+\beta')_{M+N}(\beta'-\alpha)_M(\beta-\alpha)_N}
\ee
for $M,N=0,1,2,\ldots$.

\cite{SrivastavaJMAA107}
\begin{multline}
F^{1:2;2}_{1:1;1}
\left[
\begin{array}{ccc}
\beta+\beta'-\alpha: & \beta,\gamma: & \beta',\gamma': \\
\beta+\beta': & \delta & \delta' \\
\end{array}\mid 1,1
\right]
\\
=
\frac{\Gamma(\delta)\Gamma(\delta')\Gamma(\alpha-\beta-\gamma+\delta)\Gamma(\alpha-\beta'-\gamma'+\delta')}
{\Gamma(\delta-\gamma)\Gamma(\delta'-\gamma')\Gamma(\alpha-\beta+\delta)\Gamma(\alpha-\beta'+\delta')}
F^{1:2;2}_{1:1;1}
\left[
\begin{array}{ccc}
\alpha: & \beta',\gamma: & \beta,\gamma': \\
\beta+\beta': & \alpha-\beta+\delta & \alpha-\beta'+\delta' \\
\end{array}\mid 1,1
\right]
.
\end{multline}

\cite{SrivastavaJMAA107}
\begin{multline}
F^{1:r;u}_{1:s;v}
\left[
\begin{array}{ccc}
\alpha: & (a_r): & (c_u): \\
\gamma : & (b_s): & (d_v); \\
\end{array}\mid x,y
\right]
=
\sum_{n=0}^\infty
\frac{(\alpha)_n(\gamma-\alpha)_n\prod_{j=1}^r(a_j)_n\prod_{j=1}^u(c_j)_n}
{(\gamma+n-1)_n(\gamma)_{2n}\prod_{j=1}^s(b_j)_n\prod_{l=1}^v(d_l)_n}
\,\frac{(xy)^n}{n!}
\\
\times
\,_{r+1}F_{s+1}\left(\begin{array}{cc}(a_r)+n, & \alpha+n \\
(b_s)+n, &\gamma+2n\end{array}
\mid x\right)
\,_{u+1}F_{v+1}\left(\begin{array}{cc}
(c_u)+n,& \alpha+n;\\
(d_v)+n,& \gamma+2n;\\
\end{array}\mid y\right)
.
\end{multline}

\cite{GasperPJM56}
\begin{multline}
_3F_2\left(\begin{array}{c}-n,n+a,b\\ c,d\end{array}\mid 1\right)
\,
_3F_2\left(\begin{array}{c}-n,n+a,e\\ c,f\end{array}\mid 1\right)
\\
=
\frac{(-)^n(a-c+1)_n}{(c)_n}
F^{2:2;2}_{2:1;1}\left[\begin{array}{ccc}-n,n+a:& b,e;&d-b,f-e\\d,f:&c;&a-c+1\end{array}\mid 1,1\right]
.
\end{multline}

\cite{GasperPJM56}
\begin{multline}
F^{2:2;2}_{2:1;1}\left[\begin{array}{ccc}a,b:& -x,y+e;&-y,x+d\\d,e:&c;&b\end{array}\mid 1,1\right]
\\
=
\frac{(d-a)_x(e-a)_y}{(d)_x(e)_y}
F^{2:2;1}_{2:1;0}\left[\begin{array}{ccc}a,-x:& 1+a-c,-y;&c-b\\
c,1+a-d-x:&1+a-e-y;&-\end{array}\mid 1,1\right]
\end{multline}
if $x,y=0,1,\ldots$

\cite{GasperPJM56}
\begin{multline}
F^{2:2;2}_{2:1;1}
\left[\begin{array}{ccc}a,b:& -x,y+e;&-y,x+d\\d,e:&c;&c'\end{array}\mid 1,1\right]
=
\sum_{r=0}^{\min(x,y)}
\frac{(a)_r(b)_r(a+b-c-c'+1)_r}{r!(c)_r(c')_r}
\,
\frac{(-x)_r(-y)_r}{(d)_r(e)_r}
\\
\times
\,
_3F_2\left(\begin{array}{c}r+a,r+b,r-x\\ r+c,r+c\end{array}\mid 1\right)
\,
_3F_2\left(\begin{array}{c}r+a,r+b,r-y\\ r+c',r+e\end{array}\mid 1\right)
\end{multline}
if $x,y=0,1,\ldots$

\subsection{The Confluent Hypergeometric Function}
\subsection{The Meijer G-Function}

\cite{AdamchikJCAM64}
\begin{multline}
G_{p+1,q+2}^{m+1,n+1}\left(z\mid\begin{array}{c}a,a_p\\a,b_q,b\end{array}\right)
=
(-)^{a-b}G_{p,q+1}^{m+1,n}
\left(z\mid\begin{array}{c}a_p\\b_q,b\end{array}\right)
\\
-
(-1)^{a-b}\sum_{k=1}^{a-b}\Res_{s=k-a}
\left[\frac{\Gamma(b+s)\prod_{i=1}^m\Gamma(b_i+s)\prod_{i=1}^n\Gamma(1-a_i-s)}
{\prod_{i=n+1}^p\Gamma(a_i+s)\prod_{i=m+1,q}\Gamma(1-b_i-s)}z^{-s}\right]
\end{multline}
if $a-b>0$.

\cite{AdamchikJCAM64}
\be
G_{p+1,q+2}^{m+1,n+1}\left(z\mid\begin{array}{c}a,a_p\\a,b_q,b\end{array}\right)
=
(-)^{a-b}G_{p,q+1}^{m+1,n}
\left(z\mid\begin{array}{c}a_p\\b_q,b\end{array}\right)
\ee
if $a-b\le 0$.

\cite{AdamchikJCAM64}
\be
G_{p+2,q+1}^{m,n+1}\left(z\mid\begin{array}{c}a,a_p,b\\b_q,b\end{array}\right)
=
(-)^{a-b}G_{p+1,q}^{m,n}
\left(z\mid\begin{array}{c}a_p,a\\b_q\end{array}\right)
\ee
if $a-b$ is an integer.

\cite{MilgramANE18}
\begin{multline}
G_{p+1,q+1}^{m,n+1}\left(z\mid\begin{array}{l} a,a_p \\ b_q,b \end{array}\right)
=
(-)^r G_{p+1,q+1}^{m+1,n}\left(z\mid\begin{array}{l} a_p,a \\ b,b_q \end{array}\right)
\\
-
\frac{(-)^r\prod_{j=1}^m\Gamma(b_j-b)\sum_{j=1}^n\Gamma(1+b-a_j)z^b}
{\prod_{j=m+1}^q \Gamma(1+b-b_j)\prod_{j=n+1}^p\Gamma(a_j-b)\Gamma(r)}
{}_{p+1}F_q\left(
\begin{array}{l}1+b-a_p,1-r\\ 1+b-b_q
\end{array}
\mid (-)^{p-m-n}z
\right)
\end{multline}
where $a-b=r$ is integer, $q\ge m$. This corrects \cite[(5.4.7)]{Luke1969Book}

\cite{MilgramANE18}
\begin{multline}
G_{p+1,q+1}^{m+1,k}\left(z\mid\begin{array}{l} a,a_p \\ a,b_q \end{array}\right)
=
(-)^r G_{p,q}^{m,k}\left(z\mid\begin{array}{l} a_p \\ b_q \end{array}\right)
\\
+
\frac{\prod_{j=1}^m\Gamma(b_j-a)\sum_{j=1}^n\Gamma(1+a-a_j)z^a}
{\prod_{j=m+1}^q \Gamma(1+a-b_j)\prod_{j=k+1}^p\Gamma(a_j-a)\Gamma(r)}
{}_{p+1}F_q\left(
\begin{array}{l}1+a-a_p,1\\ 1+a-b_q
\end{array}
\mid (-)^{p-m-k}z
\right)
\end{multline}
where $a_k=a+r$ is integer. This corrects \cite[(5.4.6)]{Luke1969Book}

\cite{MilgramANE18}
\begin{multline}
G_{2,4}^{4,1}\left(-1/t\mid\begin{array}{l} j,2 \\ 3/2,1,1,\epsilon; \end{array}\right)
=2\pi\Gamma(\epsilon-3/2)(-1/t)^{3/2}{}_1F_2(1/2;3/2,5/2-\epsilon;1/t)
\\
+\Gamma(3/2-\epsilon)\Gamma(\epsilon)\Gamma(1-\epsilon)/(1-\epsilon)(1-t)^\epsilon
{}_1F_2(\epsilon-1;\epsilon-1/2,\epsilon;1/t)
\\
-\surd \pi \Gamma(\epsilon-1)/t[\psi(1/2)+\psi(\epsilon-1)+\log(-t)]
\\
-2\surd \pi \Gamma(\epsilon-1)/(2-\epsilon)/t^2 {}_2F_3(1,1;3/2,3-\epsilon,2;1/t).
\\
\end{multline}

\subsection{The MacRobert E-function}
\subsection{Riemann and Hurwitz zeta functions}

\cite[(2.4.1)]{Titchmarsh}
\begin{equation}
\zeta(s)=\frac{1}{\Gamma(s)}\int_0^\infty \frac{x^{s-1}}{e^x-1}dx,\quad \sigma >1.
\end{equation}

\cite{SrivastavaZAA19,AdamchikAnal18,PilehroodDMTCS12,PilehroodEJC18}
\be
\zeta(3)=\frac{5}{2}\sum_{k\ge 1}\frac{(-)^{k-1}}{k^3{2k \choose k}}
=
-\frac{4\pi^2}{7}\sum_{k\ge 0}\frac{\zeta(2k)}{(2k+1)(2k+2)2^{2k}}.
\ee

\cite{PilehroodDMTCS12,PilehroodEJC18,LupasRoGer2000}
\be
\zeta(3) = \frac{1}{2}\sum_{n\ge 1}\frac{(-)^{n-1}(205n^2-160n+32)}{n^5\binom{2n}{n}^5}.
\ee

\cite{LupasRoGer2000}
\be
\zeta(3)= \frac{3752}{3125}+\frac{1}{625}\sum_{k=1}^\infty (-)^{k-1}\frac{P(k)}{k^5\binom{2k}{k}^5(2k+1)^5},
\ee
where $P(k)=14760k^4+28010k^3+19505k^2+5920k+671$.

\cite{LupasRoGer2000}
\be
\zeta(3,x+1)=\frac52 \sum_{k=0}^\infty \frac{(-)^{k-1}}{k^3\binom{2k+x}{k}\binom{k+x}{k}}(1-\frac{2x(3k+2x)}{5(k+x)^2}).
\ee

\cite{PilehroodDMTCS10}
\be
\zeta(4)=\sum_{n=1}^\infty \frac{1}{n!^4}\left(
\frac{4n+1}{2n^2}L_n+\frac{7n^3}{4}L_{n-1}
\right)
\ee
where $L_0=0$, $L_1=1/3$ and
\be
4(4n+3)(4n+5)L_{n+1}+2(n+1)^3(6n^3+9n^2+5n+1)L_n
-n^6(n+1)^3L_{n-1}=0.
\ee

\cite{PilehroodDMTCS12,PilehroodEJC18}
\be
\zeta(4)=\frac{36}{17}\sum_{k\ge 1} \frac{1}{k^4\binom{2k}{k}}.
\ee

\cite{PilehroodDMTCS12,PilehroodEJC18}
\be
\zeta(5)=2\sum_{k\ge 1} \frac{(-)^{k+1}}{k^5\binom{2k}{k}}
-\frac{5}{2}\sum_{k\ge 1}\frac{(-)^{k+1}}{k^3\binom{2k}{k}}\sum_{j=1}^{k-1}\frac{1}{j^2}.
\ee

\cite{PilehroodDMTCS12,PilehroodEJC18}
\be
\zeta(7)=\frac{5}{2}\sum_{k\ge 1} \frac{(-)^{k+1}}{k^7\binom{2k}{k}}
+\frac{25}{2}\sum_{k\ge 1}\frac{(-)^{k+1}}{k^3\binom{2k}{k}}\sum_{j=1}^{k-1}\frac{1}{j^4}.
\ee

\cite{PilehroodEJC18}
\be
\zeta(2n+3)=
\frac52 \sum_{k=1}^\infty \frac{(-)^{k+1}}{k^3\binom{2k}{k}}
(-)^ne_n^{(2)}(k)
+2 \sum_{j=1}^n \sum_{k=1}^\infty \frac{(-)^{k+1}}{k^{2j+3}\binom{2k}{k}}
(-)^{n-j}e_{n-j}^{(2)}(k)
\ee
where $e_r^{(s)}(k)\equiv \sum_{1\le j_1< j_2\cdots < j_1\le k-1}
(j_1j_2\cdots j_r)^{-s}$.

\cite{PilehroodEJC18}
\be
\zeta(4n+3)=
\frac52 \sum_{j=0}^n\sum_{k=1}^\infty \frac{(-)^{k+1}}{k^{3j+3}\binom{2k}{k}}
\sum_{r=0}^{n-j}4^r h_{n-j-r}^{(4)}(k)e_r^{(4)}(k)
\ee
where $h_r^{(s)}(k)\equiv [t^r] \prod_{j=1}^{k-1}[1-j^{-s}t)^{-1}$.

\cite{Silagarxiv1207}
\be
\zeta(2n+1) = \frac{(-)^n \pi^{2n+1}}{4[1-2^{-2n-1}](2n)!}
\int_0^1 \frac{E_{2n}(u)}{\sin\pi u} du.
\ee

\cite[A073002]{sloane}
\be
\zeta'(2)=-\frac{\pi^2}{6}[\gamma +\ln(2\pi)-12 \ln A]
\ee
where $A$ is the constant of \eqref{eq.A}.

\cite{MillerJCAM100}
\be
\zeta'(2k)=(-)^{k+1}\frac{(2\pi)^{2k}}{2(2k)!}
\{2k\zeta'(1-2k)-[\psi(2k)-\log(2\pi)]B_{2k}\}.
\ee

\cite{MillerJCAM100}
\begin{multline}
\zeta'(1-2k,p/q)=
\frac{[\psi(2k)-\log(2\pi q)]B_{2k}(p/q)}{2k}
-\frac{[\psi(2k)-\log(2\pi)]B_{2k}(p/q)}{q^{2k}2k}
\\
+
(-)^{k+1}\frac{\pi}{(2\pi q)^{2k}}
\sum_{n=1}^{q-1}\sin(\frac{2\pi p n}{q})\psi^{(2k-1)}(n/q)
\\
+
(-)^{k+1}\frac{2(2k-1)!}{(2\pi q)^{2k}}
\sum_{n=1}^{q-1}\cos(\frac{2\pi p n}{q})\zeta'(2k,n/q)
+\frac{\zeta'(1-2k)}{q^{2k}}.
\end{multline}

\cite{GeddesCM}
\be
\zeta(z,q) = \frac{v^x}{\Gamma(z)}\int_0^\infty \frac{x^{z-1}\exp(-qvx)}{1-exp(-vx)}dx.
\ee

\cite{MillerJCAM100}
\be
\zeta'(1-2k,1/2)=
-\frac{B_{2k}\log 2}{k4^k}
-
\frac{(2^{2k-1}-1)\zeta'(1-2k)}{2^{2k-1}}.
\ee

\cite{MillerJCAM100}
\begin{multline}
\zeta'(1-2k,1/3)=
-\frac{(9^k-1)B_{2k}\pi}{\sqrt{3}(3^{2k-1}-1)8k}
-\frac{B_{2k}\log 3}{3^{2k-1}4k}
\\
+
(-)^k\frac{\psi^{(2k-1)}(1/3)}{2\sqrt{3}(6\pi)^{2k-1}}
-
\frac{(3^{2k-1}-1)\zeta'(1-2k)}{2\times 3^{2k-1}}.
\end{multline}

\cite{PilehroodDMTCS10}
\begin{multline}
\sum_{k=1}^\infty \frac{k}{k^4-x^2k^2-y^4}
=
\sum_{n,m=0}^\infty \binom{n+m}{n}\zeta(2n+4m+3)x^{2n}
y^{4m}
\\
=
\frac{1}{2}\sum_{n=1}^\infty \frac{(-)^{n-1}r(n)}{n\binom{2n}{n}}
\frac{\prod_{m=1}^{n-1}( (m^2-x^2)^2+4y^4}
{\prod_{m=n}^{2n}(m^4-x^2m^2-y^4)}
\end{multline}
where
\be
r(n)=205n^6-160n^5+(32-62x^2)n^4+40x^2n^3+(x^4-8x^2-25y^4)n^2+10y^4n+y^4(x^2-2).
\ee

\cite{AdamchikAnal18}
\be
\sum_{k\ge 1}\frac{(-)^k}{k}[\zeta(nk)-1] = \log\left(
\prod_{j=0}^{n-1}\Gamma[2-(-)^{(2j+1)/n}]
\right)
.
\ee

\cite{AdamchikAnal18}
\be
\sum_{k=2}^\infty (-)^k[\zeta(k)-1]k^n = -1+\frac{1-2^{n+1}}{n+1}B_{n+1}
-\sum_{k=1}^n(-)^k k!\zeta(k+1)S(n+1,k+1),
\ee
with $S$ the Stirling numbers of the second kind.

\cite{AdamchikAnal18}
\be
\sum_{k=2}^\infty [\zeta(k)-1]k^n = 1
+\sum_{k=1}^n k!\zeta(k+1)S(n+1,k+1),
\ee

\cite{Adamchik}
\be
-\sum_{k=1}^\infty\frac{k}{4^{2k}}\zeta(2k+1)=G-1.
\ee

\cite{SrivastavaZAA19}
\be
\sum_{k\ge 0}\frac{\zeta(2k)}{(2k+1)4^{2k}} = -\frac{G}{\pi} -\frac{1}{4}\log 2.
\ee

\cite{SrivastavaZAA19}
\be
\sum_{k\ge 0}\frac{\zeta(2k)}{2k+1}\left(\frac34\right)^{2k} = \frac{G}{3\pi} -\frac{1}{4}\log 2.
\ee

\cite{Adamchik}
\be
\frac{1}{16}\sum_{k=1}^\infty\frac{3^k-1}{4^k}(k+1)\zeta(k+2)=G.
\ee

\cite{AdamchikAnal18}
\be
\sum_{k=1}^\infty [\zeta(2k+1)-1]=\frac14.
\ee

\cite{AdamchikAnal18}
\be
\sum_{k=1}^\infty \frac{\zeta(2k+1)-1}{k+2}=\frac{2}{3}-\frac{\gamma}{2}+\log 2+6\zeta'(-1).
\ee

\cite{AdamchikAnal18}
\be
\sum_{k=2}^\infty \frac{\zeta(k)}{(k+1)(k+2)}=\frac{1-\gamma}{6}-2\zeta'(-1).
\ee

\cite{AdamchikAnal18}
\be
\sum_{k=2}^\infty \frac{k^2}{k+1}[\zeta(k)-1]=\frac{3}{2}-\frac{\gamma}{2}+\frac{\pi^2}{6}-\frac{1}{2}\log(2\pi).
\ee

\cite{AdamchikAnal18}\cite[A256919]{sloane}
\be
\sum_{k=2}^\infty [\zeta(4k)-1]=\frac{7}{8}-\frac{\pi}{4}\coth \pi
\ee

\cite{Jasinskydec20}\cite[A339529]{sloane}
\be
\sum_{k=1}^\infty [\zeta(6k)-1]=\frac{11}{12}-\frac{\pi}{2\sqrt 3}\tanh \frac{\sqrt{3}\pi}{2}.
\ee

\cite{Jasinskydec20}
\be
\sum_{k=0}^\infty [\zeta(6k+2)-1]=-\frac{1}{4}+\frac{\pi}{2\sqrt 3}\tanh \frac{\sqrt{3}\pi}{2}.
\ee

\cite{Jasinskydec20}
\be
\sum_{k=0}^\infty [\zeta(6k+4)-1]=\frac{1}{12}.
\ee

\cite{MillerJCAM100}
\be
4\sum_{k=0}^\infty \frac{1-\zeta(2k)}{(2k+1)3^{2k}}=\log(192) -
\pi\frac{2\surd 3}{9}+\frac{\surd 3}{3\pi}\psi^{(1)}(1/3).
\ee

\cite{AdamchikAnal18}
\be
\sum_{k=1}^\infty [\zeta(4k)-1]z^{4k} = \frac{3z^4-1}{2(z^4-1)}
-\frac{\pi z}{4}[\cot(\pi z)+\coth(\pi z)],\quad |z|<2.
\ee

\cite{AdamchikAnal18}
\be
\sum_{k=1}^\infty [\zeta(2k)-1]\sin k = -\frac{1}{2}\cot(1/2)
+\frac{\pi}{2}\frac{\sin(1/2)\sin[2\pi \cos (1/2)]-\cos(1/2)\sinh[2\pi\sin(1/2)]}{\cos[2\pi\cos(1/2)]-\cosh[2\pi\sin(1/2)]}
\ee

\cite{AdamchikAnal18}
\be
\sum_{k=1}^\infty {p+k \choose k}\zeta(p+k+1,a)z^k=
\frac{(-)^p}{p!}[\psi^{(p)}(a)-\psi^{(p)}(a-z)]
.
\ee

\cite{AdamchikAnal18}
\be
\sum_{k=1}^\infty \frac{t^k}{k^2}\zeta(2k)=\log[\pi\surd t \csc(\pi \surd t].
\ee

\cite{ChoiJMAA231}
\be
\sum_{k\ge 1}\frac{\zeta(2k+1)}{(k+1)2^{2k}} = -2-\frac13\log 2 -\gamma +12 \log A,
\ee
where $A$ is Glaisher's constant.

\cite{ChoiJMAA231}
\be
\sum_{k\ge 2}(-)^k\frac{\zeta(k)-1}{(k+1)2^k} = \frac74
+\log\frac{2^{19/12}}{3^2} +\frac{\gamma }{4} -3 \log A.
\ee

\cite{ChoiJMAA231}
\be
\sum_{k\ge 2}(-)^k\frac{\zeta(k)-1}{k+1}(\frac32)^{k+1} = \frac{19}{8}
+\log\frac{2^{-17/24}}{5} +\frac{9\gamma }{8} -\frac32 \log A.
\ee

\cite{ChoiJMAA231}
\be
\sum_{k\ge 2}(-)^k\frac{\zeta(k)}{2^{2k}(k+1)} =
1+\frac{\gamma}{8}-\frac{G}{\pi}-\frac12\log\frac{2\pi A^9}{\Gamma^2(\frac14)}
\ee
and others of that type.

\cite{PilehroodDMTCS12}
\be
\sum_{n=0}^\infty \zeta(2n+2)a^{2n}=\sum_{k\ge 1}\frac{1}{k^2-a^2}=
3\sum_{k\ge 1}\frac{1}{\binom{2k}{k}(k^2-a^2)}\prod_{m=1}^{k-1}\frac{m^2-4a^2}{m^2-a^2}.
\ee

\cite{PilehroodEJC18}
\be
\sum_{n=0}^\infty \zeta(2n+2,\alpha)a^{2n}=
\sum_{n=0}^\infty \frac{(3n+3\alpha-2)(1+2a)_{n-1}(1-2a)_{n-1}}{
n\binom{2n}{n}(\alpha_a)_n(\alpha-a)_n}.
\ee
\cite{PilehroodEJC18}
\be
\sum_{n=0}^\infty \zeta(2n+2,\alpha)a^{2n}=
\sum_{n=1}^\infty \frac{p(n)}{n\binom{2n}{n}} \frac{\prod_{m=1}^{n-1}(m^2-4a^2)}
{\prod_{m=n-1}^{2n-1}((m+\alpha)^2-a^2)}
\ee
where $p(n)=2(2n-1)((2n-1+\alpha)^2-a^2+(5n+2\alpha-2)((n+\alpha-1)^2-a^2)$.

\cite{PilehroodEJC18}
\begin{multline}
\sum_{k=1}^\infty k \zeta(2k+1)a^{2k-2}=
\frac12
\sum_{n=1}^\infty \frac{(-)^{n-1}p(n)}{n\binom{2n}{n}^5
(n^2-a^2)^2(n^2-a^2/4)^2
} \prod_{m=1}^{n-1}\frac{1-4a^2/m^2}{(1-a^2/(m+n)^2)^2}
\end{multline}
where $p(n)=2(2n-1)(3n-1)(4n^2-a^2)^2+(n^2-a^2)^2(13n^2-2a^2)$.

\cite{PilehroodDMTCS12,PilehroodEJC18}
\be
\sum_{n=0}^\infty \zeta(2n+3)a^{2n}=\sum_{k\ge 1}\frac{1}{k(k^2-a^2)}=
\frac{1}{2}\sum_{k\ge 1}\frac{(-)^{k-1}}{k^3\binom{2k}{k}}\frac{5k^2-a^2}{k^2-a^2}\prod_{m=1}^{k-1}(1-a^2/m^2).
\ee

\cite{ReynoldsSciA37}
\be
\zeta(n,u)=(-)^n \frac{\psi^{(n-1)}(u)}{(n-1)!}.
\ee

\cite{Adamchik}
\be
\frac{1}{8}\sum_{k=2}^\infty\frac{k}{2^k}\zeta(k+1,3/4)=G.
\ee

\cite{Adamchik}
\be
-\frac{1}{8}\sum_{k=2}^\infty\frac{k}{2^k}\zeta(k+1,5/4)=G-1.
\ee

\cite{PilehroodEJC18}
\be
\sum_{n\ge 0}a^{2n}\zeta(2n+2,\alpha)=\sum_{n\ge 0}\frac{(3n+2\alpha-2)(1+2a)_{n-1}(1-2a)_{n-1}}{n\binom{2n}{n}(\alpha+a)_n(\alpha-a)_n}.
\ee

\cite{PilehroodEJC18}
\begin{multline}
\sum_{n\ge 0}(n+1)a^{2n}\zeta(2n+3,\alpha)=
\frac12 \sum_{n\ge 1}\frac{(-)^{n-2}(5n^2+6n(\alpha-2)+2(\alpha-1)^2-2a^2)}{n\binom{2n}{n}}
\\ \times
\frac{(n-1)!^2(1+2a)_{n-1}(1-2a)_{n-1}}{(\alpha+a)^2_n(\alpha-a)^2_n}
\end{multline}

\cite{PilehroodEJC18}
\begin{multline}
\sum_{k\ge 1}k\zeta(2k+1,\alpha)a^{2k-2}
=\frac12 \sum_{n\ge 1}\frac{(-)^{n-1}p_\alpha(n)}{n\binom{2n}{n}^5((n+\alpha-1)^2-a^2)^2((n+\alpha-1)^2-a^2/4)^2}
\\ \times
\prod_{m=1}^{n-1}\frac{1-4a^2/m^2}{((1+\frac{\alpha-1}{m+n})^2-\frac{a^2}{(m+n)^2})^2},
\end{multline}
where $p$ are polynomials in $n$, $a$ and $\alpha$.

\cite{ChoiJMAA231}
\be
\sum_{k\ge 2} (-1)^k\frac{\zeta(k)}{k(k+1)}z^{k+1}
=[\log(2\pi)-1]\frac{z}{2}+(\gamma-1)\frac{z^2}{2}
+z\log\Gamma(1+z)-\log G(1+z),
\ee
where $G$ is the reciprocal of the double Gamma function.

\cite{ChoiJMAA231}
\begin{multline}
\sum_{k=2}^\infty (-)^k\zeta(k,a) \frac{z^{k+1}}{k+1}
=
[a-\frac12-\frac12\log(2\pi)]z
+\frac12[1-\psi(a)]z^2\\
+(1-a)\log\Gamma(z+a)+\log G(z+a)
+(a-1)\log\Gamma(z)-\log G(a),\quad |z|<|a|,
\end{multline}
where $G$ is the reciprocal of the double Gamma function.

\cite{ChoiJMAA231}
\begin{multline}
\sum_{k=2}^\infty \zeta(k,a) \frac{z^{k+1}}{k+1}
=
[a-\frac12-\frac12\log(2\pi)]z
-\frac12[1-\psi(a)]z^2\\
+(a-1)\log\Gamma(z-a)-\log G(a-z)
+(1-a)\log\Gamma(z)+\log G(a),\quad |z|<|a|,
\end{multline}

\cite{AdamchikICSAC}
\be
\zeta'(1,p/q)-\zeta'(1,1-p/q)=
\pi \cot\frac{\pi p}{q}[\log(2\pi q)+\gamma]
-2\pi\sum_{j=1}^{q-1}\log\left(\Gamma(j/q)\right)\sin\frac{2\pi jp}{q}.
\ee

\cite{AdamchikICSAC}
\be
\zeta'(-n,x)+(-)^n\zeta'(-n,1-x)=
\pi i\frac{B_{n+1}(x)}{n+1}
+e^{-\pi i n/2} \frac{n!}{(2\pi)^n} \Li_{n+1}(e^{2\pi i x}).
\ee
where $n$ is a positive integer and $0<x<1$.

\cite{EspinosaRJ6}
\be
\int_0^1 \sin(2\pi q)\zeta(z,q)dq
=
\frac{(2\pi)^z}{4\Gamma(z)}
\csc\frac{z\pi}{2}.
\ee

\cite{EspinosaRJ6}
\be
\int_0^1 \sin(2k\pi q)\zeta(z,q)dq
=
\frac{(2\pi)^z k^{z-1}}{4\Gamma(z)}
\csc\frac{z\pi}{2}.
\ee

\cite{EspinosaRJ6}
\be
\int_0^1 \cos(2k\pi q)\zeta(z,q)dq
=
\frac{(2\pi)^z k^{z-1}}{4\Gamma(z)}
\sec\frac{z\pi}{2}.
\ee

\cite{EspinosaRJ6}
\be
\int_0^1 \zeta(z',q) \zeta(z,q)dq
=
-\zeta(z+z'-1)B(1-z,1-z')
\frac{\cos\frac{\pi(z-z')}{2}}{\cos\frac{\pi(z+z')}{2}}
.
\ee

\cite{EspinosaRJ6}
\be
\int_0^1 q^n \zeta(z,q)dq
=
-n!\sum_{j=1}^n\frac{\zeta(z-j)}{(z-j)_j(n-j+1)!}
.
\ee

\cite{EspinosaRJ6}
\be
\int_0^1 \ln(\sin \pi q) \zeta(z,q)dq
=
-\frac{\Gamma(1-z)}{(2\pi)^{1-z}}\sin\frac{\pi z}{2}\zeta(2-z)
.
\ee

\cite{GuilleraRJ16}
\be
\Phi(z,s,u)=2^{-s}\left[
\Phi(z^2,s,u/2)+z\Phi(z^2,s,(u+1)/2)
\right].
\ee

\cite{CrandallEM3}
\be
\zeta(r,s)\equiv \sum_{m<n}\frac{1}{n^sm^r}.
\ee

\cite{CrandallEM3}
\be
\zeta(r,s)+\zeta(s,r)=\zeta(r)\zeta(s)-\zeta(r+s).
\ee

\cite{CrandallEM3}
\be
\zeta(r,s) = -\frac{1}{2}\zeta(r+s)+\sum_{j=1,j\mathrm{odd}}
(\binom{j-1}{s-1}+\binom{j-1}{r-1})\zeta(j)\zeta(r+s-j).
\ee
for $r$ even and $s$ odd.

\cite{CrandallEM3}
\be
\zeta(1,s) = \frac{s}{2}\zeta(s+1)-\frac{1}{2}\sum_{j=1}^{s-2}
\zeta(j+1)\zeta(s-j).
\ee

\cite{CrandallEM3}
\be
\zeta(r,s) = -\frac{1}{2}\zeta(r+s)+\sum_{k=3,k\mathrm{odd}}^\infty
\Phi_k\sum_{j=0,j\mathrm{even}}^{k-1}\binom{k}{j}\eta(r-j)\eta(s-k+j),
\ee
where
\be
\Phi_k\equiv -\frac{2}{\pi}\sum_{d=1,d\mathrm{odd}}^{k-2}(-1)^{(d-1)/2}\frac{\pi^d}{d!}\zeta(k-d+1),
\ee
and
\be
\eta(s)=(1-2^{1-s})\zeta(s).
\ee

\cite{PilehroodDMTCS12}
\be
\sum_{n\ge 0}\beta(2n+1)a^{2n}=\sum_{n\ge 0}\frac{\binom{2n}{n}}{4^{2n+1}(2n+1)}\frac{3(2n+1)^2+a^2}{(2n+1)^2-a^2}
\prod_{m=0}^{n-1}(1-a^2/(2m+1)^2).
\ee

\cite{PilehroodDMTCS12}
\be
\sum_{n\ge 0}\beta(2n+2)a^{2n}=\frac{1}{8}\sum_{n\ge 1}\frac{16^n(12n^2-8n+1+a^2)}{\binom{4n}{2n}\binom{2n}{n}((2n-1)^2-a^2)}
\prod_{m=1}^{n-1}\left(1-\frac{a^2}{4m^2}\right).
\ee

\cite{PilehroodDMTCS12}
\be
\sum_{n\ge 0}\beta(2n+2)a^{2n}=\frac{1}{16}\sum_{n\ge 1}\frac{256^n p(n,a)}{n^3(2n-1)^3\binom{6n}{3n}\binom{6n}{4n}\binom{3n}{2n}
}
\frac{\prod_{m=1}^{n-1}(1-a^2/(4m^2))}{\prod_{m=n}^{3n}(1-a^2/(2m-1)^2)}
\ee
where
$$
p\equiv 580n^4-764n^3+344n^2-61n+15/4-a^2(18n^2-15n+4-a^2/4).
$$

\cite{PilehroodDMTCS12}
\be
\beta(2)=G=\frac{1}{64}\sum_{n\ge 1} \frac{(-)^{n-1}256^nq(n)}
{\binom{8n}{4n}^2\binom{2n}{n}n^3(2n-1)(4n-1)^2(4n-3)^2},
\ee
where $q(n) = 419840n^6-915456n^5+782848n^4-332800n^3+73256n^2-7800n+315$.

\cite{LupasRoGer2000}
\be
\beta(2)=G=\sum_{k\ge 1} (-)^{k-1} \frac{2^{8k}(40k^2-24k+3)}{k^3\binom{4k}{2k}^2\binom{2k}{k}64(2k-1)}.
\ee

\cite{CampbellIJMM16}
\be
\zeta(n)\prod_{k=2}^\infty\prod_{1\le j\le k,(j,k)=1}
\zeta(jm+kn)^{1/k} = \prod_p (\frac{1-p^{-(m+n)}}{1-p^{-n}})^{1/(1-p^{-m})}
\ee
if $\Re n, \Re m+n>1$,
which is a generalization of
\be
\sum_{k=1}^\infty \varphi(k)\frac{1}{k}\log \zeta(kn)=\sum_p \frac{1}{p^n-1},\, \Re n>1
\ee
where $\varphi$ is the Euler totient function, and $p$ are the primes.

\cite{CampbellIJMM16}
\be
\zeta(n)\zeta(m+n)\prod_{k=2}^\infty\prod_{1\le j\le k,(j,k)=1}
\zeta(jm+kn)^{1/k} = \prod_p (\frac{1}{1-p^{-n}})^{1/(1-p^{-m})}
\ee
if $\Re n,\Re m>1$.

\cite{CampbellIJMM16}
\be
\prod_{k=1}^\infty
\zeta(km)^{k\sum'(\alpha\beta)^{-1}} = \prod_p (\frac{1}{1-p^{-m}})^{1/(1-p^{-m})}
\ee
if $\Re m>1$, where $\sum'$ is over all positive integers $\alpha, \beta$
which are relatively prime and whose sum is $k$.

\subsection{Polylogarithms}
\cite{BarbieriNC11A}
\be
\Li_2(\frac{1}{y+i\epsilon})
=\Li_2(\frac{1}{y}-i\epsilon)
=-\Li_2(y)-\frac12 \log^2(-y-i\epsilon)-\zeta(2).
\ee
\be
\Re \Li_2(\frac{1}{y})
=-\Li_2(y)-\frac12 \log^2 y+2\zeta(2);\quad
\Im \Li_2(\frac{1}{y+i\epsilon})
=\pi\log y.
\ee

\cite{CoppoJNT150}
Let
\be
\alpha_k(s)=\frac{1}{\Gamma(s)}
\int_0^{\infty} \frac{e^{-t}}{1-e^{-2t}}\Li_k(1-e^{-2t})t^{s-1}dt,
\ee
\be
\beta(s)=\frac{1}{\Gamma(s)}
\int_0^{\infty} \frac{e^{-t}}{1-e^{-2t}}\Li_k(\frac{1-e^{-2t}}{2})t^{s-1}dt,
\ee
then
\be
\alpha_k(m+1)=\sum_{n=1}^\infty \frac{2^{2n-1}}{\binom{2n}{n}n^{k+1}}
P_m(
O_n^{(1)},\ldots
O_n^{(j)},\ldots
O_n^{(m)}
\ee
for $k\ge 1$ and
\be
\beta_k(m+1)=\sum_{n=1}^\infty \frac{2^{n-1}}{\binom{2n}{n}n^{k+1}}
P_m(
O_n^{(1)},\ldots
O_n^{(j)},\ldots
O_n^{(m)}
\label{eq.Bell}
\ee
for $k\ge 0$, where $O_n\equiv O_n{(1)}=\sum_{j=1}^n 1/(2j-1)$, 
$O_n^{(s)}=\sum_{k=1}^n 1/(2k-1)^s$ and
the modified Bell polynomials are
\be
P_m(x_1,\ldots x_m)= \sum_{k_1+2k_2+\cdots mk_m=m}
\frac{1}{k_1!k_2!\cdots k_m!} 
(x_1/1)^{k_1}
(x_2/2)^{k_2}
\cdots
(x_m/m)^{k_m}.
\ee

\cite{MezoRJ46}
\be
\int_0^1\Li_2(1-\sec^2 (\pi x))dx = -\pi^2/3.
\ee

\cite{MezoRJ46}
\be
\int_0^1\Li_2(\sec(\pi x))dx = i\frac{\pi}{2}\log 2.
\ee

\cite{MezoRJ46}
\be
\int_0^1\Li_3(\sin^2(\pi x))dx = \frac43 \log^3 2 -\frac13 \pi^2\log 2+2\zeta(3).
\ee

\subsection{Bernoulli Numbers and Polynomials}

\cite{GuilleraRJ16}
\be
B_m(x)=\sum_{n=0}^m\frac{1}{n+1}\sum_{k=0}^n(-1)^k\binom{n}{k}(x+k)^m.
\ee

\cite{Qiarxiv1401}
\be
B_n(x)=\sum_{i=0}^n (-1)^i \frac{\binom{n+1}{i+1}}{\binom{n+i}{i}}\left\{\begin{array}{c}n+i \\ i\end{array}\right\}.
\ee

\cite{NetoJIS19}
\be
B_n(m)=\sum_{l=0}^m (-1)^l \frac{l!}{l+1} \left\{\begin{array}{c}n+m \\ l+m\end{array}\right\}.
\ee

\cite{ChenJNT124}
\begin{multline}
\sum_{k=0}^m
\frac{\binom{m}{k}}{\binom{n+k+s}{s}}
y^{m-k}B_{n+k+s}(x)
=
\sum_{k=0}^n
\frac{\binom{n}{k}}{\binom{m+k+s}{s}}
(-y)^{n-k}
B_{m+k+s}(x+y)
\\
+\sum_{j=0}^{s-1}\sum_{i=0}^{s-1-j}
\binom{s-1-j}{i}\binom{s-1}{j}
\frac{(-1)^{n+1+i}sy^{m+n+s-j}B_j(x)}{(m+n+1+i)\binom{m+n+l}{n}}
\end{multline}
and
\be
\sum_{k=0}^m
\binom{m}{k}\binom{n+k}{s}
y^{m-k}B_{n+k-s}(x)
=
\sum_{k=0}^n
\binom{n}{k}\binom{m+k}{s}
(-y)^{n-k}
B_{m+k-s}(x+y)
.
\ee

\cite{ChenJNT124,SinghJIS19}
\be
\sum_{j=0}^k \binom{k+1}{j}(k+j+1)B_{k+j}=0,\quad k\ge 1.
\ee

\cite{CrandallEM3}
\be
B_s(t)B_r(t)=
\sum_{j>0,j\equiv r+s\bmod 2} \frac{1}{j}
\left(r\binom{s}{j-r}+s\binom{r}{j-s}\right) B_{r+s-j}B_j(t)+\frac{1}{2}
((-1)^{r+s}+1)
\frac{(-)^rB_{r+s}}{\binom{r+s}{s}}
\ee

\cite{ZhaoDM281}
\be
\sum_{k=j}^{n+1}\frac{(-)^k}{k}
\begin{Bmatrix}{c}n \\ k-1\end{Bmatrix}
\begin{bmatrix}{c}k \\ j\end{bmatrix}
=\frac{(-)^j}{n+1}\binom{n+1}{j}B_{n-j+1}.
\ee

\cite{GuilleraRJ16}
\be
E_m(x)=\sum_{n=0}^m\frac{1}{2^n}\sum_{k=0}^n(-1)^k\binom{n}{k}(x+k)^m.
\ee

\bibliographystyle{amsplain}
\bibliography{all}

\end{document}